\documentclass[10pt]{amsart}
\usepackage{latexsym}   % use all LaTeX fonts
\usepackage{amssymb}    % use all AMS fonts
\usepackage{amsmath}    % include some AMS-LaTeX functionality
\usepackage{amsthm}     % AMS style theorem and proof environments
\usepackage[all]{xy}
\usepackage{hyperref}
\newtheorem{definition}{Definition}

\newtheorem{example}{Example}

\newtheorem{remark}{Remark}

\newcommand{\pa}{\partial}\newcommand{\al}{\alpha}
\newcommand{\be}{\beta}

\newcommand{\Ga}{\Gamma}\newcommand{\del}{\delta}
\newcommand{\Del}{\Delta}\newcommand{\ep}{\epsilon}
\newcommand{\vaep}{\varepsilon}
\newcommand{\la}{\lambda}\newcommand{\om}{\omega}
\newcommand{\Om}{\Omega}\newcommand{\na}{\nabla}
\newcommand{\si}{\sigma}\newcommand{\ti}{\tilde}

\renewcommand{\thefootnote}

\title[   {\it The Method} of Archimedes in the Geometry of Quadrics]
{   {\it The Method} of Archimedes in the Geometry of Quadrics}

\author[  Ion I. Dinc\u{a}]{  Ion I. Dinc\u{a}}
\address{Slobozia 2, {\href{http://news.bbc.co.uk/2/hi/europe/680383.stm}
{\c{T}\u{a}nd\u{a}rei}}, Ialomi\c{t}a, 8454, Rom\^{a}nia}
\email{ionidinca@yahoo.com}

\thanks{ 
Supported by University of Notre Dame du Lac}

\begin{document}

\begin{abstract}
Confocal quadrics capture (encode) and geometrize spectral
properties of symmetric operators. Certain metric-projective
properties of confocal quadrics (most of them established in the
first half of the XIX$^{\mathrm{th}}$ century) {\it carry out}
(stick and transfer) by rolling to and influence surfaces {\it
applicable} (isometric) to quadrics and surfaces geometrically
linked to these, thus providing a wealth of integrable systems and
projective transformations of their solutions. We shall mainly
follow Bianchi's discussion of deformations (through bending) of
quadrics. Interestingly enough, {\it The Method} of Archimedes
(lost for $7$ centuries and rediscovered in the same year as
Bianchi's discovery (1906), so unknown to Bianchi) applies {\it
word by word in both spirit and the letter} and may provide the
key to generalizations in other settings. Basically we have

\begin{center}
{\it Wisdom of the wise ancients+flat connection form=integrable
systems.}
\end{center}
\begin{center}
{\it Wisdom of the wise ancients+algebraic computations of the
Bianchi Permutability Theorem=optimal discretization of integrable
systems.}
\end{center}

A four page appendix with the gist of the theory of deformations
of quadrics is included as

\S\ \ref{sec:gist} for the convenience of the rushed reader.

\end{abstract}
\maketitle \pagenumbering{roman}

$\ \ \ \ \ \ \ \ \ \ \ \ \ \ \ \ \ \ \ \ \ \ \ \ \ \ \ \ \ \ \ \ \
${\it 'It's only after we've lost everything that we're free to do
anything.'}

$\ \ \ \ \ \ \ \ \ \ \ \ \ \ \ \ \ \ \ \ \ \ \ \ \ \ \ \ \ \ \ \ \
\ \ \ \ \ \ \ \ \ \ \ \ \ \ \ \ \ \ \ \ \ \ \ \ \ \ \ \ \ \ \ \ \
\ \ \ \ \ \ \ \ \ \ \ \ \ \ \ \ \ \ \ \ \ \ \ \ $ Tyler Durden

\tableofcontents \pagenumbering{arabic}

\section{Introduction}\label{sec:introduction}\setcounter{equation}{0}

\subsection{Acknowledgements and full disclosure of the method of finding {\it The Method}}
\noindent

\noindent I am much indebted to the University of Notre Dame du
Lac for academic support during the graduate program (in
particular I am grateful that although my teaching duties began
with the second year and mainly required four tutorials per week,
most semesters I was required to hold only two; one semester even
only one) and to my Advisor Professor Brian Smyth for useful
discussions from Fall 2002 until Fall 2005 (he clocked two hours a
week during the academic year and mostly four a week during the
Summers; thus I had no choice but to present at the blackboard
mostly the opinions of the classical geometers). Also I would like
to thank my fellow graduate students (in particular to my
officemates Daniel Cibotaru, Florin Dumitrescu and Daniel Jackson)
for patiently listening and useful discussions.

I would also like to thank the Faculty of the Mathematics
Department from Bucharest University for the free high quality
undergraduate mathematical education (I was also granted a modest
financial support).

Special thanks to Daniel Cibotaru for an enlightening discussion
that took place in Fall 2005, when I was preparing a quiz for
Calculus B and I noticed a behavior which suggested a smooth ride
in a car with elliptic wheels and on a sinusoidal road (the
textbook was stating that the arc-length of
$(\cos(s),\sqrt{2}\sin(s))$ and that of $(s,\sin(s))$ are the
same; this provides rolling). I knew that the elliptic wheel has
to be attached to the car at one of the foci before even trying to
complete the computations, because replacing the circular wheel
with the elliptic wheel was the $1$-dimensional metric-projective
analogy to replacing the pseudo-sphere with an arbitrary quadric
for the $2$-dimensional theory of isometric deformations of
quadrics by means of rolling (we shall henceforth use the
classical denomination {\it deformation} instead of {\it isometric
deformation}). Thus I have stated my views to Daniel (clarifying
them as I was talking to him) and I felt compelled to do the
analytic confirmation and a Maple moving picture; the picture
suggested that the wheel rotates faster when the focus is closer
to the point of tangency, which led me to a book of history of
mathematics and the first encounter with the wisdom of the wise
ancients, who have mastered {\it 'the science of the sections
occurring in cones'}. It is this how I was able to give a possible
answer to the question which was bugging me for some time, namely
{\it Why one can do computations for quadrics?}; this was done
only after I read and understood the early history of conics (up
to including Newton), in the process simplifying the
$2$-dimensional discussion (I had all the facts, including the
metric-projective denomination, but I did not put them together up
to that point).

I thus kept in high esteem and added the early history of conics
(including Archimedes' quote from {\it The Method}) to the notes
in late December 2005-early January 2006; also I considered their
quotation in these notes to be in the correct place because I was
using a theorem of Menelaus in a place where Bianchi did not use
it. At that time I saw in those quotes a confirmation of the fact
that Bianchi was able to do computations for quadrics, so nothing
more than a weak analogy; I was looking at them as everybody does:
statements of historical importance. Unfortunately in what it
seems to be my biggest blunder so far I reduced the early history
of quadrics to less than a page by the end of January 2006 (thus
removing Archimedes' quote in the process) and I did not read the
early part of the notes up to August 2006 with the same frequency
and intensity I did before January 2006 (I used to do that on a
bad day, get my spirits up and do some work in the meat of the
notes).

In mid-August 2006 I drank a beer with Daniel and I stated my
hopes that with Levi-Civita's intuitive notion of parallel
transport taken as an obvious axiom and the methods developed by
Apollonius in {\it Conics} one can provide a synthetic proof of
the theory of deformations of quadrics, thus reducing it from a
proof Leibnitz would understand to a proof the wise ancients would
understand; again these things became more clear to me as I was
stating them to Daniel. After this discussion I noticed that at
some point I was stating in that version of these notes {\it 'All
necessary identities of the moving part boil down to valid
relevant algebraic identities of the static part; in turn these
algebraic identities naturally appear within the static part and
admit simple geometric interpretations.'} The {\it 'naturally'}
appeared from an analytic point of view (changing the ruling
family), but at some point my memory failed in my favor, as I did
not remember the change of ruling family, but just the reflection
property. But if the reflection property naturally appeared, why
would I continue to say {\it 'and admit simple geometric
interpretations'}, as the natural reflection property was itself
the simple geometric interpretation (and in trying to see why it
is natural I found the rolling of a surface on both sides of
another applicable surface also applies in the singular case of
the Bianchi-Lie ansatz: if one rolls the seed on both sides of the
considered quadric, then one gets immediately the reflection
property). Thus I tried to understand what I wrote; in doing so I
realized that the reflection property {\it as a natural geometric
approach} also works at the level of the Bianchi-Lie ansatz {\it
by means of an intuitive non-rigorous argument}.

So this was the setting in late August 2006 when I had a
discussion with the Chair (according to The Oracle in the {\it
'Matrix Trilogy':  'We cannot see beyond the choices we do not
understand'}); I believe I even mentioned to the Chair that again
I have ideas (between January 2006 and August 2006 I just run on
fumes, trying to complete a program clearly and previously
established); on that day the idea was only a fancier point of
view on the $2$-dimensional discussion, so it would not bring new
results in that setting. On the very same afternoon I went home
and I decided to bring the notes to the earlier January 2006
organization (I believe I also met Florin on campus on the way
home; he is interested in super parallel transport and I believe I
informed him on that occasion of Levi Civita's simple point of
view). I remembered Archimedes' quote and the fact that I liked
it, so first thing I did was to put it back in the notes. By
reasons of internal organization of these notes Archimedes' quote
stood for the first time alone, separate of its intended use in
the earlier version of the notes. I wanted to see if it is
correctly organized, so I read the statement above it, then
Archimedes' quote, but I never got to the statement below: as soon
as I read it something stuck to my mind and I had to investigate.
This is how I made the connection between Archimedes and the
classical geometers of the XIX$^{\mathrm{th}}$-XX$^{\mathrm{th}}$
century.

Now even if I saw this simplification (which may have helped me
look closer at Archimedes' quote) before making the connection
with Archimedes, Archimedes' point of view is still the one which
must have priority, as {\it the first one}; on top of that
Archimedes took it up a notch. I did not expect to get what I was
looking for before the first attempt to seriously look at the
ideas of the wise ancients as statements of current usefulness
instead of statements of historical usefulness. So I got the first
interpretation of Archimedes' quote by accident (it was separated
from its initial place in the earlier version of the notes because
I added the comment about Archimedes being generous {\it in
disclosing the full method of investigation}; had I've put it back
to its original place, then probably I would have not made the
connection, as I read and re-read this quote for one month in late
December 2005- early January 2006 thinking nothing else but the
situation at hand, namely Archimedes' specific applications of his
method). Further and maybe more irrationally than rationally I
almost saw the {\it actual} interpretation because {\it I just
obeyed the command shown in front of my eyes}; thus I further
investigated (as Jules would say not all signs are of a divine
nature, but by use of common sense, Fonzie cool and simple logic
the false ones can be easily removed as such) and it was clearly
revealed to me that the big fuss about what I saw and Bianchi did
not see, about what I did and Bianchi did not do would be decided
by giving the credit (as usually) to {\it a third neutral party}.
What hit me was the fact that Archimedes told me that there is no
spot of trouble in using the singular case {\it full Bianchi-Lie
ansatz} just as the non-singular case {\it Theorem I of Bianchi's
theory of deformations of quadrics}: thus basically Archimedes
solved the most difficult step of the theory of deformations of
quadrics 2200 years before anybody got any inkling about such
things. What also hit me was the power of generality of {\it The
Method} (the only way Archimedes could have accomplished this
deed), so I tried it immediately on the higher dimensional problem
and instantaneously I got {\it the first confirmation that the
higher dimensional problem is amenable to an attack strategy}:
this is my first important application of Archimedes' ideas.

Again {\it The Arnold Principle} worked out just fine (luckily
just on time for me; I felt uneasy about criticizing Bianchi and
praising myself instead; I believe even Bianchi would have
approved of Archimedes' opinion on the matter), just as {\it 'the
inconceivable effectiveness of mathematics in natural sciences'}.

Of course I was glad that it worked out this way (confounded
bitter-sweet-dark-comedy-{\it
'Twilight-Zone'-'dazed-and-confused'} to be more precise), but
several days later I remembered stories with mathematicians seeing
nonexistent hidden messages or talking to the aliens, so I
realized that there may be another rational explanation: thus I
emailed an alien (Daniel) to confirm at least my intention of
seriously reading the wise ancients and Daniel gracefully
confirmed (although a little bit confounded; when he inquired me
later for more explanations and in order to cool him down I made a
bet with a beer towards implementing the ideas of the wise
ancients in what I am doing; this bet is to be discounted because
I already knew the outcome and Daniel had already proven me right:
the key to the problem at hand was sitting on the counter the
whole time; one must only have to look straight at it while having
fresh in the memory the Bianchi-Lie ansatz in order to see it).

Archimedes is not only the best mathematician of the antiquity and
his opinions are not only of historical importance, but useful in
the study of current problems. Thus I would like to also thank
Archimedes for his choice of well crafted words (I had an early
interpretation, but it still did not fully comply with Archimedes'
written point of view, so I had to further refine it until, I
believe, Archimedes was completely satisfied).

The miracle of the {\it Archimedes-Bianchi-Lie} (ABL) method is
that Archimedes put {\it The Method} in the Bianchi-Lie ansatz,
previously used in a hurry as a strange defining axiom and rapidly
discarded as being a nonsense with no other usefulness. Bianchi
and Lie could not use the full method (although probably they
would have approved its use in the exact wording of Archimedes),
because it was a nonsense probably even according to their own
beliefs. Thus Bianchi in (\cite{B1},\S\ 374) points out a singular
behavior and takes part of it as a defining axiom. Moser in
\cite{M1} uses cautious wording and provides a number of examples;
it may be the case that Moser points out Archimedes' method, so
the rabbit hole goes even deeper (but becomes slimmer) and in a
full circle the study of the geometry of quadrics returns to its
true beginnings and not just to the beginning of the
XIX$^{\mathrm{th}}$ century: integrable system means nothing more
and nothing less than cases when one can apply Archimedes' method
of integration. It is only the French and their abstract
tendencies (see Arnold \cite{A1}) that led Darboux to dare state
such nonsense as a surface of revolution being deformed to a line,
but this was still not enough: one needs a line being {\it an
actual} deformation of a surface of revolution.

Basically {\it The Method} of Archimedes states that $P\Rightarrow
Q$ still give elegant, relevant, useful and valid consequences
$Q$, even if $P$ is a statement not accepted as being true by the
standards of proof of the time. Such an example of Q is the
statement {\it 'E pur si move!'} or the correct exponent obtained
by the scaling trick to find the needed exponent for various
inequalities involving Sobolev spaces before proving those
inequalities. Or in general all tricks used to see if a statement
is true ({\it reductio ad absurdum}), without attacking directly
the full proof of the validity of that statement (usually very
difficult); when the tricks do not point out an error it means
that one must muster courage to attack the full proof. Moreover
the whole initial effort of trying to disprove $P$ must always be
kept in mind, instead of being thrown away once its apparent
usefulness has been exhausted. Most of the time $P$ will turn out
to be true, even if not according to the standards of proof of the
time. Note that simple logic says $\neg Q\Rightarrow \neg P$; when
it is frustrating that $Q$ is headstrong enough to remain true,
that is actually a good sign that $P$ is true, just as happened
with Euclid's fifth postulate.

Luckily Archimedes had faith that P is true, even though he could
not prove it and he tried {\it The Method} on numerous examples to
confirm his beliefs and intuition, before sending a letter to
Erathostenes of the University of Alexandria. To this day
Archimedes' sentence P is not proven {\it 'such as the ancients
required'}, but instead and independently of Archimedes allowed as
a foundational axiom by the founding fathers of Calculus Leibnitz
and Newton (although the buzz around Archimedes' palimpsest is
that Archimedes may have already performed the trick).

Archimedes does not mind stating that he uses an unproven sentence
$P$, as he had already secured his back and reputation with enough
examples pointing to its validity; thus we have a great mind (he
coined the wording with the Master's perfection) being honest and
that is all that a novice needed 2300 years later. I have tried
initially to pay only partial attention to the Master's opinion
only to see later that the Master's opinion was perfect and I was
wrong in interpreting it.

The beauty of {\it The Method} of Archimedes is that it works in
some cases with same efficiency even if P is not accepted as true
by the standards of proof of the time and has a very small chance
of being accepted as true by future generations: this is the
message I got from Archimedes. So Archimedes told me to assume
that a ruling on a confocal quadric is a deformation of the
considered quadric and not to pay attention to such a {\it
'ridiculous'} statement, but only to its consequences. Before
Archimedes I justified the reflection property of the considered
distributions in the tangent bundle of the considered quadric (a
simple geometric procedure, but it does not appear naturally from
a geometric point of view, as Bianchi did not think about it) with
changing the ruling family of the confocal quadric (and thus
naturally from an analytic point of view); also the rigid motion
provided by the Ivory affinity naturally appeared from an analytic
point of view, but not from a geometric point of view, as even I
could not think about using ridiculous statements (at some point
and while trying to find vacuum solitons I thought about using
this rigid motion to see the applicability correspondence between
a ruling on a confocal quadric and a region of the considered
quadric: the corresponding region is a ruling on the considered
quadric, but the planes of the facets are just tangent planes to
the considered quadric along the ruling, each counted with an
$\infty$ multiplicity, so I dropped this nonsense immediately). At
some point the stronger version of Darboux's statement was present
in these notes as versions equivalent to it and obtained by
rolling, but that just because I had a vague idea about Darboux's
statement; when challenged to remove the stronger statement I knew
that I had to ask Monsieur Darboux about his opinion on the
matter.

Thus I saw the geometric properties of the moving picture as
consequences of the geometric properties of the static picture,
but the later naturally appeared only from an analytic point of
view (changing the ruling family and the use of a certain rigid
motion as to whose appearance I did not even dare question). After
Archimedes the geometric properties of the static picture are
consequences of the geometric properties of the moving picture and
thus naturally appear from a geometric point of view (see again
Arnold \cite{A1}).

Principles of a general nature transcend singularities and remain
valid even for singularities; moreover such principles recorded
for singularities are the only important information needed to
confirm the validity of the initial principles of a general
nature; such ideas can be found in the current literature of
integrable systems (like Hilbert-Riemann or scattering data) or
maybe in other domains and Archimedes' {\it The Method} is the
first example (so Archimedes basically proved that the balance
with fixed fulcrum survives integration and the only important
thing to do is to find the correct balance configuration at the
infinitesimal level by use of simple Euclidean properties of the
studied object: {\it he transferred and stuck slices (the
equivalent of infinitesimal leaves or facets) at their center with
$\infty$ multiplicity at the fixed left end of the balance}; this
is {\it precisely} the above mentioned nonsense of Bianchi-Lie's
and thus Archimedes' principle of the balance surviving
integration also applies to the theory of deformations of
quadrics).

Certain properties of certain structures remain valid even when
the {\it subjacent} (supporting) structure collapses to
nothingness, to ridiculous, to irrational (even the method of
finding the link to {\it The Method} of Archimedes is an
application of the same method, since an irrational hypothesis
(signs are shown to me) has a rational and perfectly founded
conclusion), thus providing not a paradox or a false ridiculous
statement, but the simplest explanation of those structures.

The area of a segment of length $f(x)$ is of course $0$, but if
that segment is a limit of thinning rectangles of same height
$f(x)$, something sticks to it and it has a memory containing more
information: the area is $f(x)dx$, an infinitesimal which by
itself and for all practical purposes is of course $0$, but it is
definitely something different from $0$. So this is the centuries
old story, told and retold to high-school kids in some countries
(and to college level kids in others).

Archimedes not only said essentially the same story more than two
thousand years ago, {\it but he also gave its geometric
description} and it is the power of geometry that makes
Archimedes' ideas amenable to generalizations in another settings,
just as {\it 'the inconceivable effectiveness of mathematics in
natural sciences' }, usually realized by a model supported on a
simple geometric picture. Bianchi's contemporaries acclaimed {\it
'Qui la geometria vive!'} and yet they did not know the whole
story.

The beauty of the ABL method is that it gives almost for free more
information from the beginning, so I believe it to be the
essential ingredient to supporting generalizations in other
settings (especially in higher dimensions); in dimension $2$ is
just a fancier reformulation of the old point of view. Of course
the information already given in dimension $2$ by the old point of
view may be generalized with a heroic effort to higher dimensions
without Archimedes' point of view, but Archimedes' point of view
is the natural one, it appearing from a {\it geometric} point of
view instead of an {\it analytic} point of view. And it is over
the course of history that {\it not} the geometric picture
subjacent to physical facts (like gravity and light), {\it but}
its analytic interpretation changed. On top of that it bypasses
complicated computations since it gives for free {\it all}
algebraic identities of the static picture necessary for the
completion of the differential identities required by the moving
picture, instead of trying to complete first the later and in
doing so see what are the former. To do this trick we need to find
the correct discretization of the said differential identities.

Consider the theory of deformations of $\mathbf{H}^n(\mathbb{R})$
in $\mathbb{R}^{2n-1}$ from Tenenblat-Terng \cite{TT1}.

From the application of the ABL method one gets immediately the
$\mathbf{O}_{n-1}(\mathbb{R})\times\mathbf{O}_{n-1}(\mathbb{R})$
symmetry due to the symmetries (in the normal bundle) of the
rolling and to the symmetry of the tangency configuration. But
these symmetries never use the fact that we work with
$\mathbf{H}^n(\mathbb{R})$, so are valid for the deformations of a
general $n$-dimensional real quadric in $\mathbb{R}^{2n-1}$
(Berger, Bryant, Griffiths \cite{BB}), should such a theory of
deformation be developed in this case. For my birthday in Spring
2006 the true gift was not the congratulatory phone call from
Romania, but finding the article of Berger, Bryant and Griffiths
(I believed it to be true, began with Cartan \cite{C2} and
followed the trail, only to find that the most difficult part
(namely the existence of the seeds of the theory) had already been
proven). Although it is not clear to me if Berger, Bryant and
Griffiths proved the flat normal bundle property, it should follow
as a simple intuitive consequence of the {\it initial value data}
of the differential system subjacent to the B\"{a}cklund
transformation (counting the leaves; their numbers are the same as
for Tenenblat-Terng since we just generalize their point of view
from a metric-projective point of view), since no functional
information is allowed in the normal bundle except a finite
$\mathbf{O}_{n-1}(\mathbb{R})$ dimensionality.

Note that the Bianchi Permutability Theorem from Terng \cite{T2}
is the correct discretization of the differential equations
required by existence of deformations of
$\mathbf{H}^n(\mathbb{R})$ in $\mathbb{R}^{2n-1}$, that is the
proof {\it 'such as the ancients required'}. The metric-projective
properties of this picture can be assumed to be true at the level
of the static picture for deformations of $n$-dimensional real
quadrics in $\mathbb{R}^{2n-1}$ by means of the ABL method and get
the needed algebraic information (it seems that the necessary
algebraic identities required by the existence of the B
transformation are encoded by the second iterated tangency
configuration and those required by the existence of the Bianchi
Permutability Theorem are encoded by the third iteration of the
tangency configuration; no further algebraic identities are
needed).

Most continuous groups of symmetries and symmetric spaces are
given by quadratic equations; some of their metric properties may
not be due to the big group of symmetries, but to the quadratic
definition: we call such properties {\it metric-projective}. Of
course to deform such things one needs a surrounding space with a
big group of symmetries and certain simplifying assumptions.

According to Terng-Thorbergsson \cite{TT2}: {\it '...while
extending the theory of sub-manifolds to ambient spaces more
general than space forms proves quite difficult if one tries to
use the same approach as for the space forms, at least for
symmetric spaces it has proved possible to develop an elegant
theory based on focal structure that reduces to the classical
theory in the case of space forms.'}

Thus should a theory of deformations of higher dimensional general
quadrics be developed in the Euclidean space, the only important
thing remaining to do is the construction of the confocal family
of a quadric in a symmetric space; probably this can be easily
done so as to satisfy the usual classical metric-projective
properties.

As to the undergraduate education at Bucharest University in what
concerns the prerequisites from my higher mathematical education
needed to read the classical geometers I wish to thank Professor
Kostake Teleman for the first year Geometry course (where it was
discussed almost all necessary prerequisites for the static
picture for quadrics including vectorial, affine and projective
spaces and their transformations; bilinear and multilinear forms;
eigenvalues and eigenvectors; diagonalization of symmetric real
and hermitian matrices; Sylvester's theorem; metric, affine and
projective classification of real quadrics; projective
classification of complex quadrics; real confocal quadrics;
B\'{e}zout's theorem for algebraic curves; etc), Assistant Victor
Vuletescu for being both the Instructor and the Assistant for the
first semester of the second year Geometry course (where the
classical geometry of curves and hyper-surfaces was summarily
covered), Professor Martin Jurchescu for the second year Complex
Analysis course (there was also the Assistant I do not remember
the name, but I met him at Notre Dame first in Summer 2004; I
distinctly remember he did a problem with a Schwartz reflection
symmetry and a problem on the exam was almost instantaneously
solved with a rotational $\mathbb{Z}_3$ symmetry plus the Schwartz
reflection symmetry and a Cauchy continuity at the origin), for
the third year Calculus on Manifolds course (when it was clearly
revealed to me by Guillemin-Pollack that any abstract manifold can
appear as a real Euclidean manifold in Poincar\'{e}'s
interpretation) and for the fourth year Multivariable Complex
Analysis course (the Romanian Dimitrie Pompeiu's generalization of
Cauchy's point of view is just a fancier reformulation for one
variable, but essential for more than one variable); unfortunately
he got sick that year; Assistant George Marinescu took over the
second semester of the course and completed it with the proof of
Nirenberg-Newlander. Now the teaching method in Romania is on the
French model: students go to lectures (which are not discussion
style, but rather the lecturer putting a whole Rudin, Ahlfors or
Schwartz on the blackboard out of memory, of course in one's own
interpretation, just as a conductor sometimes allows oneself
freedom in interpreting the message of the composer (meantime the
audience keeps quiet) or just like QT allowed himself changing the
time-line of {\it 'Pulp Fiction'}; the favor is asked back from
the student at the end of the semester) and tutorials for three
quarters of the semester and have only tough {\it written} final
examinations in the exams session of a quarter of a semester
(roughly one exam per week) and somehow easier {\it written}
examinations in a special exams session when students may choose
at the pleasure of their own heart to challenge at most two of the
first grades; as a consequence I still remember the page of the
Complex Analysis course I was reading when, being frustrated with
what I was reading, I decided together with my roommate to throw
water bags from the dorm window hoping that the frustration will
{\it stick to the water bags and transfer} to the poor passers who
happened to walk below our window and conversely.

At that time my undergraduate Advisor Professor Stere Ianu\c{s}
was doing harmonic maps, so he gave me a book of Eells-Rato on
harmonic maps with symmetries and positive curvature targets to
write the BS thesis (so all research should begin by writing some
notes to oneself where one explains in detail already known
results in an area of consistent research); in doing so I also
almost understood the article of Eells-Sampson (it is considered
the founding article of the branch of harmonic maps; using the
heat flow it provided existence of harmonic maps in negative
curvature targets). For my MS thesis at Bucharest University I
remember reading Lawson's minimal surfaces, sub-manifolds and
currents as examples of harmonic maps (I also met with Monsieur
Lawson in a course of Spin Geometry); as I was leaving Romania
Professor Stere Ianu\c{s} intended for me to read something about
loop groups (after trying a stint with isoparametric
sub-manifolds; I believe I saw something but by the time I
investigated it {\it 'I lost it, Lou!'}); I found out later that
these have something to do with deformations of quadrics. Now
harmonic maps with positive curvature targets (for example spheres
or ellipsoids) are very difficult to study; consequently one must
assume initially lots of symmetries to simplify the problem. The
Ph. D. thesis of Smith in early 1980's did just that (note also
that that of Wood's showed that the Schwartz reflection symmetry
is valid for harmonic maps in dimension $2$): CMC surfaces have
harmonic Gau\ss\ map into the unit sphere $\mathbb{S}^2$; if one
assumes this Gau\ss\ map to have rotational $\mathbb{S}^1$
symmetry, then one gets Delaunay's CMC surfaces of revolution
(this observation on Smith's results for $2$-dimensional target is
due to Calabi). Note that if one assumes the rotational
$\mathbb{S}^1$ symmetry only at the level of the linear element,
then one gets more CMC surfaces than those of Delaunay's: such
surfaces were found by Smyth in 1987, are complete, always have an
umbilic and outer dihedral symmetry (so the Gau\ss\ map does not
have rotational symmetry). Recently Dorfmeister, Pedit, Wu and
collaborators have analyzed CMC surfaces with the loop groups
tools and in particular have generalized these Smyth surfaces to
various other settings.

Finally on a very short list of people which had the most benefic
influence on my mathematical training by going beyond what was
required of them I would like to thank my mother Floarea
Dinc\u{a}, {\it \^{I}nv\u{a}\c{t}\u{a}torul Rom\^{a}n} incarnated
for me as my primary school teacher Didi\c{t}a Staicu, my
7$^{\mathrm{th}}$ grade mathematics teacher Lucia Velicu, my
formerly high-school mathematics teacher (currently Professor
Dumitru Popa of Ovidius University Constan\c{t}a), the former
mathematics Inspector Gheorghe Popescu of the Ialomi\c{t}a County,
the {\it Romanian Education Foundation} of Stanford, CA and the
weekly {\it Academia Ca\c{t}avencu}. Further this last part of the
acknowledgements is definitely biased and incomplete, it being
designed to suit a restricted {\it educational} point of view;
besides that it has a personal touch (according to Tyler Durden
{\it 'You want to make an omelet, you got to break some eggs'} and
it is the mouth of the sinner (the undersigned) which spells out
the truth: one can recognize the power of generality of
mathematics when the same model describes a-priori unrelated
situations).

It is because my mother has paid a subscription to the monthly
{\it Gazeta Matematic\u{a}} ever since I was in the
5$^{\mathrm{th}}$ grade that I had later things to munch on. She
was my best Advisor although she had no education beyond primary
school and I did not always agree with her findings (note however
that I being a minor and she being my legal tutor she was my only
Advisor ever who could {\it legally} choose and did chose notably
in Summer 1990 to impose on me one's opinion first and foremost of
mine). It is from her that I learned the useful lesson of {\it
standing and defending my ground 'lest The Heavens will fall!'}
(just as should God meet with the Satori sword from {\it 'Kill
Bill'} (designed according to the Romanian Gheorghe Zamfir's
pan-flute instructions mind you), God Himself will be cut, should
I smell that God is trying to pull a fast one on me I will hit God
Himself with a confirmation delivery letter ({\it So Help Me
God!}), because the receipt of Uncle Sam's Postal Services is the
alien's only friend on Uncle Sam's territory): one of our
neighbors' son rose to the rank of party leader of the Buz\u{a}u
county; consequently our neighbor decided to encroach on our
territory and my mother took the problem to a county court of law,
where of course she lost. Her lack of education (her father denied
his girls education beyond primary school, as he considered that
{\it B\u{a}r\u{a}gan} (the plain with a mind of its own and where
the harsh Siberian winter wind roams free; on a good rainy year it
provides the bulk of the Romanian grain which in turn used to
decide the price of grain in Europe until WWII; on a sequence of
bad dry years it almost starved to death its poorest inhabitants;
note also that it gave to the Romanians their first unifier {\it
Mihai Viteazu}, Flower's child of Floarea in the lost City of
Flocks situated 10 miles east of \c{T}\u{a}nd\u{a}rei and from
where \c{T}\u{a}nd\u{a}rei probably {\it split}) would provide the
remaining needed education; it is only after the change to
communist regime in Romania that two of her younger sisters
completed at least high school education) may be the possible
explanation of the fact that all her four children (boys) got a
college education (beginning with August 1975 she was the head of
the family and consequently worked low paying jobs as unqualified
worker at a local brick factory until she retired in Summer 1992
on a low pension on account of health reasons and worked in our
privately-owned agriculture during the Springs). In Summer 1995 I
was almost defeated by a string of bad grades; just like
\c{S}tefan cel Mare's mother put some faith into him when he was
almost defeated by the Turks, she put some faith into me, I went
back and I aced the PDE exam (it was not a perfect score, but
luckily for me the Instructor decided to round the score upwards
before seeing my student ID and my other grades); looking back I
can honestly say that no formal education was needed after that
exam, as I have not seen since then interesting tricks, methods
and tools after being hit that year with infinite dimensional
spaces (Banach, Hilbert and locally convex) and the essence of
symmetries (Lie Algebras); the {\it Fourier transform} words were
not uttered until its natural existence Shwartz's space of
tempered distribution was properly introduced. Now the custom and
only publicly recognized form of corruption in Romania is to bribe
{\it The Godfather} (Comptroller) with 30$\%$ of the train ticket
price; the briber (usually student or other small income frequent
train user social categories, because one needs both motivation
and experience; the new student user is helped by The experienced
Godfather's nonchalance in receiving the bribe and The beginner
Godfather by the experienced student's nonchalance in offering the
bribe) risks a civil fine up to $5-10$ times the price of the
ticket if caught ticket-less by The Uber-Comptroller, that is the
inspector of The Godfather's work; it seems that even The Csar of
Russia initially refused the help of The Romanian Army in the
1877-1878 Russo-Turkish war on account that a nation which does
not have the decency of buying a ticket before boarding a train
cannot be trusted as an ally in waging war; I took full advantage
of this custom as a student (in the Summers and early Falls I was
doing agriculture at \c{T}\u{a}nd\u{a}rei and during the academic
year I was carrying food to Bucharest every fortnight) and never
paid a fine; the only time when I came close to paying one (The
Uber-Comptroller confiscated our student ID's) The Wolf of
\c{T}\u{a}nd\u{a}rei happened to pass in the neighborhood and
solved the problem by producing a railway worker ID and shaming
The Uber-Comptroller into total submission; usually when The
Uber-Comptroller meets a group of students he also meets
resistance and consequently after difficult negotiations he can
bestow at most a fine fine to the whole group; thus {\it the power
lies within the group, at least five'o'em} (in which case one gets
almost even).

It is because of Didi\c{t}a Staicu that I discovered the thrill of
{\it reading and writing mathematics}: just like Jules she had an
{\it 'ultimate piercing'} look and would not let yours escape once
caught and would pierce and shake the very foundations of your
soul; consequently I still have trouble looking her in the eyes
when we meet on the street (that may be also due to my not so
proud to be of deeds).

It is because of Lucia Velicu that I let aside the life of
mischief and tormenting my teachers and our neighbors from the
Strachina neighborhood (I was first of my class, but my deeds were
still best described by Ion Creang\u{a}'s {\it 'Memories of my
childhood'}) and began doing serious mathematics: she organized a
weekly off-curriculum mathematics club. I remember my first idea
to generalize $\sqrt{2^4}=2^2$ was $\sqrt{2^n}=2^{n-2}$ (it is
customarily to remember mostly unpleasant situations although the
name itself explains: square root of $a$ is the root of $x$
squared equals $a$; even the Babilonians knew Newton's method on
how to approximate square roots thousands of years before Newton
got any linking about this method). Lucia asked us in Fall 1986 to
buy a book of problems; luckily my colleague and good friend
C\u{a}t\u{a}lin lent me some money (first and last time I touched
money not belonging to me was at age of 8 when I was playing
hide-and-seek with C\u{a}t\u{a}lin and his sister in their home; I
hid near the porcelain pig and consequently some of their money
{\it stuck to my sticky hands and transferred to my pockets};
luckily his father chose just to pull my ears the next morning
when I met C\u{a}t\u{a}lin as usually on our way to school and
without bringing my deed to the knowledge of my mother) and I
bought the last book from the bookstore; I read many answers to
the geometry problems from the answer section of that book on the
evening before the county competition in Spring 1987, {\it as they
were explained at large} (the problem with the current
mathematician is that {\it one's mathematics is difficult}; note
however that that difficulty may also arise from assuming the
student to be familiar with unfamiliar background, letting the
student fill in {\it 'the details'} by oneself as an unnecessary
homework assignment (thus virtually {\it rebuilding in a costly
masochistic fashion} the already built theory instead of {\it
investing a fraction of the same time} in understanding the same,
namely the full method of investigation) or from failing to
motivate and reinforce by means of simple examples the purpose of
the investigation (according to Number One of {\it 'Austin Powers
Goldmember': 'Throw me a ... bone!'}) and asking at each step of
the investigation the motivation of that step: taking definitions
and other things for granted produces {\it the sleep of ration});
I distinctly remember Euler's 9 points. Thus in Spring 1987 I
barely missed the national mathematics competition, but as a
consequence I acceded to county mathematics programs beginning
with Fall 1987. In Spring 1988 I got almost the perfect score for
the geometry examination of the national mathematics competition;
in the algebraic examination I remember guessing a solution of an
equation which after some change of variable (skilful handling at
that time) became a quadratic equation in $x^2$, so it had four
solutions.

It is because of Professor Dumitru Popa that I continued doing
serious mathematics: he took me under his wing at a bad time in my
life (I had just miserably failed the Romanian Literature and
Grammar high-school entry examination; on top of that The Wolf
caught me cheating first time in my life on the first English test
in high-school (in a full circle last time I cheated was on the
last test in high-school; it was French and I was caught for the
second time only); these may have also motivated me to pursue
mathematics with renewed vigor in order to redeem myself); also he
imposed his will on other teachers and always arranged that I
should have at least one week off school before mathematics
competitions (or, as a recognition of my efforts, a whole month in
1989 after returning from a national program in mathematics), a
thing I sorely missed during my last two years in high-school (in
Fall 1991 {\it I declined participation to an intercounty
mathematics competition on account of not being prepared, as no
due one week notice was given to me and instead I was fighting for
passing grades at subjects of study such as geography technology}
(drawing maps and memorizing names of Romania's mountains)). It
does not matter that the group of three-four students preparing
for the mathematics competition were in fact frustrated with the
mathematics problems and instead were playing ping-pong or various
games with chalk or coins during class time only to show diligence
before and during recesses, when we were inspected to show
progress on our work, because that work was enough (an impromptu
inspection even caught us in the act of throwing darts at the
door, as the door opened to let the inspector of our work do one's
duty). And even playing with coins can have benefic consequences:
probably \c{T}i\c{t}eica was playing with the change while
drinking beer with friends and thus noticed the problem with four
equal coins (a law of a general nature, probably amenable to
interesting generalizations for conics or even in higher
dimensions; it actually has something to do with {\it Barbilian
geometries}).

It is because of former mathematics Inspector Gheorghe Popescu of
the Ialomi\c{t}a County that I acceded to special national
programs in mathematics in Romania: due to a certain sequence of
events he chose to settle a dispute by choosing {\it a third
neutral resolution}: he bet on my geometrical abilities in Spring
1989, in the process bending some rules.

Interestingly enough in Spring 1989 I misunderstood a difficult
algebraic problem which took most of the competitors most of their
time to solve; I gave in few minutes a faulty solution and as a
consequence I had most of the remaining time to spend on the other
three problems. The remaining problem with an algebraic character
was similar to a problem which had been previously discussed by
Silviu Ciuperc\u{a} (formerly mathematics teacher at the
\c{T}\u{a}nd\u{a}rei high-school and currently turned politician)
at a county mathematics program in early 1989 (unfortunately the
last such county mathematics program ended {\it by necessity} on
the very same day with the communist regime in Romania). For one
of the geometry problems I was asked to provide a synthetic proof
to a metric-projective result of Apollonius (the plane locus of
points with constant ratio of distances to two fixed points of the
plane is a circle; when the ratio is $1$ the circle passes through
the line at $\infty$ and becomes a line in the Euclidean plane). I
knew \c{T}i\c{t}eica's book of Euclidean geometry almost by heart,
since for 4-6 hours of classes a day I was spending 9 hours on the
train and waiting in the Slobozia train station hallway (I could
have travelled around the equator with 40 miles a day for almost
four years). I had to turn back 16 years later to Apollonius and
his contemporaries, in the process discovering Arnold \cite{A1} in
early 2006 and realizing that I was right in my convictions for
the last three years.

It is because the {\it Romanian Education Foundation} paid my
TOEFL and GRE general test fees that I found later a lender to
borrow the money for the GRE mathematics test fee and thus I
acceded to a Ph. D. program in mathematics in USA (I got this REF
sponsorship tip from a fellow student one year senior to me;
unfortunately their first letter was {\it 'O scrisoare
pierdut\u{a}'} (A lost letter, as many other letters coming to
Romania from abroad; some even reach Hungary on account of
Bucharest being its capital city and thus Budapest being the
capital of Romania) and did not reach my dorm address; after
contacting them again the second one reached my permanent address,
but due to the delay and not being able to find a lender by early
Summer 1996 I had to postpone the application process with one
year).

It is because the weekly {\it Academia Ca\c{t}avencu} paid the
airline ticket that I was able to show myself at such a Ph. D.
program (I got this sponsorship tip from a fellow student, who
with the same sponsorship was able to show himself at a summer
school in Italy in 1997; on top of that their motto is {\it 'Our
readers are brighter than their readers'}, so I pressured them to
deliver on their ideas). That would not have probably been
necessary, as I had other options (including offers to borrow
money from some of my fellow students and Faculty of the
Mathematics Department of the Bucharest University), but being in
debt with a smaller amount of money gave me peace of mind.

\subsection{General remarks}
\noindent

\noindent As much as possible and always keeping in mind an
optimal mixture of clarity, convenience and simplicity we shall
choose our working definitions, notions, style of presentation,
etc in decreasing order of importance beginning from the classical
ones, followed by the current ones, followed by our choices, but
the later always not far removed from the context. Most likely the
later are not the best choice: for example {\it adjugate} was {\it
adjoint} and before that {\it dual}.

Take for example the use of the classical denomination {\it
applicable} instead of the current one {\it isometric}; the former
acknowledges that the surfaces naturally exist in a surrounding
space, while the later pays attention only to the induced linear
elements. For the same reason we use {\it deformation} instead of
{\it isometric deformation}, as there is no such later thing (the
induced linear element stays the same, so the abstract surface
itself (obtained by deletion of the surrounding space) does not
undergo any deformation). Another choice would be {\it deformation
through bending}, but this implies continuity of the process and
{\it deformation through warping} is already too complicated (see
Sabitov \cite{S1} or Spivak (\cite{S3}, Vol {\bf 5})). For the
same reason we use {\it Gau\ss\ -Weingarten equations} instead of
{\it fundamental equations of immersion}, since they are
consequences of the existence of a sub-manifold in a surrounding
space and not necessary equations required for the existence of an
isometric immersion of an abstract sub-manifold in another one.
The main point of the theory of deformations of quadrics is that
we always assume the existence of an initial ({\it seed})
deformation in a surrounding space.

In the spirit of the classical geometers we shall preserve the
thread of exposition of facts, without (numbered) theorems,
lemmas, etc.

Due to the frequent use of certain key words we shall define and
use abbreviated notations.

Due to the frequent use of  the reference Bianchi (\cite{B2},Vol
{\bf 4},(122,...)), it will be replaced with Bianchi (122,...).
Note also that Bianchi (\cite{B1},Vol 2, Part 1) contains three
chapters about deformations of quadrics; it may be the last
revised version of Bianchi's.

References are restricted to the ones I am aware of being of
direct interest here. There is a huge literature on surfaces and
integrable systems and for example Rogers-Schieff \cite{RS1} has a
consistent list of references, from where one can probably get in
at most one more step to much of the current literature in this
area. As usual, unwarranted absence of credit does not imply claim
of credit, but only that the particular fact is considered common
knowledge or that I am not aware of such correlations of names and
facts. A further hinder in this matter is the fact that the
classical geometers and their results are not considered
mainstream mathematics today; seldom there are monographs in
various branches of geometry which give short descriptions of the
results of the classical geometers in their respective historical
perspective sections. Consequently the current geometer trusts
only part of the results of the classical geometers and safely
assumes that the only logical conclusion that can be drawn from
the fact that notions such as rolling, the use of imaginaries, etc
are not present in the current literature is that there is
something fundamentally wrong with these; I had my share of doubts
on the use of imaginaries and to this day I still believe such an
explanation may be possible, but until {\it I find such a valid
argument} I will side with the classical geometers. Note that a
big part of the classical references is freely available online:
for example Bianchi (\cite{B2},Vol 2) and an earlier one volume
version from 1894, Darboux \cite{D1} and Eisenhart \cite{E2} can
be found in{\href{http://www.hti.umich.edu/u/umhistmath/} {The
University of Michigan Historical Mathematics Collection}}.
Unfortunately from Archimedes' collected works (written in the
XIX$^\mathrm{th}$ century) is missing exactly the needed piece,
but again {\href{http://www.gutenberg.org/etext/7825} {{\it The
Method}}} is freely available on the web.

There is no replacement of Bianchi's works and the present notes
do not intend this; they may be used in a parallel reading, as
there are many other ideas and points of view in Bianchi's works
that do not appear here. While a simpler proof is sometimes
useful, ideas from the original proof can be used in other
settings; this is another point of these notes.

The sine-Gordon equation $\om_{uv}=\sin\om$ appears as the Gau\ss\
equation of the {\it Gau\ss-Codazzi-Mainardi} (GCM) equations (the
integrability conditions of the {\it Gau\ss-Weingarten} (GW)
equations) for {\it constant Gau\ss\ curvature} (CGC) $-1$
surfaces, where $(u,v)$ are {\it Tchebyshev coordinates}
(coordinates for which $u$ is the arc-length parameter on the
$v=$const curves and conversely; that is the taylor's problem with
a piece of fabric) which further are asymptotic and $\om$ is the
angle between the asymptotic directions; the Codazzi-Mainardi
equations are identically satisfied in this setting. Note that the
complete CGC $-1$ surface can be realized as a connected component
of a particular quadric (the {\it space-like pseudo-sphere}) in
the Lorentz space of signature $(2,1)$. The {\it space-like}
refers to the induced linear element being positive definite
(equivalently the normal is time-like; the pseudo-sphere itself is
seen as time-like from the origin, since all normals pass through
the origin).

While the sine-Gordon equation received much attention in recent
years, I am not aware of any similar treatment of deformations of
general quadrics: Moser discusses in \cite{M1} only the geodesic
flow on quadrics and points out analytic similarities between it
and other integrable systems (like the {\it Korteweg-de Vries}
(KdV) equation), while Burstall in the monograph \cite{B4} on {\it
isothermic} surfaces (surfaces with isothermic lines of curvature)
does not discuss the link to deformations of quadrics. In
\cite{TU1} Terng and Uhlenbeck discuss, via a highly technical Lie
algebraic formalism, a natural generalization of the classical
{\it B\"{a}cklund} (B) transformation of CGC $-1$ surfaces: the B
transformation of the {\it Zakharov, \v{S}abat, Ablowitz, Kaup,
Newell, Segur} (ZS-AKNS) $\mathbf{sl}_n(\mathbb{C})$ hierarchy.
Although they do not specifically mention Bianchi's {\it
B\"{a}cklund} (B) transformation of deformations of general
quadrics, there are similarities; thus it may be the case that
their approach contains Bianchi's B transformation as a particular
case, but only from an analytical point of view and without the
simpler geometric picture. Note however that hierarchies must be
constructed at the analytic level, as at this point any geometric
meaning is above me; here again Terng-Uhlenbeck \cite{TU1} may
already have done that.

The {\it B\"{a}cklund} (B) transformation of local solutions of an
integrable system is a procedure of generating {\it leaves} (new
local solutions) from a {\it seed} (given local solution). This
procedure requires the integration of a simpler auxiliary
differential system associated to the given local solution.
Sometimes this is called an {\it auto}-B transformation, following
that a {\it general} B transformation may be a procedure of
generating local solutions of an integrable system from local
solutions of another integrable system, again requiring the
integration of a simpler auxiliary differential system (there are
also algebraic procedures with the same purpose as the general B
transformation, classically known as {\it Hazzidakis} (H)
transformations, but currently as {\it Miura} transformations).
Since we use only the auto-B transformation, we shall call it B
transformation.

By necessity (all quadrics are equivalent from a complex
projective point of view) we shall consider the complexification
$(\mathbb{C}^n,<,>)$ of the Euclidean space $(\mathbb{R}^n,<,>):\
\ \ <x,y>:=x^Ty,\ |x|^2:=x^Tx,\ x,y\in\mathbb{C}^n$. {\it
Isotropic} (null) vectors are those vectors of length $0:\
|v|^2=0$; since most vectors are not isotropic we call {\it
vector} a vector presumed non-isotropic and we shall only
emphasize {\it isotropic} when the vector is assumed to be
isotropic. This terminology will also apply in other settings: for
example we call {\it quadric} a non-degenerate quadric, following
that we shall only emphasize {\it degenerate} when the quadric is
assumed to be degenerate. Standard geometric formulae for
sub-manifolds in $\mathbb{R}^n$ remain valid (with their usual
denomination) almost everywhere. If the Euclidean product on
$\mathbb{C}^n$ induces a non-degenerate product on the
complexification of a tangent space (this complexification is
assumed to be a proper subspace of $\mathbb{C}^n$) of a general
sub-manifold in $\mathbb{C}^n$, then one can find an orthogonal
complement spanned by an orthonormal frame. Since the GW and the
{\it Gau\ss-Codazzi-Mainardi-Ricci} (GCMR) equations for a
sub-manifold $x\subset\mathbb{R}^n$ assume only the non-degeneracy
of the linear element of $x$ (and thus the existence of a moving
orthonormal frame in the normal bundle), they are still valid and
sufficient to describe the geometry of general sub-manifolds in
$\mathbb{C}^n$ almost everywhere. For example the Lorentz space of
signature $(2,1)$ can be realized as the {\it totally real
subspace} $\mathbb{R}^2\times i\mathbb{R}\subset\mathbb{C}^3$ (the
induced scalar product is real (valued) and non-degenerate). Note
that there are many definitions of totally real subspaces, some
involving a hermitian product, but all definitions coincide. All
$3$-dimensional totally real affine subspaces of $\mathbb{C}^3$
are obtained by applying a rigid motion
$(R,t)\in\mathbf{O}_3(\mathbb{C})\ltimes \mathbb{C}^3,\
(R,t)x:=Rx+t,\ x\in\mathbb{C}^3$ to one of the Lorentz spaces
$\mathbb{R}^3,\ \mathbb{R}^2\times i\mathbb{R},\
\mathbb{R}\times(i\mathbb{R})^2,\ (i\mathbb{R})^3$.

We call {\it surface} any sub-manifold of $\mathbb{C}^3\simeq
\mathbb{R}^6$ of real dimension $2,3$ or $4$ such that all its
{\it complexified tangent spaces} (called henceforth {\it tangent
spaces}) have complex dimension $2$. In this case the distribution
$T\cap iT$ formed by intersecting real tangent spaces $T$ with
$iT$ has constant real dimension (which must respectively be $0,\
2$ or $4$). This distribution is integrable (by its own
definition); by an application of Frobenius one obtains leaves; on
leaves the Nijenhuis condition holds and by Nirenberg-Newlander
this is equivalent (locally) to prescribing
$x:D\rightarrow\mathbb{C}^3$, where $D$ is a domain of
$\mathbb{R}^2$ or $\mathbb{C}\times\mathbb{R}$ or $\mathbb{C}^2$
such that $dx\times\wedge dx\neq 0$ (the $\times$ applies to the
vector structure and the $\wedge$ to the exterior form structure,
so the order does not matter: $\times\wedge=\wedge\times$; note
also that $dx\times\wedge dx=2dx\times dx$, but we also need at
some point $dx\times\wedge dy$). The GW and GCM equations suffice
to describe the geometry of most surfaces almost everywhere and
surfaces such defined are the natural completion of the usual real
surfaces, although the usual coordinates used in geometry
(asymptotic, conjugate systems, orthogonal, principal) may not be
the ones which clearly split into purely real and purely complex
(for example asymptotic coordinates on a real surface of positive
Gau\ss\ curvature). The change of coordinates still holds, even
outside bi-holomorphic change of the complex coordinates and
diffeomorphism of the real ones, but if the real and complex
coordinates are mixed, they lose the character of being purely
real or complex and the new parametrization may not be holomorphic
in some variables. For example if we have a surface
$x=x(z,t)\subset\mathbb{C}^3,\ (z,t)\in
D\subseteq\mathbb{C}\times\mathbb{R}$ and consider two complex
functions $u=u(z,t),\ v=v(z,t)$ of $z$ and $t$ with
$J:=\frac{du\wedge dv}{dz\wedge dt}=
\begin{vmatrix}\frac{\pa u}{\pa z}&\frac{\pa u}{\pa t}\\
\frac{\pa v}{\pa z}&\frac{\pa v}{\pa t}\end{vmatrix}\neq 0$, then
one can invert $z=z(u,v),\ t=t(u,v)$ at least formally and at the
infinitesimal level by the usual calculus rule of taking the
inverse of the Jacobian: $dz=\frac{\pa z}{\pa u}du+\frac{\pa
z}{\pa v}dv,\ dt=\frac{\pa t}{\pa u}du+\frac{\pa t}{\pa v}dv$.
Although the holomorphic dependence on a variable is lost, the
$2$-dimensionality character remains clear. After using formal
coordinates of geometric meaning to simplify certain computations,
one can return to the purely real and complex coordinates, just as
the use of $z,\bar z$ parametrization on real surfaces. Note that
although on such surfaces the vector fields $\pa_z,\ \pa_{\bar z}$
do not admit integral curves (lines of coordinates), statements
about infinitesimal behavior of such lines remain valid, so one
can assume such lines to exist and derive corresponding results;
thus for all practical purposes we can assume that coordinates
descend upon lines of coordinates on the surface. From a practical
point of view usually one needs only surfaces of real dimension 2
or 4; real $3$-dimensional surfaces usually appear in the
intermediary step of obtaining real $2$-dimensional totally real
surfaces from real $4$-dimensional surfaces; they are true real
$3$-dimensional when appear as real $1$-dimensional families of
holomorphic curves and this family cannot be extended to a
holomorphic $1$-dimensional family; surfaces of revolution and
ruled surfaces should provide the simplest examples.

This goes a little bit against the intuition, as for a real
surface in the real Euclidean space $\mathbb{R}^6$ the second
fundamental form is vector valued and we also need the Ricci
equations. But keep in mind the fact that the Euclidean scalar
product on $\mathbb{C}^3$ is completely different: it has a real
and an imaginary part, both being indefinite scalar products of
signature $(3,3)$ on $\mathbb{R}^6$. If we begin with a real
analytic surface, then one can easily complexify not only the
tangent bundle, as one does currently, but also the surface and
the surrounding space; thus the asymptotic coordinates (complex
conjugate $z,\bar z$) on a real analytic surface of positive
Gau\ss\ curvature will admit integral curves on the corresponding
surface parameterized by $z,\ w$ and the initial real surface is
individuated by the non-real condition $w=\bar z$ (which implies
$z=\bar w$; these two functionally independent conditions can be
considered for computations instead of two real functionally
independent conditions imposed on the four real parameters of the
surface in order to get the real surface). Thus one can easily
derive lines of length $0$ (isotropic curves) on real analytic
surfaces; this was my first encounter with imaginaries when I read
Eisenhart \cite{E1} in Fall 2002, but I ignored the use of
imaginaries for one more semester. The {\it surfaces of
translation} (sum of two curves) with generating curves being
isotropic conjugate are precisely the minimal surfaces and this
fact is due, according to Lie, to Monge. When the curves are
holomorphic and one also restricts the parameter on one to be the
conjugate of the parameter on the other one gets {\it all real
minimal surfaces}; this is one step before the {\it Weierstrass
representation formula} (so Weierstrass essentially used the
standard parametrization of the isotropic cone), recently
generalized to other settings.

It is unclear to me at this point if the classical geometers
considered only analytic surfaces; anyway except for the results
where they specifically use the power series expansion the
remaining ones remain valid in a more general setting. By the end
of the XIX$^{\mathrm{th}}$ century existence of solutions of
differential systems in more general settings appeared (for
example Picard's existence of solutions of ODE's), so the
classical geometers were aware of constructions more general than
analytic ones.

As a graduate student being involved in a program of the study of
the classical geometry of surfaces I discovered over the course of
roughly one year the beautiful (although difficult to read at the
first cursory incursion; probably due also to the fact that for
example I discovered in the Library first Vol {\bf 2}, then Vol
{\bf 3,4} and then Vol {\bf 1} of Darboux \cite{D1}) Bianchi
\cite{B1},\cite{B2} and Darboux \cite{D1},\cite{D11} and I was
introduced to the theory of deformations of quadrics, one of the
crowning achievements of the classical geometers. The main
obstacle in understanding the classical geometers is the
consistent use of imaginary numbers, as these numbers have been
dropped from current courses dealing with the geometry of
sub-manifolds in scalar product spaces. The current approach is to
split this geometry into Riemannian and Lorentz, but this comes at
a price: geometric connections between the two are lost and
similarities remain only at the analytic level. Also results about
{\it non-totally real} sub-manifolds (sub-manifolds of
$\mathbb{C}^n$ which do not lie in totally real subspaces; however
they may have real linear element) are lost. Consider for example
the statement from Darboux (\cite{D11},\S\ 169): {\it 'A surface
of revolution can not only be deformed to a line, but further to
an isotropic line'} (thus a $2$-dimensional linear element can
lose not only one dimension, but further two). Such a statement
may seem false if one adheres to the strict current rule that
applicability implies same dimension as manifolds, but if one
replaces 'surface', 'line' and 'isotropic line' with the
$2$-dimensional collection of their tangent planes (with the
points of tangency highlighted), it begins to make sense at the
differential level, at least from a dimensional point of view. Now
the condition that a $2$-dimensional family of facets (a {\it
facet} is a pair of a point and a plane passing through that
point) is the collection of the tangent planes of a sub-manifold
(leaf) does not distinguish between the cases when this
sub-manifold is $0,\ 1$ or $2$-dimensional.

The deformation of a CGC $-1$ surface to a line explains the
generation of the {\it real pseudo-sphere} (surface of revolution
with a tractrix as generating curve) and {\it Dini surfaces}
(helicoids with same type of generating curve) from a line via a B
transformation; at the level of the sine-Gordon equation this is
equivalent to generating (via a B transformation) $1$-soliton
solutions of the sine-Gordon equation from the vacuum soliton.

Note the metric-projective nature of the {\it tractrix}: it is the
real curve individuated by the requirement that its real tangents
meet a given line at constant distance from the tangency point; it
is obtain by slowly pulling a cart along a line if initially the
cart axis is perpendicular on the line.

The deformation of a CGC $-1$ surface to an isotropic line
provides a simple point of view on the B transformation, due to
Lie and apparent at the level of {\it confocal} ('with same foci')
pseudo-spheres.  Lie's ansatz is very simple: because the B
transformation is of a general nature (independent of the shape of
the seed), it must exist {\it not} only when the seed is a CGC
$-1$ in $\mathbb{R}^3$, {\it but} also when the seed coincides
with the pseudo-sphere, in which case the leaves degenerate to
isotropic (imaginary) rulings on a confocal pseudo-sphere (note
that there is no formulation of the Lie ansatz at the level of the
sine-Gordon equation, since it deals with imaginary seed). This
simplification allowed Bianchi to generalize the B transformation
of CGC $-1$ surfaces to the B transformation of deformations of
quadrics and it allows us to simplify the denomination {\it B
transformation of deformations of quadrics} to {\it B
transformation of quadrics}.

Thus the Bianchi-Lie ansatz states that the leaves (B transforms
of the seed) are rulings of a ruling family on a confocal quadric
when the seed coincides with the considered quadric.

In the famous letter {\it The Method}, sent to Eratosthenes of the
Library-University of Alexandria and lost from the
XIII$^{\mathrm{th}}$ century until 1906 (the year of Bianchi's
discovery, so it was unknown to Bianchi and Lie), Archimedes
states: {\it '... certain things first became clear to me by a
mechanical method, although they had to be proved by geometry
afterwards because their investigation by the said method did not
furnish an actual proof. But it is of course easier, when we have
previously acquired, by the method, some knowledge of the
questions, to supply the proof than it is to find it without any
previous knowledge'}.

Note that the method at the level of points of facets was known to
Bianchi and Lie; however, they had never used the full method, at
the level of the planes of the facets too (this is where
Archimedes enters into the picture), so Archimedes' part is the
essential ingredient for generalizations in other settings
(especially in higher dimensions), since it contains more
information. Thus if one assumes Theorem I of Bianchi's theory of
deformations of quadrics a-priori to be true and to be the
metric-projective generalization of Lie's approach, then by
applying the Bianchi-Lie ansatz (the {\it 'mechanical method'} in
question) both at the level of points of facets and planes of
facets one obtains in a natural and geometric way the necessary
algebraic identities needed to prove Theorem I. The fact that
Archimedes' opinion is the main ingredient in supplying the proof
of Theorem I is confirmed by the fact that Bianchi's original
proof is much more complicated, because it was initially
discovered first at the analytic level and for the composition of
two B transformations, so Bianchi had to retrace the steps back to
the basic B transformation at the analytic level and guess the
geometric picture; it took him 6 years of work on several
particular cases to do so (of course Bianchi simultaneously worked
on many other projects during this time) and thus it lacks the
{\it 'of course easier'} ingredient.

Of course  Bianchi's proof is commendable for the fact that it
solved an open problem at that time; just like Archimedes Bianchi
showed that it works on several examples (again in each case
Bianchi did computations from scratch, as was the standard of
proof of the time), so one can safely deduct that it works for all
quadrics (deformations of certain particular quadrics were
previously found, so Bianchi's method was not seen as necessary in
those cases).

However, without the power of a professional with good knowledge
of the tools of the trade, the faith in finding a result at the
end of long computations and the geometric intuition of Bianchi's,
this could have been an open problem even today, as few years
after its completion this type of problems fell into forgetfulness
for more than 50 years. It is the word superhuman ({\it
'sovrumana'}) that was chosen by Bianchi's contemporaries to
describe his efforts and geometric intuition.

We shall call this method the {\it Archimedes-Bianchi-Lie} (ABL)
method.

Note that although Archimedes and Bianchi-Lie have dealt with
different problems, their approach was the same: $P\Rightarrow Q$
with $P$ being either {\it 'The area of a segment of parabola is
an infinite sum of areas of lines'} or {\it 'A line is a
deformation of a quadric'}. Such sentences $P$ were not accepted
as true according to the standard of proof of the times, but they
had valid relevant consequences $Q$ which elegantly solved
problems not solvable with other methods of those times; if only
for this fact the standard of proof should change so as to accept
these statements as true.

At some point I thought the {\it 'mechanical method'} in question
to be the rolling of the seed back on the considered quadric (used
by Bianchi in his arguments and clearly a mechanical method),
which gives a natural geometric explanation of the necessity of
the existence of the {\it rigid motion provided by the Ivory
affinity} (RMPIA) (this is how Bianchi found the {\it
applicability correspondence provided by the Ivory affinity}
(ACPIA)) and of the reflection in the tangent bundle of the
considered quadric (one can roll a surface on both sides of one of
its deformations). These were implicit in Bianchi (122,\S\ 11) and
well known to me, but Archimedes was still not satisfied, because
these were not the {\it 'certain things'} which {\it 'first became
clear'} and they were geometric arguments which remained to be
confirmed not by geometry, but by an analytical calculus (the
double reduction ad absurdum, an analytical calculus today by
means of the $\del-\ep$ trick, was a geometric argument during
Archimedes' times). Thus Archimedes further required that to the
rolling of the seed back on the considered quadric one must also
apply the ansatz that the seed coincides with the given quadric,
in order to get the a-priori nonsensical {\it 'certain things'}
which {\it 'first became clear'}, only {\it 'to be proved by
geometry afterwards because their investigation by the said method
did not furnish an actual proof'}.

The classical geometers frequently used simple geometric arguments
to justify statements as Darboux's from above (which provides
itself a simple geometric picture of singular solutions of a
certain differential system): Darboux in \cite{D1},\cite{D11}
provides such arguments for almost any analytic computation; they
are always shorter and more intuitive. However, without Cartan's
simplifying notation (which led him to exterior calculus), their
computations of moving frames took pages (as required by the
standard of proof of the time) and thus the investigation of
higher dimensional problems lagged behind; for example Bianchi
takes time in (122,\S\ 11) to write on a full page all components
of the RMPIA, although this explicit formulation is never used
(however the components of the rotation look strangely close to
certain quadratic quantities with direct usefulness).

Another feature of the standard of proof of the time (namely
explicit solutions) had a positive influence. It is one of the
possible explanations of the success of the classical geometers in
the study of integrable systems, as these are amenable to
algebraic formulations of their solutions after a simpler
auxiliary differential system associated to a given local solution
of the initial integrable system is presumed solved. The algebraic
formulation of solutions is possible due to the {\it Bianchi
Permutability Theorem} (BPT), initially developed by Bianchi for
the B transformation of CGC $-1$ surfaces in (\cite{B2},Vol {\bf
5},(46)), but actively pursued by him in various settings
throughout his career (there are also other results and methods,
like the one mentioned in the introduction of Terng-Uhlenbeck
\cite{TU1} and due to Darboux).

With the development of calculus experts realized that most
differential systems do not exhibit such behavior and they focused
their attention on the qualitative behavior. Integrable systems
were revived in the second part of the XX$^{\mathrm{th}}$ century,
due mostly to the KdV equation and other problems from physics
(see Rogers-Schieff \cite{RS1}).

Because of the BPT Bianchi was able to put a local group action
structure on the space of local solutions of an integrable system;
thus he 'added' (via a non-linear commutative algebraic procedure,
hence the {\it permutability} denomination) two solutions and
obtained a new one. Although the classical geometers did not
further investigate this group structure (see Terng-Uhlenbeck
\cite{TU1}), Bianchi was aware of the group structure at least for
{\it $n$-solitons} of the sine-Gordon equation. The {\it
$n$-solitons} of an integrable system are obtained by the iterated
application $n$ times of the B transformation to the vacuum
soliton (practically one applies the BPT $n-1$ times to the vacuum
soliton; Bianchi calls iterations of the BPT {\it moving
M\"{o}bius configurations}). In the case of the sine-Gordon
equation the above group can be identified with its (global)
action on the vacuum soliton; its elements admit formulations as
rational functions of exponential functions (see (\cite{B2},Vol
{\bf 5},(46,\S\ VI))). When the initial value data of the
auxiliary differential system are taken into account, the
permutability loses its commutativity connotation and the BPT for
the B transformation of the pseudo-sphere is explained by two
different ways of factoring quadratic terms in the group of
rational maps on the unit sphere
$\mathbb{S}^2=\mathbb{C}\cup\{\infty\}$; the factors are
M\"{o}bius transformations (see Terng-Uhlenbeck \cite{TU1}).

These notes began with the purpose of understanding the theory of
deformations of quadrics; in the process the proofs were reduced
and simplified, thus we hope making it easier available to a wider
audience. For example the proof of the existence of the B
transformation of quadrics and of the ACPIA does not significantly
surpass both in length and difficulty the actual statement and its
prerequisites. The statements (with no major modifications) are
due to Bianchi and their beauty, clarity and simplicity provided
the motivation to go through Bianchi's original proofs. There are
other classical monographs on the theory of deformations of
quadrics, but they are incomplete or the approach is different
than Bianchi's: for example Calapso \cite{C1} or the last chapter
of Eisenhart \cite{E2}.

While all main ideas are due to the classical geometers
Archimedes, Bianchi, Darboux, Lie, etc, at the computational level
there is something new: a complete, simplified and unified
treatment of all quadrics in $\mathbb{C}^3$. For example Bianchi
takes a quadric in $\mathbb{C}^3$, intersects it with various
totally real subspaces to get real $2$-dimensional (in)definite
linear elements and then builds the B transformation in various
totally real subspaces separately for each case (providing enough
supporting examples according to the required standard of proof of
the time). We provide a unified discussion of the B transformation
and postpone the necessary splitting into totally real cases, as
an application, to the very end (a similar discussion takes place
in Terng-Uhlenbeck \cite{TU1}). Because the initial surface
(having $2,3$ or $4$ real dimensions) in the considered quadric or
its deformations may be non-totally real, the most general
configurations are obtained. In particular if the linear element
is real, then the surfaces must have real dimension $2$, so we
obtain deformations (which may be non-totally real) of totally
real surfaces in quadrics (a surface in a quadric and having real
linear element must be totally real).

In what follows we shall deal with quantities which may not be
defined everywhere; instead of always pointing out the domain
where the computations are true we shall mostly just assume them
to be true on their domain of existence, without any further
details (we discuss only the local theory). We shall require as
many derivatives as needed; in fact since most of the current
literature on integrable systems deals with solutions in the
distribution sense, all computations can be interpreted in the
distribution sense (there may be generalizations even in other
settings, like $p$-adic (non-Archimedean) calculus, but at this
point I am not familiar with such proceedings; note however that
the correct interpretation of calculus in Leibnitz's notation
gives for free the correct discretization). A global theory (if
developed) must include singularities (along points or curves) and
thus must probably be developed along lines similar to other
integrable systems. While it is essentially the trend of the
current geometers to prove global results, most such results are
about rigidity and for global deformations of complete surfaces
such results are few (see Spivak (\cite{S3}, vol {\bf 5})); thus
one can safely assume that in order to have a consistent set of
objects to study singularities must be allowed. In most current
theories singularities are to be avoided; here they provide
degeneracies, due to which certain general differential equations
(with no method of integration) degenerate to ones that can be
easier integrated. Again one can guess here the effect of
Archimedes' method and probably a geometric solution will be
either easier to obtain or, in the contrary case, the geometric
explanation of the solution obtained by analytic computations will
easily surpass and simplify the analytic computations.

As to the $n$-solitons of integrable systems known so far, they
involve difficult formulae, usually involving (hyperbolic)
trigonometric, exponential and even elliptic functions; the
geometric projective picture should serve as an easier explanation
of these strange laws governing these solitons, just like Jacobi's
simple geometric interpretation on confocal quadrics of Abel's
Theorem (see Darboux (\cite{D1},\S\ 464) and (\cite{D11}, Note
I)).

After the application of the ABL method the computations naturally
split into {\it static} (algebraic) computations on quadrics
confocal to the given one and {\it moving} (differential)
computations in its tangent bundle. The flat connection form (a
feature common to all integrable systems and always coming in a
spectral $1$-dimensional family) is provided by {\it rolling}: two
rollable (applicable) surfaces can be rolled (applied) one onto
the other, such that at any instant they meet tangentially and
with same differential at the tangency point. Conversely, if two
surfaces can be rolled one onto the other, then they must be
applicable with the applicability correspondence being given by
the correspondence between points whose tangent spaces coincide at
some instant in the rolling. Note that it is less known today that
Levi-Civita introduced the {\it parallel transport} denomination
since it is the transfer by rolling of the usual Euclidean
parallel transport in a plane, as the surface is rolled only along
the given curve (along which the parallel transport needs to be
done) on a fixed plane (thus as a plane rolls along one of its
lines on a surface it describes a geodesic on the surface and all
geodesics are born this way). The local minimizing property of
geodesics is immediate from the geometric definition, since the
shortest distance becomes at the infinitesimal level heading in
the direction of shortest distance and by rolling this heading
sticks (transfers) from the plane to the surface and conversely.
An ant on an ellipsoid would do the same as what people did from
the beginnings of mankind: it would choose a point far away as
reference or it would look at the heavens (tangent space) and it
would thus linearize its problem, rolling the heavens with it
(changing the point of reference) as it proceeds on its journey.

Currently one would be inclined to take these geometric
definitions which inspired analytic consequences as consequences
of the analytic consequences (taken currently as definitions,
instead of as confirmations), just as a trained mathematician
would rather give an $\mathbf{L}^2$ proof of the discrete
Cauchy-Schwartz inequality instead of the longer (but at
high-school level) proof using induction or just as the Calculus
II undergraduate student would give a L'Hospital proof to
$\lim_{x\rightarrow 0}\frac{\sin(x)}{x}=1$.

Just as higher dimensional calculus is essentially $1$-dimensional
calculus in all directions plus the Clairaut integrability
condition (commuting of mixed derivatives or the d-wedge
condition), higher dimensional rolling (including of sub-manifolds
in higher dimensional symmetric spaces) is essentially
$1$-dimensional rolling in all possible directions plus the same
integrability condition.

As we have only $1$-dimensional time to our imagination, it is
somewhat difficult to imagine the $2$-dimensional rolling of two
applicable surfaces, so we have to imagine $1$-dimensional rolling
along all possible corresponding curves (which will thus admit an
arc-length correspondence; since a correspondence of two surfaces
which induces arc-length correspondence of arbitrary corresponding
curves must be isometry, existence of rolling implies isometry).
Since Leibnitz's notation does not have a memory of dimension
(which is dealt with by another condition, namely {\it the
independence condition}), it is the tool best suited to attack the
rolling problem and the easy $1$-dimensional rolling problem
easily introduces the notation needed for higher dimensional
rolling problems.

By application of the ABL method certain valid relevant algebraic
identities naturally appear (from a geometric point of view)
within the static part.

All necessary identities of the moving part boil down to these
valid relevant algebraic identities of the static part (other
terms of the differential identities dissolve into nothingness
because of the flat connection form). A possible justification of
this behavior (and thus the miracle of the ABL method) appears
when one considers {\it discrete deformations of quadrics} (DDQ);
in 1996 Bobenko and Pinkall have introduced in
\cite{BP1},\cite{BP2} discrete isothermic and CGC $-1$ surfaces,
although Newton was the first one to observe the discretization of
integrable systems and of their solutions when he provided a
mostly Euclidean explanation of Kepler's Laws of planetary motion.
Discrete surfaces must obey two simple principles: on one hand
their determination from the discrete Cauchy data involves only
algebraic computations (precisely the algebraic computations
required by the BPT; thus the last important step of the theory of
deformations of quadrics becomes a basic axiom for DDQ) and on the
other hand, just as their name suggests, discrete surfaces better
approximate their corresponding surfaces as the discretization
becomes finer; therefore they algebraically encode the
differential structure of their corresponding surfaces (and thus
the most important step of DDQ provides the building blocks of the
theory of deformations of quadrics). It is thus no accident that
Bianchi had in fact developed discrete versions of various
integrable systems and of their solutions when he developed moving
M\"{o}bius configurations (iterations of the BPT) a hundred years
ago, but without being aware of this fact. One can even develop
the discrete B transformation for discrete solutions of discrete
integrable systems! (the simplest trick to do this is to preserve
Leibnitz's notation but to interpret it differently). The
infinitesimal becomes the discrete finite by iterating a
$\mathbb{Z}\times\mathbb{Z}$ lattice of B transformations on an
initial seed and picking a seed point on that seed and its
corresponding $\mathbb{Z}\times\mathbb{Z}$ heirs; the finite
becomes the infinitesimal by refining the
$\mathbb{Z}\times\mathbb{Z}$ lattice of a DDQ (while preserving
the subjacent curves of the discrete Cauchy data; since they must
become asymptotes they have prescribed torsion by Enneper's). The
discrete version of the Clairaut equality $(f_x)_y=(f_y)_x$ for
$f=f(x,y)$ implies the BPT and as one refines the
$\mathbb{Z}\times\mathbb{Z}$ lattice of B transforms the BPT
becomes infinitesimal, so it precisely describes the above
equality! Thus the theory of DDQ is the proof {\it 'such as the
ancients required'}, provided that all its algebraic computations
are put in a synthetic language.

Why one would do such a thing? Basically Archimedes discovered in
{\it The Method} integral calculus for some conics and quadrics of
revolution (he discovered almost all important formulae;  probably
there are other formulae for quadrics which can be discovered with
the same method) while Apollonius essentially discovered
differential calculus for conics in {\it Conics} (part of the
second half was available in Latin beginning unfortunately only
with 1710; among the essentially differential notions introduced
by Apollonius was the evolute of a conic). Up to Kepler some more
results were added; these allowed him to do the discrete
computations for his Laws. Note that basically all methods and
results up to Leibnitz and Newton are useful not only from a
historical point of view, but also for the fact that these hint an
integrable system behavior (take for example Roberval's simple
Euclidean arguments in computing tangents by kinematic arguments
and the area of the cycloid; it is probable due to the fact that
nature was able to do computations for the cycloid using only
simple mechanical methods and simple Euclidean arguments that it
chose it as the solution for the {\it brachistochrone} and {\it
tautochrone} problems (the resolution of the former among others
by Newton was the first instance of {\it Calculus of Variations}
or equivalently infinite dimensional calculus) and the masterly
unification of the pendulum with the tautochrone by Huygens when
he showed that the evolute of a cycloid is a cycloid). Beginning
with Leibnitz and Newton calculus was the main theoretical tool,
although discrete mathematics remained the main practical tool up
to including today (because measurements, observations and
computer memories are discrete). Beginning with the second part of
the XX$^{\mathrm{th}}$ century the use of optimal discrete
mathematics for problems with not so obvious solutions began to be
again important, as more difficult problems (such as integrable
systems) arose; this was the main motivation of Bobenko and
Pinkall for their 1996 articles. Thus a synthetic language based
on axioms and counting their use is a way to obtain not only the
optimal proof, but also to implement the optimal discrete
mathematics.

The static part requires basic knowledge of linear algebra and
projective geometry; we use for example the Menelaus theorem
(which is essentially a co-cycle theorem) to prove the
associativity of the group action whose commutative 'addition' law
is provided by the BPT; see also Arnold's observation from
\cite{A1} about the link between the Jacobi identity
(infinitesimal associativity, in itself a co-cycle theorem) and
the intersection of the heights of a triangle. Note that the
Jacobi identity for the R-matrix theory (Lie bracket deformation
$[x,y]_R:=[Rx,y]+[x,Ry]$ with $R$ linear of the Lie bracket $[,]$;
see Burstall, Ferus, Pedit and Pinkall (\cite{BFPP},\S\ 3)) may be
the possible analogy of the Menelaus theorem, as I have noticed
the $R$-matrix behavior (and consequently I was able to simplify
the notation of the algebraic computations of the BPT) before
seeing the Menelaus theorem. What are the chances of reading a
book of history of mathematics, see a theorem, remember to have
seen it in grade school and immediately realize that it is what
one needs to answer the current questions?

In the first issue of the {\it Bulletin} on 2006 I read something
about the Arnold-Jordan canonical form of matrix and I wanted to
see if Arnold had already done what I call the symmetric Jordan
canonical form; this is how I discovered \cite{A1}.

So the R-matrix theory is probably the necessary tool in higher
dimensions and the Jacobi identity of $[,]_R$ (which imposes
certain conditions on R) is essentially a generalization of the
Menelaus theorem in other settings.

For the moving part we shall only use rudiments of the absolute
(Ricci) calculus of tensors, working mostly with $(0,1)$ and
$(0,2)$ tensors. Since we would like $d^2$ to mean tensorial
second derivative, we shall use the notation $d\wedge$ (d-wedge)
for alternative (exterior) derivative. Take for example a plane
$V\subset\mathbb{C}^3$ with normal $0\neq N\in\mathbb{C}^3,\
N\perp V$; if $|N|^2\neq 0$, then one can normalize $N$ to
$|N|=1$; otherwise $N$ is an isotropic vector: $|N|=0$, so $N\in
V$ and $<,>$ restricted to $V$ is degenerate; thus the second
fundamental form of most surfaces $x\subset\mathbb{C}^3$ is
defined almost everywhere as $N^Td^2x$ (Cartan never uses the
notation $d$ with the meaning $d^2=0$, as he probably wanted to
preserve Leibnitz's notation:
$\frac{d^2y}{dx^2}=\frac{d}{dx}(\frac{dy}{dx})$; see (\cite{C2},
Ch III,\S\ 26)).

Most classical integrable systems were discovered while studying
quadrics; thus they naturally admit, at least in particular cases,
(geometric) links to quadrics. For example Darboux was studying in
\cite{D3}-\cite{D5} rolling {\it congruences} generated by
isotropic rulings of quadrics when he found the D transformation
for special isothermic surfaces (a {\it congruence} in
$\mathbb{C}^3$ is a $2$-dimensional family of objects, presumed
lines unless otherwise stated, but they can also be arbitrary
curves or even surfaces; following Pl\"{u}cker a $3$-dimensional
family of objects in $\mathbb{C}^3$ is called {\it complex}). As
the quadric rolls on its applicable surface in a $2$-dimensional
fashion (assume that the applicable surface is not ruled, in which
case the rulings correspond under the applicability to a ruling
family of the quadric and the rolling takes place in a
$1$-dimensional fashion), one if its isotropic rulings describes a
$2$-dimensional family of isotropic lines; hence the denomination
'rolling congruence'. The natural question thus arose if this
transformation exists for general isothermic surfaces and Darboux
provided an affirmative answer. The auxiliary differential system
corresponding to the D transformation for isothermic surfaces
admits a quadratic prime integral, which may explain the link
between special isothermic surfaces and quadrics; for the link
between isothermic surfaces and quadrics we probably need quadrics
in higher dimensional spaces.

Based on the classical and current literature (see for example
Moser \cite{M1} or Rogers-Schieff \cite{RS1}), we believe that any
integrable system is bound to be linked at least in particular
cases to quadrics, even in infinite dimensional spaces, where a
spectral theory is well developed but without the corresponding
geometrization; in fact confocal quadrics as particular confocal
cyclides and {\it Dupin cyclides} (surfaces with circular lines of
curvature) clearly appear in Darboux \cite{D11} as important
examples and the notion of isoparametric sub-manifolds in finite
dimensional spaces has already been generalized to isoparametric
sub-manifolds in Hilbert spaces by Thorbergsson and collaborators.
Note also that cyclides in pentaspherical coordinates have
elliptic coordinates just like confocal quadrics (see Darboux
(\cite{D1},\S\ 437 and Livre IV, Ch XII)).

We also believe that the above formulation (split the computations
into static and moving by means of the ABL method) is generous
enough to be amenable to generalizations in other settings; in
particular we believe that the B transformation of deformations of
space forms in space forms and with flat normal bundle (Tenenblat,
Terng, Uhlenbeck and collaborators) can be generalized to the B
transformation of higher dimensional quadrics in symmetric spaces
and with flat normal bundle.

\subsection{Short historical perspective on quadrics}
\noindent

\noindent The interested reader is referred for more details to
Boyer \cite{B3} or the website \cite{W2} (our main sources for
this subsection up to the beginning of the XIX$^{\mathrm{th}}$
century), to G. Scorza's and G. Fubini's eulogies in part 1 of
Volume 1 of Bianchi \cite{B2} and to R. Calapso's introduction in
part 1 of Volume 4 of Bianchi \cite{B2} (our main sources for this
subsection in what concerns the history up to 1912 of the theory
of deformations of quadrics).

Further this perspective (and its references) is definitely biased
and incomplete, it being designed to suit a restricted point of
view; besides that it has a personal touch ({\it 'You want to make
an omelet, you got to break some eggs'}). Problems of the
XIX$^{\mathrm{th}}$ century such as Monge's sphere (for any
quadric there is a sphere and a moving orthonormal frame centered
on that sphere such that the hyperplanes of the orthonormal frame
are tangent to the given quadric) and Chasles's number of 3264
conics tangent to five given conics will not be discussed in
detail here.

The primary organization of this subsection is the chronological
one; the secondary one is focused on the subject at hand.

Quadrics are the simplest non-linear geometrical objects and thus
they arose early on the interest of geometers: the ancient Greek
geometry culminated with the work of Apollonius, Archimedes,
Euclid, Pappus and others in particular on the geometry of conics.

Archimedes is credited to have been close to a theory of integral
calculus: to find areas and volumes he used at the infinitesimal
level an ingenious balancing argument and simple Euclidean
arguments (the famous letter {\it The Method}, sent to
Eratosthenes of the Library-University of Alexandria and lost from
the XIII$^{\mathrm{th}}$ century until 1906); to confirm the
validity of the formulae thus obtained he used a double reductio
ad absurdum, as was the standard of proof at the time. The
discovery of {\it The Method} is important because it changed
Archimedes' reputation from being unclear on the details and
motivations of the methods employed in proving his results
(presumably with the purpose of enhancing his own reputation) to
being very generous in {\it disclosing the full method of
investigation} to Eratosthenes, {\it 'a capable scholar and a
prominent teacher of philosophy'}, so with the stated purpose that
the later should further investigate it and teach the curious
trick to his students.

To find the area (volume) of an object Archimedes cut it into
parallel lower dimensional slices (segments for area and plane
figures for volume), then he placed all slices with their center
of mass at one end of a balance and a simpler object at the other
end; if equilibrium is obtained the area (volume) of the initial
object can be found. The equilibrium is checked by comparing
slices of the original object with corresponding slices of the
simpler object and {\it not} only simple Euclidean geometric
arguments are needed, {\it but} the choice of {\it the correct
balance with fixed fulcrum} and of {\it the simpler Euclidean
object} is revealed by the simple geometric properties of the
slices of the initial object. So basically {\it The Method} of
Archimedes takes into account the {\it discrete} version of
computing the {\it mass} and {\it center of mass} of an object
from the knowledge of the {\it masses} and {\it centers of mass}
of its component pieces and assumes it also to apply as a limiting
case to the {\it infinite} (non-discrete) version of the same
statement; the use of simple Euclidean statements about conics
available at the time was realized with a balance and thus
interesting results about conics and quadrics of revolution were
added to the bag. Nowadays this would be proof enough (the slice
theorem or Cavalieri's Principle for parallel slices), but
Archimedes was unsure of the power of his method: he would set out
to prove its validity by thickening the slices to regain the
original measure (so he would go from the singular boundary of the
space of solutions to his discrete principles back into the meat
of the original space of solutions, where one can wiggle and still
remain in the meat; this is what actually happens with the theory
of deformations of quadrics in this notes and possibly with {\it
all} integrable systems: we are unsure of the power of Archimedes'
method and provide further explanation which may be considered
unnecessary in the future) and a double reductio ad absurdum, as
was the standard of proof at the time. Note that all through his
method pervasive is the use of the balance with fixed fulcrum and
left end where facets are {\it transferred and stuck with their
center with $\infty$ multiplicity} (the variable right end of the
balance is dictated by some of the variable mass of the object at
the fixed left end, as it transfers and sticks to the apparently
wrong unit of measure of length); since the balance works with any
infinitesimal configuration, it will work with the integral
version. One can look at the good metric properties of the Ivory
affinity as being the required prerequisites of simple Euclidean
relations by means of which one puts the quadric at the right end
of the balance and rulings on a confocal quadric as leaves at the
left one. Now the existence of seed provides a way to deform the
quadric, that is a way to rearrange the object at the right end of
the balance: the balance remains pervasive through this process
and allows the rearrangement of the facets at the other end in a
manner consistent with the initial arrangement (that is admitting
leaves) so the balance wins and its principles (clearly valid at
the infinitesimal level by means of the simple Euclidean
properties) remain valid after integration. The integration
process (proven by Archimedes to be correct by means of a double
reduction ad absurdum) is realized by means of flat connection
form: facets fit nicely at the infinitesimal level on the right;
as a consequence facets fit nicely on the left. The discrete
version of the theory of deformations of quadrics can be actually
put by means of a $\del-\ep$ argument (error estimate) to the form
of a double reduction ad absurdum.

With this method Archimedes managed to compute the area and center
of mass of the region bounded by a chord on a parabola, volumes of
segments of quadrics of revolution, centers of mass of hemispheres
and of the paraboloid of revolution. He also used a geometric
argument reminiscent of what we call today 'trapezoid rule' to
find the area bounded by a chord on a parabola; the argument does
not work for ellipses or hyperbolas, as they cannot be realized as
graphs of quadratic functions; however he was able to find the
area of the ellipse.

Apollonius is mostly known for his {\it Conics}, which began as a
project of collecting the results about conics from Euclid's {\it
Elements} and other sources of the time, but stated in a more
general setting and completed with corresponding proofs; in the
process it became a subject of interest on its own and completely
replaced all other sources, remaining the standard for almost two
millenniums, just as Euclid's {\it Elements}. While it seems that
Apollonius and his contemporaries were familiar with foci of
conics and some of their properties, these were considered and
referred to only indirectly at that time. It is Apollonius who
first introduced the {\it hyperbola} (something is in excess),
{\it parabola} (something is precisely) and {\it ellipse}
(something is missing) denominations took over since then in many
other domains (such as literature).

According to Boyer (\cite{B3}, IX,\S\ 5): {\it 'If survival is a
measure of quality, the Elements of Euclid and the Conics of
Appolonius were clearly the best works in their fields'}. It is
chiefly because of the {\it Conics} that Apollonius was known as
{\it 'The Great Geometer'} of the antiquity. Understanding this
treatise was regarded as the ultimate test an accomplished
geometer should pass; this may be the possible explanation of the
fact that the second half of the {\it Conics} has survived only
partially and in arabic (Jafar Muhammad ibn Musa ibn Shakir's
version) until 1710, when, belated and apparently useful only from
a historical point of view, Halley's first Latin translation
appeared. Note however that this translation could have been done
almost a century earlier, as the manuscript was available in one
of Europe's libraries, but administrative reasons have prevented
it.

According to a {\it 'Nova'} program on Archimedes' palimpsest
calculus could have probably been discovered at least one century
earlier with {\it The Mehod} of Archimedes available at the time;
consequently the small step for {\it 'a man' } in the
XX$^{\mathrm{th}}$ century could have been on Mars (Florin pointed
out also the other side of the coin, namely humankind could have
self-destructed in the XIX$^{\mathrm{th}}$ century).

Note also that the complete version of Apollonius's {\it Conics}
could have helped in early XVII$^{\mathrm{th}}$ century: among
other essentially differential notions first introduced by him
using only synthetic computations was the {\it evolute} of a
conic. When a line rolls on the planar {\it evolute} a point of
that line will describe an {\it involute}; thus the orthogonal
family of a family of plane lines are involutes; in our situation
the given conic is among the involutes. This picture admits
generalization to curves in space by bending the plane of the
curve with a cut made along the curve (although this is still not
the correct picture, as one needs two copies of the original
plane; one can make a contraption with few tangent lines attached
to a thin metal sheet (a one sided tubular neighborhood of a plane
curve)): the normals to a involute surface envelope focal surfaces
(evolutes), the vector fields induced by normals on an evolute
surface admit integral curves which are geodesics on the evolute
surface; conversely a $1$-dimensional family of geodesics on a
surface makes that surface into an evolute of another
$1$-dimensional family of involute surfaces and each geodesic and
its involutes (curves) is in the situation of a bent plane
picture. Another possible generalization would be from a
$1$-dimensional family of plane lines to an $(n-1)$-dimensional
family of lines in $\mathbb{C}^n$; it still admits $\le n-1$
evolutes (focal hyper-surfaces or caustics), but it would not
admit in general involutes, as the condition that such an
$(n-1)$-dimensional family of lines are normals to a hyper-surface
is special.

Pappus introduced the {\it focus-directrix} property of conics (in
the plane geometry all conics are the locus of points such that
the ratio of the distances to a given point and a given line is
constant; the given point becomes focus and the given line
directrix), but he is mostly known for a mechanical method of
finding volumes of solids of revolution, which requires only
knowledge of the area and center of mass of the generating plane
figure; as a limiting case one can obtain the area of a surface of
revolution from the length and center of mass of the generating
plane curve. Replacing the generating plane figure with its center
of mass and having as mass the area of the plane figure, the
moment about the axis of revolution does not change; as a limiting
case one can obtain the area of a surface of revolution from the
length and center of mass of the generating plane curve (note
again that the part of Archimede's method concerning the mass and
center of mass provides the correct explanation if the slices are
the positions of the original plane figure (curve) as it rotates
as instructed).

At that time they endured the criticism of their opponents as to
the usefulness of their results (see page xi of {\it Preface} to
T. L. Heath \cite{H1} and the frontispiece); according to
Apollonius {\it 'They are worthy of acceptance for the sake of the
demonstrations themselves, in the same way as we accept many other
things in mathematics for this and for no other reason'}.

Up to the XVII$^{\mathrm{th}}$ century the works of the wise
ancients were required reading and their open questions
constituted the main driving force in the development of
mathematics; according to Leibnitz  {\it 'He who understands
Archimedes and Apollonius will admire less the achievements of the
foremost men of later times' }*. \footnote{* This quote (as most
others of this subsection) was introduced in the early 2006
version of the notes, but at that time I thought it to seriously
apply only up to the XVIII$^{\mathrm{th}}$ century, since Leibnitz
was not aware of Archimedes' simple point of view from {\it The
Method} and also I was not aware of the range of applications of
Archimedes' ideas. {\it 'That's a dam' shame!'}}

It is in trying to come closer to God that people chose Euclid as
a model; thus Euclid's {\it Elements} were required reading for
any important man regardless of his preoccupations and it is no
accident that in {\it The Declaration of Independence} of 1776 one
can easily see the simplicity of the proof beginning with simple
axioms and easily drawing the corresponding conclusions, just as
the balance is a universal principle (any action will have as
effect a corresponding opposite reaction equal in size in a
certain measurement).

Tartaglia (The Stammerer; he had his face sliced open by The
French in his childhood) thought initially that the trajectory of
a cannon ball is a broken segment; later on he changed his mind to
a parabola. One could easily assume that Tartaglia was a
na\"{\i}ve mathematician, only that de discovered the formula for
the roots of a cubic equation, a thing above me without the help
of Google. From the beginnings of mankind and probably up to its
ends mathematicians are seeking sponsors and most sponsors are
interested in measuring for taxing purposes (this is how geometry
appeared in ancient Egypt and Gau\ss\ initially set out to make a
geodetic survey of Hanover before he had his enlightening
remarkable discussion; it is probably because of always noticing
the same simplifications when compiling discrete measurements of
angles and distances that he got the idea) and warfare so as to
enhance and preserve the right of taxation: the first successful
algebra book in Middle Aged Europe was successful due to its
appendix in which it presented the method of {\it cooking the
books}, so the first million is questionable not only for
successful businessmen.

Using Tycho Brahe's detailed astronomical observations, Kepler
found an important application:

{\it I The orbit of a planet is an ellipse with the Sun situated
in one of its foci.}

{\it II A line joining the planet and the Sun sweeps out equal
areas in equal times as the planet describes its orbit.}

{\it III For any two planets the ratio of the squares of their
periods is the same as the ratio of the cubes of the mean radii of
their orbits.}

Newton was able to explain this by the theory of universal
gravitational attraction, thus discovering probably one of the
first integrable system, for which Kepler's laws are conserved
quantities. One could actually recognize the first integrable
systems in some of the results of the wise ancients in the form of
laws of a general nature. For example Apollonius's result that the
tangent to a hyperbola cuts the asymptotes at points equidistant
from the point of tangency is clearly a metric-projective result
which easily individuates hyperbolas as the only real curves with
this property; it remains valid if one replaces the asymptotes
with another hyperbola having the same asymptotes (some such
hyperbolas are homothetic to the given one; the whole family is
the pencil of conics generated by the asymptotes and the line at
$\infty$ counted twice; thus all such conics are tangent to the
asymptotes at their corresponding point at $\infty$) and thus any
line cuts a hyperbola and its asymptotes equidistantly (another
result of Apollonius, which thus implies the previous statement).
This can be easily generalized (by taking the slopes of the
asymptotes to be conjugate purely imaginary) to homothetic
ellipses (which have imaginary asymptotes); moreover if one
considers the $y$ coordinate purely imaginary such that the
asymptotes have parametric equation $[t\ \pm
it]^T\in\mathbb{R}\times(i\mathbb{R})$, then the above mentioned
equidistant segments also preserve same value if the subjacent
lines remain tangent to the same pseudo-circle (which thus becomes
(isotropic) pseudo-tractrix); the equidistance gives the correct
explanation of the a-priori counterintuitive fact that a plane
with positive definite induced linear element in
$\mathbb{R}^2\times(i\mathbb{R})$ intersects any light (null) cone
along a circle; one can apply this plane to a plane in
$\mathbb{R}^3$ and actually see the circle. Darboux especially has
the frequent custom of intersecting the tangent planes of a
(presumably real) surface with a sphere of zero radius and
obtaining circles.

Note however that Apollonius's result is not fully
metric-projective, as it is not true for the pencil of homothetic
parabolas (with equation $\la y^2+xz=0$; thus it is the pencil
generated by the conic $xz=0$ and the conic $y^2=0$); this can be
easily seen by considering lines passing through the origin.
Moreover the orthogonal curves of such pencils of conics are
(excepting the homothetic parabolas) seldom pencils of homothetic
conics. Thus one can safely deduct that an interesting
metric-projective notion valid for all conics must involve more
than a pencil (every pencil of conics contains two degenerate
ones) and the fact that the singular conic being a double line
must be the line at $\infty$ for the trick to work will become
clear later, since the information about the Euclidean scalar
product is projectively encoded by a quadric on this line (thus
two points).

Note that the pseudo-tractrix $\Psi=\Psi(s):=[\sinh(s)\
i\cosh(s)]^T$ and the tractrix $T=T(s):=[e^{-s}\
\int_0^s\sqrt{1-e^{-2u}}du]^T$ have the same metric definition, in
the sense that one can roll the pseudo-tractrix on the tractrix
such that the isotropic line $\Psi+\dot\Psi$ is captured and
transported by the rolling, thus becoming the line $T+\dot T$. The
rolling must be a curve
$(R,t)=(R(s),t(s))\subset\mathbf{O}_2(\mathbb{C})\ltimes\mathbb{C}^2$
such that $(R,t)(\Psi,d\Psi):=(R\Psi+t,Rd\psi)=(T,dT)$. Now the
rotation $R$ is determined up to a reflection in the tangent
bundle (the symmetry of the normal bundle) by $Rd\Psi=dT$; the
translation $t:=T-R\Psi$ satisfies $dt=-dR\Psi$. It thus becomes
straightforward that $R(\Psi+\la\dot\Psi)+t=T+\la\dot T,\ \forall
\la=\la(s)\subset\mathbb{C}$.

As an indication of the standard of proof of the time (or
according to the historian who must cover all possibilities, in
order to avoid controversy), Newton kept the use of calculus at a
minimum, using, in the spirit of Euclid, mostly geometric
arguments. A big part of his {\it 'Principia...'} deals with
properties of conics, his chosen argument for the Pappus four-line
locus problem being {\it 'not an analytical calculus but a
geometrical composition, such as the ancients required'}. Again
the historian says that this may be an indirect attack on
Descartes and his methods, who had previously solved this problem
of Pappus, in the process revolutionizing mathematics as we know
it today by the introduction of coordinates, although again the
historian says that Apollonius was the first one to have
introduced coordinates and formulae for changes of coordinates
(invariant computations) at a synthetic level.

According to Newton {\it 'If I have seen farther than Descartes,
it is because I have stood on the shoulders of giants'}, as he
probably referred to his predecessors (up to the first half of the
XX$^{\mathrm{th}}$ century one can find in the classical geometers
direct references and pride in the lineage up to including the
wise ancients).

Currently calculus as developed by Newton in {\it 'The Method of
Fluxions'} is mostly replaced by Leibnitz's elegant notation, so
Newton remained with the reputation among physicists and Leibnitz
among mathematicians (because in physics one needs computations in
functions depending on the time variable; a mathematician may find
oneself in trouble when one needs the introduction of the time
variable, since experiment and observation come foremost and
analytic confirmations later). The split that occurred when
Newton, following {\it' the bitter quarrel between adherents of
the two men' } decided to delete the reference to Leibnitz in the
third edition of his {\it 'Principia...'} was to remain on and off
for centuries and even present today (see Arnold \cite{A1}).

The fact that Newton's methods were not the correct ones (no
interesting mathematics happened in England until the
XIX$^{\mathrm{th}}$ century) is interpreted by the historian as an
indication of the superiority of Leibnitz's notation (which was
widely used in Europe).

But if it is the duty of the historian to appreciate the value of
a book by its influence, it is the duty of the current
mathematician to read the books of historical importance pertinent
to his restricted field of expertise and to appreciate not only
the value of the set of notations used, but also the value of the
ideas and motivations behind those notations. It is not excluded
that behind Newton's ideas one can see a symplectic structure and
his struggle to define it, as such structures naturally appear
from physics and Hamilton (the same with the hamiltonian and the
four quaternions $1,i,j,k$, named after the four groups of four
Roman soldiers each who guarded at some point St. Peter) was the
first important scientist in England after Newton to produce
significant research. The precursor of Cartan's exterior calculus
and higher dimensional moving frames technique is clearly present
in the classical French geometers.

The standard of proof for a differential problem decided by Newton
was to remain for two centuries: consistent use of Euclidean
geometrical arguments and finding solutions in closed form (that
is involving mostly algebraic and differential computations;
quadratures are used only if necessary). By referring to simple
Euclidean geometric arguments not only one can easily see what is
the course of investigation that must be followed and complicated
analytic computations are easily bypassed, but criticism as to the
soundness of the argument is avoided (for example Kepler found the
recent discovery of logarithms by Napier a great help for his
complicated computations; when advised not to trust logarithms as
being a contraption with possible hidden contradictions he set out
to give a Euclidean proof of the soundness of logarithms and their
properties).

By the mid-XIX$^{\mathrm{th}}$ century a crisis has reached
Euclid: the discovery of non-Euclidean geometries (although this
particular rupture was produced by renouncing a single axiom). It
seems that Gau\ss\ had previously discovered such contraptions,
but did not dare stir controversy (as it is the case with most
others contraptions of the XIX$^{\mathrm{th}}$ century). This
particular controversy was somehow closed when Beltrami showed
that the hyperbolic $2$-dimensional geometry can be realized as
the real pseudo-sphere, so ultimately the $3$-dimensional
Euclidean space contained the crisis in $2$-dimensions (one can
easily generalize this to a piece of $\mathbf{H}^n(\mathbb{R})$
being realized as a 'real pseudo-sphere' in $\mathbb{R}^{2n-1}$;
see Tenenblat-Terng \cite{TT1}; it is Bianchi who showed first
that space forms of same curvature and dimension are isometric).

However more abstract linear elements appeared: beginning with an
object in a surrounding Euclidean space, one records the induced
linear element, one deletes the surrounding ambient space and then
one studies only the properties of the linear elements with
properties as those obtained as above (Riemann). But one can
safely assume the existence of the surrounding Euclidean space in
order to get a good picture, as the only serious obstruction to
the existence of the Euclidean surrounding space for these
abstract linear elements was the difficulty of the corresponding
differential equations; thus one can work at the level of the
linear element with {\it parameterizations} instead of with {\it
maps}.

This is a major rupture with Euclid's ideas, as Euclid was
confined only to the infinitesimal world; however even in this
case general laws and structures were inspired and best explained
by a simple Euclidean picture at the infinitesimal level.
Consequently most of them began the long difficult journey back to
the finite world.

Although by mid-XX$^{\mathrm{th}}$ century this crisis of abstract
notions gave signs of being resolved, the current trends in
teaching geometry have yet to recognize this fact (see Arnold
\cite{A1}).

One cannot understand the classical geometers of the
XIX$^{\mathrm{th}}$-XX$^{\mathrm{th}}$ centuries if one pays
attention only to the two fundamental forms (the linear element
and the shape) and to the normal connection (the geometry of the
infinitesimally close surrounding space), since, in the spirit of
Euclid, most geometric constructions of the classical geometers
take place in a surrounding space and mostly at a respectable
finite distance from the considered surface and the analytic
version of such computations is meaningless at best at the first
look (for this reason Bianchi had to struggle for 6 years to find
the theory of deformations of quadrics, since equations initially
discovered for particular cases behave better at the analytic
level, so one naturally assumes that the same thing will happen
for arbitrary quadrics: he had to guess the geometric picture at
the analytic level).

Of course one can probably translate the B transformation to a
language of abstract bundles over abstract surfaces, using the
$\na$ notation for the Levi-Civita connection, the symmetric shape
operator instead of the second fundamental form, jets, etc, but
that will not make it more illuminating, as even for representing
such things we need to go back to the $3$-dimensional space and
draw a picture; moreover a quadric is naturally given by a
quadratic equation in a surrounding space and its confocal family
is not naturally defined in a bundle.

Ever since Cartan has shown that the geometry of sub-manifolds in
symmetric spaces can be attacked only with Leibnitz's notation and
with moving frames, I did not see any other simpler point of view.
One needs not the notation $x_\ast$ for $x$ an immersion of a
manifold in a symmetric space, when one can use instead $dx$ or
$x^{-1}dx$, etc for $x$ being the actual sub-manifold in the
symmetric space (or one of its representatives points and thus
allowed to move only on the sub-manifold). Consequently we need
local parametrizations of $x$ and not local maps on $x$.

Note however that for arbitrary $2$-dimensional quadrics we cannot
use orthonormal frames, as other frames (rulings) of a projective
nature are naturally chosen by the static picture; probably for
higher dimensional problems we need frames which are projective on
quadrics and orthonormal in the normal bundle.

In fact one can find in the classical geometers extensive
elementary treatments of various ideas and deep results a-priori
known to us to be of a recent nature, like comparison geometry,
geodesic flow and its rigidity (which can be extended so as to
include a-priori non-arc-length correspondence of geodesics or
just at the infinitesimal level of the flow), minimal surfaces,
Dupin's cyclides as precursors of isoparametric sub-manifolds,
etc: current monographs should at least mention these facts if
they choose not to cover them. (Un)knowingly ignoring already
known and established important facts will have as only effect the
perpetual rediscovery of same and ubiquitous simple principles,
disguised in an apparent aura of genius, modernity, rigorousness
and firm footing. The very elementary, intuitive, unsound
character and a-priori unwarranted generality of ideas and results
of the classical geometers constitutes in fact the essence of the
broad application of their ideas and results: decades before
Einstein's theory of (special) relativity the classical geometers
investigated and found interesting applications and properties of
imaginaries in geometry, only to be apologetic of the fact that
they did not see physical applications other than to apply various
tricks and restrict the whole discussion back to the real domain.
Darboux especially has the habit of applying a homothety when
discussing non-zero CGC surfaces and reducing the argument to CGC
$-1$ surfaces. Minding essentially did the same for geodesic
triangles on spheres and found hyperbolic trigonometric formulae
for geodesic triangles on pseudo-spheres (thus non-compact
continuous groups of symmetries and corresponding symmetric spaces
as counterparts of the compact ones); see (Darboux (\cite{D1},\S\
666)). Current curvature rigidity results have in general analytic
Riemanian manifolds with big continuous groups of symmetries and
other important extra-structures as end objects of investigation,
but the techniques are clearly split between positive and negative
curvature, so one can apply complexification and possibly restrict
to interesting totally real cases or methods of investigation.

Einstein's theory of general relativity as it appears in
Levi-Civita \cite{L1} is still very down to earth (he never
applied the abstractions, as the ideas of {\it 'Gau\ss, Riemann,
Christoffel and Ricci'} used by him were still very down to
earth). It is in trying to explain Einstein's ideas that
Levi-Civita came up with his connection and parallel transport, as
the surrounding Euclidean space was non-existent, so Levi-Civita
stuck the Euclidean picture to a sub-manifold and transferred it
to the abstract space; another motivation of Levi-Civita's was
mechanical problems. All interesting problems have elementary
formulations (at least in particular cases, such as the special
theory of relativity being the infinitesimal flat picture of the
general theory of relativity; when one tries to integrate this
infinitesimal picture one {\it naturally} gets the general theory)
and these elementary formulations should be presented first within
their simpler frame as a motivation for further generalizations;
further each abstract notion should come with an interesting
generalization for which the previous simpler apparatus was not
generous enough to resolve (like the Lebesque integral being the
natural generalization of the Riemann integral and built around
the necessary natural axiom that limits commute with integrals).
Most of the time the correct tools needed to attack the more
general problem will turn out to be just the correct
generalizations of the elementary tools used to attack the
elementary problem; otherwise probably the general problem is a
non-interesting one. There are rare exceptions and infinite
dimensional spaces is one of them; however even here there are few
basic results (in Banach, Hilbert and locally convex spaces) and
most interesting applications boil down to these basic results and
a-priori estimates; the interesting ones which do not yet boil
down to these are still currently shown to actually do; thus
methods to be taught to students should mainly focus on juicy
consistent applications of these principles of a general nature
instead of repeating dry abstract foundational issues for most
part of the course and to present few elementary cases (usually
involving the Euclidean space, (pseudo)-spheres and complex,
quaternionic and Cayley projective (hyperbolic) spaces and
amenable to elementary resolutions, like Desargues's theorem
proving that there is no $\mathbf{P}^n(\mathbb{C}a)$ for $n>2$ or
Pappus's theorem not holding in $\mathbf{P}^n(\mathbb{H})$) as
examples of culminations of the abstract theories. I should
mention that I was not the best (at mathematics) of my generation
in my county; one colleague and friend went to the International
Mathematics Olympiad (he is from the same city as the only other
ialomi\c{t}ean to have ever gone to an IMO), but he declared to me
one year in college: {\it 'Ion, I do not see the ideas. Where are
the ideas?'}. I did not know at the time what he was talking
about, as I saw simple multi-linear and {\it funktorial}
operations and I was hoping to later see their interesting
applications, but if he would have asked me a simple mechanical
problem like the investigation of a water bag as it was thrown
from the dorm window at the time of the questioning, then I would
have understood him, as I would not have been able to apply
Euler's equations (equivalent to the GW or their integrability
condition; I have a vague memory at this point about this trick
which can be generalized to infinite dimensional Lie groups), not
being aware of the existence of such equations. Of course the
natural question arises as to why one would study such an earthly
mundane motion instead of abstract purified heavenly mathematics
and the answer lies in the merry-go-around experiment of
Einstein's: the change of frames must be valid also in the theory
of relativity because nature does not favor any frames above
others and therein lies the trick.

Not to mention again Archimedes' method of integration, for whose
application one must completely renounce the current rigorous
definition of isometry, a fact even the classical geometers of the
late XIX$^{\mathrm{th}}$-early XX$^{\mathrm{th}}$ century were
already uneasy to do and as a consequence they renounced it with
too much restraint (they have the excuse of not having access to
Archimedes' method being taught in high-school, as Archimedes is
one of the handful whose words have such power so at to overturn
deep rooted trends (even the mighty Earth itself with the proper
choice of balance) and up to date Archimedes is the only one whose
words rewrote and may still rewrite books of history of
mathematics if the buzz around the palimpsest turns out to be
true). An interpretation is not rigorous but partial if it removes
important applications; it is not a simplification if it removes
the very foundation of the original elementary intuitive concept
(most of the times and with use of proper toys and experiments
available at the level of grade and high-school) and replaces it
with a concept only some college graduates would understand. Just
as Archimedes dealt with a raw unsound principle at his time, one
must probably carefully tread through all the current sound
foundations of geometry and allow larger non-rigorous archaic
unsound interpretations, at least initially at the level of
elementary intuitive a-priori geometric arguments (again just like
Archimedes) and in the direct interpretation as intended by their
true discoverers before partial current interpretations (if any)
were found, people were lost in translation and consequently the
heritage of our ancestors was lost. The best free minds of the
past had potential no different than that of the best free minds
of our times and they delivered magnificently on the problems at
hand. Chern was the best geometer of the second half of the
XX$^{\mathrm{th}}$ century (in fact it is Chern who raised
Geometry from almost complete oblivion to its current status), but
that happened because according to his own account he stood on the
shoulders of the giants among others Blaschke and Cartan and one
can find in the classical geometers proper justification for each
new abstract notion and distillation by means of juicy
applications (thus one must first at least cursory read Blaschke
and Cartan and only later Chern; it is locally more expensive but
on the long term is the optimal path; there is no way to fully
understand and appreciate the value of Chern's contributions
otherwise). Just as a human being passes thorough different stages
in one's mother womb (egg-invertebrate-fish-amphibian-mammal-full
developed human), a student in mathematics should be first and
foremost rather trained in the history of mathematics as it was
naturally developed and first and foremost taught rather each
important change and why it came after difficult {\it 'lupte
seculare'} (centuries' worth struggles) and the motivations behind
the scene: a full human being cannot be a full human being before
spending some time as a fish, an amphibian or a mammal; otherwise
it will lack and thus will have trouble seeing when it is
presented in front of one's eyes respectively a spine, the four
limbs or the warm blood; if chance has it that the student
actually recognizes the simpler stage of the missed development,
then the student has no choice but to fall back at that stage
(memories of one's childhood) and redo everything from scratch.
The Romanians have a strange sense of humor and consequently
corresponding jokes; I can only think of one at this juncture:
{\it A dentist from abroad pulled a tooth of a Romanian travelling
abroad during Ceau\c{s}escu's regime through his ... (aka not
through his mouth) because the Romanian was supposed to keep his
mouth shut while abroad}. I would say that Archimedes' balance
should be a tool to be kept handy by all dentists, as it can move
the mighty Earth itself with proper choice of fulcrum.

The fact that quadrics behave well {\it with respect to} (wrt)
quadratic laws and conserved quantities (like the energy or the
linear element for applicable surfaces) seems to be a tautology to
the na\"{i}ve, but to the not so na\"{i}ve the difference is
clear: 'quadratic' in quadrics is static, while 'quadratic' in
quadratic laws and conserved quantities is moving (differential).
But what is differential? According to Leibnitz and Newton it is
Euclid's world (infinitesimal or finite) in motion. While at the
infinitesimal level all things are simple and Euclidean arguments
suffice, in the finite world Euclidean arguments apply to much
less, as in the process of zooming out from the infinitesimal
world to the finite world simple objects accumulate to objects out
of the range of Euclidean geometry. Miraculously, for quadrics
certain configurations in the Euclidean infinitesimal world
(static) survive the trip to the finite world (moving) and the
na\"{i}ve turns out to be right, in the process reaping another
success: the very same survival of the configuration from the
static picture gives for free a quadratic law or conserved
quantity. Conversely, a quadratic law or conserved quantity
chooses the configuration which will survive the trip: since the
law or conserved quantity is of a general nature, it must be true
for any moving configuration at all instances (infinitesimals) in
the moving picture and thus its very nature lies in the static
picture. Once a configuration is chosen by a quadratic law or
conserved quantity, it is very likely that it will be chosen by
many others, since nature favors simplicity and thrift even more
than man. These are some of the reasons for which the wise
ancients were able to find so many interesting properties of
conics without use of calculus and for which these properties can
be found in Kepler's Laws and Newton's {\it 'Principia...'}. But
why 'quadrics' and what are 'certain configurations'? Consider an
object with as many symmetries as possible: the more symmetries,
the more chances are that that object will satisfy more laws, not
only in the static Euclidean picture, but also in the moving
picture. The first natural choice is a linear subspace and this
leads to differentials, as they live in tangent spaces (the only
objects common to both the infinitesimal and finite world and by
means of which one jumps from one to the other). The next choice
must be non-linear and simple: thus it is the sphere. In both
cases Euclid provided a consistent body of results in the static
picture and low dimensions. Now consider the moving picture: the
tangency configuration enters into our consideration with
differentials, as they live in tangent subspaces. Since 'tangency'
is a projective notion, one can replace the sphere with a
projective equivalent: an a-priori general quadric (it is known
since Apollonius that the tangents to conics can be constructed
with ruler and compass; Pascal showed that the compass can be
dropped and according to the historian there are no other curves
with this property). But keep in mind quadratic laws and conserved
quantities: these are metric properties, so we need
metric-projective tangency configurations of quadrics. Thus
appears the family of quadrics confocal to a given one: its
definition is mixed metric-projective and depends on a spectral
parameter $z$.

Consider a simple $1$-dimensional example: the rolling of an
ellipse on a plane curve; this will provide a simplified analogy
to the process of generalization of the B transformation from the
pseudo-sphere to quadrics. At any instant the ellipse and the
curve meet tangentially and with same differential at the tangency
point. Since any two curves can be rolled one onto the other
(there are no integrability restrictions on the rotation of the
rolling), the problem is of a much too general nature to be
useful, so we have to make an ansatz. If the ellipse is a circle,
the natural ansatz is that the center should preserve constant $0$
elevation, in which case the curve of rolling must be a horizontal
line. With this in mind we go back to the case of the ellipse,
when the center of the circle is replaced not with the center of
the ellipse, but with one of the two foci, as they are the heirs
of the metric-projective properties of the center of the circle.

Consider $e_1:=[1\ 0]^T,\ e_2:=[0\ 1]^T$, the focus
$F:=\sqrt{b-a}e_2$ of the ellipse $E:\
\frac{(e_1^TE)^2}{a-z}+\frac{(e_2^TE)^2}{b-z}=1$, the point
$E(s):=\sqrt{a-z}\cos\frac{s}{\sqrt{a-z}}e_1+\sqrt{b-z}\sin\frac{s}{\sqrt{a-z}}e_2$
and the rolling $(R_u,t)=(R_u(s),t(s))$ ($R_u$ is the rotation of
decreasing angle $u=u(s)$ and $t$ is a translation) of $E$ on the
curve $c(s):=c^1(s)e_1+ c^2(s)e_2:\ \ (R_u,t)E:=R_uE+t=c,\
R_udE=dc$. We need the sum of the elevation of the a-priori
inclined segment $R_uG,\ G:=F-E$ and $c^2$ to be $0$; since
$|R_u(s)G(s)|=|G(s)|=\sqrt{b-z}-\sqrt{b-a}\sin\frac{s}{\sqrt{a-z}}$,
the obvious choice of $c^2(s)$ (for which there is no inclination,
just as for the circle case) is thus revealed in the static
picture: $\ c^2(s):=-|G(s)|$. From $|dc|^2=|dE|^2$ we get
$c^1(s)=s$; now the rotation $R_u$ is uniquely determined by
$R_udE=dc$ and the translation is $t:=c-R_uE$. This is the only
choice, as the imposed $0$ elevation ansatz boils down to a first
order differential equation and thus a second preserved quantity
is revealed. Conversely, if we impose the non-inclination ansatz
$G=|G|R_{-u}e_2$, then we get the same solution. The rolled focus
$(R_u,t)F$ describes the horizontal line $R_uG+c=se_1$ uniformly
in $s$; again this provides another conserved quantity which
generalizes the circle situation. Differentiating the
non-inclination ansatz and using $R_udR_{-u}=-R_{\frac{\pi}{2}}du$
we get $|G|due_1=-R_udE-d|G|e_2 =-dc+dc^2e_2=-dse_1$, so
$ds=-|G|du$. Note also that the non-inclination ansatz is similar
to the polar coordinates at $F$: if the point $E$ is seen under
the increasing angle $\theta$ from the focus $F$ such that
$\theta=0$ for $s=-\frac{\pi\sqrt{a-z}}{2}$, then
$G=|G|R_{\theta}e_2$, so $u=-\theta$. As $E$ rolls on $c$ it is
not at equilibrium wrt the gravitational influence (constant and
pointing in the $-e_2$ direction), but this can be corrected by
considering a second counter-balancing rolling ellipse attached to
$E$ by a rigid bar connecting their foci (thus of length
$(2k+1)\sqrt{a-z}\pi,\ k\in\mathbb{N}$). Because of the
non-inclination ansatz, the length of this rigid bar remains
constant during rolling. There are further preserved quantities of
the system of two ellipses which generalize the circle situation.
The system (possibly with a weight attached at the middle point of
the rigid bar) is in equilibrium wrt the gravitational influence:
each ellipse pulls the rigid bar in opposite directions and with
force a constant multiple of
$\frac{\frac{d}{ds}|G(s)|}{1+(\frac{d}{ds}|G(s)|)^2}$. The
gravitational energy of the system is preserved, as
$G(s+(2k+1)\sqrt{a-z}\pi)=G(-s),\
\frac{F^TG(-s)}{|G(-s)|}=-\frac{F^TG(s)}{|G(s)|}$, so the heights
of the centers of the two ellipses are equal in absolute value and
of opposite signs.

From $|G(s)|+|G(-s)|=2\sqrt{b-z}$ we get
$\frac{|G(s)||G(-s)|}{a-z}=\frac{4(b-z)-|G(s)|^2-|G(-s)|^2}{2(a-z)}
=1+(\frac{d}{ds}|G|)^2$, or $\frac{1}{2|G|^2}|\frac{d}{ds}E|^2-
\frac{\sqrt{b-z}}{a-z}\frac{1}{|G|}=-\frac{1}{2(a-z)}$;
differentiating this we get
$\frac{d^2}{ds^2}|G|=\frac{1}{|G|}(1+(\frac{d}{ds}|G|)^2)-\frac{\sqrt{b-z}}{a-z}$.
The parameter $s$ behaves well for the following reason: with a
new parameter $\hat s:=\frac{s}{\sqrt{a-z}}$ independent of $z$ an
affine correspondence (Ivory affinity) with good metric properties
is established between confocal ellipses:
$E=\begin{bmatrix}\sqrt{1-za^{-1}}&0\\0&\sqrt{1-zb^{-1}}\end{bmatrix}
(\sqrt{a}\cos\hat se_1 +\sqrt{b}\sin\hat se_2)$. When $z$ varies
we get a hyperbola, as an orthogonal trajectory of the family of
confocal ellipses.

Consider the variables $a,b,z,s,\theta$ linked by the algebraic
relation $G=|G|R_{\theta}e_2$, but otherwise independent; we have
$R_{-\theta}dG=-|G|d\theta e_1+d|G|e_2,\ |G|d\theta=
\frac{-1}{|G|}dG^TR_{\frac{\pi}{2}}G=ds-sd\log\sqrt{a-z}
-\frac{\sqrt{a-z}\sqrt{b-z}}{\sqrt{b-a}}
\cos\frac{s}{\sqrt{a-z}}d\log\sqrt{\frac{a-z}{b-z}}$. Thus the
differential system $G=|G|R_{\theta}e_2,\ ds=|G|d\theta$ has as
consequence the fact that $d(a-z)=0$ {\it if and only if} (iff)
$d(b-z)=0$ iff $s$ is independent of $a-z,b-z$ (in which case $s$
can be taken as a universal variable; the system admits only the
free parameter $z$, so $d(b-a)=0$ and the rolling problem holds);
if $d(a-z)\neq 0$, then one can functionally solve for
$s=s(a-z,b-z,\frac{d(a-z)}{d(b-z)})$.

Note that $0,z<a<b,\ s\in\mathbb{R}$ for the real case, but the
computations are valid over complex numbers, so one can play with
the order of the real numbers $0,a,b,z$ and also choose $s\in
i\mathbb{R};\
(\sqrt{\ep_1}\mathbb{R})\times(\sqrt{\ep_2}\mathbb{R}),\
\ep_1,\ep_2:=\pm 1$ instead of $\mathbb{R}^2$;
$(\cosh(s),\sinh(s))=(\cos(is),-i\sin(is))$ instead of
$(\cos(s),\sin(s))$ and changing the Euclidean scalar product via
the linear transformation $[x\ y]^T\mapsto[x\ iy]^T$ or by
multiplication by $i$.

Kepler's second law can be stated:
$2=|G|^2\frac{d\theta}{dt}(=|G|\frac{ds}{dt})$; thus $t$ is the
new universal variable (independent of $a-z,b-z$) and
$s=s(t,a-z,b-z),\ \theta=\theta(t,a-z,b-z)$; now
$\frac{d^2}{dt^2}E=\frac{ds}{dt}\frac{d}{ds}
(\frac{ds}{dt}\frac{d}{ds}E)=\frac{4}{|G|}\frac{d}{ds}
(\frac{-1}{|G|}R_{\theta}(e_1+\frac{d}{ds}|G|e_2))=4(1+(\frac{d}{ds}|G|)^2
-|G|\frac{d^2}{ds^2}|G|)\frac{R_{\theta}e_2}{|G|^3}
=\frac{4\sqrt{b-z}}{a-z}\frac{1}{|G|^2}\frac{G}{|G|}$ and thus
agrees with the inverse square distance law of universal
gravitational attraction, on the condition that the mass of the
sun $\frac{\sqrt{b-z}}{a-z}$ is independent of $a-z,b-z$. But this
is where Kepler's third law states what is needed: the mean radius
of the orbit is $\sqrt{b-z}$; due to the metric properties of the
foci this is true for any probability measure on the ellipse $E$
having the symmetry which exchanges the foci; in particular the
period of the orbit is
$\int_0^Tdt=\frac{1}{2}\int_0^{2\pi\sqrt{a-z}}|G|ds=\pi\sqrt{a-z}\sqrt{b-z}$.
Note the conservation of energy
$\frac{1}{2}|\frac{d}{dt}E|^2-\frac{4\sqrt{b-z}}{a-z}\frac{1}{|G|}=-\frac{2}{a-z}$
and the fact that the mass of the sun $\frac{\sqrt{b-z}}{a-z}$ is
the maximum of the curvature of the ellipse $E$, obtained for
$s\in\pi(\frac{1}{2}+\mathbb{Z})$*. \footnote{* I was reminded by
one of my Calculus III students that the area element for a graph
in polar coordinates is $\frac{1}{2}r^2d\theta$ instead of
$r^2d\theta$; thus the area element is
$d\wedge(\frac{1}{2}r^2d\theta)=rdr\wedge d\theta$.}

Thus Newton's differential system can be interpreted as the second
in a hierarchy of differential systems whose first one is the
rolling problem and the properties of all such differential
systems emerge from the static picture.

Most known integrable systems come in hierarchies.

Although any two plane curves can be rolled one onto the other,
there exists a notion of {\it osculating rolling} of two space
curves of same curvature (and thus the rolling is restricted)
which naturally appears as a degeneration of rolling surfaces. The
tangent bundle of a plane curve is identified with a mostly
multiple cover of a certain region of the plane of the curve; if
one deforms this mostly multiple covered flat region, then we get
the tangent surface (which is flat) of a space curve having the
same curvature as the initial plane curve. One can roll this flat
surface on the plane such that the two curves osculate and roll
one onto the other. If the space curve is plane, then it must be
congruent to the original one. In our case if one rolls the
ellipse on the other side of the same (fixed) ellipse (at any
instant the two ellipses reflect in the common tangent line), then
its foci describe circles centered at the other foci, just as in
the circle case. The same statement remains true for hyperbolas
and parabolas (parabolas have a finite focus and a focus at
$\infty$; thus a circle has its centers on the line at $\infty$
and becomes a line; the other rolling focus on the line at
$\infty$ describes the line at $\infty$). For the ellipse we have
the normal direction $\hat
N:=\frac{\sqrt{b-z}}{\sqrt{a-z}}\cos\frac{s}{\sqrt{a-z}}e_1+
\sin\frac{s}{\sqrt{a-z}}e_2,\ \hat N^TdE=0,\ \hat N^TG=-|G|$ and
the rolling focus $F':=F-\frac{2\hat N^TG}{|\hat N|^2}\hat N$
satisfies $|F'+F|^2=4|F+\frac{|G|}{|\hat N|^2}\hat N|^2
=4(|F|^2+\frac{|G(s)||G(-s)|}{|\hat N|^2})=4\frac{b-z}{b-a}$.

The case when conics osculate and roll on space curves having the
same curvature is discussed in \S\ 6.5 (B transformations of
singular quadrics).

While intuitively the metric-projective character of the confocal
family of the ellipse is clear, it lacks the advantage of
generality, as it must be separately stated for the parabola (the
confocal family of the ellipse also includes hyperbolas for
$a<z<b$) and for imaginary conics. The metric-projective
definition thus must be independent of the chosen conic and it was
found by Pl\"{u}cker. Earlier Poncelet introduced the imaginary
points $[1,i,0],\ [1,-i,0]$ at $\infty$ to account for the fact
that two circles must intersect at four points, just as any other
two conics in general position (Apollonius and B\'{e}zout); now
Cayley's absolute $C(\infty):=\{[1,i,0],\ [1,-i,0]\}$ encodes the
metric structure of $(\mathbb{C}^2,<,>)$, as rigid motions and
homotheties are the only {\it homographies} (projective
transformations) which preserve $C(\infty)$. Note that $C(\infty)$
is the quadric $(x^1)^2+(x^2)^2=0$ in the line at $\infty:\
[x^1,x^2,0],\ (x^1,x^2)\neq(0,0)$ of $\mathbb{CP}^2$. The ellipse
is given by the equation $\frac{(x^1)^2}{a}+\frac{(x^2)^2}{b}=1,\
a,b>0$, the hyperbola by the same equation with $a>0>b$ and the
parabola by the equation with $\frac{(x^1)^2}{a}=2x^2,\ a>0$; all
conics are given by the equation $x^TAx=0,\ x=[x^1\ x^2\
1]^T\simeq[x^1,x^2,1]$ with certain matrices
$A=A^T\in\mathbf{GL}_3(\mathbb{C})$. One can define the {\it
adjugate} of the conic $x:\ \ x^TAx=0$ to be the conic $x^*:\ \
{x^*}^TA^{-1}x^*=0$; if $y=Hx$ for
$H\in\mathbf{PGL_3(\mathbb{C})}$, then $y^*=(H^T)^{-1}x^*$ and the
'adjugate' denomination becomes clear. $C(\infty)$ is by
definition a self-adjugate quadric in the line at $\infty$ (this
is equivalent to bringing an arbitrary non-degenerate scalar
product $<x,y>:=x^TQy,\ Q=Q^T\in\mathbf{GL}_2(\mathbb{C})$ to the
canonical form $<x,y>:=x^Ty$).

The family $\{x_z\}_{z\in\mathbb{C}}$ of conics confocal to the
given conic $x$ is {\it the adjugate of the pencil of conics
generated by the adjugate of the initial conic and $C(\infty)$}:
$x_z^T(A^{-1}-zI_{1,2})^{-1}x_z=0,\ I_{1,2}:=\mathrm{diag}[1\ 1\
0]$; thus $x_0=x$ and the mixed metric-projective character
becomes clear.

For $z:=\infty$ we obtain the line at $\infty$ as a singular conic
$x_{\infty}$ of the confocal family; it contains as a singular set
the quadric $C(\infty)$. Similarly for $z$ inverses of non-zero
eigenvalues of the quadratic part of the conic ($z:=a$ or $b$ for
ellipses or hyperbolas and $z:=a$ for parabola) we obtain singular
conics $x_z$ of the confocal family (lines); again these contain
as singular sets other quadrics: the two foci (the parabola has a
finite focus and a focus on the line at $\infty$). Thus the family
of confocal conics can be interpreted as the adjugate of the
pencil of conics generated by the adjugate of the quadric
consisting of the two foci and the quadric $C(\infty)$.

Pl\"{u}cker's characterization of confocal quadrics is: for
ellipses and hyperbolas each conic of the family is tangent to the
four lines joining the two foci and $C(\infty)$ (the two other
points of intersection of these four lines are imaginary foci);
for parabolas each conic of the family is tangent to the three
lines joining the foci and $C(\infty)$. This characterization
follows from simple projective arguments.

An easy consequence of this metric-projective definition is that
if $H$ is a homography such that the line containing
$H^{-1}(C(\infty))$ intersects the line containing $C(\infty)$ at
a point not in $C(\infty)$, then $H$ takes the family of confocal
quadrics with foci among $H^{-1}(C(\infty))$ to the family of
confocal quadrics with foci among $H(C(\infty))$, while preserving
all metric-projective properties and thus all integrable systems
whose integrability boils down to these metric-projective
properties (see also Darboux (\cite{D1},\S\ 603)).

Stepping on the path opened by Euler half a century earlier
(namely differential geometry of surfaces, dynamics of a rigid
body and fluids), Dupin, Ivory and Lam\'{e} met in the beginning
of the XIX$^{\mathrm{th}}$ century at confocal quadrics (some of
the properties of confocal quadrics were also known to Monge, who
integrated the equation of lines of curvature on the ellipsoid, in
the process introducing the notion of {\it umbilic}; it is less
known since Darboux that there are 12 umbilics on the general
ellipsoid, among which four are the usual real ones; if one
excludes the importance of the remaining imaginary $8$, then one
must also exclude the D transformation of special isothermic
surfaces and thus Bianchi's theory of deformations of quadrics).
They were not primarily interested in the geometrical properties
of confocal quadrics, but in natural phenomena amenable to a
mathematical formulation at the time: (optical properties of)
light, gravity and heat.

Dupin studied light and its reflection in surfaces; geometrically
this amounts to normal congruences. These congruences remain
normal after reflection and refraction (Malus) in surfaces; a
surface normal to the normal congruence is cut by the {\it flats}
(developables) of the congruence along lines of curvature.
Developables are important for the following reason: the normal
congruence is tangent to its two focal surfaces (caustics) and
thus it induces direction fields on them; the integral curves of
these direction fields (geodesics) give the developables (each
developable is generated from the lines of the congruence by
keeping constant a parameter of an integral curve and varying the
parameter of the other or equivalently they are the positions of
the line in a plane as that plane rolls along the line on the
focal surface).

Ivory studied the gravitational attraction of ellipsoids; in
particular he transformed the problem of the gravitational
attraction a homogeneous ellipsoid exercises upon an exterior
point to the problem of the gravitational attraction a confocal
homogeneous ellipsoid exercises upon a corresponding interior
point; his work on the ellipsoidal equilibrium configuration of
self-gravitating fluids was an extension of that of Laplace
(Laplace kept Ivory in high esteem for his critical commentary of
{\it M\`{e}canique c\'{e}leste}) and influenced the results of
Jacobi and Liouville on the geodesic flow.

Lam\'{e} was interested in heat and confocal ellipsoids provide an
example of level sets of the temperature for the stable heat
equation with a singular ellipsoid of the family (the convex hull
of the singular set which is ellipse) being the source of constant
temperature; since the gradient of the temperature is normal on
its level sets, this provides an easy physical justification of
Lam\'{e}'s result. The actual problem Lam\'{e} was interested is
more general: {\it To find all triply orthogonal families formed
by isothermic surfaces}. Bonnet completely solved it and showed
that, excluding developables, only quadrics can appear as surfaces
in such families (see Darboux (\cite{D1},\S\ 513),(\cite{D11},Ch
III)).

The Ivory affinity between two confocal quadrics can be easily
seen now: the temperature flow from a point on the first quadric
cuts the second quadric at the point which gives the Ivory
affinity.

Consider the Euclidean scalar product on $\mathbb{C}^n:\ \ \
<x,y>:=x^Ty,\ |x|^2:=x^Tx$. A {\it quadric} $x\subset\mathbb{C}^n$
is given by the equation $Q(x)=0$, where
$Q(x):=\begin{bmatrix}x\\1\end{bmatrix}^T
\begin{bmatrix}A&B\\B^T&C\end{bmatrix}\begin{bmatrix}x\\1\end{bmatrix}=x^T(Ax+2B)+C,\
A=A^T\in\mathbf{M}_n(\mathbb{C}),\ B\in\mathbb{C}^n,\
C\in\mathbb{C},\ \begin{vmatrix}A&B\\B^T&C\end{vmatrix}\neq 0$;
$x$ is called {\it degenerate quadric} if
$\begin{vmatrix}A&B\\B^T&C\end{vmatrix}=0$. A quadric $x$ belongs
to either of the types: {\it quadric with center} (QC) for
$\ker(A)=0,\ Q(-A^{-1}B)\neq 0$, the center being $-A^{-1}B$; {\it
quadric without center} (QWC) for $\ker(A)=\mathbb{C}v,\ v^TB\neq
0,\ |v|^2\neq 0$; {\it isotropic quadric without center} (IQWC)
for $\ker(A)=\mathbb{C}v,\ v^TB\neq 0,\ |v|^2=0$. After applying
suitable rigid motions, most quadrics are {\it diagonal} QC ($A$
is diagonal, $B=0, C=-1$).

One can complete a quadric $x$ to a family of confocal quadrics
$\{x_z\}_{z\in\mathbb{C}}$: the adjugate of the pencil of quadrics
generated by the adjugate of $x$ and the adjugate of Cayley's
absolute (taken to be self-adjugate). Thus
$x_z\subset\mathbb{C}^n$ is given by $Q_z(x_z)=0$, where
$Q_z(x_z):=\begin{bmatrix}x_z\\1\end{bmatrix}^T(\begin{bmatrix}A&B\\B^T&C\end{bmatrix}^{-1}-z
\begin{bmatrix}I_n&0\\0^T&0\end{bmatrix})^{-1}\begin{bmatrix}x_z\\1\end{bmatrix}$, so $x=x_0$.

The condition that a {\it nonsingular quadric} ($z$ is not inverse
of non-zero eigenvalue of $A$ and $z\neq\infty$) of a {\it general
confocal family} ($x$ does not have continuous rotational
symmetries; the obvious choice for the origin is clear for QC; a
similar choice will appear later for (I)QWC) contains a point
$x\in\mathbb{C}^n$ is a polynomial of degree $n$ in $z$ which has
in general distinct roots. Therefore through any point
$x\in\mathbb{C}^n$ pass in general $n$ quadrics of the confocal
family: $x\in x_{z_j},\ j=1,...,n$; along the higher co-dimension
set of $\mathbb{C}^n$ where normals are isotropic pass less than
$n$ confocal quadrics. With certain modifications the statement
remains true for singular quadrics of the family or non-general
quadrics. Replacing the cartesian coordinates of $x$ with
$(z_j)_{j=1,...,n}$ gives elliptic coordinates on most of
$\mathbb{C}^n$ (Lam\'{e}). For non-general quadrics one can
complete the set of such $(z_j)_{j=1,...,n-k}$ with $k$
coordinates of the space obtained by factoring the continuous
group of rotational symmetries with the subgroup fixing a general
point of the quadric, since all confocal quadrics will admit the
same continuous group of rotational symmetries. Thus all
non-general quadrics are general in a non-isotropic subspace of
$\mathbb{C}^n$ which intersects orthogonally all trajectories of
points of $\mathbb{C}^n$ under the action of a continuous group of
symmetries (lines passing through centers of spheres for spheres
and planes passing through the axis of revolution for quadrics of
revolution for $n=3$) and conversely from general quadrics one can
obtain non-general ones in higher dimensions by acting on the
surrounding space of the initial general quadric with a continuous
group of rotational symmetries of the bigger space so as to fill
it. The sense in which this statement can be interpreted will be
explained later; intuitively a quadric $x$ is general iff all
eigenvalues of $A$ have geometric multiplicity $1$ (including the
metrically most degenerate case $\mathrm{spec}(A)=\{0\}$ which
thus must be an IQWC). Thus the use of elliptic coordinates
(completed with other coordinates for non-general quadrics) is
essential in the {\it metric classification} of all quadrics and
in studying integrable systems related to them (as shown since
Jacobi and Liouville, most known integrable systems boil down one
way or another to the use of elliptic coordinates).

Up to homographies all possible degenerate quadrics appear as
singular quadrics of a certain confocal family (one can find a
non-degenerate quadric in a hyperplane such that the pencil
generated by the adjugates (in a sense that will become clear
later) is mostly formed by non-degenerate quadrics and consider
the later quadric as Cayley's absolute).

The classical theorems about confocal quadrics are as follows:

{\it (Lam\'{e}) If $x$ is a general quadric, then the normal to
$x_z$ is proportional to $\pa_zx_z$ and the vectors
$(\pa_z|_{z=z_j}x_{z_j})_{j=1,...,n}$ are orthogonal; thus the
family of quadrics confocal to a given general quadric forms an
orthogonal system.}

{\it (Dupin) General confocal quadrics cut each other along
principal sub-manifolds; thus they admit lines of curvature
parametrization (the elliptic  coordinates).}

Thus Dupin's indicatrix (essentially the shape) is encoded by the
confocal family.

If the quadric is not general, then one can still complete the
confocal family to an orthogonal family; for example confocal
quadrics in $\mathbb{C}^3$ and of revolution around an axis can be
completed with the family of planes passing through the axis of
revolution. For confocal (pseudo-)spheres the indeterminacy is
even higher, because the principal lines are not defined; however
the level sets of the standard spherical coordinates (planes
through an axis and circular cones with the same axis) appear as a
natural choice.

This result subsists for surfaces applicable to quadrics only at
the level of the second fundamental form missing mixt terms ({\it
conjugate} system of coordinates; in this matter Dupin preserved
the denomination chosen by Apollonius since at the level of Dupin
indicatrix the two coincide). The fact that Apollonius was onto
something is confirmed by the simple projective explanation of
conjugate systems which preserves the property of normal conjugate
diameters (corresponding to principal directions on a surface):
the ruled surface generated by conjugate directions along a curve
on a surface is developable (a projective notion, since asymptotic
directions coincide on such surfaces and osculating planes of
asymtotes are tangent to the surface) and any two pairs of
conjugate diameters of a conic are in harmonic ratio; moreover
Apollonius introduced the tangent to a conic with center among
others as the parallel to a conjugate diameter through the given
point (which thus becomes the endpoint of a diameter). Conversely
if we bend a plane on a surface along a curve, then the rulings on
the resulting developable are the conjugate directions of the
curve. Thus the statement 'triply orthogonal system of quadrics'
becomes 'triply conjugate system of surfaces containing a family
of surfaces applicable to quadrics'; triply conjugate system of
surfaces means surfaces cutting each other along conjugate lines
(see Bianchi (\cite{B2},Vol {\bf 4},(143),(146))). Because the
principal lines on (pseudo-)spheres are indeterminate, the triply
conjugate system of surfaces containing a family of CGC surfaces
are actually triply orthogonal (the common conjugate system on a
surface and on the unit sphere is realized by the lines of
curvatures on the given surface and the correspondence is given by
the Gau\ss\ map). Thus Calapso studied deformations of quadrics by
imposing the natural condition that a conjugate system on the
initial quadric is also a conjugate system on one of its
deformations; by this {\it infinitesimal analytic} method
(although more complicated than Bianchi's simple geometric
picture) he managed to take Bianchi's results up a notch in 1912.

{\it (Ivory) There exists an affine transformation of
$\mathbb{C}^n$ which takes $x_0$ to $x_z$; it preserves the
lengths of segments between confocal quadrics: with
$V_0^1:=x_z^1-x_0^0,\ V_1^0:=x_z^0-x_0^1$ we have
$|V_0^1|^2=|V_1^0|^2$ for pairs of points $(x_0^0,x_z^0),\
(x_0^1,x_z^1)$ corresponding on $(x_0,x_z)$ under this affine
transformation}.

Such an affine transformation is called {\it the Ivory affinity}
because it is friendlier than an usual affine transformation. In
the same vein we have the result mostly due to Bianchi (see
(122,\S\ 10)):

{\it (Bianchi I) The Ivory affinity preserves lengths of segments
on rulings (segments on tangent cones) of confocal quadrics, the
tangency configuration (TC), angles between segments (between
confocal quadrics) and rulings, between rulings and between polar
rulings.

These can be written as: $|w_z|^2=|w_0|^2$ if $x_0^0+tw_0\in x_0,\
x_z^0+tw_z\in x_z,\ \forall t\in \mathbb{C}$ and $w_0,\ w_z$
correspond under the Ivory affinity;
$(V_0^1)^T\pa_z|_{z=0}x_z^0=(V_1^0)^T\pa_z|_{z=0}x_z^1$, thus
$x_z^1\in T_{x_0^0}x_0$ iff $x_z^0\in T_{x_0^1}x_0;\
(V_0^1)^Tw_0^0= -(V_1^0)^Tw_z^0,\ (w_0^0)^Tw_z^1=(w_z^0)^Tw_0^1,\
(w_z^0)^T\hat w_z^0=(w_0^0)^T\hat w_0^0$ if $x_0^0+sw_0^0+t\hat
w_0^0\in x_0,\ \forall s,t\in \mathbb{C}$.

Thus for $x_0^0,\ x_0^1$ points on a quadric
$x_0\subset\mathbb{C}^{2n-1}$ and $W_0^0,\ W_0^1$ maximal sets of
mutually polar rulings at  $x_0^0,\ x_0^1$ (each containing $n-1$
rulings), with $F:=[V_0^1\ W_0^0\ W_z^1]$ we have
$F^TF=((0\leftrightarrow z)\circ F)^T((0\leftrightarrow z)\circ
F)$ and there exists a RMPIA $(R_0^1,t_0^1)$ such that
$R_0^1F=(0\leftrightarrow z)\circ F,\
t_0^1=x_z^0-R_0^1x_0^0(=x_0^1-R_0^1x_z^1)$; moreover
$(R_0^1,t_0^1)$ is unique precisely when $\det F\neq 0$.}

The preservation of lengths of segments on rulings under the Ivory
affinity between confocal quadrics ($n=3$) give 'almost'
Tchebyshev coordinates which are further asymptotic; they are
'almost' because in general we cannot change these coordinates to
make the rulings of unit length. Thus we get (up to rigid motions)
the family of confocal quadrics if we deform the doubly ruled
fabric of an initial quadric, changing the angles between the
ruling families while preserving their lengths. This construction
of the articulated hyperbolic paraboloid is due to Henrici so some
of these results may have already been folklore at the time (see
Bianchi \cite{B1}).

In 1838 Jacobi integrated the equation of the geodesics on an
ellipsoid; geometrically his result can be stated:

{\it (Jacobi) The tangent lines to a geodesic on $x_0$ remain
tangent to $n-2$ other confocal quadrics.}

This is the first instance when the simple metric-projective
definition of confocal quadrics has as consequence a deep
metric-projective property, which hints at the fact that the
confocal family of a quadric influences surfaces applicable to
that quadric. Chasles proved its converse:

{\it (Chasles) The common tangents to $n-1$ confocal quadrics form
a normal congruence and envelope geodesics on the $n-1$ confocal
quadrics.}

Thus confocal quadrics not only {\it infinitesimally} cut each
other orthogonally (Dupin and Lam\'{e}) but also {\it finitely},
that is as seen by the viewer from any point of the space
(Darboux). This is another principle: {\it infinitesimal} laws for
confocal quadrics admit corresponding {\it finite} laws (and
according to Archimedes the converse is true); the theory of
deformations of quadrics is one of its facets and its
discretization realizes the complete unification of the {\it
infinitesimal} and {\it finite}.

Consider the lines $x_{z_1}+sv,\ s\in \mathbb{C}$ tangent to the
$n-1$ quadrics $x_{z_j},\ j=1,...,n-1$ at points
$x_{z_j}:=x_{z_1}+s_{z_j}v$. As we shall see later, this implies
that the normals along the common tangents are orthonormal:
$N_{z_j}^TN_{z_k}=\del_{jk},\ k=1,...,n-1$; thus $0,\ \{z_j\}_j$
can be considered eigenvalues with corresponding eigenvectors $v,\
\{N_{z_j}\}_j$ for a certain symmetric operator which undergoes an
isospectral deformation. The approach to integrable systems by
means of isospectral deformation of a linear operator was first
developed by Lax for KdV and it has since emerged into a powerful
technique used to study most known integrable systems. Moser
established in \cite{M1} and via the above considered symmetric
operator links between confocal quadrics and other integrable
systems, thus returning to the basics: confocal quadrics capture
(encode) and geometrize spectral properties of symmetric
operators.

While there are many transformations of $\mathbb{C}^n$ which take
a family of confocal quadrics to another one, a transformation
which preserves the previous properties of confocal quadrics must
take lines to lines, so it must be a homography. The converse is
mostly true and due to Bianchi (see (\cite{B2},Vol {\bf 5},(90))):

{\it (Bianchi II) For any general homography $H$ of
$\mathbb{C}^n$, there are two unique families of confocal quadrics
in $\mathbb{C}^n$ such that $H$ takes one family to the other.
Moreover, it takes principal directions ($\pa_z|_{z=z_j}x,\
j=1,...,n$, where $\mathbb{C}^n\ni x\in x_{z_j}$) to principal
directions (Lam\'{e}) and geodesics to geodesics (Chasles). Thus
$H$ induces a linear map on each tangent space of $\mathbb{C}^n:\
\ dH(x):T_x\mathbb{C}^n\mapsto T_{H(x)}\mathbb{C}^n$ and the
singular value decomposition (SVD) of $dH(x)$ gives the principal
directions in both the domain and the range of $dH(x)$. Moreover
the Ivory affinities between $(x_0,x_z)$ and $(H(x_0),H(x_z))$
(and thus the RMPIA) correspond under $H$; also $H$ preserves the
isospectral deformation character of the above considered
symmetric operator.}

Note that $A=A^T\in\mathbf{M}_n(\mathbb{C})$ cannot always be
diagonalized via a conjugation with a rotation
$R\in\mathbf{O}_n(\mathbb{C})$; we can hope at best to obtain a
{\it symmetric Jordan} (SJ) form (a form as close as possible to a
diagonal form and depending on as few as possible constants). Thus
the SVD does not always reduce $A\in\mathbf{GL}_n(\mathbb{C})$ to
a diagonal matrix: $A=R_1DR_2^T,\ R_1,\
R_2\in\mathbf{O}_n(\mathbb{C})$ and $D$ is SJ ($R_1,R_2,D$ can be
mostly found from $A^TA=R_2D^2R_2^T,\ AA^T=R_1D^2R_1^T$). Note
also that the two families of confocal quadrics must have the same
continuous group of rotational symmetries.

The natural question about the existence of the SJ form appeared
since I was using identities of the form
$A\sqrt{I_n-zA}=\sqrt{I_n-zA}A,\
A=A^T\in\mathbf{M}_n(\mathbb{C}),\
z\in\mathbb{C}\setminus\mathrm{spec}(A)^{-1}$.

Note that according to Apollonius any line cuts two homothetic QC
along two equal segments thus individuated (including the singular
cases when the homothety is $0,\infty$, so the QC becomes a
quadratic cone or the hyperplane at $\infty$). Since this is not a
true metric-projective property, the natural question arises if
the result still stands by replacing the pencil of homothetic QC
with the family of confocal quadrics. The answer is no, but
because of the symmetry of the TC the next natural question arises
wether a segment between two points on a quadric and whose
subjacent line is tangent to a confocal quadric is taken by the
Ivory affinity to a ray of light trajectory (that is reflecting in
the initial quadric) and tangent to the given confocal quadric at
two points. The answer is not only is yes, but also its converse
(a ray of light tangent to a quadric and reflecting in a confocal
quadric is still tangent to the given quadric; moreover the Ivory
affinity takes the broken segment into a single one) is valid;
applying Ivory's Theorem and the symmetry of the TC we get an
intuitive a-priori idea about rays of light tangent to a quadric
and reflecting on a confocal one (note that any such ray comes
with a dual friend of same length obtained from the Ivory
affinity). In this sense we have a result due to Chasles-Darboux
about a ray of light whose trajectory is formed by polygonal
segments on lines circumscribed (tangent) to $n-1$ given confocal
quadrics and bouncing (reflecting) sequentially on another given
sequence of quadrics $Q_1,Q_2,...$ of the same confocal family as
the initial one: if one such trajectory is closed, then all are;
moreover all such trajectories have same length (see Darboux
(\cite{D1},\S\ 465)). Note again the principle of correspondence
between {\it infinitesimal} and {\it finite} by means of {\it
discretization}: the rolling of an $(n-1)$ dimensional space along
one of its lines on a quadric by means of which one can obtain
Jacobi's geodesics becomes at the discrete level {\it precisely}
the construction and motivation of Darboux's interpretation (due
to the good metric properties of the Ivory affinity as pointed out
by Ivory himself). One can replace the ray of light with a
billiard ball; billiards on confocal quadrics is a current subject
of research in integrable systems (due mostly to Serge
Tabachnikov). Of course even such things are preserved by the
homography of Bianchi II and one can investigate them first in the
complex setting and then restrict, as an application, to totally
real cases*. \footnote{* At this point the fact that the Ivory
affinity provides the correct explanation is pure speculation;
Darboux uses differential computations in elliptic coordinates;
should the trick with the Ivory affinity work one can use instead
algebraic computations.}

Delaunay proved in 1841 that if we roll a conic on an line and
then rotate around the line the curve described by its foci
(focus), then we get all CMC (minimal) surfaces of revolution.

In 1859 the French Academy posed the problem:
\newline
\begin{center}
{\it To find all surfaces applicable to a given one}.
\newline
\end{center}
It became the driving force which led to the flourishing of the
classical differential geometry in the second half of the
XIX$^{\mathrm{th}}$ century and its profound study by illustrious
geometers led to interesting results (see Bianchi
\cite{B1},\cite{B2}, Darboux \cite{D1}, Eisenhart
\cite{E1},\cite{E2}, Sabitov \cite{S1} and its references for
results up to 1990's or Spivak (\cite{S3}, Vol {\bf 5})). Today it
is still an open problem in its full generality, but basic
familiar results like the Gau\ss-Bonnet theorem and the
Codazzi-Mainardi equations (independently discovered also by
Peterson) were first communicated to the French Academy. A list
(most likely incomplete) of the winners of the prize includes
Bianchi, Bonnet, Guichard, Weingarten.

Weingarten proved in 1863 that the focal surfaces of the
congruence of normals of a surface having a functional
relationship between its curvatures are applicable to surfaces of
revolution (which depend only on the functional relationship) and
conversely, the tangents to geodesics corresponding to meridians
on a surface applicable to a surface of revolution (meridians on
the given surface) form a normal congruence, whose normal surfaces
have a functional relationship between their curvatures. Such
normal surfaces are called {\it Weingarten} (W). Note that in the
correspondence established between the two focal surfaces by the
normal W congruence to the meridians on a focal surface correspond
parallels on the other, since on any normal congruence the
developables of different families cut each other perpendicularly
(Dupin and Malus). Since the two focal surfaces as above are
applicable to surfaces of revolution, Darboux's statement about
the deformation of a surface of revolution to a line has a simple
geometric explanation: as a focal surface rolls on its applicable
surface of revolution, the other focal surface (called {\it
complementary} surface) rolls on the axis of revolution, with the
arc-length of the meridians corresponding to the arc-length of the
axis of revolution (each meridian on the second focal surface
rolls on the axis of revolution precisely when the first focal
surface rolls on the applicable surface of revolution along a
parallel). Or: as the first focal surface is deformed to its
corresponding surface of revolution, the complementary surface is
deformed to the axis of revolution, each meridian becoming the
axis of revolution.

Ribaucour proved that the second fundamental forms of focal
surfaces as above are proportional; therefore Bianchi calls a
congruence on whose focal surfaces the second fundamental forms
are proportional a {\it Weingarten} (W) congruence. This notion is
invariant under homographies, since involves only tangency: the
asymptotic directions (a projective notion) correspond.

W congruences play an important r\^{o}le in the infinitesimal
deformation problem (Darboux proved that the later generates W
congruences and Guichard proved the converse); thus the
infinitesimal deformation problem is essentially projective.

As an example Weingarten classified the W surfaces with
$2(r_2-r_1)=\sin(2(r_2+r_1))$, in which case the lines of
curvature correspond via the Gau\ss\ map to geodesic confocal
ellipses and hyperbolas on the unit sphere. By 1894 Darboux proved
that the {\it evolutes} (focal surfaces of the normal congruence)
of such W surfaces are applicable to imaginary paraboloids of
revolution, completed the study of such W surfaces by replacing
$\sin$ with $\sinh$ (in which case the paraboloids become real)
and provided a simple geometrical construction of such W surfaces:
they are surfaces of translation with the two generating curves
having constant complex conjugate torsions (or opposite real);
their evolutes give all surfaces applicable to (imaginary)
paraboloids of revolution (see (\cite{D1},\S\ 745-\S\ 751)).
Therefore Darboux posed the natural question (still open in its
full generality): {\it 'To find all algebraic curves of constant
torsion'} (see (\cite{D1},\S\ 39 and Note 4 of Vol {\bf 4})).

B\"{a}cklund constructed in 1883 a transformation for CGC $-1$
surfaces (this point of view is due, according to Bianchi
(\cite{B1},\S\ 374), to Lie): the tangent spaces to the unit
pseudo-sphere $x_0$ cut a confocal pseudo-sphere $x_z$ along
circles, thus highlighting a circle in each tangent space of
$x_0$. Each point of the circle, the segment joining it with the
origin of the tangent space and one of the (imaginary) rulings on
$x_z$ passing through that point determine a facet. We have thus
highlighted a $3$-dimensional integrable distribution of facets:
its leaves are the ruling families on $x_z$. If we roll the
distribution while rolling $x_0$ on an non-rigidly applicable
surface $x$ (called {\it seed}), it turns out that the
integrability condition of the rolled distribution is always
satisfied (we have complete integrability), so the integrability
of the rolled distribution does not depend on the shape of the
seed. The rolled distribution is obtained as follows: each facet
of the original distribution corresponds to a point on $x_0$; we
act on that facet with the rigid motion of the rolling
corresponding to the highlighted point of $x_0$ in order to obtain
the corresponding facet of the rolled distribution. Note that the
facets of the $3$-dimensional distribution are differently
re-distributed into $2$-dimensional families of facets as tangent
planes to leaves when the shape of the seed changes, but
principles and properties independent of the shape of the seed
remain valid for facets even in the singular picture: this is
Archimedes' contribution to the Bianchi-Lie ansatz. The leaves of
the rolled distribution (called the B transforms of $x$, denoted
$B_z(x)$ and whose determination requires the integration of a
Ricatti equation) are applicable to the pseudo-sphere. Moreover
the seed and any leaf are the focal surfaces of a W congruence, so
the inversion of the B transformation has a simple geometric
explanation (the seed and leaf exchange places). For simplicity we
ignore in the notation $B_z(x)$ the dependence on the initial
value of the Ricatti equation.

The transformation originally constructed by B\"{a}cklund states
that if a CGC $-1$ seed $x$ in general position and an angle
$0<\theta<\frac{\pi}{2}$ are given, then the $3$-dimensional
distribution formed by facets with centers on circles of radius
$\sin\theta$ in tangent spaces of $x$, of inclination $\theta$ to
these and passing through the origin of these is integrable;
moreover the leaves are CGC $-1$ surfaces. Lie's point of view
becomes clear now: because the B transformation is of a general
nature (independent of the shape of $x$), it must exist (at least
as a limiting case) also in the case when $x$ in not in general
position, but it actually coincides with the pseudo-sphere. In
this case the $1$-dimensional family of non-degenerated leaves
(surfaces) degenerates to a $1$-dimensional family of degenerated
leaves (isotropic rulings on confocal pseudo-sphere) and thus the
true nature of the B transformation is revealed at the static
level of confocal pseudo-spheres.

The B transformation when the confocal pseudo-sphere is the
isotropic cone (sphere of zero radius or equivalently
$\theta=\frac{\pi}{2}$) was constructed even earlier (1879, upon
some results of Ribaucour from 1870) by Bianchi in his PhD thesis
and named the complementary transformation, since the seed and the
leaf are the focal surfaces of a normal W congruence. Ribaucour
did not look at his results as a method of generating CGC $-1$
surfaces from a given one (transformation) and thus Bianchi is
credited with the idea of using transformations in the study of
surfaces: one can obtain families of CGC $-1$ surfaces depending
on arbitrary many constants by iterating the complementary
transformation (each constant being introduced by the integration
of a Ricatti equation).

Lie proved the inversion of the complementary transformation; thus
its iteration is realized by quadratures since a Ricatti equation
degenerates to a linear one once a solution is known. He also
discovered the $1$-dimensional (spectral) family of CGC $-1$
surfaces, but this was only {\it infinitesimally} determined (the
fundamental forms are known), with no explicit procedure of
constructing it without the integration of Ricatti equations and
quadratures; moreover this Lie transformation cannot be iterated.
According to a formula recently discovered by Sym in \cite{S4},
one can replace the quadrature step with a derivative in the
spectral parameter; however the most difficult part (the
integration of a family of Ricatti equations) remains. According
to Bianchi the B transformation is a conjugation of a
complementary transformation with a Lie transformation.

Bianchi proved in 1892 the BPT for the B transformation; it is
with the BPT that the B transformation can be iterated using only
algebraic computations (after first integrating a  $1$-dimensional
family of Ricatti equations and taking some derivatives to find
the rolling).

Up to Weingarten's prize winning article of 1896 (\cite{W1}),
Darboux's intrinsic Monge-Amp\`{e}re equation of isometric
immersion of an abstract linear element already appeared in
conjunction with the deformation problem, but Weingarten made the
essential remark that the problem already assumes the existence of
a given surface (see Darboux (\cite{D1},vol {\bf 4}, ch XIII and
XIV) and Eisenhart (\cite{E1}, ch X)). However, few applications
of Weingarten's method have been since found; most of them give
the already found results for surfaces of revolution or of
Goursat's deformations of certain imaginary quadrics. Probably one
must modify Weingarten's approach so as to take into consideration
not only the existence of the given surface, but also that of a
deformation ({\it seed}). The only question from an audience of
graduate students after a talk about Weingarten's method in Spring
2004 (where I mentioned, just as Darboux does, an equation of
Euler's in conjunction with the propagation of the sun, clearly an
example of integrable systems) was where am I heading with this
program. Looking back I would say that that is an important
question; the answer would still elude me for more than two years.
It has been an open question ever since the classical times to
find all surfaces for which an interesting theory of deformations
can be build, just as for quadrics (see \S\ 5.17 of R. Calapso's
introduction to Vol {\bf 4}, part 1 of Bianchi's {\it Opere}).

In 1899 Guichard considered Q(W)C of revolution around the focal
axis: when such a quadric rolls on an applicable surface, its foci
(focus) describe(s) CMC (minimal) surfaces. Note that since a
surface of revolution can be rolled on a line with the arc-length
of the meridian corresponding to the arc-length of the line,
Guichard's result is a generalization of Delaunay's. Note also
that according to Bonnet when the unit sphere rolls on an
applicable surface its center describes a CMC $\frac{1}{2}$
surface, so Guichard's result is a metric-projective
generalization of Bonnet's result (compare with the rolling
problem for the ellipse). Note also that when a Q(W)C of
revolution around the focal axis rolls on itself (at each instant
the fixed quadric and the rolling quadric reflect in the common
tangent space) its foci (focus) describe(s) CMC (minimal)
surfaces: spheres (planes), so Guichard proved that this property
is essentially metric.

With Guichard's result the race to find the deformations of
general quadrics was on.

It ended in 1906 with Bianchi's discovery; in the process the
differential geometry underwent a fundamental change through the
results of geometers involved in the study of quadrics: mostly
Bianchi, Calapso, Darboux and Guichard, but a list of geometers
with results related to quadrics and not mentioned so far is
longer: Chieffi, Peterson, \c{T}i\c{t}eica (aka Tzitzeica), etc.
In particular Peterson discovered certain surfaces (later on
extensively studied by Bianchi and thus called {\it Bianchi
surfaces}; note however that The Arnold Principle applies here
only partially, as Bianchi gave due credit to Peterson in
\cite{B1} and he never calls them Bianchi surfaces; also
Peterson's work became widely known in western Europe only after
1903 (see Darboux (\cite{D1},\S\ 111))) and \c{T}i\c{t}eica (see
\cite{T3}) discovered the affine sphere (which led Blaschke to the
development of the affine geometry) while studying certain
deformations of quadrics. The {\it 'illustre geometra rumeno'}
\c{T}i\c{t}eica was a student of Darboux's (he got only a Master
degree in France) and the founder of the modern geometry school
from Romania a couple of decades after its independence, but in
Romania he is known even to the non-mathematician, mostly due to
his high-school level book of problems in Euclidean geometry
(which is the standard even today). Upon his untimely death from a
cerebral vascular accident in 1938 the two main pretenders to the
position of Chair of the Geometry Desk at the Mathematics
Department of the Bucharest University were Dan Barbilian and
Gheorghe Vr\^{a}nceanu. Dan Barbilian (mostly known for the {\it
Barbilian geometries}, a generalization of Cayley's point of view)
is also known in Romania as Ion Barbu, the pseudonym under which
he published successful poetry; again his poetry is required
subject of study in Romanian high-schools, so kids are tormented
to comment on his messed up poetry (what other type of poetry a
mathematician could create?). One could easily understand why the
apparently messed up moving picture {\it 'Pulp Fiction'} (the name
itself explains: juicy raw rough difficult to believe underground
story) is generous enough to admit interpretations in other
settings (for this reason {\it 'Pulp Fiction'} is actually a good
model of {\it 'Generic Essence of Pure Reality'}), because besides
the natural flow of logical deductions ({\it dialogue,
environment} and even {\it lyrics} further the story, in itself a
patch of more than two decades' worth of pop culture; the facets
of the infinitesimal picture (facets and bits of the story) fit
perfectly and consistently integrate (by means of a flat
connection form, namely QT) to the finite world thus providing the
correct explanation of the moving picture as a whole) and episodes
comparable with the elegance and simplicity of the proofs of the
wise ancients (but not apparent at the first cursory incursion:
just like good medicine and sometimes the bare truth tastes
bitter, {\it the simplest possible explanation is always the true
one}, no matter how confounding it is), Butch explaining the
character of his day ({\it 'without a doubt the single weirdest
... day of my life'}), whose chopper was his current chopper ({\it
'I had to crash that Honda ... might've broke my nose'}) and when
inquired {\it 'Who's Zed?'} he replies with the apparently wrong
answer {\it 'Zed's dead, baby, Zed's dead'} (according to The Wolf
{\it 'time is a factor'} when we have {\it 'a body in a car minus
a head in a garage'} (no need for parenthesis, change of order of
words or commas, as the whole sentence is associative, {\it
'court'} and self-contained and on top of that the garage can be
counted with either $\pm$), but the inquirer does not know it yet)
is the actual scientist of the XVII$^{\mathrm{th}}$ century in the
first few days after being hit with the apple on his head, when
the whole world collapsed with that apple and must then be
rearranged anew: {\it interesting things are simple, already exist
and can at most reveal themselves}. {\it 'Dead'} ain't no friend I
ever heard of, just as {\it 'What ain't no country I ever heard
of'}; should such a country exist the natural question {\it 'Do
they speak English in What?'} immediately arises. First he does
not believe and tries to disprove himself, but in doing so he
already changed the laws of the solar system and sees no error yet
and everything is {\it 'Flowers on The Walls'} ({\it 'Playing
Solitaire till dawn with a deck of 51'} is a good picture of
distributions of facets; if the croupier is fast enough to throw
the last card before the first one touches the casino table and a
snapshot is taken at that instant, then the perfect picture of a
space distribution of facets appears; also if all cards are
counted with their centers in one point and with 51 multiplicity,
then a perfect picture of Archimedes' method appears); then it
turns out that things are not so simple and rosy, so within
minutes (a longer time for the scientist) the first {\it shadow of
a doubt} appears; he tries to fend it off immediately ({\it 'It's
good to see you, I must go, I know I look a fright'}), but that
turns out not to be so easy, so he feels the bounds of his own
sanity within the reach of his bound hands and only then, in order
to preserve one's sanity, he lets go of the old views as
unimportant details, changes his own beliefs and unties his hands
with a supreme effort and ultimately rides chop-chop out of town
with the eery {\it 'Twilight Zone'} theme omnipresent and
providing the natural smooth change in time-line and on the first
available chopper ({\it 'Grace'}), not before paying to the Caesar
of the land (always bearing a Latin name; while bestowing {\it
'oak'} on friends also the prefix {\it 'Uncle'}) the debt he owed,
just as everybody else (for example Jules) also does before
leaving town, so as to be {\it 'Kool and The Gang'}. But then he
remains with fear of God for the remainder of his life, as he
cannot find other rational explanation for the falling of the
apple {\it precisely} when he was studying the movement of the
moon as it was explained by his predecessors (or equivalently,
while turning back for the forgotten heritage of his predecessors,
studying an unusual strange object at hand and finding {\it
without even thinking} its immediate and natural usefulness {\it
precisely} when the toast ripened and was thus released, {\it The
Brothers Startled} and one of them got smoke-toasted, thus
fulfilling in the process the previously established {\it 'grease
spot'} denomination as a result of his beliefs in a {\it 'freak
occurrence'}; most of the {\it hyperbolic} denominations stick to
their object and {\it parabolically} fulfill their purpose, even
if Butch apparently has no meaning in English). So he finds it
likely that God has created the world in a small number of days
making thrift use only of the ruler, the compass, the balance and
other simple mechanical methods (it seems that even God Himself
took some recreational breaks from this endeavor to play some
dice; in doing so some of the dice spilled and stuck to this
endeavor and the resulting ambiguity had to be decided among
others by means of {\it 'eenie, meeny, miney, moe'}) and provides
the explanation along these lines, passing through several
interpretations until he finds the one which he believes finally
fits the moving picture, just as Jules's {\it last minute}
interpretation of {\it Ezekiel 12:57} in the actual {\it bony
situation} (the story where some characters {\it 'get into
character'} and {\it 'have character'} returns in a full circle to
its true beginnings). Playing with the dice and things sticking
and transferring to their apparently wrong place is common not
only to DNA (where time and nature take over the resulting
ambiguity and decide its rightfulness) or to human language (all
languages are inherently ambiguous and continuously changing; as a
consequence the closest relative in Romanian to the English {\it
kitchen}, the French {\it cuisine} and the Italian {\it cucina} is
{\it cocina} (pig sty; most of the time the correct description of
my kitchen and probably the same explanation {\it parabolically}
applied in the Middle Ages), a fact revealed to me on a trip to
Chicago in Summer 2005 when I saw an Italian restaurant and began
laughing instead of in the first year of the study of French and
English), but in all aspects and facets of reality. When
instructed {\it 'Say something!'} or in the morning {\it 'Say
goodnight, Raquel!'} the instructed obeys the command shown in
front of one's eyes and replies ({\it 'quite rationally, in my
opinion'}) {\it 'Something!'} or {\it 'Goodnight, Raquel!'}. There
is another structured and mathematical story with a {\it Jack
Rabbit's Slim hole} (it is mostly things that shock that are not
forgotten, probably as a side effect of evolution): the {\it
'Matrix Trilogy'} is a good model of the {\it infinitesimal
world}, where Delaunay's surface (closely resembling UEFA's
Delaunay trophy; strange things appear in front of one's eyes when
one asks {\it Google, show me!} instead of {\it Open, Sesame!})
clearly appears as the trail of a bullet to make the point out of
it and to provide the link to {\it gas dynamics}; also some
information sticks (probably due to entropy) and a-priori
inexplicably consistently transfers from The One to Agent Smith.
Consequently I have a vague idea at this point about the
mathematical definition of {\it equivariancy}, but I remember a
talk where the speaker approached a member of the audience in the
first row, winked (always winking with the right eye is the trick;
understanding this trick puts the {\it I} in {\it Robot} and by
entropy conversely) and identified oneself: {\it 'I am your friend
and you are my friend; thus I treat your other friends right and
in return you will do the same'}; {\it consequently} the TC can be
interpreted as something to do with tangent cones (as seen from
$x_z$) taken by RMPIA's to cones of tangents (as seen from $x_0$),
so the RMPIA may be the wink by means of which {\it cones} and
{\it tangents} for quadrics befriended and sealed their friendship
with the symmetry of the B transformation (that is Archimedes'
balance) and the BPT. And by the way, {\it 'Antwan should've ...
better known better'} because Archimedes {\it 'is The One that
says 'Bad ...' '} (Marsellus himself had to learn the hard way
that one cannot defeat Archimedes directly, so he {\it 'took
advantage of a poorly guarded fortification which he had seen
during diplomatic negotiations'}; it may be a possible explanation
of Archimedes' unfortunate premature death or so they say {\it
'around the campfire'}): getting {\it 'the technique down, I don't
tickle or nothing'} is useless if one loses the focus and purpose
of the idea: {\it 'sticking'} facets with their centers at a point
and with an (close to) $\infty$ multiplicity may provide the
required {\it 'charming ... pig'} as the correct substitution to
the initial unattractive one. Also even if the Dutch decided to
put mayonnaise on French fries ({\it 'I've seen them do it, man.
They ... drown them in that ...'}), the choice of {\it 'ketchup'}
should be the natural one (there is a Romanian equivalent of the
joke with three tomatoes: when the junior one is warned by a
senior one to move from the road to the sidewalk because of the
danger posed by cars it replies {\it etee, fleo\c{s}c!}; this can
be translated as either {\it 'eat my shorts!'} or {\it 'eat my
splashhh!'}; on the other hand the Romanians also have another
saying: {\it 'Orice pasere pre limba ei piere'}; this can be
translated as either {\it 'Every bird dies singing its song'} or,
by entropy and not paying corresponding due attention to the cat,
{\it 'Every bird dies because of its song'}). Then he realizes
that when God said {\it 'Let there be light!'} it meant that God
decided to disclose His Plans to man ({\it 'at the request of'}
man, since out of his own vanity (the devil's favorite sin) man
chose the apple and later on went on A Mountain, repented and made
{\it A Covenant}), so he begins looking in religious texts
available at the time in hope of finding hidden messages from the
wise ancients. I would say that this is the most logical and
scientific step a deeply religious scientist of those times would
have made, since God or a Librarian of Alexandria or a combination
thereof, knowing that the age of the wise ancients was coming to
an end, would have most likely chosen a religious text to
partially hide and thus protect the knowledge in order to better
preserve it for the next generations. Also at that time the wise
ancients were regarded with the proper due respect (even during
his lifetime there was a general consensus that the wise ancients
reached a level of development of scientific investigation much
higher than the current one; in particular this became clear with
the almost complete version of Apollonius's {\it Conics} of 1710)
and it was still believed (as inferred from the available
literature) that an important part of their results was missing.
It turns out that indeed the Word of God (wisdom: note that
Archimedes was also hit with the apple; as a consequence he would
run naked on the streets shouting like a person out of one's mind)
was hidden in the Word of God (religious text), so {\it reading
between the lines of a religious text of the time was indeed the
best way to reach the wisdom of the wise ancients}. Indeed the
finding is simple: to the inquiry {\it 'You 'Flock of Seagulls'
... why don't you tell my man Vince here where you got the ... hid
at'} the answer is of course the simple and logical {\it 'It's in
the cupboard'}, but the final delivery is a little bit more
complicated: {\it 'I cannot give this case to you because it don't
belong to me'}. Vr\^{a}nceanu (probably the most famous Romanian
geometer) had been a student of Levi-Civita's and by that time
even Cartan had already built on his ideas. The next oral story
still circulates today at Bucharest University: it seems that
Barbilian formally filed his request before Vr\^{a}nceanu; after
hearing about Vr\^{a}nceanu formally filing the request (at that
time Vr\^{a}nceanu was Professor at Cern\u{a}u\c{t}i University),
Barbilian warned the Committee charged to choose \c{T}i\c{t}eica's
successor that the only way he would not become Chair and would
not contest the findings of such a Committee is for Vr\^{a}nceanu
to become Chair, so the Committee had no choice but to put
Vr\^{a}nceanu Chair; this position was held until his retirement
in 1970.

Bianchi generalized the B transformation of the pseudo-sphere to
the B transformation of quadrics: the same statement remains true
if in Lie's interpretation 'pseudo-sphere' is replaced with
'quadric' and 'circle' with 'conic'. The completion of the study
of the B transformation was realized by Calapso in 1912 by taking
into consideration the singular case $z=\infty$ for diagonal QC
(see \cite{C1}).

A simple heuristic argument in favor of the existence of the B
transformation for general quadrics goes as follows: the B
transformation for CGC $-1$ surfaces is of a general nature
(independent of the shape of the seed), so the cancellations of
the shape of the seed in the integrability condition of the
considered rolled distribution (equivalent to the complete
integrability condition of the Ricatti equation) occur because of
general equations (the GCM equations of the seed and of the
pseudo-sphere) coupled with algebraic consequences of the TC
(apparent in the static picture of confocal pseudo-spheres): these
are valid for general quadrics. Thus the fact that a general
quadric has less symmetries does not influence the existence of
the B transformation.

However, because of the smaller group of symmetries, the
applicability correspondence on CGC surfaces (which can be found
in $\infty ^3$ ways) does not have an easy generalization to
quadrics. Without knowledge of the applicability correspondence it
is very difficult in general to check if two surfaces are
applicable, the test being similar to and much older than the
Cartan-Ambrose-Singer theorem: basically one must find an
arc-length correspondence between same level curves of the Gau\ss\
curvatures such that their orthogonal trajectories also admit an
arc-length correspondence (see Eisenhart (\cite{E1},\S\ 136)).
This test answers a problem posed by Minding in 1839: {\it 'To
find a necessary and sufficient condition that two surfaces be
applicable'} (see Eisenhart \cite{E1}). Once the Minding problem
was solved, the problem posed by the French Academy naturally
appeared.

A first attempt is to prove that the correspondence between the
seed and the leaf given by the W congruence is the applicability
correspondence: the differential of the leaf depends on the shape
(but not on its derivatives) of the seed, but in the linear
element of the leaf the Gau\ss\ theorem does not cancel all the
shape of the seed. Thus to follow the general method and find the
applicability correspondence directly between the seed and the
leaf (or between a surface in the initial quadric and the leaf) is
out of question. According to Bianchi's own account from (122,\S\
IV), should the applicability correspondence exist in general, it
can be presumed to be of a general nature and thus independent of
the shape of the seed: therefore the answer lies not in the
picture with the seed and leaf, but in the static picture with
confocal quadrics and we have to roll back (or one can apply the
ABL method). Thus to find the applicability correspondence one
must look for a natural correspondence between confocal quadrics:
the Ivory affinity provides such a correspondence and proving that
the ACPIA is valid remained a matter of computations. This
applicability correspondence is realized as follows: if we roll
the seed $x^0$ on the applicable surface $x_0^0\subseteq x_0$,
then the point on the rolled B transform $x^1=B_z(x^0)$
corresponding to the point of tangency of $x_0^0$ and the rolled
surface $x^0$ is $x_z^1\in x_z$. Now $x_z^1$ is taken by the Ivory
affinity between $x_0,x_z$ into a point $x_0^1$ on the initial
quadric $x_0$: this is the applicability correspondence. Although
the applicability correspondence is restricted, the structure
becomes richer (a quadric with less symmetries has less
degeneracies; for example a segment of B transformations
degenerates to the complementary transformation for a quadric of
revolution or the geodesic flow on an ellipsoid degenerates to the
one on a sphere).

The necessary algebraic relations of the static picture do not
naturally appear in the static picture from a geometric point of
view, but are naturally chosen by the moving picture. For example
if we assume the B transformation to exist and we roll the given
quadric on a side of the seed, then we get two families of leaves,
corresponding to the two ruling families on the confocal quadric.
If we roll the initial quadric on the opposite face of the seed,
then we get the same two families of leaves (otherwise we would
have three or four ruling families on the confocal quadric), so
the two rolled distributions reflect in the tangent bundle of the
seed and this property is independent of the shape of the seed.
According to the ABL method, the reflection in the tangent bundle
of the seed property must be true also in the static picture, when
the seed coincides with the initial quadric and the leaves are the
ruling families on the confocal quadric. If one applies the same
ABL method to the ACPIA, then one naturally obtains the RMPIA.

Now it becomes clear why the initial attempt to prove the
applicability correspondence directly between the seed and the
leaf will be unsuccessful: the leaf may be applicable to a surface
$x_0^1\subseteq x_0$ different from $x_0^0$. For example by
applying the B transformation to a real seed $x^0$ applicable to
the real ellipsoid $x_0^0\subseteq x_0$ (and $z\in\mathbb{R}$
chosen so that the confocal quadric $x_z$ intersects
$\mathbb{R}^3$ along a hyperboloid with one sheet), the leafs are
real surfaces, but they are not applicable to the given ellipsoid:
they are applicable to another totally real region of $x_0$ (the
ACPIA becomes {\it ideal}, in Peterson's denomination). To get
back real surfaces applicable to the given ellipsoid we need to
apply the B transformation to the leaf also. This is the reason
why the B transformation for the sphere was developed much later,
as the first application of the B transformation to real CGC $1$
surfaces gives imaginary surfaces and one needs iteration of the B
transformation to get back real CGC $1$ surfaces.

Such is the simplicity of the answer to a question which eluded
geometers for a long time! Bianchi's contemporaries acclaim {\it
'Qui la geometria vive!'} and even in the 1950's, when Bianchi's
{\it Opere} were published, geometers still considered it one of
the most important discoveries of his career and the culmination
of the golden age of classical differential  geometry. It is very
sad that the undergraduate student hears today only of Bianchi's
identities (appearing in a four page article; it seems that Ricci
had previously discovered them while doing computations for his
newly discovered absolute calculus of tensors contraption, but he
did not see them important so as to publish them): G. Fubini, in
his eulogy appearing in part 1 of Bianchi (\cite{B2},Vol {\bf 1})
mentions Bianchi's identities only in a sentence on page 60,
crammed among four other results (including the list of spaces
from Thurston's geometrization program), but mentions in a
consistent footnote on page 70 a beautiful theorem where Bianchi
described all normal W congruences by putting an Euclidean
configuration into motion and under certain tangency requirements
(so basically {\it integrable systems} means nothing more and
nothing less than {\it Euclid in motion}; as to what Euclidean
configurations are chosen to behave well in the moving picture is
decided by {\it Archimedes' method of integration}).

There are tree main theorems in the theory of deformations of
quadrics, all due to Bianchi (see (122,\S\ I-XII)):

-the existence and inversion of the B transformation of quadrics
and the ACPIA:

{\it Theorem I

Every surface $x^0\subset \mathbb{C}^3$ applicable to a surface
$x_0^0\subseteq x_0$ ($x_0$ being a quadric) appears as a focal
surface of a $2$-dimensional family of W congruences, whose other
focal surfaces $x^1=B_{z}(x^0)$ are applicable, via the Ivory
affinity, to surfaces $x_0^1$ in the same quadric $x_0$. The
determination of these surfaces requires the integration of a
family of Riccatti equations depending on the parameter $z$.
Moreover, if we compose the inverse of the RMPIA with the rolling
of $x_0^0$ on $x^0$, then we obtain the rolling of $x_0^1$ on
$x^1$ and $x^0$ reveals itself as a $B_z$ transform of $x^1$;}

-the Bianchi Permutability Theorem:

{\it Theorem II

If $x^1=B_{z_1}(x^0),\ x^2=B_{z_2}(x^0)$, then one can find only
by algebraic computations a surface
$B_{z_2}(x^1)=x^3=B_{z_1}(x^2)$; thus $B_{z_2}\circ B_{z_1}
=B_{z_1}\circ B_{z_2}$ and once all $B$ transforms of the seed
$x^0$ are found, the $B$ transformation can be iterated using only
algebraic computations;}

-the Hazzidakis transformation:

{\it Theorem VIII

If a surface $x^0\subset\mathbb{C}^3$ is applicable to a surface
$x_0^0\subseteq x_0$ and the homography $H$ of Bianchi II takes
the confocal family $x_z$ to another confocal family $\ti x_{\ti
z},\ \ti z=\ti z(z),\ \ti z(0)=0$, then one infinitesimally knows
a surface $\ti x^0=H(x^0)$, called the Hazzidakis (H) transform of
$x^0$ and applicable to a surface $\ti x_0^0\subseteq\ti x_0$.
Moreover the H transformation commutes with the B transformation
($H\circ B_z=B_{\ti z}\circ H$) and the $B_{\ti z}(\ti x^0)$
transforms can be algebraically recovered from the knowledge of
$\ti x^0$ and $B_z(x^0)$}.

If we apply Theorem VIII to QC of revolution around the focal
axis, then this induces some transformation of CMC surfaces
described by the foci (Guichard). This transformation was
discovered earlier by Hazzidakis and is the reason why Bianchi
still calls the transformation for general quadrics Hazzidakis.
The current denomination for Hazzidakis type (algebraic)
transformations of local solutions of same or different integrable
systems is {\it Miura} transformations, as Miura discovered a
transformation between local solutions of the KdV equation and
local solutions of the {\it modified} KdV (mKdV) equation (see
Roger-Schieff (\cite{RS1}).

An important tool in proving Theorem VIII is the notion of {\it
surfaces conjugate in deformation}: pairs of non-flat
non-homothetic surfaces on which the asymptotic coordinates and
all {\it virtual asymptotic coordinates} (coordinates which can
become asymptotic on an applicable surface) correspond. Thus for
any surface applicable on the first one one infinitesimally knows
a surface applicable to the other. An immediate consequence of the
Codazzi-Mainardi equations, the change of Christoffel symbols for
the change of coordinates and the equation of geodesics is the
fact that 'conjugate in deformation' is equivalent to
correspondence of asymptotic coordinates and geodesics (and we
know that this is realized by the homography of Bianchi II);
conversely Bianchi in (\cite{B2},Vol {\bf 5},(87),(90)) and
Servant in {\cite{S11} proved that only deformations of quadrics
can appear as surfaces conjugate in deformation (see also Darboux
(\cite{D1},\S\ 603)); Bianchi proved the statement for quadrics of
revolution and Servant for the remaining quadrics.

The notion of surfaces conjugate in deformation has another
important application (Bianchi (122,\S\ XI,\S\ 56)): if we
consider $x_z$ for a particular value of $z$ as the Cayley's
absolute, then the Euclidean confocal family of $x_0$ remains
confocal family in this new geometry. Lines become geodesics (so
remain lines) and all properties of confocal quadrics are
preserved; in particular Euclidean geodesics and asymptotic
coordinates on $x_0$ remain geodesics and asymptotic coordinates
on $x_0$ in the new geometry, so all virtual asymptotic
coordinates on $x_0$ in the new geometry still correspond to
virtual asymptotic coordinates of $x_0$ in the initial Euclidean
geometry. Thus deformations of quadrics in other geometries are
reduced to deformations of quadrics in the Euclidean space.

With this construction the deformation of all quadrics in
$\mathbf{O}_3(\mathbb{C})$ except the flat Clifford torus can be
reduced to the deformation of quadrics in $\mathbb{C}^3$; all flat
surfaces in $\mathbf{O}_3(\mathbb{C})$ appear as rotations of
rollings of flat surfaces in $\mathbb{C}^3$.

Darboux split the prize of the French Academy between Bianchi and
Guichard for solving the problem for quadrics (Guichard had also
other results, including some transformations G of deformations of
general quadrics) and asked if there are any relationships between
the transformations found by Bianchi and Guichard. Bianchi showed
that certain G transformations are compositions of two B
transformations of different ruling families $B_z\circ B'_z$ with
finite $z$; Calapso showed in 1912 that the remaining G
transformations are compositions of two B transformations
$B_{\infty}\circ B'_{\infty}$ (see \cite{C1}). Thus most
transformations of deformations of quadrics were reduced to
compositions of B and H transformations. The problem was solved in
the sense that solutions depending on arbitrarily many constants
have been found: beginning with a seed surface we integrate a
$1$-dimensional family of Ricatti equations; after that only
algebraic computations are required to produce solutions (BPT).
Bianchi communicated to Calapso (see  \cite{C1}) that the
$B_{\infty}$ transforms are generated by some trivial
infinitesimal deformations (isotropic translations); therefore we
can find them by quadratures.

It thus becomes apparent that certain metric-projective properties
of confocal quadrics (most of them established in the first half
of the XIX$^{\mathrm{th}}$ century) {\it carry out} (stick and
transfer) by rolling to and influence deformations of quadrics and
surfaces geometrically linked to these, therefore providing a
wealth of integrable systems and projective transformations of
their solutions (the B transformation boils down to a Ricatti
equation, the cross-ratio of whose four solutions is constant),
directly linked to quadrics or closely related to these: see Moser
\cite{M1} for similarities between quadrics and other integrable
systems, Bianchi (\cite{B2},Vol {\bf 4},(108)) or Burstall
\cite{B4} and collaborators for isothermic surfaces, Bianchi
surfaces (first studied, according to Bianchi (\cite{B1},\S\
294-\S\ 295), by Peterson), Rogers-Schieff \cite{RS1} for
\c{T}i\c{t}eica surfaces, etc. For example the sine-Gordon
equation (with its many real forms) describes many types of
surfaces (with positive definite metric or not) in spaces of
constant curvature (with positive definite metric or not) and all
these are generated just by a particular type of confocal
quadrics: the complex spheres. In particular the deformation of
the (pseudo-)sphere is equivalent to the deformation of Darboux
quadrics (quadrics tangent at one point to $C(\infty)$; they
cannot be realized as quadrics in $\mathbb{R}^3$, although their
deformations can: see Darboux \cite{D2}), of paraboloids and of
quadrics of revolution (see Eisenhart (\cite{E2},\S\ 144));
deformations of general quadrics are put in relation with special
coordinates on the (pseudo-)sphere (Eisenhart \cite{E2}).

There is literature of deformations of $2$-dimensional quadrics in
higher dimensional spaces including up to the 1930's, due mostly
to Calapso and Guichard.

Cartan (see \cite{C2}) studied in 1918 deformations of space forms
in space forms; in particular he proved that there are no
deformations of the hyperbolic space form
$\mathbf{H}^n(\mathbb{R})$ in $\mathbb{R}^{n+p},\ p<n-1$ and all
deformations of $\mathbf{H}^n(\mathbb{R})$ in $\mathbb{R}^{2n-1}$
appear only with flat normal bundle. In 1972 Moore (see
\cite{TT1}) introduced generalized Tchebyshev coordinates on
deformations of $\mathbf{H}^n(\mathbb{R})$ in $\mathbb{R}^{2n-1}$.
Upon a suggestion of S.S. Chern and using Cartan-Moore, Tenenblat
and Terng (see \cite{TT1}) developed in 1980 the B transformation
for deformations of $\mathbf{H}^n(\mathbb{R})$ in
$\mathbb{R}^{2n-1}$ (and Terng developed in \cite{T2} the BPT for
such a B transformation). Note that a simple corollary of the ABL
method is the
$\mathbf{O}_{n-1}(\mathbb{R})\times\mathbf{O}_{n-1}(\mathbb{R})$
symmetry of the isoclinic facets for the Tenenblat-Terng B
transformation of deformations of $\mathbf{H}^n(\mathbb{R})$ in
$\mathbb{R}^{2n-1}$ (the rolling can be realized with the symmetry
$\mathbf{O}_{n-1}(\mathbb{R})$ of the normal bundle and one takes
into consideration the symmetry of the tangency configuration).
More recently Tenenblat, Terng, Uhlenbeck and collaborators (see
\cite{TU1} and its references for earlier results and names) have
developed in particular the B transformation of space forms (with
flat normal bundle) in space forms and in general have adapted a
machinery, Lie algebraic in character, to deal with such B and
Darboux transformations (different from the D transformation of
isothermic surfaces) for integrable systems, obtaining for example
variants of the Bianchi Permutability Theorem and $n$-soliton
formulae. Bobenko and Pinkall have developed in 1996 (see
\cite{BP1},\cite{BP2}) a theory of discrete isothermic and CGC
$-1$ surfaces, at the heart of which lies the BPT and Bianchi's
notion of M\"{o}bius configurations (composed with a Lie contact
transformation for discrete isothermic surfaces). While Bianchi
developed in (\cite{B2},Vol {\bf 5},(117)) moving M\"{o}bius
configurations for {\it B\"{a}cklund} (B) transformations of
general surfaces (the two focal surfaces of a W congruence are
considered B transforms one of the other), he was not aware of the
fact that pieces of discrete surfaces are subsets of M\"{o}bius
configurations. Probably discrete mathematics (present since the
wise ancient's times; it is by discrete mathematics that
Apollonius and Kepler completed their computations) was not seen
important, as there were no computers and much need for it.

Since most of Cartan's arguments are projective, the natural
question arises if a rich family of deformations in
$\mathbb{R}^{2n-1}$ of space-like regions in $n$-dimensional
quadrics exists. Such regions can appear as the ones in
$n$-dimensional quadrics of the Lorentz space of signature $(n,1)$
along which the normal is time-like (thus this linear element is
incomplete except for the pseudo-sphere, since there are isotropic
normal directions) or as quadrics in $\mathbb{R}^{n+1}$. Berger,
Bryant and Griffiths provided an affirmative answer in 1983 (see
\cite{BB}); in fact such regions of quadrics are the only
$n$-dimensional Riemannian manifolds whose family of deformations
in $\mathbb{R}^{2n-1}$ is as rich as possible: it depends on
$n(n-1)$ functions of one variable. Once this question was
settled, the natural problem of generalizing Tenenblat-Terng's B
transformation to such deformations appears. Note that the
$\mathbf{O}_{n-1}(\mathbb{R})\times\mathbf{O}_{n-1}(\mathbb{R})$
symmetry deduced from the ABL method remains true for general
quadrics, since it has to do with the normal bundle. One can study
a more general question by allowing the linear element of the
quadric to be indefinite, every component of its signature being
of course smaller than or equal to the corresponding component of
the signature of the ambient Lorentz space of dimension  $2n-1$,
but once the issues are clarified for a general quadric all other
cases should follow immediately.

\subsection{Statement of results}
\noindent

\noindent First a short history of the theory of deformations of
quadrics from my point of view; since at this point I have a vague
memory about its early part, the only sources are the Library's
stamped due dates (assumed to be roughly half an year after the
day books were borrowed), the talks I gave and some of my emails.
I became aware in early Fall 2002 of Eisenhart \cite{E1}; except
for the first and last chapters and most of the exercises (which
are sort of a short list of important facts from Bianchi \cite{B1}
and Darboux \cite{D1}) it was in general accessible. By early
Spring 2003 (after the qualifier) I also found Eisenhart \cite{E2}
and the theory of deformations of quadrics began to loom large, as
it always appeared in the last chapters and thus clearly beyond my
level of expertise. The find of Eisenhart \cite{E2} brought to my
attention isothermic surfaces (and their {\it Darboux} (D)
transformation) and W congruences (and their B transformation) as
important objects of investigation, but that book is very
technical and does not allow a discontinuous reading (it builds an
impressive apparatus from the beginning and provides precise
proofs later just by sending the reader to previously proved
results; by the time the reader reads the applications the
impressive apparatus is already forgotten): after a serious
reading of the first chapters I just jumped to a cursory reading
of the remaining ones and from the theory of deformations of
quadrics I was not able to understand much. During this time I was
also involved in studying some current literature, but the current
techniques and tools clearly appeared insufficient (I was involved
in a parallel reading of the classical literature and the current
one up to Summer 2005, after which time I mostly dropped the later
from my readings). By mid Spring 2003 I also found Darboux
(\cite{D1}, Vol {\bf 2}) and I realized that it is much cheaper to
read Darboux than to solve the problems of Eisenhart \cite{E1}, at
least until one gets used with the techniques (whole sections of
Darboux's are compressed as exercises in Eisenhart). In late
Spring 2003 the article Bianchi (\cite{B2},Vol {\bf 4},(108)) was
brought to my attention. Since at that time I was involved in a
dual classical-current literature reading, I searched on the same
afternoon {\it 'isothermic surfaces'} on Altavista (for some time
I use Google). In the internet age the best research is done by
use of search machines on the internet; of course there are also
leads which {\it apparently} have less to do with mathematics, but
in the first 20 entries one also finds the relevant literature
(this is how I found Moser \cite{M1} and later Arnold \cite{A1}).
Today it is very difficult to find Burstall \cite{B4} among the
first hits, mainly due to the size of the literature on isothermic
surfaces. After a cursory reading of Bianchi (\cite{B2},Vol {\bf
4},(108)) and the last section of Burstall \cite{B4}, I realized
again the importance of the theory of deformations of quadrics;
Bianchi seemed to provide in (\cite{B2},Vol {\bf 4},(108,\S\ 15))
an approach simpler than Eisenhart's and also the correspondence
between the D transformation of {\it special isothermic surfaces}
(the closure of {\it constant mean curvature} (CMC) surfaces under
conformal transformations of $\mathbb{C}^3$) and deformations of
quadrics provided the necessary bridge to show that quadrics put
their hands in a lot of cookie jars. The most successful
expedition in the Library turned out to be when I searched Bianchi
as author and I found his {\it Opere}; on the same day I grabbed a
few of his volumes (including the two parts Vol {\bf 4} concerning
the theory of deformations of quadrics), Vol {\bf 2} of Bianchi
\cite{B2} and Vol {\bf 3,4} of Darboux \cite{D1}. I do not believe
there has been more booty in an expedition, except probably when
the Crusaders first broke Constantinopole's walls, bringing
economic decline and thus easying the rise of the Ottoman Empire
and the fall of but the last remnants of the Eastern Roman Empire
two hundred years later . It is on the same day that I read and
understood Bianchi's theorems of the theory of deformations of
quadrics; understanding their proofs would take me one more year
and repeated unsuccessful attempts of scaling the walls of the
fortification (I was still trying at the time to use current tools
in finding alternative simpler proofs or new results, but I was
never leaving the infinitesimal neighborhood of the seed; in fact
my purpose at that time was to provide an explanation according to
the current tools (loop groups, solitons, etc) of the theory of
deformation of quadrics; as I have delved deeper into the
classical geometers I found other interesting things to do). Since
Guichard's result appeared to be easier in Bianchi's
interpretation than in that of Eisenhart \cite{E2}, by late Fall
2003 I had a good understanding of it; in Spring 2004 I decided to
seriously read Weingarten's method (one of the last chapters of
Eisenhart \cite{E1}), since it involved a particular case of
deformations of quadrics. At that time I already had a good
understanding of the rolling problem (including the trick with
infinitesimal deformations of the rotation of the rolling
generating the corresponding pairs of applicable surfaces) and a
small result of cyclic systems in space forms. I distinctly
remember considering the next point of view on the rolling and the
deformation problem: imagine a piece of paper; one deforms it
while preserving a virtual copy of the initial planar position
such that a given planar curve remains planar and nonsingular in
the initial plane (the crumpled piece of paper is allowed to
develop singularities outside a small tubular infinitesimal
neighborhood of the curve); the same trick but in higher
dimensions would be used in the deformation problem. However too
many variables were introduced and it seemed to me to always boil
down to the same difference of the two GCM equations (again I was
never leaving the infinitesimal neighborhood of the seed, but the
difference of the two GCM equations got a simple geometric
explanation). And even understanding deformations of sheets of
paper is still an open question: I have a vague memory of a talk
where a sheet of paper was deformed while keeping a linear
boundary fixed (two strong opposite laws): {\it parabolas} of
singularities developed near the linear boundary and in a {\it
fractal} (fractured) pattern, so I can honestly say that {\it
deformation into parabolic patterns} naturally appears in
mathematics as in every facet of reality and nature; when such a
pattern appears it usually comes in abundance and singular manner
(thus finite energies are discharged and distributed on
infinitesimal domains creating $\infty$ multiplicity), so as to
equilibrate the strong opposing forces without breaking the usual
smooth rules (the smooth part of the problem notices nothing
interesting, but according to Archimedes the simplest explanation
possible and thus the truth lies in the singular behavior; when in
doubt nature always {\it runs home to mama} and chooses quadrics
in order to simplify and expedite the non-linear problem). There
is a famous joke with a famous current mathematician stating that
he studies quadratic equations (to the dismay of the inquirer): it
is in studying iterations of quadratic functions over complex
numbers that the most famous fractal patterns naturally appear.

Being involved with teaching a course in Summer 2004, I cut down
on research and I concentrated instead on reading Bianchi's
proofs, which by that time I almost understood; as soon as I went
through his proofs I decided to write a set of notes to myself.
Still due to that Summer 2004 there is a {\it 'The Tonight Show'}
show of Jay Leno's to whom I may have thank for my current
research. So Jay Leno decided one night to invite a guest which
was a Professor of Mathematics; that fella came riding a car with
square wheels; of course it was a bumpy ride, but the car was
lifted and fitted on a special road: on that road the car had a
smooth ride. The road looked like a cycloid (I had just covered
the cycloid in the Calculus II course I was teaching at the time);
thus I wanted to see if that is true; it turns out that it uses
the most difficult and unnatural to remember integral formula for
trigonometric functions, namely
$\int\sec(s)ds=\ln(\sec(s)+\tan(s))+\mathcal{C}$, because the road
is given by $(\ln(\sec(s)+\tan(s)),-\sec(s)),\ 0\le
s\le\frac{\pi}{4}$ and then it extends by reflections in the
mirrors $x=0,x=\ln(1+\sqrt{2})$. So after the fact one can easily
see that the road cannot be the cycloid, as it turns at $90^\circ$
(as it should) instead of at $180^\circ$, but doing the
computations and a Maple moving picture (I applied the rolling I
already knew to curves) helped me later immediately do the Maple
moving picture for the elliptical wheel. While preparing the final
exam for the course taught in Summer 2004 and working in the same
time on the notes I noticed the first simplification in Bianchi's
proofs, which was something of the form
$du_1=u_{1u_0}du_0+u_{1v_0}dv_0+u_{1v_1}dv_1$ for
$u_1=u_1(u_0,v_0,v_1)$; so Bianchi was doing computations with the
{\it right hand side} (rhs) and I did them with the {\it left hand
side} (lhs) (the sinister side of the equation; also I thought at
some point that Bianchi was a member of numerous Academies having
to do with silviculture, since {\it forestiere} can be interpreted
not only as foreign, but also as something having to do with a
forest; thus {\it foreign} has a Latin root, a fact previously
unknown to me; having good knowledge of Romanian and none of
Italian I realized however the correct use of {\it forestiere}, as
opposed to {\it nostrane}; not to talk about {\it concorso
banditi} being probably an {\it admission contest} in a respected
elite group, such as {\it Scuola Normale Superiore}). As a
consequence by late August 2004 I saw the relations
(\ref{eq:simpas}) and the one line proof of the ACPIA and I
decided to stay on the theory of deformations of quadrics, which I
did since then. I remember being very optimistic in Fall 2004
about finishing a simplification of Bianchi's computations and
thus having the correct tools at hand to attack the higher
dimensional problem (by that time I already found Bobenko-Pinkall
\cite{BP1},\cite{BP2} and Tenenblat-Terng \cite{TT1}). It is then
that I realized that by considering imaginaries one can put
together all QC in a single case (Bianchi was using
parametrization by trigonometric functions for real ellipsoids and
hyperboloids with two sheets). However, I got bogged down in a
project of proving the algebraic computations of the BPT for QC
(\ref{eq:qc1}); after completing the difficult computations and
going further to prove the cross-ratio property I realized that
the cross-ratio property can be used to simplify the algebraic
computations of the BPT. By Spring 2005 I headed to the H
transformation; thus I followed Bianchi's example of a Darboux
quadric and realized that all quadrics can be treated in the tree
cases QC, QWC and IQWC and that the last two were similar. In late
Spring 2005 in trying to discuss the BPT for the $B'_z$
transformation I had to redo all computations for the other ruling
family: it is this how I found the simple geometric proofs of the
necessary algebraic relations (\ref{eq:simpas}) (the simple
geometric interpretations were clear enough from Summer 2004, but
their proof was not). Then I have jumped however to the most
difficult part of the algebraic discussion of the static picture
in higher dimensions (namely Bianchi II) without having the SJ
canonical form; just as in the Romanian {\it M\u{a}n\u{a}stirea
Arge\c{s}ului} story what was built one day was destroyed over
night. It is then that I realized the need of the SJ canonical
form, since even the proofs of the other classical theorems about
confocal quadrics would succumb if $\sqrt{I_n-zA}$ does not
commute with $A$ (this is an obvious statement for $z$ close to
$0$ but even the existence, definition and good name of
$\sqrt{I_n-zA}$ is questionable for $z$ away from $0$). Few
examples in dimension $n=4,5$ showed the proper attack strategy: a
global one, since the eigenvalue-by-eigenvalue strategy for
non-isotropic eigenvalues was the story with destroying what was
previously built. Once the attack strategy established everything
fell into place and all possible information was squeezed and
completely used (as expected). With the SJ canonical form in place
the proof of Bianchi II was amenable to attack, but again false
leads took their toll (even at this point there is the question of
a diagonalization process). By the end of Fall 2005 the proof of
the M\"{o}bius configuration $\mathcal{M}_3$ was attacked; as a
consequence the simpler notation of the M\"{o}bius configuration
$\mathcal{M}_2$ revealed itself and consequently the algebraic
computations of the M\"{o}bius configurations $\mathcal{M}_n$ also
immediately revealed themselves (there is still the question of
validating those formulae). For the most degenerate quadric from a
metric point of view IQWC (\ref{eq:iqwc2}) these computations boil
down to a determinant formulation, a mixture of Vandermonde
determinant and homographies (separate linearity in rulings); from
this example one can deduce the invariant algebraic formula for
$\mathcal{M}_n$ and assume that it is true for all other quadrics.
With use of Menelaus's theorem a simple geometric argument for all
cases is immediate. Still due to the identities involved in the
argument with Menelaus's theorem a proof of the fact that the
cross-ratio property is the consequence of the second iterated
tangency configuration was immediate; as a consequence what up to
that point was a non-rigorous argument on top of which a
complicated identity with $u_1:=\infty$ still had to be checked
became rigorous and the checking became unnecessary. All specific
algebraic computations were done with Maple; consequently I knew I
was on the wrong path (or I was using the wrong notation) when my
computer was crunching numbers and did not even want to shut down
when I was pressing the corresponding button and I knew I was on
the correct one when the computer was polite enough to confirm the
result in a few seconds for $n=3$ and up to 10 minutes for $n=6$
for QC (\ref{eq:qc1}); moreover the computer clearly showed me
many tricks on many other occasions and verified most of algebraic
identities appearing in these notes. In Spring 2006 with the
algebraic computations of the higher dimensional static picture
completed (of course it was and is still missing the algebraic
computations specifically needed for the theory of deformations of
quadrics, namely those obtained by the ABL method) I decided to
investigate Cartan \cite{C2}; up to that point it was unclear to
me if for example $\mathbf{H}^n(\mathbb{R})$ can appear as a
sub-manifold of $\mathbb{C}^{n+p}$ with $p<n-1$ other than the
pseudo-sphere itself. It is after seeing Cartan's mostly
projective arguments that I realized that the dimension $2n-1$
cannot be lowered even if one relaxes the restriction on not using
imaginaries. The strategy is clear: just like Lie find the leaves
when the seed $\mathbf{H}^n(\mathbb{R})$ actually coincides with
the actual pseudo-sphere in
$\mathbb{C}^{n+1}\subset\mathbb{C}^{2n-1}$, then just like Bianchi
take their metric-projective generalization (thus the leaves must
be cones in the tangent cone including possibly the most
degenerate case of lines) and in trying to prove the higher
dimensional theory of deformation of quadrics the necessary
algebraic identities of the static picture will reveal themselves.
However, just like the sine-Gordon equation is not best suited for
the Lie ansatz (I was calling it at that time approach), similarly
the Tenenblat-Terng approach is not best suited for the Lie ansatz
since it assumes real seed: one gets only $1$-solitons similarly
to getting $1$-solitons from the axis of the tractrix for CGC $-1$
surfaces. The find of the Berger-Bryant-Griffiths \cite{BB}
article not only confirmed the program, but also removed its most
difficult part (I should mention that this article was first
mentioned to me in a course on the isometric embedding problem in
Fall 2004, but I did not hear the word {\it quadric} on that
occasion, just the fact that it deals with some difficult
particular cases where the dimension of the surrounding space is
below $\frac{n(n+1)}{2}$). However any hope of understanding that
article without serious prerequisites of Algebraic Geometry is
unfounded in my opinion ({\it 'I lost it, Lou!'} when I was hit
with that Serre's stuff taken for granted), so a simpler proof of
the weaker part (the program I needed and intended to prove) that
$n$-dimensional totally real quadrics admit rich families of
deformations in $2n-1$ Lorentz spaces of corresponding signatures
and with flat normal bundles should be provided first (the use of
elliptic coordinates on the quadric and lines of curvatures on its
deformation should render it a perfect example of Cartan's
calculus, moving orthonormal frames and counting certain numbers
to provide solutions of the corresponding differential system and
the dimensionality of the solutions space). Just like Terng's
generalized sine-Gordon equation, use of proper coordinates should
reduce the problem of existence of seeds and of the B
transformation to a Frobenius (complete) integrability, thus
avoiding the use of the more complicated Cartan-K\"{a}hler
theorem. The Ivory affinity still provides the applicability
correspondence and the $n,2n-1=(n-1)+1+(n-1)$ dimensions fit
perfectly to uniquely provide in general the rigid motion taking a
facet centered on the confocal quadric (but not tangent to it; it
even in general lives outside the space $\mathbb{C}^{n+1}$ where
the quadric lives, so the easy generalization of the RMPIA to
higher dimensions is not the correct one) to the tangent facet
corresponding by the Ivory affinity and does the reverse trick to
the other pair of facets. The discussion of totally real cases
should close the $2$-dimensional discussion in Euclidean
surrounding space, but again the stronger statement that the only
surfaces in quadrics and having real linear element are totally
real took its toll and the discussion is incomplete (only the most
degenerate cases are completed; I believe there is a simpler
unifying approach using rolling, but that still eludes me at this
point). So this was the setting when I found Archimedes' opinion
on the matter: this provided further insight into the higher
dimensional discussion and a method to find {\it all} necessary
algebraic identities needed for the moving picture from the static
picture, without actually trying to complete the computations of
the moving picture in order to see what are those algebraic
relations. Once those algebraic relations are found again their
geometric proof should be easiest.

In \S\ \ref{sec:algpren} we are interested in bringing a quadric
to the canonical form (depending on as few constants as possible),
finding its confocal family and proving the classical theorems
(only Bianchi II presents some challenge). The computations are on
occasion lengthy, although they always require only elementary
linear algebra notions. While part of these results with proofs
are present in the quoted bibliography, we believe that the
unified complete treatment appears here for the first time. For
example Bianchi in (\cite{B2},Vol {\bf 5},(90)) proves Bianchi II
with $n>3$ only for diagonal QC; also although most (if not all)
of the totally real quadrics of $\mathbb{C}^3$ were known to the
classical geometers, we have not found a reference with an
exhaustive and organized discussion.

First we bring a symmetric complex matrix, via conjugation by a
rotation $R\in\mathbf{O}_n(\mathbb{C})$, to a form depending on as
few as possible constants (SJ form); it inherits the properties of
matrices in Jordan canonical form. Such a SJ matrix is formed by
SJ diagonal blocks $aI_p+J_p,\ p\ge 1,\ J_1:=0=0_{1,1}$; the SJ
blocks $J_p,\ p\ge 2$ are formed by augmenting either (according
to $p$ being either even or odd) of the blocks
$J_2:=\frac{1}{2}\begin{bmatrix}1&-i\\-i&-1\end{bmatrix},\
J_3:=\frac{1}{\sqrt{2}}\begin{bmatrix}0&0&1\\0&0&-i\\1&-i&0\end{bmatrix}$,
diagonally on the upper left and until the dimension $p$ is
reached, with blocks
$\begin{bmatrix}0_{2,2}&M\\M^T&0_{2,2}\end{bmatrix},\
M:=\frac{1}{2}\begin{bmatrix}1&i\\-i&1\end{bmatrix}$ (the right
lower $0_{2,2}$ block is superimposed by addition on the left
upper $2\times 2$ block of $J_p$ in order to get $J_{p+2}$). The
normalization $f_1:=\frac{e_1-ie_2}{\sqrt{2}}$ instead of
$f_1:=e_1-ie_2$ for the first standard isotropic vector (the
$e_j,\ j=1,...,n$ are the standard basis: $e_j^Te_k=\del_{j,k}$)
is preferred in order to avoid cumbersome powers of $2$ when
working with SJ blocks.

Although a SJ matrix as defined above efficiently encodes all
metric invariants of a symmetric complex matrix $A$ (its
eigenvalues $\la$ and orders of cyclic vectors of $A-\la I_n$), we
prefer when convenient to work with the {\it SJ type} (the block
$J_p$ is replaced with a polynomial $P(J_p)$ in $J_p$ with
$P'(0)\neq 0$), since it is closed under needed algebraic
operations; the stronger SJ notion would require a rotation after
each such algebraic operation.

For $z=\infty$  or inverses of non-zero eigenvalues of $A$ we have
singular quadrics of the confocal family; the singular set
$\mathcal{S}(x_z)\subset x_z$ for $x_z$ singular is a lower
dimensional quadric. These singular sets  play an important
r\^{o}le in the proof of Bianchi II, in the course of which we
provide a classification of $\mathbf{PGL}_n(\mathbb{C})$ suited to
it.

We close \S\ \ref{sec:algpren} with a parametrization of confocal
IQWC (this allows us to consider them as metrically degenerated
QWC and thus to treat later all (I)QWC as a single case) and with
a discussion of totally real quadrics.

In \S\ \ref{sec:algpre2} we restrict our attention to quadrics in
$\mathbb{C}^3$; in \S\ \ref{subsec:algpre21} we present a list of
all confocal quadrics in $\mathbb{C}^3$, their ruling families
parametrization invariant under the Ivory affinity and
homographies of Bianchi II. Bianchi in (\cite{B1},\S\ 425) and
(122,\S\ 51-\S\ 57) uses only the homographies which take QC
(\ref{eq:qc1}) to QC (\ref{eq:qc1}) and QC (\ref{eq:qc2}) to QWC
(\ref{eq:qwc1}), although the general geometric method from
(\cite{B2},Vol {\bf 5},(90)) applies to all cases.

In \S\ \ref{subsec:algpre22} we discuss totally real quadrics; the
part of the argument where one proves that a surface in a quadric
and having real (valued) linear element must be a totally real
quadric may be replaced in \S\ \ref{subsec:deformations6} with a
simpler unifying argument using rolling.

While currently it is well known that there is a single type of
quadrics from a complex projective point of view, two types of
quadrics from a complex affine point of view (namely with or
without center), two types of real quadrics from a real projective
point of view (namely ruled or not ruled), 5 types of real
quadrics from a real affine point of view (namely three with
center: connected ruled, connected not ruled and disconnected not
ruled; two without center: ruled and not ruled) and 12 types of
real quadrics from a real metric point of view (all subcases of
the previous 5), it is less known since the classical times that
there are 11 types of quadrics from a complex metric point of view
(all subcases of the 7 cases QC (\ref{eq:qc1})-(\ref{eq:qc3}),
QWC(\ref{eq:qwc1})-(\ref{eq:qwc2}), IQWC
(\ref{eq:iqwc1})-(\ref{eq:iqwc2})) and numerous types of totally
real quadrics. The number of types of totally real quadrics is
even since multiplication by $i$ changes the signature, both of
the linear element of the totally real surface and of the ambient
totally real space, so it is enough to count the 12 types of real
quadrics plus the types of quadrics in
$\mathbb{R}^2\times(i\mathbb{R})$ and double this number.

A simple heuristic argument in favor of counting all types of
totally real linear elements of a complex metric type which is not
a complex sphere is to double the number of components of singular
quadrics of the confocal family (including the singular quadric at
infinity), since each type of totally real linear element comes
with its totally real confocal family, which contains singular
linear elements. The change in the type of linear element occurs
along the singular quadric; thus the change of the type of linear
element due to the multiplication by $i$ occurs at infinity. Since
most families of confocal quadrics are generated by finite conics
$c$ (the family being the adjugate of the pencil generated by the
adjugate of $c$ and the adjugate of $C(\infty)$), it is expected
that all information about such a family is recorded by the pair
$(c,C(\infty))$.

Take for example the real confocal family of a general ellipsoid;
the change in the real confocal family from ellipsoids to
hyperboloids with one sheet occurs along a real singular quadric
of the real family: a mostly doubly covered real principal plane.
In that real singular quadric the singular part (covered once) is
an ellipse; it separates the interior convex region (parameterized
by complex conjugate parameters and thus having imaginary rulings)
from the exterior concave region (parameterized by real parameters
and thus having real rulings); the two parametrizations coincide
on the ellipse. Therefore the ellipsoids of the real confocal
family admit parametrization by complex conjugate parameters and
have imaginary rulings, while the hyperboloids with one sheet of
the real confocal family admit parametrization by real parameters
and have real rulings.

Note that not only the types of totally real quadrics are distinct
(one cannot be obtained from another by means of a complex rigid
motion), but also their linear elements are distinct, since two
quadrics with same linear element can be rolled one onto the other
and the rolling with the determinant of its rotation being $1$
must be a rigid motion. Thus all real symmetries of the linear
element of a totally real quadric are in correspondence with rigid
motions preserving the quadric; the complex symmetries are
accounted by changing the type of the totally real quadric within
the same complex metric type of quadric and they must preserve the
character of real or imaginary rulings.

Since the B transformation is essentially projective, we must use
projective coordinates (either asymptotic or conjugate systems)
either on $x_0^0$ or on $x^0$ (the sine-Gordon equation
$\om_{uv}=\sin\om$ for deformations of the pseudo-sphere is
obtained by considering the net of arc-length asymptotes on
$x^0$). Due to the essential r\^{o}le played by the Ivory affinity
in the B transformation of quadrics, the use of rulings on $x_0$
(asymptotic coordinates on $x_0^0$) is more appropriate.

Beginning with \S\ \ref{subsec:algpre23} we present the algebraic
identities necessary for the theory of deformations of quadrics.
While they and their proof are simple and on occasion elementary,
the natural character of these relations is definitely not
obvious; these relations and their necessity become obvious once
certain natural ansatzs are made in the moving part (such as the
ABL method, the ACPIA or ruled surfaces for an a-priori valid
BPT).

The algebraic identities (\ref{eq:simpas}), their equivalency to
the TC and their simple geometric interpretation do not appear in
Bianchi's account of deformations of quadrics, although he uses
equivalent algebraic identities (for example (29) in (122,\S\ 4))
and the place of use is the same: bring the differential system of
the B transformation to a Ricatti equation, check the complete
integrability of this Ricatti equation and prove the ACPIA.
Bianchi's approach to the ACPIA in (122,\S\ 7,\S\ 8,\S\ 9) is to
take $dv_1$ from the (\ref{eq:dify2}), $du_1$ from (\ref{eq:du1})
and to show that
$|dx^1|^2-|dx_z^1|_{v_1=\mathrm{const}}|^2=|dx_0^1|^2-|dx_z^1|_{v_1=\mathrm{const}}|^2$
as symmetric quadratic forms in $du_0,\ dv_0$; in (122,\S\ 31,\S\
32) he takes $du_1$ from (\ref{eq:du1}), imposes the ACPIA and
obtains the Ricatti equation (\ref{eq:dify2}) (in both instances
he uses the Gau\ss\ theorem to cancel the influence of the shape
of $x^0$ with the influence of the shape of $x_0^0$). Bianchi's
choice to prove $|dx^1|^2-|dx_z^1|_{v_1=\mathrm{const}}|^2
=|dx_0^1|^2-|dx_z^1|_{v_1=\mathrm{const}}|^2$ instead of just
$|dx^1|^2=|dx_0^1|^2$ becomes clear for the following reasons:

(I) For a function $f=f(u_0,v_0,v_1(u_0,v_0)),\ du_0\wedge
dv_0\neq 0$ he considered for example the explicit derivative
$\pa_{u_0}f=\frac{\pa f}{\pa u_0}$ of $f$ wrt $u_0$ as opposed to
the implicit derivative $\pa_{v_1}fv_{1u_0}$ of $f$ wrt $u_0$; the
total derivative is
$f_{u_0}=\frac{df}{du_0}=\pa_{u_0}f+\pa_{v_1}fv_{1u_0}=
f_{u_0}|_{v_1=\mathrm{const}}+\pa_{v_1}fv_{1u_0}$.

(II) When $x^0$ coincides with $x_0^0$ the leaves of the
considered integrable distribution are the rulings
$v_1=\mathrm{const}$, so all formulae already found for $x^1$ in
general position apply to $x_z^1$ by considering
$v_1=\mathrm{const}$ and replacing the shape of $x^0$ with the
shape  of $x_0^0$.

To simplify the notation we shall use $f_u$ with the meaning of
$\pa_uf$ when the explicit dependence of $f$ on $u$ is clear
(other variables appearing in the definition of $f$ do not depend
on $u$); for example in the above situation we assume first
$du_0\wedge dv_0\wedge dv_1\neq 0$, derive
$df=f_{u_0}du_0+f_{v_0}dv_0+f_{v_1}dv_1$ and after restricting to
a surface $v_1=v_1(u_0,v_0)$ we have further
$dv_1=v_{1u_0}du_0+v_{1v_0}dv_0$.

We propose (in \S\ \ref{sec:deformations}) a different approach:
show that $|dx^1|^2$ and $|dx_0^1|^2$ are equal as symmetric
quadratic forms in $du_1,\ dv_1$; in doing so the computations are
reduced almost to a tautology and the necessary relations
(\ref{eq:simpas}) naturally appear from an analytic point of view;
they will be consistently used throughout the remaining sections.
Note that it is very likely that Bianchi saw the simpler proof
using only $du_1,\ dv_1$ and without expanding these as $1$-forms
in $du_0,\ dv_0$, since the needed rudiments of Cartan's exterior
calculus and of Ricci's absolute differential calculus of tensors
were present at the time for almost a decade, but he probably felt
more secure in providing an explanation according to the standards
of proof of the time; it took another decade for Ricci's absolute
differential calculus of tensors and several more decades for
Cartan's exterior calculus to become mainstream mathematics.

The simplest proof of the ACPIA is not almost a tautology, but an
actual tautology; it is essentially due to Bianchi, although he
preferred the security of an analytic confirmation to the power of
his geometric arguments. To prove the inversion of the B
transformation in (122,\S\ 11), Bianchi defines the rolling of
$x_0^1$ on $x^1$ as the composition of the RMPIA with the rolling
of $x_0^0$ on $x^0$; he also checks that the $v_1$ ruling goes to
the correct place (if a family of rigid motions takes a surface
and its tangent spaces to another surface and its tangent spaces,
it is not a-priori a rolling, as it may rotate the tangent spaces
into themselves). In doing so he used the ACPIA and found a simple
formula for the rolling $x_0^1$ on $x^1$, but he actually found
the simplest proof for the ACPIA, since existence of rolling
implies applicability correspondence (one needs only check the
complete integrability of the Ricatti equation (\ref{eq:dify2}),
since we need the leaf $x^1$ to exist).

While Bianchi knew in (122,\S\ 11) that by changing the ruling
family on $x_z^1$ the rotation of the RMPIA must acquire $-1$
determinant (and thus must be a rotation of determinant $1$
composed with a reflection), he was not aware of the simple
observation that if the required reflection takes place in the
tangent space of $x_0^0$, then the rotation does not otherwise
change at all (to change the rulings $u_0\leftrightarrow v_0$ he
uses instead a reflection $r_{e_2}$ in the $e_2^Tx=0$ plane for
$x_0^0$ being QWC (\ref{eq:qwc1})). In fact initially we were also
not aware of this observation and thus provided an analytic
confirmation of (\ref{eq:simpas}); it is only after considering
the algebraic computations of the BPT that we realized that a good
accounting of all RMPIA's for all rulings is needed. Using this
observation a simple geometric proof of (\ref{eq:simpas}) is
immediate; as a result we have gained much respect for what
Bianchi calls {\it 'elementary properties of the Ivory affinity'}
and consistent use of symmetries is carried through the remaining
sections*. \footnote{* After noticing the link to {\it The Method}
of Archimedes this observation lost most of its initial purpose.}
We introduce four quadratic quantities $\Del^{\pm},\ \Del'^{\pm}$
which are roughly the angles between the normal fields $m^1_0,\
m'^1_0$ of the considered distributions and the rulings of $x_0$
at $x_0^0$ with the denominator parts removed; the normal fields
$m^1_0,\ m'^1_0$ themselves have the usual definition (involving a
cross product) with the denominator part removed; further we
derive some useful algebraic identities. Note that these analytic
computations are not necessary for the B transformation for finite
$z$, but we need them for the case $z=\infty$, since in that case
the geometric picture disappears and we can only use analytic
computations.

Bianchi derives the properties of the {\it second iteration of the
tangency configuration} (SITC) and proves the BPT for ruled
deformations of QWC (\ref{eq:qwc1}) in (122,\S\ 21-\S\ 25); as a
corollary he proves in (122,\S\ 26) the beautiful cross-ratio
property (the BPT is essentially a cross-ratio theorem; in fact it
roughly boils down to the product of a cross-ratio property and
same cross-ratio property composed with a simple symmetry). We
propose a change in Bianchi's approach: consider the cross ratio
property from the beginning, to be proven together with the other
algebraic properties of the SITC; this interlacing simplifies the
computations of the SITC (we need this simplification especially
for QC (\ref{eq:qc1}), where even with the use of symmetries
computations will become involved at some point). Bianchi uses
directly (\ref{eq:v1v2}) and a counterpart as a homography between
$u_0,u_3$, which is inconvenient because it depends quadratically
on $v_1$. With moderate effort and use of simple symmetries we
propose to replace (\ref{eq:v1v2}) and its counterpart with a
single relation (\ref{eq:v1v21}), which is separately linear in
all variables, just like the cross-ratio property
(\ref{eq:cross}). The fact that both (\ref{eq:v1v21}) and
(\ref{eq:cross}) must be true imposes the homography
(\ref{eq:phi2c}) which reveals the necessary 'algebraic
computations' of the BPT. Our formula (\ref{eq:phi2c}) for QWC
(\ref{eq:qwc1}) coincides with (III) in Bianchi (122,\S\ 23); as
for the remaining quadrics Bianchi just claims that such formulae
exist by similar computations. He is not far from the truth, since
any projective property (like the cross-ratio property) which is
valid for a particular quadric is valid for all other quadrics
(all quadrics are projectively equivalent), but since in building
the confocal family of a given quadric Cayley's absolute
$C(\infty)$ plays an important r\^{o}le, it remains to prove that
the cross-ratio property is also independent of a choice of a
conic in a hyperplane of $\mathbb{CP}^3$ as Cayley's absolute
(probably a purely synthetic proof unifying all cases exists). It
is this what we did; in fact we did more: we found the necessary
'algebraic computations' of the BPT for all quadrics. It may be
the case that Bianchi, who first proved the BPT for the
pseudo-sphere in 1892 (in which case the proof is easier), has
noticed by chance (or looking on purpose for projective
properties) the cross-ratio property. Once the cross-ratio
property for pseudo-spheres is found it is natural to conjecture
that it is true for all quadrics and indeed it is just a matter of
computations to prove its validity for general quadrics. Bianchi
(122,\S\ 25) assumes the use of the algebraic identity
(\ref{eq:calib}) in the proof of the BPT without explicitly
stating it: {\it 'Mais toutes ces conditions, ... , ne peuvent
\^{e}tre remplies que si la deformation se r\'{e}duit \`{a} un
mouvement invariable, ...'}. Its use in the proof of the BPT can
be easily justified with a differential symmetry (it is probably
what Bianchi had in mind), but to  err on the safe side we prefer
to use only algebraic symmetries. In fact should we accept
differential symmetries just as easily as algebraic symmetries,
the BPT would require no proof beyond algebraic computations,
since it is essentially a statement about differential symmetries.

We conclude \S\ \ref{sec:algpre2} with further iterations of the
TC: M\"{o}bius configurations $\mathcal{M}_n$ (the TC is
$\mathcal{M}_1$ and the SITC is $\mathcal{M}_2$). Again Bianchi's
computations from (\cite{B2},Vol {\bf 5},(117)) proves their
existence and assures that they can be obtained only from
algebraic computations; we provide a simple geometric argument for
their existence (using for example the Menelaus theorem) and
provide the actual algebraic computations. As chance has it, I was
reading at leisure Boyer \cite{B3} during my office hours around
the days I was trying to find the algebraic computations
supporting the existence of $\mathcal{M}_3$; it is then that I
found the Menelaus theorem and I realized that it was what I
needed. The algebraic identity (\ref{eq:calib}) naturally appears
and plays an important r\^{o}le in the analytic computations of
M\"{o}bius configurations, thus closing the list of {\it all}
algebraic identities necessary for the moving part. M\"{o}bius
configurations will play an important r\^{o}le in the discussion
of DDQ and thus in creating computer images of deformations of
quadrics.

Beginning with \S\ \ref{sec:rolling} we present basic results and
formulae about the moving part. While most of these results are
due to the classical geometers, the compact notation makes them
more comprehensible and other issues become more clear.

On occasion there are small observations which were not apparent
to the classical geometers; for example Bianchi did not notice
that fact that pairs applicable surfaces in $\mathbb{C}^3$ are in
correspondence with infinitesimal deformations of surfaces in
$\mathbf{O}_3(\mathbb{C})$ (which become the rotations of the
rollings), although he was in possession of that differential
equation and of other geometric interpretations.

In \S\ \ref{sec:deformations}, \ref{sec:particular} we present the
important issues (in my opinion at this time) of the moving part
for deformations of quadrics in $\mathbb{C}^3$ and respectively in
higher dimensions. Note that many important issues from Bianchi
are not covered; for example the {\it permanent conjugate system}
(the conjugate system common to the surface in the considered
quadric and to its seed deformation; it is also the permanent
conjugate system on the leaf). The permanent conjugate system is
best suited to describe the B transformation at the analytic level
(this was Bianchi's original approach to the theory of
deformations of quadrics: he began with Guichard's result and
found the analytic formulation of the B transformation for
quadrics of revolution and diagonal paraboloids; take for example
the article (\cite{B2}, Vol {\bf 4},(107)) of 1905 where
everything was in place except the static picture of confocal
quadrics and the Ivory affinity). Also it is best suited to
describe triply conjugate systems of surfaces containing a family
of deformations of quadrics and it may be the correct tool in
attacking the higher dimensional problem of conjugate systems
containing a family of deformations of quadrics (orthogonal
systems are clearly present in the current literature as examples
of integrable systems; see for example \cite{TU2}).

\section{Algebraic preparatives for quadrics in
$\mathbb{C}^n$}\label{sec:algpren} \setcounter{equation}{0}

\subsection{Confocal quadrics in canonical form}
\noindent

\noindent Consider a non-degenerate symmetric bilinear form
(scalar product) on $\mathbb{C}^n:\ \ \ \ <x,y>:=x^TQy,\
Q=Q^T\in\mathbf{GL}_n(\mathbb{C})$. Since the map
$\mathbf{GL}_n(\mathbb{C})\ni A\mapsto A^TA$ has kernel
$\mathbf{O}_n(\mathbb{C})$, it has $\frac{n(n+1)}{2}$-dimensional
connected image ($\mathbf{GL}_n(\mathbb{C})$ is connected), so
$Q=A^TA$ for some $A\in\mathbf{GL}_n(\mathbb{C})$ (the set of
symmetric matrices in $\mathbf{GL}_n(\mathbb{C})$ is connected);
equivalently $A$ is the product of row (and thus $A^T$ the product
of same column) operations which bring $Q$ to the reduced echelon
form. Under the linear transformation $\mathbb{C}^n\ni x\mapsto
Ax$, the scalar product can be considered the Euclidean one:
$<x,y>:=x^Ty,\ |x|^2:=x^Tx$. Consider the standard basis
$\{e_j\}_{j=1,...,n},\ e_j^Te_k=\del_{jk}$ of $\mathbb{C}^n$ and
$f_j:=\frac{e_{2j-1}-ie_{2j}}{\sqrt{2}},\ j=1,...,[\frac{n}{2}],\
I_{j,k}:=\Sigma_{l=j}^ke_le_l^T$.

Consider a quadric $x\subset \mathbb{C}^n:\ \
Q(x):=x^T(Ax+2B)+C=0,\ A=A^T\in\mathbf{M}_n(\mathbb{C}),\
B\in\mathbb{C}^n,\ C\in\mathbb{C},\
\begin{vmatrix}A&B\\B^T&C\end{vmatrix}\neq 0$.

If $\ker(A)=0$, then $\begin{vmatrix}A&B\\B^T&C\end{vmatrix}
\begin{vmatrix}I_n&-A^{-1}B\\0_{1,n}&1\end{vmatrix}=
\begin{vmatrix}A&0_{n,1}\\B^T&Q(-A^{-1}B)\end{vmatrix}$,
so we have $Q(-A^{-1}B)\neq 0$. These are QC, since they have
center of symmetry $o$ (which remains center of symmetry under
rigid motions and homotheties): $x\in Q\Leftrightarrow 2o-x\in Q$,
equivalent to $Q(2o)=Q(0),\ Ao+B=0$ and  $o:=-A^{-1}B$ satisfies
both.

If $M=M^T\in\mathbf{M}_n(\mathbb{C}),\ \ker(M)=F\mathbb{C}^j\oplus
E\mathbb{C}^k,\ F:=[f_1...f_j],\ E:=[e_{2j+1}...e_{2j+k}]$, then
$M$ restricted to $\bar F\mathbb{C}^j\oplus
I_{2j+k+1,n}\mathbb{C}^n$ and co-restricted to
$M(\mathbb{C}^n)=F\mathbb{C}^j\oplus I_{2j+k+1,n}\mathbb{C}^n$ is
invertible. Note $F^TF=0,\ \bar F^TF=I_j,\ \bar FF^T+F\bar
F^T=I_{1,2j},\ I_{2j+1,n}F=0$ (the last $n-2j$ rows of $F$ are
$0$). Denote $M^{\circ -1}$ the inverse of the restricted and
co-restricted $M$, extended as $0$ on $\bar F\mathbb{C}^j\oplus
E\mathbb{C}^k:\ M^{\circ -1}:=(M+\bar F\bar
F^T+EE^T)^{-1}-FF^T-EE^T,\ MM^{\circ -1}=I_n-\bar FF^T-EE^T,\
MM^{\circ -1}M=M,\ M^{\circ -1}MM^{\circ -1}=M^{\circ -1}$ and
$Mv=e^{-a}v,\ v\in\mathbb{C}^n\setminus\{0\} \Leftrightarrow
M^{\circ -1}v=e^av$. As we shall see later, the kernel of any
symmetric matrix can be brought to the form above by applying
rotations $R\in\mathbf{O}_n(\mathbb{C})$, so we can define
$M^{\circ -1}:=R^T(RMR^T)^{\circ -1}R$ for general
$M=M^T\in\mathbf{M}_n(\mathbb{C}),\
\ker(RMR^T)=F\mathbb{C}^k\oplus E\mathbb{C}^p$. Although $M^{\circ
-1}$ may depend of the choice of $R$, we still have $MM^{\circ
-1}M=M,\ M^{\circ -1}MM^{\circ -1}=M^{\circ -1}$ and $Mv=e^{-a}v,\
v\in\mathbb{C}^n\setminus\{0\}\Leftrightarrow M^{\circ -1}v=e^av$.

If $\ker(A)=\mathbb{C}v,\ v\in\mathbb{C}^n\setminus\{0\}$, then
$v^TB\neq 0$ (otherwise $[v^T\
0]^T\in\ker\begin{bmatrix}A&B\\B^T&C\end{bmatrix}$). Such $Q$ does
not have center and is thus (I)QWC according to $|v|^2(=0)\neq 0$.
However, it has another set of interest, the vertex $o$: for QWC
it satisfies $A^{\circ -1}(Ao+B)=0,\ Q(o)=0$, equivalent to
$o=-A^{\circ -1}B+ \frac{B^TA^{\circ -1}B-C}{2v^TB}v$ ($o$ is a
well defined point and remains vertex under rigid motions:
$o\circ(R,t)=(R,t)\circ o$ and homotheties). A similar formula for
IQWC fails, since $A^{\circ -1}$ (and thus $o$) depend on the
rotation $R$ with $\ker(RAR^T)=\mathbb{C}f_1$, but since $A^{\circ
-1}$ always has isotropic $1$-dimensional kernel (in fact as we
shall see later, if $|w|^2=0,\ w^Tv\neq 0$, then we can make
$\ker(A^{\circ -1})=\mathbb{C}w$), we can define the vertex of an
IQWC to be the co-dimension $1$ sub-manifold of $Q$ along which
normal directions are isotropic:
$\mathcal{O}:=\{o\in\mathbb{C}^n|\ |Ao+B|^2=Q(o)=0\}$.
$\mathcal{O}$ remains vertex under rigid motions and homotheties;
a canonical choice of a vertex $o\in\mathcal{O}$ will appear
later.

If $\dim(\ker(A))\ge 2$, then
$\begin{vmatrix}A&B\\B^T&C\end{vmatrix}=0$.

We would like to reduce $x$, by applying rigid motions
$(R,t)\in\mathbf{O}_n(\mathbb{C})\ltimes\mathbb{C}^n,\
(R,t)x:=Rx+t$, to the canonical form (its equation depending on as
few constants as possible); equivalently prescribing a coset of
the action of rigid motions on quadrics. Once we prescribe the
rotation $R$ such that $RAR^T$ depends on as few constants as
possible (SJ form) the translation will be $-o$, so as to bring
the center (vertex) to the origin.

Complete a maximal set of orthonormal eigenvectors of $A$ to an
orthonormal basis of $\mathbb{C}^n$ and consider the vectors of
this basis as the rows of the rotation $R$. $RA R^T$ is decomposed
into two blocks on the diagonal, one block being a diagonal matrix
and the other a symmetric matrix with isotropic eigenspaces; we
have thus reduced the problem to symmetric matrices $A$ with
isotropic eigenspaces.

Consider an isotropic subspace $W\mathbb{C}^p\subset\mathbb{C}^n,\
W:=[w_1...w_p]$ of rank $p$; since $W^TW=0$ we have $w_1^TW=0$.
Because of the transitive action of $\mathbf{O}_n(\mathbb{C})$ on
the isotropic cone $|x|^2=0$, we can make $w_1=f_1$; replacing
$WI_{2,p}$ with $W-f_1(\bar f_1^TW)=I_{3,n}W$ we make
$WI_{2,p}=I_{3,n}WI_{2,p}$; by induction hypothesis we can make
$W=F:=[f_1...f_p]$.

Just as for the diagonalization process for symmetric matrices
with non-isotropic eigenvectors, we would like to find rotations
which preserve isotropic subspaces of $\mathbb{C}^n$ (the
translations are the subspaces themselves). A rotation $R$
preserving $F\mathbb{C}^p$ must satisfy $RF=FL,\
L\in\mathbf{GL}_p(\mathbb{C})$; with
$R':=I_{1,2p}+I_{2p+1,n}RI_{2p+1,n}\in\mathbf{M}_n(\mathbb{C}),\
M:=-\bar F^TR^TI_{2p+1,n}\in\mathbf{M}_{p,n}(\mathbb{C}),\ N:=\bar
F^TR^T(F\bar F^T+\frac{1}{2}I_{2p+1,n})R\bar
F\in\mathbf{M}_p(\mathbb{C})$ we have
$R'R'^T=I_n-I_{2p+1,n}(FL\bar F^TR^T+R\bar
FL^TF^T)I_{2p+1,n}=I_n,\ N+N^T=\bar F^TR^TR\bar F=0,\
R^TF=FL^{-1},\ LMR'=-L\bar F^TR^T(I_n-F\bar F^T-\bar
FF^T)RI_{2p+1,n}=\bar F^TRI_{2p+1,n},\ L(N-\frac{1}{2}MM^T)=L\bar
F^TR^TF\bar F^TR\bar F=\bar F^TR\bar F$, so $R=FL\bar F^T+\bar
F(L^T)^{-1}F^T+FL(N-\frac{1}{2}MM^T)F^T+FLMR'-M^TF^T
+I_{2p+1,n}R'=:R(L,M,N,R')$. Conversely, given
$L\in\mathbf{GL}_p(\mathbb{C}),\ MI_{2p+1,n}=
M\in\mathbf{M}_{p,n}(\mathbb{C}),\
-N^T=N\in\mathbf{M}_p(\mathbb{C}),\
R'\in\mathbf{O}_n(\mathbb{C}),\ R'I_{1,2p}=I_{1,2p}$ we have
$R(L,M,N,R')\in\mathbf{O}_n(\mathbb{C}),\ R(L,M,N,R')F=FL$. If
$R_j:=R(L_j,M_j,N_j,R'_j),\ j=1,2$, then $R_1R_2=R(L,M,N,R'),\
L:=L_1L_2,\ M:=M_2R'^T_1+L_2^{-1}M_1,\
N:=N_2+L_2^{-1}N_1(L_2^T)^{-1}-\frac{1}{2}(L_2^{-1}M_1R'_1M_2^T-(L_2^{-1}M_1R'_1M_2^T)^T),\
R':=R'_1R'_2,\ I_n=R(I_p,0,0,I_n),\
R^{-1}(L,M,N,R')=R(L^{-1},-LMR',-LNL^T,R'^T)$.

We have the decomposition $R(L,M,N,R')=R(L)R(M)R(N)R(R'),\\
R(L):=FL\bar F^T+\bar F(L^T)^{-1}F^T+ I_{2p+1,n},\
R(M):=I_n+FM-(I_n+\frac{1}{2}FM)M^TF^T,\\ R(N):=I_n+FNF^T,\
R(R'):=R'$ and the relations $R(L_1)R(L_2)=R(L_1L_2),\\
R(M)R(L)=R(L)R(L^{-1}M),\ R(N)R(L)=R(L)R(L^{-1}N(L^T)^{-1}),\
R(R')R(L)=R(L)R(R'),\\
R(M_1)R(M_2)=R(M_1+M_2)R(\frac{1}{2}(M_2M_1^T-M_1M_2^T)),\
R(N)R(M)=R(M)R(N),\\ R(R')R(M)=R(MR'^T)R(R'),\
R(N_1)R(N_2)=R(N_1+N_2),\ R(R')R(N)=R(N)R(R'),\\ R(R'_1)R(R'_2)
=R(R'_1R'_2)$.

At the level of the Lie algebra $\mathbf{o}_n(\mathbb{C})$ we have
$r\in\mathbf{o}_n(\mathbb{C}),\ rF=Fl,\
l\in\mathbf{M}_p(\mathbb{C})\Leftrightarrow r=r(l,m,n,r'):=Fl\bar
F^T- \bar Fl^TF^T+Fm-m^TF^T+FnF^T+r',\
mI_{2p+1,n}=m\in\mathbf{M}_{p,n}(\mathbb{C}),\
n\in\mathbf{o}_p(\mathbb{C}),\
I_{2p+1,n}r'=r'\in\mathbf{o}_n(\mathbb{C}),\
[r(l_1,m_1,n_1,r'_1),r(l_2,m_2,n_2,r'_2)]=r([l_1,l_2],l_1m_2-l_2m_1+m_1r'_2-m_2r'_1,l_1n_2
+n_2l_1^T-(l_2n_1+n_1l_2^T)+m_2m_1^T-m_1m_2^T,[r'_1,r'_2]),\
R(e^l)=e^{r(l)},\ R(m)=e^{r(m)},\ R(n)=e^{r(n)},\
R(e^{r'})=e^{r(r')}$. Note $r(m)^3=r(n)^2=0$; this allows simple
integral formulations of $M(t),\ N(t)$: if $R(t):=e^{rt}$, then
$rdt=R^TdR,\ ldt=\bar F^TR^TdRF=L^{-1}dL,\\
r'dt=I_{2p+1,n}R^TdRI_{2p+1,n}=R'^TdR',\ mdt=-\bar
F^TR^TdRI_{2p+1,n}=L^{-1}d(LM)R',\\ n=\bar F^TR^TdR\bar
F=dN+(lN+Nl^T)dt+\frac{1}{2}((L^{-1}d(LM)M^T)^T-L^{-1}d(LM)M^T)$,
so $L(t)=e^{lt},\ R'(t)=e^{r't},\
M(t)=\int_0^te^{-(t-s)l}me^{-sr'}ds,\\ N(t)=\int_0^te^{-sl}
(\frac{1}{2}\int_0^{t-s}(me^{-r'u}m^Te^{-l^Tu}-e^{-lu}me^{r'u}m^T)du+n)e^{-sl^T}ds$.

{\it Example

$R\in\mathbf{O}_n(\mathbb{C})$ which preserve $f_1^Tx=0$ are of
the form $R=R(c)R(v)R(R'),\ R(c):=e^{ic}f_1\bar f_1^T+e^{-ic}\bar
f_1f_1^T+I_{3,n},\ c\in\mathbb{C},\
R(v):=I_n+f_1v^T-(\frac{f_1v^T}{2}+I_n)vf_1^T,\
v=I_{3,n}v\in\mathbb{C}^n,\ R(R'):=R'\in\mathbf{O}_n(\mathbb{C}),\
I_{1,2}R'=I_{1,2}$. Further $R(c_1)R(c_2)=R(c_1+c_2),\
R(v)R(c)=R(e^{-ic}v)R(c),\ R(R')R(c)=R(c)R(R'),\
R(v_1)R(v_2)=R(v_1+v_2),\ R(R')R(v)=R(R'v)R(R')$ and
$\{R(c)R(v)\bar f_1|\ c,\ v\}=\{x\in\mathbb{C}^n|\ |x|^2=0,\
f_1^Tx\neq 0\}$.}

In what follows up to obtaining the SJ form we shall always have
$k,l\in\{1,...,m\},\ j\in\{2,...,m\}$. Consider $m$ orthogonal
isotropic subspaces
$W_k\mathbb{C}^{j_k-j_{k-1}=:p_k}\subset\mathbb{C}^n,\
W_k:=[w_{j_{k-1}+1}...w_{j_k}]$ of rank $p_k,\ j_0:=0,\ j_m:=p,\
W:=[W_1...W_m]$ of rank $p,\ W^TW=0$. By applying rigid motions
$R\in\mathbf{O}_n(\mathbb{C})$ we would like to make
$W_k=F_k:=[f_{j_{k-1}+1}...f_{j_k}]$ and thus $W=F:=[F_1...F_m]$.
As above we can make $W_1=F_1$ and $W=F_1(\bar
F_1^TW)+I_{2p_1+1,n}W$. We have orthogonal isotropic subspaces
$I_{2p_1+1,n}W_j\mathbb{C}^{p_j}\subset I_{2p_1+1,n}
\mathbb{C}^n\simeq\mathbb{C}^{n-2p_1},\ I_{2p_1+1,n}W$ of rank
$p-p_1$; by the induction hypothesis there exists
$R'_1\in\mathbf{O}_n(\mathbb{C}),\ R'_1I_{1,2p_1}=I_{1,2p_1}$
which makes $I_{2p_1+1,n}W=FI_{p_1+1,p}=I_{2p_1+1,n}F$. We can
finish by choosing $M_1:=-\bar F_1^TW\bar F^TI_{2p_1+1,n}:\
R(M_1)W=F$.

Consider the distinct eigenvalues $a_k$ of
$A=A^T\in\mathbf{M}_n(\mathbb{C})$ and their isotropic eigenspaces
$\ker(A-a_kI_n)=W_k\mathbb{C}^{p_k},\ W_k$ as above. Thus we can
make $W=F$ and with $a:=\mathrm{diag}[a_1I_{p_1}...a_mI_{p_m}],\
\mathcal{A}'^T=\mathcal{A}':=\bar F^TA\bar
F\in\mathbf{M}_p(\mathbb{C}),\
\mathcal{H}I_{2p+1,n}=\mathcal{H}:=\bar
F^TAI_{2p+1,n}\in\mathbf{M}_{p,n}(\mathbb{C}),\
\mathcal{A}I_{2p+1,n}=\mathcal{A}^T=\mathcal{A}:=I_{2p+1,n}AI_{2p+1,n}\in\mathbf{M}_n(\mathbb{C})$
we have $A=Fa\bar F^T+\bar
FaF^T+F\mathcal{A}'F^T+F\mathcal{H}+\mathcal{H}^TF^T+
\mathcal{A}$. Since $A-a_kI_n$ has kernel $F_k\mathbb{C}^{p_k}$
and eigenspaces $F_l\mathbb{C}^{p_l}$ for eigenvalues $a_l-a_k,\
l\neq k,\ A-a_kI_n$ restricted to subspaces of $(I_n-F\bar
F^T)\mathbb{C}^n=\bar F\mathbb{C}^p+I_{2p+1,n}\mathbb{C}^n$ and
co-restricted to its image is invertible.

In order to further simplify $A$ we would like to consider
rotations $R$ which preserve each $F_k\mathbb{C}^{p_k}$. As above
these have the form $R(L,M,N,R')$, but we further require that
$RF_k=F_kL_k$ for each $k$ (there are no other restrictions on
$M,N,R'$), so a block decomposition of $L,M,N$ appears:
$L=\mathrm{diag}[L_1...L_m]\in\mathbf{GL}_p(\mathbb{C}),\
M=[M_1^T...M_m^T]^T\in\mathbf{M}_{p,n}(\mathbb{C}),\
-N^T=N=(N_{lk})_{l,k} \in\mathbf{M}_p(\mathbb{C}),\
R'\in\mathbf{O}_n(\mathbb{C}),\ R'I_{1,2p}=I_{1,2p}$.

Under such an $R$ we have $aL=La,\ Fa\bar F^T+\bar FaF^T$ is
preserved and $\mathcal{A},\ \mathcal{H},\ \mathcal{A}'$
respectively become $R(\mathcal{A}):=R'\mathcal{A}R'^T,\
R(\mathcal{H}):=L(\mathcal{H}R'^T+MR(\mathcal{A})-aM),\
R(\mathcal{A}'):=L(\mathcal{A}'+Na-aN-\frac{1}{2}((MR(\mathcal{A})-aM)M^T+
M(MR(\mathcal{A})-aM)^T))L^T+R(\mathcal{H})(LM)^T+(LM)R(\mathcal{H})^T$.
In particular $R(M)(\mathcal{H})=\mathcal{H}+M\mathcal{A}-aM,\
R(N)(\mathcal{A}')=\mathcal{A}'+Na-aN$; if further
$M\mathcal{A}-aM=0$, then $R(L,M)(\mathcal{H})=L\mathcal{H},\
R(L,M)(\mathcal{A}')=L(\mathcal{A}'+\mathcal{H}M^T+M\mathcal{H}^T)L^T$.

By induction hypothesis we can prescribe $R'$ so that
$R(\mathcal{A})$ is SJ (note that $\det(A-\la
I_n)=\prod_{k=1}^m(a_k-\la)^{2p_k}\det(I_{1,2p}+\mathcal{A}-\la
I_{2p+1,n})$ and since all eigenvalues and eigenspaces of $A$ are
accounted for, eigenvalues of the restriction and co-restriction
of $\mathcal{A}$ to $\mathbb{C}^{n-2p}\simeq
I_{2p+1,n}\mathbb{C}^n$ are among the eigenvalues of $A$, while
its eigenspaces are $\subseteq I_{2p+1,n}\mathbb{C}^n$ and thus
disjoint from the eigenspaces of $A$, which are $\subset
I_{1,2p}\mathbb{C}^n$).

Let $\mathcal{H}=:[\mathcal{H}_1^T...\mathcal{H}_m^T]^T,\
\mathcal{H}_kI_{2p+1,n}=\mathcal{H}_k\in\mathbf{M}_{p_k,n}$; if
$\ker(I_{1,2p}+\mathcal{A}-a_1I_{2p+1,n})=F_{m+1}\mathbb{C}^{p_{m+1}}\oplus
E_1\mathbb{C}^{q_1},\ F_{m+1}:=[f_{p+1}...f_{p+p_{m+1}}],\
E_1:=[e_{2(p+p_{m+1})+1}...e_{2(p+p_{m+1})+q_1}]$, then take
$M_1:=-\mathcal{H}_1(I_{1,2p}+\mathcal{A}-a_1I_{2p+1,n})^{\circ
-1}$, which makes $\mathcal{H}_1=\mathcal{H}_1(F_{m+1}\bar
F_{m+1}^T+E_1E_1^T)=:\mathcal{H}'_1[\bar F_{m+1}\ E_1]^T,\\
\mathcal{H}'_1\in\mathbf{M}_{p_1,r_1:=p_{m+1}+q_1}$; we can apply
the same simplification to $\mathcal{H}_2,...,\mathcal{H}_m$. With
$\mathcal{A}'=:(\mathcal{A}'_{kl})_{kl},\\
\mathcal{A}'_{kl}\in\mathbf{M}_{p_k,p_l}(\mathbb{C})$ take
$N_{kl}:=\frac{\mathcal{A}'_{kl}}{a_k-a_l},\ k\neq l$, which makes
$\mathcal{A}'_{kl}=0,\ k\neq l$. In order to preserve the
properties so far obtained we can further arbitrarily prescribe
$L$ and $M$ with the restriction $M\mathcal{A}-aM=0$: for such $M$
we have $M_1=:M'_1[F_{m+1}\ E_1]^T,\ M'_1\in\mathbf{M}_{p_1,r_1}$
and similar formulae for $M_j$; as
$\ker(I_{1,2p}+\mathcal{A}-a_kI_{2p+1,n}),\
\ker(I_{1,2p}+\mathcal{A}-a_lI_{2p+1,n}),\ k\neq l$ are subspaces
of non-isotropic orthogonal subspaces of $\mathbb{C}^n$ we have
$\mathcal{H}_kM^T_l=0,\ k\neq l$. We still have to simplify
$\mathcal{H}'_k,\ \mathcal{A}'_{kk}$; by induction it is enough to
prescribe $L_1,\ M'_1$ so as to simplify $\mathcal{H}'_1,\
\mathcal{A}'_{11}$. We have
$\mathcal{H}_j(F_{m+1}\mathbb{C}^{p_{m+1}}\oplus
E_1\mathbb{C}^{q_1})=0$ (again because eigenspaces for different
eigenvalues of $\mathcal{A}$ are supported in orthogonal
non-isotropic subspaces of $\mathbb{C}^n$); since $A-a_1I_n$
restricted to $F_{m+1}\mathbb{C}^{p_{m+1}}\oplus
E_1\mathbb{C}^{q_1}$ and co-restricted to its image
$F_1\mathcal{H}'_1\mathbb{C}^{r_1}$ is invertible,
$\mathcal{H}'_1$ must have maximal rank $r_1\le p_1$. Take $L_1$
the product of the row operations which bring $\mathcal{H}'_1$ to
its reduced echelon form $\mathcal{H}'_1=[e_1...e_{r_1}]$; we can
further prescribe $L_1$ with the restriction
$L_1\mathcal{H}'_1=\mathcal{H}'_1$ (equivalently $L_1$ is the
product of further arbitrary row operations not involving the
pivot rows). Now choose $[M'_1\ 0_{p_1,p_1-r_1}]:=
-I_{r_1+1,p_1}\mathcal{A}'_{11}I_{1,r_1}-\frac{1}{2}I_{1,r_1}\mathcal{A}'_{11}I_{1,r_1}$,
which makes
$\mathcal{A}'_{11}=I_{r_1+1,p_1}\mathcal{A}'_{11}I_{r_1+1,p_1}$;
since $A-a_1I_n$ restricted to $\bar
F_1I_{r_1+1,p_1}\mathbb{C}^{p_1}$ and co-restricted to its image
$F_1\mathcal{A}_{11}'$ is invertible, $\mathcal{A}_{11}'$ must
have maximal rank $p_1-r_1$. Further choose $L_1$ the product of
row (and thus $L_1^T$ the product of same column) operations which
bring $\mathcal{A}_{11}'$ to the reduced echelon form
$I_{r_1+1,p_1}$. We now permute the coordinates of $\mathbb{C}^n$
in order to follow each eigenvector of $A$ on the flag structure
$A,\ \mathcal{A},...$ along its eigenvalue:

{\it Any matrix $A=A^T\in\mathbf{M}_n(\mathbb{C})$ can be reduced,
via conjugation by a rotation $R\in\mathbf{O}_n(\mathbb{C})$, to a
SJ matrix: a block diagonal matrix formed by SJ blocks $aI_p+J_p,\
J_1:=0,\ J_p:=V_pI'_p\bar V_p^T,\ p\ge 2$, where
$I'_p:=\Sigma_{j=1}^{p-1}e_je_{j+1}^T\in
\mathbf{M}_p(\mathbb{C}),\ p\ge 1$ is the usual upper diagonal in
a Jordan block, $\hat I_p:=\Sigma_{j=1}^pe_je_{p+1-j}^T,\
V_{2p}:=[F\ \bar F\hat I_p],\ V_{2p+1}:=[F\ e_{2p+1}\ \bar F\hat
I_p],\ F:=[f_1...f_p]$. Moreover $\bar V_p^TV_p=I_p,\ V_p\hat
I_p=\bar V_p$, so $J_p^T=J_p,\ J_p^{\circ -1}=\bar J_p$; $J_p$ has
$\bar f_1$ as cyclic vector of order $p$, so $J_p$ has only
eigenvalue $0$ with only eigendirection $f_1$ and minimal
polynomial $x^p$.

Define the type of a SJ matrix $A$ to be the set of all matrices
$M$ with same block diagonal decomposition as that of $A$, each
block $J_p$ being replaced with a polynomial $P(J_p)$ in $J_p$,
$P'(0)\neq 0$. Then two matrices of the same type commute, a
matrix with same block diagonal decomposition as that of a SJ
matrix $A$ and which commutes with $A$ has the type of $A$
(without the restriction $P'(0)\neq 0$), $M^{-1}$ has the type of
$M$ when it exists, $f(I_n-zM)$ has the type of $M$ if $f$ is
analytic near $1$ and
$z\in\mathbb{C}\setminus(\mathrm{spec}(M))^{-1}$ is in the domain
of convergence of $f$ (we are interested mostly in
$\sqrt{I_n-zM}$).}

Although a SJ matrix as defined above efficiently encodes all
metric invariants of a symmetric complex matrix, we prefer when
convenient to work with the SJ type, since it is closed under
certain algebraic operations; the stronger SJ notion would require
a rotation after each algebraic operation.

Given $\sum_{j=1}^{n-1}b_jJ_n^j,\ b_1\neq 0$ we know that there is
$R\in\mathbf{O}_n(\mathbb{C})$ such that
$R(\sum_{j=1}^{n-1}b_jJ_n^j)R^T=J_n$; moreover $RJ_nR^T=J_n,\
R\in\mathbf{O}_n(\mathbb{C})\Rightarrow R=\pm I_n$. With $S:=\bar
V_n^TRV_n$ we have $\hat I_nS^T\hat I_nS=I_n$ and we need
$S(\sum_{j=1}^{n-1}b_jI'^j_n)\hat I_nS^T\hat I_n=I'_n$. If
$S:=b_1^{-\frac{1+n}{2}}\sum_{j=1}^nb_1^je_je_j^T$, then $\hat
I_nS^T\hat I_nS=I_n,\ b_1SI'_n\hat I_nS^T\hat I_n=I'_n$, so we can
assume $b_1=1$. With $\mathcal{K}:=(k_1,...,k_{n-1}),\
k_1,...,k_{n-1}\in\{0,...,n-1\},\
s_1(\mathcal{K}):=\sum_{j=1}^{n-1}k_j,\
s_2(\mathcal{K}):=\sum_{j=1}^{n-1}jk_j,\
a_{kj}:=\sum_{s_1(\mathcal{K})=k,s_2(\mathcal{K})=j}(^{\ \ \ \
k}_{k_1,...,k_{n-1}}) b_2^{k_2}...b_{n-1}^{k_{n-1}}$ we have
$(a_{kj})_{k,j=1,...,n-1}-I_{n-1}$ strictly upper triangular and
$(\sum_{j=1}^{n-1}b_jI'^j_n)^k=\sum_{j=k}^{n-1} a_{kj}I'^j_n,\
k=1,...,n-1$. Thus $\sum_{j=k}^{n-1}a_{kj}SI'^j_n\hat I_nS^T\hat
I_n=I'^k_n,\ k=1,...,n-1$ and $SI'^k_n\hat I_nS^T\hat
I_n=\sum_{j=k}^{n-1}b_{kj}I'^j_n$, where
$(b_{kj})_{k,j=1,...,n-1}-I_{n-1}$ is strictly upper triangular
and $\sum_{j=k}^qa_{kj}b_{jq}=\del_{kq}$. Now
$\\\sum_{j=1}^{n-k}(Se_j)(Se_{n+1-k-j})^T=
\sum_{j=k}^{n-1}b_{kj}\sum_{q=1}^{n-j}e_qe_{n+1-q-j}^T$. By
induction after decreasing $k=n-1,...,1$, $S-I_n$ is strictly
upper triangular and $s_{jk}:=e_j^TSe_k$ can be found by induction
after decreasing $k$ from $s_{j\ n-k}=b_{k\
n-j}-\sum_{l=1}^{n-k-j} s_{1\ l+1}s_{j\ n-k-l},\ k=1,...,n-1,\
j=1,...,n-k$.

We are interested in the rotations which preserve the SJ form of
$A$; since $\ker(A-a_jI_n)^k$ increases with $k$ and stabilizes to
the subspace of $\mathbb{C}^n$ on which the SJ blocks
corresponding to the eigenvalue $a_j$ of $A$ are supported, such a
rotation preserves these subspaces, so it admits a diagonal block
decomposition of rotations of such subspaces. We have thus reduced
the problem to $A$ having only eigenvalue $a_1=0:\
A=\bigoplus_{j=1}^kJ_{p_j},\ p_1\ge...\ge p_k\ge 1,\
\sum_{j=1}^kp_j=n$ and we need $R\in\mathbf{O}_n(\mathbb{C})$ with
$RAR^T=A$. If $k=1,\ p_1=2$, then $n=2k,\ R=R(L)R(N)=FL\bar
F^T+\bar FLF^T+FLNF^T,\ F:=[f_1...f_k],\
L\in\mathbf{O}_k(\mathbb{C}),\ -N^T=N\in\mathbf{M}_k(\mathbb{C})$.
If $k=2,\ p_1=3,\ p_2=2$, then let $A$ be formed of $k_1$ blocks
$J_3$ and $k_2$ blocks $J_2,\ 3k_1+2k_2=n$. Let
$F_1,E_1\in\mathbf{M}_{n,k_1}(\mathbb{C}),\
F_2\in\mathbf{M}_{n,k_2}(\mathbb{C})$ respectively be the matrices
with columns the $f_1$'s of $J_3$'s, $e_3$'s of $J_3$'s, $f_1$'s
of $J_2$'s, considered as vectors in $\mathbb{C}^n$; thus
$A=F_1E_1^T+E_1F_1^T+F_2F_2^T$. Using $\ker(RA^j)=\ker(A^jR)$ and
$AR=RA$ we get $R=F_1L_1\bar F_1^T+(F_2-F_1M_1^T)L_2\bar F_2^T+
(E_1+F_2M_1+F_1(N_1-\frac{1}{2}M_1^TM_1))L_1E_1^T+(\bar
F_2-E_1M_1^T+F_2(N_3-\frac{1}{2}M_1M_1^T)
-F_1(M_1^TN_3+N_1M_1^T+M_2^T))L_2F_2^T+(\bar F_1+\bar
F_2M_1+E_1(N_1-\frac{1}{2}M_1^TM_1)+F_2M_2+F_1(N_2-\frac{1}{2}(M_2^TM_1+M_1^TM_2
+(N_1-\frac{1}{2}M_1^TM_1)^T(N_1-\frac{1}{2}M_1^TM_1))))L_1F_1^T,\
L_j\in\mathbf{O}_{k_j}(\mathbb{C}),\
M_j\in\mathbf{M}_{k_2,k_1}(\mathbb{C}),\
-N_j^T=N_j\in\mathbf{M}_{k_1}(\mathbb{C}),\ j=1,2,\
-N_3^T=N_3\in\mathbf{M}_{k_2}(\mathbb{C})$ arbitrary. While there
is a precise determination of $R$ for $p_1\ge 4$ by similar
computations, it increases in complexity (also the relations among
such rotations). Again the picture at the level of the Lie algebra
should be simpler: $rA-Ar=0,\ r\in\mathbf{o}_n(\mathbb{C})$ and it
becomes clear that for $A=J_n$ such an $r$ must be a polynomial in
$J_n$ (thus $r^T=r$ and $r=0$).

Returning to the quadric $x:\ \ Q(x)=x^T(Ax+2B)+C=0$, once $A$ is
SJ and $\ker(A)$ spanned by $\emptyset,\ e_n$ respectively for
Q(W)C, choose the translation $t:=-o$ and we have respectively
canonical Q(W)C: $B=0,-e_n,\ C=-1,0$ (still denote by $A$ the
matrix $-Q(o)^{-1}A,\ \frac{-1}{e_n^TB}A$, although another
rotation may be required to make the new $A$ SJ). Therefore for
Q(W)C the rigid motions preserving the canonical form are
rotations preserving $A$, so $\mathrm{spec}(A)$ is preserved by
such rotations. For IQWC once $A$ is SJ, $\ker(A)=\mathbb{C}f_1$
and the SJ block for the eigenvalue $0$ of $A$ is $J_p$, choose
the translation $t:=-o$ to get $B=e^{-ic}\bar f,\ C=0$. If we
conjugate $A$ with a rotation $R\in\mathbf{O}_n(\mathbb{C}),\
RI_{p+1,n}=I_{p+1,n},\ R(J_p\oplus
I_{p+1,n})R^T=(e^{2ic_1}J_p)\oplus I_{p+1,n}$ (we use the $\oplus$
to denote the blocks in a block diagonal decomposition; such a
rotation $R$ and $-R+2I_{p+1,n}$ are the only ones with this
property), then with $c_1:=-\frac{c}{p+1}$ the vector $B$ becomes
$-\bar f_1$ and $A$ becomes $e^{-2ic_1}RAR^T$: this is the
canonical form (of course, after a further rotation preserving the
first $p$ coordinates to bring the remaining part of $A$ to the SJ
form). Two canonical forms differ only by a rotation and the
corresponding $A$'s must differ in the first $p$ coordinates by a
rotation as above with $c_1\in\frac{\pi}{p+1}\mathbb{Z}$, in which
case the spectra of these $A$'s differ by a $(p+1)^{\mathrm{th}}$
root of unity. Thus for canonical IQWC $\mathrm{spec}(A)$ is not
preserved by rotations preserving the canonical form, but the
collection of spectra $\{e^{\frac{2\pi
ij}{p+1}}\mathrm{spec}(A)\}_{j=0,...,p}$ is; for simplicity we
shall make a choice of $\mathrm{spec}(A)$ and work with it.

We have the diagonal Q(W)C respectively for
$A=\Sigma_{j=1}^na_j^{-1}e_je_j^T,\
A=\Sigma_{j=1}^{n-1}a_j^{-1}e_je_j^T$; the diagonal IQWC come in
different flavors, according to the block of $f_1:\
A=J_p+\Sigma_{j=p+1}^na_j^{-1}e_je_j^T$; in particular if $A=J_n$,
then $\mathrm{spec}(A)=\{0\}$ is unambiguous. Thus general
quadrics are those for which all eigenvalues have geometric
multiplicity $1$; equivalently each eigenvalue has an only
corresponding SJ block.

With $R_z:=I_n-zA$ the family of quadrics $\{x_z\}_z$ confocal to
$x_0$ is given by
$Q_z(x_z)=x_z^TAR_z^{-1}x_z+2(R_z^{-1}B)^Tx_z+C+zB^TR_z^{-1}B=0$.

Consider the IQWC $x^T(Ax+2B)+C=0,\ \ker(A)=\mathbb{C}v,\
|v|^2=0,\ v\in\mathbb{C}^n\setminus\{0\}$ and $o\in\mathcal{O}$.
Applying the translation $-o$ and a rotation $R,\ Rv=f_1,\
R(Ao+B)=-e^{-ic}\bar f_1,\ e^{-ic}:=-v^TB$ we can make
$\ker(A)=\mathbb{C}f_1,\ B=-e^{-ic}\bar f_1,\ C=0$. If we begin
with a different point in $\mathcal{O}$, we get a similar
situation; the difference between the two is a rigid motion
$(R,t)$ which preserves the form of $A,\ B,\ C:\
\ker(RAR^T)=\mathbb{C}f_1,\ e^{-ic}R\bar
f_1+RAR^Tt=e^{-i(c+c_1)}\bar f_1,\ Q(-R^Tt)=0$, so
$Rf_1=e^{ic_1}f_1$ and $t=e^{-i(c+c_1)}(I_n-f_1\bar
f_1^T)RA^{\circ -1}R^T\bar f_1$ is obtained by plugging into the
formula $o=-A^{\circ -1}B+\frac{B^TA^{\circ -1}B-C}{2v^TB}v$ the
values $RA^{\circ -1}R^T,\ -e^{-i(c+c_1)}\bar f_1,\ 0,\ f_1$
respectively instead of $A^{\circ -1},\ B,\ C,\ v$. In this
simplified setting it is easy to see that the formula of $o$ is
well defined up to rotations of the type $R(c_1)R(R')$ (which also
fix $\mathbb{C}\bar f_1$); the remaining type of rotations $R(v)$
give a parametrization of $\mathcal{O}:\
\mathcal{O}=\{e^{ic}(I_n-f_1\bar f_1^TR(v)^T) (R(v)AR(v)^T)^{\circ
-1}R(v)\bar f_1|\ v=I_{3,n}v\in\mathbb{C}^n\}$.

There is a unifying point of view in the projective space: the
action of the group of rigid motions and homotheties
$\mathbb{C}^n\ni x\mapsto (e^a(R,t))x:=e^a(Rx+t),\
(R,t)\in\mathbf{O}_n(\mathbb{C})\ltimes\mathbb{C}^n,\ a\in
\mathbb{C}$ can be extended to
$\mathbb{CP}^n\supseteq[\mathbb{C}^{nT},1]\simeq \mathbb{C}^n$ as
$(e^a(R,t))[x^T,x^{n+1}]:=[(\begin{bmatrix}R&t\\0^T&e^{-a}\end{bmatrix}
\begin{bmatrix}x\\x^{n+1}\end{bmatrix})^T]=[(Rx+x^{n+1}t)^T,e^{-a}x^{n+1}]$ and is
transitive on Cayley's absolute $C(\infty):=\{[Y^T,0]|\
Y\in\mathbb{C}^n\setminus\{0\},\ |Y|^2=0\}$ (circle at infinity
for $n=3$); in fact it is the maximal subgroup of
$\mathbf{PGL}_{n}(\mathbb{C})$ which preserves $C(\infty)$. Any
quadric in $\mathbb{CP}^{n-1}=[(\mathbb{C}^n\setminus\{0\})^T,0]$
can be brought, via a homography, to $C(\infty)$, thus confirming
that any symmetric non-degenerate bilinear form on $\mathbb{C}^n$
is isometric to the Euclidean scalar product.

We have $\begin{bmatrix}A&B\\B^T&C\end{bmatrix}^{-1}
=\begin{bmatrix}A^{\circ -1}&b\\b^T&C\end{bmatrix}$, where $A$ is
SJ and respectively for canonical (I)Q(W)C $\\\ker(A)=0,\
\mathbb{C}e_n,\ \mathbb{C}f_1,\ B=0,-e_n,-\bar f_1,\
b=0,-e_n,-f_1,\ C=-1,0,0$.
If $z\in\mathbb{C},\ \det R_z\neq 0$, then $\begin{bmatrix}AR_z^{-1}&R_z^{-1}B\\
B^TR_z^{-1}&C+zB^TR_z^{-1}B\end{bmatrix}=\begin{bmatrix}A^{\circ -1}-zI_n&b\\
b^T&C\end{bmatrix}^{-1}$.

A quadric $[x^T,x^{n+1}]\subset\mathbb{CP}^n$ satisfies
$\begin{bmatrix}x\\x^{n+1}\end{bmatrix}^T
\begin{bmatrix}A&B\\B^T&C\end{bmatrix}\begin{bmatrix}x\\x^{n+1}\end{bmatrix}=0,\ A^T=A,\
\begin{vmatrix}A&B\\B^T&C\end{vmatrix}\neq 0$ and its family of confocal quadrics
$\begin{bmatrix}x_z\\x_z^{n+1}\end{bmatrix}^T(\begin{bmatrix}A&B\\B^T&C\end{bmatrix}^{-1}-z
\begin{bmatrix}I_n&0\\0^T&0\end{bmatrix})^{-1}\begin{bmatrix}x_z\\x_z^{n+1}\end{bmatrix}=0$.

If we define the adjugate of the quadric
$\begin{bmatrix}x\\x^{n+1}\end{bmatrix}^T
\begin{bmatrix}A&B\\B^T&C\end{bmatrix}\begin{bmatrix}x\\x^{n+1}\end{bmatrix}=0$ to be the quadric
${\begin{bmatrix}x\\x^{n+1}\end{bmatrix}
^*}^T\begin{bmatrix}A&B\\B^T&C\end{bmatrix}^{-1}
\begin{bmatrix}x\\x^{n+1}\end{bmatrix}^*=0$, then
$(\begin{bmatrix}R&t\\0^T&e^{-a}\end{bmatrix}\begin{bmatrix}x\\x^{n+1}\end{bmatrix})^*=
(\begin{bmatrix}R&t\\0^T&e^{-a}\end{bmatrix}^T)^{-1}\begin{bmatrix}x\\x^{n+1}\end{bmatrix}^*$;
the family of quadrics confocal to $[x^T,x^{n+1}]$ is the adjugate
of the pencil generated by the adjugate of $[x^T,x^{n+1}]$ and
$C(\infty)$.

Conversely, given a quadric $Q\subset\mathbb{CP}^n$, a hyperplane
$P\subset\mathbb{CP}^n$ and a quadric $C(\infty)\subset P$, we can
consider the Euclidean metric introduced by $C(\infty)$ on
$\mathbb{CP}^n\setminus P;\ P$ becomes the hyperplane
$\mathbb{CP}^{n-1}$ at $\infty$ and we can recover the
construction above. If $P$ is not tangent to $Q$, then $Q$ is QC;
if $P$ is tangent to $Q$ (that is $P\cap Q$ is a quadratic cone),
then $Q$ is (I)QWC; if $C(\infty)$ does not pass through the point
of tangency of $P$ and $Q$ (that is the vertex of the cone $P\cap
Q)$, then $Q$ is QWC, otherwise $Q$ is IQWC.

If $\ker(A)=0$ or is not isotropic, $A$ SJ, then we can define
$\sqrt{A}$ to be a symmetric matrix of the same type as $A:\
\sqrt{e^aI_p+J_p}=\sqrt{e^a}\Sigma_{j=0}^{p-1}e^{-ja}
(^{\frac{1}{2}}_j)J_p^j,\ p\ge 1$; (the sign of $\sqrt{.}$ does
not matter, because any quadric is symmetric wrt the principal
spaces, but in order to keep a strict account of these symmetries
(to be used for other reasons), we shall always choose
$\sqrt{re^{i\theta}}:=\sqrt{r}e^{i\frac{\theta}{2}},\ r\ge 0,\
-\pi<\theta\le\pi$ and we shall avoid simplifications of the form
$\sqrt{a}\sqrt{b}=\pm\sqrt{ab}\simeq\sqrt{ab}$). Consider the unit
sphere $X\subset \mathbb{C}^n,\ |X|^2=1$ and the equilateral
paraboloid $Z\subset\mathbb{C}^n ,\ Z^T(I_{1,n-1}Z-2e_n)=0$; then
we have parametrization $x_0=(\sqrt{A})^{-1}X$ of QC and
$x_0=(\sqrt{A+e_ne_n^T})^{-1}Z$ of QWC. A natural parametrization
of IQWC will appear later; because of this we should rather
require that $A^{\circ -1}$ is SJ.

For $z$ inverses of nonzero eigenvalues of $A$ (which make $R_z$
singular) we have finite singular quadrics; for $z=\infty$ we have
the singular quadric at $\infty$; all are obtained as sets by
letting $z$ tend to the singular value in the closure of
$x_z\simeq[x_z^T,1]\subset\mathbb{CP}^n$. If $A$ is of SJ type
with $\ker(R_a)$ spanned by $f_1,...,f_k,e_{2k+1},...,e_p$, then
by L'Hospital $x_a=\{x\in \ker(R_a)||x|^2=0\}\times
I_{p+1,n}\mathbb{C}^n,\ x_{\infty}=\mathbb{CP}^{n-1}=[X^T,0]\cup
C(\infty)$; for example if $A=aI_n$, then $x_a$ is the isotropic
cone $|x|^2=0$. Singular quadrics have extra projective structure:
given a quadric $Q\subset \mathbb{C}^n$ and $p\in \mathbb{C}^n$
the (quadratic) tangent cone of $Q$ at $p$ is
$C_pQ:=\{v\in\mathbb{C}^n|Q(p+sv)=0$ has a double root $s\}$ and
degenerates to the singular quadratic cone $T_pQ$ if $p\in Q$.
Thus $\mathbb{C}^n$ can be seen as a quadratic cone bundle over
itself, factorized by the relation of equivalence $p\sim q$ iff
$p=q$ as points in $\mathbb{C}^n$. The factor set has singular set
$Q$ (which can again be seen as quadratic cone bundle over itself
factorized by the same relation of equivalence, but this time the
factor set has no singular points, since all cones $Q\cap T_pQ,\
p\in Q$ are projectively equivalent). Similarly singular confocal
quadrics can be seen as quadratic cone bundles over themselves,
factorized by equivalence relationships whose singular sets are
lower dimensional quadrics and this (singular) projective
structure of singular quadrics of the confocal family is the limit
of the (non-singular) projective structure of non-singular
quadrics of the confocal family; the whole confocal family can be
recovered from the knowledge of the singular set of certain finite
singular quadrics of the family (including those for which the
eigenvalue has geometric multiplicity $1$). The singular set of
$x_{\infty}$ is $\mathcal{S}(x_{\infty}):=C(\infty)$; for $z=a_1,\
A$ of the SJ type with $\ker(R_{a_1})\neq 0$ the singular set
$x:=\mathcal{S}(x_{a_1})$ of the singular quadric $x_{a_1}$
satisfies: $x^T(AR_{a_1}^{\circ -1}x+2R_{a_1}^{\circ
-1}B)+C+a_1B^TR_{a_1}^{\circ -1}B=0,\ x=R_{a_1}R_{a_1}^{\circ
-1}x$. For a family of confocal quadrics, the spectrum is the set
of singular values of $z$ (thus it is not well defined, since it
depends on the choice of the initial quadric); the relative
spectrum is the set of differences of values in the spectrum (thus
it is well defined).

\subsection{Proofs of Lam\'{e}, Dupin, Ivory, Bianchi I, Jacobi
and Chasles theorems; a ray of light}
\noindent

\noindent Assume that the quadric $x_0$ is in the canonical form:
$A$ SJ, $\ker(A)=0,\ \mathbb{C}e_n,\ \mathbb{C}f_1,\ B=0,-e_n,\
-\bar f_1,\ C=-1,0,0$ respectively for (I)Q(W)C. The Ivory
affinity for QC should preserve the center $0$ of the quadrics
$x_z$, so it is given by a linear transformation $x_z=M(z)x_0$,
where $M(z)\in\mathbf{GL}_n(\mathbb{C})$ satisfies
$M(z)^TAR_z^{-1}M(z)=A$; the simplest choice $M(z):=\sqrt{R_z}$ is
the correct one.

For (I)QWC we make the ansatz $x_z=\sqrt{R_z}x_0+C(z)$ and further
require it to satisfy Lam\'{e}'s:
$0=-2dx_w^T\pa_z|_{z=w}x_z=dx_0^T\sqrt{R_w}(A(\sqrt{R_w})^{-1}x_0-2C'(w))$.
Using $dx_0^T(Ax_0+B)=0$ we get
$Ax_0-2\sqrt{R_w}C'(w)=c(w,x_0)(Ax_0+B)$; multiplying on the left
with $v^T,\ 0\neq v\in\ker(A)$ we get
$c(w,x_0)=c(w)=-\frac{2v^T\sqrt{R_w}C'(w)}{v^TB}$; making $x_0=0$
we get $2\sqrt{R_w}C'(w)=-c(w)B$, so $(1-c(w))Ax_0=0$ and since
the quadric $x_0$ does not lie in $\ker(A)$, we get $c(w)=1$ and
$C(z)=(-\frac{1}{2}\int_0^z(\sqrt{R_w})^{-1}dw)B (=\frac{z}{2}e_n$
for QWC; for IQWC it is the Taylor series of
$\frac{1}{2}\int_0^z(\sqrt{1-w})^{-1}dw$ at $z=0$ with each
monomial $z^{k+1}$ replaced by $-z^{k+1}J_p^kB$, where $J_p$ is
the block of $f_1$ and thus a polynomial of degree $p$ in $z$).
Note $AC(z)+(I_n-\sqrt{R_z})B=0=(I_n+\sqrt{R_z})C(z)+zB$. Applying
$d$ to $Q_z(x_z)=0$ we get $dx_z^TR_z^{-1}(Ax_z+B)=0$, so the unit
normal $N_z$ is proportional to $\hat N_z:=-2\pa_zx_z$. If
$\mathbb{C}^n\ni x\in x_{z_1},x_{z_2}$, then $\hat
N_{z_j}=R_{z_j}^{-1}(Ax+B)$; using $R_z^{-1}-I_n=zAR_z^{-1},\
z_1R_{z_1}^{-1}-z_2R_{z_2}^{-1}=(z_1-z_2)R_{z_1}^{-1}R_{z_2}^{-1}$
we get $0=Q_{z_1}(x)-Q_{z_2}(x)=(z_1-z_2)\hat N_{z_1}^T\hat
N_{z_2}$. The polynomial equation  $Q_z(x)=0$ has degree $n$ in
$z$ and it has multiple roots iff $0=\pa_zQ_z(x)=|\hat N_z|^2$.

Ivory  becomes
$|V_0^1|^2=|x_0^0+x_0^1-C(z)|^2-2(x_0^0)^T(I_n+\sqrt{R_z})x_0^1+zC=|V_1^0|^2$;
Bianchi I becomes for lengths of rulings: if $w_0^TAw_0=w_0^T\hat
N_0=0,\ w_z=\sqrt{R_z}w_0$, then
$w_z^Tw_z=|w_0|^2-zw_0^TAw_0=|w_0|^2$; for the symmetry of the TC:
$(V_0^1)^T\hat N_0^0
=(x_0^0)^TA\sqrt{R_z}x_0^1-B^T(x_z^0+x_z^1-C(z))+C=(V_1^0)^T\hat
N_0^1$; for angles between segments and rulings:
$(V_0^1)^Tw_0^0+(V_1^0)^Tw_z^0=-z(\hat N_0^0)^Tw_0^0=0$; for
angles between rulings:
$(w_0^0)^Tw_z^1=(w_0^0)^T\sqrt{R_z}w_0^1=(w_z^0)^Tw_0^1$; for
angles between polar rulings: $(w_z^0)^T\hat w_z^0=(w_0^0)^T\hat
w_0^0-z(w_0^0)^TA\hat w_0^0 =(w_0^0)^T\hat w_0^0$.

For Jacobi's result a geodesic in affine parametrization must
satisfy $Q(x_0)=\dot x_0^T\hat N_0=\ddot x_0^T dx_0=0$, so
$0=(\dot x_0^T\hat N_0)^{\cdot}=\ddot x_0^T\hat N_0+\dot
x_0^T\dot{\hat N}_0$, or $\ddot x_0=-\frac{\dot x_0^T\dot{\hat
N}_0}{|\hat N_0|^2}\hat N_0$. The tangent line at $x_0(t)$ is
tangent to a confocal quadric $x_{z(t)}$ iff the quadratic
equation in $s$: $\frac{1}{2}Q_z(x_0+s\dot x_0)=0$ has a double
solution, so we need its discriminant $T_0:=(\dot
x_0^TR_z^{-1}\hat N_0)^2-\dot x_0^TAR_z^{-1}\dot x_0Q_z(x_0)$ to
be 0. With $a,\ b,\ c,\ d\in\mathbb{C}^n$ and $[a,b]:=ba^T-ab^T,\
[[a,b],c]:=(ba^T-ab^T)c$ we have
$[a,b]\bullet[c,d]:=\frac{1}{2}\mathrm{tr}([a,b]^T[c,d])=a^T[b,[c,d]]=a^Tcb^Td-a^Tdb^Tc$.
Using also $Q_z(x_0)=z\hat N_0^TR_z^{-1}\hat N_0$ we get $[\dot
x_0,\hat N_0]\bullet[R_z^{-1}\hat N_0,R_z^{-1}\dot x_0] +|\dot
x_0|^2\hat N_0^TR_z^{-1}\hat N_0=0$. After eliminating the
denominator, $T_0$ becomes a polynomial of degree $n-1$ in $z$
with the coefficient of the {\it highest order term} (hot) being
$|\dot x_0|^2\neq 0$, $0$ as one root and in general having
distinct roots; thus Jacobi's result is proven if $\dot z=0$.
Dotting $T_0=0$ we get $\dot T_0=T_1+\dot zT_2=0$ where $T_1:=
[\dot x_0,\dot{\hat N}_0]\bullet[R_z^{-1}\hat N_0,R_z^{-1}\dot
x_0]+ [\dot x_0,\hat N_0]\bullet[R_z^{-1}\dot{\hat
N}_0,R_z^{-1}\dot x_0] +2|\dot x_0|^2\hat N_0^TR_z^{-1}\dot{\hat
N}_0=2(\hat N_0^TR_z^{-1}\dot x_0\dot x_0^TAR_z^{-1}\dot x_0 -\dot
x_0^T(R_z^{-1}-I_n)\dot x_0\dot x_0^TAR_z^{-1}\hat N_0)=0$ and
$T_2:=\pa_zT_0=-[\dot x_0,\hat N_0]\bullet([AR_z^{-2}\hat
N_0,R_z^{-1}\dot x_0] +[R_z^{-1}\hat N_0,AR_z^{-2}\dot x_0])-|\dot
x_0|^2\hat N_0^TAR_z^{-2}\hat N_0\neq 0$ in general; thus $\dot
z=0$. We have $T_2=0$ precisely when $z$ is a multiple root for
the considered polynomial, but in this case if $\dot z\neq 0$,
then similarly we have $\dot T_2=T_3+\dot zT_4,\ T_3:=-[\dot
x_0,\dot{\hat N}_0]\bullet ([AR_z^{-2}\hat N_0,R_z^{-1}\dot
x_0]+[R_z^{-1}\hat N_0,AR_z^{-2}\dot x_0]) -[\dot x_0,\hat
N_0]\bullet([AR_z^{-2}\dot{\hat N}_0,R_z^{-1}\dot x_0]
+[R_z^{-1}\dot{\hat N}_0,AR_z^{-2}\dot x_0])-2|\dot x_0|^2\hat
N_0^TAR_z^{-2}\dot{\hat N}_0=\\=-2([\dot x_0,\hat
N_0]\bullet([A^2R_z^{-2}\dot x_0,R_z^{-1}\dot x_0] +[AR_z^{-1}\dot
x_0,AR_z^{-2}\dot x_0])+|\dot x_0|^2\hat N_0^TA^2R_z^{-2}\dot
x_0)=0$ (use $zAR_z^{-j}=R_z^{-j}-R_z^{-j+1},\ j=1,2$), so
$T_4:=\pa_zT_2=0$ and so on, which implies the fact that $z$ has
multiplicity $n-1$, so the geodesic is a ruling on $x_0$,
contradicting $\dot z\neq 0$. Note that we have proven that the
roots of $T_0,T_2,T_4,...$ or of any linear combination thereof
satisfy the same property: they remain constant along the
geodesic.

For diagonal QC ($A=\Sigma_{j=1}^na_j^{-1}e_je_j^T$) we have
$T_0=\Sigma_{j<k}\frac{(x_0^j\dot x_0^k-x_0^k\dot x_0^j)^2}
{(a_j-z)(a_k-z)}-\Sigma_j\frac{(\dot x_0^j)^2}{a_j-z}$; for
diagonal QWC ($A=\Sigma_{j=1}^{n-1}a_j^{-1}e_je_j^T$) we have
$T_0=\Sigma_{j<k}\frac{(x_0^j\dot x_0^k-x_0^k\dot
x_0^j)^2}{(a_j-z)(a_k-z)} +\Sigma_j\frac{2x_0^j\dot
x_0^jx_0^n-(2x_0^n-a_j)(\dot x_0^j)^2}{a_j-z}-|\dot x_0|^2$.

The part of Chasles's result involving geodesics can be easily
read from the proof of Jacobi's result. Let the $n-1$ confocal
quadrics be $x_{z_0:=0},\ x_{z_1},...,x_{z_{n-2}}$ and
$z\in\{z_0,z_1,...,z_{n-2}\}$. Consider the line $x_0+sv$ tangent
to the $n-1$ quadrics at points $x_z:=x_0+s_zv$. We thus have
$(v^TR_z^{-1}\hat N_0)^2 -zv^TAR_z^{-1}v\hat N_0^TR_z^{-1}\hat
N_0=0,\ s_z=-\frac{v^TR_z^{-1}\hat N_0}{v^TAR_z^{-1}v},\
x_z=x_0+s_zv
=\frac{[[v,x_0],AR_z^{-1}v]-v^TR_z^{-1}Bv}{v^TAR_z^{-1}v}$ and
$zv^TAR_z^{-1}v\hat N_0^T\hat N_z =z\hat
N_0^TR_z^{-1}(A[[v,x_0],AR_z^{-1}v]+[[Av,B],R_z^{-1}v])= v^T\hat
N_0\hat N_0^TR_z^{-1}v-\\-(v^TR_z^{-1}\hat
N_0)^2+zv^TAR_z^{-1}v\hat N_0^TR_z^{-1}\hat N_0=0$.

Therefore $N_0^TN_z=0$ along the common tangents and by symmetry
this is true for any pair of normals of the $n-1$ confocal
quadrics (the choice of the initial quadric $x_0$ is arbitrary).
The remaining part is due to Dupin and Malus (see Darboux
(\cite{D1},\S\ 441)): a congruence of lines in $\mathbb{C}^n$
which admits $n-1$ focal sub-manifolds whose normals are
orthonormal is normal. We would like to find a sub-manifold
$y:=x_0+sv$ orthogonal to unit $v:\ 0=v^Tdy=v^Tdx_0+ds$. Imposing
the compatibility condition $d\wedge$ we thus need $0=dx_0^T\wedge
dv$. From $dx_0=(\Sigma_zN_zN_z^T+vv^T)dx_0$ we thus need to prove
that $dx_0^T\Sigma_{z\neq z_0}N_zN_z^T\wedge dv=0$ and
$0=dx_z^TN_z=(dx_0+s_zdv)^TN_z$ finishes the proof.

Note that $n-1$ confocal quadrics $Q_{z_1},...,Q_{z_{n-1}}$ are
seen orthogonally from any point $p$ on any lines tangent to all
$n-1$ quadrics, that is the cones through $p$ tangent to the given
quadrics intersect along the principal directions at that point;
moreover all common tangent to $Q_{z_1},...,Q_{z_{n-1}}$ and
passing through $p$ are generated by reflection in the principal
hyperplanes of the initial common tangent through $p$ (Darboux
(\cite{D1},\S\ 465)).

\subsection{Proof of Bianchi II}
\noindent

\noindent For Bianchi II consider
$\mathbb{C}^n\simeq[\mathbb{C}^{nT},1]\subset\mathbb{CP}^n$ and a
homography $x\mapsto H(x)$ is given by $H(x):=\frac{H'x+h_1}
{h_2^Tx+a},\
H:=\begin{bmatrix}H'&h_1\\h_2^T&a\end{bmatrix}\in\mathbf{GL}_{n+1}(\mathbb{C}),\
H'\in\mathbf{M}_n(\mathbb{C}),\ h_1,h_2\in\mathbb{C}^n,\
a\in\mathbb{C}$ being determined up to multiplication by a
non-zero constant. Suppose $H$ takes a family of confocal quadrics
$x_z$ given by $\begin{bmatrix}A&B\\B^T&C\end{bmatrix}$ to another
one $\ti x_{\ti z}:=H(x_z),\ \ti z(0)=0$ given by
$\begin{bmatrix}\ti A&\ti B\\\ti B^T&\ti C\end{bmatrix}$. We
exclude rigid motions and homotheties (for which
$H(C(\infty))=C(\infty)$) as they do not metrically change the
family $x_z$; as we shall see later $H$ cannot be an affine
transformation of $\mathbb{C}^n$ (that is $h_2\neq 0$) and further
$\det(H\begin{bmatrix}I_n&0\\0^T&0\end{bmatrix}H^T-\la\begin{bmatrix}I_n&0\\0^T&0\end{bmatrix})=
\begin{vmatrix}H'H'^T-\la I_n&H'h_2\\(H'h_2)^T&|h_2|^2\end{vmatrix}\neq 0$ for some
$\la\in\mathbb{C}$ (and thus for all but $\le n\
\la\in\mathbb{C}$). Let $\ti H:=H^{-1}$; in the following we
consider tilde quantities as counterparts of non-tilde quantities
(and thus are defined to have same properties); all equalities
remain valid if tilde and non-tilde quantities are exchanged.

By applying rigid motions in both the domain and range of $H$ we
make $\ker(A)=0,\ \mathbb{C}e_1,\ \mathbb{C}f_1,\ B=0,\ -e_1,\
-\bar f_1,\ b=0,-e_1,\ -f_1,\ C=-1,0,0$ respectively for (I)Q(W)C.
As $\begin{bmatrix}A&B\\B^T&C\end{bmatrix}^{-1}
=\begin{bmatrix}A^{\circ -1}&b\\b^T&C\end{bmatrix}$ we have $H\begin{bmatrix}A^{\circ -1}-zI_n&b\\
b^T&C\end{bmatrix}H^T=\la(z)^{-1}\begin{bmatrix}\ti A^{\circ -1}-\ti z(z)I_n&\ti b\\\ti b^T&\ti C
\end{bmatrix}$; thus the singular values of $z,\ \ti z$ correspond; in particular since
$H^{-1}(\mathbb{CP}^{n-1})\neq\mathbb{CP}^{n-1},\
H^{-1}(C(\infty)),\ H(C(\infty))$ are respectively singular sets
of finite singular quadrics of the families $x_z,\ \ti x_{\ti z}$
whose corresponding eigenvalues $a_1,\ \ti a_1$ thus must have
geometric multiplicity $1$. We need $H\begin{bmatrix}A^{\circ
-1}&b\\b^T&C\end{bmatrix}H^T =\la(0)^{-1}\begin{bmatrix}\ti
A^{\circ -1}&\ti b\\\ti b^T&\ti C\end{bmatrix},\
H\begin{bmatrix}I_n&0\\0^T&0\end{bmatrix}H^T=\frac{\la(z)^{-1}-\la(0)^{-1}}{-z}
\begin{bmatrix}\ti A^{\circ -1}&\ti b\\\ti b^T&\ti C\end{bmatrix}
+\frac{\la(z)^{-1}\ti
z(z)}{z}\begin{bmatrix}I_n&0\\0^T&0\end{bmatrix}$; thus
$\\\frac{\la(z)^{-1}-\la(0)^{-1}}{-z},\ \frac{\la(z)^{-1}\ti
z(z)}{z}$ are constants. Under the homothety $e^t\ A,\ z$ become
$e^{-2t}A,\ e^{2t}z$ if $x_z$ are QC and $e^{-t}A,\ e^tz$ if $x_z$
are (I)QWC (when $\begin{bmatrix}A^{\circ
-1}-zI_n&b\\b^T&0\end{bmatrix}$ becomes
$e^t\begin{bmatrix}e^t(A^{\circ -1}-zI_n)&b\\b^T&0\end{bmatrix}$).
If $\la(z)$ is constant, then $H(C(\infty))=C(\infty)$; otherwise
$\la(z)^{-1}=e^s(z-a_1),\ \ti z(z)=a_1\ti
a_1(a_1^{-1}+(z-a_1)^{-1})$ and by applying homotheties in either
the range, domain of $H$ or both we make $\ti a_1=a_1^{-1}$;
further fixing the non-zero constant upon which $H$ depends we
make $s=0$ and $\begin{bmatrix}\ti A^{\circ -1}-\ti a_1I_n&\ti
b\\\ti b^T&\ti C\end{bmatrix}
=-H\begin{bmatrix}I_n&0\\0^T&0\end{bmatrix}H^T,\
H\begin{bmatrix}A^{\circ -1}-a_1I_n&b\\b^T&C\end{bmatrix}H^T
=-\begin{bmatrix}I_n&0\\0^T&0\end{bmatrix}$. To preserve these we
exclude further homotheties in either the range or domain of $H$.

Applying rigid motions in both the domain and range of $H$, we can
suppose that $A$ is of the SJ type, $\ker(A)=\mathbb{C}V_pe_1,\
B=-\bar V_pe_1,\ b=-V_pe_1,\ C=-\sqrt{1-B^Tb}$ and further the
block of $V_pe_1$ in $A$ is $-J_p$ (the definitions of $V_p,\ J_p$
as $V_{2k}:=[f_1...f_k\ \bar f_k...\bar f_1],\ k\ge 1,\
V_{2k+1}:=[f_1...f_k\ e_{2k+1}\ \bar f_k...\bar f_1],\ k\ge 0,\
J_p:=V_pI'_{1,p-1}\bar V_p^T,\ p\ge 1$ can be extended to
$V_0:=\emptyset,\ J_0:=\emptyset$ with the convention that
matrices whose either dimension is $0$ are empty and are excluded
from the block decomposition, $V_0e_1=0=0_{n,1}$ as the span of
$\emptyset$). Further we require that the block of $a_1,\
\mathcal{S}(x_{a_1})= H^{-1}(C(\infty))$ in $A^{\circ -1}$ is
$a_1I_q-J_q,\ q\ge 1$, so $A^{\circ -1}=
\begin{bmatrix}-\bar J_p&0&0\\0&a_1I_q-J_q&0\\0&0&a_1I_r-A'^{-1}\end{bmatrix}
=:(-\bar J_p)\oplus(a_1I_q-J_q)\oplus(a_1I_r-A'^{-1}),\ p,r\ge 0,\
q\ge 1,\ p+q+r=n$ and
$A'^{-1}:=\bigoplus_{j=2}^m((a_1-a_j)I_{s_j}+\bigoplus_kJ_{s_{jk}})
\in\mathbf{GL}_r(\mathbb{C}),\ \Sigma_ks_{jk}=s_j>0,\ j=2,...,m,\
\Sigma_js_j=r$. To simplify the notation we drop the $\oplus$ for
zero terms: for example $AJ_q=(a_1I_q-J_q)^{-1}J_q$.

Quadrics whose eigenspaces of all non-zero eigenvalues of $A$ have
dimension $\ge 2$ do not enjoy the property of their confocal
family being taken to another one; for example QC with
$A=a^{-1}I_n$, whose only finite singular quadric (namely the
isotropic cone $|x|^2=0$) is not projectively equivalent to
$\mathbb{CP}^{n-1}$) or IQWC with $A=J_n$ (in which case there are
no finite singular quadrics).

Conversely, given $x_z$ as above, then since
$\mathcal{S}(x_{a_1})$ is projectively equivalent to $C(\infty)$,
as a quadric in a hyperplane of $\mathbb{CP}^n$, there exists $H$
homography (unique up to rigid motions in the domain of $H$ which
fix $\mathcal{S}(x_{a_1})$ and rigid motions and homotheties in
the range of $H$) such that $H(\mathcal{S}(x_{a_1}))=C(\infty)$.
$H$ takes $x_z$ to another family of confocal quadrics $\ti x_{\ti
z}$ (since $\ti x_{\ti z}$ is uniquely determined by $C(\infty)$
and $\mathcal{S}(\ti x_{\ti a_1}):=H(C(\infty))$).

Since $\mathcal{S}(x_{a_1})=H^{-1}(C(\infty)):\
|H'x+h_1|^2=h_2^Tx+a=0$ we have $a=0,\ h_2=e^{c_1}V_qe_1$. If
$q=1$, then $V_qe_1=e_{p+1},\ x_z\cap\{e_{p+1}^Tx=0\}$ is a
quadric in $e_{p+1}^Tx=0$ and is taken by $H$ into $\ti x_{\ti
z}\cap\mathbb{CP}^{n-1}$, which thus must be a quadric and
therefore $\ti x_{\ti z}$ are QC; if $q\ge 2$, then
$x_z\cap\{(V_qe_1)^Tx=0\}$ is a degenerate quadric in
$(V_qe_1)^Tx=0$ (a quadratic cone with vertex $[(V_qe_1)^T,0]\in
C(\infty)$, since the quadric $x_z$ is tangent to $(V_qe_1)^Tx=0$
at $[(V_qe_1)^T,0]$) and is taken by $H$ into $\ti x_{\ti
z}\cap\mathbb{CP}^{n-1}$, which thus must be a degenerate quadric
and therefore $\ti x_{\ti z}$ are (I)QWC, according to
$[(V_qe_1)^T,0]\notin\mathcal{S}(x_{a_1})
(\in\mathcal{S}(x_{a_1}))\Leftrightarrow (V_qe_1)^TR_{a_1}^{\circ
-1}(V_qe_1)\neq 0\ (=0) \Leftrightarrow(V_qe_1)^T\bar
J_q(V_qe_1)\neq 0\ (=0)\Leftrightarrow q=2\ (q>2)$. Moreover
vertices of importance to the characterization of quadrics
transform under $H$ as follows: $H[(V_qe_1)^T,0]=[(V_{\ti
p}e_1)^T,0]\ (=[0^T,1])$ for $\ti p>0\ (\ti p=0)$, $q=\ti p+1$ and
$H'$ splits as $H'=I_{1,n-r}H'I_{1,n-r}+H'',\ 0_{n-r,n-r}\oplus
H''\simeq H'' \in\mathbf{GL}_{r}(\mathbb{C}),\ \ti H''=H''^{-1}$.
As $H''H''^T=\bigoplus_{j=2}^m((\ti a_1-\ti
a_j)I_{s_j}+\bigoplus_kJ_{s_{jk}}),\ H''^TH''=
\bigoplus_{j=2}^m\bigoplus_k((a_1-a_j)I_{s_j}+J_{s_{jk}})^{-1},\
H''$ is formed by blocks of the form
$R_1(\sqrt{(a_1-a_j)^{-1}}I_{s_{jk}}+J_{s_{jk}})R_2^T$, where the
rotations $R_1,\ R_2\in\mathbf{O}_{s_{jk}}(\mathbb{C})$ can be
found from a diagonalization process:
$R_1(2\sqrt{(a_1-a_j)^{-1}}J_{s_{jk}}+J^2_{s_{jk}})R_1^T=
R_2((\sqrt{(a_1-a_j)^{-1}}I_{s_{jk}}+J_{s_{jk}})^{-2}-(a_1-a_j)I_{s_{jk}})R_2^T=J_{s_{jk}}$.
As expected, excluding the singular values $a_1,\ \infty$ of $z$
and the corresponding singular values $\infty,\ \ti a_1$ of $\ti
z$, the remaining singular values of $z,\ \ti z$ not only
correspond ($\ti a_j=\ti a_1+(a_j-a_1)^{-1},\ j=2,...,m$), but
also their singular quadrics are projectively equivalent (with the
projective structure induced by their singular sets as previously
described).

We shall now provide an analytic confirmation of these facts: we
have $H'\ti H'+h_1\ti h_2^T=I_n,\ H'\ti h_1=\ti H'^Th_2=0,\
h_2^T\ti h_1=1$. Multiplying the first relation on the left by
$\ker(H'^T)^T$ and on the right by $\ker(\ti H')$ we get
$\ker(H'^T)=\mathbb{C}\ti h_2,\ \ker(\ti H')=\mathbb{C}h_1$; also
$\ti a_1I_n-\ti A^{\circ -1}=H'H'^T,\ H'h_2=-\ti b,\ |h_2|^2=-\ti
C;\ 0=H'H'^T\ti h_2=\ti h_2-\ti A^{\circ -1}\ti h_2$, so $\ti
h_2=e^{\ti c_1}V_{\ti q}e_1;\ -\ti H'\ti b=\ti H'H'h_2=h_2+\ti
C\ti h_1,\ -\ti H'^T\ti H'\ti b=\ti C\ti H'^T\ti h_1=0$, so $\ti
h_1=h_2=e_{p+1}$ for $\ti x_{\ti z}$ QC; $|h_1|^2bb^T=\ti
H'((H'\ti H')^TH'\ti H' -I_n)\ti H'^T=(A^{\circ
-1}-a_1I_n)H'^TH'(A^{\circ -1}-a_1I_n)+(A^{\circ -1}-a_1I_n)$.
Multiplying it on the left by $R_{a_1}^{\circ -1}A$, on the right
by $AR_{a_1}^{\circ -1}$ and using $H'^TH'b=0$ we get
$R_{a_1}^{\circ-1}R_{a_1}H'^TH'R_{a_1}R_{a_1}^{\circ-1}=
-R_{a_1}^{\circ-1}AR_{a_1}R_{a_1}^{\circ-1}$; further multiplying
on the left by $I_n-h_2\ti h_1^T$ and on the right by $I_n-\ti
h_1h_2^T$ we get $H'^TH'=-(I_n-h_2\ti h_1^T)R_{a_1}^{\circ-1}A
R_{a_1}R_{a_1}^{\circ-1}(I_n-\ti h_1h_2^T)$. If $x_z$ are (I)QWC,
then $0=(I_n-h_1\ti h_2^T-H'\ti H')\ti H'^TB=(\ti
H'^T+h_1b^T-a_1H')B=h_1+(\ti H'^T-a_1H')B$, so
$H'^Th_1=(a_1H'^TH'+h_2\ti h_1^T-I_n)B,\
R_{a_1}H'^Th_1=-R_{a_1}(a_1AR_{a_1}^{\circ-1}+I_n)B =-B,\
-|h_1|^2=-|(a_1H'-\ti
H'^T)B|^2=a_1B^T(I_n-a_1H'^TH')B=\\a_1B^TR_{a_1}^{\circ -1}B;\
-\ti H'^TR_{a_1}^{\circ -1}B=\ti H'^TR_{a_1}^{\circ
-1}R_{a_1}H'^Th_1=\ti H'^TH'^Th_1$, so $h_1=-a_1B^TR_{a_1}^{\circ
-1}B\ti h_2-\ti H'^TR_{a_1}^{\circ -1}B$.

The equalities
$R_{a_1}^{\circ-1}R_{a_1}H'^TH'R_{a_1}R_{a_1}^{\circ-1}=
-R_{a_1}^{\circ-1}AR_{a_1}R_{a_1}^{\circ-1},\ R_{a_1}H'^Th_1=-B,\
|h_1|^2=-C-a_1B^TR_{a_1}^{\circ -1}B$ can also be found by
accounting the two equations of
$H^{-1}(C(\infty))=\mathcal{S}(x_{a_1}):\ |H'x+h_1|^2=h_2^Tx=0$
and $x^T(AR_{a_1}^{\circ -1}x+2R_{a_1}^{\circ
-1}B)+a_1B^TR_{a_1}^{\circ -1}B+C=0,\ x=R_{a_1}R_{a_1}^{\circ
-1}x$.

Rewriting the above: $h_2=e^{c_1}V_qe_1,\ \ker(H')=\mathbb{C}\ti
h_1,\ \ker(H'^T)= \mathbb{C}V_{\ti q}e_1,\ H'V_pe_1=Ch_1+e^{\ti
c_1}V_{\ti q}e_1,\ H'V_qe_1=e^{-c_1}V_{\ti p}e_1,\
(I_n-V_qe_1e_1^T\bar V_q^T)H'^Th_1=(I_p+a_1J_p)^{-1}\bar V_pe_1,\
|h_1|^2=(-a_1)^p,\ H'H'^T=(\ti a_1I_{\ti p}+\bar J_{\ti p}) \oplus
J_{\ti q}\oplus\ti A'^{-1},\ (I_n-V_qe_1e_1^T\bar
V_q^T)H'^TH'(I_n-\bar V_qe_1e_1^TV_q^T)
=(I_p+a_1J_p)^{-1}J_p\oplus\bar J_q\oplus A',\ h_1=(-a_1)^pe^{\ti
c_1}V_{\ti q}e_1+ \ti H'^T(I_p+a_1J_p)^{-1}\bar V_pe_1$.

We have
$\ker((H'H'^T)^j)\subseteq\ker(H'^T(H'H'^T)^j)\subseteq\ker((H'H'^T)^{j+1}),\
\ker((H'^TH')^j)\subseteq\\\ker(H'(H'^TH')^j)\subseteq\ker((H'^TH')^{j+1}),\
(H'^T(H'H'^T)^j)^T=H'(H'^TH')^j,\ j\ge 0$, so
$\\\mathrm{dim}(\ker((H'H'^T)^j))\le
\mathrm{dim}(\ker((H'^TH')^{j+1}))\le\mathrm{dim}(\ker((H'H'^T)^{j+2})),\
j\ge 0$. For $M=M^T\in\mathbf{M}_n(\mathbb{C})$ the function $0\le
j\mapsto\mathrm{dim}(\ker(M^j))$ is increasing and stabilizes to
the sum of the dimensions of the SJ blocks of $\ker(RMR^T)$ for
$j\ge k,\ RMR^T$ being the SJ form of $M$ and $k$ the maximum
dimension of the SJ blocks of $\ker(RMR^T)$. Since
$\ker(H'^TH')=\mathbb{C}\ti h_1\oplus\mathbb{C}b,\ |\ti h_1|^2\neq
0$ we have $\mathrm{dim}(\ker((H'^TH')^j))=p+1,\ j\ge
p+\del_{p0}$; since $\mathrm{dim}(\ker((H'H'^T)^j))=\ti q,\ j\ge
\ti q$ we have $\ti q=p+1$.

If $\ti x_{\ti z}$ are QC, then $\ti p=0,\ q=1,\
\ker(H'^T)=\mathbb{C}V_{p+1}e_1,\ \ker(H')= \mathbb{C}e_{p+1}$, so
$H'H'^T=J_{p+1}\oplus \ti A'^{-1},\
H'^TH'=(I_p+a_1J_p)^{-1}J_p\oplus A'$; moreover since
$H'(H'^TH')^{p+1}=(H'H'^T)^{p+1}H',\ H'$ splits as
$H'=I_{1,p+1}H'I_{1,p+1}+H'',\ H''\in\mathbf{GL}_r(\mathbb{C})$.
With $V_{p+1}^TH'V_p=:[0_{p,1}\ \hat I_pU]^T,\
U\in\mathbf{GL}_p(\mathbb{C})$ we have $U^T\hat I_pU=\hat I_p,\
UI'_p\hat I_pU^T=(I_p+a_1I'_p)^{-1}I'_p\hat I_p,\ U^Te_p=e^{\ti
c_1}e_p$; the second equality can be replaced by $I'_pU
=(I_p+a_1I'_p)UI'_p$. With $R:=V_p\hat I_pU^TV_p^T+I_{p+1,n}$ we
have $R\in\mathbf{O}_n(\mathbb{C}),\
J_p=R(I_p+a_1J_p)^{-1}J_pR^T$, so $R$ appears from a
diagonalization process and thus exists and is mostly unique; we
can do better and actually compute $U$. With
$e_j^TUe_k=:(-a_1)^{k-j}u_{jk}$ we have $ \ti c_1=0,\ U-I_p$
strictly upper triangular and

$$(I)\ \ \ \Sigma_{s=0}^ru_{(p+1-j)-s\ (p+1-j)}u_{j+s\ j+r}=\del_{r0},\ j,j+r=1,...,p,$$
$$(II)\ \ \ u_{j+1\ j+r+1}+u_{j+1\ j+r}-u_{j\ j+r}=0,\ j,j+r=1,...,p-1.$$

$u'_{j\ j+r}:=(^{\frac{p}{2}-j}_r),\ j,j+r>0$ satisfy (II) and (I):
$\Sigma_{s=0}^ru'_{(p+1-j)-s\ (p+1-j)}u'_{j+s\ j+r}=\\
(^{\frac{p}{2}-j-1}_r)\Sigma_{s=0}^r(-1)^s(^r_s)=\del_{r0}$.
Suppose $u_{j\ j+s}=u'_{j\ j+s}$ for $j,j+s=1,...,p,\ s=0,...,r-1$
and for $s=r$ and some $j=1,...,p-r$; from (II) it is true for
$j,j+s=1,...,p,\ s=0,...,r$. Thus in order to prove that
$u_{jk}=u'_{jk},\ j,k=1,...,p$ we need: $u_{j\ j+s}=u'_{j\ j+s}$
for $j,j+s,r+1=1,...,p,\ s=0,...,r-1\Rightarrow u_{j\ j+r}=u'_{j\
j+r}$ for some $j=1,...,p-r$. For $p=2k$ if $r=2m-1$, then let
$j=k-m+1$ in (I) to get $u_{k-m+1\ k+m}=u'_{k-m+1\ k+m}$; if
$r=2m$, then let $j=k-m$ in (I) and (II) to get $u_{k-m\
k+m}=u'_{k-m\ k+m}$; for $p=2k-1$ if $r=2m-1$, then let $j=k-m$ in
(I) and (II) to get $u_{k-m+1\ k+m}=u'_{k-m+1\ k+m}$; if $r=2m$,
then let $j=k-m$ in (I) to get $u_{k-m\ k+m}=u'_{k-m\ k+m}$. Thus
$U=U(p,a_1):=\Sigma_{j=1}^pe_je_j^T(I_p-
a_1I'_p)^{\frac{p}{2}-j},\
U^{-1}=\Sigma_{j=1}^pe_je_j^T(I_p+a_1I'_p)^{\frac{p}{2}-j}$ (we
have
$\\\Sigma_{r=0}^l(-1)^{l-r}(^{\frac{p}{2}-j}_r)(^{\frac{p}{2}-j-r}_{l-r})=
\Sigma_{r=0}^l(^{\frac{p}{2}-j}_r)(^{j-\frac{p}{2}+l-1}_{l-r})=\del_{l0}$
as the coefficient of $x^l$ in
$(1+x)^{\frac{p}{2}-j}(1+x)^{j-\frac{p}{2}+l-1}=(1+x)^{l-1}$). Now
$H'=V_{p+1}[\hat I_p\ 0_{p,1}]^TU^TV_p^T+H'',\ h_1=V_{p+1}(I_{p+1}
-a_1I'_{p+1})^{\frac{p}{2}}e_{p+1},\\ \ti H'=V_p[\hat
I_p\Sigma_{j=1}^pa_1^j(^{\frac{p}{2}}_j)e_j\ \ \ U\hat
I_p]V_{p+1}^T+H''^{-1}$. Note that $e_j^T(U\hat I_p)e_k
=a_1^{p+1-j-k}(^{\frac{p}{2}-k}_{p+1-j-k}),\ j,k=1,...,p$; for
$k=0$ we get the first column $e_j^T(\hat
I_p\Sigma_{l=1}^pa_1^l(^{\frac{p}{2}}_l)e_l)
=a_1^{p+1-j}(^{\frac{p}{2}}_{p+1-j}),\ j=1,...,p$ of $V_p^T\hat
I_p\ti H'\hat I_{p+1}V_{p+1}$ and for $j=0$ we get  the first $p$
entries $e_k^T(I_{p+1}-a_1I'_{p+1})^{\frac{p}{2}}e_{p+1}=\\
(-a_1)^{p+1-k}(^{\frac{p}{2}}_{p+1-k}),\ k=1,...,p$ of
$V_{p+1}^T\hat I_{p+1}h_1$.

If $x_z,\ \ti x_{\ti z}$ are (I)QWC, then $p,\ \ti p>0,\
H'H'^TV_{\ti p}=V_{\ti p}(\ti a_1I_{\ti p}+ I'^T_{\ti p}),\
H'^TH'V_{\ti p+1}I_{1,\ti p}=V_{\ti p+1}(e_1v^T+I'^T_{\ti p}),\
v^T:=e_{\ti p+1}^TV_{\ti p+1}^TH'^TH'V_{\ti p+1}I_{1,\ti p}$. By
induction after $j\ge 0$ we have $H'V_{\ti p+1}(e_1v^T+I'^T_{\ti
p+1})^je_1 =e^{-c_1} V_{\ti p}(\ti a_1I_{\ti p}+I'^T_{\ti
p})^je_1,\ j=0,...,\ti p,\ V_{\ti p+1}(e_1v^T+I'^T_{\ti
p+1})^{j+1}e_1=e^{-c_1}H'^TV_{\ti p} (\ti a_1I_{\ti p}+I'^T_{\ti
p})^je_1,\ j=0,...,\ti p-1$; thus $H'^TV_{\ti p}=V_{\ti p+1}(\bar
V_{\ti p+1}^TH'^TV_{\ti p})$ and $\ti h_1=V_{\ti p+1}((-\ti
a_1)^{\ti p}e^{c_1}e_1 +\bar V_{\ti p+1}^TH'^TV_{\ti p}(I_{\ti
p}+\ti a_1I'_{\ti p})^{-1}e_{\ti p}),\\
H'^TH'=(I_p+a_1J_p)^{-1}J_p\oplus(I_{\ti p+1}-h_2\ti h_1^T)\bar
J_{\ti p+1}(I_{\ti p+1}- \ti h_1h_2^T)\oplus A'$. Since
$H'(H'^TH')^{p+1}=(H'H'^T)^{p+1}H'$, we have $V_{\ti p}^TH'V_p=0,\
I_rH'V_p=0,\ V_{p+1}^TH'I_r=0$; further since
$H'(H'^TH')=(H'H'^T)H'$, we have $V_{p+1}^TH'V_{\ti
p+1}\in\mathbb{C}e_{\ti p+1}\ti h_1^TV_{\ti p+1}$. As
$V_{p+1}^TH'V_{\ti p+1}e_1=e^{-c_1}V_{p+1}^TV_{\ti p}e_1=0$ we
have $V_{p+1}^TH'V_{\ti p+1}=0$; now $V_{p+1}^TH'V_p=:[0_{p,1}\
\hat I_pU]^T,\ U\in\mathbf{GL}_p(\mathbb{C})$ must be as before,
so $\ti c_1=0,\ U=U(p,a_1)$. As $\ti H'H'+\ti h_1h_2^T=I_n$ we
have $(V_{\ti p+1}^T\ti H'V_{\ti p})\hat I_{\ti p} (V_{\ti
p}^TH'I_r)=0,\ (I_r\ti H'I_r)(I_rH'V_{\ti p+1})=0,\ (I_r\ti
H'I_r)(I_rH'I_r)=I_r$, so $V_{\ti p}^TH'I_r=0,\ I_rH'V_{\ti
p+1}=0$ and finally $\\H'=V_{p+1}[\hat I_p\
0_{p,1}]^TU(p,a_1)^TV_p^T+V_{\ti p}[\hat I_{\ti
p}\Sigma_{j=1}^{\ti p}\ti a_1^j (^{\frac{\ti p}{2}}_j)e_j\ \ \
U(\ti p,\ti a_1)\hat I_{\ti p}]V_{\ti p+1}^T+H'',\ \ti h_1=V_{\ti
p+1}(I_{\ti p+1}-\ti a_1I'_{\ti p+1})^{\frac{\ti p}{2}}e_{\ti
p+1}$.

A homography $H$ does not have the property of taking a family of
confocal quadrics to another iff it is an affine transformation
with $H(C(\infty))\neq C(\infty)$ or $\begin{vmatrix}H'H'^T-\la
I_n&H'h_2\\(H'h_2)^T&|h_2|^2
\end{vmatrix}=0,\ \forall\la\in\mathbb{C}$; in the later case $H^{-1}(C(\infty)):\
|H'x+h_1|^2=h_2^Tx+a=0$ (and thus $H(C(\infty))$) cannot appear as
the singular set of a singular quadric of a family of confocal
quadrics. By applying a rigid motion in the domain of $H$ and
fixing the constant upon which $H$ depends we make $a=0,\ h_2=f_1\
(=e_1)$ according to being isotropic (or not) and bring
$H^{-1}(C(\infty))$ to the form of (I)Q(W)C in $f_1^Tx=0\
(e_1^Tx=0):\ \ker(A)$ is spanned by $f_1,\ \{f_1,e_3\},\
\{f_1,f_2\}\ (e_1,\ \{e_1,e_3\},\ \{e_1,f_2\}),\ B=0,e_3,-\bar
f_2,\ C=-1,0,0$ or a quadric of the form $f_1^Tx=x^T(Ax-2\bar
f_1)=0,\ A=A^T,\ \ker(A)=I_{1,2}\mathbb{C}^n$. As we have seen,
$H^{-1}(C(\infty))$ must be of this last type; thus $|h_1|^2=0,\
(I_n-f_1\bar f_1^T)H'^Th_1=-e^c\bar f_1,\ H'^TH'f_1=c_1f_1,\
|H'f_1|^2=0$; by applying a rotation in the range of $H$ we make
$h_1=e^cf_1$, so $H'^Tf_1=-\bar f_1+c_2f_1,\ f_1^TH'f_1=-1$. Since
$\{R(h)|h=I_{3,n}h\in\mathbb{C}^n\}$ acts transitively on
$\{x\in\mathbb{C}^n|\ |x|^2=0,\ f_1^Tx=-1\}$, by applying a
rotation $R(h)$ in the range of $H$ we further make $H'f_1=-\bar
f_1$, so $H'^T\bar f_1=-c_1f_1,\ H'=(f_1\bar f_1^T+\bar f_1f_1^T
+I_{3,n})H'(f_1\bar f_1^T+\bar f_1f_1^T+I_{3,n})=h_3f_1^T-\bar
f_1\bar f_1^T+I_{3,n}H'',\ h_3\in\mathbb{C}^n,\
H''\in\mathbf{GL}_n(\mathbb{C}),\
I_{3,n}H''=H''I_{3,n}=H''-I_{1,2}$;  further applying the
translation $-h_3$ in the range of $H$ we make $h_3=0,\
H'H'^T=I_{3,n}H''H''^T,\ \begin{vmatrix}H'H'^T-\la
I_n&H'h_2\\(H'h_2)^T&|h_2|^2\end{vmatrix} =\begin{vmatrix}I_n&-\la
f_1\\0^T&1\end{vmatrix}
\begin{vmatrix}I_{3,n}H''H''^T-\la I_n&-\bar f_1\\-\bar f_1^T&0\end{vmatrix}
\begin{vmatrix}I_n&0\\-\la f_1^T&1\end{vmatrix}
=\begin{vmatrix}I_{3,n}H''H''^T-\la I_{3,n}&-\bar f_1\\-\bar f_1^T&0\end{vmatrix}=0$.

Summing up:

{\it Any homography
$H:=\begin{bmatrix}H'&h_1\\h_2^T&a\end{bmatrix}\in\mathbf{PGL}_n(\mathbb{C}),\
H'\in\mathbf{M}_n(\mathbb{C}),\ h_1,\ h_2\in\mathbb{C}^n,\
a\in\mathbb{C}$, with inverse $(H^{-1}=)\ti H:=\begin{bmatrix}\ti
H'&\ti h_1\\\ti h_2^T&\ti a\end{bmatrix}$ can be brought, via
rigid motions and homotheties in both its domain and range, to one
of the two cases below:

(I) (i) $H=\ti H=I_{n+1}$; (ii) $a=\ti a=1,\ h_1=\ti h_1=h_2=\ti
h_2=0,\ \ti
H'^{-1}=H'\in\mathbf{GL}_n(\mathbb{C})\setminus\{\mathbb{C}I_n\}$
SJ; (iii) $a=\ti a=0,\ e^{-c}h_1=h_2=f_1,\ \ti h_1=e^{c}\ti
h_2=\bar f_1,\ H'=-\bar J_2+I_{3,n}H'',\ \ti
H'=-J_2+I_{3,n}H''^{-1},\ c\in\mathbb{C},\
H''\in\mathbf{GL}_n(\mathbb{C})$ SJ, $I_{3,n}H''=H''-I_{1,2}$.

(II) $a=\ti a=0$ and $H$ takes confocal (I)Q(W)C $x_z$ to confocal
(I)Q(W)C $\ti x_{\ti z},\ \ti a_1:=a_1^{-1},\ \ti z:=\ti
a_1+(z-a_1)^{-1}$. By symmetry all the following remain true if
tilde and non-tilde quantities are exchanged: with
$U(p,a_1):=\Sigma_{j=1}^pe_je_j^T(I_p-a_1I'_{p})^{\frac{p}{2}-j}$
we have $h_1=V_{p+1}(I_{p+1}-a_1I'_{p+1})^{\frac{p}{2}}e_{p+1},\
h_2=V_{\ti p+1}e_1,\ H'=V_{p+1}[\hat I_p\
0_{p,1}]^TU(p,a_1)^TV_p^T+ V_{\ti p}[(\hat I_{\ti
p}\Sigma_{j=1}^{\ti p}\ti a_1^j(^{\frac{\ti p}{2}}_j)e_j)\ U(\ti
p,\ti a_1)\hat I_{\ti p}]V_{\ti p+1}^T+H'',\
H''\in\mathbf{GL}_{n-p-\ti p+1}(\mathbb{C}),\ \ti H'':=H''^{-1},\
H''H''^T= \bigoplus_{j=2}^m((\ti a_1-\ti
a_j)I_{s_j}+\bigoplus_kJ_{s_{jk}}),\ H''^TH''
=\bigoplus_{j=2}^m\bigoplus_k((a_1-a_j)I_{s_{jk}}+J_{s_{jk}})^{-1},\
\Sigma_ks_{jk}=s_j>0,\ j=2,...,m,\ \Sigma_js_j=n-p-\ti p-1,\ \ti
a_j=\ti a_1+(a_j-a_1)^{-1},\ j=2,...,m$.

Thus $H''$ is formed by blocks of the form
$R_1(\sqrt{(a_1-a_j)^{-1}}I_{s_{jk}}+J_{s_{jk}})R_2^T$, where the
rotations $R_1,\ R_2\in\mathbf{O}_{s_{jk}}(\mathbb{C})$ can be
found from a diagonalization process:
$J_{s_{jk}}=R_1(2\sqrt{(a_1-a_j)^{-1}}J_{s_{jk}}+J^2_{s_{jk}})R_1^T=
R_2((\sqrt{(a_1-a_j)^{-1}}I_{s_{jk}}+J_{s_{jk}})^{-2}-(a_1-a_j)I_{s_{jk}})R_2^T$.

Also $A^{\circ -1}=\bar J_p\oplus(a_1I_{\ti p+1}-J_{\ti
p+1})\oplus\bigoplus_{j=2}^m\bigoplus_k(a_jI_{s_{jk}}-
J_{s_{jk}}),\ B=-\bar V_pe_1,\ C=-\del_{p0}$.

Conversely, given a quadric in a hyperplane of $\mathbb{C}^n$,
excluding quadrics which belong, after a rigid motion, to the type
$f_1^Tx=x^T(Ax-2\bar f_1)=0,\ A$ SJ,
$\ker(A)=I_{1,2}\mathbb{C}^n$, one can consider it to be
$\mathcal{S}(x_{a_1})$ and complete it to a family of confocal
quadrics $x_z$; then one can recover $H,\ \ti H$ and $\ti x_{\ti
z}$ as above.}

\subsection{Parametrization of confocal isotropic quadrics without center}
\noindent

\noindent Consider the equilateral QWC $Z^T(I_{1,n-1}Z-2e_n)=0$
and a canonical IQWC $x_0^T(Ax_0-2\bar f_1)=0,\
\ker(A)=\mathbb{C}f_1,\ A=J_p\oplus ...$ SJ. We are looking for a
linear map $L\in\mathbf{GL}_n(\mathbb{C})$ such that $x_0=LZ$,
equivalently $L^TAL=e^{2a}I_{1,n-1},\ L^T\bar f_1=e^{2a}e_n$.
Replacing $L$ with $L(e^{-a}I_{1,n-1}+e^{-2a}e_ne_n^T)$ we can
make $a=0$. Thus $Le_n=f_1,\ L^T(A+\bar f_1\bar f_1^T)L=I_n$, so
$L^{-1}=R^T\sqrt{A+\bar f_1\bar f_1^T},\ R^TR=I_n$ with
$Re_n=\sqrt{A+\bar f_1\bar f_1^T}f_1$ (note that $Re_n$ has, as
required, length $1$). Once $R\in\mathbf{O}_n(\mathbb{C})$ with
the above property is found, $L$ thus defined satisfies $L^T\bar
f_1=e_n$ and thus $L^TAL=I_{1,n-1}$. $L$ with the above properties
is unique up to rotations fixing $e_n$ in its domain and a
canonical choice reveals itself:
$I_{1,n-1}L^{-1}\sqrt{R_z}LI_{1,n-1}=\sqrt{I_{1,n-1}-zL^TA^2L}=:
\sqrt{R'_z}-e_ne_n^T,\ A':=L^TA^2L,\
\ker(A')=\mathbb{C}e_n\oplus\mathbb{C}L^{-1}A^{\circ -1}f_1
=\mathbb{C}e_n\oplus\mathbb{C}L^Tf_1$; choose $R$ which makes $A'$
SJ. We have $LL^T=A^{\circ -1}+f_1f_1^T,\
L^{-1}\sqrt{R_z}Le_n=e_n$, so
$I_{1,n-1}\sqrt{R'_z}=I_{1,n-1}L^{-1}\sqrt{R_z}L=L^TA\sqrt{R_z}L$.

Consider now the ruling families parametrization on $Z$ for
$n=2k+1:\ u_j:=\bar f_j^TZ,\ v_j:=f_j^TZ,\ j=1,...,k$; thus
$Z=Z(u_j,v_j)=\Sigma_{j=1}^k(u_jf_j+v_j\bar f_j+ u_jv_je_n)$; if
we let $u_k=v_k$, then the $2k^{\mathrm{th}}$ coordinate is $0$
and we obtain the ruling families parametrization on $Z$ for
$n=2k$. Let
$e_n^TL^{-1}\sqrt{R_z}L=:\Sigma_{j=1}^k(u_j(z)f_j+v_j(z)\bar
f_j)^T+e_n^T$, $u_{jz}:=u_j+u_j(z),\ v_{jz}:=v_j+v_j(z),\
Z^z=Z^z(u_j,v_j):=Z(u_{jz},v_{jz}),\
Z^z(0):=Z^z(0_j,0_j)=I_{1,n-1}L^T\sqrt{R_z}(L^T)^{-1}e_n+\\
\frac{1}{2}|I_{1,n-1}L^T\sqrt{R_z}(L^T)^{-1}e_n|^2e_n$. We have
$I_{1,n-1}\sqrt{R'_z}Z^z(0)=L^T\sqrt{R_z}AA^{\circ
-1}\sqrt{R_z}(L^T)^{-1}e_n= \\L^T(I_n-\sqrt{R_z}-zA)\bar f_1,\
e_n^TZ^z(0)=\frac{1}{2}\bar f_1^T\sqrt{R_z}A^{\circ
-1}\sqrt{R_z}\bar f_1;\ -2\pa_z(\sqrt{R_z}A^{\circ
-1}\sqrt{R_z})=\sqrt{R_z}(R_z^{-1}AA^{\circ -1} +A^{\circ
-1}AR_z^{-1})\sqrt{R_z}=2I_n-(\sqrt{R_z})^{-1}\bar f_1f_1^T-
f_1\bar f_1^T(\sqrt{R_z})^{-1}$, so $e_n^TZ^z(0)=\bar f_1^TC(z)$;
further $Z^z=Z+(e_nZ^TI_{1,n-1}+I_n)Z^z(0)$ and
$L^{-1}x_z-\sqrt{R'_z}Z^z+2zI_{1,n-1}\pa_z|_{z=0}Z^z(0)=
L^{-1}C(z)-\sqrt{R'_z}Z^z(0)-zL^TA(L^T)^{-1}e_n=(L^{-1}-e_n\bar
f_1^T)C(z)-L^T(I_n-\sqrt{R_z}-zA)\bar f_1-zL^TA\bar f_1=0$ (use
$L^{-1}=L^TA+e_n\bar f_1^T,\ AC(z)-I_n-\sqrt{R_z})\bar f_1=0$).
Thus $x_z=L\sqrt{R'_z}Z^z+z\bar f_1$; with $f:=e_n\ (\bar f_1)$
respectively for (I)QWC and $A':=A,\ L:=(\sqrt{A+e_ne_n^T})^{-1}$
for QWC we have a parametrization for (I)QWC:
$x_z=L\sqrt{R'_z}Z^z+z(1-\frac{|f|^2}{2})(L^T)^{-1}e_n$. The TC
becomes: $(V_0^1)^T\hat N_0^0=Z_0^TI_{1,n-1}\sqrt{R'_z}Z_1
-(Z_0+Z_1)^T(I_{1,n-1}Z^z(0)+e_n)-e_n^TZ^z(0)-\frac{|f|^2}{2}z=0$,
so IQWC can be considered metrically degenerated QWC.

\subsection{Totally real quadrics in $\mathbb{C}^n$}
\noindent

\noindent Note that if $V\subseteq\mathbb{C}^n$ is a real
$k$-dimensional subspace (a-priori $0\le k\le 2n$) such that the
Euclidean scalar product on $\mathbb{C}^n$ induces a
non-degenerate real (valued) scalar product on $V$, then the
orthogonal basis of $V$ (formed by vectors of length $\pm 1$)
provided by Sylvester's is, by complexification, an orthogonal
base of a (complex) $k$-dimensional subspace of $\mathbb{C}^n$, so
$k\le n$.

Consider a real $n$-dimensional subspace $V\subset\mathbb{C}^n$ as
above and a basis $\{v_j\}_{j=1,...,n}$ of $V$ with
$v_j^Tv_k=\ep_{j}\del_{jk},\ \ep_j:=\pm 1,\ j=1,...,n,\
\mathcal{E}:=\mathrm{diag}[\ep_1\ ...\ \ep_n]$. Now $R:=[v_1\ ...\
v_n]\sqrt{\mathcal{E}}\in\mathbf{O}_n(\mathbb{C})$ and
$R^TV=\sqrt{\mathcal{E}}\mathbb{R}^n$ satisfies
$\overline{x}=\mathcal{E}x,\ \forall x\in R^TV$.

We are interested in the quadrics which intersect a totally real
subspace $\sqrt{\mathcal{E}}\mathbb{R}^n$ along an
$(n-1)$-dimensional totally real quadric and bringing such
quadrics to the canonical form (a form depending on as few as
possible constants) by means of rigid motions preserving
$\sqrt{\mathcal{E}}\mathbb{R}^n$. The rigid motions preserving
$\sqrt{\mathcal{E}}\mathbb{R}^n$ are $(R,t),\
R\in\mathbf{O}_n(\mathbb{C}),\ \bar R=\mathcal{E}R\mathcal{E},\
t\in\sqrt{\mathcal{E}}\mathbb{R}^n$; the homotheties preserving
$\sqrt{\mathcal{E}}\mathbb{R}^n$ are $e^{a},\ a\in\mathbb{R}$.

However in order to take advantage of the SJ form it is more
convenient to intersect a quadric in canonical form with a totally
real affine space $(R_0,t_0)\sqrt{\mathcal{E}}\mathbb{R}^n,\
(R_0,t_0)\in\mathbf{O}_n(\mathbb{C})\ltimes\mathbb{C}^n$ and put
the condition that the sub-manifold thus obtained has dimension
$n-1$. Note that $(R_0,t_0)$ is uniquely determined modulo
compositions on the right with rigid motions preserving
$\sqrt{\mathcal{E}}\mathbb{R}^n$, so $R:=\bar
R_0\mathcal{E}R_0^T,\ t:=t_0-R^T\bar t_0$ are uniquely determined.
Note $R\in\mathbf{O}_n(\mathbb{C}),\ \bar R^T=R,\ \det
R=\det\mathcal{E}$; conversely given $R$ as above one can find
$R_0$ as above (which is uniquely determined modulo multiplication
on the right with rotations which preserve
$\sqrt{\mathcal{E}}\mathbb{R}^n$). This is a standard dimensional
argument: since the map $\mathbf{O}_n(\mathbb{C})\ni
R_0\mapsto\bar R_0\mathcal{E}R_0^T
\in\{R\in\mathbf{O}_n(\mathbb{C})|\bar R^T=R,\ \det
R=\det\mathcal{E}\}$ is continuous with connected codomain of real
dimension $\frac{n(n-1)}{2}$ and kernel of the same real
dimension, it is onto. Similarly because of the real decomposition
$\mathbb{C}^n\ni v=\frac{v+\mathcal{E}\bar
v}{2}+\frac{v-\mathcal{E}\bar
v}{2}\in\sqrt{\mathcal{E}}\mathbb{R}^n\oplus
i\sqrt{\mathcal{E}}\mathbb{R}^n$, from $t$ with $\bar t=-Rt$ one
can obtain $t_0$ as above (which is unique modulo
$R_0\sqrt{\mathcal{E}}\mathbb{R}^n$). Thus we need
$\begin{bmatrix}R_0&t_0\\0&1\end{bmatrix}^T
\begin{bmatrix}A&B\\B^T&C\end{bmatrix}\begin{bmatrix}R_0&t_0\\0&1\end{bmatrix}=
e^{2c}\begin{bmatrix}\mathcal{E}&0\\0&1\end{bmatrix}
\begin{bmatrix}\bar R_0&\bar t_0\\0&1\end{bmatrix}^T
\begin{bmatrix}\bar A&\bar B\\\bar B^T&\bar C\end{bmatrix}
\begin{bmatrix}\bar R_0&\bar t_0\\0&1\end{bmatrix}
\begin{bmatrix}\mathcal{E}&0\\0&1\end{bmatrix}$;
if such a condition is not satisfied, then applying the rigid
motion $(R_0,t_0)^{-1}$ and using the real decomposition
$\mathbb{C}^n=\sqrt{\mathcal{E}}\mathbb{R}^n\oplus
i\sqrt{\mathcal{E}}\mathbb{R}^n$ wrt
$\sqrt{\mathcal{E}}\mathbb{R}^n$ the considered sub-manifold
becomes the intersection of two (possibly empty or degenerate)
quadrics in $\sqrt{\mathcal{E}}\mathbb{R}^n$ and thus it has
dimension at most $n-2$.

Applying the compatibility condition conjugation to the above
relation we get $c\in i\mathbb{R}$; multiplying
$\begin{bmatrix}A&B\\B^T&C\end{bmatrix}$ with $e^{-c}$ we make
$c=0$ (note however that by doing this $B$ and $C$ will change
from their usual values for quadrics in canonical form).

Now $\bar A=RAR^T,\ \bar B=R(At+B),\ \bar C=t^T(B+R^T\bar B)+C$;
since $\bar J_p$ is obtained from $J_p$ by applying reflections in
the $e_{2j}^Tx=0$ hyperplanes (which take $f_j$ to $\bar f_j$),
multiplying $R$ on the left with the product $r$ of these
reflections we obtain the rotation $rR$ preserving the SJ form of
a matrix and with a determinant and real condition. Thus the set
of eigenvalues of $A$ for SJ blocks of the same dimension is
closed under conjugation, so eigenvalues of $A$, besides being
either real or complex conjugate have corresponding SJ blocks of
the same dimension. Moreover $R$ exchanges the subspaces of
$\mathbb{C}^n$ on which the SJ blocks of the eigenvalues $a_j$ and
$\bar a_j$ are supported, although it may not preserve the
subspaces corresponding to the SJ blocks of same dimension, so $R$
admits the corresponding block diagonal decomposition. If the
corresponding blocks of $\mathcal{E}$ and of $R$ have the same
determinant, then one can choose $R_0$ with the same block
decomposition, following that the blocks on which the determinants
differ will be grouped by two to get a block for $R_0$; the
simplest choice of $R_0$ will be the one giving so far the form
closest to the canonical form; further simplifications will appear
below. For the blocks $\bigoplus_{k=1}^{s_j}J_{p_k}$ corresponding
to eigenvalue $a_j\in\mathbb{R}$ we have the corresponding block
$R_j$ of $R$ with $R_j(\bigoplus_{k=1}^{s_j}J_{p_k})R_j^T
=\bigoplus_{k=1}^{s_j}\bar J_{p_k}$ (which must still satisfy the
reality condition $\bar R_j^T=R_j$). For complex conjugate
eigenvalues $a_j,\bar a_j\in\mathbb{C}\setminus\mathbb{R}$ we have
the corresponding SJ block $(a_jI_l+\bigoplus_{k=1}^{s_j}J_{p_k})
\bigoplus(\bar a_jI_l+\bigoplus_{k=1}^{s_j}J_{p_k})$ of dimension
$2l$; by considering $R(A-a_jI_n)^kR^T,\ k\ge 1$ we get
$R_j=\begin{bmatrix}0&Q_j\\\bar Q_j^T&0\end{bmatrix},\
Q_j\in\mathbf{O}_l(\mathbb{C})$ so we have reduced the problem to
finding $Q_l\in\mathbf{O}_l(\mathbb{C})$ with
$Q_l(\bigoplus_{k=1}^{s_j}J_{p_k})Q_l^T=\bigoplus_{k=1}^{s_j}J_{p_k}$
(without a reality condition).

Now we can already infer that most quadrics do not contain totally
real quadrics and that totally real quadrics for which
$\mathrm{spec}(A)\neq\{0\}$ belong to totally real confocal
families obtained for $z\in\mathbb{R}$ (the necessity follows from
the requirement on the eigenvalues of $A$; the sufficiency is
clear from the definition of the confocal family). Note also that
only for diagonal QC with $A=-\mathcal{E}|A|$ we have empty
totally real quadrics, but this can be corrected by considering
pairs $(\mathcal{E},-\mathcal{E})$ instead of only $\mathcal{E}$.
This pairing is natural, since multiplication by $i$ exchanges
both $\mathcal{E}$ to $-\mathcal{E}$ and the signature of the
induced linear element of the totally quadric. According to
Sylvester's this linear element has signature
$(\frac{\mathrm{tr}(\mathcal{E}+|\mathcal{E}|)}{2},
\frac{\mathrm{tr}(\mathcal{E}-|\mathcal{E}|)}{2}-1)$ along the
regions where the normal is time-like and
$(\frac{\mathrm{tr}(\mathcal{E}+|\mathcal{E}|)}{2}-1,
\frac{\mathrm{tr}(\mathcal{E}-|\mathcal{E}|)}{2})$ along the
regions where the normal is space-like.

For QC we have $t=0,\ C=\pm 1$; by a real homothety one can make
$C=-1$.

For QWC we have $Re_n=\ep e_n,\ \ep:=\pm 1,\ B=-e^{-c}e_n,\
e^{2c}=\ep,\ t=0$.

For IQWC if the SJ block of $f_1$ is $J_p$, then $R$ admits a
corresponding block diagonal decomposition and on the first $p$
coordinates $R$ is $\pm$ the product of the reflections which take
$J_p$ to $\bar J_p$, so $Rf_1=\ep\bar f_1,\ R\bar f_1=\ep f_1,\
\ep:=\pm 1,\ B=-e^{-c}\bar f_1,\ e^{2c}=\ep,\ t=0$.

Therefore in all cases one can choose $t_0:=0$ and the ambient
totally real affine space passes through $0$.

In order to get a classification of totally real quadrics and
their confocal families we need to exclude cases obtained as above
and which coincide. Since by complexification of two rigidly
applicable totally real quadrics we obtain two rigidly applicable
(complex) quadrics, two such totally real quadrics must belong to
the same complex metric type; moreover the rigid motion by means
of which the above rigidly applicability correspondence is
realized must be a rigid motion which preserves the canonical form
and conversely a rigid motion of the later type takes totally real
quadrics to totally real quadrics. Such a rigid motion cannot
preserve a totally real affine $n$-dimensional space, since a
quadric intersects such a space in at most a totally real quadric.
The group of rigid motions preserving the canonical form of a
quadric acts transitively on each type of totally real quadrics in
that quadric. Thus the last step in finding the canonical form of
a totally real quadric is revealed from rotations preserving the
SJ form: among all such rotations $R'$ choose one which makes
$\bar R'RR'^T$ as simple as possible and then continue with the
program described above to find the simplest $R_0$.

Thus totally real quadrics naturally appear in canonical form not
in totally real subspaces $\sqrt{\mathcal{E}}\mathbb{R}^n$, but in
totally real subspaces $R_0\sqrt{\mathcal{E}}\mathbb{R}^n$. If one
takes this to be the definition of the canonical form of totally
real quadrics, then with the exception of finding the signature
$\mathcal{E}$ one needs no work from a metric point of view beyond
finding the canonical form of the corresponding (complex) quadric
and checking the conditions on the eigenvalues and their
corresponding SJ blocks.

\section{Algebraic preparatives for quadrics in
$\mathbb{C}^3$}\label{sec:algpre2}\setcounter{equation}{0}

\subsection{Canonical confocal quadrics and homographies of Bianchi II}\label{subsec:algpre21}
\noindent

\noindent In what follows we shall restrict our attention to
quadrics in $\mathbb{C}^3$; we shall choose the ruling families
parametrization invariant under the Ivory affinity. Intersections
of isotropic rulings in $\mathbb{C}^3$ give all umbilics (because
the first and second fundamental forms are proportional); also
because of the preservation of lengths of segments on rulings
under the Ivory affinity umbilics are preserved by the Ivory
affinity. Excluding (pseudo-)spheres, the remaining four points of
intersections of isotropic rulings (counting multiplicities) are
situated on the circle at $\infty ,\ C(\infty)$ (B\'{e}zout). When
the four points are distinct we get general diagonal Q(W)C:
$x_z\cap C(\infty)=x_z\pitchfork C(\infty)$; for higher
multiplicities the quadrics are tangent to $C(\infty):\ x_z\cap
C(\infty)= (x_z\cap^TC(\infty))\cup(x_z\pitchfork C(\infty))$. If
there are two points with double multiplicities, then we get
diagonal Q(W)C and of revolution; in all other cases we get
Darboux quadrics ($x_z\cap^TC(\infty)$ consists of one point);
they cannot be realized as quadrics in $\mathbb{R}^3$, but if we
intersect them with certain totally real spaces, then we get
positive definite linear element so we get real deformations of
these quadrics. If $x_z\cap\mathbb{CP}^2$ is non-degenerate, then
$x_z$ are QC; otherwise $x_z\cap\mathbb{CP}^2$ is a cone whose
vertex belonging or not to $C(\infty)$ decides that $x_z$ are IQWC
or QWC.

Let $X\subset\mathbb{C}^3$ be the unit sphere $|X|^2=1$ and
$u:=\frac{e_3^TX+1}{f_1^TX},\ v:=\frac{e_3^TX-1}{f_1^TX},\ u\neq
v\in \mathbb{C}\cup\infty$ be the parametrization of $X$ by its
ruling families: $X=X(u,v)=\frac{-uvf_1+2\bar f_1+(u+v)e_3}{u-v}$.
$X$ has asymptotic isotropic cone $u=v=\tau:\ \infty Y(\tau),\
Y(\tau):=-\tau^2f_1+2\bar f_1+2\tau e_3$. Consider
$\si_1:(u,v)\mapsto (iu,iv)$ and the involutions $\si_1^2,\
\si_2,\si_3,\si_4,\si_5:(u,v)\mapsto(-u,-v),\ (v,u),\
(2u^{-1},2v^{-1}),\ (\bar u,\bar v),\
(\frac{2-\sqrt{2}u}{\sqrt{2}+u},\frac{2-\sqrt{2}v}{\sqrt{2}+v})$;
note the relations $\si_1\si_3=\si_3\si_1^{-1},\
\si_1\si_4=\si_4\si_1^{-1},\ \si_1\si_5=\si_1^2(\si_5\si_1)^2,\
\si_3\si_5=\si_5\si_1^2$. The group generated by reflections in
the planes $e_1^Tx,e_2^Tx,e_3^Tx=0$ corresponds to the group
generated by $\si_2\si_3,\ \si_1^2\si_2\si_3,\ \si_1^2\si_2$ and
the group of permutations of coordinates generated by
$e_1\leftrightarrow e_2,\ e_2\leftrightarrow e_3,\
e_3\leftrightarrow e_1$ corresponds to the group generated by
$\si_1\si_2\si_3,\ \si_1\si_2\si_3\si_5\si_1,\ \si_2\si_5$. Let
$Z\subset \mathbb{C}^3$ be the equilateral paraboloid
$Z^T(I_{1,2}Z-2e_3)=0$ and $u:=\bar f_1^TZ,\ v:=f_1^TZ,\ u,\
v\in\mathbb{C}$ be a parametrization of $Z$ by its ruling
families: $Z=Z(u,v)=uf_1+v\bar f_1+uve_3$. We have:
\begin{eqnarray}\label{eq:XY}
f_1\times\bar f_1=ie_3,\ \bar f_1\times e_3=i\bar f_1,\ e_3\times
f_1=if_1;\ Y(u)=\frac{u}{2}(Y'(u)+\overline{Y'(-2\bar
u^{-1})})=-\frac{u^2}{2}\overline{Y(-2\bar u^{-1})};\nonumber\\
X(u,v)=\frac{1}{u-v}Y(u)-\frac{1}{2}Y'(u)
=\frac{1}{u-v}Y(v)+\frac{1}{2}Y'(v);\ X\circ\si_2=-X,\
\bar X=X\circ(\si_1^2\si_2\si_3\si_4);\nonumber\\
Z\circ\si_2=(e_2^Tx\mapsto -e_2^Tx)\circ Z,\ \bar Z=Z\circ(\si_2\si_4);\
Y(u)^TY(v)=-2(u-v)^2,\nonumber\\Y(u)\times Y'(u)=2iY(u),\ Y(u)\times Y(v)=-2i(u-v)^2X(u,v);\
|dX|^2=\frac{4du\odot dv}{(u-v)^2},\nonumber\\
X_{uv}=-\frac{2}{(u-v)^2}X;\ I_3=-\frac{Y(u)Y(v)^T+Y(v)Y(u)^T}{2(u-v)^2}+X(u,v)X(u,v)^T.
\end{eqnarray}
(we use the notation $\phi\odot\varphi:=\frac{\phi\otimes\varphi+\varphi\otimes\phi}{2}$
for $\phi,\varphi$ $(0,1)$ tensors).
Note also:
\begin{eqnarray}\label{eq:detA}
(Ax)\times(Ay)=(\det A)(A^T)^{-1}(x\times y),\  A\in\mathbf{GL}_3(\mathbb{C}),\ x,y\in\mathbb{C}^3
\end{eqnarray}
of which we shall make consistent use (on the standard basis it
states that entries of $(\det A)(A^T)^{-1}$ are co-factors of $A$;
the formula remains valid for $A\in\mathbf{M}_3(\mathbb{C})$, if
one replaces $(\det A)(A^T)^{-1}$ with the adjugate of $A$ or for
the cross product of $n-1$ vectors in $\mathbb{C}^n$).

Under the rotation $R\in\mathbf{O}_3(\mathbb{C}),\ \det R=1$, the
rulings of the unit sphere $X$ transform as:
$(u,v)\mapsto(\frac{e_3^TRY(u)}{f_1^TRY(u)},\frac{e_3^TRY(v)}{f_1^TRY(v)})=
(\frac{au+b}{cu+d},\frac{av+b}{cv+d})$, where
$R=\frac{d^2}{2(ad-bc)}Y(\frac{b}{d})f_1^T-\frac{c^2}{ad-bc}Y(\frac{a}{c})\bar
f_1^T+
X(\frac{a}{c},\frac{b}{d})e_3^T=\frac{a^2}{2(ad-bc)}f_1Y(-\frac{b}{a})^T-
\frac{c^2}{ad-bc}\bar
f_1Y(-\frac{d}{c})^T-e_3X(\frac{b}{a},\frac{d}{c})^T$; if $\det
R=-1$ compose with $\si_2$. For $Z\ R$ must satisfy $Re_3=e_3$; if
$\det R=1$, then $R=e^{ic}f_1\bar f_1^T+e^{-ic}\bar
f_1f_1^T+e_3e_3^T$ and $(u,v)\mapsto (e^{ic}u,e^{-ic}v)$;
otherwise compose with $\si_2$. We have the unit pseudo-sphere
$|iX|^2=-1,\ iX\subset\mathbb{R}^2\times
i\mathbb{R}\Leftrightarrow X\subset
(i\mathbb{R})^2\times\mathbb{R}\Leftrightarrow u\bar v=2;\
|d(iX)|^2=\frac{8dv\odot d\bar v}{(2-|v|^2)^2}$ and the real unit
sphere $X\subset\mathbb{R}^3,\ |X|^2=1 \Leftrightarrow u\bar
v=-2;\ |dX|^2=\frac{8dv\odot d\bar v}{(2+|v|^2)^2}$. The usual
linear element of the (pseudo-)sphere $(\frac{4dv\odot d\bar
v}{(1-|v|^2)^2})\ \frac{4dv\odot d\bar v}{(1+|v|^2)^2}$ is
obtained by replacing $v$ with $\sqrt2 v$, or equivalently for
$f_1:=e_1-ie_2,\ Y(v):=-v^2f_1+\bar f_1+2ve_3$, but we prefer to
normalize $f_1^T\bar f_1=1$ since the computations involving SJ
matrices would otherwise acquire as cumbersome constants powers of
$2$. The real equilateral paraboloid of revolution
$Z\subset\mathbb{R}^3$ is obtained for $u=\bar v$; the real
equilateral hyperbolic paraboloid $\mathrm{diag}[1\ \ i\ \ 1]
Z\subset\mathbb{R}^3$ is obtained for $u=\bar u,\ v=\bar v$.

With $\ep(z,a):=\frac{\sqrt{1-za^{-1}}(\sqrt{a^{-1}})^{-1}}
{\sqrt{a-z}}=\pm 1$ the confocal families of diagonal Q(W)C and of
the Darboux quadrics, the ruling families parametrization
invariant under the Ivory affinity, isotropic rulings, umbilics
and ruling parameters (if any) which give $x_z\cap^TC(\infty)$ can
be put as: QC:(\ref{eq:qc1})-(\ref{eq:qc3}),
QWC:(\ref{eq:qwc1})-(\ref{eq:qwc2}),
IQWC:(\ref{eq:iqwc1})-(\ref{eq:iqwc2}):
\begin{eqnarray}\label{eq:qc1}
x_z^TAR_z^{-1}x_z-1=0,\ A^{-1}:=\Sigma_{j=1}^3a_je_je_j^T;\
x_z(u,v):=\sqrt{R_z}(\sqrt A)^{-1}X(u,v);\ |\sqrt{R_z}(\sqrt A)^{-1}Y(v)|^2=\nonumber\\
0\Leftrightarrow\frac{v}{\sqrt{2}}=\pm\sqrt{\frac{a_1-a_3}{a_1-a_2}}\pm
i\sqrt{\frac{a_2-a_3}{a_2-a_1}}\ (=0,\infty\ \mathrm{if}\
a_1=a_2);\
\sqrt{R_z}(\sqrt A)^{-1}(\pm\sqrt{\frac{a_2-a_1}{a_2-a_3}}e_2\pm\nonumber\\
\sqrt{\frac{a_3-a_1}{a_3-a_2}}e_3),\ \sqrt{R_z}(\sqrt
A)^{-1}(\pm\sqrt{\frac{a_1-a_2}{a_1-a_3}}e_1
\pm\sqrt{\frac{a_3-a_2}{a_3-a_1}}e_3),\
\sqrt{R_z}(\sqrt A)^{-1}(\pm\sqrt{\frac{a_1-a_3}{a_1-a_2}}e_1\pm\nonumber\\
\sqrt{\frac{a_2-a_3}{a_2-a_1}}e_2)\ (\pm\sqrt{R_z}(\sqrt
A)^{-1}e_3\ if\ a_1=a_2);\ \mathrm{none}\ (x_z(0,0),\
x_z(\infty,\infty)\ \mathrm{if}\ a_1=a_2).
\end{eqnarray}
\begin{eqnarray}\label{eq:qc2}
x_z^TAR_z^{-1}x_z-1=0,\ A^{-1}:=(a_1I_2+J_2)\oplus a_2I_1;\
x_z(u,v):=\sqrt{R_z}(\sqrt A)^{-1}X(u,v)=\nonumber\\
(\frac{\sqrt{a_1-z}}{\ep(z,a_1)}(I_2+\frac{J_2}{2(a_1-z)})\oplus
\frac{\sqrt{a_2-z}}{\ep(z,a_2)}I_1)X(u,v);\
|\sqrt{R_z}(\sqrt A)^{-1}Y(v)|^2=0\Leftrightarrow\frac{v}{\sqrt{2}}=\infty,\nonumber\\
\pm\frac{1}{\sqrt{a_1-a_2}};\
\pm\sqrt{R_z}(\sqrt A)^{-1}(\frac{f_1}{\sqrt{a_1-a_2}}\pm e_3),\
\pm\sqrt{R_z}(\sqrt A)^{-1}\frac{f_1+2(a_1-a_2)\bar f_1}{2\sqrt{a_1-a_2}}\nonumber\\
(\mathrm{none}\ \mathrm{if}\ a_1=a_2);\ x_z(\infty,\infty).
\end{eqnarray}
\begin{eqnarray}\label{eq:qc3}
x_z^TAR_z^{-1}x_z-1=0,\ A^{-1}:=a_1I_3+J_3;\ x_z(u,v):=\sqrt{R_z}(\sqrt A)^{-1}X(u,v)
=\frac{\sqrt{a_1-z}}{\ep(z,a_1)}(I_3+\nonumber\\
\frac{J_3}{2(a_1-z)}-\frac{J_3^2}{8(a_1-z)^2})X(u,v);\
|\sqrt{R_z}(\sqrt A)^{-1}Y(v)|^2=0\Leftrightarrow v=\infty,\ 0;\
\pm\sqrt{R_z}(\sqrt A)^{-1}e_3;\nonumber\\ x_z(\infty,\infty).
\end{eqnarray}
\begin{eqnarray}\label{eq:qwc1}
x_z^T(AR_z^{-1}x_z-2e_3)+z=0,\ A^{\circ-1}:=a_1I_1\oplus a_2I_1;\
x_z(u,v):=\sqrt{R_z}LZ(u,v)+\frac{z}{2}e_3,\nonumber\\
L:=(\sqrt{A\oplus I_1})^{-1};\
|x_{zu}(u,\pm\sqrt{\frac{a_2-a_1}{2}})|^2=|x_{zv}(\pm\sqrt{\frac{a_2-a_1}{2}},v)|^2=0;\
\sqrt{R_z}L(\pm\sqrt{a_2-a_1}e_1+\nonumber\\\frac{a_2-a_1+z}{2}e_3),\
\sqrt{R_z}L(\pm\sqrt{a_1-a_2}e_2+\frac{a_1-a_2+z}{2}e_3);\
\mathrm{none}\ (x_z(0,\infty),\ x_z(\infty,0)\
\mathrm{if}\nonumber\\ a_1=a_2).
\end{eqnarray}
\begin{eqnarray}\label{eq:qwc2}
x_z^T(AR_z^{-1}x_z-2e_3)+z=0,\ A^{\circ -1}:=a_1I_2+J_2;\
x_z(u,v):=\sqrt{R_z}LZ(u,v)+\frac{z}{2}e_3,\nonumber\\
L:=(\sqrt{A\oplus I_1})^{-1},\ \sqrt{R_z}L=\frac{\sqrt{a_1-z}}{\ep(z,a_1)}(I_2+
\frac{J_2}{2(a_1-z)})\oplus I_1;\ |x_{zu}(u,0)|^2=\nonumber\\
|x_{zv}(\pm i,v)|^2=0;\ \pm i\sqrt{a_1-z}f_1+\frac{z}{2}e_3;\
x_z(\infty,0).
\end{eqnarray}
\begin{eqnarray}\label{eq:iqwc1}
(x_z-z\bar f_1)^T(AR_z^{-1}(x_z-z\bar f_1)-2\bar f_1)=0,\
A^{\circ -1}:=\bar J_2\oplus a_1I_1;\ x_z(u,v):=\nonumber\\
L\sqrt{I_3-zL^TA^2L}Z(u_z,v_z)+z\bar f_1,\
L:=f_1e_3^T+\bar f_1e_1^T-\sqrt{a_1}e_3e_2^T,\nonumber\\
u(z)=v(z):=-\frac{z}{2\sqrt{2}};\
|x_{zu}(u,\infty)|^2=|x_{zv}(\infty,v)|^2
=|x_{zu}(u,\frac{a_1}{2\sqrt{2}})|^2=\nonumber\\
|x_{zv}(\frac{a_1}{2\sqrt{2}},v)|^2=0;\
\frac{(a_1-z)^2}{8}f_1+\frac{a_1+z}{2}\bar f_1;\
x_z(\infty,\infty).
\end{eqnarray}
\begin{eqnarray}\label{eq:iqwc2}
(x_z-z\bar f_1)^T(AR_z^{-1}(x_z-z\bar f_1)-2\bar f_1)=0,\
A^{\circ -1}:=\bar J_3;\ x_z(u,v):=\nonumber\\
L\sqrt{I_3-zL^TA^2L}Z(u_z,v_z)+z\bar f_1,\
L:=J_3+\bar J_3^2,\ u(z):=-\frac{z}{2},\nonumber\\
v(z):=-\frac{z^2}{8};\ |x_{zu}(u,\infty)|^2=|x_{zv}(\infty,v)|^2
=|x_{zu}(u,0)|^2=0;\ \mathrm{none};\ x_z(\infty,\infty).
\end{eqnarray}

Let $f:=e_3\ (:=\bar f_1)$ respectively for (I)QWC, $u(z)=v(z)=0$
for QWC, $u_z:=u+u(z),\ v_z:=v+v(z)$. Note:
\begin{eqnarray}\label{eq:no}
L^{-1}=L^TA+e_3f^T,\ L^Tf=e_3,\ L^{-1}\sqrt{R_z}L=\sqrt{R'_z}+e_3(u(z)f_1+v(z)\bar f_1)^T,\nonumber\\
A':=L^TA^2L=A\ \mathrm{for\ QWC},\ =a_1^{-1}e_2e_2^T\ \mathrm{for\
IQWC}\ (\ref{eq:iqwc1}),\
=\bar J_3^2\ \mathrm{for\ IQWC}\ (\ref{eq:iqwc2});\nonumber\\
(L^TL)^{-1}e_3=e_3\ \mathrm{for\ QWC},\ =e_1\ \mathrm{for\ IQWC\
(\ref{eq:iqwc1})},\
=f_1\ \mathrm{for\ IQWC\ (\ref{eq:iqwc2})};\nonumber\\
L^{-1}x_z(u,v)=\sqrt{R'_z}Z(u_z,v_z)+z(1-\frac{|f|^2}{2})(L^TL)^{-1}e_3,\
L^T\hat N_0=I_{1,2}Z-e_3.
\end{eqnarray}
Note that the last relation of (\ref{eq:XY}) has a metric
projective symmetric equivalent for any quadric:
\begin{eqnarray}\label{eq:basic}
x_{zu}x_{zv}^T+x_{zv}x_{zu}^T-\frac{x_{zuv}x_z^T+x_zx_{zuv}^T}{2}
=-\frac{2}{\mathcal{B}}(zI_3-A^{\circ -1})\ \mathrm{for\ QC},\nonumber\\
x_{zu}x_{zv}^T+x_{zv}x_{zu}^T-(x_{zuv}x_z^T+x_zx_{zuv}^T)
=-\frac{2}{\mathcal{B}}(zI_3-A^{\circ -1})\ \mathrm{for\ (I)QWC},\nonumber\\
\mathcal{B}:=(u-v)^2\ \mathrm{for\ QC},\ :=2\ \mathrm{for\
(I)QWC}.
\end{eqnarray}
For QC it is just the last relation of (\ref{eq:XY}) after using
the next to last relation of (\ref{eq:XY}) and multiplying on both
left and right with $(\sqrt{R_z})^{-1}\sqrt{A}$. For (I)QWC after
multiplying on the left with $L^{-1}$ and on the right with
$(L^T)^{-1}$ it becomes $(\sqrt{R'_z}f_1+v_ze_3)(\sqrt{R'_z}\bar
f_1+u_ze_3)^T+(\sqrt{R'_z}\bar
f_1+u_ze_3)(\sqrt{R'_z}f_1+v_ze_3)^T-e_3(u_z\sqrt{R'_z}f_1+v_z\sqrt{R'_z}\bar
f_1+u_zv_ze_3+z(1-\frac{|f|^2}{2})(L^TL)^{-1}e_3)^T-(u_z\sqrt{R'_z}f_1+v_z\sqrt{R'_z}\bar
f_1+u_zv_ze_3+z(1-\frac{|f|^2}{2})(L^TL)^{-1}e_3)e_3^T=-z(L^TL)^{-1}+L^{-1}A^{\circ
-1}(L^T)^{-1}$; this boils down to
$-L^TA^2L=(1-\frac{|f|^2}{2})((L^TL)^{-1}e_3e_3^T+e_3e_3^T(L^TL)^{-1})-(L^TL)^{-1}
(=-I_{1,2}(L^TL)^{-1}I_{1,2}),\ I_{1,2}=L^{-1}A^{\circ
-1}(L^T)^{-1}$ which are straightforward. Note that
(\ref{eq:basic}) admits easy generalizations to quadrics in
$\mathbb{C}^n$; for the (I)QWC case the relevant parametrization
and necessary identities have already been introduced; for the QC
case the only important issue is finding the equivalent of the
last relation of (\ref{eq:XY}). Again just as for (I)QWC a
parametrization of the unit sphere
$X_{2n+1}\subset\mathbb{C}^{2n+1}$ comes first, following that
$X_{2n}:=X_{2n+1}\cap\{x\in\mathbb{C}^{2n+1}|e_{2n+1}^Tx=0\}$.
With $X(u,v)*:\mathbb{C}^{2n-1}\rightarrow\mathbb{C}^{2n+1},\
X*V:=[\frac{2-uv}{\sqrt{2}(u-v)}\ i\frac{2+uv}{\sqrt{2}(u-v)}\
\frac{u+v}{u-v}V^T]^T$ we have by induction
$X_{2n+1}(u^1,v^1,...,u^n,v^n)=X(u^n,v^n)*X_{2n-1}(u^1,v^1,...,u^{n-1},v^{n-1})$.

Note that all quadrics $\{x_z\}_{z\in\mathbb{C}}$ confocal to the
given one $x_0$ and all quadrics homothetic to $x_0$ are of the
$x_0$ type. Since modulo rigid motions there is a single IQWC
(\ref{eq:iqwc2}), all quadrics $x_z$ confocal to the IQWC
(\ref{eq:iqwc2}) $x_0$ or homothetic to $x_0$ are rigidly
applicable to $x_0$. This can be confirmed analytically as
follows: the rigid motion $(I_3-\frac{z}{2}(f_1e_3^T-e_3f_1^T)
-\frac{z^2}{8}e_3e_3^T,\frac{z^3}{16}f_1-\frac{3z}{2}\bar
f_1-\frac{3z^2}{8}e_3)
\in\mathbf{O}_3(\mathbb{C})\ltimes\mathbb{C}^3$ brings $x_z(u,v)$
to the canonical form $x_0(u-z,v)$; also the rigid motion
$(R(c),0)\in\mathbf{O}_3(\mathbb{C})\ltimes\mathbb{C}^3$ brings
$e^{2ic}x_0(u,v)$ to the canonical form $x_0(e^{ic}u,e^{2ic}v)$.

$x_z\cap^TC(\infty)$ can be seen after a homography (homographies
preserve tangency and its order); we consider only $z=0$. Under
the homography $\mathbb{C}^3\ni x:=[x^1\ x^2\ x^3\
1]\mapsto[\frac{f_1^Tx}{\bar f_1^Tx}\ \frac{e_3^Tx}{\bar f_1^Tx}\
\frac{1}{\bar f_1^Tx}\ 1],\ C(\infty)$ becomes $c(t):=-2[t^2\ t\
0]$ and $[f_1^T,\ 0]$ becomes $0$. The curves $(u^{-1},v):=
(O(t^2),-\sqrt{a_1}t+O(t^2)),\
(u^{-1},v^{-1}):=(-i\sqrt{\frac{2}{a_1}}t+O(t^2),
i\sqrt{\frac{2}{a_1}}t+O(t^2))$ respectively on the transformed
(\ref{eq:qwc2}), (\ref{eq:iqwc1}) meet $c(t)$ tangentially at
$t=0$; the curves $(u^{-1},v^{-1}):=
(t-\frac{3}{a_1}t^2+O(t^3),t-\frac{3}{a_1}t^2+O(t^3)),\
(u^{-1},v^{-1}):=(-2t+O(t^3),-2t^2+O(t^3))$ respectively on the
transformed (\ref{eq:qc3}), (\ref{eq:iqwc2}) osculate $c(t)$ at
$t=0$ up to order $2$; the curves
$(u^{-1},v^{-1}):=(t-\frac{t^3}{a_1}+O(t^4),
t-\frac{t^3}{a_1}+O(t^4))$ on the transformed (\ref{eq:qc2}) for
$a_2=a_1$ osculate $c(t)$ at $t=0$ up to order $3$.

According to Eisenhart (\cite{E1},\S\ 151), Goursat integrated the
equations for the deformations of the paraboloids (\ref{eq:qwc2})
and (\ref{eq:iqwc2}); also Darboux reduced in \cite{D2} the
deformations of the Darboux quadrics to the deformations of the
sphere.

For QC (\ref{eq:qc1}) if $a_j,\ j=1,2,3$ are distinct, then $x_z$
meets $C(\infty),\ \mathcal{S}(x_{a_j}),\ j=1,2,3$ only
transversally; if $a_1=a_2\neq a_3$, then $x_z$ meets $C(\infty),\
\mathcal{S}(x_{a_3})$ only tangentially at $[f_1^T,0],\ [\bar
f_1^T,0]$ and $\mathcal{S}(x_{a_1})$ only transversally at
$\pm\sqrt{a_3-z}e_3$. For QC (\ref{eq:qc2}) if $a_1\neq a_2$, then
$x_z$ meets $C(\infty),\ \mathcal{S}(x_{a_2})$ tangentially at
$[f_1^T,0]$ and transversally at two more points each; further
$x_z$ meets $\mathcal{S}(x_{a_1})$ only transversally; if
$a_1=a_2$, then $x_z$ meets $C(\infty),\ \mathcal{S}(x_{a_1})$
only tangentially at $[f_1^T,0]$. For QC (\ref{eq:qc3}) $x_z$
meets $C(\infty),\ \mathcal{S}(x_{a_1})$ tangentially at
$[f_1^T,0]$ and transversally at one respectively two more points.
For QWC (\ref{eq:qwc1}) if $a_1\neq a_2$, then $x_z$ meets
$C(\infty)$ only transversally; further it meets
$\mathcal{S}(x_{a_j}),\ j=1,2$ tangentially at $[e_3^T,0]$ and
transversally at two more points each; if $a_1=a_2$, then $x_z$
meets $C(\infty)$ only tangentially at $[f_1^T,0],\ [\bar
f_1^T,0]$ and $\mathcal{S}(x_{a_1})$ only transversally at
$\frac{z}{2}e_3,\ [e_3^T,0]$. For QWC (\ref{eq:qwc2}) $x_z$ meets
$C(\infty),\ \mathcal{S}(x_{a_1})$ tangentially at $[f_1^T,0]$,
respectively $[e_3^T,0]$ and transversally at two more points
each. For QC (\ref{eq:iqwc1}) $x_z$ meets $C(\infty),\
\mathcal{S}(x_{a_1})$ tangentially at $[f_1^T,0]$ and
transversally at two respectively one more points. For QC
(\ref{eq:iqwc2}) $x_z$ meets $C(\infty)$ tangentially at
$[f_1^T,0]$ and transversally at one more point. The points of
intersection with finite singular quadrics are situated in
$\mathbb{CP}^2$ iff they are tangency points.

For $n=3$ Bianchi has a beautiful geometric proof of Bianchi II
(see also Darboux (\cite{D1},\S\ 603)): all coordinates on a
surface whose definition involves only projective invariants
(tangency of different orders; for example tangent planes of a
surface, osculating planes of a curve) are preserved by
homographies. In particular asymptotes and conjugate systems are
preserved; since on developables the asymptotes coincide,
homographies take developables to developables. There is a
developable circumscribed to two given arbitrary curves and the
isotropic developable with generating curve $c(u)$ is the
developable $c(u)+vY(f(u))$ circumscribed (tangent) to $c,\
C(\infty)$, so $c'(u)^TY(f(u))=0$; thus it is an algebraic surface
of degree four if $c$ is a conic. Now nonsingular quadrics $x_z$
intersect $\mathcal{S}(x_{a_1})$ at four points (B\'{e}zout) and
the isotropic rulings of the developable through those points are
all the isotropic rulings on $x_z$; by inspection all quadrics are
uniquely determined by their isotropic rulings (most quadrics have
three distinct isotropic rulings of one ruling family; as a
quadric is uniquely determined by three rulings of a ruling family
the statement follows in these cases). As a quadric varies in its
confocal family its umbilics describe the conics which are the
singular sets of singular quadrics of the family and its isotropic
rulings isotropic developables generated by these conics. Thus $H$
takes the isotropic developable of the conic
$\mathcal{S}(x_{a_1})$ to the isotropic developable of the conic
$\mathcal{S}(\ti x_{\ti a_1}):=H(C(\infty))$ and $\ti x_{\ti z}$
is uniquely recovered. For example for QC (\ref{eq:qc1}) if
$\frac{v}{\sqrt{2}}=\ep_1
(\sqrt{\frac{a_1-a_3}{a_1-a_2}}+i\ep_2\sqrt{\frac{a_2-a_3}{a_2-a_1}}),\
c(z):=\sqrt{R_z}(\sqrt
A)^{-1}\ep_3(\sqrt{\frac{a_1-a_3}{a_1-a_2}}e_1+\
\ep_4\sqrt{\frac{a_2-a_3}{a_2-a_1}}e_2),\ \ep_j=\pm 1,\
j=1,...,4$, then $c'(z)^T\sqrt{R_z}(\sqrt A)^{-1}Y(v)
=i\ep_3\sqrt{2}\sqrt{\frac{a_1-a_3}{a_1-a_2}}\sqrt{\frac{a_2-a_3}{a_2-a_1}}
(\sqrt{\frac{a_1-a_3}{a_1-a_2}}+i\ep_2\sqrt{\frac{a_2-a_3}{a_2-a_1}})(\ep_2-\ep_4)=0$
for $\ep_2=\ep_4$.

Let $H\in\mathbf{PGL}_4(\mathbb{C})$ be a homography of $\mathbb{C}^3\ni
x\simeq [x^T,1]\in\mathbb{CP}^3$, with $\ti x:=\frac{I_{1,3}Hx}{e_4^THx}$.
We shall always have $\ti x_{\ti z}(\ti u,\ti v):=H(x_z(u,v)),\
H:=\begin{bmatrix}H'&h_1\\h_2^T&0\end{bmatrix},\ \ti H:=H^{-1}=\begin{bmatrix}\ti H'&\ti h_1\\
\ti h_2^T&0\end{bmatrix}$; the rulings $\ti u,\ \ti v$ may change
with $z$ (this corresponds to symmetries of $\ti x_{\ti z}$ wrt
the principal planes).

If $H'=\Sigma_{j=2}^3\frac{1}{\sqrt{a_1-a_j}}e_je_j^T,\ \ti
H'=\Sigma_{j=2}^3\sqrt{a_1-a_j}e_je_j^T,\ h_1=h_2=\ti h_1=\ti
h_2=e_1$, then $x_z(u,v),\\ \ti x_{\ti z}(\ti u,\ti v)$ are QC of
type (\ref{eq:qc1}), $\ti a_1:=a_1^{-1},\ \ti a_j:=\ti
a_1+(a_j-a_1)^{-1},\ j=2,3,\ \ti z:=\ti a_1+(z-a_1)^{-1}$; $H$
takes $x_z\cap\{e_1^Tx=0\}= \{\sqrt{R_z}(\sqrt
A)^{-1}X(2v^{-1},v)\}$ to $\ti x_{\ti
z}\cap\mathbb{CP}^2=\{[(\sqrt{\ti R_{\ti z}}(\sqrt{\ti
A})^{-1}Y(\ti v))^T,0]\}$, so $(\ti u,\ti
v)=(i\ep_3\frac{2-i\ep_2\sqrt{2}u}{\sqrt{2}+i\ep_2u},
-i\ep_3\frac{2+i\ep_2\sqrt{2}v}{\sqrt{2}-i\ep_2v}),\
\ep_j:=\frac{i\ep(\ti z,\ti a_1)\ep(z,a_j) \sqrt{a_j-z}\sqrt{\ti
a_1-\ti z}}{\ep(\ti z,\ti a_j)\sqrt{a_1-a_j}\sqrt{\ti a_j-\ti
z}},\ j=2,3$. If $a_2\neq a_3$, then $H(x_z\pitchfork
C(\infty))=\ti x_{\ti z}\pitchfork\mathcal{S}(\ti x_{\ti a_1}),\
H(x_z\pitchfork\mathcal{S}(x_{a_1}))=\ti x_{\ti z}\pitchfork
C(\infty),\ H(x_z\pitchfork\mathcal{S}(x_{a_j}))=\ti x_{\ti
z}\pitchfork\mathcal{S}(\ti x_{\ti a_j}),\ j=2,3$; if $a_2=a_3$,
then $H(x_z\cap^TC(\infty))=\ti x_{\ti z}\cap^TC(\infty),\
H(x_z\pitchfork\mathcal{S}(x_{a_2}))=\ti x_{\ti
z}\pitchfork\mathcal{S}(\ti x_{\ti a_2})$. Under the homothety
$\frac{1}{\sqrt{(a_1-a_2)(a_1-a_3)}}$ and $x^2\leftrightarrow x^3$
in the range of $H$, $\ti a_j$ become $a_j-w,\ j=1,3,\
w:=a_1-\frac{(a_1-a_2)(a_1-a_3)}{a_1}$, so the new $H$ takes the
family of confocal QC (\ref{eq:qc1}) into itself.

If $H'=i\sqrt{\ti a_2-\ti a_1}((a_2-a_1)f_1\bar
f_1^T+\frac{1}{2}f_1f_1^T- \frac{1}{a_2-a_1}\bar f_1f_1^T),\ \ti
H'=\frac{-i}{\sqrt{\ti a_2-\ti a_1}} (\frac{1}{a_2-a_1}f_1\bar
f_1^T+\frac{1}{2}f_1f_1^T-(a_2-a_1)\bar f_1f_1^T),\ h_1=h_2=\ti
h_1=\ti h_2=e_3$, then $x_z(u,v),\ \ti x_{\ti z}(\ti u,\ti v)$ are
QC of the type (\ref{eq:qc2}), $\ti a_2:=a_2^{-1},\ \ti a_1:=\ti
a_2+(a_1-a_2)^{-1},\ \ti z:=\ti a_2+(z-a_2)^{-1};\ H$ takes the
cone $x_z\cap\{f_1^Tx=0\}=\{\sqrt{R_z}(\sqrt A)^{-1}X(\infty,v)\}
\cup\{\sqrt{R_z}(\sqrt A)^{-1}X(u,\infty)\}$ (with vertex
$x_z\cap^TC(\infty))$) to the cone $\ti x_{\ti
z}\cap\{f_1^Tx=0\}=\{\sqrt{\ti R_{\ti z}}(\sqrt{\ti
A})^{-1}X(\infty,\ti v)\} \cup\{\sqrt{\ti R_{\ti z}}(\sqrt{\ti
A})^{-1}X(\ti u,\infty)\}$ with same vertex, so $(\ti u,\ti
v)=(-\ep_1u,\ep_1v),\ \ep_1:=\frac{i\ep(z,a_1)\ep(z,a_2)\ep(\ti
z,\ti a_1)\sqrt{a_1-z}} {\sqrt{\ti a_2-\ti
a_1}\sqrt{a_2-z}\sqrt{\ti a_1-\ti z}}$ if
$\frac{\sqrt{a_2-z}\sqrt{\ti a_2-\ti z}}{\ep(z,a_2)\ep(\ti z,\ti
a_2)}=1$ and composed with $\si_2$ otherwise; also
$H(x_{z}\cap^TC(\infty))= \ti x_{\ti z}\cap^T\mathcal{S}(\ti
x_{\ti a_2}),\ H(x_{z}\pitchfork C(\infty))= \ti x_{\ti
z}\pitchfork\mathcal{S}(\ti x_{\ti a_2}),\
H(x_z\cap^T\mathcal{S}(x_{a_2}))=\ti x_{\ti z}\cap^TC(\infty),\
H(x_z\pitchfork\mathcal{S}(x_{a_2}))= \ti x_{\ti z}\pitchfork
C(\infty),\ H(x_z\pitchfork\mathcal{S}(x_{a_1}))= \ti x_{\ti
z}\pitchfork\mathcal{S}(\ti x_{\ti a_1})$.

If $H'=ie_3(\frac{\ti a_1}{2}f_1-\bar f_1)^T+\sqrt{\ti a_1-\ti
a_2}e_2e_3^T,\ \ti H'=if_1e_3^T+\frac{1}{\sqrt{\ti a_1-\ti
a_2}}e_3e_2^T,\ h_1=\ti h_2=e_1,\ h_2=if_1,\ \ti h_1=-i(\bar
f_1+\frac{a_1}{2}f_1)$, then $x_z(u,v)$ are QC of type
(\ref{eq:qc2}), $\ti x_{\ti z}(\ti u,\ti v)$ are QWC of type
(\ref{eq:qwc1}): $\ti a_1:=a_1^{-1},\ \ti a_2:=\ti
a_1+(a_2-a_1)^{-1},\ \ti z:= \ti a_1+(z-a_1)^{-1};\ H$ takes the
cone $x_z\cap\{f_1^Tx=0\} =\{\sqrt{R_z}(\sqrt
A)^{-1}X(\infty,v)\}\cup\{\sqrt{R_z}(\sqrt A)^{-1}X(u,\infty)\}$
(with vertex $x_z\cap^TC(\infty)\notin\mathcal{S}(x_{a_1})$) to
the cone $\ti x_{\ti z}\cap \mathbb{CP}^2=\{[\ti x_{\ti
z}(\infty,\ti v)^T,0]\}\cup \{[\ti x_{\ti z}(\ti
u,\infty)^T,0]\}=\{[(\sqrt{\ti R_{\ti z}}L\bar f_1+2\ti
ve_3)^T,0]\}\cup \{[(\sqrt{\ti R_{\ti z}}L\bar f_1+2\ti
ue_3)^T,0]\}$ with vertex $[e_3^T,0]\notin C(\infty)$, so $(\ti
u,\ti v)=(\frac{-i\ep_1}{\sqrt{2}}u, \frac{i\ep_1}{\sqrt{2}}v),\
\ep_1:=\frac{\sqrt{a_1-z}\sqrt{\ti a_1-\ti z}} {\ep(z,a_1)\ep(\ti
z,\ti a_1)}$ if $\frac{-i\ep(z,a_2)\ep(\ti z,\ti a_1)\sqrt{\ti
a_2-\ti z}} {\sqrt{\ti a_1-\ti a_2}\sqrt{\ti a_1-\ti
z}\sqrt{a_2-z}}=\ep(\ti z,\ti a_2)$ and composed with $\si_2$
otherwise; also $H(x_{z}\cap^TC(\infty))=\ti x_{\ti
z}\cap^T\mathcal{S}(\ti x_{\ti a_1}),\ H(x_{z}\pitchfork
C(\infty))=\ti x_{\ti z}\pitchfork\mathcal{S}(\ti x_{\ti a_1}),\
H(x_z\pitchfork\mathcal{S}(x_{a_1}))=\ti x_{\ti z}\pitchfork
C(\infty),\ H(x_z\cap^T\mathcal{S}(x_{a_2}))=\ti x_{\ti
z}\cap^T\mathcal{S}(\ti x_{\ti a_2}),\
H(x_z\pitchfork\mathcal{S}(x_{a_2}))=\ti x_{\ti
z}\pitchfork\mathcal{S}(\ti x_{\ti a_2})$.

If $H'=i(f_1\bar f_1^T+\bar f_1(-\ti a_1f_1+e_3)^T),\ \ti
H'=-i(f_1\bar f_1^T+e_3f_1^T),\ h_1=\ti h_2=e_3,\ \ti h_1=i(\bar
f_1+\ti a_1e_3),\ h_2=-if_1$, then $x_z(u,v)$ are QC of type
(\ref{eq:qc3}), $\ti x_{\ti z}(\ti u,\ti v)$ are IQWC of type
(\ref{eq:iqwc1}): $\ti a_1:=a_1^{-1},\ \ti z:=\ti
a_1+(z-a_1)^{-1};\ H$ takes the cone
$x_z\cap\{f_1^Tx=0\}=\{\sqrt{R_z}(\sqrt A)^{-1}X(\infty,v)\}
\cup\{\sqrt{R_z}(\sqrt A)^{-1}X(u,\infty)\}$ with vertex
$x_z\cap^TC(\infty)\in\mathcal{S}(x_{a_1})$ to the cone $\ti
x_{\ti z}\cap\mathbb{CP}^2=\{[\ti x_{\ti z}(\infty,\ti
v)^T,0]\}\cup\{[\ti x_{\ti z}(\ti u,\infty)^T,0]\}=\{[(\bar
f_1+(\sqrt{2}\ti v-\frac{\ti z}{2})f_1+ i\frac{\sqrt{\ti a_1-\ti
z}}{\ep(\ti z,\ti a_1)}e_3)^T,0]\} \cup\{[(\bar f_1+(\sqrt{2}\ti
u-\frac{\ti z}{2})f_1 -i\frac{\sqrt{\ti a_1-\ti z}}{\ep(\ti z,\ti
a_1)}e_3)^T,0]\}$ with vertex $\ti x_{\ti z}\cap^TC(\infty)\in
C(\infty)$, so $(\ti u,\ti v)=(\frac{1}{\sqrt{2}}(\frac{\ti
a_1}{2}-u), \frac{1}{\sqrt{2}}(\frac{\ti a_1}{2}-v))$ if
$\frac{\sqrt{a_1-z}\sqrt{\ti a_1-\ti z}} {\ep(z,a_1)\ep(\ti z,\ti
a_1)}=1$ or composed with $\si_2$ otherwise; also
$H(x_{z}\cap^TC(\infty))=\ti x_{\ti z}\cap^T\mathcal{S}(\ti x_{\ti
a_1}),\ H(x_{z}\pitchfork C(\infty))=\ti x_{\ti
z}\pitchfork\mathcal{S}(\ti x_{\ti a_1}),\
H(x_z\cap^T\mathcal{S}(x_{a_1}))=\ti x_{\ti z}\cap^TC(\infty),\
H(x_z\pitchfork\mathcal{S}(x_{a_1}))=\ti x_{\ti z}\pitchfork
C(\infty)$.

If $H'=i(f_1e_3^T+e_3(\frac{\ti a_1}{2}f_1-\bar f_1)^T),\ \ti
h_1=-i(\frac{\ti a_1}{2}f_1+\bar f_1),\ h_2=if_1$ and for $\ti
H',\ h_1,\ \ti h_2$ replace $\ti a_1$ with $a_1$, then $x_z(u,v),\
\ti x_{\ti z}(\ti u,\ti v)$ are QWC of type (\ref{eq:qwc2}), $\ti
a_1:=a_1^{-1},\ \ti z:=\ti a_1+(z-a_1)^{-1};\ H$ takes the cone
$x_z\cap\{f_1^Tx=0\}=
\{\frac{\sqrt{a_1-z}}{\ep(z,a_1)}uf_1+\frac{z}{2}e_3\}\cup\{[(\frac{\sqrt{a_1-z}}{\ep(z,a_1)}f_1
+ve_3)^T,0]\}$ with vertex
$x_z\cap^TC(\infty)\notin\mathcal{S}(x_{a_1})$ to  the cone $\ti
x_{\ti z}\cap \mathbb{CP}^2=\{[(\frac{\sqrt{\ti a_1-\ti
z}}{\ep(\ti z,\ti a_1)}(\bar f_1+ \frac{f_1}{2(\ti a_1-\ti
z)})+\ti ue_3)^T,0]\}\cup\{[(\frac{\sqrt{\ti a_1-\ti z}} {\ep(\ti
z,\ti a_1)}f_1+\ti ve_3)^T,0]\}$ with vertex $[e_3^T,0]\notin
C(\infty)$, so $(\ti u,\ti v)=(\ep_1u,\frac{-\ep_1}{v}),\
\ep_1:=\frac{\sqrt{a_1-z}\sqrt{\ti a_1-\ti z}} {\ep(z,a_1)\ep(\ti
z,\ti a_1)}$; also $H(x_{z}\cap^TC(\infty))=\ti x_{\ti z}\cap^T
\mathcal{S}(\ti x_{\ti a_1}),\ H(x_{z}\pitchfork C(\infty))=\ti
x_{\ti z}\pitchfork \mathcal{S}(\ti x_{\ti a_1}),\
H(x_z\cap^T\mathcal{S}(x_{a_1}))=\ti x_{\ti z}\cap^TC(\infty),\
H(x_z\pitchfork\mathcal{S}(x_{a_1}))=\ti x_{\ti z}\pitchfork
C(\infty)$.

Homographies with $H'=-\bar J_2\oplus e^aI_1,\ \ti H'=-J_2\oplus
e^{-a}I_1,\ e^{-c}h_1=h_2=f_1,\ \ti h_1=e^c\ti h_2=\bar f_1$ or
$H=\begin{bmatrix}H'&0\\0^T&1\end{bmatrix},\ e^cI_3\neq
H'\in\mathbf{GL}_3(\mathbb{C})$ SJ do not have the property of
taking a family of confocal quadrics to another one and up to
rigid motions and homotheties in both the range and domain of $H$
all homographies $H\in\mathbf{PGL}_4(\mathbb{C})$ can be brought
to one of the cases above or to the identity. Quadrics whose
confocal family cannot be taken to another one by homographies are
the (pseudo-)spheres, QC (\ref{eq:qc2}) and QWC (\ref{eq:qwc1})
for $a_1=a_2$, IQWC (\ref{eq:iqwc2}).

\subsection{Totally real confocal quadrics in $\mathbb{C}^3$ and the corresponding
homographies of Bianchi II}\label{subsec:algpre22} \noindent

\noindent With $\mathcal{A}:=i\det(\sqrt{A})^{-1}$ for QC,
$:=i\det L$ for (I)QWC ($=\sqrt{a_1}$ for IQWC (\ref{eq:iqwc1}),
$=i$ for IQWC (\ref{eq:iqwc2})) we have the Gau\ss\ curvature
$K(x_0)=\frac{-1}{\mathcal{A}^2|\hat N_0|^4}$ of $x_0$.

We are interested in surfaces $x_0^0$ in the quadric $x_0$ having
real linear element (that is $|dx_0^0|^2=|d\bar x_0^0|^2$). By
Sylvester's such surfaces must have real dimension $2$, so they
are obtained by imposing two functionally independent real
relations between the four real parameters of $x_0$. We shall see
in \S\ \ref{subsec:rolling2} that the only such surfaces in
spheres are the totally real (pseudo-)spheres. The totally real
unit spheres $X(u,v)$ are obtained as follows:
$X\subset\mathbb{R}^3$ has signature $(2,0)$ for $v=-2\bar
u^{-1},\ u\in\mathbb{C}$;
$X\subset(i\mathbb{R})^2\times\mathbb{R}$ has signature $(0,2)$
for $v=2\bar u^{-1},\ u\in\mathbb{C}$; $X\subset\mathbb{R}^2\times
i\mathbb{R}$ has signature $(1,1)$ for
$(u,v)=(\sqrt{2}e^{is},\sqrt{2}e^{it}),\ s,t\in\mathbb{R}$; the
totally real unit pseudo-spheres $iX(u,v)$ are obtained through
multiplication by $i$, when the components of both the signature
of $X$ and of the ambient totally real space are exchanged.

For the remaining quadrics $K(x_0^0)=\overline{K(x_0^0)}$, or
equivalently

$$(I)\ \ \mathcal{A}|\hat N_0^0|^2=\ep\overline{\mathcal{A}|\hat N_0^0|^2},\ \ep:=\pm 1$$

provides a non-vacuous real relation between the four real
parameters of $x_0$. Keeping account of $\hat N_0^0=Ax_0^0+B$ and
applying $d$ to (I) we get $(dx_0^0)^T(\mathcal{A}A\hat
N_0^0)=\ep\overline{(dx_0^0)^T(\mathcal{A}A\hat N_0^0)}$. Now
Sylvester's allows us to find an adapted moving frame
$R[\sqrt{\ep_1}e_1\ \ \sqrt{\ep_2}e_2\ \ e_3],\ \ep_j:=\pm 1,\
j=1,2$ chosen so that $\mathrm{diag}[\ep_1\ \ \ep_2]$ is the
signature of $|dx_0^0|^2,\ R\subset\mathbf{O}_3(\mathbb{C}),\
R^{-1}dx_0^0=\sqrt{\ep_1}\psi_1e_1+\sqrt{\ep_2}\psi_2e_2$ and
$\psi_j,\ j=1,2$ are real (valued) independent $1$-forms. If
$R^{-1}\mathcal{A}A\hat
N_0^0=:\sqrt{\ep_1}v_1e_1+\sqrt{\ep_2}v_2e_2+v_3e_3,\
v_j\in\mathbb{C}$, then we get $\ep_1(v_1-\ep\bar
v_1)\psi_1+\ep_2(v_2-\ep\bar v_2)\psi_2=0$, so $v_j=\ep\bar v_j,\
v_j^2=\bar v_j^2,\ j=1,2$. This implies that the tangential
component $\frac{(\hat N_0^0\times(\mathcal{A}A\hat
N_0^0))\times\hat N_0^0} {|\hat N_0^0|^2}$ of the vector
$\mathcal{A}A\hat N_0^0$ has real length; using also (I) we
obtain:

$$(II)\ \ \mathcal{A}^3|V_2|^2=
\ep\overline{\mathcal{A}^3|V_2|^2},\ V_2:=\hat N_0^0\times A\hat N_0^0.$$

Similarly to the step (I) $\Rightarrow$ (II), one can apply $d$ to
(II) and get
$\mathcal{A}^7|V_3|^2=\ep\overline{\mathcal{A}^7|V_3|^2},\
V_3:=\hat N_0^0\times A(A(V_2\times\hat N_0^0)-V_2\times A\hat
N_0^0)$; in fact one can get a whole hierarchy of consequences,
but (I)$\&$(II) will be sufficient for our purposes, because they
are functionally independent for all quadrics which are not
spheres. We still need to satisfy

$$(III)\ \ |dx_0^0|^2-|d\bar x_0^0|^2=0.$$

Since we already have the ruling families parametrization $(u,v)$
on $x_0$, it is more convenient to work with the parametrization
$(u,v,\bar u,\bar v)$ or closely related ones $(U,V,\bar U,\bar
V)$, $U,\ V$ being simple rational functions of $u$ and $v$ with
$\frac{dU\wedge dV}{du\wedge dv}\neq 0$ (all functions below will
be rational). (I)$\&$(II) will be equivalent to two real relations
$f_1(V,\bar V)=f_2(U,\bar U)=0$ or will imply $f(U,V,\bar V)=0$
(and thus its conjugate). If one can solve for $U=F(V,\bar V)$,
then this and its conjugate will be functionally independent and
thus will imply (I)$\&$(II). However, in this case (III) further
imposes two real conditions, namely the condition that
$|dx_0^0|^2-|d\bar x_0^0|^2$ is degenerate as a symmetric bilinear
form in $dV,\ d\bar V$ and the condition that the coefficient of
$dV\odot d\bar V$ in $|dx_0^0|^2-|d\bar x_0^0|^2$ is $0$; these
conditions must not impose a real condition between $V$ and $\bar
V$. For most cases this program becomes tedious, so other
simplifications (usually considering only the hot's in $V,\ \bar
V$) may be required.

One can try to find another functional relationship between
$U,V,\bar U,\bar V$, which is easier to obtain than the third
functional relationship of the hierarchy described above or than
the method described above; hopefully the condition that
(I)$\&$(II) and this new functional relationship are functionally
dependent will provide other useful information. An example of
such a functional relationship is as follows: (I) gives one of the
a-priori independent variables $U,V,\bar U,\bar V$ (for example
$\bar U$) as a function of the other three variables; replacing
this into (III) further provides a functional relationship between
the remaining three variables, namely the condition that
$|dx_0^0|^2-|d\bar x_0^0|^2$ is degenerate as a symmetric bilinear
form in the differentials of the three remaining variables; this
condition is valid since the real $2$-dimensionality of $x_0^0$
implies that $|dx_0^0|^2-|d\bar x_0^0|^2$, as a symmetric bilinear
form in the differentials of the three remaining variables, admits
two linearly independent orthogonal isotropic vectors, namely
$\pa_s,\ \pa_t$ for $s,t$ real parametrization of $x_0^0$.
Although finding this functional relationship is more difficult
than finding the ones imposed by (I)$\&$(II), a nice feature
appears: this functional relationship always factors as a product
of $|\hat N_0^0|^2\neq 0$ and another term being $0$, because if
$|\hat N_0^0|^2=0$, then $x_0^0$ is a complex holomorphic curve
with linear element extended to isotropic one on $x_0^0$ (it is
real $2$-dimensional, but not real valued), so it has degenerate
real $3$-dimensional linear element. Unfortunately in all cases
except IQWC (\ref{eq:iqwc2}) this functional relationship is the
consequence $f(U,V,\bar V)=0$ of (I)$\&$(II) (for IQWC
(\ref{eq:iqwc2}) it is the consequence $f_1(u,\bar u)=0$ of
(I)$\&$(II) if we begin with $\bar v$ as a function of $u,\ \bar
u,\ v$), so it is useful only as a confirmation of the already
completed computations.

We shall prove, after a study case by case, that all surfaces in
quadrics and having real linear element are totally real, so most
quadrics do not contain such surfaces; a simpler unifying argument
using rolling will appear in \S\ 6.6.

For QC we have $\mathcal{A}=i\det(\sqrt{A})^{-1},\
x_0^0=(\sqrt{A})^{-1}X,\ \hat N_0^0=\sqrt{A}X,\  |\hat
N_0^0|^2=X^TAX,\ |V_2|^2=\det A(X\times AX)^TA^{-1}(X\times AX),\
|dx_0^0|^2=dX^TA^{-1}dX$.

For (I)QWC we have $\mathcal{A}=i\det L,\ x_0^0=LZ,\ \hat
N_0^0=(L^T)^{-1}(I_{1,2}Z-e_3),\ |\hat N_0^0|^2
=(I_{1,2}Z-e_3)^T(L^TL)^{-1}(I_{1,2}Z-e_3),\ |V_2|^2=(\det
L)^{-2}((I_{1,2}Z-e_3)\times
(A'Z-(L^TL)^{-1}e_3+|f|^2e_3))^T(L^TL)\\
((I_{1,2}Z-e_3)\times(A'Z-(L^TL)^{-1}e_3+|f|^2e_3)),\
|dx_0^0|^2=dZ^T(L^TL)dZ,\ L^TL=A^{\circ -1}+e_3e_3^T\ \mathrm{for\
QWC},\ =e_3e_1^T+e_1e_3^T+a_1e_2e_2^T\ \mathrm{for\ IQWC\
(\ref{eq:iqwc1})},\ =J_3^2+\bar J_3\ \mathrm{for\ IQWC\
(\ref{eq:iqwc2})},\ (L^TL)^{-1}=A+e_3e_3^T\ \mathrm{for\ QWC},\
=e_3e_1^T+e_1e_3^T+a_1^{-1}e_2e_2^T\ \mathrm{for\ IQWC\
(\ref{eq:iqwc1})},\ =\bar J_3^2+J_3\ \mathrm{for\ IQWC\
(\ref{eq:iqwc2})}$.

Since the most metrically degenerate a complex metric type of a
quadric is, the less totally real types it contains and also
keeping account of the homographies of Bianchi II, the discussion
will take place as follows: IQWC (\ref{eq:iqwc2}), IQWC
(\ref{eq:iqwc1}) $\&$ QC (\ref{eq:qc3}), QWC (\ref{eq:qwc2}), QWC
(\ref{eq:qwc1}) $\&$ QC (\ref{eq:qc2}), QC (\ref{eq:qc1}).

For $\ep,\ep_1:=\pm 1$ note $\frac{\sqrt{-\ep}}{\sqrt{\ep}}=i\ep,\
\frac{\sqrt{-\ep_1\sqrt{-\ep}}}{\sqrt{\ep_1\sqrt{-\ep}}}=-i\ep\ep_1$,
etc, of which we shall make consistent use below.

For IQWC (\ref{eq:iqwc2}) (I)$\&$(II) become $u^2-2v=-\ep(\bar
u^2-2\bar v),\ 3u^2+2v=-\ep(3\bar u^2+2\bar v)$, or $u^2=-\ep\bar
u^2,\ v=-\ep\bar v$. Keeping account of $|dx_0^0|^2=2du\odot
d(uv)+dv^2$, (III) will be satisfied.

Thus we have four quadrics $x_0^0$ in $x_0=x_0(u,v)=uvf_1+u\bar
f_1+ve_3$, given by
$(u,v)=(\sqrt{\ep_1\sqrt{-\ep}}s,\sqrt{-\ep}t),\ \ep_1:=\pm 1,\
s,t\in\mathbb{R}$. Since $\ep_1$ in $|dx_0^0|^2$ can be absorbed
by $(s,t)\mapsto(\ep_1s,\ep_1t)$, there are two types of real
linear element of IQWC({\ref{eq:iqwc2}) from a real point of view
(namely $-\ep(2ds\odot d(st)+dt^2)$); from a complex point of view
there is only one, since further $-\ep$ can be absorbed by
$(s,t)\mapsto(\sqrt[4]{-\ep}s,\sqrt{-\ep}t)$. All symmetries of
the linear element are realized by rigid motions composed when
necessary with a multiplication by $i$:
$iR(\frac{\pi}{4})x_0(u,v)=x_0(\sqrt{i}u,iv)$; $\ep_1$ appears in
the linear element of $x_0$ because of the rotation
$(iR(\frac{\pi}{4}))^2=-R(\frac{\pi}{2})$ which generates the
group of rigid motions preserving the canonical form of $x_0:\
-R(\frac{\pi}{2})x_0(u,v)=x_0(-iu,-v)$. Note that the IQWC
(\ref{eq:iqwc2}) $x_0$ is the only quadric such that $ix_0$
differs from $x_0$ by a rigid motion; for all other quadrics
$ix_0$ will be a quadric of the same type, but which is different
up to rigid motions from $x_0$. The signature of $|dx_0^0|^2$ is
$(1,1)$ for $2\ep_1t<s^2$ and $\mathrm{diag}[-\ep\ \ -\ep]$ for
the other inequality; the change of signature of $|dx_0^0|^2$
occurs along the curves where the normal directions are isotropic.

With $\mathcal{E}:=\mathrm{diag}[1\ -1\ -\ep],\ R_0:=\frac{f_1\bar
f_1^T}{\sqrt{\ep_1\sqrt{-\ep}}} +\sqrt{\ep_1\sqrt{-\ep}}\bar
f_1f_1^T+e_3e_3^T\in\mathbf{O}_3(\mathbb{R})$ we have
$\mathcal{E}f_1=\bar f_1$, so $\mathcal{E}\bar f_1=f_1$; further
$\overline{R_0^Tx_0^0}=\mathcal{E}R_0^Tx_0^0$ and thus all four
surfaces $x_0^0$ are totally real.

Note that the rigid motion $(I_3-\frac{z}{2}(f_1e_3^T-e_3f_1^T)
-\frac{z^2}{8}e_3e_3^T,\frac{z^3}{16}f_1-\frac{3z}{2}\bar
f_1-\frac{3z^2}{8}e_3)
\in\mathbf{O}_3(\mathbb{C})\ltimes\mathbb{C}^3$ brings $x_z(u,v)$
to the canonical form $x_0((u-z),v)$, so $x_z^0$ has real linear
element iff $(u-z,v)$ are as $(u,v)$ above; in particular both
$x_0^0,\ x_z^0$ have real linear element iff
$z\in\sqrt{\ep_1\sqrt{-\ep}}\mathbb{R}$, so each of the four
surfaces comes with its totally real confocal family. Since
$iR(\frac{\pi}{4})x_z(u,v)=x_{\sqrt{i}z}(\sqrt{i}u,iv)$,
multiplication by $i$ and the change of $\ep_1$ has same effect on
the totally real confocal family as that on $x_0^0$; further since
$\overline{R_0^Tx_z^0}=\mathcal{E}R_0^Tx_z^0$ we conclude that
$x_0^0$ and $x_z^0$ stay in the same totally real space. Note that
the Ivory affinity preserves each totally real confocal family and
the collection of these four totally real confocal families meets
each of $\{x_{(\sqrt{i})^jz}\}_{j=0,...7},\ z\in\mathbb{R}^*$
fixed just in one of their four totally real surfaces.

Consider the ruling families parametrization $x_z(\hat u,\hat
v):=\frac{z}{2}(-\hat u\hat vf_1+ (1-\frac{\hat u}{z})2\bar
f_1+(\hat u+\hat v)e_3)$; note that it differs from the usual
ruling families parametrization on $x_z$ invariant under the Ivory
affinity (it is obtained by applying to the later the
transformation
$(u,v)\mapsto(\frac{z}{2}-u,\frac{2}{z}(v-\frac{z^2}{8}))$); for
totally real surface $x_0^0$ and
$z\in\sqrt{\ep_1\sqrt{-\ep}}\mathbb{R}$ we have $(\hat u,\hat v)=
(\sqrt{\ep_1\sqrt{-\ep}}\hat
s,\frac{\sqrt{-\ep}}{\sqrt{\ep_1\sqrt{-\ep}}}\hat t),\ \hat s,\
\hat t\in\mathbb{R}$. Thus $[x_z^0(\hat u,\hat
v)^T,1]=[\frac{z}{|z|}(-\hat u\hat vf_1^T +(1-\frac{\hat
u}{z})2\bar f_1^T+(\hat u+\hat v)e_3^T),\frac{2}{|z|}]$; as
$z\rightarrow\infty$ in its real line and $\hat u,\hat v$ are
fixed we get $[x_{\infty}^0(\hat u,\hat
v)^T,0]=[\sqrt{\ep_1\sqrt{-\ep}}(-\hat u\hat vf_1^T+2\bar f_1^T
+(\hat u+\hat v)e_3^T),0]=[(\ep\ep_1\hat s\hat tf_1^T+2\bar f_1^T
+\sqrt{-\ep}(\ep_1\hat s+\hat t)e_3^T)R_0^T,0]$, which is a double
cover of the concave region of the plane $[(\mathbb{R}\times
i\mathbb{R}\times \sqrt{-\ep}\mathbb{R})R_0^T,0]$ at $\infty$ of
the totally real space $R_0(\mathbb{R}\times
i\mathbb{R}\times\sqrt{-\ep}\mathbb{R})$. The change of $\ep$ in
the linear element because of multiplication by $i$ thus becomes
evident at the infinite singular quadric, where also the change of
$\ep_1$ by the rotation $-R(\frac{\pi}{2})$ can be seen. The
change of $\ep$ occurs because all real lines
$\la\mathbb{R}\cup\{\infty\}\subset\mathbb{C}\cup\{\infty\},\
\la\in\mathbb{C}\setminus\{0\}$ meet at $\infty$ and $0$; when
$z=\infty$ and $\hat u,\hat v$ are finite in their corresponding
real lines, $u,v$ are infinite in their corresponding real lines
(although $du=-d\hat u$ and $dv$ for $d\hat v=0$ may be finite).
In this case $\hat u,\hat v$ change when changing their real
lines, but $u, v$ don't change, so the change in the linear
element is continuous from a projective point of view. Note that
the above considered lines meet also at $0$, which can be
considered a locus of change of linear element ($u=v=du=dv=0$),
but $0$ is not the natural choice, being a point instead of a
surface naturally appearing as a limiting case of a family of
surfaces.

In the diagram below continuous lines are quadrics of the
respective confocal family and with real rulings; to get all four
totally real confocal families one must apply to the diagram below
a rotation $-R(\frac{\pi}{2})$ and superimpose the two.

\begin{center}
$\xymatrix@!0{\{v\in x_{\infty}=\mathbb{P}(\mathbb{R}\times i\mathbb{R}\times\mathbb{R})|
|v|^2>0\}\ar@/_/@{-}[dd]\ar@/^/@{<-}[dd]^{z}&&&&&&&&
\{v\in x_{\infty}=\mathbb{P}(\mathbb{R}\times(i\mathbb{R})^2)|
|v|^2<0\}\ar@/_/@{-}[dd]\ar@/^/@{<-}[dd]^{\sqrt{i}z}\\
&&\ar@{~>}[rrrr]_{iR(\frac{\pi}{4})}&&&&&&\\
&&&&&&&&}$
\end{center}

For IQWC (\ref{eq:iqwc1}), with $U:=\frac{u+v}{\sqrt{2}},\
V:=\frac{u-v}{\sqrt{2}},\ a:=\sqrt{a_1}$ we have $\mathcal{A}=a,\
|dx_0^0|^2=-a^2dV^2+dU\odot d(U^2-V^2)$. (I)$\&$(II) become
$2aU+\frac{V^2}{a}=\ep(2\bar a\bar U+\frac{\bar V^2}{\bar a}),\
2V^2(\frac{U}{a^3}+\frac{1}{a})-a^3=\ep(2\bar V^2(\frac{\bar
U}{\bar a^3}+\frac{1}{\bar a})-\bar a^3)$, which imply
$(\frac{V^2}{a^4}-\frac{\bar V^2}{\bar a^4}) (2aU+a^3(2-\frac{\bar
V^2}{\bar a^4}))-(a^3-\ep\bar a^3) (1-\frac{\bar V^2}{\bar
a^4})^2=0$. If $a^3-\ep\bar a^3\neq 0$, then one can solve for $U$
as a rational function of $V,\ \bar V$; however in this case the
condition that the coefficient of $dV\odot d\bar V$ in
$|dx_0^0|^2-|d\bar x_0^0|^2$ is $0$ cannot be realized without
imposing a real relation between $V$ and $\bar V$. This can be
seen by considering the hot's in $V,\ \bar V$: with
$\la:=\frac{V^2}{a^4}-\frac{\bar V^2}{\bar a^4}$ we have $\la
U=\frac{V^2\bar V^2}{2|a|^4\bar a^2}-\frac{\ep\bar V^4}{2|a|^2\bar
a^4}+\mathrm{lot},\ \la^2dU= (\frac{\ep\bar a^3-a^3}{|a|^{10}\bar
a^3}V\bar V^4+\mathrm{lot})dV+(\frac{V^4\bar
V}{|a|^8a^2}-\frac{2\ep V^2\bar V^3}{|a|^{10}}+\frac{\ep\bar
V^5}{|a|^2\bar a^8}+\mathrm{lot})d\bar V$, so $\la^5((2U\frac{\pa
U}{\pa V}-V)\frac{\pa U}{\pa\bar V} -(2\bar U\frac{\pa\bar
U}{\pa\bar V}-\bar V)\frac{\pa\bar U}{\pa V})$ is a polynomial of
degree $14$ in $V,\ \bar V$, with the coefficient of $V\bar
V^{13}$ being a nonzero multiple of $\ep\bar a^3-a^3$ (it is
obtained by considering only $\la^5U\frac{\pa U}{\pa V}\frac{\pa
U}{\pa\bar V}$). Thus $\ep\bar a^3-a^3=0$ and (III) further
imposes $\frac{V^2}{a^4}-\frac{\bar V^2}{\bar a^4}=0$; using again
(I) we get $\bar V=\ep_1\frac{\bar a^2}{a^2}V,\ \ep_1:=\pm 1,\
\bar U=\frac{\ep a}{\bar a}U$ and (III) is satisfied.

With $\mu:=\frac{\ep\bar
a}{a}=\frac{\ep|a|^2}{a^2}=\frac{\ep|a_1|}{a_1}$ we have
$\mu^3=1$, so there are $6$ types of IQWC (\ref{eq:iqwc1})
$x_0=x_0(u,v)=uvf_1+\frac{u+v}{\sqrt{2}}\bar
f_1+i\sqrt{\ep\mu^2|a_1|}\frac{u-v}{\sqrt{2}}e_3$ containing
surfaces $x_0^0$ with real linear element; each contains two
surfaces $x_0^0$, given by $(u,v)=(\sqrt{\mu}s,\sqrt{\mu}t),\
s,t\in\mathbb{R}$ and $v=\mu\bar u,\ u\in\mathbb{C}$. In the first
case the sign $(\sqrt{\mu})^3$ in $|dx_0^0|^2$ can be absorbed by
$(s,t)\mapsto(-s,-t)$; in the second case $\mu$ in $|dx_0^0|^2$
can be absorbed by $u\mapsto\mu^{-1}u$, so we have four types of
real linear element of IQWC (\ref{eq:iqwc1}) from a real point of
view (namely $-\ep|a_1|\frac{(d(s-t))^2}{2}+\sqrt{2}d(st)\odot
d(s+t)$ and $-\ep|a_1|\frac{(d(u-\bar
u))^2}{2}+\sqrt{2}d|u|^2\odot d(u+\bar u)$); from a complex point
of view we have only two types, since further $\ep$ can be
absorbed by $a_1\mapsto -a_1$. We have $-R(\pi)x_0(u,v)=x_0(v,u)$;
the rigid motions preserving the canonical form account for the
choice of $\mu$ as a cubic root of unity: $R(\frac{2\pi
k}{3})x_0(u,v)=(-R(\pi))^j\hat x_0(-\sqrt{\hat\mu}u,
-\sqrt{\hat\mu}v),\ \sqrt{\hat\mu}\bar
u=\hat\mu\overline{\sqrt{\hat\mu}u}$, where $k=1,2,\ \mu:=1,\
\hat\mu:=e^{\frac{2\pi ik}{3}},\ \sqrt{\hat
a_1}:=(-1)^{j+1}(\sqrt{\hat\mu})^{-1}\sqrt{a_1},\ \hat\ep:=\ep$;
$j$ is chosen so that $\sqrt{\hat a_1}$ is a square root.
Multiplication by $i$ changes $\ep$ while preserving $\mu$; it is
enough to choose $\mu:=1,\ a_1>0:\ iR(-\frac{\pi}{2})x_0(u,v)=\hat
x_0(-u,-v))$, where $\hat\mu:=1,\ \hat a_1:=-a_1$.

In the first case $|dx_0^0|^2$ has signature $(1,1)$ for
$2\sqrt{2}\ep|a_1|(\sqrt{\mu})^3(s+t)>-(s-t)^2$, respectively
$\mathrm{diag}[-\ep\ -\ep]$ for the other inequality; in the
second case $|dx_0^0|^2$ has signature $(1,1)$ for
$2\sqrt{2}\ep|a_1|(\mu u+\bar\mu\bar u)<-(\mu u-\bar\mu\bar u)^2$,
respectively $\mathrm{diag}[\ep\ \ep]$ for the other inequality.
We can make in what follows $\mu=1,\ a_1=\ep|a_1|$. We have $\bar
x_0^0=\mathcal{E}x_0^0$, where in the first case
$\mathcal{E}:=\mathrm{diag}[1\ -1\ -\ep]$ and in the second case
$\mathcal{E}:=\mathrm{diag}[1\ -1\ \ep]$; thus $x_0^0$ are totally
real.

Since $x_z(u,v)-z\bar f_1=(-R(\pi))^j\hat
x_0(u-\frac{z}{2\sqrt{2}},v-\frac{z}{2\sqrt{2}}),\ \sqrt{\hat
a_1}:=(-1)^j\sqrt{a_1}\sqrt{1-za_1^{-1}}$ is the canonical form of
$x_z$, $x_z^0$ is totally real iff $\hat a_1\in\mathbb{R}$; in
particular both $x_0^0,\ x_z^0$ are totally real iff
$z\in\mathbb{R}$; in both cases $\hat\mu=\mu=1$, but $\hat\ep$ may
be different than $\ep$ (precisely when $1-za_1^{-1}<0$). In both
cases since $\bar x_z^0=\hat{\mathcal{E}}x_z^0$ we conclude that
$x_z^0$ stays in the same totally real space as that of $x_0^0$ if
$1-za_1^{-1}>0$ and stays in a totally real space which differs
from that of $x_0^0$ by a multiplication with $i$ in the third
coordinate if $1-za_1^{-1}<0$; the Ivory affinity preserves each
totally real confocal family. For $z=a_1$ we have the singular
totally real quadric $x_{a_1}^0=\mathbb{R}\times
i\mathbb{R}=\{m_1f_1+m_2\bar f_1|m_1,\ m_2\in\mathbb{R}\}$; the
concave domain given by the inequality $(m_2-a_1)^2-2m_1>0$ is
realized as a singular totally real quadric with linear element of
the first case type
$(x_{a_1}^0(s,t)=stf_1+\frac{s+t}{\sqrt{2}}\bar f_1 +a_1\bar f_1,\
s,t\in\mathbb{R}$) and the convex domain given by the other
inequality is realized as a singular totally real quadric with
linear element of the second case type ($x_{a_1}^0(u,\bar u)=
|u|^2f_1+\frac{u+\bar u}{\sqrt{2}}\bar f_1+a_1\bar f_1,\
u\in\mathbb{C}$); both domains meet in the singular part
$\mathcal{S}(x_{a_1}^0)(m_2):=\frac{(m_2-a_1)^2}{2}f_1+m_2\bar
f_1$ for $s=t=u=\bar u:=\pm\frac{m_2-a_1}{\sqrt{2}}$. The
concavity (convexity) can be easily seen because of the
parametrization by two real (complex conjugate) parameters $s,t$
($u,\ \bar u$). Note again that the change of $\mu$ due to the
group generated by $R(\frac{2\pi}{3})$ or the multiplication by
$i$ has same effect on the confocal family as that on $x_0^0$; it
is enough to consider $\mu:=1,\ z\in \mathbb{R}$, when $\hat
z=-\sqrt{\hat\mu}z$.

The homography of Bianchi II previously considered and its
composition with a homothety $i$ takes the totally real IQWC
(\ref{eq:iqwc1}) $x_0^0(u,v)$ to the QC (\ref{eq:qc3}) $\ti
x_0^0(u,v):=\pm\sqrt{\ep_2}(\sqrt{a_1(I_3+a_1J_3)^{-1}})^{-1}
X(\frac{a_1}{2}-\sqrt{2}u,\frac{a_1}{2}-\sqrt{2}v),\ \ep_2:=\pm
1$, so it has the linear element of the image $|d\ti
x_0^0|^2=\frac{2\ep_2}{(u-v)^4}(a_1(d(u-v))^2-2\sqrt{2}(vdu
-udv)\odot d(u-v)+2a_1^{-1}(u-v)^2du\odot dv)$, which is real only
for $\mu:=1$; for $\mu\neq 1$ it is real after an imaginary
homothety. We have $z\in\mathbb{R}\Rightarrow\ti z\in\mathbb{R}$,
so these homographies take totally real confocal quadrics to
totally real confocal quadrics; as we shall see below all totally
real QC (\ref{eq:qc3}) appear as above. The discussion of the
infinite singular quadric of IQWC (\ref{eq:iqwc1}) (where either
$\ep$ or the two cases change) is unnecessary at this point, since
it will correspond via the above considered homographies to the
discussion of the finite singular quadric for QC (\ref{eq:qc3});
also at that time the discussion of the corresponding infinite
singular quadric will not be necessary, it corresponding to the
already completed discussion of the finite singular quadric for
IQWC (\ref{eq:iqwc1}). In the diagram below $a_1=|a_1|$,
(dis)continuous lines represent the quadrics of the respective
confocal families having (imaginary) real rulings; the spectral
parameter $z\ (-z)$ increases with the arrow.

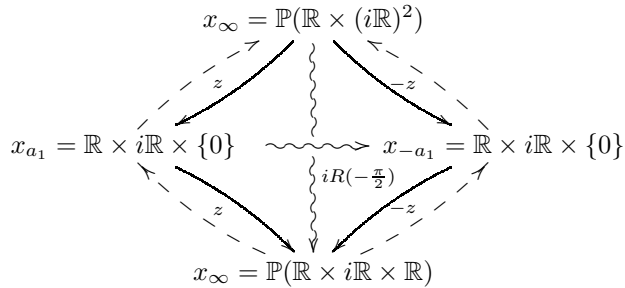
\begin{center}
$\xymatrix@!0{&&&x_{\infty}=\mathbb{P}(\mathbb{R}\times(i\mathbb{R})^2)
\ar@/_/@{<--}[ddlll]\ar@/^/@{->}[ddlll]_{z}
\ar@/^/@{<--}[ddrrr]_{-z}\ar@/_/@{->}[ddrrr]\ar@{~}[dd]&&&\\
&&&&&&\\
x_{a_1}=\mathbb{R}\times i\mathbb{R}\times \{0\}\ar@/_/@{<--}[ddrrr]^{z}
\ar@/^/@{->}[ddrrr]&&\ar@{~>}[rr]&\ar@{~>}[dd]^<<<<{iR(-\frac{\pi}{2})}&&&
x_{-a_1}=\mathbb{R}\times i\mathbb{R}\times \{0\}\ar@/^/@{<--}[ddlll]_{-z}
\ar@/_/@{->}[ddlll]\\
&&&&&&\\
&&&x_{\infty}=\mathbb{P}(\mathbb{R}\times i\mathbb{R}\times\mathbb{R})&&&}$
\end{center}

For QC (\ref{eq:qc3}) with $a:=(\sqrt{a_1^{-1}})^{-1},\
V:=a^4(u-v)^2,\ U:=2\frac{1-a^2(u+v)}{V}+\frac{1}{2}$ we have
$\mathcal{A}=ia^3,\ |dx_0^0|^2=\frac{a^2}{V}((\frac{d(UV)}{2}
-\frac{dV}{4} +\frac{dV}{V})^2-\frac{U}{2V}dV^2)$. (I)$\&$(II)
become $aU=-\ep\bar a\bar U,\
a((U+\frac{2}{V}-\frac{1}{2})^2-2\frac{U}{V})=-\ep\bar a((\bar
U+\frac{2}{\bar V}-\frac{1}{2})^2-2\frac{\bar U}{\bar V})$, which
imply $f(U,V,\bar V):=(1+\frac{\ep a}{\bar a})(\frac{\ep a}{\bar
a}U^2 +2(\frac{1}{V} -\frac{1}{\bar V})U
+(\frac{2}{V}-\frac{1}{2})^2)-2(\frac{1}{V} -\frac{1}{\bar
V})(U+\frac{2}{V}+\frac{2}{\bar V}-1)=0$. If $1+\frac{\ep a}{\bar
a}\neq 0$, then $\frac{\pa f}{\pa U}\neq 0$ and one can solve for
$dU$ from $df=0$; the condition that the coefficient of $dV\odot
d\bar V$ in $|dx_0^0|^2-|d\bar x_0^0|^2$ is $0$ becomes the cubic
equation $\frac{\ep a}{\bar a}(U\frac{\pa f}{\pa U}(\frac{\pa
f}{\pa V}-\frac{\pa f}{\pa\bar V})+(V-\bar V)\frac{\pa f}{\pa
V}\frac{\pa f}{\pa\bar V}-\frac{\pa f}{\pa U}\frac{\pa f}{\pa\bar
V}(\frac{2}{V}-\frac{1}{2}))-\frac{\pa f}{\pa U}\frac{\pa f}{\pa
V}(\frac{2}{\bar V}-\frac{1}{2})=0$ in $U$, which together with
the quadratic equation $f(U,V,\bar V)=0$ in $U$ imposes a real
relation between $V$ and $\bar V$.

Thus $\bar a+\ep a=0$, when (III) implies $\bar V=V$ (the
condition $\bar V\neq V$ imposes $U=1-\frac{2}{V}-\frac{2}{\bar
V}=\bar U$ and a functional relationship between $V$ and $\bar
V$); finally $a_1=-\ep|a_1|,\ (u,v)=(s,t),\ s,t\in\mathbb{R}$ or
$v=\bar u,\ u\in\mathbb{C}$.

We have two  types of QC (\ref{eq:qc3})
$x_0=x_0(u,v)=\frac{\sqrt{|a_1|}}{\sqrt{-\ep}(u-v)}(-(u+\frac{\ep}{2|a_1|})
(v+\frac{\ep}{2|a_1|})f_1+2\bar f_1+(u+v-\frac{\ep}{|a_1|})e_3)$
containing surfaces with real linear element; each contains two
such surfaces $x_0^0$. Thus we have four types of real linear
element of QC (\ref{eq:qc3}) from a real point of view (namely
$8\frac{(tds-sdt)\odot d(s-t)}{(s-t)^4}-4\ep|a_1|\frac{ds\odot
dt}{(s-t)^2},\ 8\frac{(\bar udu-ud\bar u)\odot d(u-\bar
u)}{(u-\bar u)^4}-4\ep|a_1|\frac{du\odot d\bar u}{(u-\bar u)^2}$);
from a complex point of view we have only two, since $\ep$ can be
absorbed by $a_1\mapsto -a_1$. We have
$ix_0^0(u,v)=-\ep(I_3-2e_3e_3^T)\hat x_0^0(-u,-v),\\ \hat
a_1:=-a_1$, so $\hat\ep:=-\ep$.

In the first case $|dx_0^0|^2$ has signature $(1,1)$ for
$4(1+\ep|a_1|(s+t))>-|a_1|^2(s-t)^2$, respectively
$\mathrm{diag}[-\ep\ -\ep]$ for the other inequality; in the
second case $|dx_0^0|^2$ has signature $(1,1)$ for
$4(1+\ep|a_1|(u+\bar u))<-|a_1|^2(u-\bar u)^2$, respectively
$\mathrm{diag}[\ep\ \ep]$ for the other inequality.

With $\ep':=\mathrm{diag}[-\ep\ \ep\ -\ep]$ in the first case or
$:=\mathrm{diag}[\ep\ -\ep\ \ep]$ in the second case we have $\bar
x_0^0=\ep'x_0^0$, so $x_0^0$ are totally real.

Now $x_z$ contains totally real surfaces iff $a_1-z\in\mathbb{R}$;
in particular both $x_0^0,\ x_z^0$ are totally real iff
$a_1,z\in\mathbb{R}$, in which case $x_z^0(u,v)=\hat x_0^0(u,v)$
for $1-za_1^{-1}>0$ and $=\ep\hat x_0^0(u,v)$ for the other
inequality, where $\hat a_1:=a_1-z$. In both cases, since $\bar
x_z^0=\hat\ep'x_z^0$, we conclude that $x_z^0$ stays in the same
totally real space as that of $x_0^0$ if $1-za_1^{-1}>0$ and stays
in a totally real space which differs from that of $x_0^0$ by a
multiplication with $i$ for the other inequality; the Ivory
affinity preserves each confocal family.

For QWC (\ref{eq:qwc2}) with $a:=\sqrt{a_1^{-1}},\ U:=uv,\ V:=v^2$
we have $\mathcal{A}=ia^{-2},\
|dx_0^0|^2=(d(uv))^2+dv^2+\frac{2}{a^2}du\odot
dv=dU^2+\frac{dV^2}{4V}
+\frac{1}{a^2}dV\odot(\frac{dU}{V}-\frac{UdV}{2V^2})$. (I)$\&$(II)
become $a^2(2a^2U-a^4V+1)=-\ep\bar a^2(2\bar a^2\bar U-\bar
a^4\bar V+1),\ a^{12}(2U- a^2V(3+a^4V))=-\ep\bar a^{12}(2\bar U
-\bar a^2\bar V(3+\bar a^4\bar V))$, which imply $2 a^4(a^8-\bar
a^8)U-a^{14}V(3+a^4V) -\ep\bar a^{14}\bar V(3+\bar a^4\bar V)
+a^2\bar a^8(a^4V-1)+\ep\bar a^2a^8(\bar a^4\bar V-1)=0$. If
$a^8-\bar a^8\neq 0$, then (III) cannot be satisfied without
imposing at least a real condition between $V$ and $\bar V$; thus
$a^8=\bar a^8$, from which also follows $a^2(a^4V+1)^2+\ep\bar
a^2(\bar a^4\bar V+1)^2=0$ and (I)$\&$(II) are equivalent, when
$x_0^0$ exists, to

$$(I)'\ \ a^2(2a^2U-a^4V+1)=-\ep\bar a^2(2\bar a^2\bar U-\bar a^4\bar V+1),$$
$$(II)'\ \ a^2(a^4V+1)^2=-\ep\bar a^2(\bar a^4\bar V+1)^2.$$

Now $|dx_0^0|^2=(d(uv))^2+dv^2+\frac{2}{a^2}du\odot
dv=dU^2+\frac{dV^2}{4V}
+\frac{1}{a^2}dV\odot(\frac{dU}{V}-\frac{UdV}{2V^2})$. Applying
$d$ to (I)'$\&$(II)' we get $-\ep d\bar V=\frac{a^6(a^4V+1)}{\bar
a^6(\bar a^4\bar V+1)}dV,\ -\ep d\bar U=\frac{a^4}{\bar
a^4}dU+\frac{a^6(\bar a^4\bar V-a^4V)} {2\bar a^4(\bar a^4\bar
V+1)}dV$, so (III) becomes the condition that a symmetric bilinear
form in $dU,\ dV$ is identically $0$. In particular the condition
that the coefficient of $dU\odot dV$ is $0$ gives $(a^4V-\bar
a^4\bar V)((a^4V-1)(\bar a^4\bar V-1)-2)=0$, which together with
(II)' implies $a^4V=\bar a^4\bar V$, so finally $\bar a_1=-\ep
a_1,\bar v=-\ep\ep_1v,\ \bar u=\ep_1u,\ \ep_1:=\pm 1$ and (III) is
satisfied.

Thus we have two types of QWC (\ref{eq:qwc2})
$x_0=x_0(u,v)=(\sqrt{a_1^{-1}})^{-1}((u+\frac{v}{2a_1})f_1+v\bar
f_1)+uve_3$ containing surfaces $x_0^0$ with real element (given
by the choice of $\ep$ in $a_1=(\sqrt{-\ep})^{-1}\ep_2|a_1|,\
\ep_2:=\pm 1$); each contains two surfaces $x_0^0$ (given by the
choice of $\ep_1$ in
$(u,v)=(\sqrt{\ep_1}s,\sqrt{-\ep}\sqrt{\ep_1}t),\
s,t\in\mathbb{R}$). The transformation $u\mapsto -u$ changes
$\ep_2$ in $|dx_0^0|^2$, so we have from a real point of view four
different types of real linear element of QWC (\ref{eq:qwc2})
(namely $|dx_0^0|^2=-\ep((d(st))^2+\ep_1dt^2+2|a_1|ds\odot dt)$);
from a complex point of view we have only one type, since
$(s,t)\mapsto (is,-it)$ changes $\ep_1$ and $(a_1,s,t)\mapsto
(ia_1,is,-it)$ changes $\ep$. All symmetries of the linear element
are realized by rigid motions composed when necessary with a
multiplication by $i$ or with a simple complex affine
transformation: $R(\pi)x_0(u,v)=x_0(-u,-v),\ (-i\ep\ep_2f_1\bar
f_1^T+i\ep\ep_2\bar f_1f_1^T-e_3e_3^T)x_0(u,v)=\hat x_0(-u,v)$,
where $\hat a_1:=-a_1$,
$i(\frac{\sqrt{\ep_2\sqrt{\ep}}}{\sqrt{\ep_2\sqrt{-\ep}}}f_1\bar
f_1^T+\frac{\sqrt{\ep_2\sqrt{-\ep}}}{\sqrt{\ep_2\sqrt{\ep}}}\bar
f_1f_1^T+\ep e_3e_3^T)x_0(u,v)=\hat x_0(\ep u,iv)$, where $\hat
a_1:=i\ep a_1$, $\mathrm{diag}[i\ \ i\ \ -1]
(\frac{\sqrt{\ep\ep_2\sqrt{-\ep}}}{\sqrt{\ep_2\sqrt{-\ep}}}f_1\bar
f_1^T
+\frac{\sqrt{\ep_2\sqrt{-\ep}}}{\sqrt{\ep\ep_2\sqrt{-\ep}}}\bar
f_1f_1^T+\ep e_3e_3^T) x_0(u,v)=\hat x_0(i\ep u,iv)$, where $\hat
a_1:=\ep a_1$.

Now $|dx_0^0|^2$ has signature $(1,1)$ for
$\ep_1t^2-2\ep_1\ep_2st|a_1|-|a_1|^2<0$ and $\mathrm{diag}[-\ep\
-\ep]$ for the other inequality. With $\ep':=\mathrm{diag}[1\ \
-1\ \ -\ep],\ R_0:=\sqrt{\frac{\ep_1\ep_2}{\sqrt{-\ep}}}f_1\bar
f_1^T+ (\sqrt{\frac{\ep_1\ep_2}{\sqrt{-\ep}}})^{-1}\bar
f_1f_1^T+e_3e_3^T$ we have $\overline{R_0^Tx_0^0}
=\ep'R_0^Tx_0^0$, so all $x_0^0$ are totally real.

Since $x_z(u,v)-\frac{z}{2}e_3=R(\pi)^j\hat x_0(u,v)$, where $\hat
a_1:=a_1-z,\ j:=0,1$ is the canonical form of $x_z$, $x_z^0$ is
totally real iff $a_1-z$ is as $a_1$ above; in particular both
$x_0^0,\ x_z^0$ are totally real iff $\bar z=-\ep z$. In this case
$j=0$ for $\hat\ep_2=\ep_2$ or $\ep\ep_2=1$ and $j=1$ otherwise;
since $\overline{R_0^Tx_z^0}=\ep'R_0^Tx_z^0$ for
$\hat\ep_2=\ep_2$, it is natural to require that the whole
confocal family should stay in the same totally real space
$R_0(\mathrm{R})\times i\mathbb{R}\times(\sqrt{-\ep}\mathbb{R}))$;
this will suffice to conclude $\hat\ep_1=\ep_1\ep_2\hat\ep_2$. The
rigid motions considered above and the multiplication by $i$ have
the same effect on the totally real confocal family $x_z^0$ as
that on $x_0^0$ ($\ep_2$ is not replaced with $\hat\ep_2$ in the
definition of these rigid motions); however the affine
transformation considered above cannot have the same effect (as we
know that the only affine transformations which take confocal
quadrics to confocal quadrics are homotheties and rigid motions):
it involves a translation $-\ep ze_3$ on the rhs term and $\ep_2$
is replaced with $\hat\ep_2$. Note that the Ivory affinity does
not preserve the totally real confocal family $x_z^0$; for
$\ep_2\hat\ep_2=-1$ it mixes the cases $\ep_2:=1$ and $\ep_2:=-1$,
but according to the above affine transformation (or variants
thereof). When $z$ varies with $\bar z=-\ep z$ and $\hat a_1$
changes sign, $\hat\ep_1$ changes. Again this can be explained if
one considers the singular quadric $x_{a_1}$.

Consider the ruling families parametrization $x_z^0(\hat u,\hat
v):=(\hat u+\frac{\hat v}{2})f_1 +(a_1-z)\hat v\bar f_1+(\hat
u\hat v+\frac{z}{2})e_3$; note that it differs from the usual
ruling families parametrization on $x_z^0$ invariant under the
Ivory affinity (it is obtained by applying to the later the
transformation $(u,v)\mapsto
(\sqrt{1-za_1^{-1}}(\sqrt{a_1^{-1}})^{-1}u,\\
\frac{(\sqrt{a_1^{-1}})^{-1}}{a_1\sqrt{1-za_1^{-1}}}v)$, which, as
$z\rightarrow a_1$, contracts $u$ and expands $v$), but in a way
such that $\hat u,\ \hat v$ vary in real lines, if $\bar z=-\ep z$
and $u,\ v$ vary in their corresponding real lines. Moreover, the
real lines of $\hat u,\ \hat v$ are independent of the sign of
$1-za_1^{-1}$, since $\ep_1$ changes sign with $1-za_1^{-1}$.
Multiplying by $i$ and applying further rigid motions, if
necessary, it is enough to consider only the case $a_1>0,\
\ep_1:=1$, when we have $(\hat u,\hat v)=(\hat s,\hat t),\ \hat
s,\hat t\in\mathbb{R}$. In this case $x_{a_1}^0(\hat s,\hat
t)=(\hat s+\frac{\hat t}{2})f_1+(\hat s\hat t+\frac{a_1}{2})e_3$
appears as a singular case for both types of totally real linear
element $|dx_z^0|^2$ with $1-za_1^{-1}>0$ and $1-za_1^{-1}<0$, it
has degenerate (isotropic) linear element $(d(\hat s\hat t))^2$,
covers twice the concave domain
$\{m_1f_1+m_2e_3|m_1,m_2\in\mathbb{R},\ m_1^2-2m_2+a_1>0\}$ and
meets $\mathcal{S}(x_{a_1})=\{m_1f_1+m_2e_3|m_1,m_2\in\mathbb{C},\
m_1^2-2m_2+a_1=0\}$ for $\hat s=\frac{\hat t}{2}$. Although $u,v$
change their real lines for $z\neq a_1$, they do not change it for
$z=a_1$, when $u=\infty,\ v=0,\ uv=\hat u\hat v,\ du=\infty$ if
$d\hat u\neq 0$ and $dv=0$ if $d\hat v\neq \infty$.

The homography of Bianchi II previously considered already takes
totally real confocal quadrics QWC (\ref{eq:qwc2}) to same type of
quadrics: $\ti x_0^0(u,v)=\sqrt{a_1^{-1}}(uf_1-\frac{1}{v}\bar
f_1)- \frac{1}{2\sqrt{a_1^{-1}}v}f_1-\frac{u}{v}e_3$; it and its
composition with a homothety $i$ are the only such homographies.
We thus have $\ti a_1=a_1^{-1},\ \ti u=u,\ \ti v=-v^{-1}$; for
$(\sqrt{a_1})^{-1}=-\sqrt{a_1^{-1}}$ or multiplication by $i$
further rigid motions may be required.

Again the discussion of the infinite singular quadric of QWC
(\ref{eq:qwc2}) (where $\ep$ may change) is not necessary, since
it corresponds via the above considered homographies to the
already completed discussion of the finite singular quadric of QWC
(\ref{eq:qwc2}).

For QWC (\ref{eq:qwc1}) with $a:=(\sqrt{a_1^{-1}})^{-1},\
b:=-i(\sqrt{a_2^{-1}})^{-1},\ \mathcal{U}:=\frac{u+v}{\sqrt{2}},\
\mathcal{V}:=\frac{u-v}{\sqrt{2}}$ we have
$x_0^0=a\mathcal{U}e_1+b\mathcal{V}e_2+\frac{\mathcal{U}^2-\mathcal{V}^2}{2}e_3,\
\hat N_0^0=\frac{\mathcal{U}}{a}e_1-\frac{\mathcal{V}}{b}e_2-e_3,\
A\hat
N_0^0=\frac{\mathcal{U}}{a^3}e_1+\frac{\mathcal{V}}{b^3}e_2,\
V_2=\frac{\mathcal{V}}{b^3}e_1-\frac{\mathcal{U}}{a^3}e_2
+\frac{a^2+b^2}{a^3b^3}\mathcal{U}\mathcal{V}e_3,\
\mathcal{A}=-ab$. With $(U,V):=(\mathcal{U}^2,\mathcal{V}^2)$
(I)$\&$(II) become $\frac{b}{a}U+\frac{a}{b}V+ab=\ep(\frac{\bar
b}{\bar a}\bar U +\frac{\bar a}{\bar b}\bar V+\bar a\bar b),\
\frac{b^3}{a^3}U+\frac{a^3}{b^3}V
+\frac{(a^2+b^2)^2}{a^3b^3}UV=\ep(\frac{\bar b^3}{\bar a^3}\bar U
+\frac{\bar a^3}{\bar b^3}\bar V+\frac{(\bar a^2+\bar b^2)^2}{\bar
a^3\bar b^3} \bar U\bar V)$, which imply
$\frac{b}{a}U(\frac{b^2}{a^2}-\frac{\bar b^2}{\bar a^2}
+\frac{(a^2+b^2)^2}{a^3b^3}\frac{a}{b}V-\frac{(\bar a^2+\bar
b^2)^2}{\bar a^3\bar b^3} \frac{\bar a}{\bar b}\bar
V)+\frac{a^3}{b^3}V-\ep\frac{\bar a^3}{\bar b^3}\bar V
-(\frac{\bar b^2}{\bar a^2}+\frac{(\bar a^2+\bar b^2)^2}{\bar
a^3\bar b^3} \frac{\bar a}{\bar b}\bar V)(\frac{a}{b}V
+ab-\ep(\frac{\bar a}{\bar b}\bar V+\bar a\bar b))=0$. If
$a^2+b^2=0$, then also $ab=\ep\bar a\bar b$; since both $a^{-1}$
and $-ib^{-1}$ are square roots, we get $b=-ia$ and $a^2=-\ep\bar
a^2$; now (I) becomes $\bar U-\bar V=-\ep(U-V)$, or $\bar u\bar
v=-\ep uv$. Keeping account of
$|dx_0^0|^2=\frac{1}{4}(a^2\frac{dU^2}{U}
+b^2\frac{dV^2}{V}+(dU-dV)^2)$, (III) factors as $d\log\frac{\bar
v}{v}\odot d\log\frac{\bar v}{u}=0$, so we have two types of QWC
(\ref{eq:qwc1}) with $a_1=a_2$ containing surfaces $x_0^0$ with
real linear element (given by $\ep$ in
$a_1=\ep_2(\sqrt{-\ep})^{-1}|a_1|,\ \ep_2:=\pm 1$); each type
contains three surfaces $x_0^0$, according to
$(u,v)=(\sqrt{-\ep}e^{ic}s,e^{-ic}t),\ s,t\in\mathbb{R}$ or
$v=e^{-2c}\ep_1\sqrt{-\ep}\bar u,\ u\in\mathbb{C},\ \ep_1:=\pm 1$;
in both cases $c\in\mathbb{R}$ is a constant. In the first case
$s\mapsto -s$ changes $\ep_2$ in $|dx_0^0|^2$; in the second case
$c$ can be absorbed by $u\mapsto e^{-c}u$, so we have $6$ types of
real linear element of QWC (\ref{eq:qwc1}) with $a_1=a_2$ (namely
$2|a_1|ds\odot dt-\ep(d(st))^2$ or $2\ep_1\ep_2|a_1|du\odot d\bar
u -\ep(d|u|^2)^2$); from a complex point of view we have only two
types, since the transformation $(a_1,u)\mapsto(i\ep a_1,i\ep u)$
changes $\ep$ in the first case or $\ep,\ \ep_1$ in the second one
and $a_1\mapsto -a_1$ changes $\ep_2$ in the second case. All
symmetries of the linear element are realized by rigid motions
composed when necessary with a multiplication by $i$ or a simple
complex affine transformation: $R(-c)$ in the first case or
$R(ic)$ in the second one makes $c=0,\
(I_3-2e_2e_2^T)x_0(u,v)=x_0(v,u),\
(\frac{\sqrt{-\ep_2\sqrt{-\ep}}}{\sqrt{\ep_2\sqrt{-\ep}}}f_1\bar
f_1^T
+\frac{\sqrt{\ep_2\sqrt{-\ep}}}{\sqrt{-\ep_2\sqrt{-\ep}}}\bar
f_1f_1^T -e_3e_3^T)x_0(u,v)=\hat x_0(-u,v)$, where $\hat
a_1:=-a_1$, so it changes $\ep_2$ in the first case and both
$\ep_1,\ \ep_2$ in the second case;
$ix_0(u,v)=(-i\frac{\sqrt{\ep_2\sqrt{-\ep}}}{\sqrt{\ep_2\sqrt{\ep}}}f_1\bar
f_1^T +i\frac{\sqrt{\ep_2\sqrt{\ep}}}{\sqrt{\ep_2\sqrt{-\ep}}}\bar
f_1f_1^T+\ep e_3e_3^T)\hat x_0(i\ep u,v)$, where $\hat a_1:=i\ep
a_1$ so multiplication by $i$ changes, up to a rigid motion, $\ep$
in the first case and $\ep,\ \ep_1$ in the second one; also
$\mathrm{diag}[i\ i\ 1]x_0^0(u,v)=\hat
x_0^0(\ep\ep_2u,\ep\ep_2v),\ \hat a_1:=-a_1$, so it changes
$\ep_2$.

In the first case $|dx_0^0|^2$ has signature $(1,1)$ for
$2\ep\ep_2st<|a_1|$ (respectively $\mathrm{diag}[-\ep\ \ -\ep]$
for the other inequality); in the second case we have signature
$(1,1)$ for $2\ep\ep_1\ep_2e^{-2c}|u|^2>|a_1|$ (note that we need
in this case $\ep\ep_1\ep_2=1$), respectively
$\mathrm{diag}[\ep_1\ep_2\ \ \ep_1\ep_2]$  for the other
inequality. In the first case with
$R_0:=\sqrt{\ep_2\sqrt{-\ep}}f_1\bar f_1^T
+\frac{1}{\sqrt{\ep_2\sqrt{-\ep}}}\bar f_1f_1^T+e_3e_3^T,\
\ep':=\mathrm{diag}[1\ \ -1\ \ -\ep]$ we have
$\overline{R_0^Tx_0^0}=\ep'R_0^Tx_0^0$; in the second case with
$\ep':=\mathrm{diag}[\ep_1\ep_2\ \ \ep_1\ep_2\ \ -\ep]$ we have
$\bar x_0^0=\ep'x_0^0$, so all surfaces $x_0^0$ are totally real.

Since $x_z(u,v)-\frac{z}{2}e_3=\hat x_0(u,v),\ \hat a_1:=a_1-z$ is
the canonical form of $x_z$, $x_z^0$ is totally real iff $\hat a$
is as $a$ above; in particular both $x_0^0,\ x_z^0$ are totally
real iff $z\in\sqrt{-\ep}\mathbb{R}$ (although $\ep_2$ may
change).

If $a^2+b^2\neq 0$, then with $\la:=\frac{(a^2+b^2)^2}{a^3b^3}$
the consequence of (I)$\&$(II) can be written as
$(\la\frac{b}{a}U-\bar\la\frac{\bar a}{\bar b}\bar V+
\frac{a^2}{b^2}-\frac{\bar b^2}{\bar
a^2})(\frac{b^2}{a^2}-\frac{\bar b^2}{\bar a^2}
+\la\frac{a}{b}V-\bar\la\frac{\bar a}{\bar b}\bar
V)-(\bar\la-\ep\la) \frac{(\bar V(\bar a^2+\bar b^2)+\bar
b^4)^2}{\bar a\bar b^5}=0$. If $\bar\la\neq\ep\la$, then one can
solve the previous for $U$ as a rational function of $V,\ \bar V$,
but (III) cannot be satisfied without imposing at least a real
relation between $V$ and $\bar V$. This can be seen by looking at
the hot's in $U,\ V,\ \bar V$; with
$\mu:=\la\frac{a}{b}V-\bar\la\frac{\bar a}{\bar b}\bar V$ we have:
$\la\frac{b}{a}U\mu=\bar\la\frac{\bar a}{\bar b}\bar V(\mu
+(\bar\la-\ep\la)\frac{\bar a}{\bar b}\bar V)+\mathrm{lot},\
\la\frac{b}{a}\mu^2dU=\bar\la\frac{\bar a}{\bar b}(\mu d(\bar
V(\mu+(\bar\la-\ep\la)\frac{\bar a}{\bar b}\bar V))-\bar V(\mu
+(\bar\la-\ep\la)\frac{\bar a}{\bar b}\bar V)d\mu)+\mathrm{lot}$.
If we ignore the lot's, the hot of the numerator of the
coefficient of $dV\odot d\bar V$ in $|dx_0^0|^2-|d\bar x_0^0|^2$
is a homogeneous polynomial of degree $4$ in $V$ and $\bar V$,
namely $(\mu^2\frac{\pa U}{\pa V}) (\mu^2\frac{\pa U}{\pa\bar
V})-\mu^2(\mu^2\frac{\pa U}{\pa V}) -(\mu^2\frac{\pa\bar U}{\pa
V})(\mu^2\frac{\pa\bar U}{\pa\bar V}) +\mu^2(\mu^2\frac{\pa\bar
U}{\pa V})$; the coefficient of $\bar V^4$ is $-\ep(a^2+b^2)
(\bar\la-\ep\la)\frac{a\bar a^5\bar\la^3}{b^3\bar b^5}\neq 0$, so
(III) imposes a functional relationship between $V$ and $\bar V$.
Thus $\bar\la=\ep\la$, when we have either $\frac{b}{a}U-
\ep\frac{\bar a}{\bar b}\bar V+\frac{1}{\la}
(\frac{a^2}{b^2}-\frac{\bar b^2}{\bar a^2})=0$ or
$\frac{b^2}{a^2}- \frac{\bar b^2}{\bar
a^2}+\la(\frac{a}{b}V-\ep\frac{\bar a}{\bar b}\bar V)=0$. In the
later case (I)$\&$(II) are equivalent to

$$(I)'\ \ \frac{a}{b}V-\ep\frac{\bar a}{\bar b}\bar V
=\frac{1}{\la}(\frac{\bar b^2}{\bar a^2}-\frac{b^2}{a^2}),$$
$$(II)'\ \ \frac{b}{a}U-\ep\frac{\bar b}{\bar a}\bar U
=\frac{1}{\la}(\frac{\bar a^2}{\bar b^2}-\frac{a^2}{b^2}).$$

Applying $d$ to (I)'$\&$(II)' we get $\frac{a}{b}dV- \ep\frac{\bar
a}{\bar b}d\bar V=\frac{b}{a}dU-\ep\frac{\bar b}{\bar a}d\bar
U=0$; now (III) imposes $\frac{a}{\bar a}=\ep\frac{\bar b}{b}=\pm
1$; finally we get $a_1=\ep_1|a_1|,\ a_2=-\ep\ep_1|a_2|,\
\ep_1:=\pm 1,\ (u,v)=(\sqrt{\ep_2}s,\sqrt{\ep_2}t),\
s,t\in\mathbb{R}$ or $v=\ep_2\bar u,\ u\in\mathbb{C},\ \ep_2:=\pm
1$.

In the former case $\frac{b}{a}U-\ep\frac{\bar a}{\bar b}\bar
V+\frac{1}{\la} (\frac{a^2}{b^2}-\frac{\bar b^2}{\bar a^2})=0$ and
its conjugate are equivalent to (I)$\&$(II). Differentiating these
we get $\frac{b}{a}dU- \ep\frac{\bar a}{\bar b}d\bar V=\frac{\bar
b}{\bar a}d\bar U- \ep\frac{a}{b}dV=0$; now (III) imposes $\bar
a^2=\ep b^2$; finally we get $a_2=-\ep\bar a_1,\
a_1\in\mathbb{C}\setminus\{0\}$ (note that $a_1\in\mathbb{R}
\setminus\{0\}$ is allowed),
$(u,v)=(\sqrt{\ep_1\sqrt{\ep}}s,\sqrt{-\ep_1\sqrt{\ep}}t),\
s,t\in\mathbb{R},\ \ep_1:=\pm 1$.

For QC (\ref{eq:qc2}) with $a:=(\sqrt{a_1^{-1}})^{-1},\
b:=(\sqrt{a_2^{-1}})^{-1},\ V:=(u-v)^{-2},\ U:=(u+v)^2V$
(I)$\&$(II) become
$b(U(1-\frac{a^2}{b^2})-1+\frac{4V}{a^2})=-\ep\bar b(\bar
U(1-\frac{\bar a^2}{\bar b^2})-1+\frac{4\bar V}{\bar a^2}),\
b(U(1-\frac{a^2}{b^2})((1-\frac{a^2}{b^2})(U-1)
+(3-\frac{a^2}{b^2})\frac{4V}{a^2})+\frac{b^2}{a^2}(\frac{4V}{a^2})^2)=
-\ep\bar b(\bar U(1-\frac{\bar a^2}{\bar b^2})( (1-\frac{\bar
a^2}{\bar b^2})(\bar U-1)+(3-\frac{\bar a^2}{\bar b^2})\frac{4\bar
V}{\bar a^2})+\frac{\bar b^2}{\bar a^2}(\frac{4\bar V}{\bar
a^2})^2)$; with $\mathcal{V}:=\frac{4V}{a^2},\
\mathcal{E}:=1-\frac{a^2}{b^2},\
\mathcal{U}:=U\mathcal{E}-1+\mathcal{V}$ these imply $(1+\frac{\ep
b}{\bar b})(\frac{\ep b}{\bar b}\mathcal{U}^2+
(\mathcal{E}(\mathcal{V}-1)-\overline{\mathcal{E}}(\overline{\mathcal{V}}-1))\mathcal{U}+
1-\mathcal{E}((1-\mathcal{V})^2+\mathcal{V}^2))-
(\mathcal{E}(\mathcal{V}-1)-\overline{\mathcal{E}}(\overline{\mathcal{V}}-1))\mathcal{U}+
\mathcal{E}((1-\mathcal{V})^2+\mathcal{V}^2)
-\overline{\mathcal{E}}((1-\overline{\mathcal{V}})^2+\overline{\mathcal{V}}^2)
=0$; further $|dx_0^0|^2=\frac{b^2}{4U}dU^2-\frac{a^2}{2V}dU\odot
dV +\frac{4V-a^2(1-U)}{4V^2}dV^2$.

If $a^2=b^2$, then (I)$\&$(II) imply $\bar a=-\ep a,\ \bar V=V$;
now (III) becomes $\frac{V}{U}(d\frac{U}{V})^2-\frac{\bar V}{\bar
U}(d\frac{\bar U}{\bar V})^2=0$, or $(du+dv)^2-(d\bar u+d\bar
v)^2=0$ and finally $a_1=a_2=-\ep|a_1|,\
(u,v)=(\sqrt{\ep_1}(s+ic), \sqrt{\ep_1}(t+ic)),\ s,t\in\mathbb{R}$
or $v=\ep_1\bar u-i\sqrt{\ep_1}c,\ u\in\mathbb{C},\
c\in\mathbb{R}$ constant, $\ep_1:=\pm 1$.

If $a^2\neq b^2$, then we have
$|dx_0^0|^2=\frac{a^2(d\mathcal{U}-d\mathcal{V})^2}
{4\mathcal{E}(1-\mathcal{E})(\mathcal{U}+1-\mathcal{V})}
-\frac{a^2}{2\mathcal{E}\mathcal{V}}d(\mathcal{U}-\mathcal{V})\odot
d\mathcal{V}+\frac{a^2}{4\mathcal{V}^2}(\mathcal{V}-1
+\frac{\mathcal{U}+1-\mathcal{V}}{\mathcal{E}})d\mathcal{V}^2$.

\subsection{Two algebraic consequences of the tangency configuration}\label{subsec:algpre23}
\noindent

\noindent Choose a surface $x_0^0$ in the quadric $x_0$ (that is
$(u_0,v_0)=(u_0(s,t),v_0(s,t))$, where $s,t$ are independent
variables (real or complex, but not necessarily both of the same
type) such that $u_0,\ v_0$ are still independent; thus
$T_{x_0^0}x_0^0=T_{x_0^0}x_0$). $T_{x_0^0}x_0$ cuts $x_z$ along a
conic; for every point $x_z^1$ on this conic we have the two
facets at $x_z^1$ spanned by $V_0^1:=x_z^1-x_0^0$ and one of the
rulings of $x_z$ at $x_z^1$. The two distributions of facets
$\mathcal{D}^1,\ \mathcal{D'}^1$ are integrable, their leaves
being the ruling families on $x_z$ (we shall see in \S\
\ref{subsec:deformations1} that the complete integrability depends
only on the linear element of $x_0^0$). Choose
$m_0^1:=\mathcal{B}_1x_{zu_1}^1\times V_0^1$ a normal field of the
distribution $\mathcal{D}^1$ with leaves $x_{zu_1}^1$ (and
similarly $m'^1_0:=\mathcal{B}_1x_{zv_1}^1\times V_0^1$ by
considering the other ruling family on $x_z^1$). Note that
$m^1_0=-i\det(\sqrt{R_z}(\sqrt A)^{-1}) (\sqrt{R_z})^{-1}\sqrt
AY(v_1)-x_0^0\times\sqrt{R_z}(\sqrt A)^{-1}Y(v_1)$ for QC,
$=2x_{zu_1}^1\times(x_z^1(0,v_1)-x_0^0)$ for (I)QWC ($m'^1_0=
-i\det(\sqrt{R_z}(\sqrt A)^{-1})(\sqrt{R_z})^{-1}\sqrt AY(u_1)
+x_0^0\times\sqrt{R_z}(\sqrt A)^{-1}Y(u_1)$ for QC,
$=2x_{zv_1}^1\times(x_z^1(u_1,0)-x_0^1)$ for (I)QWC) depends
quadratically on $v_1\ (u_1)$, fact which will make the
differential system of the B transformation a Ricatti equation.

The next two algebraic results (essentially due to Bianchi; see
(122)) play a fundamental r\^{o}le in the theory of deformations
of quadrics:

{\it If $x_z^1\in T_{x_0^0}x_0$, then:

I The change in the linear element from $x_z^1$ to $x_0^1$ is four
times the product of the orthogonal projections of the
differentials of the rulings of $x_z^1$ on the normal of $x_0$ at
$x_0^0$.

II The facets at $x_z^1$ spanned by $V_0^1$ and one of the rulings
of $x_z^1$ reflect in $T_{x_0^0}x_0$; therefore the distributions
$\mathcal{D}^1,\ \mathcal{D'}^1$ reflect in $Tx_0^0$.}

This can be put as: if $(V_0^1)^T\hat N_0^0=0$, then:
\begin{eqnarray}\label{eq:simpas}
4(x_{zu_1}^1)^TN_0^0(N_0^0)^Tx_{zv_1}^1du_1dv_1=|dx_z^1|^2-|dx_0^1|^2,\
(x_{zv_1}^1)^T(I_3-2N_0^0(N_0^0)^T)m_0^1=0.
\end{eqnarray}
Since $|dx_z|^2-|dx_0|^2=-zdx_0^TAdx_0=-z|dX|^2$ for QC,
$=-z|I_{1,2}dZ|^2$ for (I)QWC, these boil down to:
$\mathcal{B}_1(x_{zu_1}^1)^T\hat N_0^0(\hat
N_0^0)^Tx_{zv_1}^1=-z|\hat N_0^0|^2,\ (x_{zv_1}^1\times
x_{zu_1}^1)^T\hat N_0^0(\hat N_0^0)^TV_0^1=(\hat N_0^0)^T
(x_{zu_1}^1(x_{zv_1}^1)^T+x_{zv_1}^T(x_{zu_1}^1)^T)(V_0^1\times\hat
N_0^0)$.

The relation
\begin{eqnarray}\label{eq:baseb}
(x_{zv_1}^1(x_{zu_1}^1)^T+x_{zu_1}^1(x_{zv_1}^1)^T)\hat N_0^0=
x_{zu_1v_1}^1(V_0^1)^T\hat
N_0^0-\frac{2}{\mathcal{B}_1}(V_0^1+z\hat N_0^0)
\end{eqnarray}
follows immediately from (\ref{eq:basic}).

Multiplying (\ref{eq:baseb}) with $(\hat N_0^0)^T$ on the left we
get the fact that the first relation of (\ref{eq:simpas}) for
(I)QWC is equivalent to the TC $(V_0^1)^T\hat N_0^0=0$; for QC it
is equivalent to the TC $(V_0^1)^T\hat N_0^0=0$ or to
$X_0^T\sqrt{R_z}X_1=-1$, which can be brought to the TC by
$u_0\leftrightarrow v_0$. Multiplying (\ref{eq:baseb}) with
$(V_0^1\times\hat N_0^0)^T$ on the left we get the fact that the
second relation of (\ref{eq:simpas}) is equivalent to the TC
$(V_0^1)^T\hat N_0^0=0$ or to $0=(\hat N_0^0)^T(V_{0u_1}^1\times
V_0^1)_{v_1}=(\hat N_0^0)^T(\frac{m_0^1}{\mathcal{B}_1})_{v_1}$.
If $(V_0^1)^T\hat N_0^0\neq 0$ and the facets reflect in
$T_{x_0^0}x_0$, then from $(V_0^1)^Tm'^1_0=0$ we get $(\hat
N_0^0)^Tm_0^1=0$; thus $m_0^1,\ m'^1_0$ are multiples of $\hat
N_z^1$ and $(\hat N_0^0)^T\hat N_z^1=0$; essentially by Chasles's
result we get $(V_0^1)^T\hat N_0^0=0$, a contradiction. Therefore
the two algebraic consequences of the TC are actually equivalent
to the TC.

Note:
\begin{eqnarray}\label{eq:symrel}
(V_0^1)^T\hat N_0^0=X_0^T\sqrt{R_z}X_1-1\ \mathrm{for\ QC},\nonumber\\
=Z_0^TI_{1,2}\sqrt{R'_z}Z_1-f^T(x_z^0+x_z^1-C(z))=Z_0^TI_{1,2}\sqrt{R'_z}Z_1-
\nonumber\\-(Z_0+Z_1)^T(u(z)f_1+v(z)\bar f_1+e_3)-u(z)v(z)-\frac{|f|^2}{2}z\ \mathrm{for\ (I)QWC}.
\end{eqnarray}
If we let $u_1=v_1=0$ in (\ref{eq:symrel}) for (I)QWC we get
$(\hat N_0^0)^TC(z)+f^Tx_z^0=(\hat N_0^0)^Tx_0^0(=u_0v_0)$.

Consider now the TC $(V_0^1)^T\hat N_0^0=0$; one can choose $v_1$
arbitrarily and define
\begin{eqnarray}\label{eq:u2}
u_1:=v_1+\frac{Y(v_1)^T\sqrt{R_z}X_0}{1-\frac{1}{2}Y'(v_1)^T\sqrt{R_z}X_0}\
\mathrm{for\ QC},\
u_1:=\frac{v_1(\hat N_0^0)^Tx_{zv}(0,0)-f^Tx_z^0}{v_1-(\hat N_0^0)^Tx_{zu}(0,0)}\
\mathrm{for\ (I)QWC};\nonumber\\
\end{eqnarray}
it is a ratio of two functions separately linear in $u_0,v_0,v_1$
and thus a homography is established between $u_0,v_0,u_1,v_1$.

If we let $u_1=v_1=0$ in (\ref{eq:baseb}) for (I)QWC and multiply
with $(\hat N_0^0)^T$ on the left we get:
\begin{eqnarray}\label{eq:n(i)qwc}
\frac{z}{2}|\hat N_0^0|^2=f^Tx_z^0-x_{zu}(0,0)^T\hat N_0^0(\hat N_0^0)^Tx_{zv}(0,0);
\end{eqnarray}
note also for (I)QWC
\begin{eqnarray}\label{eq:nor}
(\hat N_0^0)^T(u_1x_{zu}(0,0)+ v_1x_{zv}(0,0))=(\hat N_0^0)^TV_0^1+f^Tx_z^0+u_1v_1,
\end{eqnarray}
which is equivalent to (\ref{eq:u2}) in the TC case. The algebraic
relation
\begin{eqnarray}\label{eq:intalg}
(N_0^0)^T(2zm_0^1+m_0^1\times m_{0v_1}^1)=0
\end{eqnarray}
will appear as the complete integrability condition of the Ricatti
equation subjacent to the B transformation. Using
(\ref{eq:simpas}) this becomes:
$0=\frac{z(m_0^1)^Tx_{zv_1}^1}{(N_0^0)^Tx_{zv_1}^1}-\mathcal{B}_1(x_{zu_1}^1)^TN_0^0(V_0^1)^T
m_{0v_1}^1=\\=\frac{z(V_0^1)^T(\mathcal{B}_1x_{zv_1}^1\times
x_{zu_1}^1+(\mathcal{B}_1x_{zu_1}^1\times
V_0^1)_{v_1})}{(N_0^0)^Tx_{zv_1}^1}$, which is straightforward.
Replacing $(m_0^1,v_1)$ with $(m'^1_0,u_1)$ we get a similar
relation.

\subsection{The Ivory affinity provides a rigid motion}\label{subsec:algpre24}
\noindent

\noindent
From the fact that the Ivory affinity preserves lengths and angles
of segments between (and rulings on) confocal quadrics, we shall
find that if $x_0^0(u_0,v_0),\ x_0^1(u_1,v_1)$ are on the
quadric $x_0$ and choosing rulings $w_0^j,\ w_z^j$ at each of
the points $x_0^j,\ x_z^j$, then there is an rigid motion
$(R_0^1,t_0^1)$ of $\mathbb{C}^3$ which takes $x_0^0$ to
$x_z^0$, $x_z^1$ to $x_0^1$, $w_0^0$ to $w_z^0$ and $w_z^1$
to $w_0^1$:
\begin{eqnarray}\label{eq:compy2}
(R_0^1,t_0^1)x_0^0=x_z^0,\ (R_0^1,t_0^1)x_z^1=x_0^1,\
R_0^1w_0^0=w_z^0,\  R_0^1w_z^1=w_0^1.
\end{eqnarray}
The first two conditions can be replaced by
$R_0^1V_0^1=-V_1^0,\ V_0^1:=x_z^1-x_0^0,\ V_1^0:=x_z^0-x_0^1$.
\begin{center}
$\xymatrix@!0{&&x_0^0\ar@{-}[drdr]\ar@/_/@{-}[rr]^{x_0}\ar@{~>}[dd]_{(R_0^1,t_0^1)}&&
x_0^1\ar@{<~}[dd]^{(R_0^1,t_0^1)}&&\\
\ar@{-}[urr]^{w_0^0}&&&\ar[dl]^>>>>{V_1^0}\ar[dr]^>>>>>{V_0^1}&&&\ar@{-}[ull]_{w_0^1}\\&&
x_z^0\ar@{-}'[ur][urur]\ar@/_/@{-}[rr]_{x_z}&&x_z^1&\\\ar@{-}[urr]^{w_z^0}&&&
&&&\ar@{-}[ull]_{w_z^1}}$
\end{center}
Such an $R_0^1$ can be found because: $[(x_z^1-x_0^0)\ \ w_0^0\
\ w_z^1]^T[(x_z^1-x_0^0)\ \ w_0^0\ \ w_z^1]=\\=\begin{bmatrix}
|x_z^1-x_0^0|^2&(w_0^0)^T(x_z^1-x_0^0)&(w_z^1)^T(x_z^1-x_0^0)\\
(w_0^0)^T(x_z^1-x_0^0)&|w_0^0|^2&(w_0^0)^Tw_z^1\\
(w_z^1)^T(x_z^1-x_0^0)&(w_0^0)^Tw_z^1&|w_z^1|^2
\end{bmatrix}=[(x_0^1-x_z^0)\ \ w_z^0\ \ w_0^1]^T[(x_0^1-x_z^0)\ \ w_z^0\ \ w_0^1]$.
To fix the ideas we shall henceforth make a choice:
$w_0^j:=x_{0u_j}^j,\ j=0,1$; the rotations (in the case of the TC
$(V_0^1)^T\hat N_0^0=0$) when changing the ruling family on
$x_0^1$, on $x_0^0$ or on both are obtained, according to the
second relation of (\ref{eq:simpas}), by first reflecting in
$T_{x_0^0}x_0$ or further in $T_{x_0^1}x_0$ or both; conversely,
the second relation of (\ref{eq:simpas}) is a consequence of the
existence of the RMPIA. Because of this
$x_{0v_1}^1=R_0^1(I_3-2N_0^0(N_0^0)^T)x_{zv_1}^1$; multiplying
this on the left with $(x_{0u_1}^1)^T$ we get the first relation
of (\ref{eq:simpas}). Note that $(R_0^1,t_0^1)$ depends only on
$(v_0,v_1)$: with $F:=[\mathcal{B}_0x_{0u_0}^0\ \ \
\mathcal{B}_1x_{zu_1}^1\ \ \ (x_0^0(0,v_0)-x_z^1(0,v_1))]$ we
have:
\begin{eqnarray}\label{eq:r12c}
R_0^1=((0\leftrightarrow z)\circ F)F^{-1};\ t_0^1=x_z^0(0,v_0)-R_0^1x_0^0(0,v_0).
\end{eqnarray}
Let $\Del^-=\Del^-(z,v_0,v_1):=-\mathcal{B}_0(m_0^1)^Tx_{0u_0}^0,\
\Del^+=\Del^+(z,u_0,v_1) :=\mathcal{B}_0(m_0^1)^Tx_{0v_0}^0$;
similarly we can consider $\Del'^{\pm}$ by changing the ruling
family on $x_z^1$ (replace $(m_0^1,v_1)$ with $(m'^1_0,u_1)$).

We have
$\frac{-\Del^-}{4\mathcal{A}}=\frac{1}{4}Y(v_0)^T(\det\sqrt{R_z}(\sqrt{R_z})^{-1}
+\sqrt{R_z})Y(v_1)$ for QC,
$=-i((f_1+v_0e_3)\times(\sqrt{R'_z}f_1+v_{z1}e_3))^T(v_{z1}\sqrt{R'_z}\bar
f_1 -v_0\bar
f_1+z(1-\frac{|f|^2}{2})(L^TL)^{-1}e_3)=(1+\det\sqrt{R'_z})v_0v_1-
\bar
f_1^T\sqrt{R'_z}f_1((v_0+v_1+v(z))^2-2v_0v_1)-zf_1^T(L^TL)^{-1}e_3
(v_0+v_1+v(z))+z\frac{|f|^2}{2}f_1^T\sqrt{R'_z}f_1=v_0v_1(1+\det\sqrt{R'_z})
-(v_0^2+v_1^2)\bar f_1^T\sqrt{R'_z}f_1
+\frac{z}{2}f_1^T\sqrt{R'_z}f_1$ for QWC,
$=-\frac{1+\det\sqrt{R'_z}}{2}(v_0-v_1)^2+v(z)(1-\det\sqrt{R'_z})
(v_0+v_1)+\frac{v(z)^2}{2}(3-\det\sqrt{R'_z})$ for IQWC
(\ref{eq:iqwc1}), $=-(v_0-v_1)^2-2v(z)(v_0+v_1)-v(z)^2$ for IQWC
(\ref{eq:iqwc2}). We have $\det F=\Del^-=(v_0\leftrightarrow
v_1)\circ\Del^-= \det((v_0\leftrightarrow v_1)\circ
F)=\det((0\leftrightarrow z)\circ F)$, so $\det R_0^1=1$.

Differentiating $(V_0^1)^T\hat N_0^0=(V_1^0)^T\hat N_0^1$ wrt
$u_0$ we get:
\begin{eqnarray}\label{eq:dV12}
(V_0^1)^T\hat N_{0u_0}^0=(x_{zu_0}^0)^T\hat N_0^1.
\end{eqnarray}
Consider now the TC $(V_0^1)^T\hat N_0^0=0$; because
$R_0^1\frac{m_0^1}{|m_0^1|}= \pm\frac{\hat N_0^1}{|\hat N_0^1|}$
we have $R_0^1V_0^1 =\mathcal{A}\hat N_0^1\times R_0^1\hat N_0^0$
(multiply it on the left with $(R_0^1(x_{0u_0}^0\times x_0^0))^T$,
use (\ref{eq:dV12}), $x_{0u}\times x_{0}= \mathcal{A}\hat N_{0u},\
\hat N_0^Tx_0=1$ for QC and $x_0^T\hat N_0=uv,\ x_{0u}\times x_0=
-v\mathcal{A}(\hat N_0-u\hat N_{0u})$ for(I)QWC). Thus:
\begin{eqnarray}\label{eq:V21c}
R_0^1V_0^1=\mathcal{A}\hat N_0^1\times R_0^1\hat N_0^0,\
R_0^1m_0^1=\mathcal{A}\mathcal{B}_1\hat N_0^1(\hat N_0^0)^Tx_{zu_1}^1.
\end{eqnarray}
If we change the ruling family on $x_0^0$ or on $x_z^1$ or on
both, then the sign of the rhs changes in the second or the first
or both relations of (\ref{eq:V21c}). Using (\ref{eq:simpas}) and
(\ref{eq:V21c}) we thus have:
\begin{eqnarray}\label{eq:deltasi2}
\Del^-=-\mathcal{A}\mathcal{B}_0\mathcal{B}_1(x_{zu_0}^0)^T
\hat N_0^1(\hat N_0^0)^Tx_{zu_1}^1,\nonumber\\
\Del^+=-\mathcal{A}\mathcal{B}_0\mathcal{B}_1(x_{zv_0}^0)^T
\hat N_0^1(\hat N_0^0)^Tx_{zu_1}^1,\nonumber\\
\Del^+\Del^-=-z\mathcal{A}^2\mathcal{B}_0\mathcal{B}_1^2
|\hat N_0^1|^2((x_{zu_1}^1)^T\hat N_0^0)^2.
\end{eqnarray}
By changing the ruling family on $x_z^1$ we get similar formulae
for $\Del'^\pm$ (the sign of the rhs does not change),
$\Del'^+\Del'^-,\ \Del^-\Del'^-,\ \Del^+\Del'^+$. With
$\mathcal{N}:=\frac{1}{\mathcal{A}\mathcal{B}|\hat N_0|^2}$ the
first relation of (\ref{eq:simpas}) and (\ref{eq:deltasi2}) imply
\begin{eqnarray}\label{eq:del+-+-}
\Del^-\Del'^+=\frac{z^2}{\mathcal{N}_0\mathcal{N}_1}=\Del^+\Del'^-;
\end{eqnarray}
conversely $\Del^-\Del'^+-\Del^+\Del'^-=0$ is equivalent to
$0=\frac{(m^1_0)^Tx_{0v_0}^0(m'^1_0)^Tx_{0u_0}^0
-(m^1_0)^Tx_{0u_0}^0(m'^1_0)^Tx_{0v_0}^0}{(N_0^0)^T(x_{0u_0}^0\times
x_{0v_0}^0)} =(m'^1_0\times m^1_0)^TN_0^0=-\mathcal{B}_1^2
(N_z^1)^T(x_{zu_1}^1\times
x_{zv_1}^1)(V_0^1)^TN_z^1(V_0^1)^TN_0^0$, so it is equivalent
either to the TC or to $(V_0^1)^TN_z^1=0$.

Suppose $\{m_j\}_{j=1,2,3}$ are functions of the independent
variables $\{y_k\}_{k=1,...,4}$ (that is $dy_1\wedge...\wedge
dy_4\neq 0,\ dm_j\wedge dy_1\wedge...\wedge dy_4=0$). If we can
restrict $y_4$ to $dy_1\wedge dy_2\wedge dy_3\neq 0,\
dy_4=\sum_{j=1}^3m_jdy_j$, then imposing the compatibility
condition $d\wedge$ to the last relation we get
\begin{eqnarray}\label{eq:wet}
\pa_{y_4}\log\frac{m_{j+1}}{m_j}+\frac{1}{m_j}\pa_{y_j}\log m_{j+1}
-\frac{1}{m_{j+1}}\pa_{y_{j+1}}\log m_j=0,\ j=\circlearrowleft_1^3.
\end{eqnarray}

Applying $d$ to the TC $(V_0^1)^T\hat N_0^0=0$, using the first
relation of (\ref{eq:simpas}), (\ref{eq:dV12}) and
(\ref{eq:deltasi2}) we get
\begin{eqnarray}\label{eq:du1}
du_1=-\frac{\mathcal{N}_0}{z}(\Del'^-du_0+\Del'^+dv_0)-\frac{\Del'^+}{\Del^+}dv_1.
\end{eqnarray}
Applying (\ref{eq:wet}) to
$(y_1,y_2,y_3,y_4,m_1,m_2,m_3):=(u_0,v_0,v_1,u_1,
-\frac{\mathcal{N}_0\Del'^-}{z},-\frac{\mathcal{N}_0\Del'^+}{z},-\frac{\Del'^+}{\Del^+}),\
j=1,3$ and to
$(y_1,y_2,y_3,y_4,m_1,m_2,m_3):=(v_0,u_0,v_1,u_1,-\frac{z}{\mathcal{N}_1\Del^-},
-\frac{\mathcal{N}_0\Del'^-}{z},-\frac{\Del'^-}{\Del^-}),\ j=1$ we
get the algebraic relations
\begin{eqnarray}\label{eq:alg}
\frac{\mathcal{N}_0}{z}\pa_{u_1}\log\frac{\Del'^+}{\Del'^-}=
\frac{1}{\Del'^-}\pa_{u_0}\log(\mathcal{N}_0\Del'^+)-
\frac{1}{\Del'^+}\pa_{v_0}\log(\mathcal{N}_0\Del'^-)
=\frac{1}{\Del'^-}\pa_{u_0}\log\frac{\Del'^+}{\Del^+}=\nonumber\\
\frac{\mathcal{N}_0}{z}\pa_{u_1}\log(\mathcal{N}_1\Del'^+)
-\frac{1}{\Del'^+}\pa_{v_0}\log(\mathcal{N}_0\Del'^-).
\end{eqnarray}
Note that while explicit derivatives $\pa_{\cdot}$ obey usual
calculus rules, one must avoid the use of consequences of the TC
inside the $\pa_{\cdot}$ sign; for example
$\pa_{u_0}\log\frac{\Del^-}{\Del'^-}=0\neq\pa_{u_0}\log\frac{\Del^+}{\Del'^+}$.
Since the first relation of (\ref{eq:alg}) does not depend on
$v_1$, it is valid for independent $u_0,v_0,u_1$.

Note
$((u_0,v_0)\leftrightarrow(u_1,v_1))\circ(\Del^-,\Del^+,\Del'^-,\Del'^+)
=(\Del^-,\Del'^-,\Del^+,\Del'^+)$; using (\ref{eq:del+-+-}) the
relation (\ref{eq:du1}) (and thus (\ref{eq:alg})) has the
symmetries $(u_0,v_0,\Del^+,\Del'^-,\mathcal{N}_0)\leftrightarrow
(u_1,v_1,\Del'^-,\Del^+,\mathcal{N}_1),\\
(u_0,\Del^+,\Del'^+)\leftrightarrow(v_0,\Del^-,\Del'^-),\
(u_1,\Del'^-,\Del'^+)\leftrightarrow(v_1,\Del^-,\Del^+)$. Under
these symmetries and compositions thereof, there are no other
algebraic identities besides (\ref{eq:alg}) and its immediate
consequences (like
$\frac{1}{\Del'^-}\pa_{u_0}\log\frac{\Del'^+}{\Del^+}=
-\frac{1}{\Del'^+}\pa_{v_0}\log\frac{\Del'^-}{\Del^-}$) which
appear from (\ref{eq:du1}) and its equivalent versions.

The TC has also other symmetries: for QC $u_0\leftrightarrow v_0\
(u_1\leftrightarrow v_1)$ and $\sqrt{R_z}\leftrightarrow
-\sqrt{R_z}$ (which reflect $x_0^0 (x_0^1)$ in the origin and
restore the TC: $\sqrt{R_z}X(u,v)=-\sqrt{R_z}X(v,u)$); for QWC
(\ref{eq:qwc1}) and IQWC (\ref{eq:iqwc1}) $u_0\leftrightarrow
v_0,\ (u_1\leftrightarrow v_1)$ and $\sqrt{R'_z}\leftrightarrow
r_{e_2}\sqrt{R'_z} =\sqrt{R'_z}r_{e_2},\ r_{e_2}:=f_1f_1^T+\bar
f_1\bar f_1^T+e_3e_3^T =I_3-2e_2e_2^T$ (which reflect $x_0^0\
(x_0^1)$ in the plane $e_2^Tx=0$ for QWC (\ref{eq:qwc1}) or in the
plane $e_3^Tx=0$ for IQWC (\ref{eq:iqwc1}) and restore the TC:
$\sqrt{R'_z}I_{1,2}Z(u,v)=r_{e_2}\sqrt{R'_z}I_{1,2}Z(v,u)$).

For QC
$\Del^+=((u_0,\sqrt{R_z})\leftrightarrow(v_0,-\sqrt{R_z}))\circ\Del^-,\
\Del'^+=((u_0,u_1)\leftrightarrow(v_0,v_1))\circ\Del^-$; for QWC
(\ref{eq:qwc1}) and IQWC (\ref{eq:iqwc1})
$\Del^+=((u_0,\sqrt{R'_z})\leftrightarrow(v_0,r_{e_2}\sqrt{R'_z}))\circ\Del^-,\
\Del'^+=((u_0,u_1)\leftrightarrow(v_0,v_1))\circ\Del^-$. For QWC
(\ref{eq:qwc2}) and IQWC (\ref{eq:iqwc2}) the change
$\sqrt{R'_z}\leftrightarrow\sqrt{R'_z}r_{e_2}$ for changing the
ruling family on $x_0^1$ or $f_1\leftrightarrow r_{e_2}f_1=\bar
f_1$ for changing the ruling family on $x_0^0$ still has some
effect similar to QWC (\ref{eq:qwc1}) and IQWC (\ref{eq:iqwc1}),
but since $\sqrt{R'_z}r_{e_2}\neq r_{e_2}\sqrt{R'_z}$ new formulae
appear:
$\frac{-\Del^+}{4\mathcal{A}}=v_1(u_0(1-\det\sqrt{R'_z})-v_1\bar
f_1^T\sqrt{R'_z}\bar f_1) +\frac{z}{2}\bar f_1^T\sqrt{R'_z}f_1,\
\frac{-\Del'^+}{4\mathcal{A}}=u_0u_1(1+\det\sqrt{R'_z})
-(u_0^2+u_1^2)\bar f_1^T\sqrt{R'_z}f_1+\frac{z}{2}\bar
f_1^T\sqrt{R'_z}\bar f_1$ for QWC (\ref{eq:qwc2}),
$\frac{-\Del^+}{4\mathcal{A}}=\frac{z}{2}u_0(u_0-z)-zv_1+\frac{z^3}{8},\
\frac{-\Del'^+}{4\mathcal{A}}=\frac{z^2}{4}-(u_0-u_1)^2$ for IQWC
(\ref{eq:iqwc2}).

\subsection{Second iteration of the tangency configuration}\label{subsec:algpre25}
\noindent

\noindent Consider the SITC: $x_{z_1}^0,\ x_{z_2}^3\in
T_{x_0^1}x_0$. In the pencil of planes containing $x_{z_2}^0,\
x_{z_1}^3$ choose a certain one tangent to $x_0$ at $x_0^2$ (in
general there are two choices; for the other choice we must
consider the other ruling families on $x_0^2,\ x_0^3$). To have a
consistent notation we shall order $0<z_1<z_2$ and the position of
$j,k$ in $V_k^j,\ (R_k^j,t_k^j)$ is increasing wrt this order: for
example $V_0^3:=x_{z_2}^3-x_{z_1}^0\in T_{x_0^1}x_0$ since there
is no other incidence relationship involving the upper indices
$0,3$ and lower indices $0,z_1$ or $0,z_2;\
R_3^0(x_{z_2u_0}^0,x_{z_1u_3}^3,V_3^0)=(x_{z_1u_0}^0,x_{z_2u_3}^3,-V_0^3),\
(R_3^0)^{-1}=R_0^3$. Note:
\begin{eqnarray}\label{eq:config}
V_0^3=V_1^3-V_1^0,\ V_3^0=V_2^0-V_2^3,\ V_1^2=V_0^2-V_0^1,\
V_2^1=V_3^1-V_3^2,\ (V_0^1)^T\hat N_0^0 =(V_0^2)^T\hat
N_0^0=\nonumber\\ (V_1^0)^T\hat N_0^1= (V_1^3)^T\hat
N_0^1=(V_2^0)^T\hat N_0^2=(V_2^3)^T\hat N_0^2=(V_3^1)^T\hat
N_0^3=(V_3^2)^T\hat N_0^3=0.
\end{eqnarray}
In particular the SITC has the group of symmetries of the square
generated by the flip $(0\leftrightarrow 1,2\leftrightarrow 3)$
and the rotation $(0\rightarrow 1\rightarrow 3\rightarrow 2
\rightarrow 0,\ z_1\leftrightarrow z_2)$ (consider the square with
vertices $0,1,2,3$ and sides labelled $z_1,\ z_2$ such that both
the sums of opposite vertices and the labels of opposite sides are
equal; $0,1,2,3$ do not yet respectively correspond, up to some
rigid motions, to
$(R_0^1,t_0^1)x_0^0,x_0^1,(R_3^0,t_3^0)x_0^2,(R_3^1,t_3^1)x_0^3$,
but the symmetries can be seen at the level of (\ref{eq:config})).

The next result (proved by Bianchi in (122,\S\ 26)) plays a
fundamental r\^{o}le in the proof of the BPT:

{\it If $x_{z_1}^0,\ x_{z_2}^3\in T_{x_0^1}x_0$ and for a certain
choice of $x_0^2$ such that $x_{z_2}^0,\ x_{z_1}^3\in
T_{x_0^2}x_0$, then the rulings $x_{0u_1}^1,R_3^0x_{0u_2}^2$ cut
the line $[x_{z_1}^0\ x_{z_2}^3]$ with cross ratio
$\frac{z_1}{z_2}$.}
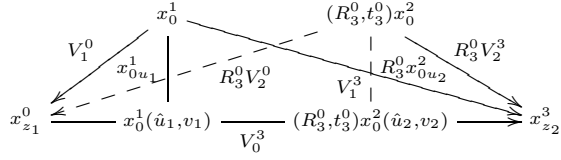
\begin{center}
$\xymatrix{&
^{x_0^1}\ar[dl]_{V_1^0}\ar@{-}[d]_{x_{0u_1}^1}\ar[drr]_<<<<<<<<<<<<<<<<<<<<<<<{V_1^3}&
^{(R_3^0,t_3^0)x_0^2}\ar@{-->}[dll]^<<<<<<<<<<<<{R_3^0V_2^0}\ar@{--}[d]^{R_3^0x_{0u_2}^2}
\ar[dr]^{R_3^0V_2^3}&\\^{x_{z_1}^0}\ar@{-}[r]&^{x_0^1(\hat
u_1,v_1)}\ar@{-}[r]_<<<<<{V_0^3}&^{(R_3^0,t_3^0)x_0^2(\hat
u_2,v_2)}\ar[r]&^{x_{z_2}^3}}$
\end{center}
From (\ref{eq:u2}) we get:
$\frac{Y(v_1)^T\sqrt{R_{z_1}}X_0}{1-\frac{1}{2}Y'(v_1)^T\sqrt{R_{z_1}}X_0}=u_1-v_1=
\frac{Y(v_1)^T\sqrt{R_{z_2}}X_3}{1-\frac{1}{2}Y'(v_1)^T\sqrt{R_{z_2}}X_3}$
for QC, $\frac{v_1x_{z_1v}(0,0)^T\hat N_0^0-f^Tx_{z_1}^0}{v_1-
x_{z_1u}(0,0)^T\hat N_0^0}=u_1=\frac{v_1x_{z_2v}(0,0)^T\hat
N_0^3-f^Tx_{z_2}^3} {v_1-x_{z_2u}(0,0)^T\hat N_0^3}$ for (I)QWC
and other similar equations, obtained by $(0\leftrightarrow
3,1\leftrightarrow 2)$. These become:
\begin{eqnarray}\label{eq:v1v2}
\frac{1}{2}Y(v_1)^T(\sqrt{R_{z_2}}X_3-\sqrt{R_{z_1}}X_0-i(\sqrt{R_{z_2}}X_3)
\times(\sqrt{R_{z_1}}X_0))=0\ \mathrm{for\ QC},\nonumber\\
v_1^2(x_{z_1v}(0,0)^T\hat N_0^0-x_{z_2v}(0,0)^T\hat N_0^3)+v_1(f^TV_0^3
-x_{z_1v}(0,0)^T\hat N_0^0x_{z_2u}(0,0)^T\hat N_0^3+\nonumber\\
x_{z_1u}(0,0)^T\hat N_0^0x_{z_2v}(0,0)^T\hat
N_0^3)-x_{z_1u}(0,0)^T\hat N_0^0f^Tx_{z_2}^3 +x_{z_2u}(0,0)^T\hat
N_0^3f^Tx_{z_1}^0=0\nonumber\\ \mathrm{for\ (I)QWC}.
\end{eqnarray}
When $v_j,\ j=0,1,3$ are fixed $u_0,\ u_3$ vary preserving
(\ref{eq:v1v2}), so a homography (invariant under
$(0\leftrightarrow 3,1\leftrightarrow 2)$) is established between
the rulings $x_{z_1u_0}^0,x_{z_2u_3}^3$. The quadric $Q^{301}\
(Q^{032})$ generated by $x_{z_2u_3}^3,x_{z_1u_0}^0,x_{0u_1}^1\
(x_{z_2u_0}^0,x_{z_1u_3}^3,x_{0u_2}^2)$ cuts $x_0$ along
$x_{0u_1}^1\ (x_{0u_2}^2)$ (and possibly other sets we exclude
from our discussion). The points on $x_{z_1u_0}^0,x_{z_2u_3}^3\
(x_{z_2u_0}^0,x_{z_1u_3}^3)$ corresponding under the homography
lie on the other ruling family of $Q^{301}\ (Q^{032})$, so
$(R_0^3,t_0^3)Q^{301} =Q^{032}$ (if a line $l_1$ intersects two
non-co-planar lines $l_2,l_3$, then there is a $2$-dimensional
family of quadrics having $l_2,l_3$ in a ruling family and $l_1$
in the other ruling family; each quadric is uniquely determined by
the homography it establishes between $l_2,l_3$ and it corresponds
to the choices of an arbitrary line $l_4$ at a fixed point on
$l_1$; $l_4$ will belong to the first ruling family), which is
equivalent to the rulings $x_{z_1u_0}^0,\ x_{z_2u_3}^3$ cutting
the line $[x_0^1\ (R_3^0,t_3^0)x_0^2]$ (again with cross-ratio
$\frac{z_1}{z_2}$, because four rulings of a ruling family on a
quadric cut any ruling of the other ruling family with same cross
ratio). This is equivalent to $R_1^0R_3^0V_2^0\in T_{x_0^0}x_0,\
R_1^3R_3^0V_2^3\in T_{x_0^3}x_0$, from which follow the co-cycle
relations
\begin{eqnarray}\label{eq:cocy}
(R_0^1,t_0^1)(R_2^0,t_2^0)=(R_3^0,t_3^0)=(R_3^1,t_3^1)(R_2^3,t_2^3),\
(R_1^0,t_1^0)(R_3^1,t_3^1)=(R_2^1,t_2^1)=(R_2^0,t_2^0)(R_3^2,t_3^2).\nonumber\\
\end{eqnarray}
To prove this let
$(R_0,t_0):=(R_1^0,t_1^0)(R_3^0,t_3^0)(R_0^2,t_0^2),\
(R_1,t_1):=(R_0^1,t_0^1)(R_2^1,t_2^1)(R_1^3,t_1^3),\
(R_2,t_2):=(R_3^2,t_3^2)(R_1^2,t_1^2)(R_2^0,t_2^0),\
(R_3,t_3):=(R_2^3,t_2^3)(R_0^3,t_0^3)(R_3^1,t_3^1)$; note:
$(R_j,t_j)(x_0^j,x_{0u_j}^j)=(x_0^j,x_{0u_j}^j),\ j=0,...,3$, so
$R_j$ are rotations around $x_{0u_j}^j,\ j=0,...,3$. But
$R_0V_0^2\in T_{x_0^0}x_0,\ R_3^{-1}V_3^2\in T_{x_0^3}x_0$; thus
$R_0(T_{x_0^0}x_0)=T_{x_0^0}x_0,\ R_3(T_{x_0^3}x_0)=T_{x_0^3}x_0$;
since $\det R_0=\det R_3=1$ we get $R_0=R_3=I$, so $t_0=t_3=0$.
From the first four relations of (\ref{eq:config}) we get
$-R_0^1(I-R_0)V_0^2+R_3^1(I-R_3^{-1})V_3^2+(I-R_1^{-1})R_0^1V_1^2=0$,
so $R_1=I$ and
$R_2^{-1}=R_0^2R_1^0R_1R_0^1R_0R_2^0R_3^2R_3R_2^3=I$, so
$t_1=t_2=0$. Another way to prove $(R_1,t_1)=(R_2,t_2)=(I_3,0)$
once we know $(R_0,t_0)=(R_3,t_3)=(I_3,0)$ is to apply to the SITC
the symmetry $(0\rightarrow 1\rightarrow 3\rightarrow 2
\rightarrow 0,\ z_1\leftrightarrow z_2)$.

Applying the rigid motions $(R_0^3,t_0^3),\ (R_1^0,t_1^0),\
(R_1^3,t_1^3)$ to the configuration
$\\((R_0^1,t_0^1)(x_0^0,x_{0u_0}^0),(x_0^1,x_{0u_1}^1),(R_3^0,t_3^0)(x_0^2,x_{0u_2}^2),
(R_3^1,t_3^1)(x_0^3,x_{0u_3}^3))$ we get the configurations
$\\((R_0^2,t_0^2)(x_0^0,x_{0u_0}^0),(R_0^3,t_0^3)(x_0^1,x_{0u_1}^1),(x_0^2,x_{0u_2}^2),
(R_3^2,t_3^2)(x_0^3,x_{0u_3}^3))$,
$\\((x_0^0,x_{0u_0}^0),(R_1^0,t_1^0)(x_0^1,x_{0u_1}^1),(R_2^0,t_2^0)(x_0^2,x_{0u_2}^2),
(R_2^1,t_2^1)(x_0^3,x_{0u_3}^3))$,
$\\((R_1^2,t_1^2)(x_0^0,x_{0u_0}^0),(R_1^3,t_1^3)(x_0^1,x_{0u_1}^1),(R_2^3,t_2^3)(x_0^2,x_{0u_2}^2),
(x_0^3,x_{0u_3}^3))$.

These rigid motions realize the symmetries $(0\leftrightarrow
3,1\leftrightarrow 2), (0\leftrightarrow 1,2\leftrightarrow
3),(0\leftrightarrow 2,1\leftrightarrow 3)$; the symmetry
$(0\leftrightarrow 3,z_1\leftrightarrow z_2)$ is realized by the
identity (at the level of the square it is the flip in the
diagonal $[1\ 2]$).

Because of (\ref{eq:n(i)qwc}) the discriminant of the quadratic
equation (\ref{eq:v1v2}) which defines $v_1$ is
$\mathbf{\Del}:=(X_0^T\sqrt{R_{z_1}}\sqrt{R_{z_2}}X_3-1)^2
-z_1z_2|\hat N_0^0|^2|\hat N_0^3|^2$ for QC,
$:=(x_{z_1u}(0,0)^T\hat N_0^0x_{z_2v}(0,0)^T\hat N_0^3+\\
x_{z_1v}(0,0)^T\hat N_0^0x_{z_2u}(0,0)^T\hat
N_0^3-f^T(x_{z_1}^0+x_{z_2}^3))^2 -z_1z_2|\hat N_0^0|^2|\hat
N_0^3|^2=(Z_0^T\sqrt{R'_{z_1}}\sqrt{R'_{z_2}}I_{1,2}Z_3
-(Z_0+Z_3)^T(u(z_1+z_2)f_1+v(z_1+\ep z_2)\bar
f_1+e_3)-(u(z_1)-u(z_2))(v(z_1)-v(z_2))
-\frac{|f|^2}{2}(z_1+z_2))^2-z_1z_2|\hat N_0^0|^2|\hat N_0^3|^2$
for (I)QWC, where $\ep:=-1$ for IQWC (\ref{eq:iqwc2}) and $:=1$
otherwise. Thus $(0\leftrightarrow 3,1\leftrightarrow
2)\circ\mathbf{\Del}=\mathbf{\Del}$ (note also that the first
squared terms are obtained by polarization* from
(\ref{eq:symrel})) and we have: \footnote{* This polarization is
clear at the analytic level from the ansatz $z_1=z_2$, but it may
also have another geometric explanation as in higher dimensions
probably another discriminant must also behave well. For this
reason the computations of Jacobi's higher dimensional version and
its corresponding discretization {\it ray of light} must be
similar: there also the discriminant behaved well because a
solution (namely $0$) of a quadratic equation was known; as a
consequence other solutions of a polynomial of degree $n-1$ with
$0$ among its roots also behaved well. Thus one can consider the
polarization of the next geometric statement: a line tangent to
$x_0$ at $x_0^1$ cuts a confocal quadric $x_z$ at the known point
$x_z^0$ and thus at another computable point $x_z^3$; now one {\it
geometrically} polarizes by {\it wiggling} the points
$x_0^0,x_0^3$ and the $z$ (separate $z_1$'s and $z_2$'s stick to
each of $x_0^0,x_0^3$).}
\begin{eqnarray}\label{eq:v1v21}
\frac{1}{2}Y'(v_1)^T(\sqrt{R_{z_2}}X_3-\sqrt{R_{z_1}}X_0
-i(\sqrt{R_{z_2}}X_3)\times(\sqrt{R_{z_1}}X_0))=\pm\sqrt{\mathbf{\Del}}=\nonumber\\
\pm\frac{1}{2}Y'(v_2)^T(\sqrt{R_{z_2}}X_0-\sqrt{R_{z_1}}X_3
-i(\sqrt{R_{z_2}}X_0)\times(\sqrt{R_{z_1}}X_3))\ \mathrm{for\ QC},\nonumber\\
2v_1(x_{z_1v}(0,0)^T\hat N_0^0-x_{z_2v}(0,0)^T\hat N_0^3)+f^TV_0^3
-x_{z_1v}(0,0)^T\hat N_0^0x_{z_2u}(0,0)^T\hat N_0^3+\nonumber\\
x_{z_1u}(0,0)^T\hat N_0^0x_{z_2v}(0,0)^T\hat
N_0^3=\pm\sqrt{\mathbf{\Del}}=
\pm(2v_2(x_{z_1v}(0,0)^T\hat N_0^3-x_{z_2v}(0,0)^T\hat N_0^0)+f^TV_3^0-\nonumber\\
x_{z_1v}(0,0)^T\hat N_0^3x_{z_2u}(0,0)^T\hat N_0^0+
x_{z_1u}(0,0)^T\hat N_0^3x_{z_2v}(0,0)^T\hat N_0^0)\ \mathrm{for\
(I)QWC}.
\end{eqnarray}
For QC (\ref{eq:v1v2}) has the symmetry $(u_j\leftrightarrow v_j,\
j=0,1,3)$ (since reflection in the origin preserves the SITC), so
the lhs of (\ref{eq:v1v21}) is skew wrt $(u_j\leftrightarrow v_j,\
j=0,1,3)$; for (I)QWC the lhs of (\ref{eq:v1v21}) is skew wrt
$u_1\leftrightarrow v_1$ (use $x_{z_ju}(0,0)\leftrightarrow
x_{z_jv}(0,0),\ j=1,2$ and (\ref{eq:nor})).

We first choose the $-$ sign in the rhs of (\ref{eq:v1v21}) (this
is the choice of $x_0^2$). Replacing $u_0,\ u_3$ from
(\ref{eq:u2}), (\ref{eq:v1v21}) becomes a quadratic polynomial of
$u_1$ being identically $0$; on account of the homography between
the rulings $x_{z_1u_0}^0,\  x_{z_2u_3}^3$ it is enough to give
$u_1$ a particular value (the relation thus found between $v_j,\
j=0,...,3$ will guarantee that the required quadratic polynomial
in $u_1$ is identically $0$). If we let $u_1:=\infty$, then we
obtain for all quadrics except (\ref{eq:qc1}) with $a_1\neq a_2$ a
homography (as a polynomial it has degree $2$) between $v_j,\
j=0,...,3:\ \phi=\phi(z_1,z_2,v_0,v_1,v_2,v_3)=0$ having the
required symmetries of the SITC; for (\ref{eq:qc1}) with $a_1\neq
a_2$ we obtain a polynomial of degree $6:\
\varphi=\varphi(z_1,z_2,v_0,v_1,v_2,v_3)=0$ and missing some
symmetries. Symmetrizing $\varphi$ we get a homography (as a
polynomial it has degree $4$) with all the required symmetries
($\varphi':=\varphi-\varphi\circ(v_0\leftrightarrow
v_3,v_1\leftrightarrow v_2)$ is still valid and gains the symmetry
$(v_0\leftrightarrow v_3,v_1\leftrightarrow v_2)$; finally
$\phi:=\frac{\varphi'-\varphi'\circ(v_0\leftrightarrow
v_1,v_2\leftrightarrow v_3)}{v_0v_2-v_1v_3}$ is a homography with
all the required symmetries). Now one can easily factor $\varphi$
as the product of $\phi$ and a polynomial of degree $2$ (linear in
$v_0,v_1,v_3$) by accounting the terms containing $v_0v_1v_2v_3$;
the remaining terms also satisfy it.

The symmetrization process to find (\ref{eq:phi2c}) for QC
(\ref{eq:qc1}) with $a_1\neq a_2$ is complicated; once
(\ref{eq:phi2c}) is known, the factorization of $\varphi$ as above
is easily achieved. Therefore we propose a change in Bianchi's
approach: while he first completes the discussion of the SITC \&
BPT and then proves the cross-ratio property as a corollary, we
insert the cross-ratio property in the discussion of the SITC,
thus easily arriving at (\ref{eq:phi2c})-(\ref{eq:phi2i1});
further the above approach by choosing $u_1:=\infty$ (and
factorization of $\varphi$ for QC (\ref{eq:qc1}) with $a_1\neq
a_2$) turns out to be unnecessary.

Consider the cross ratio property: the ruling $x_{0u_1}^1$ cuts
the line $[x_{z_1}^0\ x_{z_2}^3]$ at $x_0(\hat
u_1,v_1):=x_0^1+\frac{\hat u_1-u_1}{\hat
u_1-v_1}(u_1-v_1)x_{0u_1}^1$ for QC or $x_0^1+(\hat
u_1-u_1)x_{0u_1}^1$ for (I)QWC; therefore $(x_0(\hat
u_1,v_1)-x_{z_1}^0)\times(x_{z_1}^0-x_{z_2}^3)=0$, or multiplying
on the left with $(x_{z_1}^0)^T:\ x_0(\hat
u_1,v_1)^T(x_{z_1}^0\times x_{z_2}^3)=0$. Thus $x_0(\hat
u_1,v_1)=\frac{(x_{0u_1}^1\times x_0^1)\times(x_{z_1}^0\times
x_{z_2}^3)} {(x_{0u_1}^1)^T(x_{z_1}^0\times x_{z_2}^3)}$ and
$(x_0(\hat u_1,v_1)-x_{z_1}^0)\div (x_0(\hat
u_1,v_1)-x_{z_2}^3)=\frac{(x_{z_1}^0\times
x_{z_2}^3)^T(x_{0u_1}^1\times V_1^0)} {(x_{z_1}^0\times
x_{z_2}^3)^T(x_{0u_1}^1\times V_1^3)}=\frac{(\hat
N_0^1)^T(x_{0u_1}^1\times V_1^0)} {(\hat N_0^1)^T(x_{0u_1}^1\times
V_1^3)}=^{(\ref{eq:V21c})}\frac{(x_{z_1u_1}^1)^T\hat N_0^0}
{(x_{z_2u_1}^1)^T\hat
N_0^3}(=^{(\ref{eq:u2})}\frac{1-\frac{1}{2}Y'(v_1)^T\sqrt{R_{z_1}}X_0}
{1-\frac{1}{2}Y'(v_1)^T\sqrt{R_{z_2}}X_3}$ for QC). We thus need:
\begin{eqnarray}\label{eq:cross}
\frac{(1-\frac{1}{2}Y'(v_1)^T\sqrt{R_{z_1}}X_0)(1-\frac{1}{2}Y'(v_2)^T\sqrt{R_{z_1}}X_3)}
{(1-\frac{1}{2}Y'(v_1)^T\sqrt{R_{z_2}}X_3)(1-\frac{1}{2}Y'(v_2)^T\sqrt{R_{z_2}}X_0)}=
\frac{z_1}{z_2}\ \mathrm{for\ QC},\nonumber\\
\frac{(x_{z_1u_1}^1)^T\hat N_0^0(x_{z_1u_2}^2)^T\hat N_0^3}
{(x_{z_2u_1}^1)^T\hat N_0^3(x_{z_2u_2}^2)^T\hat N_0^0}=\frac{z_1}{z_2}\ \mathrm{for\ (I)QWC}.
\end{eqnarray}

We have to prove that the cross-ratio property (\ref{eq:cross}) or
its counterpart obtained by choosing the other ruling family on
$x_0^2$ is equivalent to the SITC. Let
$\varphi_+:=\frac{(x_{z_1u_1}^1)^T\hat N_0^0(x_{z_1u_2}^2)^T\hat
N_0^3} {(x_{z_2u_1}^1)^T\hat N_0^3(x_{z_2u_2}^2)^T\hat N_0^0},\
\varphi^{\pm}_+:=\frac{(x_{z_1u_1}^1)^T\hat
N_0^0}{(x_{z_2u_1}^1)^T\hat N_0^3}\pm \frac{(x_{z_1v_1}^1)^T\hat
N_0^0}{(x_{z_2v_1}^1)^T\hat N_0^3},\
\varphi_-:=\frac{(x_{z_1u_1}^1)^T\hat N_0^0(x_{z_1v_2}^2)^T\hat
N_0^3} {(x_{z_2u_1}^1)^T\hat N_0^3(x_{z_2v_2}^2)^T\hat N_0^0},\
\varphi^{\pm}_-:=(0\leftrightarrow 3,1\leftrightarrow
2)\circ\varphi^{\pm}_+$. By (\ref{eq:simpas}) we have
$\frac{(x_{z_1v_1}^1)^T\hat N_0^0(x_{z_1v_2}^2)^T\hat N_0^3}
{(x_{z_2v_1}^1)^T\hat N_0^3(x_{z_2v_2}^2)^T\hat
N_0^0}=\frac{z_1^2}{z_2^2\varphi_+},\ \frac{(x_{z_1v_1}^1)^T\hat
N_0^0(x_{z_1u_2}^2)^T\hat N_0^3} {(x_{z_2v_1}^1)^T\hat
N_0^3(x_{z_2u_2}^2)^T\hat N_0^0}=\frac{z_1^2}{z_2^2\varphi_-}$, so
we need to prove
$\frac{2z_1}{z_2}=\frac{(\varphi^+_+\varphi^+_-)\pm(\varphi^-_+\varphi^-_-)}{2}
(=\varphi_++\frac{z_1^2}{z_2^2\varphi_+}$ for the $+$ sign or
$\varphi_-+\frac{z_1^2}{z_2^2\varphi_-}$ for the $-$ sign). For QC
we have $\frac{-z_2|\hat
N_0^3|^2\varphi^+_+}{2}=X_0^T\sqrt{R_{z_1}}\sqrt{R_{z_2}}X_3-1=
\frac{-z_2|\hat N_0^0|^2\varphi^+_-}{2},\ \frac{z_2|\hat
N_0^3|^2\varphi^-_+}{2}=
-iX_1^T((\sqrt{R_{z_2}}X_3)\times\sqrt{R_{z_1}}X_0)=$lhs of
(\ref{eq:v1v21}), $\frac{z_2|\hat
N_0^0|^2\varphi^-_-}{2}=-iX_2^T((\sqrt{R_{z_2}}X_0)\times\sqrt{R_{z_1}}X_3)=$rhs
of (\ref{eq:v1v21}) with the $+$ sign; for (I)QWC and using
(\ref{eq:nor}) we have $\frac{-z_2|\hat
N_0^3|^2\varphi^+_+}{2}=x_{z_1u}(0,0)^T\hat
N_0^0x_{z_2v}(0,0)^T\hat N_0^3+ x_{z_1v}(0,0)^T\hat
N_0^0x_{z_2u}(0,0)^T\hat N_0^3
-f^T(x_{z_1}^0+x_{z_2}^3)=\frac{-z_2|\hat
N_0^0|^2\varphi^+_-}{2},\ \frac{z_2|\hat
N_0^3|^2\varphi^-_+}{2}=$lhs of (\ref{eq:v1v21}), $\frac{z_2|\hat
N_0^0|^2\varphi^-_-}{2}=$rhs of (\ref{eq:v1v21}) with the $+$
sign, so $\varphi^+_+\varphi^+_-=\frac{4\mathbf{\Del}}{z_2^2|\hat
N_0^0|^2|\hat N_0^3|^2}+\frac{4z_1}{z_2},\
\varphi^-_+\varphi^-_-=\pm\frac{4\mathbf{\Del}}{z_2^2|\hat
N_0^0|^2|\hat N_0^3|^2}$, where $\mathbf{\Del}$ is the
discriminant of the quadratic equation (\ref{eq:v1v2}) which
defines $v_1$.

Note that (\ref{eq:v1v21}) can be interpreted as a homography
between $u_0,u_3:\ C_1u_0u_3+C_2u_0+C_3u_3+C_4=0,\
C_j=C_j(z_1,z_2,v_0,v_1,v_2,v_3),\ j=1,...,4$ and (\ref{eq:cross})
as the same homography between $u_0,u_3:\
D_1u_0u_3+D_2u_0+D_3u_3+D_4=0,\ D_j=D_j(z_1,z_2,v_0,v_1,v_2,v_3),\
j=1,...,4$; thus $C_j,D_j$ must be proportional: $D_j=\eta C_j,\
j=1,...,4,\ \eta=(1\leftrightarrow 2, z_1\leftrightarrow
z_2)\circ\eta$; this condition establishes the sought homographies
between $v_j,\ j=0,...,3$. For QC accounting the coefficients of
$u_1,\ u_2$ from $\frac{D_2}{C_2}=\frac{D_3}{C_3}$ we get
$\eta=\frac{v_0v_3z_1f_1^T\sqrt{R_{z_2}}f_1(1+e_3^T\sqrt{R_{z_2}}e_3)-
2z_2\bar f_1^T\sqrt{R_{z_1}}f_1(1-e_3^T\sqrt{R_{z_1}}e_3)}
{v_0v_3f_1^T\sqrt{R_{z_2}}f_1(1-e_3^T\sqrt{R_{z_1}}e_3) -2\bar
f_1^T\sqrt{R_{z_1}}f_1(1+e_3^T\sqrt{R_{z_2}}e_3)}(=\frac{1-e_3^T\sqrt{R_{z_1}}e_3}
{1+e_3^T\sqrt{R_{z_2}}e_3}z_2$ for QC
(\ref{eq:qc2}),(\ref{eq:qc3})); for all (I)QWC except QWC
(\ref{eq:qwc2})
$\eta=\frac{D_1}{C_1}=\frac{z_1f_1^T\sqrt{R'_{z_2}}f_1} {2\bar
f_1^T\sqrt{R'_{z_1}}f_1}$; for QWC (\ref{eq:qwc2}) the homography
between $v_j,\ j=0,...,3$ is obtained directly from
(\ref{eq:cross}); taking $u_1=\infty$ in (\ref{eq:v1v21}) we
obtain the same homography times a constant multiple of $v_1$.

We shall use the notation $f(z_1\wedge z_2)$ for $f(z_1,z_2)$ if
$f(z_2,z_1)=-f(z_1,z_2)$ and $f(z_1\odot z_2)$ for $f(z_1,z_2)$ if
$f(z_2,z_1)=f(z_1,z_2)$.

For QC let
$M_2(z):=z^{-1}((\mathrm{tr}\sqrt{R_z})I_3-\sqrt{R_z})(I_3+\sqrt{R_z})
=z^{-1}(\mathrm{tr}(\det\sqrt{R_z}(\sqrt{R_z})^{-1}+\sqrt{R_z})I_3
-(\det\sqrt{R_z}(\sqrt{R_z})^{-1}+\sqrt{R_z})),\ M_3(z_1\wedge
z_2):=((\mathrm{tr}(\sqrt{R_{z_1}})I_3-\sqrt{R_{z_1}})(\mathrm{tr}
(\sqrt{R_{z_2}})I_3-\sqrt{R_{z_2}})-\\-\mathrm{tr}(\sqrt{R_{z_1}}\sqrt{R_{z_2}})I_3
+\sqrt{R_{z_1}}\sqrt{R_{z_2}})(\sqrt{R_{z_2}}-\sqrt{R_{z_1}})^{-1},\
m_2(z):=e_3^TM_2(z)e_3,\ m_3(z_1\wedge
z_2):=e_3^TA^{-1}M_3(z_1\wedge z_2)e_3(=(\frac{1}{m_2(z_1)}
-\frac{1}{m_2(z_2)})^{-1}$ if $\mathrm{spec}(A)=\{a_1^{-1}\}$),
$a:=\bar f_1^TA\bar f_1e_3^TAe_3$ for QC
(\ref{eq:qc1}),(\ref{eq:qc2}), $:=a_1^{-3}$ for QC (\ref{eq:qc3}),
$\del:=1$ for QC (\ref{eq:qc1}) and $:=0$ otherwise.

For (I)QWC let $M_3(z_1\wedge
z_2):=2((\mathrm{tr}(\sqrt{R'_{z_1}})I_3-\sqrt{R'_{z_1}})(\mathrm{tr}
(\sqrt{R'_{z_2}})I_3-\sqrt{R'_{z_2}})-\mathrm{tr}(\sqrt{R'_{z_1}}\sqrt{R'_{z_2}})I_3
+\sqrt{R'_{z_1}}\sqrt{R'_{z_2}})(z_1-z_2)^{-1},\
m_2(z):=2z^{-1}(\mathrm{tr}(\sqrt{R'_z})-1),\ m_3(z_1\wedge z_2):=
e_3^TM_3(z_1\wedge
z_2)e_3(=(\frac{1}{m_2(z_1)}-\frac{1}{m_2(z_2)})^{-1}$ for IQWC),
$a:=f_1^TA'f_1,\ k:=\frac{-1}{2\sqrt{2}}$ for IQWC
(\ref{eq:iqwc1}).

The homography between $v_j,\ j=0,...,3$ is:
\begin{eqnarray}\label{eq:phi2c}
\frac{v_0v_3+v_1v_2}{m_3(z_1\wedge z_2)}
+\frac{v_0v_2+v_1v_3}{m_2(z_2)}-\frac{v_0v_1+v_2v_3}{m_2(z_1)}
+a\frac{4+\del v_0v_1v_2v_3}{m_2(z_1)m_2(z_2)m_3(z_1\wedge z_2)}=0\nonumber\\ \mathrm{for\ QC\
(\ref{eq:qc1}),(\ref{eq:qc2}),\ QWC\ (\ref{eq:qwc1}),(\ref{eq:qwc2})};
\end{eqnarray}
\begin{eqnarray}\label{eq:phi2c3}
\frac{v_0v_3+v_1v_2}{m_3(z_1\wedge z_2)}
+\frac{v_0v_2+v_1v_3}{m_2(z_2)}-\frac{v_0v_1+v_2v_3}{m_2(z_1)}
-2a\frac{v_0+v_1+v_2+v_3}{m_2(z_1)m_2(z_2)m_3(z_1\wedge z_2)}
+(\frac{1}{m_3(z_1\wedge z_2)^2}+\nonumber\\
\frac{1}{m_2(z_1)m_2(z_2)})\frac{4a^2}{m_2(z_1)m_2(z_2)m_3(z_1\wedge
z_2)}=0\ \mathrm{for\ QC\ (\ref{eq:qc3}),\ IQWC\
(\ref{eq:iqwc2})};
\end{eqnarray}
\begin{eqnarray}\label{eq:phi2i1}
\frac{v_0v_3+v_1v_2}{m_3(z_1\wedge z_2)}
+\frac{v_0v_2+v_1v_3}{m_2(z_2)}-\frac{v_0v_1+v_2v_3}{m_2(z_1)}
-4k\frac{v_0+v_1+v_2+v_3+k(z_1+z_2)}{a_1m_2(z_1)m_2(z_2)m_3(z_1\wedge z_2)}-\nonumber\\
\frac{k^2z_1z_2}{m_3(z_1\wedge z_2)}=0\ \mathrm{for\ IQWC\
(\ref{eq:iqwc1})}.
\end{eqnarray}
Note that the above homography for QC can be written as
$$X(v_0,v_3)^T(M_2(z_1)-M_2(z_2))(M_2(z_1)+M_2(z_2))^{-1}X(v_1,v_2)=1.$$
The relations (\ref{eq:v1v21}) and (\ref{eq:cross}) can be
interpreted as the same homography between $v_0,v_3$, thus
establishing a homography between $u_0,v_1,v_2,u_3$ analogue to
(\ref{eq:phi2c})-(\ref{eq:phi2i1}). For QC we just perform the
change $(u_0\leftrightarrow v_0,u_3\leftrightarrow v_3,
\sqrt{R_{z_j}}\leftrightarrow -\sqrt{R_{z_j}})$; for QWC
(\ref{eq:qwc1}) and IQWC (\ref{eq:iqwc1}) the change
$(u_0\leftrightarrow v_0,u_3\leftrightarrow
v_3,\sqrt{R'_{z_j}}\leftrightarrow r_{e_2}\sqrt{R'_{z_j}})$. For
QWC (\ref{eq:qwc2}) we obtain a new homography ($m_2^-,\ m_3$ have
same definition with $\sqrt{R'_{z_j}}$ replaced by
$\sqrt{R'_{z_j}}r_{e_2}$ and it is the reason why $m_3$ acquires
the $-$ sign; $a^-:=\bar f_1^TA'\bar f_1$):
\begin{eqnarray}\label{eq:phi2wv0v3}
-\frac{2v_1v_2}{m_3(z_1\wedge z_2)}
+\frac{u_0v_2+v_1u_3}{m_2^-(z_2)}-\frac{u_0v_1+v_2u_3}{m_2^-(z_1)}
-\frac{2a^-}{m_2^-(z_1)m_2^-(z_2)m_3(z_1\wedge z_2)}=0\nonumber\\ \mathrm{for\ QWC\
(\ref{eq:qwc2})};
\end{eqnarray}
for IQWC (\ref{eq:iqwc2}) we must perform $(u_j\leftrightarrow
v_j,\ j=1,2)$ in (\ref{eq:v1v21}) and consider the other ruling
families on $x_0^1,\ x_0^2$ in (\ref{eq:cross}); the new
(\ref{eq:v1v21}), (\ref{eq:cross}) are valid and can be
interpreted as the same homography between $u_0,u_3$, thus
establishing the homography ($m_2^-,\ m_3$ have same definition
with $\sqrt{R'_{z_j}}$ replaced by $r_{e_2}\sqrt{R'_{z_j}}$):
\begin{eqnarray}\label{eq:phi2i2v0v3}
-\frac{2u_1u_2}{m_3(z_1\wedge z_2)}
+\frac{v_0u_2+u_1v_3}{m_2^-(z_2)}-\frac{v_0u_1+u_2v_3}{m_2^-(z_1)}+2\frac{v_0+v_3}{m_3(z_1\wedge z_2)}
+2(\frac{1}{m_2(z_1)}+\frac{1}{m_2(z_2)})\frac{u_1+u_2}{m_3(z_1\wedge z_2)}-\nonumber\\
\frac{4}{m_3(z_1\wedge
z_2)^3}-\frac{8}{m_2(z_1)m_2(z_2)m_3(z_1\wedge z_2)}=0\
\mathrm{for\ IQWC\ (\ref{eq:iqwc2})}.
\end{eqnarray}
Because of the symmetry $(0\leftrightarrow 1,2\leftrightarrow 3)$
of the SITC, the homography between $v_0,u_1,u_2,v_3$ is valid
with a similar formula: for QC we just perform the change
$(u_1\leftrightarrow v_1,u_2\leftrightarrow v_2,
\sqrt{R_{z_j}}\leftrightarrow -\sqrt{R_{z_j}})$; for QWC
(\ref{eq:qwc1}) and IQWC (\ref{eq:iqwc1}) the change
$(u_1\leftrightarrow v_1,u_2\leftrightarrow
v_2,\sqrt{R'_{z_j}}\leftrightarrow r_{e_2}\sqrt{R'_{z_j}})$;
(\ref{eq:phi2wv0v3}) is valid for $(u_0,v_1,v_2,u_3)
\rightarrow(u_1,v_0,v_3,u_2)$ and (\ref{eq:phi2i2v0v3}) is valid
for $(v_0,u_1,u_2,v_3)\rightarrow(v_1,u_0,u_3,v_2)$.

Because of the first relation of (\ref{eq:simpas}),
(\ref{eq:cross}) remains valid if we consider  the other ruling
families on $x_0^1,\ x_0^2$ and since $(u_j\leftrightarrow v_j,\
j=0,...,3) \circ(\ref{eq:v1v21}),\ (u_j\leftrightarrow v_j,\
j=1,2)\circ(\ref{eq:v1v21})$ remain valid with the same sign
respectively for QC, (I)QWC, the new (\ref{eq:v1v21}),
(\ref{eq:cross}) can be interpreted as the same homography between
$v_0,\ v_3$; this establishes homographies between $u_j,\
j=0,...,3$ analogue to (\ref{eq:phi2c})-(\ref{eq:phi2i1}): for QC,
QWC (\ref{eq:qwc1}) and IQWC (\ref{eq:iqwc1}) we perform just the
change $u_j\leftrightarrow v_j,\ j=0,...,3$ (because of the
involutory character of the change in $\sqrt{R_{z_j}}\
(\sqrt{R'_{z_j}})$ and the change being done twice: once for
$(u_0\leftrightarrow v_0,u_3\leftrightarrow v_3)$ and once for
$(u_1\leftrightarrow v_1, u_2\leftrightarrow v_2)$); for the
remaining quadrics we get a new homography ($m_2,m_3$ remain the
same, since replacing $\sqrt{R'_{z_j}}$ by
$r_{e_2}\sqrt{R'_{z_j}}r_{e_2}$ does not change either;
$a^-:=-f_1^TA'f_1$ for IQWC (\ref{eq:iqwc2})); for IQWC
(\ref{eq:iqwc2}) it is obtained directly from the new
(\ref{eq:cross}):
\begin{eqnarray}\label{eq:phi2wu}
\frac{u_0u_3+u_1u_2}{m_3(z_1\wedge z_2)}
+\frac{u_0u_2+u_1u_3}{m_2(z_2)}-\frac{u_0u_1+u_2u_3}{m_2(z_1)}
+\frac{4a^-}{m_2(z_1)m_2(z_2)m_3(z_1\wedge z_2)}=0\nonumber\\
\mathrm{for\ QWC\ (\ref{eq:qwc2})\ \&\ \ IQWC\ (\ref{eq:iqwc2})}.
\end{eqnarray}
For the other choice of the sign in (\ref{eq:v1v21}) perform in
(\ref{eq:v1v21}) $u_j\leftrightarrow v_j,\ j=0,2,3\
(u_2\leftrightarrow v_2)$ respectively for QC ((I)QWC) and
consider in (\ref{eq:cross}) the other ruling family on $x_0^2$;
again the new (\ref{eq:v1v21}), (\ref{eq:cross}) can be
interpreted as the same homography between $u_0,\ v_3$; this
establishes a homography between $v_0,v_1,u_2,u_3$ analogue to
(\ref{eq:phi2c})-(\ref{eq:phi2i1}): for QC we just perform the
change $(u_2\leftrightarrow v_2,u_3\leftrightarrow v_3,
\sqrt{R_{z_2}}\leftrightarrow -\sqrt{R_{z_2}})$; for QWC
(\ref{eq:qwc1}) and IQWC (\ref{eq:iqwc1}) the change
$(u_2\leftrightarrow v_2,u_3\leftrightarrow v_3,
\sqrt{R'_{z_2}}\leftrightarrow r_{e_2}\sqrt{R'_{z_2}})$. For the
remaining quadrics we obtain a new homography ($m_3^-$ has same
definition as $m_3$ with $\sqrt{R'_{z_2}}$ replaced by
$\sqrt{R'_{z_2}}r_{e_2}$ for QWC (\ref{eq:qwc2}) and
$r_{e_2}\sqrt{R'_{z_2}}$ for IQWC (\ref{eq:iqwc2}), but it loses
its skew-symmetry):
\begin{eqnarray}\label{eq:phi2wu0v3}
\frac{v_0u_3+v_1u_2}{m_3^-(z_1\odot z_2)}
+\frac{v_0u_2+v_1u_3}{m_2^-(z_2)}-\frac{2v_0v_1}{m_2(z_1)}
+\frac{2a^-}{m_2(z_1)m_2^-(z_2)m_3^-(z_1\odot z_2)}=0\nonumber\\
\mathrm{for\ QWC\ (\ref{eq:qwc2})};
\end{eqnarray}
\begin{eqnarray}\label{eq:phi2iu0v3}
\frac{v_0u_3+v_1u_2}{m_3^-(z_1\odot z_2)}+\frac{v_0u_2+v_1u_3}{m_2^-(z_2)}
-\frac{2u_2u_3}{m_2(z_1)}+2\frac{v_0+v_1}{m_2(z_1)}
+2(\frac{2}{m_2(z_2)}+\frac{1}{m_2(z_1)})\frac{u_2+u_3}{m_2(z_1)}-\nonumber\\
\frac{2}{m_2(z_1)^3}-(\frac{2}{m_2(z_2)}+\frac{1}{m_2(z_1)})^2\frac{2}{m_2(z_1)}=0\
\mathrm{for\ IQWC\ (\ref{eq:iqwc2})}.
\end{eqnarray}
Choosing $u_1:=\infty$ in (\ref{eq:v1v21}) with the $+$ sign we
get the homography between $v_0,v_1,u_2,u_3$ for QC, QWC
(\ref{eq:qwc1}), IQWC (\ref{eq:iqwc1}) (this fact follows from the
similar statement for the homography between $v_0,v_1,v_2,v_3$)
and IQWC (\ref{eq:iqwc2}); for QWC (\ref{eq:qwc2}) we get the left
side of (\ref{eq:phi2wu0v3}) times a constant multiple of
$v_0(\bar f_1^T\sqrt{R'_{z_2}}f_1)^2-v_1\bar
f_1^T\sqrt{R'_{z_1}}f_1$.

Further we can consider in this case the corresponding
homographies between the other ruling families, but since the
homography between $v_0,v_1,u_2,u_3$ is invariant only under the
symmetry $(0\leftrightarrow 1,2\leftrightarrow 3)$ of the SITC,
all other symmetries of the SITC generate the sought homographies:
$(0\leftrightarrow 3, 1\leftrightarrow 2),\ (0\leftrightarrow
3,z_1\leftrightarrow z_2),\ (1\leftrightarrow 2,z_1\leftrightarrow
z_2)$ respectively generate the homographies between
$u_0,u_1,v_2,v_3;\ u_0,v_1,u_2,v_3;\ v_0,u_1,v_2,u_3$.

Thus the rulings of $x_0^0,x_0^3$ and $x_0^1,x_0^2$ involved in a
homography are simultaneously same or different; this dichotomy
gives the two choices of $x_0^2$ and can be interpreted as a
$\mathbb{Z}_2$ co-cycle: if we assign the signature $\ep:=1\ (-1)$
to $v\ (u)$ rulings, then $\ep_0\ep_1\ep_2\ep_3=1$.

If the homographies (\ref{eq:phi2c})-(\ref{eq:phi2i1}) are written
as $\phi=\phi(v_0,(v_1,z_1)\wedge(v_2,z_2),v_3)=0$, then the
algebraic relation:
\begin{eqnarray}\label{eq:calib}
\frac{z_1\mathcal{A}(\phi\phi_{v_2v_3}-\phi_{v_2}\phi_{v_3})}{\Del^-(z_1,v_0,v_1)}=
\frac{1}{m_2(z_1)m_2(z_2)m_3(z_1\wedge z_2)}
\end{eqnarray}
for independent variables $v_j,\ j=0,...,3$ plays an important
r\^{o}le in the analytic computations of the M\"{o}bius
configurations and in the proof of the BPT; it remains valid if we
perform changes $u_j\leftrightarrow v_j$ for some $j$'s according
to the previous recipes (in which case $\Del^-$ may change to
$\Del^+$ or $\Del'^\pm;\ m_2$ to $m_2^-;\ m_3$ to $-m_3$ or $\pm
m_3^-$). Bianchi (122,\S 25) assumes its use in the proof of the
BPT without explicitly stating it.

The homography H of Bianchi II preserves the Ivory affinity, the
RMPIA and the (SI)TC; therefore all results of the last three
subsections. The preservation of the SITC is suggested by the form
of the homographies (\ref{eq:phi2c})-(\ref{eq:phi2iu0v3}) and can
be confirmed analytically if one takes into consideration the
change of $z_k,\ k=1,2,\ A,\ u_j,\ v_j,\ j=0,...,3$ under $H$, but
an obvious geometric argument can replace these analytic
computations.

\subsection{M\"{o}bius configurations and an application of the Menelaus theorem}\label{subsec:algpre26}
\noindent

\noindent

A tetrahedron consists of $2^2$ points and $2^2$ planes, each
point (plane) belonging to (containing) $2+1$ planes (points).
M\"{o}bius considered configurations $\mathcal{M}_3$ of two
tetrahedra inscribed one into the other, that is configurations of
$2^3$ points and $2^3$ planes, each point (plane) belonging to
(containing) $3+1$ planes (points); two $\mathcal{M}_3$
configurations inscribed one into the other gives raise to a
configuration of $2^4$ points and $2^4$ planes, each point (plane)
belonging to (containing) $4+1$ planes (points). Therefore Bianchi
(\cite{B2},Vol {\bf 5},(117)) calls a configuration of $2^n$
points and $2^n$ planes such that any point (plane) belongs to
(contains) $n+1$ planes (points) a M\"{o}bius configuration
$\mathcal{M}_n$. We can consider the TC as $\mathcal{M}_1$, the
SITC as $\mathcal{M}_2$ and further iterates as $\mathcal{M}_n$'s.
Let $n\ge 2,\ k=1,...,n,\ j,j_1,j_2,...=0,...,2^n-1,\ x_0^0\in
x_0,\ z_k\in\mathbb{C},\ x_{z_k}^{2^{k-1}}\in x_{z_k}\cap
T_{x_0^0}x_0$; choose rulings $\hat v_0(=v_0$ or $u_0$) at $x_0^0$
and $\hat v_{2^{k-1}}(=u_{2^{k-1}}$ or $v_{2^{k-1}}$) at
$x_0^{2^{k-1}}$; as before we can order $0<z_1<...<z_n$. We would
like to get by repeated application of the SITC well defined
points $x_0^j\in x_0$ and rulings $\hat v_j$ at those points (and
thus RMPIA's $(R_{j_1}^{j_2},t_{j_1}^{j_2}),\
\log_2|j_1-j_2|\in\mathbb{N}$; the Ivory affinities are between
$x_0,\ x_{z_{\log_2|j_1-j_2|+1}})$. The rulings are well defined
(the signature $\ep_j=\ep_0^{s-1}\ep_{k_1}...\ep_{k_s}$ for
$0<k_1<...<k_s,\ j=\Sigma_{p=1}^s2^{k_p-1}$ is independent of the
path from $x_0^0$ to $x_0^j$). To define $x_0^j$ apply repeatedly
the SITC to a path $x_0^{j_0},x_0^{j_1},...x_0^{j_s},\
0=j_0<j_1<...<j_s=j,\ \log_2(j_p-j_{p-1})\in\mathbb{N},\
p=1,...,s$; this definition must be independent of the chosen path
(assume this for the moment); $\mathcal{M}_n$ is formed by facets
$(R_{j_1}^{j_0},t_{j_1}^{j_0})\circ(R_{j_2}^{j_1},t_{j_2}^{j_1})\circ...
\circ(R_{j_s}^{j_{s-1}},t_{j_s}^{j_{s-1}})(x_0^j,dx_0^j)$ (because
of the co-cycle relations (\ref{eq:cocy}) these are independent of
the sequence $j_0,j_1,j_2,...,j_s$ as above). If the rulings are
not prescribed, then for $n\ge 2$ we can choose $\hat v_0:=v_0$
(since considering other rulings on all $x_0^0,\
x_{z_k}^{2^{k-1}}$ does not change $\mathcal{M}_n$), so we have
$\frac{2^n}{2}=2^{n-1}$ choices of $\mathcal{M}_n$ (among the
$2^n$ choices of rulings on $x_{z_k}^{2^{k-1}}$ we have to remove
$\frac{1}{2}$ of them, since the dichotomy 'simultaneously same or
different' allows to consider the other rulings on all
$x_{z_k}^{2^{k-1}}$); for $n=1$ we have two choices of
$\mathcal{M}_1$.

We need now prove that the definition of $x_0^j$ does not depend
on the path from  $x_0^0$ to $x_0^j$; it is enough to consider
$n=3$, since by induction if all $x_0^j,\ j=0,...,2^{n+1}-2$ are
well defined, then we can uniquely define $x_0^{2^{n+1}-1}$ from
the quadrilateral $x_0^{2^{n-1}-1},
x_0^{2^n-1},x_0^{2^n+2^{n-1}-1},x_0^{2^{n+1}-1}$; if we separately
add each $x_0^{2^{n+1}-2^{k-1}-1},\ k=1,...,n-1$ to this
quadrilateral and complete it to an $\mathcal{M}_3$, then
$x_0^{2^{n+1}-1}$ is well defined for each $\mathcal{M}_3$
obtained, so $x_0^{2^{n+1}-1}$ is well defined for all paths from
$x_0^0$ to $x_0^{2^{n+1}-1}$. For $n=3:\ k=1,2,3,\
z_k\in\mathbb{C},\ x_0^0\in x_0,\ x_{z_k}^{2^{k-1}}\in x_{z_k}\cap
T_{x_0^0}x_0$ and consider the rulings $v_0,\ v_{2^{k-1}}$;
applying the SITC three times we get well defined facets
$(x_0^j,dx_0^j),\ j=0,...,6$. To simplify the notation replace
$x_0^0$ with $0,\ x_{z_1}^1$ with $1$, etc. If the rulings
$v_1,v_4,v_2$ respectively cut the lines $[3\ 5],\ [5\ 6],\ [6\
3]$ at the points $1_3^5,\ 4_5^6,\ 2_6^3$, then
$((1_3^5-3)\div(1_3^5-5))((4_5^6-5)\div(4_5^6-6))((2_6^3-6)\div(2_6^3-3))
=\frac{(x_{z_2u_1}^1)^T\hat N_0^3}{(x_{z_3u_1}^1)^T\hat N_0^5}
\frac{(x_{z_1u_4}^4)^T\hat N_0^5}{(x_{z_2u_4}^4)^T\hat N_0^6}
\frac{(x_{z_3u_2}^2)^T\hat N_0^6}{(x_{z_1u_2}^2)^T\hat N_0^3}=
\frac{(x_{z_2u_1}^1)^T\hat N_0^3}{(x_{z_1u_2}^2)^T\hat N_0^3}
\frac{(x_{z_1u_4}^4)^T\hat N_0^5}{(x_{z_3u_1}^1)^T\hat N_0^5}
\frac{(x_{z_3u_2}^2)^T\hat N_0^6}{(x_{z_2u_4}^4)^T\hat
N_0^6}=^{(\ref{eq:cross})}
\frac{z_2}{z_1}\frac{(x_{z_1u_1}^1)^T\hat
N_0^0}{(x_{z_2u_2}^2)^T\hat N_0^0}
\times\frac{z_1}{z_3}\frac{(x_{z_3u_4}^4)^T\hat
N_0^0}{(x_{z_1u_1}^1)^T\hat N_0^0}
\frac{z_3}{z_2}\frac{(x_{z_2u_2}^2)^T\hat
N_0^0}{(x_{z_3u_4}^4)^T\hat N_0^0}=1$. By the converse of the
Menelaus theorem  (which is essentially a co-cycle theorem) the
points $1_3^5,\ 4_5^6,\ 2_6^3$ in the plane $[3\ 5\ 6]$ are
co-linear. Consider the points $7_3^5,\ 7_5^6,\ 7_6^3$
respectively on the lines $[3\ 5],\ [5\ 6],\ [6\ 3]$ such that
they and the previous points cut these lines with the required
cross-ratios $\frac{z_j}{z_k}$; again by the converse of the
Menelaus theorem $7_3^5,\ 7_5^6,\ 7_6^3$ lie on a line $l_7$; a
similar statement holds for the other rulings $u_1,u_2,u_4$ and we
obtain another line $l'_7$: take $7:=l_7\cap l'_7$. By symmetry we
need now only prove that $7$ is in the facet of $5$, since in this
case the rulings $v_5,v_6$ will cut the line $[4\ 7]$ with same
cross-ratio as the one the rulings $l_7,v_4$ cut the line $[5\ 6]\
(v_4,v_5,v_6,l_7$ will be rulings of the same ruling family on a
quadric), so $7$ completes $4,5,6$ to an SITC quadrilateral. We
have $\frac{(x_{z_1u_0}^0)^T\hat N_0^1}{(x_{z_2u_0}^0)^T\hat
N_0^2}- \frac{(x_{z_1v_0}^0)^T\hat N_0^1}{(x_{z_2v_0}^0)^T\hat
N_0^2}=^{(\ref{eq:V21c})}\frac{(V_0^1)^TR^0_2\hat N_0^2(\hat
N_0^0)^T(x_{0v_0}^0\times x_{0u_0}^0)} {\mathcal{A}|\hat
N_0^0|^2(x_{z_2u_0}^0)^T\hat N_0^2(x_{z_2v_0}^0)^T\hat N_0^2}
=^{(\ref{eq:simpas})}-\frac{2(V_0^1)^TR^0_2\hat N_0^2}{z_2|\hat
N_0^2|^2}$; since $7-5$ is a multiple of
$\frac{z_3}{z_2}(\frac{(x_{z_2u_1}^1)^T\hat
N_0^3}{(x_{z_3u_1}^1)^T\hat N_0^5}-\frac{(x_{z_2v_1}^1)^T\hat
N_0^3}{(x_{z_3v_1}^1)^T\hat N_0^5})(6-4)-
\frac{z_1}{z_2}(\frac{(x_{z_2u_4}^4)^T\hat
N_0^6}{(x_{z_1u_4}^4)^T\hat N_0^5}-\frac{(x_{z_2v_4}^4)^T\hat
N_0^6}{(x_{z_1v_4}^4)^T\hat N_0^5})(3-1)$, the condition
$0=(7-5)^TR_1^0R_5^1\hat N_0^5$ becomes $0=(\hat
N_0^5)^T(R_1^5R_0^1R_4^0 V_4^6(R_1^5V_1^3)^T
-R_1^5R_0^1R_1^0V_1^3(R_4^5V_4^6)^T)\hat N_0^5=^{(\ref{eq:cocy})}
(\hat N_0^5)^T(R_4^5V_4^6(R_1^5V_1^3)^T
-R_1^5V_1^3(R_4^5V_4^6)^T)\hat N_0^5$.

By an obvious geometric argument and induction on $n,\
\mathcal{M}_n$ has the symmetries of an $n$-dimensional cube, the
static $n$-soliton $\hat v_{2^n-1}$ is a homography of $\hat v_0,\
\hat v_{2^{k-1}}:\ \hat v_{2^n-1}=\\\mathcal{S}_n(\hat v_0,(\hat
v_1,z_1),... ,(\hat v_{2^{n-1}},z_n))$; moreover symmetries of the
$n$-dimensional cube which do not induce symmetries of this
homography generate valid homographies between other rulings. Note
that $\hat v_0$ does not appear in $\hat v_{2^n-1}$ for $n$ odd:
take for example $n=3$. To fix the ideas consider $\hat v_j:=v_j,\
j=0,...,7$ satisfying the homographies imposed by the SITC; if we
replace $v_0,v_3,v_5,v_6$ with $u_0,u_3,u_5,u_6$ we have similar
homographies imposed by the SITC. Let now $v_0$ vary and keep
$u_0,v_1,v_2,v_4$ constant; thus $u_3,u_5$ are constant; because
of the homography between $v_1,u_3,u_5,v_7$ we conclude that $v_7$
is constant.

The 'addition' of static solitons now becomes clear: $\hat
v_{j_1}\oplus_{\hat v_0}\hat v_{j_2} =\hat v_{j_1+j_2}$, where
$1\le j_1,j_2\le 2^n-1$ have nonzero digits in base $2$ on
different positions. The SITC defines the law of 'addition' and
proves its commutativity on static $1$-solitons: $\hat
v_1\oplus_{\hat v_0}\hat v_2=\hat v_2\oplus_{\hat v_0}\hat
v_1=\hat v_3$; the symmetry of the TC provides an inverse. The
existence of $\mathcal{M}_3$ assures the associativity of the
'addition'. Once these are checked on static $1$-solitons $\hat
v_{2^{k-1}}$, they are true for static $n$-solitons (the existence
of $\mathcal{M}_2,\ \mathcal{M}_3$ assures the existence of
$\mathcal{M}_n$).

An analytic confirmation of the above properties is more
complicated, but all symmetries are revealed in the implicit form
(homography) $\phi_n(\hat v_0,\bigwedge_{k=1}^n(\hat
v_{2^{k-1}},z_k), \hat v_{2^n-1})=0$ of $\mathcal{S}_n:\\
\mathcal{S}_n(\hat v_0,\bigodot_{k=1}^n(\hat v_{2^{k-1}},z_k))
=-\frac{\phi_n(\hat v_0,\bigwedge_{k=1}^n(\hat
v_{2^{k-1}},z_k),0)} {\pa_v\phi_n(\hat v_0,\bigwedge_{k=1}^n(\hat
v_{2^{k-1}},z_k),v)}$. We use the notation
$f(\bigodot_{k=1}^nx_k)$ for $\\f(x_1,...,x_n)$ if
$f(x_{\si(1)},...,x_{\si(n)})=f(x_1,...,x_n)$ and
$f(\bigwedge_{k=1}^nx_k)$ for $f(x_1,...,x_n)$ if
$\\f(x_{\si(1)},...,x_{\si(n)})=(-1)^{\ep(\si)}f(x_1,...,x_n)$ for
$\si$ permutation of $1,...,n$; although $\phi_n$ has
skew-symmetries (and thus $\mathcal{S}_n$ has symmetries) only
when separately permuting $u$ rulings or $v$ rulings, we still
maintain for simplicity the $\bigwedge\ (\bigodot)$ notation for
all rulings $\hat v_{2^{k-1}},\ k=1,...,n$. Also $\phi_n$ (and
thus $\mathcal{S}_n$) depends on the signatures $\ep_0,\
\ep_{2^{k-1}},\ k=1,...,n$; only when strictly necessary we shall
use the notation $\phi_n^{\vaep},\
\vaep:=2^{n-1}(1-\ep_0)+\sum_{k=1}^n(1-\ep_{2^{k-1}})2^{k-2}$ to
point out the explicit dependence on rulings. Since for $n$ odd
$\hat v_{2^n-1}$ does not depend on $\hat v_0$ we drop in this
case $\hat v_0$ from $\phi_{n+1}(...)$.

Since IQWC (\ref{eq:iqwc2}) is metrically the most degenerate,
(\ref{eq:phi2wu}) provides the simplest formulae for $\phi_n$;
thus it is natural to study it first and conjecture relations
valid for all quadrics and rulings. We have
$\phi_2(u_0,\bigwedge_{k=1,2}(u_{2^{k-1}},z_k),u_3)=\begin{vmatrix}(u_1-u_0)& z_1(u_3-u_1)\\
(u_2-u_0)& z_2(u_3-u_2)\end{vmatrix}+
\frac{1}{4}\begin{vmatrix}z_1&z_1^2\\z_2&z_2^2\end{vmatrix},\\
\phi_{2m+1}(u_0,\bigwedge_{k=1}^{2m+1}(u_{2^{k-1}},z_k),u_{2^{2m+1}-1})=
\sum_{j=1}^{2m+1}(-1)^{j-1}\phi_{2m}(u_0,\bigwedge_{k=1,\ k\neq
j}^{2m+1}(u_{2^{k-1}},z_k),u_{2^{2m+1}-1}),\\
\phi_{2m+2}(u_0,\bigwedge_{k=1}^{2m+2}(u_{2^{k-1}},z_k),u_{2^{2m+2}-1})=\\
(\prod_{j=1}^{2m+2}z_j)\sum_{j=1}^{2m+2}\frac{(-1)^{j-1}(u_{2^{j-1}}-u_0)}{z_j}
\phi_{2m+1}(u_0,\bigwedge_{k=1,\ k\neq
j}^{2m+2}(u_{2^{k-1}},z_k),u_{{2^{2m+2}-1}}),\ m\ge 1$. Thus
$\phi_n$ is the linear combination of $2$ determinants
respectively with coefficients $1,\ \frac{1}{4}$; $\phi_{2m+1}$ is
obtained from $\phi_{2m}$ by replacing $u_{2^{2m}-1}$ with
$u_{2^{2m+1}-1}$, then padding its determinants with a first
column of $1$'s and with a last row of the same type as the
previous ($z_{2m+1},\ u_{2^{2m}}$ replace $z_j,\ u_{2^{j-1}}$);
$\phi_{2m+2}$ is obtained from $\phi_{2m+1}$ by multiplying the
$j^{\mathrm{th}}$ row of each of its determinants with $z_j$,
replacing $u_{2^{2m+1}-1}$ with $u_{2^{2m+2}-1}$, then padding its
determinants with a first column of
$(u_{2^{j-1}}-u_0)_{j=1,...,2m+1}$ and with a last row of the same
type as the previous ($z_{2m+2},\ u_{2^{2m+1}}$ replace $z_j,\
u_{2^{j-1}}$); for example
$\phi_3(u_0,\bigwedge_{k=1}^3(u_{2^{k-1}},z_k),u_7)
=\begin{vmatrix}1&u_1&z_1(u_7-u_1)\\1&u_2&z_2(u_7-u_2)\\
1&u_4&z_3(u_7-u_4)\end{vmatrix}+
\frac{1}{4}\begin{vmatrix}1&z_1&z_1^2\\1&z_2&z_2^2\\1&z_3&z_3^2\end{vmatrix},\\
\phi_4(u_0,\bigwedge_{k=1}^4(u_{2^{k-1}},z_k),u_{15})
=\begin{vmatrix}(u_1-u_0)&z_1&z_1u_1&z_1^2(u_{15}-u_1)\\
(u_2-u_0)&z_2&z_2u_2&z_2^2(u_{15}-u_2)\\
(u_4-u_0)&z_3&z_3u_4&z_3^2(u_{15}-u_4)\\
(u_8-u_0)&z_4&z_4u_8&z_4^2(u_{15}-u_8)\end{vmatrix}
+\frac{1}{4}\begin{vmatrix}(u_1-u_0)&z_1&z_1^2&z_1^3\\
(u_2-u_0)&z_2&z_2^2&z_2^3\\(u_4-u_0)&z_3&z_3^2&z_3^3\\
(u_8-u_0)&z_4&z_4^2&z_4^3\end{vmatrix}$. Note that the separate
linearity in $u_0,\ u_{2^{k-1}},\ u_{2^n-1}$ behavior mixes with
the Vandermonde determinant in $z_k$ behavior, so powers
$\lessapprox\frac{n}{2}$ in $z_k$ suffice to guarantee
non-singular determinants.

The determinant formulation of $\phi_n$ is proven by induction
after $n\ge 3$: consider the determinants
$D^j_n(u_0,\bigwedge_{k=1}^n(u_{2^{k-1}},z_k)),\ j=1,2,3$, where
$\phi_n(u_0,\bigwedge_{k=1}^n(u_{2^{k-1}},z_k),u_{2^n-1})=:\\u_{2^n-1}D_n^1(u_0,
\bigwedge_{k=1}^n(u_{2^{k-1}},z_k))+
D_n^2(u_0,\bigwedge_{k=1}^n(u_{2^{k-1}},z_k))
+\frac{1}{4}D_n^3(u_0,\bigwedge_{k=1}^n(u_{2^{k-1}},z_k))$. We
have
$\\D_{n+1}^1(u_0,\bigwedge_{k=1}^{n+1}(u_{2^{k-1}},z_k))((u_{2^n}-u_0)^2
-\frac{1}{4}z_{n+1}^2)^{[\frac{n}{2}]}(\frac{1+(-1)^n}{2}\prod_{k=1}^nz_k+\frac{1-(-1)^n}{2})
=\\D_n^1(u_{2^n},\bigwedge_{k=1}^n(u_{2^n+2^{k-1}},z_k))\prod_{k=1}^n
D_2^1(u_0,(u_{2^{k-1}},z_k)\wedge(u_{2^n},z_{n+1})),\ n\ge 2$ and
a similar relation obtained by replacing $D_{n+1}^1(...),\
D_n^1(...)$ with $D_{n+1}^2(...)+\frac{1}{4}D_{n+1}^3(...),\
D_n^2(...)+\frac{1}{4}D_n^3(...)$; since
$\phi_n(u_{2^n},\bigwedge_{k=1}^n(u_{2^n+2^{k-1}},z_k),u_{2^{n+1}-1})=0$
this will suffice. Just as the Vandermonde determinant $V_{n+1}$
in $z_k,\ k=1,...,n+1$ can be interpreted as a polynomial of
degree $n$ in $z_{n+1}$, with roots $z_k,\ k=1,...,n$ and $V_n$ as
the dominant coefficient, each of the two relations between the
functionally independent $u_0,\ u_{2^{k-1}},\ z_k,\ k=1,...,n+1$
can be interpreted as polynomials in $z_n$ respectively of degrees
$[\frac{n+1}{2}],\ [\frac{n+1}{2}]+1$, with roots among $z_k,\
k=1,...,n-1,n+1$ and whose dominant coefficients are $0$. The
condition that the dominant coefficients of the polynomials in
$z_n$ are $0$ can be interpreted in turn as two polynomials in
$z_{n-1}$ of degrees $[\frac{n+1}{2}]-1,\ [\frac{n+1}{2}]$, with
roots among $z_k,\ k=1,...,n-2,n+1$ and whose dominant
coefficients are $0$ (the dominant coefficients are linear in
$u_{2^{n-2}}$; the coefficients of $u_{2^{n-2}}$ being $0$ are the
relations that must be initially proved, but with $n-2$ replacing
$n$; the other coefficients being $0$ are also valid relations
proven by backward induction with step two). To prove the fact
that roots of the considered polynomials are among $z_k,\
k=1,...,n-1,n+1$ or $k=1,...,n-2,n+1$ we need $z_j=z_k\Rightarrow
u_{2^{j-1}}=u_{2^{k-1}}$; take for example $z_1=z_2$; then either
$u_1=u_2$ and $u_0,\ u_3$ are functionally independent, or
$u_0=u_3$ and $u_1,\ u_2$ are functionally independent; in the
second case the $n^{\mathrm{th}}$ soliton $u_{2^n-1}$ degenerates
to the $(n-2)^{\mathrm{th}}$ soliton $u_{2^n-4}$ and its formula
is thus valid; optionally one can apply a symmetry of the $n$
dimensional cube and exchange the functionally independent $u_1,\
u_2$ with $u_0,\ u_3$ and thus still get $u_1=u_2$.

For the other quadrics we cannot expect Vandermonde determinants
since $z_k$ do not appear polynomially, but some of the behavior
of the previous determinants remains; in particular since
$\frac{\Del'^+(z_{n+1},u_0,u_{2^n})}{4\mathcal{A}}=(u_{2^n}-u_0)^2
-\frac{1}{4}z_{n+1}^2$ the use of (\ref{eq:calib}) and backward
induction by step two naturally appears.

The algebraic relation
\begin{eqnarray}\label{eq:mob}
\phi_n(\hat v_{2^n},\bigwedge_{k=1}^n
(-\frac{\phi_2(\hat v_0,\bigwedge_{j=k,n+1}(\hat v_{2^{j-1}},z_j),0)}
{\pa_v\phi_2(\hat v_0,\bigwedge_{j=k,n+1}(\hat v_{2^{j-1}},z_j),v)},z_k),\hat v_{2^{n+1}-1})
\prod_{k=1}^n\pa_v\phi_2(\hat v_0,\bigwedge_{j=k,n+1}(\hat v_{2^{j-1}},z_j),v)=\nonumber\\
t_{n+1}(\frac{\hat \Del(z_{n+1},\hat v_0,\hat
v_{2^n})}{z_{n+1}\mathcal{A}})^{[\frac{n}{2}]} \phi_{n+1}(\hat
v_0,\bigwedge_{k=1}^{n+1}(\hat v_{2^{k-1}},z_k),\hat
v_{2^{n+1}-1}),\ n\ge 2
\end{eqnarray}
is valid for independent variables $\hat v_0,\ \{\hat
v_{2^{k-1}}\}_{k=1,...,n+1},\ \hat v_{2^{n+1}-1}$ and certain
constants $\\t_n=t_n(z_1,...,z_n)$ ($\hat\Del$ is one of
$\Del^{\pm},\ \Del'^{\pm}$ according to the choice of the rulings
$\hat v_0,\ \hat v_{2^n}$ and thus the use of (\ref{eq:calib})
becomes clear). Since $\hat v_{2^n+2^{k-1}}= -\frac{\phi_2(\hat
v_0,\bigwedge_{j=k,n+1}(\hat v_{2^{j-1}},z_j),0)}
{\pa_v\phi_2(\hat v_0,\bigwedge_{j=k,n+1}(\hat
v_{2^{j-1}},z_j),v)}$ (more generally for
$j:=\Sigma_{s=1}^p2^{k_s-1},\ 0<k_1<...<k_p\le n$ we have
$\phi_p(\hat v_0,\bigwedge_{s=1}^p(\hat
v_{2^{k_s-1}},z_{k_s}),\hat v_j)=0$), (\ref{eq:mob}) gives a
procedure of generating $\phi_n$ by induction after $n$.

From a practical point of view (generate $\phi_n$) we don't need
to completely prove (\ref{eq:mob}): note that the lhs of
(\ref{eq:mob}) and $\phi_{n+1}(...)$ must be simultaneously
annihilated; thus (\ref{eq:mob}) must be true with the coefficient
of $\phi_{n+1}(...)$ being an a-priori arbitrary polynomial in
$\hat v_0,\ \hat v_{2^n}$ (the remaining variables $\{\hat
v_{2^{k-1}}\}_{k=1,...,n+1},\ \hat v_{2^{n+1}-1}$ appear linearly,
since we know that $\phi_{n+1}$ must be a homography). Using few
terms of the highest degree of this polynomial (with the
simplification $\hat v_0=0$ for $n$ even, since $\hat
v_{2^{n+1}-1}$ does not depend on $\hat v_0$) and the required
skew-symmetryies of $\phi_{n+1}(...)$ wrt $\{(\hat
v_{2^{k-1}},z_k)\}_{k=1,...,n+1}$, one can find $\phi_{n+1}(...),\
t_{n+1}$ and the considered terms (which turn out to be as  they
should).

By this practical procedure we first find $\phi_n$ for $n=3,4$ and
assume it to be found for $n\ge 5$; then (\ref{eq:mob}) is proven
by backward induction after $n$ with step $2$ (the validity of
(\ref{eq:mob}) for $n$ is reduced to the validity of
(\ref{eq:mob}) for $n-2$).

For $\phi_3$ again there is a formulation for all QC:
$$X(\hat v_1,\hat v_7)^T(M_3(z_1\wedge z_2)-M_3(z_1\wedge z_3))(M_3(z_1\wedge z_2)
+M_3(z_1\wedge z_3))^{-1}X(\hat v_2,\hat v_4)=1;$$ one replaces
$\sqrt{R_{z_k}}$ with $-\sqrt{R_{z_k}}$ in $M_3$ if $\hat v_0,\
\hat v_{2^{k-1}}$ have different signatures.

We have:
\begin{eqnarray}\label{eq:phi3}
\phi_3(\bigwedge_{k=1}^3(\hat v_{2^{k-1}},z_k),\hat v_7)=\sum_{k=1}^3(-1)^{k-1}
\phi_2(\hat v_7,\bigwedge_{j=1,\ j\neq k}^3(\hat v_{2^{j-1}},z_j),\hat v_{2^{k-1}}).
\end{eqnarray}
Note that $\phi_3$ for rulings of same signature has the most
symmetries (the symmetries of a regular tetrahedron), since all
rulings are on equal footing (for example for $\phi_2$ the
relationship between $v_0,v_1$ is different from that between
$v_0,v_3$).

To prove (\ref{eq:mob}) for $n=2,\
t_3=-\prod_{k=1}^3\frac{1}{m_2(z_k)},\ \hat v_0=0$ we need
$$I(\bigwedge_{k=1}^3z_k):\
\sum_{k=1}^3(-1)^{k-1}\frac{m_2(z_k)}{m_3(\bigwedge_{j=1,j\neq k}^3z_j)}
=\prod_{k=1}^3\frac{m_2(z_k)}{m_3(\bigwedge_{j=1,j\neq k}^3z_j)},$$
$$II(\bigwedge_{k=1}^3z_k):\
\sum_{k=1}^3(-1)^{k-1}\frac{\prod_{j=1,j\neq
k}^3m_2(z_j)}{m_3(\bigwedge_{j=1,j\neq k}^3z_j)}
=4a^2\del\prod_{k=1}^3\frac{1}{m_3(\bigwedge_{j=1,j\neq
k}^3z_j)};$$ again we maintain the $\wedge$ notation for $m_3$,
although sometimes $m_3$ loses the skew-symmetry and we skip the
$-$ when it appears in $m_2^-,\ m_3^-$. For all (I)Q(W)C except
the general QC (\ref{eq:qc1}) ($A^{-1}$ having distinct
eigenvalues) and all rulings  one can provide a justification of
(\ref{eq:phi3}) without using (\ref{eq:mob}) (and thus I\& II are
consequences and do not need proof): from
(\ref{eq:phi2c})-(\ref{eq:phi2iu0v3}),(\ref{eq:phi3}) we have
$\sum_{k=1}^3(-1)^{k-1}(\phi_2(\hat v_0,\bigwedge_{j=1,\ j\neq
k}^3(\hat v_{2^{j-1}},z_j), \hat v_{7-2^{k-1}})+\phi_2(\hat
v_7,\bigwedge_{j=1,\ j\neq k}^3(\hat v_{7-2^{j-1}},z_j), \hat
v_{2^{k-1}}))=\phi_3(\bigwedge_{k=1}^3(\hat v_{2^{k-1}},z_k),\hat
v_7) +\phi_3(\bigwedge_{k=1}^3(\hat v_{7-2^{k-1}},z_k),\hat v_0)$.
By the same argument as when was proven that $\hat v_7$ does not
depend on $\hat v_0$, one can keep $\hat v_1,\hat v_2,\hat
v_4,\hat v_7$ constant and vary $\hat v_0,\hat v_3,\hat v_5,\hat
v_6$ and conversely; thus $\phi_3(\bigwedge_{k=1}^3(\hat
v_{2^{k-1}},z_k),\hat v_7)=c(z_1,z_2,z_3)$; by symmetry
$\phi_3(\bigwedge_{k=1}^3(\hat v_{7-2^{k-1}},z_k),\hat
v_0)=c(z_1,z_2,z_3)$, so $c(z_1,z_2,z_3)=0$.

The reader familiar with the proof of the BPT for the B
transformation of the sine-Gordon equation (see Bianchi
(\cite{B2},Vol {\bf 5},(46))) may recall similar cancellations and
rightfully so: we shall see later that these algebraic
computations encode all necessary information for the proof of the
BPT.

Because $\phi_2$ for QWC (\ref{eq:qwc2}), IQWC (\ref{eq:iqwc2})
are missing some quadratic terms for mixed rulings,
(\ref{eq:phi3}) and the above argument have to be slightly
modified. For example:
\begin{eqnarray}\label{eq:phi34}
\phi_3^4((v_1,z_1)\wedge(v_2,z_2),(u_4,z_3),u_7)=\phi_2^2(u_7,(v_2,z_2),(u_4,z_3),v_1)
-\phi_2^2(u_7,(v_1,z_1),(u_4,z_3),v_2)+\nonumber\\
2\phi_2^0(u_7,(v_1,z_1)\wedge(v_2,z_2),0)\
\mathrm{for\ QWC\ (\ref{eq:qwc2})},\nonumber\\
\phi_3^4((v_1,z_1)\wedge(v_2,z_2),(u_4,z_3),u_7)=\phi_2^2(v_1,(v_2,z_2),(u_4,z_3),u_7)
-\phi_2^2(v_2,(v_1,z_1),(u_4,z_3),u_7)+\nonumber\\
\phi_2^7(0,(0,z_1)\wedge(0,z_2),0)-
\phi_2^7(u_7,(v_1,z_3)\wedge(v_2,z_3),u_4)\ \mathrm{for\ IQWC\
(\ref{eq:iqwc2})};
\end{eqnarray}
further
$\phi_2^2(v_0,(v_2,z_2),(u_4,z_3),u_6)-\phi_2^2(v_0,(v_1,z_1),(u_4,z_3),u_5)
+\phi_2^2(v_3,(v_1,z_2),(u_7,z_3),u_5)-\\\phi_2^2(v_3,(v_2,z_1),(u_7,z_3),u_6)
+2\phi_2^0(v_0,(v_1,z_1)\wedge(v_2,z_2),v_3)=\phi_3^4((v_1,z_1)\wedge(v_2,z_2),(u_4,z_3),u_7)
+\phi_3^4((v_0,z_1)\wedge(v_3,z_2),(u_5,z_3),u_6)$ for QWC
(\ref{eq:qwc2}); for IQWC (\ref{eq:iqwc2}) replace
$\phi_2^0(v_0,(v_1,z_1)\wedge(v_2,z_2),v_3)$ with
$\phi_2^7(u_7,(u_6,z_1)\wedge(u_5,z_2),u_4)$. We don't need
(\ref{eq:phi34}) to find $\phi_3^4$: because of the symmetries of
the $3$-dimensional cube we have
$\phi_3^4((v_1,z_1)\wedge(v_2,z_2),(u_4,z_3),u_7)=\phi_3^6((v_1,z_3),(u_7,z_2)\wedge(u_4,z_1),v_2);\
$ for $\phi_3^6$ (\ref{eq:phi3}) and the above argument remain
valid without modification. Note also
$\\\sum_{k=1}^3(-1)^{k-1}(\phi_2^3(v_0,\bigwedge_{j=1,\ j\neq
k}^3(u_{2^{j-1}},z_j),
v_{7-2^{k-1}})+\phi_2^4(u_7,\bigwedge_{j=1,\ j\neq
k}^3(v_{7-2^{j-1}},z_j),
u_{2^{k-1}}))=\\-2\phi_3^0(\bigwedge_{k=1}^3(v_{7-2^{k-1}},z_k),v_0)$
for QWC (\ref{eq:qwc2}),
$=-2\phi_3^{15}(\bigwedge_{k=1}^3(u_{2^{k-1}},z_k),u_7)$ for IQWC
(\ref{eq:iqwc2}). The $-$ sign appears in the rhs (for all
quadrics) because $m_3(z_{k_1}\wedge z_{k_2})$ acquires a $-$ sign
for $\phi_2^3,\ \phi_2^4$.

We have
$t_n=-\frac{1}{m_2(z_n)^{n-2}}\prod_{k=1}^{n-1}\frac{1}{m_2(z_k)}$;
to get $\phi_4$ from (\ref{eq:mob}) for $n=3$ as previously
described we need the consequences $\mathrm{rhs\
II}(\bigwedge_{k=1}^3z_k)\prod_{k=1}^3\frac{1}{m_3(\bigwedge_{j=1,j\neq
k}^3z_j)}= \sum_{k=1}^3\frac{(-1)^{k-1}\mathrm{lhs\
I}(\bigwedge_{j=1,j\neq k}^4z_j)}
{m_2(z_4)m_3(\bigwedge_{j=k,4}z_j)}(=\sum_{k=1}^3\frac{(-1)^{k-1}}
{m_3(\bigwedge_{j=1,\ j\neq k}^3z_j)m_3(\bigwedge_{j=k,4}z_j)}),\
\sum_{k=1}^4\frac{(-1)^{k-1}\mathrm{rhs\ I}(\bigwedge_{j=1,j\neq
k}^4z_j)} {\prod_{j=1,\ j\neq
k}^4m_2(z_j)}=\sum_{k=1}^4\frac{(-1)^{k-1}\mathrm{lhs\ I}
(\bigwedge_{j=1,j\neq k}^4z_j)}{\prod_{j=1,\ j\neq
k}^4m_2(z_j)}\\(=0)$ of I,II; in fact one can get a whole
hierarchy of consequences of I,II for  $n\ge 5$.

\subsection{Doubly ruled M\"{o}bius configurations}\label{subsec:algpre27}
\noindent

\noindent We consider {\it doubly ruled M\"{o}bius configurations}
(DRMC) as a natural generalization of M\"{o}bius configurations,
apparent at the level of analytic computations. Let $k=1,2,3$;
consider rulings $\hat v_0(=u_0$ or $v_0),\hat
v_{2^{k-1}}(=u_{2^{k-1}}$ or $v_{2^{k-1}})\in \mathbb{C}$,
signatures $\ep_0(=1$ if $\hat v_0=v_0$ and $-1$ otherwise),
$\ep_{2^{k-1}}(=1$ if $\hat v_{2^{k-1}}=v_{2^{k-1}}$ and $-1$
otherwise) associated to these rulings and $z_k\in \mathbb{C}$.
The homographies below depend a-priori in an arbitrary fashion on
$z_k$.

Assume:

{\it I The $1^{\mathrm{st}}$ DRMC $\mathcal{M}_1$ is a homography
$\phi_1((u_0,v_0)\odot(u_1,v_1),z_1)=0$}.

{\it II Given $\phi_1((u_0,v_0)\odot(u_1,v_1),z_1)=0,\
\phi_1((u_0,v_0)\odot(u_2,v_2),z_2)=0$, there are two choices of
$(u_3,v_3)$ (according to $\ep_1\ep_2=\pm 1$) with
$\phi_1((u_3,v_3)\odot(u_1,v_1),z_2)=0,\
\phi_1((u_3,v_3)\odot(u_2,v_2),z_1)=0$}.

The $2^{\mathrm{nd}}$ DRMC $\mathcal{M}_2$ is a homography
$\phi_2(\hat v_0,(\hat v_1,z_1),(\hat v_2,z_2),\hat v_3)=0,\
\ep_3=\ep_0\ep_1\ep_2$. Moreover symmetries of the square $0123$
below which do not induce symmetries of $\phi_2$ generate valid
versions of $\phi_2$ (for simplicity we drop in the notation
$\phi_2$ the dependence on $\ep_0,\ep_1,\ep_2$).

{\it III Given $\phi_2(\hat v_0,(\hat v_2,z_2),(\hat v_4,z_3),\hat v_6)=0,\
\phi_2(\hat v_0,(\hat v_1,z_1),(\hat v_4,z_3),\hat v_5)=0,\\
\phi_2(\hat v_0,(\hat v_1,z_1),(\hat v_2,z_2),\hat v_3)=0,\
\ep_{2^{k_1-1}+2^{k_2-1}}=\ep_0\ep_{2^{k_1-1}}\ep_{2^{k_2-1}}$,
there is a unique $\hat v_7,\ \ep_7=\ep_1\ep_2\ep_4$ with
$\phi_2(\hat v_1,(\hat v_3,z_2),(\hat v_5,z_3),\hat v_7)=0,\
\phi_2(\hat v_2,(\hat v_3,z_1),(\hat v_6,z_3),\hat v_7)=0,\
\phi_2(\hat v_4,(\hat v_5,z_1),(\hat v_6,z_2),\hat v_7)=0$.}

The $3^{\mathrm{rd}}$ DRMC $\mathcal{M}_3$ is a homography
$\phi_3((\hat v_1,z_1),(\hat v_2,z_2),(\hat v_4,z_3),\hat v_7)=0,\
\ep_7=\ep_1\ep_2\ep_4$. Moreover symmetries of the cube below
which do not induce symmetries of $\phi_3$ generate valid versions
of $\phi_3$.
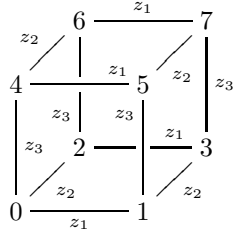
\begin{center}
$\xymatrix@!0{& 6 \ar@{-}[rr]^{z_1}\ar@{-}'[d][dd]_{z_3} & & 7 \ar@{-}[dd]^{z_3} \\
4 \ar@{-}[ur]^{z_2}\ar@{-}[rr]^<<<<<<<<<<<<{z_1}\ar@{-}[dd]^{z_3} & & 5 \ar@{-}[ur]_>>>>>>{z_2}
\ar@{-}[dd]_>>>>>>>>>>>{z_3} \\
& 2 \ar@{-}'[r][rr]^{z_1} & & 3 \\
0 \ar@{-}[rr]_{z_1}\ar@{-}[ur]_{z_2} & & 1 \ar@{-}[ur]_{z_2}}$
\end{center}
As in the previous subsection these three assumptions allows us to
generate the $n^{\mathrm{th}}$ DRMC $\mathcal{M}_n$ as a
homography $\phi_n(\hat v_0,(\hat v_1,z_1),..., (\hat
v_{2^{n-1}},z_n),\hat v_{2^n-1})=0$. Moreover symmetries of the
$n$-dimensional cube which do not induce symmetries of $\phi_n$
generate valid versions of $\phi_n$.

Thus the natural question arises if all DRMC can be reduced, via
M\"{o}bius transformations $\\u\mapsto\frac{au+b}{cu+d}$ (applied
to all rulings), to one of the cases of the previous subsection.

In particular, is the relation $\phi_3((\hat v_1,z_1),(\hat
v_2,z_2),(\hat v_4,z_3),\hat v_7)=\phi_2(\hat v_7,(\hat
v_2,z_2),(\hat v_4,z_3),\hat v_1)-\\\phi_2(\hat v_7,(\hat
v_1,z_1),(\hat v_4,z_3),\hat v_2)+\phi_2(\hat v_7,(\hat
v_1,z_1),(\hat v_2,z_2),\hat v_4)$ still valid (after proper
normalization of $\phi_2$ and with some modifications in
exceptional cases)?

\section{The rolling problem}\label{sec:rolling}\setcounter{equation}{0}

\subsection{Submanifolds in $\mathbb{C}^n$}\label{subsec:rolling1}
\noindent

\noindent We shall recall without many details the concept of
sub-manifolds of $\mathbb{C}^n$ (called henceforth sub-manifolds).
In the early development of the geometry of sub-manifolds in
scalar product spaces geometers made liberal use of imaginary
numbers, until this geometry has been divided into Lorentz and
Riemannian.

We would like to have the  GW and  the GCMR equations for a
sub-manifold $x$ in $\mathbb{C}^n$ just as the GW and GCMR
equations for a real sub-manifold in $\mathbb{R}^n$. The important
thing in creating these formulae is the existence of the normal
bundle; thus the natural definition for $k$-dimensional
sub-manifolds is that they have $k$-dimensional complex tangent
spaces (called henceforth tangent spaces), $1\le k\le n-1$ and
thus $(n-k)$-dimensional normal spaces. We exclude the case of
isotropic tangent spaces (spaces having tangency with the
isotropic cone), when the linear element degenerates and one
obtains among others  sub-Riemannian manifolds (the linear element
is real and positive semi-definite). The real dimension of the
real tangent spaces is a-priori any number between $k$ and $2k$
(we exclude the case when the dimension of the real tangent spaces
varies while the dimension of the complex tangent spaces remains
constant). The distribution $T\cap iT$ formed by intersecting real
tangent spaces $T$ with $iT$ is Frobenius integrable (by its own
definition); on its leaves the Nijenhuis condition holds, the
complex structure being induced by the surrounding space; from
Nirenberg-Newlander the leaves are holomorphic sub-manifolds. The
Frobenius theorem gives a parametrization holomorphic on these
leaves; in this parametrization it appears clearly that the
manifold is essentially $k$-dimensional, although apparently
initially could have had real dimension between $k$ and $2k$. If
it so happens that after possibly applying a change of variables
the parametrization in some of the remaining real parameters
reveals itself to be real analytic, it can be extended by
analyticity to a $k$-dimensional sub-manifold with more complex
parameters (and less real).

Take for example the IQWC (\ref{eq:iqwc2}) $x_0=uvf_1+u\bar
f_1+ve_3$ and let $x_0^0\subseteq x_0$ be a surface of real
dimension $3$ obtained by restricting $(u,v)$ to the domain
$(u,v)=(x+iy,s+ixy),\ x,y,s\in\mathbb{R}$; this is chosen so that
$K(x_0^0)=\frac{1}{(u^2-2v)^2}=\frac{1}{(x^2-y^2-2s)^2}>0$. We
know that there is a parametrization $(z,t),\ z\in\mathbb{C},\
t\in\mathbb{R}$ of $x_0^0$; since in this case $K(x_0^0)$ is
holomorphic in $z$ while real valued, we get $K(x_0^0)=f(t)$; we
can choose $t:=u^2-2v=x^2-y^2-2s$ and
$x_0^0=u\frac{u^2-t}{2}f_1+u\bar f_1+\frac{u^2-t}{2}e_3$, so we
can choose $z:=u$. In this case $T\cap iT=\mathbb{C}x_{0z}^0=
\mathbb{C}(x_{0u}^0+ux_{0v}^0)$ is real $2$-dimensional.

We need the normal bundle (and the fundamental forms of the
sub-manifold) to be holomorphic in the complex parameters. This is
done as follows: along the holomorphic leaves the tangent spaces
to the sub-manifold are spanned by the derivatives of the position
vector in the complex parameters and in the real parameters, but
because along such leaves the real parameters are set to
constants, the vectors spanning the tangent spaces are holomorphic
in the complex parameters, so we can find in the normal spaces
orthonormal frames which vary holomorphically in the complex
parameters. Thus the linear element and the second fundamental
form are holomorphic in the complex parameters and the
sub-manifold is recovered (the Gau\ss-Bonnet theorem) modulo rigid
motions of $\mathbb{C}^n$ from the knowledge of its linear element
and the second fundamental form through the integration of a
Ricatti equation (the GW equations) and a quadrature.

For two applicable sub-manifolds of $\mathbb{C}^n$ one can choose
orthonormal frames for both normal bundles and define the rotation
of the rolling as the family of orthogonal transformations which
takes the bases of the tangent spaces formed by parametric vectors
and the orthonormal frames for the normal bundles at corresponding
points one onto the other. The rotation can be arbitrarily chosen
so as to take any orthonormal frame of the first normal bundle
into any other of the second normal bundle, so we take advantage
and define it to take a holomorphic (in the complex parameters)
orthonormal frame into a holomorphic one. In this way the rotation
is holomorphic in the complex parameters (and unique modulo right
transformations through sections of $\mathbf{O}_{n-k}(\mathbb{C})$
bundles holomorphic in the complex parameters) and so will turn
out to be the translation.

For example if $x(u,v)$ is a holomorphic surface and we impose a
nonanalytic condition $f(u,\bar u)=0$ or explicitly $u(t)=g(t),\
t\in\mathbb{R}$, then the resulting surface $y(t,v):=x(u(t),v)$ is
nonanalytic in $t$, but can be extended to the holomorphic surface
$x(u,v)$.

Ruled surfaces with real ruling or real surfaces of revolution
extend by analyticity to ruled surfaces or surfaces of revolution.

However if we take for example a holomorphic ruled surface
$x(u,v):=a(u)v+b(u)$, non-analytically restrict the complex
variable $u$ to a curve $u(t)$ and then one deforms the resulting
$y(t,v):=x(u(t),v)$ (and analytically inconsistent with $u(t)$'s
non-analyticity) to a ruled surface $z(t,v)$ with the same ruling,
then $z(t,v)$ cannot be extended to a holomorphic ruled surface
applicable to the initial one.

Intersecting such a $k$-dimensional sub-manifold of $\mathbb{C}^n$
with a totally real subspace (here we use the definition of a
totally real space as being obtained after a rigid motion of
$\mathbb{C}^n$ from an $n$-dimensional real subspace of
$\mathbb{C}^n$ having real or purely imaginary entries in the
standard basis; all definitions of totally real subspaces
coincide) we get a real sub-manifold of a Lorentz space (the
imaginary part of the scalar product in $\mathbb{C}^n$
disappears), possibly analytic in some parameters.

A real analytic $k$-dimensional sub-manifold of $\mathbb{R}^n$
extends by analyticity to a holomorphic $k$-dimensional
sub-manifold of $\mathbb{C}^n$. One can iterate this procedure: a
holomorphic $k$-dimensional sub-manifold in $\mathbb{C}^n$ can be
considered as a real analytic $2k$-dimensional sub-manifold of
$\mathbb{R}^{2n}$ which can be complexified to a holomorphic
$2k$-dimensional sub-manifold in $\mathbb{C}^{2n}$, but the
Euclidean scalar product on $\mathbb{C}^n$ becomes a pair of
indefinite scalar products on $\mathbb{C}^{2n}$, so it is not
interesting from a metric point of view.

CMC(minimal) surfaces have analytic data in conformal
parametrization, so they extend to holomorphic CMC(minimal)
surfaces. The Lawson or Bianchi-Bryant-Weierstrass isometric
correspondence between minimal surfaces and CMC surfaces in space
forms always has as domain a space form of smaller curvature; for
the domain to be a space form of higher curvature we need
time-like mean curvature vector, so we need deformations in
Lorentz spaces.

According to Bianchi (\cite{B2},Vol {\bf 4},(202)), Lie noticed
that the formula developed by Monge for minimal surfaces can be
stated as: all minimal surfaces are surfaces of translation with
isotropic generating curves, that is $x=f(z)+\bar f(\bar z),\
f_z^Tf_z=0$ and this is one step from the Weierstrass formula,
only that at that time Cauchy had not yet rigorously developed the
theory of holomorphic functions.

Bianchi found in (\cite{B2},Vol {\bf 4},(108)) a conformal
correspondence between CMC surfaces in $\mathbb{R}^3$ and CMC
surfaces in $\mathbb{S}^3,\ (\mathbb{H}^3)$: if one replaces the
Euclidean metric on a ball containing two parallel CMC surfaces
(the existence of this configuration was proven by Bonnet) with
the metric induced by considering the boundary of the given ball
as the Cayley's absolute, we obtain two parallel surfaces in
$\mathbb{S}^3,\ (\mathbb{H}^3)$ (that is they are normal to a
normal congruence of geodesics) in conformal correspondence (the
conformal correspondence from the Euclidean space is preserved),
so they must have CMC.

While this conformal correspondence is not of applicability, it is
good enough to study the topology of these surfaces; for example
CMC tori (Wente, Pinkall- Sterling \cite{PS}, Bobenko's CMC tori
in terms of $\theta$-functions) or CMC compact surfaces of higher
genus (Kapouleas), including the case when these surfaces are
tangent to the boundary of the ball in discussion or when they
have singularities (for example cusps or degenerate linear
element) at this boundary.

\subsection{Rolling surfaces}\label{subsec:rolling2}
\noindent

\noindent The study of the rolling problem was initiated by
Ribacour and has been extensively pursued in Bianchi
\cite{B1},(\cite{B2},Vol {\bf 7}) and Darboux \cite{D1}.

Let $(u^1,u^2)\ni D$ with $D$ a domain in $\mathbb{R}^2,\
\mathbb{C}\times\mathbb{R}$ or $\mathbb{C}^2$ and $x:D\longmapsto
\mathbb{C}^3$ be a surface. Then $|dx|^2=dx^Tdx=g_{jk}du^j\otimes
du^k$ is the first fundamental form. Since $dx^T\wedge dx=0$,
$|dx|^2$ is a symmetric $(0,2)$ tensor. We have $d^2x^Tdx=\Gamma
_{jk,l}du^j\otimes du^k\otimes du^l$, where $\Gamma _{jk,l}$ are
the Christoffel symbols of first type. Since $d^2x$ is symmetric
we have $\Gamma _{jk,l}=\Gamma _{kj,l}$. Denote $\pa
x^2:=g^{jk}\pa _{u^j}\otimes\pa_{u^k}$, so $\pa
x^2\lrcorner|dx|^2=|dx|^2\llcorner\pa x^2=\pa_{u^j}\otimes du^j$
is the identity element of both $\lrcorner,\ \llcorner$ (with
$a,b$ tensors by $a\lrcorner b=b\llcorner a$ we contract the first
contravariant entry of $a$ with the first covariant entry of $b$;
one can stack the contractions in wait for the contravariant or
covariant entries without problems when working with symmetric
tensors).

For a function $f:D\longmapsto \mathbb{C}$ we have $\nabla
f:=df\llcorner\pa x^2=f_{u^j}g^{jk}\pa _{u^k}$ the tensor dual to
$df$, so we have $\Gamma ^l_{jk}\pa _{u^l}\otimes du^j\otimes
du^k=\na x^Td^2x$. Denote $\nabla ^2f:=d^2f-df\llcorner\na
x^Td^2x$ the second covariant derivative of $f$. We have
$0=d(|dx|^2\llcorner\pa x^2)$ so $d(\pa x^2)=-\pa x^2\lrcorner
d(|dx|^2)\llcorner\pa x^2=-\pa x^2\lrcorner d^2x^T\na x-\na
x^Td^2x\llcorner\pa x^2$.

If $N$ is a unit normal frame of the surface, then the second
fundamental form is $-dx^TdN=d^2x^TN$ and the GW equations are:
$d^2x=(d^2x)^{\top}+(d^2x)^{\bot}=dx\llcorner\na
x^Td^2x+NN^Td^2x$, or equivalently: $\na ^2x=NN^Td^2x$,
$dN=(dN)^{\top}=dx\llcorner\na x^TdN=-\na x\lrcorner d^2x^TN$.

For $\om_1,\ \om_2$ vector 1 forms in $\mathbb{C}^3$ and
$a,b\in\mathbb{C}^3$, we have $a^T\om_1\wedge b^T\om_2=((a\times
b)\times\om_1+b^T\om_1a)^T\wedge\om_2=(a\times
b)^T\om_1\times\wedge\om_2+b^T\om_1\wedge a^T\om_2$. In particular
$a^T\om\wedge b^T\om=\frac{1}{2}(a\times b)^T\om\times\wedge\om$.
Since both $\times$ and $\wedge$ are skew-symmetric, we have
$\om_1\times\wedge\om_2=\om_1\times\om_2+\om_2\times\om_1=\om_2\times\wedge\om_1$.

Consider the scalar product $<,>$ on $\mathbf{M}_3(\mathbb{C})$:
$<X,Y>:=\frac{1}{2}$tr$(X^TY)$. We have the isometry

$\al :\mathbb{C}^3\rightarrow\mathbf{o}_3(\mathbb{C}),\ \al(\begin{bmatrix}x^1\\
x^2\\
x^3
\end{bmatrix})=\begin{bmatrix}
0&-x^3&x^2\\
x^3&0&-x^1\\
-x^2&x^1&0
\end{bmatrix},\
x^Ty=<\al(x),\al(y)>=\frac{1}{2}$tr$(\al(x)^T\al(y)),\\
\al(x\times y)=[\al(x),\al(y)]=\al(\al(x)y)$.

Note
$<\al(x)^2,\al(y)>=\frac{1}{2}$tr$((\al(x)^2)^T\al(y)))=-\frac{1}{2}$tr$(\al(y)^T\al(x)^2)=0$.

Let $x\subset\mathbb{C}^3$ be a surface non-rigidly applicable to
a surface $x_0\subset\mathbb{C}^3$:
\begin{eqnarray}\label{eq:roll}
(x,dx)=(R,t)(x_0,dx_0):=(Rx_0+t,Rdx_0),
\end{eqnarray}
where $(R,t)$ is sub-manifold in
$\mathbf{O}_3(\mathbb{C})\ltimes\mathbb{C}^3$ (in general surface,
but it is a curve if $x_0,\ x$ are ruled and the rulings
correspond under the applicability). The sub-manifold $R$ gives
the rolling of $x_0$ on $x$, that is if we rigidly roll $x_0$ on
$x$ such that points corresponding under the applicability will
have the same differentials, $R$ will dictate the rotation of
$x_0$; the translation $t$ will satisfy $dt=-dRx_0$.

Therefore $\mathbf{O}_3(\mathbb{C})\ltimes\mathbb{C}^3$ acts on
$2$-dimensional integrable distributions of facets $(x_0,dx_0)$ in
$T^*\mathbb{C}^3$ as: $(R,t)(x_0,dx_0)=(Rx_0+t,Rdx_0)$ a rolling
is a sub-manifold
$(R,t)\in\mathbf{O}_3(\mathbb{C})\ltimes\mathbb{C}^3$ such that
$(R,t)(x_0,dx_0)$ is still integrable.

Consider first the case $[R^{-1}dR,\wedge R^{-1}dR]\neq 0$ ($R$ is
a surface).

Applying the compatibility condition $d\wedge$ to (\ref{eq:roll})
we get:
\begin{eqnarray}\label{eq:comp}
R^{-1}dR\wedge dx_0=0=dRR^{-1}\wedge dx.
\end{eqnarray}
Let $b$ be an right infinitesimal deformation of $R$:

$R_{\ep}:=R\exp(\ep b),\ \frac{d}{d\ep}|_{\ep
=0}(R_{\ep}^{-1}dR_{\ep})=db-[b,R^{-1}dR],\ 0=\frac{d}{d\ep}|_{\ep
=0}<R_{\ep}^{-1}dR_{\ep},R_{\ep}^{-1}dR_{\ep}>
=<R^{-1}dR,db>+<db,R^{-1}dR>$, so $db=[R^{-1}dR,\al(y)]$ for some
surface $y$ in $\mathbb{C}^3$, or $d\al^{-1}(b)=R^{-1}dRy$.
Imposing the compatibility condition $d\wedge$ we get
$R^{-1}dR\wedge dx_0=0$ for $x_0:=y+\al^{-1}(b)$ and this is just
(\ref{eq:comp}). Thus any infinitesimal deformation of a surface
$R$ is put into correspondence with pairs $(x_0,x)$ of applicable
surfaces in the Euclidean space $\mathbb{C}^3$ (this is
essentially due to Bianchi; see (\cite{B1},\S\ 494*)). Note that
$t:=-R\al^{-1}(b)$ will be the translation of the rolling of $x_0$
on $x:=Ry,\ dx=Rdx_0$. Also
$R^{-1}N_R:=\frac{[R^{-1}dR,R^{-1}dR]}{|[R^{-1}dR,R^{-1}dR]|}$
(the left Gau\ss\ map of $R$) is parallel to $\al(N_0),\
N_0:=\frac{dx_0\times\wedge dx_0}{|dx_0\times\wedge dx_0|}$ and
$N_RR^{-1}$ to $\al(N:=RN_0)$, so any statement about $R,\ x_0$
becomes a statement about $R^{-1},\ x$. This dichotomy into left
and right geometries is due to Clifford; he got interesting
geometric results  (Clifford left and right translations) of what
is basically currently know as
$\mathbf{o}_4(\mathbb{C})=\mathbf{o}_3(\mathbb{C})
\oplus\mathbf{o}_3(\mathbb{C})$; this can be lifted for all
practical purposes (locally) to the Lie groups and thus provided a
tool to attack the geometry of $3$-dimensional space forms. Sym's
formula from \cite{S4} basically uses a combination of both right
and left differentiation: $d(R^{-1}\pa_zR)=R^{-1}\pa_z(dRR^{-1})R$
for $R=R(u,v,z),\ d\cdot:=(\pa_u\cdot)du+(\pa_v\cdot)dv$ when
$dRR^{-1}$ is linear in $z$ (or more generally rational in $z$);
most surfaces of interest to the classical geometers and arising
as integrable systems (CGC, CMC) come in spectral families (Lie's
spectral family for CGC $-1$ and Bonnet's family for CMC); for
such families $dRR^{-1}$ is linear in $z$ and the family of
surfaces is
$R^{-1}\pa_zR\subset\mathbf{o}_3(\mathbb{R})\simeq\mathbb{R}^3$.

Note that any two infinitesimal isometries of $R$ will give pairs
of applicable surfaces such that the corresponding normals are
parallel.

Bianchi calls two surfaces $x,\ y\subset\mathbb{C}^3$ associate if
$dx\times\wedge dy=0$ (this has to do, as expected, with
infinitesimal deformations of both $x$ and $y$- see \S\
\ref{subsec:rolling5}).

Applying $R^{-1}d$ to (\ref{eq:roll}) we get
\begin{eqnarray}\label{eq:roll1}
R^{-1}d^2x=R^{-1}dRdx_0+d^2x_0.
\end{eqnarray}
Also since $N=RN_0$ we have:
\begin{eqnarray}\label{eq:secoi}
R^{-1}dRN_0=R^{-1}dN-dN_0.
\end{eqnarray}
Thus for a parametrization $(u,v)$ on $(x_0,\ x)$ with the same
off-diagonal second fundamental forms (equivalent to $R_vx_u=0$;
for example the common conjugate system on $(x_0,\ x)$ or $u\ (v)$
asymptotic parameter on $x_0\ (x)$) the axes of rotation of the
rotation $R$ along a parametric line are the tangents to $x_0$
along the other parametric line. Bianchi calls such coordinates
kinematically conjugate. The rotation of the rolling with the
other face of $x_0$ (or on the other face of $x$) is
$R':=R(I-2N_0N_0^T)=(I-2NN^T)R$.

Since $R^{-1}dR$ is skew-symmetric and using (\ref{eq:comp}) we
have
\begin{eqnarray}\label{eq:dx0}
dx_0^TR^{-1}dRdx_0=0.
\end{eqnarray}
For $a\in\mathbb{C}^3$ we have
$R^{-1}dRa=R^{-1}dR(a^{\bot}+a^{\top})=a^TN_0R^{-1}dRN_0-a^TR^{-1}dRN_0N_0=\om\times
a,\ \om:=N_0\times R^{-1}dRN_0=(\det R)R^{-1}(N\times
dN)-N_0\times dN_0$. Thus $R^{-1}dR=\al(\om)$ and
\begin{eqnarray}\label{eq:om}
d\wedge\om+\frac{1}{2}\om\times\wedge\om=0,\ \om\times\wedge dx_0=0.
\end{eqnarray}
Therefore any such flat connection form $\om$ in the tangent
bundle of $x_0$ is of the form $\om=R^{-1}dR,\ R^{-1}dR\wedge
dx_0=0$ and we get all surfaces $x$ applicable to $x_0$ by
integrating the Ricatti equation $R^{-1}dR=\om$ and a quadrature
for $dx=Rdx_0$. To show that $R^{-1}dR=\om$ is a Ricatti equation
consider $X(a,b):=R^{-1}c$, where $c$ is a constant unit vector,
so $0=R^{-1}dc=\om\times X+dX$. Multiplying on the left with
$Y(a)^T,\ Y(b)^T$ we get $Y(a)^T\om=2ida,\ Y(b)^T\om=2idb$; from
any 2 solutions $a,b$ of the Ricatti equation $Y(a)^T\om=2ida$ we
can find $X(a,b)$, then $R$ is uniquely determined by the
constraint $d(RX)=0$ and the initial condition (since $c$ depends
only on two constants, one adds another constant for the rotations
which preserve $c$). Note that the integrability condition is
satisfied: $2id\wedge
da=Y(a)^Td\wedge\om+\frac{i}{2}Y'(a)^T\om\wedge
Y(a)^T\om=Y(a)^T(d\wedge\om+\frac{1}{2}\om\times\wedge\om)=0$.
This Ricatti equation will also characterize the $B_{\infty}$
transformation of QC.

In this vein Darboux showed in (\cite{D1},\S\ 725) that the problem of
finding all surfaces applicable to a given one is equivalent to
finding all virtual asymptotic coordinates. If we infinitesimally know a
surface (that is the first and second fundamental forms), then we
need to integrate a Ricatti equation and a quadrature to find the
surface (for example from the GW equations we know $F^{-1}dF$, where
$F:=[N\ \frac{x_u}{|x_u|}\ N\times\frac{x_u}{|x_u|}]$).

If $x=x(u^1,u^2)$ is a given surface and $(\ti u^1,\ti u^2)$ are
virtual asymptotic coordinates, then keeping account of the change of
Christoffel symbols $\frac{\pa u^l}{\pa \ti
u^c}\ti\Ga_{ab}^c=(\frac{\pa^2u^l}{\pa\ti u^a\pa\ti u^b}+\frac{\pa
u^j}{\pa\ti u^a}\frac{\pa u^k}{\pa\ti u^b}\Ga_{jk}^l)$ and of the
Codazzi-Mainardi equations $d\log(K)=4(\ti\Ga_{12}^2d\ti
u^1+\ti\Ga_{12}^1d\ti u^2)$, we get the equations of virtual
asymptotic coordinates:
\begin{eqnarray}\label{eq:virta}
\frac{\pa^2u^l}{\pa\ti u^1\pa\ti u^2}+\frac{\pa u^j}{\pa\ti
u^1}\frac{\pa u^k}{\pa\ti u^2}\Ga_{jk}^l+\frac{1}{4}(\frac{\pa
u^j}{\pa\ti u^1}\frac{\pa u^l}{\pa\ti u^2}+\frac{\pa u^l}{\pa\ti
u^1}\frac{\pa u^j}{\pa\ti u^2})\frac{\pa\log(K)}{\pa u^j}=0,\
l=1,2.
\end{eqnarray}
The virtual asymptotic coordinates may not admit level curves; for
example real surfaces of positive Gau\ss\ curvature have imaginary
asymptotic directions and thus do not admit integral curves if the
surface is not analytic, but one can use in this case {\it
isothermal-conjugate systems} (conformal parametrization of the
second fundamental form); their level curves are integral curves
of the real and imaginary parts of the asymptotic directions.

Let $|dx_0|^2=:g_{jk}du^j\otimes du^k,\
g:=\frac{\det(|dx_0|^2)}{(du^1\wedge du^2)^2}$ and
$s:=s_{jk}du^j\otimes du^k:=(x_{jk}-x_{0jk})du^j\otimes
du^k:=N_0^T(R^{-1}d^2x-d^2x_0)$ be the difference of the second
fundamental forms of $x,\ x_0$. Note:
\begin{eqnarray}\label{eq:omjk}
\sqrt{g}\om=(s_{12}x_{0u^1}-s_{11}x_{0u^2})du^1+(s_{22}x_{0u^1}-s_{21}x_{0u^2})du^2;
\end{eqnarray}
($s_{12}=s_{21}$ is equivalent here to $\om\times\wedge dx_0=0$;
$d\wedge\om+\frac{1}{2}\om\times\wedge\om=0$ encodes the
difference of the GCM equations of $x_0$ and $x$).

From (\ref{eq:secoi}) we get the first and left second fundamental
forms of $R:\ |R^{-1}dR|^2=|\om|^2=|R^{-1}dRN_0|^2 =s\llcorner\pa
x_0^2\lrcorner s=s_{jk}g^{kl}s_{lm}du^j\otimes du^m,\
<R^{-1}N_R,d(R^{-1}dR)>=N_0^Td\om=N_0^T(dN_0\times
R^{-1}dRN_0)=(N_0\times
dN_0)^T(R^{-1}dN-dN_0)=\frac{1}{\sqrt{g}}(x_{01j}s_{2k}-
x_{02j}s_{1k})du^j\otimes du^k$. Note that the left second
fundamental form of $R$ is not symmetric; its symmetric part is
the second fundamental form of $R$ and its antisymmetric part
accounts for the curvature $\frac{1}{4}$ of
$\mathbf{O}_3(\mathbb{C})$.

Thus the common conjugate system on $(x_0,x)$ corresponds to
asymptotic coordinates on $R$ (Bianchi) and asymptotic coordinates
on $x_0\ (x)$ give a conjugate system on $R$. This is true again
for associated surfaces; in fact this property together with
parallel normals is Bianchi's definition of associated surfaces.

Gau\ss\ theorem states that $dN\times\wedge dN=2K(x)Nda$ ($da$ is
the area form). Because $R^{-1}N_R$ and $N_0$ are parallel, we
thus get $2K(x_0)\sqrt{\det(|dx_0|^2)}=\pm
2K(R)\sqrt{\det(|R^{-1}dR|^2)}$, so rotations of rollings of flat
surfaces in $\mathbb{C}^3$ are flat surfaces in $\mathbf{O}_3(\mathbb{C})$
(also surfaces associated to flat surfaces are flat).

We have the mean curvature
$H(R)=-\frac{\sqrt{g}}{\det{s}}g^{jk}(x_{01j}s_{2k}-x_{02j}s_{1k})$
and $db=\al(\om\times R^{-1}x),\\ \frac{d}{d\ep}|_{\ep
=0}(R_{\ep}^{-1}dR_{\ep})=\al(\om\times x_0),\
\frac{d}{d\ep}|_{\ep
=0}(R_{\ep}^{-1}N_{R_{\ep}})=\frac{[db,R^{-1}dR]+[R^{-1}dR,db]}{|[R^{-1}dR,R^{-1}dR]|}
-[b,R^{-1}N_R]=\al(N_0\times R^{-1}x-(x_0-R^{-1}x)\times
N_0)=\al(N_0\times x_0),\ \frac{d}{d\ep}|_{\ep
=0}<R_{\ep}^{-1}N_{R_{\ep}},d(R_{\ep}^{-1}dR_{\ep})>
=N_0^T(\om\times dx_0)=s$.

Therefore $H(R_{\ep})=H(R)+\ep\frac{g}{\det
s}(H(x)-H(x_0))+O(\ep^2)$, so infinitesimal deformations of an
isothermic $R$ preserving its mean curvature are put in
correspondence with Bonnet pairs in $\mathbb{C}^3$ (Bianchi
(\cite{B1},\S\ 494*)).

We are interested in the mean and Gau\ss\ curvatures $H(t),\ K(t)$
of the surface $t$ described by a point (which can be considered
the origin) rigidly attached to $x_0$ in this rolling. We thus
have $t=x-Rx_0$ and $R^{-1}dt=-R^{-1}dRx_0=-x_0^Tdx_0\llcorner
R^{-1}dR\na
x_0-x_0^TN_0R^{-1}dRN_0=-x_0^TN_0R^{-1}dRN_0+x_0^Tdx_0\llcorner
\na x_0^TR^{-1}dRN_0N_0$, so $|dt|^2=(s\llcorner\pa x_0^2\lrcorner
x_0^Tdx_0)^2+(x_0^TN_0)^2|s\llcorner\pa x_0^2\lrcorner
dx_0|^2=s\llcorner B\lrcorner s,\\
B:=(\na\frac{|x_0|^2}{2})^2+(x_0^TN_0)^2\pa
x_0^2=B^{jk}\pa_{u^j}\otimes\pa_{u^k}
=(x_0^Tg^{jl}\pa_lx_{0}x_0^Tg^{km}\pa_mx_{0}+(x_0^TN_0)^2g^{jk})\pa_{u^j}\otimes\pa_{u^k}$.
$N_t=\frac{1}{|x_0|}Rx_0$ is the normal to $t$ and the second
fundamental form of $t$ is: $N_t^Td^2t=\frac{1}{|x_0|}x_0^T
(R^{-1}d^2Rx_0-R^{-1}dRdx_0)=\frac{1}{|x_0|}(x_0^T
(R^{-1}d^2R-d(R^{-1}dR))x_0-x_0^TR^{-1}dRdx_0)
=\frac{1}{|x_0|}(-|dt|^2+x_0^TN_0s)$.

Using a kinematically conjugate system on $(x_0,x)$ we have:
$2H(t)=\frac{-1}{|x_0|}(2-\frac{x_0^TN_0s_{jj}B^{jj}}{s_{11}s_{22}(B^{11}B^{22}-(B^{12})^2)}),\\
K(t)=\frac{1}{|x_0|^2}(1+\frac{x_0^TN_0(x_0^TN_0-s_{jj}B^{jj})}{s_{11}s_{22}(B^{11}B^{22}
-(B^{12})^2)})$. But
$B^{11}B^{22}-(B^{12})^2=g^{-1}(x_0^TN_0)^2|x_0|^2$ and using the
Gau\ss\ equation for both $x_0$, $x$ we get in general
coordinates, with $\\n:=2(H(t)|x_0|+1)|x_0|^2x_0^TN_0,\
q:=(K(t)|x_0|^2-1)|x_0|^2x_0^TN_0$:
\begin{eqnarray}\label{eq:cmc}
(2K(x)n+x_{0jk}B^{jk})-
x_{jk}(n(-1)^{j+k}x_{0(j+1)(k+1)}g^{-1}+B^{jk})=0,\nonumber\\
(2K(x)q-x_0^TN_0-x_{0jk}B^{jk})-
x_{jk}(q(-1)^{j+k}x_{0(j+1)(k+1)}g^{-1}-B^{jk})=0.
\end{eqnarray}
We have
$|d\frac{x_0}{|x_0|}|^2=|dN_t|^2+\frac{1}{|x_0|^2}|dt|^2-\frac{2}{|x_0|}N_t^Td^2t=:d\sigma^2$.

Therefore a surface $t$ appears as the translation of a rolling if
$d\sigma^2$ has curvature $1$. This constitutes a second order PDE
(linear in the second order terms) for $\frac{1}{|x_0|}$; once a
solution is known, $x_0$ can be recovered by the integration of a
Ricatti equation. Any surface $t$ appears as the translation of a
rolling in a fashion depending on two arbitrary functions of one
variable and to the lines of curvature parametrization of $t$
corresponds on $x_0$ coordinates appearing orthogonal from the
origin (Bianchi (\cite{B1},\S\ 323)).

We are interested in totally real deformations $x$ of $x_0:\ \bar
x=\ep x,\ \ep:=\mathrm{diag}[\pm 1\ \pm 1\ \pm 1]$. A-priori $x_0$
does not have to be totally real; since $|dx|^2$ is real, so it
must be $|dx_0|^2$. Consider $x_0\subset\mathbb{C}^3$ surface with
real linear element of signature $\ep':=\mathrm{diag}[\ep_1\
\ep_2]$; by Sylvester's it must have real dimension $2$ choosing
$v_j:=\sqrt{\ep_j}e_j,\ j=1,2,3,\ V:=[v_1\ v_2\ v_3]$, ($\ep_3$ is
undetermined at this point), we have $dx_0=RV\psi,\
R\subset\mathbf{O}_3(\mathbb{C})$ where $\psi$ is $\mathbb{R}^3$
valued $1-$ form, $e_3^T\psi=0$. If $dx_0=R'V\psi'$ is a similar
decomposition, then $R'^TR\subset\mathbf{O}_3(\mathbb{C}),\
R'^TRe_3=\pm e_3,\ \overline{R'^TR}=\ep R'^TR\ep,\ \psi'=r\psi,\
r:=\ep VR_1^TRV$ is a curve of real orthogonal transformations of
$\mathbb{R}^3$ preserving the signature $\ep$, $re_3=\pm e_3$, so
we can consider $r\subset\mathbf{O}_2(\mathbb{R})$ if $\det\ep'=1$
and $r\subset\mathbf{O}_{1,1}(\mathbb{R})$ if $\det\ep'=-1$. We
have $d\bar x_0=\bar R\ep V\psi= \bar R\ep R^Tdx_0$, so $\bar R\ep
R^T$ (which does not change if $R$ changes as above) is the
rotation of the rolling of $x_0$ on $\bar x_0$. Note that if
$d\bar x_0=Rdx_0,\ R\subset\mathbf{O}_3(\mathbb{C})$, then $\bar
R^T=R$, so such surfaces $x_0$ are put in correspondence with
special right infinitesimal deformations of surfaces
$R\subset\mathbf{O}_3(\mathbb{C}),\ \bar R^T=R$: with
$y:=R^{-1}\bar x_0$ we have $d(R^{-1}\bar y)-R^{-1}d(Ry)=0$;
equivalently the special right infinitesimal deformation
$b:=\al(Ry-\bar y)$ of $R$ satisfies its requirement
$db=[R^{-1}dR,\al(y)]$. If $d\bar x_0=Rdx_0,\
R\in\mathbf{O}_3(\mathbb{C})$, then $\bar x_0=(R,t)x_0$, $(R,t)$
rigid motion with $\bar R^T=R,\ \bar t+\bar Rt=0$. Should $x_0$ be
totally real, then there would be signature $\ep$ and rigid motion
$(R_0^T,t_0)$ such that
$\overline{(R_0^T,t_0)x_0}=\ep((R_0^T,t_0)x_0)$; equivalently
$R=\bar R_0\ep R_0^T,\ \bar t_0-\ep t_0+\bar R_0^Tt=0$; once we
found $R_0$ with $R=\bar R_0\ep R_0^T$, $\bar R_0^Tt$ will satisfy
the compatibility condition $\overline{\bar R_0^Tt}=-\ep\bar
R_0^Tt$ needed in order to prescribe $t_0$. Note that $R_0,\ t_0$
can be prescribed modulo rigid motions of the corresponding
totally real space, which satisfy $\bar R_0=\ep R_0\ep,\ \bar
t_0=\ep t_0$. Now the signature $\ep'$ can be found from the
linear element of $x_0$, the last component $\ep_3$  of the
signature $\ep$ (so far undetermined) can be found from $\det
R=\det\ep$. The usual argument with the dimension of kernel and
codomain and connectedness of the codomain of the map
$\mathbf{O}_3(\mathbb{C})\ni R_0 \mapsto\bar R_0\ep R_0^T$ applies
to prove the existence of $R_0$.

Thus all totally real surfaces $x_0$ are recognized from $d\bar x_0=Rdx_0,\
R\in\mathbf{O}_3(\mathbb{C}),\ \bar R^T=R$.

Note that there are non-totally real surfaces with real linear
element. To see this note that there are curves in arc-length
parametrization (that is $c=c(s),\ s\in\mathbb{R},\ |dc|^2=\pm
ds^2$) which are not totally real. All real and real analytic
$2$-dimensional linear elements admit local isometric embedding in
totally real spaces and depending on two functions of a variable
as long as the signature of the ambient space allows the signature
of the linear element. If one takes a real analytic curve in
arc-length parametrization on the surface and one deforms it to a
general position (non-totally real), then by Cauchy-Kovalewskaia a
small portion of the surface around the curve can be deformed to
the new general position: thus is non-totally real. In general
Cauchy-Kovalewskaia fails for characteristic data either in
non-existence or non-uniqueness of solutions, but in our case it
fails in the good way: if the general position of the deformed
curve is asymptote on the new surface, it is asymptote on an
$1$-dimensional family of deformations of the considered surface
(Darboux proved that the characteristics of the deformation
problem are the asymptotes; thus it becomes clear that Bianchi's
notion of W congruence is the correct tool in studying the
deformation problem since although it does not give in general the
applicability correspondence, it gives the correspondence of
characteristics). Here Eisenhart in (\cite{E1},\S\ 139) puts the
condition that for real surfaces and deformations the curvature of
the deformed curve should be greater than or equal to the geodesic
curvature of the curve on the surface (equality is obtained when
the deformed curve becomes asymptote on the deformed surface). For
the same problem (real surface and deformation of a curve on the
surface) Darboux does not put this condition; he just mentions
instead that imaginary deformations are obtained if this condition
is not satisfied. If one takes in Darboux's imaginary deformation
case the deformed curve in discussion to be non-planar (one can
change the ambient totally space along planes), then the imaginary
surface thus obtained is non-totally real and has positive
definite linear element.

The unit normal field $N_0$ of a totally real surface $x_0\subset
V$, where $\bar v=\ep v,\ \forall v\in V$, satisfies $N_0\subset
V$ if $\ep_3=1$ (in which case the second fundamental form is
real) and $N_0\subset iV$ otherwise (in which case the second
fundamental form is purely imaginary). Thus $x_0$ has asymptotic
lines iff the sign of the Gau\ss\ curvature $K_0$ is
$-\ep_1\ep_2\ep_3$; if $K_0\equiv 0$, then the asymptotic lines
coincide (the rulings of the developable $x_0$).

Take for example surfaces $x_0\subset\mathbb{R}^2\times
i\mathbb{R},\ x\subset\mathbb{R}^3,\ (x,dx)=(R,t)(x_0,dx_0),\
(R,t)\subset\mathbf{O}_3(\mathbb{C})\ltimes\mathbb{C}^3$. Then
with $\ep:=\mathrm{diag}[1\ 1\ -1]$ we have
$iN_0\subset\mathbb{R}^2\times i\mathbb{R},\
\ep(I_3-2N_0N_0^T)\bar R^TR=I_3$. If $K_0>0$, then $x_0$ has
asymptotic lines while $x$ does not; if $K_0<0$, then $x$ has
asymptotic lines while $x_0$ does not.

We shall now prove that the only surfaces in unit complex
(pseudo-)spheres having real linear element are the totally real
(pseudo-)spheres. Assume $x\subset\mathbb{C}^3,\ |x|^2=\ep_1,\
\ep_1:=\pm 1,\ d\bar x=Rdx,\ R\subset\mathbf{O}_3(\mathbb{C}),\
\bar R^T=R$. Thus $\bar N=\ep_2 RN$, where $\ep_2:=\pm 1,\
|N|^2=1,\ N^Tdx=0$; since $N=\sqrt{\ep_1}x$, we get $\bar
x=\ep_1\ep_2 Rx$. If $\ep_1\ep_2=-1$, then one can replace $R$
with $R(I_3-2NN^T) =R(I_3-2\ep_1 xx^T)$ to get $\ep_1\ep_2=1$; if
$\ep_1\ep_2=1$, then $dRx=0$ so $0=R^{-1}dRx=\om\times x$, or
$\om\times N=0$; since $\om=\om^\top$, we get $\om=0$ and $R$ is
constant.

\subsection{Curves corresponding under the applicability}\label{subsec:rolling3}
\noindent

\noindent Consider $c_0,\ c$ curves on $x_0,\ x$ corresponding
under the applicability. We would like to see how the curvatures
and torsions transform under the applicability. If
$F_c:=[T_c:=\frac{dc}{|dc|}\ \ N_c:=\frac{(dc\times d^2c)\times
dc}{|(dc\times d^2c)\times dc|}\ \ B_c:=\frac{dc\times
d^2c}{|dc\times d^2c|}],\ \ti F_c:=[T_c\ \ N\ \ T_c\times N],\
k_c:=\frac{|dc\times d^2c|}{|dc|^3},\ \tau_c:=\frac{(dc\times
d^2c)^Td^3c}{|dc\times d^2c|^2},\ k_{nc}:=\frac{N^Td^2c}{|dc|^2},\
k_{gc}:=-\frac{(dc\times d^2c)^TN}{|dc|^3},\
\tau_{gc}:=\frac{(N\times dN)^Tdc}{|dc|^2}$ are the usual frames
along $c$, curvatures and torsions of $c$ (and similar quantities
for $c_0$), then we have $F_c^{-1}dF_c=|dc|\al([\tau_c\ 0\
k_c]^T),\ \ti F_c^{-1}d\ti F_c=|dc|\al([\tau_{gc}\ -k_{gc}\
k_{nc}]^T),\ |dc|\ti F_c^{-1}N_c=\frac{1}{k_c}\ti
F_c^{-1}dT_c=|dc|[0\ \frac{k_{nc}}{k_{c}}\ \frac{k_{gc}}{k_c}]^T,\
|dc|\ti F_c^{-1}B_c=\frac{1}{\tau_c}(d(\ti F_c^{-1}N_c)+\ti
F_c^{-1}d\ti F_c\ti F_c^{-1}N_c+|dc|k_c\ti
F_c^{-1}T_c)=\frac{1}{\tau_c}[0\ \ \
d(\frac{k_{nc}}{k_c})-|dc|\frac{\tau_{gc}k_{gc}}{k_c}\ \ \
d(\frac{k_{gc}}{k_c})+|dc|\frac{\tau_{gc}k_{nc}}{k_c}]^T$.
Therefore $|dc|\ti F_c^{-1}F_c=\begin{bmatrix}
|dc|&0&0\\
0&|dc|\frac{k_{nc}}{k_c}&\frac{1}{\tau_c}(d(\frac{k_{nc}}{k_c})-|dc|\frac{\tau_{gc}k_{gc}}{k_c})\\
0&|dc|\frac{k_{gc}}{k_c}&\frac{1}{\tau_c}(d(\frac{k_{gc}}{k_c})+|dc|\frac{\tau_{gc}k_{nc}}{k_c})
\end{bmatrix}=|dc|\begin{bmatrix}
1&0&0\\
0&\frac{k_{nc}}{k_c}&-\frac{k_{gc}}{k_c}\\
0&\frac{k_{gc}}{k_c}&\frac{k_{nc}}{k_c}
\end{bmatrix}$ and we have
$|dc|\tau_c=|dc|\tau_{gc}+\frac{k_{nc}dk_{gc}-k_{gc}dk_{nc}}{k_c^2}$
(Bonnet).

We have $\ti F_c=R|_{c_0}\ti F_{c_0}$, so $\ti F_c^{-1}d\ti
F_c=\ti F_{c_0}^{-1}d\ti F_{c_0}+\ti
F_{c_0}^{-1}R^{-1}dR|_{c_0}\ti F_{c_0}$. Therefore
\begin{eqnarray}\label{eq:tors}
k_{nc}=k_{nc_0}+(\frac{dc_0}{|dc_0|}\times N_0)^T\frac{\om|_{c_0}}{|dc_0|},\nonumber\\
k_{gc}=k_{gc_0},\
\tau_{gc}=\tau_{gc_0}-\frac{dc_0^T}{|dc_0|}\frac{\om|_{c_0}}{|dc_0|}.
\end{eqnarray}
Note that if we choose $(u,v)$ coordinates on $x$ such that
$c(u)=x(u,0)$ and $v=$const are parallel, $u=$const are geodesics,
then $\tau _{gc}=\frac{N^Tx_{uv}}{|x_u||x_v|}$, so
$\tau_{gc}=\tau_c=\sqrt{-K}$ for $c$ asymptote on $x$ (Enneper).
In particular the asymptotes of CGC surfaces are curves of
constant torsion.

Conversely, we can compute the curvatures and torsions of the
corresponding curve $c_R$ on $R$ in function of the curvatures and
torsions of $c_0,\ c$:

\noindent$|dc_R|=|c_R^{-1}dc_R|=|\om|=\sqrt{(\tau_{gc}
-\tau_{gc_0})^2+(k_{nc}-k_{nc_0})^2}|dc_0|,\
k_{nc_R}=\frac{k_{nc_0}\tau_{gc}-\tau_{gc_0}k_{nc}}{(\tau_{gc}
-\tau_{gc_0})^2+(k_{nc}-k_{nc_0})^2},\\
k_{gc_R}|dc_0|=\frac{(k_{nc}-k_{nc_0})d(\tau_{gc}-\tau_{gc_0})
-(\tau_{gc}-\tau_{gc_0})d(k_{nc}-k_{nc_0})}
{((\tau_{gc}-\tau_{gc_0})^2+(k_{nc}-k_{nc_0})^2)^{\frac{3}{2}}},\
\tau_{gc_R}=\frac{\tau_{gc_0}(\tau_{gc}-\tau_{gc_0})+k_{nc_0}(k_{nc}-k_{nc_0})}
{(\tau_{gc}-\tau_{gc_0})^2+(k_{nc}-k_{nc_0})^2}.$

\subsection{Ruled surfaces}\label{subsec:rolling4}
\noindent

\noindent If $[R^{-1}dR,\wedge R^{-1}dR]=0$, then by a change of
coordinates we can suppose $R_u=0$, so $x_{0u}=\phi_{u}a(v),\
R^{-1}dRa=0,\ |a|^2=1\ \mathrm{or}\ a=Y(v),\ \phi_{u}\neq 0$ and
$x_0(u,v)=\phi(u,v)a(v)+b(v)$, $b(v)$ arbitrary such that
$a\times(\phi a_v +b_v)\neq 0$. With the change of coordinates
$(u',v')=(\phi(u,v),v)$, we get the formula for ruled surfaces:
$x_0(u,v)=ua(v)+b(v)$ (for convenience after all these changes of
coordinates we still keep the notation $(u,v)$ for our
parametrization).

In this case $x=uRa+c$, where $c(v)=\int Rdb=Rb-\int\varphi
R(a\times b)|db|$ and $R$ is a solution of the ODE
$R^{-1}dR=\varphi|db|\al(a)$. The line of striction $b$ of $x_0$
is defined as $da^Tdb=0$ (the curve along which the rulings are
closest) can be realized after a change of the $u$ parameter. In
particular if $x_0$ is developable, then its line of striction
coincides with the line of regression. Because
$d(Ra)^TRdb=db^T(da+R^{-1}dRa)=0$, the line of striction is
invariant under deformations preserving the rulings.

We can prescribe various conditions on $c$ by choosing suitable
$\varphi$'s (see Eisenhart (\cite{E1},\S\ 144)).

In our case $N_0=\frac{a\times db}{|a\times db|}$, so we have
$k_{nc}=\frac{d^2b^T(a\times db)}{|a\times
db||db|^2}+\varphi\frac{|a\times db|}{|db|},\
\tau_{gc}=(1-\frac{a^Tdb}{|db|})\frac{(a\times
db)^T}{|db|}(\frac{a^Tdbd^2b}{|db|^3}-
(1+\frac{a^Tdb}{|db|})\frac{da}{|db|})-\varphi\frac{a^Tdb}{|db|}$.

Therefore we can prescribe all quantities except the torsion by
solving an ODE; for the torsion we need to solve two ODE's.

According to Bonnet, if two surfaces $x_0,\ x$ are applicable with
corresponding asymptotes of parameter $u$, then
$R_ux_{0u}=R_ux_{0v}=0$, so $R_u=0$, so $x_0,\ x$ are ruled with
the rulings corresponding under the applicability.

Also (see Eisenhart (\cite{E1},\S\ 144)) if $x_0,\ x$ are ruled
and applicable with rulings not corresponding under the
applicability, then choosing parametrization $(u,v)$ such that the
rulings are the curves $u=$const, $v=$const, $(u,v)$ will be
virtual asymptotes which are geodesics. Thus there is a surface on
which $(u,v)$ are asymptotes and geodesics, so $x_0,\ x$ are
applicable (with correspondence of rulings) to a quadric.

According to a theorem of Chieffi, if we have a deformation of a
ruled surface, then the ruled surface formed by the tangents (to
the geodesics corresponding to rulings) along an asymptotic line
is applicable to the initial ruled surface.

If $x_0(u_0,v_0)=u_0a(v_0)+b(v_0),\ (x,dx)=(R,t)(x_0,dx_0)$ and
$(u,v)$ is the parametrization of $x$ with $v$ asymptotic
parameter, then setting $u:=$const the ruled surface $\ti
x(w,v):=wx_{u_0}+x=R(wa+x_0)+t$ is applicable to $\ti
x_0(w,v):=wa+x_0=(w+u_0)a+b=x_0(w+u_0,v_0)$, because $d\ti x
=Rd\ti x_0+wdRx_{0u_0}=Rd\ti x_0$. Here we used the fact that
$dRx_{0u_0}=R_vx_{0u_0}dv=0$, which follows from (\ref{eq:roll1})
(because $u_0=$const, $v=$const are asymptotes on $x_0$
respectively $x$, the Gau\ss\ curvature information is recorded
only in the off-diagonal terms of the second fundamental forms of
$x_0,\ x$, so
$0=N_0^T(R^{-1}x_{u_0v}-x_{0u_0v})=R^{-1}R_vx_{0u_0}$).

Because $u=$const is an asymptote on $x$, its osculating spaces
coincide with the tangent spaces of $x$ along $u=$const; therefore
$u=$const is an asymptote also on $\ti x$.

For example, because the catenoid is applicable to the minimal
helicoid, it is applicable to the ruled surface formed by the
tangents (to the geodesics corresponding to the rulings of the
minimal helicoid) along one of its asymptotes. This example (for
surfaces applicable to the catenoid) is due to Bianchi and was the
inspiration of Chieffi's result.

Note that we have the more general result, due to Bianchi: if $x,\
x_0$ are applicable surfaces with rolling $(x,dx)=(R,t)(x_0,dx_0)$
and we have $(u_0,v_0),\ (u,v)$ parameterizations on $x_0,\ x$
with $(u_0,v)$ kinematically conjugate, then setting $u:=$const
the ruled surface $\ti x (w,v):=wx_{u_0}+x$ is applicable to the
ruled surface $\ti x_0(w,v):=wx_{0u_0}+x_0$ with the same rolling
as that of $(x_0,\ x)$ and correspondence of rulings. We have $\ti
x=R\ti x_0+t,\ d\ti x=Rd\ti x_0+wdRx_{0u_0}=Rd\ti x_0$, because
$dRx_{0u_0}=R_vx_{0u_0}dv=0$. Darboux in (\cite{D1},\S\ 932) has a
result in the same vein: if $(c_0,\ c)$ are curves of $x_0,\ x$
corresponding under the applicability, then since
$0=R^{-1}dR|_{c_0}\frac{\om|_{c_0}}{|dc_0|}$, the ruled surfaces
$c_0+w\frac{\om|_{c_0}}{|dc_0|},\ c+wR\frac{\om|_{c_0}}{|dc_0|}$
are applicable with rolling $(R,t)|_{c_0}$.

Conversely, if $(u_0,v_0),\ (u,v)$ are  a-priori  arbitrary
parameterizations on $x_0,\ x$ and a point-wise correspondence is
established between $x_0,\ x$ such that for any $u=$const the
ruled surface $\ti x(w,v):=wx_{u_0}+x$ is applicable to the ruled
surface $\ti x_0(w,v):=wx_{0u_0}+x_0$, $(\ti x,d\ti x)=(R,t)(\ti
x_0,d\ti x_0)$ with correspondence of rulings, then
$x_{u_0}dw+wdx_{u_0}+dx=d\ti x=Rd\ti
x_0=Rx_{0u_0}dw+wRdx_{0u_0}+Rdx_0$, so $x_{u_0}=Rx_{0u_0},\
wx_{u_0v}+x_v=wRx_{0u_0v}+Rx_{0v}$. Setting $w=0$ we get
$x_v=Rx_{0v}$, so  $x$ is applicable to $x_0$. Moreover we get
$x_{u_0v}=Rx_{0u_0v}$, so $R_vx_{0u_0}=0$ and thus $(u_0,v)$ are
kinematically conjugate.

In this vein according to Bianchi (122) Lie has reduced the study
of B transformations of the (pseudo-)sphere to the study of B
transformations of their asymptotes, curves of constant torsion.
Also, according to Bianchi (122,\S\ 67), Serret has integrated the
equation of ruled CGC surfaces: all are obtained by choosing the
ruling as one of the isotropic directions in the osculating plane
of a curve of constant torsion. This can be easily justified with
Chieffi's and Enneper's results; such surfaces $x$ are of the form
$x(u,v)=\frac{1}{u-v}R(v)Y(v)-\int R(v)f_1dv,\
R\subset\mathbf{O}_3(\mathbb{C}),\ R'(v)Y(v)=0$ for the CGC $1$
case and $ix$, $x$ as above for the CGC $-1$ case. In particular
if we let $u\bar v=-2\ (=2)$ for the CGC $1\ (-1)$ case, then the
linear element of $x$ ($ix$) is real positive definite; such
surfaces $x\ (ix)$ are real analytic. Further for the CGC $1\
(-1)$ case we have $d\bar x=\overline{R(v)}R(v)^Tdx\ \ \
(d\overline{(ix)}=\overline{R(v)}(I_3-2e_3e_3^T)R(v)^Td(ix))$;
should $\overline{R(v)}R(v)^T\ \ \
(\overline{R(v)}(I_3-2e_3e_3^T)R(v)^T)$ be constant,
$\overline{R(v)}$ would be holomorphic in $v$ and in $\bar v$, so
it would be a constant. Thus $x\ (ix)$ is non-totally real except
when $x\ (ix)$ is up to rigid motions the standard totally real
unit (pseudo-)sphere.

\subsection{Infinitesimal deformations, Weingarten congruences, Darboux's $12$ surfaces
and transformations of surfaces}\label{subsec:rolling5}
\noindent

\noindent (see also Sabitov \cite{S1})

Let $s_1$ be an infinitesimal deformation of the surface $s_0:\
0=\frac{d}{d\ep}|_{\ep =0}|d(s_0+\ep s_1)|^2=ds_0^Tds_1
+ds_1^Tds_0$. Therefore there is $s_2$ (the infinitesimal rotation
of the infinitesimal deformation) such that $ds_0=s_2\times ds_1$.
With $s_3:=s_0-s_2\times s_1$ we have $ds_3=s_1\times ds_2$ and
$ds_2^Tds_3+ds_3^Tds_2=0$, so there is $s_4$ such that
$ds_2=s_4\times ds_3$ and we can iterate this construction
(replace $s_k,\ k=0,1$ with $s_{2j+k}, j=1,2,...$). Darboux proved
in (\cite{D1},\S\ 883-\S\ 924) that this process ends after 12
steps ($s_j=s_{j+12}$) and found interesting relationships between
these 12 surfaces.

We have $ds_{2j+3}=s_{2j+1}\times
ds_{2j+2}=s_{2j+1}\times(s_{2j+4}\times
ds_{2j+3})=-s_{2j+1}^Ts_{2j+4}ds_{2j+3}$, so
$s_{2j+1}^Ts_{2j+4}+1=0$. Also $ds_{2j+2}=s_{2j+4}\times
ds_{2j+3}=s_{2j+4}\times(s_{2j+1}\times
ds_{2j+2})=-s_{2j+1}^Ts_{2j+4}ds_{2j+2}+ds_{2j+2}^Ts_{2j+4}s_{2j+1}$,
so $ds_{2j+2}^Ts_{2j+4}=0$.

But $ds_{2j+2}\times\wedge ds_{2j+1}=d\wedge ds_{2j}=0$, so
$s_{2j+1},\ s_{2j+2}$ are associate; in particular they have
parallel tangent planes, so $ds_{2j+1}^Ts_{2j+4}=0$, so
$s_{2j+1},\ s_{2j+4}$ are {\it polar reciprocal with respect
to the pseudo-sphere} (prwrtps).

From $s_{2j+4}^Tds_{2j+1}=s_{2j-1}^Tds_{2j+1}=0$ we get
$s_{2j+4}=\varphi s_{2j-1}$ for some function $\varphi$.
Multiplying this last relation on the left with $s_{2j+1}^T,\
s_{2j+2}^T$ we get
$s_{2j+4}=-\frac{s_{2j-1}}{s_{2j+1}^Ts_{2j-1}},\
s_{2j-1}=-\frac{s_{2j+4}}{s_{2j+2}^Ts_{2j+4}}$. Therefore
$s_{2j+3}=s_{2j}+\frac{s_{2j-3}\times
s_{2j+1}}{s_{2j-3}^Ts_{2j-1}}$ and $s_{2j-3}^Ts_{2j+3}+1=0$. Then
$s_{2j+7}=s_{2j+4}+\frac{s_{2j+1}\times
s_{2j+5}}{s_{2j+1}^Ts_{2j+3}}=-\frac{s_{2j-1}}{s_{2j+1}^Ts_{2j-1}}+
\frac{s_{2j+1}\times(s_{2j-1}^Ts_{2j+1}s_{2j+2}+s_{2j-1}\times
s_{2j+3})}{s_{2j+1}^Ts_{2j+3}s_{2j-1}^Ts_{2j+1}}
=-\frac{s_{2j-1}}{s_{2j+1}^Ts_{2j-1}}+\frac{s_{2j+1}\times
s_{2j+2}}{s_{2j}^Ts_{2j+1}}-\frac{s_{2j+3}}{s_{2j+1}^Ts_{2j+3}}
+\frac{s_{2j-1}}{s_{2j-1}^Ts_{2j+1}}=-\frac{s_{2j}}{s_{2j}^Ts_{2j+1}}=
-\frac{s_{2j}}{s_{2j}^Ts_{2j-2}}=s_{2j-5}$ and
$s_{2j+6}=-\frac{s_{2j+1}}{s_{2j+3}^Ts_{2j+1}}
=-\frac{s_{2j-11}}{s_{2j-9}^Ts_{2j-11}}=s_{2j-6}$.

The prwrtps is uniquely determined:
$s_{2j+1}=\frac{-N_{2j+4}}{s_{2j+4}^TN_{2j+4}},\
s_{2j+4}=\frac{-N_{2j+1}}{s_{2j+1}^TN_{2j+1}}$; $s_{2j+1}$ and
$s_{2j+2}$ are associate (parallel tangent planes and to
asymptotes on one surface corresponds on the other a conjugate
system). Conversely from two associate surfaces $s_1,\ s_2\
(ds_1\times\wedge ds_2=0$) we can recover $s_0$ by a quadrature
from $ds_0=s_2\times ds_1$.

We have $s_{2j+3}-s_{2j}=s_{2j+1}\times s_{2j+2}$; since
$s_{2j+2}^Tds_{2j}=s_{2j+1}^Tds_{2j+3}=0$ we see that
$s_{2j+3}-s_{2j}$ is tangent to both $s_{2j},\ s_{2j+3}$. With
$\phi:=s_{2j+2}^Ts_{2j+4}$ we have $0=ds_{2j+3}^Tds_{2j+1}(\phi
s_{2j+1} +s_{2j+2})^Ts_{2j+4}=(s_{2j+4}\times ds_{2j+3})^T((\phi
s_{2j+1}+s_{2j+2})\times ds_{2j+1})=ds_{2j+2}^T(\phi
s_{2j+1}\times ds_{2j+1}+ds_{2j})=\phi
d^2s_{2j+3}^Ts_{2j+1}-d^2s_{2j}^Ts_{2j+2}$; therefore the second
fundamental forms of $s_{2j},\ s_{2j+3}$ are proportional and
$s_{2j},\ s_{2j+3}$ are the focal surfaces of a W congruence.

Conversely Guichard proved that if $s_0, \ s_3$ are the focal
surfaces of a W congruence, then we can recover infinitesimal
deformations $s_1,(s_2)$ of $s_0,(s_3)$ such that
$s_1^Tds_3=s_2^Tds_0=0,\ s_3-s_0=s_1\times s_2$. These $s_j,\
j=0,3$ will satisfy $ds_0-s_2\times ds_1=ds_3-s_1\times ds_2,\
0=s_2^T(ds_0-s_2\times ds_1)=s_1^T(ds_3-s_1\times ds_2)
=ds_1^T(ds_0-s_2\times ds_1+ds_1^T(ds_0-s_2\times
ds_1)=ds_2^T(ds_3-s_1\times ds_2)+ds_2^T(ds_3-s_1\times ds_2)$, so
$0=ds_0-s_2\times ds_1=ds_3-s_1\times ds_2$.

To prove Guichard's result we refer $s_0,\ s_3$ to their
asymptotic coordinates $u,\ v:\ \theta_j:=|\theta_j|N_j,\
|\theta_j|^2:=\frac{|s_{ju}|}{|N_{ju}|}=\frac{|s_{jv}|}{|N_{jv}|}$,
so $ds_j=\theta_j\times(\theta_{ju}du-\theta_{jv}dv), \
N_j=\frac{\theta_j}{|\theta_j|},\ \theta_j\times \theta_{juv}=0,\
j=0$ or $3$ (Lelieuvre formulae). All infinitesimal deformations
$s_1$ of $s_0$ are given by:
$ds_1=\theta_0(\ti\theta_{0u}du-\ti\theta_{0v}dv)-\ti\theta_0(\theta_{0u}du-\theta_{0v}dv)$,
where $\ti\theta_0,\ (\theta_0)$ are scalar (vector) solutions of
the Laplace equation $\theta_{uv}=M_0(u,v)\theta$; then
$s_2:=\frac{\theta_0}{\ti\theta_0}$ satisfies $ds_0-s_2\times
ds_1=0$. We have $s_3-s_0=\rho\theta_0\times\theta_3$ for some
function $\rho$, so $\rho
d\theta_0^T(\theta_0\times\theta_3)=-\theta_0^Td(s_3-s_0)=
-(\theta_0\times\theta_3)^T(\theta_{3u}du-\theta_{3v}dv)$ and
$\rho d\theta_3^T(\theta_0\times\theta_3)=-\theta_3^Td(s_3-s_0)
=-(\theta_0\times\theta_3)^T(\theta_{0u}du-\theta_{0v}dv)$.
Therefore
$(\rho^2-1)d\theta_0^T(\theta_0\times\theta_3)=(\rho^2-1)d\theta_3^T(\theta_0\times
\theta_3)=0$ and $\rho=1$ (we can absorb the sign in $\theta_0$;
also $\theta_0\times\theta_3=0$ is impossible). Thus
$0=ds_3-ds_0-d\theta_0\times\theta_3-\theta_0\times
d\theta_3=(\theta_3-\theta_0)\times(\theta_3+\theta_0)_udu-
(\theta_3+\theta_0)\times(\theta_3-\theta_0)_vdv$ and
$(\theta_3+\theta_0)_u=\log(\ti\theta)_u(\theta_3-\theta_0),\
(\theta_3-\theta_0)_v=\log(\theta)_v(\theta_3+\theta_0)$ for some
functions $\ti\theta,\ \theta$. Differentiating these two last
relations wrt $v,\ u$ and using the relations themselves we get
$0=\theta_3(M_3-\log(\ti\theta)_{uv}-\log(\ti\theta)_u\log(\theta)_v)+
\theta_0(M_0+\log(\ti\theta)_{uv}-\log(\ti\theta)_u\log(\theta)_v)=
\theta_3(M_3-\log(\theta)_{uv}-\log(\ti\theta)_u\log(\theta)_v)-
\theta_0(M_0+\log(\theta)_{uv}-\log(\ti\theta)_u\log(\theta)_v)$,
or $(\log(\theta)-\log(\ti\theta))_{uv}=0,\
M_0=-\log(\ti\theta)_{uv}+\log(\ti\theta)_u\log(\theta)_v,\
M_3=\log(\ti\theta)_{uv}+\log(\ti\theta)_u\log(\theta)_v$.
Therefore
$\ti\theta_3:=f(v)\ti\theta=g(u)\theta=:\frac{1}{\ti\theta_0}$ for
some functions $f,\ g,\ \ti\theta_{0uv}=M_0\ti\theta_0,\
\ti\theta_{3uv}=M_3\ti\theta_3$ and
$s_1:=\frac{\theta_3}{\ti\theta_3},\
s_2:=\frac{\theta_0}{\ti\theta_0}$ are the surfaces we are looking
for. To the geometric W congruence corresponds analytically the
Moutard transformation of solutions of the Laplace equations with
functions $M_0,\ M_3$. Note that $K(s_0)=-\frac{1}{|\theta_0|^4}$,
so a congruence is W iff its focal surfaces $s_0,\ s_3$ satisfy
$K(s_0)K(s_3)=\frac{|N_0\times N_3|^4}{|s_0-s_3|^4}$ (Ribacour).

Thus the infinitesimal deformation problem becomes equivalent to W
congruences. These in turn become {\it Ribacour sphere
congruences} (2-dimensional families of spheres upon whose
envelopes the lines of curvature correspond) via Lie's contact
transformation. The envelopes of the Ribacour sphere congruences
correspond to the focal surfaces of the W congruence. Any
transformation preserving W congruences (homographies,
correlations) preserves the problem of infinitesimal deformation.

If we take any two surfaces of the 12 and which are not prwrtps
one of the other, then by taking their prwrtps we get four
surfaces among which there will be 2 who are in the relation: one
infinitesimal deformation of the other, associate or focal
surfaces of a W congruence, so we can recover the 12 surfaces.

Since a linear isomorphism of a plane preserves the cross-ratio of
four lines through the origin, a correspondence between two
surfaces preserves the cross-ratio of four vector fields. Thus if
asymptotic directions correspond under the correspondence (which
is the case for the focal surfaces of a W congruence), then all
conjugate systems correspond (because conjugate directions and
asymptotic directions have cross-ratio $1$).

On the 12 surfaces we have three systems of coordinates of
interest; each will give on each surface asymptotic or conjugate
systems with equal (tangential) invariants (that is
$\Ga_{12v}^2=\Ga_{12u}^1$, where the Christoffel symbols are taken
wrt the linear element of the surface (Gau\ss\ map)). On two
surfaces one infinitesimal deformation of the other the asymptotic
coordinates on one correspond to a conjugate system with equal
invariants on the other and the common conjugate system has equal
tangential invariants on both. On two associate surfaces the
asymptotic coordinates on one correspond to a conjugate system
with equal tangential invariants on the other and the common
conjugate system has equal invariants on both surfaces. On the
focal surfaces of a W congruence (or on a surface and its prwrtps)
the asymptotic coordinates correspond and the conjugate system
with equal invariants on one corresponds to a conjugate system
with equal tangential invariants on the other. Knowledge of a
conjugate system with equal invariants gives by a quadrature the
associate surface (unique modulo translations) upon which the
given conjugate system is still conjugate with equal invariants,
so knowledge of such systems is equivalent to the problem of
infinitesimal deformation.

Consider the congruence $s_1+t(s_2-s_1)$; it cuts its focal
surfaces and $s_1,\ s_2$ with cross-ratio $1$. The developables of
the congruence (inducing conjugate systems with equal invariants
on the focal surfaces) correspond to the common conjugate system
on $s_1,\ s_2$; therefore Koenigs considered such surfaces $s_1,\
s_2$ as transforms one of the other (a particular case is the D
transformation of isothermic surfaces). Eisenhart in \cite{E2} and
Jonas have generalized this transformation to the {\it
fundamental} (F) transformation (they do not require equal
invariants). According to Eisenhart \cite{E2} this transformation
F together with Bianchi's generalized {\it B\"{a}cklund} (B)
transformation (the focal surfaces of a W congruence are
considered as transforms one of the other) account for most if not
all of transformations of surfaces for which an interesting theory
has been developed (including for example variants of the BPT).

\subsection{Rolling congruences}\label{subsec:rolling6}
\noindent

\noindent All results of this subsection in the particular case
when the curve $c$ is a line appear in Bianchi (\cite{B2}, Vol{\bf
7},(157),(166)).

Let $c:=c(w)$ be a curve in $\mathbb{C}^3$. As $x_0$ rolls on $x$
and rigidly moving $c$ the curve $c$ will describe a congruence of
congruent curves $(R(u,v),t(u,v))c(w)$. We are interested in the
focal surfaces (envelopes) of such a congruence, that is
prescribing $w(u,v)$ such that the surface
$(R(u,v),t(u,v))c(w(u,v))$ meets tangentially the curve
$w\rightarrow(R(u,v),t(u,v))c(w)$. We thus need
$\\0=(Rdc)^T(d((R,t)c)\times\wedge d((R,t)c))
=dc^T((dc+\om\times(c-x_0))\times\wedge(dc+\om\times(c-x_0)))
=dc^T((\om\times(c-x_0))\times\wedge(\om\times(c-x_0)))
=N_0^T(\om\times\wedge\om)N_0^T(c-x_0) dc^T(c-x_0)$. Therefore the
focal surfaces are given by $N_0^T(c-x_0)=0,\ dc^T(c-x_0)=0$ and
are thus independent of the rolling (we exclude the case
$N_0^T(\om\times\wedge\om)=0$ of $x_0,\ x$ ruled with
correspondence of rulings).

Thus the focal surfaces are given by the intersection(s) of the
tangent plane of $x_0$ with $c$ in the first case and the
perpendicular(s) from $x_0$ to $c$ in the second case, since the
normals to the focal surfaces are multiples of $R(N_0\times dc)$
in the first case and of $R(c-x_0)$ in the second case, the focal
planes are the tangent plane to $c$ normal to the tangent plane of
$x_0$ in the first case and normal to $c-x_0$ in the second case.

Conversely, if a congruence $C(u,v,w):=(R(u,v),t(u,v))c(w)$ of
congruent curves has the property that on the normals of a focal
surface $F_2:=C(u,v,w_2(u,v))$ one can choose a surface $x$ such
that the tangent space of the surface $x$ is normal to the tangent
space of another focal surface $F_1:=C(u,v,w_1(u,v))$, then this
congruence is in general a congruence of rolling. If $c$ does not
admit a continuous group of symmetries (that is c does not have
constant curvature and torsion; equivalently $c$ is not a line,
circle or helix; note however that a line does not have curvature
and torsion defined geometrically; equivalently from a geometric
point of view a line can have arbitrarily prescribed curvature and
torsion), then we can recover $(R,t)$ (up to finitely many
choices, as the curve $c$ may have a finite group of symmetries)
from the knowledge of the congruence $C$. With $x_0:=R^{-1}(x-t)$
and because $F_1,\ F_2$ are focal surfaces we get
$0=d(c(w_2))^T(c(w_2)-x_0)=(R^{-1}dF_2)^T(c(w_2)-x_0)=(R^{-1}N)^T(c(w_1)-x_0)
=(R^{-1}dF_1)^T((R^{-1}N)\times d(c(w_1)))$. Since
$R^{-1}dF_j=R^{-1}dRc(w_j)+d(c(w_j))+R^{-1}dt,\ j=1,2$ we get
$0=(c(w_2)-x_0)^T(R^{-1}dRx_0+R^{-1}dt)$...

\subsection{Rolling distributions}\label{subsec:rolling7}
\noindent

\noindent Ribaucour proved that if we take a tangential congruence
of curves (each curve in a tangent space of $x_0$), then the
integrability condition of the distribution of facets formed by
the normal planes of the curves depends only on the linear element
of $x_0$, so if the distribution in the initial position was
integrable and we roll $x_0$ on an applicable surface $x$, then
the rolled distribution remains integrable. The statement includes
the degenerate case of a congruence of curves in the tangent
planes of a curve $x_0$, in which case the integrability depends
on the arc-length and curvature of $x_0$, so if we deform the
tangent surface of $x_0$, then the distribution remains integrable
if it was integrable in the initial position (see Bianchi
(\cite{B1},\S\ 359)).

In the particular case when the curves are circles and the
distribution is integrable such a congruence of circles is called
{\it cyclic system}. In particular cyclic systems with all circles
of the same radius $1$ exist when $x_0$ has CGC $-1$ and the
centers of circles are the corresponding points of $x_0$; in this
case the leaves of the distribution are Bianchi's complementary
transforms (also CGC $-1$ surfaces).

Darboux proved that all cyclic systems appear when we intersect
the tangent spaces of $x_0$ with a particular isotropic cone thus
highlighting a circle in each tangent space. The leaves of the
distribution in this case are the rulings of the isotropic cone
and in general will be the focal surfaces of the rolling
congruences of rulings (Darboux (\cite{D1},\S\ 936). When $x_0$ is
a quadric and the isotropic cones have umbilics as vertices among
the leaves of the rolled distribution are special isothermic
surfaces (which envelope the congruence of rolled isotropic
rulings).

Bianchi was able to prove in (\cite{B2},Vol {\bf 4},(108)) the
Permutability Theorem for the D transformation of isothermic
surfaces, a thing considered miraculous at the time since the
Ricatti character of the Darboux system was not evident. The
Ricatti character was recently discovered and involves Clifford
multiplication (see Burstall \cite{B4}). Bianchi proved in
(\cite{B2},Vol {\bf 4},(114)) that a D transformation for special
isothermic surfaces is a composition of two B transformations of
quadrics.

Consider an integrable $3$-dimensional distribution of facets
$(p,P)=(p,P)(u,v,w)$ in $\mathbb{C}^3$ along the surface
$x=x(u,v)$ (thus distributed as a $1$-dimensional distribution of
facets along each point of $x$). It is thus natural to inquire
when this distribution remains integrable if we roll it as $x$
rolls on any applicable surface $x_0$ (see Bianchi (\cite{B2},Vol
{\bf 4},(173))). Let $m=m(u,v,w)\bot P(u,v,w)$ be a normal field
along the distribution. We shall use the notation: $\ti
df=f_udu+f_vdv+f_wdw=df+f_wdw$ for $f=f(u,v,w)$. If the
distribution is integrable with leaves $p(u,v,w(u,v))$, then
$0=m^T\ti dp=m^Tp_wdw+m^Tdp$ holds. In particular if
$m^TN=N^T(p-c(w)x)=0$, then the integrability condition depends
only on the linear element of $x$. Imposing the integrability
condition $m^Tp_w\ti d\wedge$ (suppose $m^Tp_w\neq 0$) and using
the equation itself we get: $(m_w^Tdp-dm^Tp_w)\wedge
m^Tdp+m^Tp_wdm^T\wedge dp=0$, or:
\begin{eqnarray}\label{eq:inteco}
(dp\times p_w)^T\wedge(m\times dm)+\frac{1}{2}(m_w\times
m)^T(dp\times\wedge dp)=0.
\end{eqnarray}
If $(x_0,dx_0)=(R,t)(x,dx)$ is the rolling of $x$ on $x_0$, then
the rolled distribution will be $(R,t)(V+x,P)=(RV+x_0,RP),\
V:=p-x$, so (\ref{eq:inteco}) must be satisfied with $RV,\ Rm$
replacing $V,\ m$: $((d(V+x)+\om\times V)\times
V_w)^T\wedge(m\times (dm+\om\times m))+\frac{1}{2}(m_w\times
m)^T((d(V+x)+\om\times V)\times\wedge(d(V+x)+\om\times V))=\ =0$,
or, using (\ref{eq:inteco}): $((m\times((d(V+x)\times V_w)\times
m)-V\times(V_w\times(m\times dm)+d(V+x)\times(m_w\times
m)))^T\wedge\om+(V^T(m_w\times m)V+V_w\times((m\times V)\times
m))^T\frac{\om\times\wedge\om}{2}=0$. Using the Gau\ss\ theorem
($dN\times\wedge dN=KN\sqrt{g}du\wedge dv$) for both $x_0,\ x$ we
have
$\frac{\om\times\wedge\om}{2}=\frac{1}{2}N^T((R^{-1}dN_0-dN)\times\wedge(R^{-1}dN_0-dN))N
=-N^T(dN\times\wedge(R^{-1}dN_0-dN))N=dN^T\wedge\om N$ and we thus
need: $(m\times((d(V+x)\times V_w)\times
m)-V\times(V_w\times(m\times dm)+d(V+x)\times(m_w\times
m))+(V^T(m_w\times m)V+V_w\times((m\times V)\times
m))^TNdN)^T\wedge\om=0$ for any flat connection form $\om$. This
is equivalent to $(m\times((d(V+x)\times V_w)\times m)-
V\times(V_w\times(m\times dm)+d(V+x)\times(m_w\times m))
+(V^T(m_w\times m)V+V_w\times((m\times V)\times m))^TNdN)^Tdx=0$,
or using $dx\times dx=Nda=N\sqrt{g}du\wedge dv$:
\begin{eqnarray}\label{eq:disa}
(m^Td(V+x)(m\times V)_w-((m^Td(V+x))_w-d(m^TV_w))m\times V-m^TV_wd(m\times\ V))^Tdx-\nonumber\\
m^TV_wm^TNda+((m\times V)\times(m\times V)_w-m\times
Vm^TV_w)^TNN^Td^2x=0.
\end{eqnarray}
But (\ref{eq:disa}) is satisfied if $RV,\ Rm$ replace $V,\ m$, so
it depends only on the linear element of $x$, as expected.
Therefore the distributions with the required properties must
satisfy (\ref{eq:inteco}) and (\ref{eq:disa}) in which we ignore
the terms involving linearly the second fundamental form of $x$
(these equation can also be obtained from (\ref{eq:inteco}) (which
becomes an affine equation in the second fundamental form of $x$
by an application of the Gau\ss\ theorem) by setting all the
coefficients of the second fundamental form of $x$ to $0$). We
have $dV=V^I+V^{II},\ V^I:=d(V^TN)N+d(V^T\na x\lrcorner dx^T)\na
x\lrcorner dx,\ V^{II}:=N^Td^2x\llcorner((\na x\times N)\times
V),\ N\times dx=-\na x\lrcorner da,\ m\times V=(m^T\na xV^T\na
xN+\begin{vmatrix}m^TN&m^T\na x\\V^TN&V^T\na x\end{vmatrix}\na
x)\lrcorner\lrcorner da$. Therefore we need:
\begin{eqnarray}\label{eq:disb}
((V^I+dx)\times V_w)^T\wedge(m\times m^I)+\frac{1}{2}(m_w\times m)^T
((V^I+dx)\times \wedge(V^I+dx))-\nonumber\\
KN^T((m\times V)\times(m\times V)_w-m\times Vm^TV_w)da=0,\nonumber\\
(m^T(V^I+dx)(m\times V)_w-((m^T(V^I+dx))_w-d(m^TV_w))m\times V-\nonumber\\
m^TV_w(m\times V)^I)^Tdx-m^TV_wm^TNda=0.
\end{eqnarray}
In particular if the distribution is integrable for a particular
deformation of $x$ (which can be taken as the definition of the
distribution), then we need not consider the first equation. If we
take the intersections of the tangent planes of $x$ with an
isotropic developable (rulings are isotropic vectors, or tangent
to $C(\infty)$), then these curves are like in Ribacour's result
and the leaves of the distribution in this initial position are
the rulings (Darboux (\cite{D1},\S\ 762)).

Suppose that the leaves of the initial distribution $(p,P)$ are
curves: $0=\ti d(V+x)\times\wedge \ti
d(V+x)=V_w^T(d(V+x)\times\wedge d(V+x))m^TV_wm$. Since $m^TV_w\neq
0$ this is equivalent to
$0=V_w^T(KN^TVVda+(V^I+dx)\times\wedge(V^I+2V^{II}+dx))$...

\section{Deformations in $\mathbb{C}^3$ of quadrics in
$\mathbb{C}^3$}\label{sec:deformations}\setcounter{equation}{0}

\subsection{The Ricatti equation and the applicability correspondence provided by the Ivory affinity}
\label{subsec:deformations1} \noindent

\noindent Let $x^0\subset\mathbb{C}^3$ be non-rigidly applicable
to the surface $x_0^0$ in the quadric $x_0$:
$(x^0,dx^0)=(R_0,t_0)(x_0^0,dx_0^0)$. If the facets of the
distribution $\mathcal{D}$ become tangent spaces to surfaces
(leaves) $x^1:=(R_0,x^0)V_0^1$ when rigidly moved with $x_0^0$
when it rolls on $x^0$, then $0=(R_0m_0^1)^Tdx^1$ holds.

This will constitute a completely integrable Ricatti equation on
the parameter of the conics $T_{x_0^0}x_0\cap x_z$ (chosen to be
$v_1$), so the arbitrariness of the transformation $B_z$ taking
$x^0$ into $x^1$ will be given by the choice of a point on an
initial conic (for simplicity we ignore this $1$-dimensionality in
the notation $B_z$); 4 transforms $B_z(x^0)$ will cut these conics
with constant cross-ratio. Therefore to find all transforms
$B_z(x^0)$ with fixed $z$ we need only one solution of this
Ricatti equation, two of a linear differential equation and
further algebraic computations.

We have:
$R_0^{-1}dx^1=d(x_0^0+V_0^1)+R_0^{-1}dR_0V_0^1=dx_z^1+\om_0\times
V_0^1,\ \om_0:=N_0^0\times R_0^{-1}dR_0N_0^0$. But
$(\om_0)^{\bot}=0$ and $dx_z^1=x_{zv_1}^1dv_1+x_{zu_1}^1du_1$, so
$0=(R_0m_0^1)^Tdx^1$ becomes:
\begin{eqnarray}\label{eq:dify1}
-(V_0^1)^T(\om_0\times N_0^0)(m_0^1)^TN_0^0+(m_0^1)^Tx_{zv_1}^1dv_1=0.
\end{eqnarray}
If similarly to Bianchi's original approach (122,\S\ 32) we
postpone the proof of the complete integrability of
(\ref{eq:dify1}) until later and first try to prove the ACPIA:
$|dx^1|^2=|dx_0^1|^2$, then the necessary algebraic conditions
(\ref{eq:simpas}) appear.

Using the second relation of (\ref{eq:simpas}), (\ref{eq:dify1})
becomes:
\begin{eqnarray}\label{eq:dify3}
-(V_0^1)^T(\om_0\times N_0^0)+2(x_{zv_1}^1)^TN_0^0dv_1=0.
\end{eqnarray}
Multiplying this by $\mathcal{B}_1(x_{zu_1}^1)^TN_0^0$ and using
$-\mathcal{B}_1(x_{zu_1}^1)^TN_0^0V_0^1=\mathcal{B}_1(V_0^1\times
x_{zu_1}^1)\times N_0^0=-m_0^1\times N_0^0$ we finally get:
\begin{eqnarray}\label{eq:dify2}
(m_0^1)^T\om_0+2zdv_1=0.
\end{eqnarray}
Thus the Ricatti character of the equation appears clearly. We
have $dm_0^1=m_{0v_1}^1dv_1+\mathcal{B}_1dx_0^0\times x_{zu_1}^1$,
so $d(m_0^1)^T\wedge\om_0=dv_1\wedge(m_{0v_1}^1)^T\om_0$. Imposing
the condition $d\wedge$ on (\ref{eq:dify2}) and using the equation
itself we need:
$-(m_0^1)^T\om_0\wedge(m_{0v_1}^1)^T\om_0+2z(m_0^1)^Td\wedge\om_0=0$,
or using (\ref{eq:om}): $(N_0^0)^T(2zm_0^1+m_0^1\times
m_{0v_1}^1)(N_0^0)^T(\om_0\times\wedge\om_0)=0$. Therefore the
complete integrability is equivalent to (\ref{eq:intalg}). We get
the transformation $B_z'$ for the other ruling family by replacing
$(m_0^1,v_1)$ with $(m'^1_0,u_1)$. Since the distributions
$\mathcal{D},\ \mathcal{D'}$ reflect in $Tx_0^0$, the rolled
distributions $(R_0,t_0)\mathcal{D},\ (R_0,t_0)\mathcal{D'}$
reflect in $Tx^0$, so $B'_z(x^0)$ is just $B_z(x^0)$ when $x_0^0$
rolls on the other face of $x^0$. This geometric observation can
be confirmed analytically if one keeps account of the change of
$\om_0$ to $\om'_0:=-\om_0-2N_0^0\times dN_0^0$ when we roll
$x_0^0$ on the other face of $x^0$ and of (\ref{eq:du1}). Thus
\begin{eqnarray}\label{eq:BB'}
(m_0^1)^T\om_0+2zdv_1=0\ \Leftrightarrow\ (m'^1_0)^T\om'_0+2zdu_1=0\ \Leftrightarrow\ B_z\
\mathrm{transformation},\nonumber\\
(m_0^1)^T\om'_0+2zdv_1=0\ \Leftrightarrow\ (m'^1_0)^T\om_0+2zdu_1=0\ \Leftrightarrow\ B'_z\
\mathrm{transformation}.
\end{eqnarray}

Note that to obtain (\ref{eq:intalg}) as the complete
integrability of (\ref{eq:dify2}) we just used
$m_0^1=f(z,v_1)+g(z,v_1)\times x_0^0$. The condition that the
vectors $f,g$ satisfy (\ref{eq:intalg}) imposes certain
restrictions on the geometry of $x_0^0$ (for example if $f,g$ are
rational functions in $v_1$, then we have algebraic restriction on
the geometry of $x_0^0$, which may be over-determined and prevent
its existence); also equation (\ref{eq:dify2}) loses its geometric
meaning.

Using $R_0^{-1}dx^1=dx_z^1+\om_0\times V_0^1$ and (\ref{eq:dify3})
we have $R_0^{-1}dx^1=dx_z^1-2(x_{zv_1}^1)^TN_0^0N_0^0dv_1$ and
the ACPIA becomes:
$|dx_z^1-2(x_{zv_1}^1)^TN_0^0N_0^0dv_1|^2=|dx_0^1|^2$, which boils
down the first relation of (\ref{eq:simpas}). Note that if we
conversely impose the ACPIA, then we get (\ref{eq:dify2}) (or its
correspondent for the other ruling family).

We can try to prescribe the function $v_1(u_0,v_0)$ and solve
(\ref{eq:dify2}) for flat connection forms $\om_0$. If
$\frac{\om_0}{\mathcal{B}_0}=:(\om_{0u_0}^{u_0}x_{0u_0}^0-\om_{0u_0}^{v_0}x_{0v_0}^0)du_0+
(\om_{0v_0}^{u_0}x_{0u_0}^0-\om_{0v_0}^{v_0}x_{0v_0}^0)dv_0$, then
$0=\om_0\times\wedge dx_0^0$ is equivalent to
$\om_{0u_0}^{u_0}=\om_{0v_0}^{v_0}$ and the first line of
(\ref{eq:BB'}) becomes:
\begin{eqnarray}\label{eq:syst}
\begin{bmatrix}v_{1u_0}\\v_{1v_0}\end{bmatrix}=
\begin{bmatrix}\om_{0u_0}^{v_0}&\om_{0u_0}^{u_0}\\
\om_{0u_0}^{u_0}&\om_{0v_0}^{u_0}\end{bmatrix}
\begin{bmatrix}\frac{\Del^+}{2z}\\\frac{\Del^-}{2z}\end{bmatrix};
\begin{bmatrix}u_{1u_0}\\u_{1v_0}\end{bmatrix}=
\begin{bmatrix}-\om_{0u_0}^{v_0}&-(\om_{0u_0}^{u_0}+2\mathcal{N}_0)\\
-(\om_{0u_0}^{u_0}+2\mathcal{N}_0)&-\om_{0v_0}^{u_0}\end{bmatrix}\begin{bmatrix}\frac{\Del'^+}{2z}\\
\frac{\Del'^-}{2z}\end{bmatrix}.
\end{eqnarray}
With $\rho_0^2:=\frac{-1}{K(x_0^0)}=\mathcal{A}^2|\hat N_0^0|^4,\
g_0^0:=\frac{\det(|dx_0^0|^2)}{(du_0\wedge dv_0)^2}=
\frac{4\mathcal{A}^2|\hat N_0^0|^2}{\mathcal{B}_0^2}$ we have
$\mathcal{N}_0=\frac{1}{\mathcal{B}_0\rho_0}$ and the Christoffel
symbols of $x_0^0$ are $\Ga^1_{22}=\Ga^2_{11}=0,\
-\Ga^1_{11}=\Ga^2_{22}=\pa_{u_0}\log\mathcal{B}_0,\
\Ga^2_{12}=\pa_{u_0}\log\sqrt{g_0^0}-\Ga^1_{11}=\pa_{u_0}\log\sqrt{\rho_0},\
\Ga^1_{12}=\pa_{v_0}\log\sqrt{g_0^0}-\Ga^2_{22}=\pa_{v_0}\log\sqrt{\rho_0}$.
Keeping account of (\ref{eq:omjk}), the condition
$d\wedge\om_0+\frac{1}{2}\om_0\times\wedge\om_0=0$ becomes:
\begin{eqnarray}\label{eq:flatcon}
(\om_{0u_0}^{u_0})^2+2\mathcal{N}_0\om_{0u_0}^{u_0}=\om_{0v_0}^{u_0}\om_{0u_0}^{v_0},\nonumber\\
\pa_{v_0}\om_{0u_0}^{v_0}=\om_{0u_0}^{u_0}\pa_{u_0}\log\frac{\om_{0u_0}^{u_0}}{\mathcal{N}_0},\
\pa_{u_0}\om_{0v_0}^{u_0}=\om_{0u_0}^{u_0}\pa_{v_0}\log\frac{\om_{0u_0}^{u_0}}{\mathcal{N}_0}
\end{eqnarray}
(the difference of the GCM equations for $x^0,\ x_0^0$). In
particular if $v_{1u_0}=0$, then
$\om_{0u_0}^{u_0}=\om_{0u_0}^{v_0}=R_{0u_0}=t_{0u_0}=0$ and
$x^0=(R_0(v_0),t_0(v_0))x_0^0$ is a ruled surface, in which case
$\om_{0v_0}^{u_0}=\varphi_0(v_0),\\
\varphi_0(v_0):=\frac{2zv_{1v_0}}{\Del^-}$. The nonzero diagonal
term of the second fundamental form of $x^0$ is
$\mathcal{B}_0\sqrt{g_0^0}\varphi_0(v_0)$. Using (\ref{eq:dre}) we
see that $x^1$ is also ruled with rulings corresponding to
$x_{0u_1}^1$. Therefore as we roll a quadric on a ruled applicable
surface, the ruling families of a confocal quadric will generate
congruences decomposable in ruled surfaces applicable to the
initial quadric. Note that since $\Del^+$ is linear in $v_1$ for
IQWC (\ref{eq:iqwc2}), if $x^0$ is ruled with rulings
corresponding to $x_{0v_0}^0$, then its $B_z$ transforms can be
found by quadratures only (therefore its $B'_z$ transforms can
also be found by quadratures only).

Using first equations of (\ref{eq:syst}) and  (\ref{eq:flatcon})
we get:
\begin{eqnarray}\label{eq:flatco}
\om_{0u_0}^{u_0}=\frac{2z^2v_{1u_0}v_{1v_0}}{\mathcal{N}_0\Del^+\Del^-+zv_{1u_0}\Del^+
+zv_{1v_0}\Del^-},\
\om_{0u_0}^{v_0}=\frac{2zv_{1u_0}(\mathcal{N}_0\Del^-+zv_{1u_0})}{\mathcal{N}_0\Del^+\Del^-
+zv_{1u_0}\Del^++zv_{1v_0}\Del^-},\nonumber\\
\om_{0v_0}^{u_0}=\frac{2zv_{1v_0}(\mathcal{N}_0\Del^++zv_{1v_0})}{\mathcal{N}_0\Del^+\Del^-
+zv_{1u_0}\Del^++zv_{1v_0}\Del^-}.
\end{eqnarray}
From the compatibility condition
$(v_{1u_0})_{v_0}=(v_{1v_0})_{u_0}$ applied to the first equation
of (\ref{eq:syst}), using
$((u_1,\Del'^-,\Del'^+)\leftrightarrow(v_1,\Del^-,\Del^+))\circ$
the first relation of (\ref{eq:alg}) and the first relation of
(\ref{eq:flatcon}) (which $\om_0$ as defined by (\ref{eq:flatco})
automatically satisfies) we get the fact that the last two
relations of (\ref{eq:flatcon}) are linearly dependent with
coefficients $-\Del^+,\ \Del^-$; thus either of them must be
satisfied. The second equation of (\ref{eq:flatcon}) becomes:
$(\mathcal{N}_0\Del^-+zv_{1u_0})
(\log\frac{v_{1u_0}(\mathcal{N}_0\Del^-+zv_{1u_0})}{\mathcal{N}_0\Del^+\Del^-+zv_{1u_0}\Del^+
+zv_{1v_0}\Del^-})_{v_0}=zv_{1v_0}(\log\frac{v_{1u_0}v_{1v_0}}{\mathcal{N}_0(\mathcal{N}_0
\Del^+\Del^-+zv_{1u_0}\Del^++zv_{1v_0}\Del^-)})_{u_0}$, or using
again
$((u_1,\Del'^-,\Del'^+)\leftrightarrow(v_1,\Del^-,\Del^+))\circ$
the first relation of (\ref{eq:alg}):
\begin{eqnarray}\label{eq:new1}
z(\mathcal{N}_0\Del^++zv_{1v_0})(\log\frac{v_{1u_0}}{v_{1v_0}})_{u_0}+z(\mathcal{N}_0\Del^-+zv_{1u_0})
(\log\frac{v_{1v_0}}{v_{1u_0}})_{v_0}-\mathcal{N}_0^2\Del^+\Del^-\frac{v_{1u_0v_0}}{v_{1u_0}v_{1v_0}}
-\nonumber\\ z^2(v_{1v_0}(\log \mathcal{N}_0)_{u_0}+v_{1u_0}(\log
\mathcal{N}_0)_{v_0})
-\mathcal{N}_0\Del^+\Del^-(\frac{z}{\Del^-}\pa_{u_0}\log(\mathcal{N}_0\Del^+)-\nonumber\\
\mathcal{N}_0\pa_{v_1}\log\Del^++\frac{z}{\Del^-}(\log
\mathcal{N}_0)_{u_0} +\frac{z}{\Del^+}(\log
\mathcal{N}_0)_{v_0})=0.
\end{eqnarray}
Therefore $v_1$ must satisfy a second order PDE linear in the
second order terms and having the symmetry
$(u_0,\Del^+)\leftrightarrow(v_0,\Del^-)$ (because of this the
third equation of (\ref{eq:flatcon}) and (\ref{eq:new1}) are
equivalent). Since (\ref{eq:flatcon}) is valid with
$(\om_{0u_0}^{v_0},\om_{0u_0}^{u_0},\om_{0v_0}^{u_0})
\leftrightarrow(-\om_{0u_0}^{v_0},-(\om_{0u_0}^{u_0}+2M_0),-\om_{0v_0}^{u_0})$,
(\ref{eq:new1}) is valid for
$(u_1,\Del'^-,\Del'^+)\leftrightarrow(v_1,\Del^-,\Del^+)$.

\subsection{Inversion of the B\"{a}cklund transformation and the Weingarten congruence}
\label{subsec:deformations2}
\noindent

\noindent Let $x_0^0,\ x_0^1$ be in the TC $(V_0^1)^TN_0^0=0$ and
$v_1$ be given by (\ref{eq:dify2}), where
$(x^0,dx^0)=(R_0,t_0)(x_0^0,dx_0^0)$.

If we roll $x^0$ on $x_0^0$, then the tangent space of
$x^1:=(R_0,x^0)V_0^1$ will be applied to the space at $x_z^1$
generated by $x_{zu_1}^1,\ V_0^1$. Further applying
$(R_0^1,t_0^1)$, this space will be applied to $T_{x_0^1}x_0$ with
$(R_0^1,t_0^1)(x_z^1,x_{zu_1}^1,V_0^1)=(x_0^1,x_{0u_1}^1,-V_1^0)$.
If $(x^1,dx^1)=(R_1,t_1)(x_0^1,dx_0^1)$ is the rolling of $x_0^1$
on $x^1$ with $\det(R_1^{-1}R_0)=1$, then we have just proved
$(R_0^1,t_0^1)(R_0,t_0)^{-1}(x^1,dx^1)=(R_1,t_1)^{-1}(x^1,dx^1)$,
so
\begin{eqnarray}\label{eq:dre}
(R_0^1,t_0^1)=(R_1,t_1)^{-1}(R_0,t_0).
\end{eqnarray}
If we make the ansatz $x^0=x_0^0$, so $(R_0,t_0)=(I_3,0),\
\om_0^0=0,\ v_1=\mathrm{const}$ or
$(R_0,t_0)=(I_3-2N_0^0(N_0^0)^T,2N_0^0(N_0^0)^Tx_0^0),\
\om_0^0=-2N_0^0\times dN_0^0,\ v_1=\mathrm{const}$, then we get in
a natural geometric way all interesting algebraic identities of
the static part (including the explanation for the necessity of
the existence of the RMPIA).

Moreover, $x^0=(R_0,x^1)(-V_0^1)=(R_0,x^1)(R_0^1)^{-1}V_1^0
=(R_1,x^1)V_1^0$ and thus $x^0$ reveals itself as a $B_z$
transform of $x^1$. Therefore we can find all the B transforms
$B_z(x^1)$ for fixed $z$ only by two quadratures and further
algebraic computations. We shall see later that we don't need the
quadratures once we know the B transforms $B_z(x^0)$ for all $z$
as one can take a derivative in the spectral parameter $z$ in a
L'Hospital situation of the BPT (although Sym's formula was not
known to the classical geometers, tricks similar to it were). Also
equations (\ref{eq:syst})-(\ref{eq:new1}) admit the symmetry
$(u_0,v_0,\Del^+,\Del'^-,\om_0,\mathcal{N}_0)\leftrightarrow
(u_1,v_1,\Del'^-,\Del^+,\om_1,\mathcal{N}_1)$. We shall now
provide an analytic confirmation of the inversion of the B
transformation, equivalently of the fact that
$((u_0,v_0,\Del^+,\Del'^-,\mathcal{N}_0)\leftrightarrow
(u_1,v_1,\Del'^-,\Del^+,\mathcal{N}_1))\circ$(\ref{eq:new1}) is
still valid. This analytic confirmation will play an important
r\^{o}le in the discussion of the $B_{\infty}$ transformation,
since the easy geometric justification will not be available.

We have$\begin{bmatrix}u_{0u_1}&u_{0v_1}\\
v_{0u_1}&v_{0v_1}\end{bmatrix}=\frac{1}{J}\begin{bmatrix}v_{1v_0}&-u_{1v_0}\\
-v_{1u_0}&u_{1u_0}\end{bmatrix}=\frac{1}{J}\begin{bmatrix}v_{1v_0}&\frac{\Del'^-}{\Del^-}
(\frac{\mathcal{N}_0\Del^+}{z}+v_{1v_0})\\-v_{1u_0}&-\frac{\Del'^-}{\Del^-}
(\frac{\mathcal{N}_0\Del^-}{z}+v_{1u_0})\end{bmatrix}$, where
$J:=\frac{du_1\wedge dv_1}{du_0\wedge
dv_0}=\frac{\mathcal{N}_0\Del'^+}{z}v_{1u_0}
-\frac{\mathcal{N}_0\Del'^-}{z}v_{1v_0}$ (assume $J\neq 0$).

Now $((u_0,v_0,\Del^+,\Del'^-,\mathcal{N}_0)\leftrightarrow
(u_1,v_1,\Del'^-,\Del^+,\mathcal{N}_1))\circ$(\ref{eq:new1})
becomes

$z(\mathcal{N}_1\Del'^-+zv_{0v_1})((\log\frac{v_{0u_1}}{v_{0v_1}})_{u_0}u_{0u_1}
+(\log\frac{v_{0u_1}}{v_{0v_1}})_{v_0}v_{0u_1})
+z(\mathcal{N}_1\Del^-+zv_{0u_1})((\log\frac{v_{0v_1}}{v_{0u_1}})_{u_0}u_{0v_1}+\\
(\log\frac{v_{0v_1}}{v_{0u_1}})_{v_0}v_{0v_1})
-\mathcal{N}_1^2\Del'^-\Del^-((\log v_{0u_1})_{v_0}
+\frac{u_{0v_1}}{v_{0v_1}}(\log v_{0u_1})_{u_0})
-z^2(v_{0v_1}(\log\mathcal{N}_1)_{u_1}+v_{0u_1}(\log\mathcal{N}_1)_{v_1})
-\mathcal{N}_1\Del'^-\Del^-(\frac{z}{\Del^-}\pa_{u_1}\log(\mathcal{N}_1\Del'^-)
-\mathcal{N}_1\pa_{v_0}\log\Del'^-+\frac{z}{\Del^-}(\log\mathcal{N}_1)_{u_1}+\frac{z}{\Del'^-}
\log(\mathcal{N}_1)_{v_1})=0$. Using (\ref{eq:del+-+-}), the
coefficient of $(\log\frac{v_{0u_1}}{v_{0v_1}})_{u_0}$ becomes
$0$; the remaining part becomes

$0=-\frac{zv_{1u_0}\mathcal{N}_0\mathcal{N}_1(\Del^-)^2}{J(\mathcal{N}_0\Del^-+zv_{1u_0})}
(\frac{\mathcal{N}_0\Del'^-}{zv_{1u_0}}(\pa_{v_0}\log(\mathcal{N}_0\Del'^-)
+u_{1v_0}\pa_{u_1}\log\Del'^--(\log
v_{1u_0})_{v_0})+\frac{\Del'^-}{\Del^-}
(\pa_{v_0}\log\frac{\Del'^-}{\Del^-}+\\u_{1v_0}\pa_{u_1}\log\Del'^--v_{1v_0}\pa_{v_1}\log\Del^-))
+\frac{\mathcal{N}_1^2\Del'^-\Del^-v_{1u_0}}{J}(\frac{\mathcal{N}_0\Del'^+}{z}
((\log\mathcal{N}_0)_{v_0}+u_{1v_0}\pa_{u_1}\log\Del'^+)-\\
\frac{\mathcal{N}_0\Del'^-v_{1v_0}}{zv_{1u_0}}(\pa_{v_0}\log(\mathcal{N}_0\Del'^-)+
u_{1v_0}\pa_{u_1}\log\Del'^-+(\log\frac{v_{1v_0}}{v_{1u_0}})_{v_0})
-\frac{\mathcal{N}_0\Del^++zv_{1v_0}}{\mathcal{N}_0\Del^-+zv_{1u_0}}
(\frac{\mathcal{N}_0\Del'^+}{z}(\pa_{u_0}\log(\mathcal{N}_0\Del'^+)+\\
u_{1u_0}\pa_{u_1}\log\Del'^+)-\frac{\mathcal{N}_0\Del'^-v_{1v_0}}{zv_{1u_0}}
((\log\mathcal{N}_0)_{u_0}+u_{1u_0}\pa_{u_1}\log\Del'^--(\log\frac{v_{1u_0}}{v_{1v_0}})_{u_0})))
-\frac{z^2}{J}(u_{1u_0}(\log\mathcal{N}_1)_{u_1}-\\
v_{1u_0}(\log\mathcal{N}_1)_{v_1})-\mathcal{N}_1\Del'^-\Del^-(\frac{z}{\Del^-}\pa_{u_1}
\log(\mathcal{N}_1\Del'^-)-\mathcal{N}_1\pa_{v_0}\log\Del'^-
+\frac{z}{\Del^-}(\log\mathcal{N}_1)_{u_1}+\frac{z}{\Del'^-}\log(\mathcal{N}_1)_{v_1})$
Comparing the coefficients of
$(\log\frac{v_{1u_0}}{v_{1v_0}})_{u_0},\
(\log\frac{v_{1v_0}}{v_{1u_0}})_{v_0},\
\frac{v_{1u_0v_0}}{v_{1u_0}v_{1u_0}}$ we see that so far this is
just (\ref{eq:new1}) multiplied by
$-\frac{\mathcal{N}_0\mathcal{N}_1^2\Del^-
(\Del'^-)^2v_{1v_0}}{z^2J(\mathcal{N}_0\Del^-+zv_{1u_0})}$.
Replacing $u_{1u_0},u_{1v_0},J$ with their values (note that the
terms in the second parenthesis containing $u_{1u_0},u_{1v_0}$
cancel), the remaining part is the same multiple of
(\ref{eq:new1}) (this boils down to a quadratic polynomial in
$v_{1u_0},v_{1v_0}$ being identically $0$; no further expansion of
derivatives is needed; one uses instead the algebraic identities
(\ref{eq:del+-+-}), (\ref{eq:alg}) and its versions under
symmetries).

Another way to prove the inversion of the B transformation is to
show that $\om_1$ as defined by
$((u_0,v_0,\Del^+,\Del'^-,\mathcal{N}_0)\leftrightarrow
(u_1,v_1,\Del'^-,\Del^+,\mathcal{N}_1))\circ$(\ref{eq:flatco})
satisfies the last two equations of
$\\((u_0,v_0,\om_0,\mathcal{N}_0)\leftrightarrow
(u_1,v_1,\om_1,\mathcal{N}_1))\circ$(\ref{eq:flatcon}) (the first
equation is automatically satisfied); by symmetry it is enough to
prove only the second one. The differential of $x^1$ depends on
the shape (but not its derivatives) of $x^0$; a-priori the shape
of $x^1$ depends also on the first derivatives of the shape of
$x^0$, but it actually depends only on the shape of $x^0$. Using
the GCM equations and the first derivatives of the Gau\ss\
equation of $x^0$, the required equation is algebraic in
$u_0,v_0,u_1,v_1,\om_{0u_0}^{v_0},\om_{0u_0}^{u_0},
\pa_{u_0}\om_{0u_0}^{u_0},\pa_{v_0}\om_{0u_0}^{u_0}$; accounting
the coefficients of the four terms involving the shape of $x^0$
and its derivatives they will be identically $0$.

Note that we have provided another proof for the ACPIA by
effectively finding the rolling. This proof is the one susceptible
for generalization and the essence of the theory of deformations
of quadrics revealed itself: take a surface $x_0=x_0(u,v)$ in
$\mathbb{C}^3$ and a $3$-dimensional integrable distribution of
facets $(p,P)=(p,P)(u,v,w)$, whose integrability condition depends
only on the linear element of $x_0$. If $x^0\subset\mathbb{C}^3$
is applicable to $x_0^0:=x_0(u_0,v_0)$ with rolling
$(x^0,dx^0)=(R_0,t_0)(x_0^0,dx_0^0)$, then the rolled distribution
will still be integrable with leaves

$x^1:=(R_0(u_0,v_0),x^0(u_0,v_0))(p(u_0,v_0,w(u_0,v_0))-x^0(u_0,v_0))$.
If for each triple $f:=(u,v,w)$ we have a triple of functions $\ti
f:=(\ti u,\ti v,\ti w)$ of $(u,v,w)$, a rigid motion $(R_{f\ti
f},t_{f\ti f})$ taking $(p,P)(u,v,w)$ to $(x_0,dx_0)(\ti u,\ti v)$
and $(R_{f\ti f}R_0^{-1}dx_1)(u_1,v_1)=dx_0(u_1,v_1)$, where

$(u_1,v_1):=(\ti u(u_0,v_0,w(u_0,v_0)),\ \ti
v(u_0,v_0,w(u_0,v_0)))$, so by a change of coordinates (when
$(u_1,v_1)$ are functionally independent) we can consider the
$(u_1,v_1)$ coordinates on $x^1,\ x_0^1:=x_0(u_1,v_1)$, then $x^1$
is applicable to $x_0^1$ with explicit rolling. If further
$(R_{f\ti f},t_{f\ti f})$ takes $(x_0,dx_0)(u,v)$ to $(p,P)(\ti
u,\ti v,\ti w)$, then the relationship between $x^0,\ x^1$ is
symmetric. To find such $(R_{f\ti f},t_{f\ti f})$ we need a
correspondence between the given distribution of facets $(p,P)$
and the distribution of facets $(x_0,dx_0)$; therefore we need a
submersion between the two distributions. In the case of quadrics
this submersion factors through the leaves (a ruling family on the
confocal quadric) and becomes the Ivory affinity. It is thus
natural for this submersion to factor through the leaves of the
initial distribution (when these are curves) to a diffeomorphism
between the surface described by these curves and $x_0$. In this
case $(R_{f\ti f},t_{f\ti f})$ are not uniquely determined and it
is natural that the extra condition should require that it takes
the tangent vector to the curve (leaf) to a certain vector of the
tangent space of $x_0$ at the point corresponding under the
diffeomorphism.

However, it may be the case that only the B transformation for
quadrics satisfies such requirements. Bianchi (for example in
(\cite{B2}, Vol {\bf 4},(173)) tried generalizations (natural from
a geometric point of view) of various aspects of the B
transformation for quadrics, but most of such configurations arise
only for quadrics.

According to R. Calapso in \S\ 5.17 of the introduction to Bianchi
(\cite{B2}, Vol {\bf 4}, part 1), it may be an interesting problem
to find the surfaces for which a theory of deformations similar to
that for quadrics can be developed.

We are interested in the case when $x^1$ degenerates to a curve,
that is $0=dx^1\times\wedge dx^1$, or using (\ref{eq:dre}):
$du_1\wedge dv_1=0$, or replacing $du_1$ from (\ref{eq:du1}):
$v_{1u_0}\Del^+-v_{1v_0}\Del^-=0$, which must be adjoined to
(\ref{eq:new1}). In general $v_{1u_0}v_{1v_0}=0$ is not possible,
so we need $v_{1u_0},v_{1v_0}\neq 0$. Using (\ref{eq:new1}) we
thus need:
\begin{eqnarray}\label{eq:solit}
\mathcal{N}_0^2(\frac{\Del^-}{v_{1u_0}})^2v_{1u_0v_0}+(\mathcal{N}_0\Del^+\Del^-+zv_{1u_0}\Del^+)
(\frac{z}{\Del^-}\pa_{u_0}\log(\mathcal{N}_0\Del^+)+\frac{z}{\Del^+}\pa_{v_0}\log(\mathcal{N}_0\Del^-))
+\nonumber\\\mathcal{N}_0\Del^+\Del^-(\frac{z}{\Del^-}\pa_{u_0}\log(\mathcal{N}_0\Del^+)
-\mathcal{N}_0\pa_{v_1}\log\Del^+)=0,\
v_{1u_0}\Del^+-v_{1v_0}\Del^-=0.
\end{eqnarray}
Differentiating the second relation of (\ref{eq:solit}) wrt $u_0$
and $v_0$ and using the first one we get $v_{1u_0u_0},\\
v_{1v_0v_0},\ v_{1u_0v_0}$ in terms of $v_1,\ v_{1u_0}$ and so on;
therefore all derivatives of $v_1$ can be written in function of
$v_1,\ v_{1u_0}$, so in general (formally) such a $v_1$ depends on
two constants besides the parameter $z$. An application of
Cauchy-Kovalewskaia assures in general the existence of such
solutions. Because of the first relation of (\ref{eq:alg}) all
derivatives of $v_1$ admit the symmetry
$(u_0,\Del^+)\leftrightarrow(v_0,\Del^-)$. Once the solution $v_1$
is found, (\ref{eq:flatco}) gives the shape of the $1$-solitons.
Note that the family of $1$-solitons may depend up to rigid
motions on less than three constants; for example Dini's CGC $-1$
helicoids depend up to rigid motions only on a constant (instead
of three; in fact the axis of the tractrix as a limiting case of a
Dini helicoid is a good picture of it as a CGC $-1$ surface). Note
however that Serret's ruled surfaces generated by taking isotropic
directions in tangent planes along asymptotes of Dini's helicoids
still provide another $1$-dimensional family of ruled imaginary
$1$-solitons (again they cannot be seen at the level of the
sine-Gordon equation); this must be also the case for general
quadrics (by Chieffi's), but at that level one cannot see the
analytic behavior. If we take according to Rogers-Schieff
(\cite{RS1},(1.85)) the Dini helicoids $x=x(u,v):=
\cosh(\frac{u-v\cos(\zeta)}{\sin(\zeta)})^{-1}\sin(\zeta)(\cos(v)e_1+\sin(v)e_2)
+(u-\sin(\zeta)\tanh(\frac{u-v\cos(\zeta)}{\sin(\zeta)}))e_3$,
then it becomes clear that for $\zeta=0$ we obtain the $e_3$-axis;
for $u=\frac{c+s}{2},\ v=\ep_1\frac{c-s}{2},\ \ep_1:=\pm 1, \ c$
constant we have ruled CGC $-1$ surfaces
$S(s,t):=x+t(\pa_ux+i\ep_2\sinh(\frac{u-v\cos(\zeta)}{\sin(\zeta)})\pa_vx),\
\ep_2:=\pm 1$; on Rogers-Schieff (\cite{RS1}, fig 1.5 pag 36)
asymptotes are obtained by drawing diagonals of the rectangles
obtained by lines of curvature $u,v$.

According to a theorem of Ribacour $x^0,\ x^1$ are the focal
surfaces of a W congruence iff
$\\K(x^0)K(x^1)=\frac{\sin^4(\be)}{d^4}$, where $\be$
(respectively $d$) is the angle (the distance) between the tangent
planes at the corresponding points (the corresponding points).
This is basically a statement about tangent curves (an asymptote
on $x^0\ (x^1)$ lies in the osculating bundle of its corresponding
asymptote on $x^1\ (x^0)$ and the relation between their torsions
becomes via Enneper's formula a relation between the curvatures of
$x^0$ and $x^1$; see Bianchi (\cite{B2},Vol {\bf 5},(117))). In
our case that becomes: $K(x^0)K(x^1)=\frac{|m_0^1\times
N_0^0|^4}{|m_0^1|^4|V_0^1|^4}=\frac{((N_0^0)^Tx_{zu_1}^1)^4}{|x_{zu_1}^1\times
V_0^1|^4}=^{(\ref{eq:V21c})}\frac{1}{\mathcal{A}^4|\hat
N_0^0|^4|\hat N_0^1|^4}$, which is true since
$K(x_0^j)=\frac{-1}{\mathcal{A}^2|\hat N_0^j|^4},\ j=0,1$.

Bianchi (122) has a beautiful geometric argument: using Cheffi's
result it is enough to prove the W congruence property only for
ruled surfaces applicable to quadrics, in which case the rulings
correspond under the B transformation (use (\ref{eq:dre})). Thus
asymptotes of one system (rulings) correspond, so Ribacour's
theorem applies and therefore the asymptotes of the other system
correspond. Note that although the rulings are plane curves, they
do not have zero torsion because the Frenet trihedral is not
defined; therefore the torsion is undefined, or equivalently
arbitrarily defined so as to suit our purposes.

Again the geometric picture gives the correct answer: the
curvature of a space curve measures its bending in the normal
direction (the normal is naturally defined from a geometric point
of view in the osculating plane of the curve) and the torsion
measures its twisting (tendency to leave its osculating plane).
Since a line does not admit a well defined osculating plane, the
curvature and torsion of a line cannot be defined from a geometric
point of view, or equivalently they can be arbitrarily defined
(usually as a limiting process) so as to suit our purposes. For
example the family of helices $x=x(s)=\la\cos(\frac{s}{\sqrt{\la
a}})e_1+\la\sin(\frac{s}{\sqrt{\la
a}})e_2+\frac{\sqrt{\la(a-\la)}}{\sqrt{\la a}}se_3$ has constant
curvature $\frac{1}{a}$, so we get the curvature of the line
$se_3$ to be $\frac{1}{a}$ for $\la\rightarrow 0$ and the tangent
surfaces of the family (all deformations one of the other) further
contains as a deformation the tangent bundle of the line $se_3$.
One can perform an experiment for the real plane region between
concentric circles of radii $a,\ a+\ep$ with a cut: this can be
tightly coiled around itself; as $\ep\rightarrow 0$ we get the
required line*. \footnote{* Again this discussion took place in
Calculus III; since students did not a-priori see the fact that
the projection of helices on the $xy$ plane become circles with
smaller radius I saw myself put into the position of further
deforming the piece of paper so it becomes clear that the
projection becomes a point in the limit (I was not prepared to
attack the point picture, but my students forced me to); thus they
spelled out the truth that the circles become smaller and in the
limit we obtain a point. I then asked {\it the natural question}
about the curvature of the line without being prepared with the
correct answer; when a student spelled out the truth that the line
has curvature $\frac{1}{a}$ I knew that he was right (Bianchi
proved him right) but that the discussion went longer than I was
prepared and I had to clarify it on the next meeting by trimming
the piece of paper to a smaller one-sided tubular neighborhood of
the given circle, because the piece of paper just with the disk
removed and a cut did not leave me coil it tightly and just
pulling it apart (as one would do with a coiled wire) would break
it.} The same experiment can be done with any one-sided tubular
neighborhood of a plane curve and thus a line can have any
curvature. Of course such a statement is ridiculous, but
statements and structures whose validity depend only on the linear
element of any of such tangent surfaces may remain valid and
interesting in the limit.

Because of the correspondence of asymptotes on $x^0,\ x^1$ we can
reduce (similarly to Lie's approach for CGC surfaces) the B
transformation of deformations of quadrics to the B transformation
of their asymptotes (if a curve corresponds with preserved
arc-length with a curve on a quadric and is an asymptote on a
deformation of that quadric, then its torsion is determined by
Enneper's formula) and conversely we can extend B transformations
of such curves to B transformations of deformations of quadrics.
Similarly we can extend B transformations (of deformations of
quadrics) to B transformations of triply conjugate families of
surfaces (and containing a family of deformations of quadrics). In
fact according to Bianchi (\cite{B2},Vol {\bf 4},(143)) the family
$x^1(c):=x^1(u_1(u_0,v_0,c),v_1(u_0,v_0,c))=B_{a_1}(x^0)$ of
singular B transforms of $x^0$ gives a family of surfaces (all
applicable to that same quadric) which can be completed to a
triply conjugate system of surfaces. All such deformations of the
considered quadric have as the considered conjugate system the
permanent conjugate system, so for CGC surfaces we obtain triply
orthogonal family of surfaces containing a family of CGC surfaces.

Since $s_0:=x^0, s_3:=x^1$ are the focal surfaces of a W
congruence, there exists $s_1:=r_1R_0m_0^1$ infinitesimal
deformation of $x^0:\
0=2zr_1^{-1}(ds_1^Tdx^0+(dx^0)^Tds_1)=2z(d(\log
r_1)(m_0^1)^Tdx_0^0 +(dx_0^0)^Tm_0^1d(\log
r_1)+(dm_0^1)^Tdx_0^0+(dx_0^0)^Tdm_0^1-2(m_0^1)^T(\om_0\times
dx_0^0))=(2zd(\log r_1)-\\(m_{0v_1}^1)^T\om_0)(m_0^1)^Tdx_0^0+
(dx_0^0)^Tm_0^1(2zd(\log r_1)-(m_{0v_1}^1)^T\om_0)$ (use
$(dm_0^1)^Tdx_0^0+(dx_0^0)^Tdm_0^1=\\dv_1(m_{0v_1}^1)^Tdx_0^0+(dx_0^0)^Tm_{0v_1}^1dv_1$,
(\ref{eq:intalg}) and $0=(dx_0^0)^T(\om_0\times dx_0^0)$).
Therefore $r_1$ should be given by $2zd(\log
r_1)=(m_{0v_1}^1)^T\om_0$, which must be integrable. But
differentiating (\ref{eq:intalg}) wrt $v_1$ we get
$(N_0^0)^T(2czm_{0v_1}^1+m_0^1\times m_{0v_1v_1}^1)=0$, which
ensures $d\wedge(m_{0v_1}^1)^T\om_0=0$. Thus we can say that $
r_1^{-1}$ is an infinitesimal deformation of $v_1:\ r_1^{-1}=\del
v_1$, or $r_1^{-1}=\pa_cv_1$ if $v_1=v_1(z,u_0,v_0,c)$. Therefore
$s_1$ and Darboux's 12 surfaces can be generated by a quadrature
(or no quadrature if we know the general solution of
$(\ref{eq:dify2})$ for fixed $z$). We have
$r_1^{-1}R_0^{-1}ds_1=dx_0^0\times(\mathcal{B}_1x_{zu_1}^1-
\om_0\lrcorner\frac{(\na x_0^0)^T(2zm_0^1+m_0^1\times
m_{0v_1}^1)}{2z})$, so
$s_{-1}:=r_1R_0(-\mathcal{B}_1x_{zu_1}^1+\om_0\lrcorner\frac{(\na
x_0^0)^T(2zm_0^1+m_0^1\times m_{0v_1}^1)}{2z})$ is the prwrtps of
$s_2:=r_2R_1m_1^0$ and $0=1+s_{-1}^Ts_2=1+r_1r_2\Del^-$ (use
(\ref{eq:V21c}) and the first equation of (\ref{eq:deltasi2}) for
$R_1m_1^0$). Thus
\begin{eqnarray}\label{eq:ro12}
1+r_1r_2\Del^-=0.
\end{eqnarray}

\subsection{Bianchi Permutability Theorem and moving M\"{o}bius configurations}
\label{subsec:deformations3}
\noindent

\noindent Let $B_{z_1}(x^0)=x^1=B_{z_2}(x^3)$. The BPT states that
there is $B_{z_2}(x^0)=x^2=B_{z_1}(x^3)$; moreover
$x^3=(B_{z_1}\circ B_{z_2})(x^0)=(B_{z_2}\circ B_{z_1})(x^0)$ is
obtained from $(x^0,x^1,x^2)$ only by algebraic computations.
\begin{center}
$\xymatrix{\ar@{}[dr]|{\#}x^2\ar@{<->}[d]_{B_{z_2}}\ar@{<->}[r]^{B_{z_1}}&x^3\ar@{<->}[d]^{B_{z_2}}\\
x^0\ar@{<->}[r]_{B_{z_1}}&x^1}$
\end{center}
Let $(R_j,t_j)(x_0^j,dx_0^j)=(x^j,dx^j),\ j=0,1,3$ be the rollings
known so far. One can find $x_0^2$ according to the algebraic
computations of the SITC; the sought surface
$x^2:=(R_0,x^0)V_0^2=(R_3,x^3)V_3^2$ is recovered using only
algebraic computations (we need to take derivatives just to find
$R_0$). From the cross ratio property we get $x^3=x^0$ for
$z_1=z_2$, so the Permutability Theorem does not a-priori get rid
of the quadratures for the iteration of the transformations $B_z$
with the same $z$. However an application of L'Hospital will
suffice once we know all $B_z(x^0)$ (see Bianchi (122)).

We have $x^2=(R_1,t_1)(R_3^0,t_3^0)x_0^2=(R_2,t_2)x_0^2,\
R_2:=R_1R_3^0,\ t_2:=(R_1,t_1)t_3^0$ and we need to prove
$dx^2=R_2dx_0^2$, that is
$R_1^{-1}dR_1R_0^1V_1^2+d(R_3^0,t_3^0)x_0^2=0$. Using
(\ref{eq:dify2}) this becomes an algebraic relation between
$x_0^j,\ j=0,...,3$ and linear in the second fundamental form of
$x^1$. Using (\ref{eq:syst}) for
$\frac{\om_1}{\mathcal{B}_1}=:(\om_{1u_1}^{u_1}x_{0u_1}^1-\om_{1u_1}^{v_1}x_{0v_1}^1)du_1+
(\om_{1v_21}^{u_1}x_{0u_1}^1-\om_{1v_1}^{v_1}x_{0v_1}^1)dv_1$ we
thus need the algebraic relations:
\begin{eqnarray}\label{eq:permuw}
\mathcal{B}_1x_{0u_1}^1\times R_0^1V_1^2
+(\frac{\Del^-(z_1,v_0,v_1)}{2z_1}(R_3^0,t_3^0)_{v_0}
+\frac{\Del^-(z_2,v_1,v_3)}{2z_2}(R_3^0,t_3^0)_{v_3})x_0^2=0,\nonumber\\
-\mathcal{B}_1x_{0v_1}^1\times R_0^1V_1^2
+(\frac{\Del^+(z_1,v_0,u_1)}{2z_1}(R_3^0,t_3^0)_{v_0}
+\frac{\Del^+(z_2,u_1,v_3)}{2z_2}(R_3^0,t_3^0)_{v_3})x_0^2=0.\nonumber\\
\end{eqnarray}
If $x^j,\ j=0,...,3$ are ruled with rulings corresponding to
$x_{0u_j}^j,\ j=0,...,3$ only the first relation of
(\ref{eq:permuw}) appears when the rulings correspond to
$x_{0u_0}^0,\ x_{0v_1}^1,\ x_{0v_2}^2,\ x_{0u_3}^3$ only the
second one appears.

Therefore we have reduced the proof of the BPT to the case of
ruled surfaces (Bianchi in (122) applies the inverse of Chieffi's
result for deformations of quadrics). Suppose $x^j,\ j=0,...,3$
are ruled (with rulings corresponding to $x_{0u_j}^j,\ j=0,...,3$)
satisfying the BPT. Using (\ref{eq:syst}) we have
$2z_1dv_0=\Del^-(z_1,v_0,v_1)\varphi_1(v_1)dv_1,\
2z_2dv_3=\Del^-(z_2,v_1,v_3)\varphi_1(v_1)dv_1,\
2z_2dv_0=\Del^-(z_2,v_0,v_2)\varphi_2(v_2)dv_2,\\
2z_1dv_3=\Del^-(z_1,v_2,v_3)\varphi_2(v_2)dv_2$, so
$\frac{z_2\Del^-(z_1,v_0,v_1)}
{z_1\Del^-(z_2,v_1,v_3)}=\frac{dv_0}{dv_3}=
\frac{z_1\Del^-(z_2,v_0,v_2)} {z_2\Del^-(z_1,v_2,v_3)}$. Therefore
\begin{eqnarray}\label{eq:cad}
\frac{z_2\Del^-(z_1,v_0,v_1)} {z_1\Del^-(z_2,v_1,v_3)}=
\frac{z_1\Del^-(z_2,v_0,v_2)} {z_2\Del^-(z_1,v_2,v_3)},\nonumber\\
\frac{z_2\Del^+(z_1,v_0,u_1)} {z_1\Del^+(z_2,u_1,v_3)}=
\frac{z_1\Del^+(z_2,v_0,u_2)} {z_2\Del^+(z_1,u_2,v_3)}
\end{eqnarray}
are necessary algebraic conditions for the BPT (the second one is
similarly obtained, when the rulings correspond to $x_{0u_0}^0,\
x_{0v_1}^1,\ x_{0v_2}^2,\ x_{0u_3}^3$). Applying to
(\ref{eq:cross}) the symmetry $(0\leftrightarrow
1,2\leftrightarrow 3)$, multiplying the two relations and using
(\ref{eq:deltasi2}) we get the first relation of (\ref{eq:cad})
(the second one is similar has to do with the cross-ratio of the
intersections of the rulings $x_{0v_1}^1,\ R_3^0x_{2v_2}^2$ with
the line $[x_{z_1}^0\ x_{z_2}^3]$ or we can use the last relation
of (\ref{eq:deltasi2})).

Conversely they are also sufficient for ruled surfaces (we shall
prove the first one; the second one is similar). If we are in the
conditions of the Permutability Theorem and we denote
$\varphi_2(v_2):=\frac{2z_2}{\Del^-(z_2,v_0,v_2)}\frac{dv_0}{dv_2}=
\frac{2z_1}{\Del^-(z_1,v_2,v_3)}\frac{dv_3}{dv_2}$, then the ruled
surface $x'^2$ with the nonzero diagonal entry of the second
fundamental form given by
$\mathcal{B}_2\sqrt{g_0^2}\varphi_2(v_2)$ is a $B_{z_2}$
($B_{z_1}$) transform of $x'^1,\ (x'^3)$ (which is applicable to
$x_0^0,\ (x_0^3)$).

Let $(x^j,dx^j)=(R_j,t_j)(x_0^j,dx_0^j),\ j=0,1,3,\
(x'^k,dx'^k)=(R'_k,t'_k)(x_0^k,dx_0^k),\ k=0,2,3$ be the rollings
known so far. Then using (\ref{eq:dre}) and (\ref{eq:cocy}) we
have $(R_0,t_0)(R'_0,t'_0)^{-1}=\\
(R_1,t_1)(R_3^0,t_3^0)(R'_2,t'_2)^{-1}=
(R_3,t_3)(R'_3,t'_3)^{-1}$. If $(R_0,t_0)(R'_0,t'_0)^{-1}$ is a
rigid motion, then $x^2:=(R_0,t_0)(R'_0,t'_0)^{-1}x'^2$ will
satisfy the conditions of the BPT. For this we need the second
fundamental forms of $x^0,\ x'^0$ to be the same, or
\begin{eqnarray}\label{eq:dv3dv2}
v_0=v_0(v_1),\ v_3=v_3(v_1),\ \frac{dv_3}{dv_0}=\frac{z_1\Del^-(z_2,v_1,v_3)}
{z_2\Del^-(z_1,v_0,v_1)}\Rightarrow
\frac{dv_2}{dv_1}=\frac{z_1\Del^-(z_2,v_0,v_2)}{z_2\Del^-(z_1,v_0,v_1)}.
\end{eqnarray}
While it seems that we are finished if we apply the symmetry
$(0\leftrightarrow 1,2\leftrightarrow 3)$ to the SITC in the
previous conditions, a reduction to an algebraic identity is
safer. If the lhs of (\ref{eq:phi2c})-(\ref{eq:phi2i1}) are
$\phi=\phi(z_1,z_2,v_0,v_1,v_2,v_3)$ and we differentiate $\phi=0$
with the constraints $v_0=v_0(v_1),\ v_3=v_3(v_1),\
\frac{dv_3}{dv_0}=\frac{z_1\Del^-(z_2,v_1,v_3)}{z_2\Del^-(z_1,v_0,v_1)}$,
then we get
$\frac{dv_0}{dv_1}[\phi_{v_0}+\phi_{v_3}\frac{dv_3}{dv_0}]
+\phi_{v_1}+\phi_{v_2}\frac{dv_2}{dv_1}=0$. If (\ref{eq:dv3dv2})
were to be true, then this last relation becomes a linear function
(with coefficients algebraic functions of $v_j,\
j=0,...,3,z_1,z_2$) of the a-priori arbitrary function
$\frac{dv_0}{dv_1}(v_1)$ being $0$, which is true if
$\phi_{v_0}+\phi_{v_3}\frac{dv_3}{dv_0}=\phi_{v_1}+\phi_{v_2}\frac{dv_2}{dv_1}=0$.
Therefore we need the algebraic relations
$\frac{\phi_{v_0}}{\phi_{v_3}}=-\frac{z_1\Del^-(z_2,v_1,v_3)}{z_2\Del^-(z_1,v_0,v_1)},\
\frac{\phi_{v_1}}{\phi_{v_2}}=-\frac{z_1\Del^-(z_2,v_0,v_2)}{z_2\Del^-(z_1,v_0,v_1)}$.
If either is proven, then by applying the symmetry
$(0\leftrightarrow 1,2\leftrightarrow 3)$ (which preserves $\phi$)
we get the other one. The fact that the rhs of (\ref{eq:calib}) is
skew wrt $(0\leftrightarrow 3,z_1\leftrightarrow z_2)$ finishes
the proof. Roughly the same argument works for changes
$u_j\leftrightarrow v_j$ for some $j$'s, so we also get
$B_{z_1}\circ B'_{z_2}= B'_{z_2}\circ B_{z_1},\ B'_{z_1}\circ
B'_{z_2}=B'_{z_2}\circ B'_{z_1}$.

An application of L'Hospital algebraically solves the problem of
the iteration of the B transformation for $z_2=z_1$: let
$v=v(u_0,v_0,z,b)$ be the solution of (\ref{eq:dify2}),
$v_1=v_1(u_0,v_0,z_1,b_1)$ a particular solution,
$v_3=v_3(u_0,v_0,z_1,z,b_1,b),\ z\neq z_1$ be given by
(\ref{eq:phi2c})-(\ref{eq:phi2i1}) with $z_2,\ v_2$ replaced by
$z,\ v$ and $b=b(z)$ a differentiable function with $b(z_1)=b_1,\
\frac{db}{dz}(z_1)=:b_0$. Letting now $z\rightarrow z_1$ in
(\ref{eq:phi2c})-(\ref{eq:phi2i1}) and applying L'Hospital we get
a definite value of $v_3=v_3(u_0,v_0,z_1,b_1,b_0)$, so the
iteration of the B transformation with the same $z_1$ (done by
means of algebraic computations and derivatives) introduces a new
constant as expected. Note that in this case $x^3$ has the same
rolling as $x^0:\ (x^3,dx^3)=(R_3,t_3)(x_0^3,dx_0^3),\
(R_3(u_3,v_3),t_3(u_3,v_3))=(R_0(u_0(u_3,v_3),v_0(u_3,v_3)),t_0(u_0(u_3,v_3),v_0(u_3,v_3)))$,
so the rolling of a quadric on an applicable surface gives a
$3$-dimensional family of deformations of that quadric having the
same rolling, a fact a-priori not evident. For the $B_z\circ B'_z$
transformation there is no need of L'Hospital; in fact these are
the G transformations discovered by Guichard independently of
Bianchi; they behave better at the analytic level and by use of
permanent conjugate systems (especially for real deformations of
real ellipsoids they give directly real deformations of the same
real ellipsoid). Because of the BPT M\"{o}bius configurations are
put in motion and become moving M\"{o}bius configurations once the
$n$ leaves $\{x^{2^{k-1}}\}_{k=1,...,n}$ are B transforms of the
seed $x^0$. The static $n$-soliton $x^{2^n-1}$ of $\mathcal{M}_n$
becomes an actual (moving) soliton if the seed is the vacuum
soliton; again one can apply the L'Hospital trick and even higher
derivatives when some of $z_1,...z_n$ are not distinct; again the
rolling can be considered for a higher dimensional family of
deformations of the considered quadric.

\subsection{Quadrics conjugate in deformation and the Hazzidakis transformation}
\label{subsec:deformations4}
\noindent

\noindent Bianchi in (\cite{B1},\S\ 423) calls two non-homothetic
non-flat surfaces conjugate in deformation if a pointwise
correspondence is established between them such that the
asymptotic coordinates and all the virtual asymptotic coordinates
correspond. Because of this the coefficients of (\ref{eq:virta})
on the conjugate in deformation $x,\ \ti x$ must be the same, or
if we choose $(u^1,u^2)$ asymptotic coordinates on $x,\ \ti x$:
\begin{eqnarray}\label{eq:via}
\pa_{u^j}\log\frac{K(\ti x)}{K(x)}=2(\Ga_{jj}^j-\ti\Ga_{jj}^j)=4(\Ga_{12}^{j+1}-\ti\Ga_{12}^{j+1}),\
\Ga_{jj}^{j+1}=\ti\Ga_{jj}^{j+1},\ j=1,2.
\end{eqnarray}
The geodesic equations
$\frac{d^2u^l}{ds^2}+\Ga_{jk}^l\frac{du^j}{ds}\frac{du^k}{ds}=0,\
l=1,2$ can be put in implicit form ($u^2=\phi(u^1)$) as
$(\frac{du^1}{ds})^2(\phi''
-(\phi')^3\Ga_{22}^1+(\phi')^2(\Ga_{22}^1-2\Ga_{12}^1)-\phi'(\Ga_{11}^1-2\Ga_{12}^2)+\Ga_{11}^2)=0$
and $\frac{du^1}{ds}=0\Rightarrow \Ga_{22}^1=0$ (Eisenhart
(\cite{E1},\S\ 85)). Thus the geodesics on $x,\ \ti x$ correspond.
Conversely if on two non-flat non-homothetic surfaces $x,\ \ti x$
the asymptotic coordinates and geodesics correspond, then the
second fundamental forms are proportional with factor
$\sqrt{\frac{K(\ti x)\ti g}{K(x)g}}$; from the Codazzi-Mainardi
equations we see that the equations of virtual asymptotic
coordinates (\ref{eq:virta}) on $x,\ \ti x$ are the same. Thus
conjugation in deformation is equivalent to correspondence of
asymptotic coordinates and geodesics. Bianchi proved in
(\cite{B2},Vol {\bf 5},(87),(90)) for quadrics of revolution and
Servant in {\cite{S11} essentially for the remaining ones that
only deformations of quadrics can be conjugate in deformation;
these appear according to Bianchi II for $n=3$. We shall reproduce
below this proof.

Dini solved the problem of finding all the surfaces with geodesic
correspondence (see  Darboux (\cite{D1},\S\ 600-\S\ 603)). First
one can choose on the two surfaces corresponding orthogonal
systems. Take a particular tangent facet $T$ of $x$; directions
conjugate to isotropic directions in $T$ and to the directions
corresponding to isotropic directions in the corresponding tangent
facet $\ti T$ of $\ti x$ are in general unique and distinct; if
they are not unique, then the linear elements of $x,\ \ti x$ are
proportional; if they are not distinct (Lie pointed out this gap
in Dini's proof), then $x,\ \ti x$ are flat and do not present
interest to us. Because (\ref{eq:via}) is coordinate invariant, it
is valid for the common orthogonal system $(u,v):\
|dx|^2=Edu^2+Gdv^2,\ \ |d\ti x|^2=\ti Edu^2+\ti Gdv^2$. From
$\Ga_{jj}^j-\ti\Ga_{jj}^j=2(\Ga_{12}^{j+1}-\ti\Ga_{12}^{j+1}),\
j=1,2$ we get $\ti E=\frac{E}{VU^2},\ \ti G=\frac{G}{UV^2},\
U=U(u),\ V=V(v)$. From $\Ga_{jj}^{j+1}=\ti\Ga_{jj}^{j+1},\ j=1,2$
we get $E=U_1^2(U-V),\ G=V_1^2(U-V),\ U_1=U_1(u),\ V_1=V_1(v)$ (we
exclude $U=V=$const which leads to homothetic $x,\ \ti x$), so
$|dx|^2=(U-V)(U_1^2du^2+V_1^2dv^2),\ |d\ti
x|^2=(\frac{1}{V}-\frac{1}{U})(\frac{U_1^2}{U}du^2+
\frac{V_1^2}{V}dv^2)$. Thus both $x,\ \ti x$ have Lioville linear
element (the coordinate curves are orthogonal geodesic confocal
ellipses and hyperbolas). This type of linear element generalizes
the linear element of quadrics in the elliptic coordinates. Under
the transformation $(U,V,U_1^2,V_1^2)\rightarrow
(m(U+h),m(V+h),\frac{U_1^2}{m},\frac{V_1^2}{m}),\ m,\ h$ constants
$|dx|^2$ is preserved while $|d\ti x|^2$ becomes up to the
constant factor $\frac{1}{m^3}$:
$(\frac{1}{V+h}-\frac{1}{U+h})(\frac{U_1^2}{U+h}du^2+\frac{V_1^2}{V+h}dv^2)$.
We can assume $|dx|^2=(u-v)(U^2du^2+V^2dv^2),\ |d\ti x|^2
=(\frac{1}{v+h}-\frac{1}{u+h})(\frac{U^2}{u+h}du^2+\frac{V^2}{v+h}dv^2)$.
Since $\Ga_{jj}^j+\Ga_{12}^{j+1}=\frac{1}{2}\pa_{u^j}\log(g),\
j=1,2$, (\ref{eq:via}) are equivalent to correspondence of
geodesics and (the coordinate invariant)
$cK(x)^{\frac{3}{4}}da=K(\ti x)^{\frac{3}{4}}d\ti a$ for some
constant $c$ (which we can assume to be $1$); here $da$ and $d\ti
a$ are the area forms. The last relation becomes:
$\pa_u(\frac{V}{U}\frac{1}{u-v})-\pa_v(\frac{U}{V}\frac{1}{u-v})=
\frac{1}{\sqrt{(u+h)(v+h)}}(\pa_u(\frac{V}{U}\frac{1}{u-v}\sqrt{\frac{v+h}{u+h}})
-\pa_v(\frac{U}{V}\frac{1}{u-v}\sqrt{\frac{u+h}{v+h}}))$, or
$\frac{2}{u-v}(\frac{1}{U^2}(\frac{1}{u+h}-1)+\frac{1}{V^2}(\frac{1}{v+h}-1))
+\pa_v(\frac{1}{V^2}(\frac{1}{v+h}-1))-\pa_u(\frac{1}{U^2}(\frac{1}{u+h}-1))=0$.
Thus $\frac{1}{U^2}(\frac{1}{u+h}-1)=au^2+bu+c,\
\frac{1}{V^2}(\frac{1}{v+h}-1)=-(av^2+bv+c)$; after a
normalization we have linear elements of quadrics in elliptic
coordinates: $|dx|^2=(u-v)(\frac{u}{(u-a)(u-b)(u-c)}du^2-
\frac{v}{(v-a)(v-b)(v-c)}dv^2),\ |d\ti x|^2=(\ti u-\ti v)
(\frac{\ti u}{(\ti u-\ti a)(\ti u-\ti b)(\ti u-\ti c)}d\ti u^2-
\frac{\ti v}{(\ti v-\ti a)(\ti v-\ti b)(\ti v-\ti c)}d\ti v^2),\
(\ti u,\ti v):=\frac{1}{(b-a)(c-a)}(\frac{u}{u-a},\frac{v}{v-a}),\
(\ti a,\ti b,\ti
c):=\frac{1}{(b-a)(c-a)}(1,\frac{b}{b-a},\frac{c}{c-a})$. For a
comprehensive discussion see Eisenhart (\cite{E2},\S\ 145).

Because of the virtual asymptotic coordinates correspondence, for
any deformation of $x$ we know infinitesimally a deformation of
the H transform $\ti x=H(x)$, so we need to integrate a Ricatti
equation and a quadrature to recover it.

However the H transformation will turn out to commute with the B
transformation ($B_{\ti z}\circ H=H\circ B_z$) and $(B_{\ti
z}\circ H)(x^0)=\ti x^1=(H\circ B_z)(x^0)$ can be algebraically
recovered from $x^0,\ x^1=B_z(x^0),\ \ti x^0=H(x^0)$.

In this vein the H transformation has the same involutory
character wrt the B transformation as the Christoffel
transformation wrt the D transformation of isothermic
surfaces.

Assume first that $\ti u_j,\ti v_j$ depend respectively only on
$u_j,v_j,\ j=0,1$. Since
$\frac{\ti\om_0}{\ti{\mathcal{B}}_0}=\frac{\ti{\mathcal{N}_0}}{\mathcal{N}_0}
((\om_{0u_0}^{u_0}\ti x_{0\ti u_0}^0-(\om_{0u_0}^{v_0}\frac{\ti
v_{0v_0}}{\ti u_{0u_0}}) \ti x_{0\ti v_0}^0)d\ti
u_0+((\om_{0v_0}^{u_0}\frac{\ti u_{0u_0}}{\ti v_{0v_0}})\ti
x_{0\ti u_0}^0 -\om_{0u_0}^{u_0}\ti x_{0\ti v_0}^0)d\ti v_0)$,
$\ti v_1$ given by the homography of Bianchi II satisfies $(\ti
m_0^1)^T\ti \om_0+2\ti zd\ti v_1=0$ if the algebraic relations
\begin{eqnarray}\label{eq:alr}
\frac{\ti\Del^-d\ti u_0}{\Del^-du_0}=\frac{\ti\Del^+d\ti v_0}{\Del^+dv_0}=
\frac{\ti zd\ti v_1/\ti{\mathcal{N}}_0}{zdv_1/\mathcal{N}_0}
\end{eqnarray}
are valid. But since (\ref{eq:du1}) and
$\widetilde{(\ref{eq:du1})}$ are valid for independent $u_0,\
v_0,\ v_1$, (\ref{eq:alr}) follows immediately. If $\ti u_j,\ti
v_j$ do not depend respectively only on $u_j,v_j,\ j=0,1$, then
the argument is similar.

Moreover there is an algebraic correspondence between
infinitesimal deformations of surfaces conjugate in deformation
and thus between W congruences having surfaces conjugate in
deformation as one of their focal surfaces (see Bianchi
({\cite{B1},\S\ 427)). If $s_0,\ s_3:=s_0+as_{0u}+bs_{0v},\
a=a(u,v),\ b=b(u,v)$ are the focal surfaces of a W congruence and
$\ti s_0$ is conjugate in deformation to $s_0$, then $\ti s_0,\
\ti s_3:=\ti s_0+a\ti s_{0u}+b\ti s_{0v}$ are the  focal surfaces
of a W congruence. To prove this we can choose $(u,v)$ common
conjugate system for the focal surfaces of the W congruence, so
$b=0$. Let $s_2:=rN_0$ be the infinitesimal deformation of $s_3$
from Darboux's 12 surfaces; using (\ref{eq:via}) $\ti
s_2:=r\sqrt[4]{\frac{K(s_0)}{K(\ti s_0)}}\ti N_0$ will be the
infinitesimal deformation of $\ti s_3$ from Darboux's 12 surfaces.

\subsection{(Singular) B\"{a}cklund transformations of (singular) quadrics}
\label{subsec:deformations5} \noindent

\noindent If we let $z$ approach the singular values $0,\ \infty$
or inverses $a_j$ of nonzero eigenvalues of $A$ (under
restrictions as follows: we consider only ruled $x^0$ for $z=0$,
only QC for $z=\infty$ and only $a_j$'s such that the eigenspace
of $a_j^{-1}$ is not isotropic; these restriction will reveal
themselves later), then we get singular B transformations. If we
allow one of the $a_j$'s tend to $0$ (again under restrictions as
above), then we get B transformations of singular quadrics. We
shall call $B_0$ {\it infinitesimal B transformation}, $B_{a_j}$
{\it finite singular B transformation} and $B_{\infty}$ {\it
infinite singular B transformation}.

For $z\rightarrow 0$ we have
$\frac{\Del^-}{z}\rightarrow\infty(v_1-v_0)^2,\
\frac{\Del'^+}{z}\rightarrow\infty(u_1-u_0)^2$ and
$\frac{\Del^+}{z},\ \frac{\Del'^-}{z}$ have definite limits. If
$v_1=v_0$, then $\om_{0u_0}^{v_0}=\om_{0u_0}^{u_0}=0$ and $x^1$ is
a particular position of $x_0^0$ as it rolls on the ruled $x^0$
(with ruling $x_{0u_0}$ corresponding to $x_{0u_0}^0$). If
$v_1\neq v_0$, then $\om_{0u_0}^{u_0}=\om_{0v_0}^{u_0}=0,\
dv_1=\lim_{z\rightarrow 0}\frac{\Del^+}{2z}\om_{0u_0}^{v_0}du_0$
and $x^0,\ x^1$ are ruled with rulings corresponding to
$x_{0v_0}^0,\ x_{0u_1}^1$. Note in particular
\begin{eqnarray}\label{eq:limd}
(v_1\rightarrow v_0)\circ\lim_{z\rightarrow 0}\frac{\Del^+}{z}=-\frac{1}{\mathcal{N}_0}.
\end{eqnarray}

In this case $x^1$ is obtained from $x^0$ by the application of
Chieffi's result and finding solutions of (\ref{eq:syst}) is
equivalent to finding asymptotes of $x^0$. All these results have
a simple geometric interpretation: as we roll a quadric on an
applicable ruled surface, its ruling families describe two
congruences decomposable in ruled surfaces applicable to that
quadric. The first strata are the particular positions of the
quadric as it rolls on its deformation; the other ones are
obtained by assembling the tangents (to geodesics corresponding to
the other ruling family on the quadric) along a curvilinear
asymptote of the initial ruled surface. We can consider
$z_1\rightarrow 0,\ z_2\neq 0$ or $z_1,z_2\rightarrow 0$ in the
BPT for $x^0$ ruled.

If we let $z$ tend to a finite singular value, then all equations
and quantities of interest, infinitesimal (fundamental  forms) or
finite (surfaces, rollings) have definite limits and enjoy the
same properties; the only difference is that the conics (in the
tangent spaces of $x^0$ and which are cut by four B transforms
with constant cross ratio) become lines. Note that although the
Ivory affinity degenerates, its action from $x_0$ to $x_{a_j}$ is
well defined (double cover), by means of which we can find a
double valued inverse (or an inverse from the double cover of
$x_{a_j}$ as a set, considered as the tangent bundle of the conic
$\mathcal{S}(x_{a_j})$) and complete the computations required to
produce the RMPIA. For quadrics of revolution ($a_1=a_2$ in QC
(\ref{eq:qc1}) or in QWC (\ref{eq:qwc1})) $B_{a_1}$ becomes the
complementary transformation; therefore $x^1$ will roll on the
axis of revolution as $x^0$ rolls on $x_0^0$. The BPT remains
valid when either (or both) of $z_1,z_2$ take finite singular
values. If the eigenspace of $a_j^{-1}$ is isotropic, then one
cannot define $\sqrt{R_{a_j}}$ and thus the Ivory affinity fails.
One can try to replace the Ivory affinity with another degenerate
affine transformation which should preserve as much as possible
from the properties of the Ivory affinity, but most of these are
lost and the program cannot be concluded (for example for QC one
can preserve the symmetry of the TC by choosing a symmetric linear
transformation linear combination of $J_3,J_3^2,e_3e_3^T$).
However limiting processes like the ones used in the discussion of
totally real quadrics may provide the correct tools needed to
extend the finite singular B transformation to isotropic
eigenspaces.

For example for QC (\ref{eq:qc1}) we only consider the singular
$B_{a_1}$ transformation (the other $B_{a_j}$'s are similar). The
singular quadric
$x_{a_1}(u,v):=\sqrt{R_{a_1}}(\sqrt{A})^{-1}X(u,v)$ covers
$\mathbb{C}^2$ twice: for any $\sqrt{R_{a_1}}(\sqrt{A})^{-1}[0\ a\
b]^T\in\mathbb{C}^2$ take $v$ a solution of $[i\ -b\ a]Y(v)=0$ and
$u:=-v\frac{1+b}{1-b}$. Through any point
$x_{a_1}(u_0,v_0)\in\mathbb{C}^2$ we have two tangents
$x_{a_1}(u,v_0),\ x_{a_1}(u_0,v),\ u,v\in\mathbb{C}$ at the conic
$\frac{(x^2)^2}{a_2-a_1}+\frac{(x^3)^2}{a_3-a_1}=1$ which touch
the conic for $u= 2v_0^{-1},\ v=2u_0^{-1}$. For $z\rightarrow a_1$
we have $\sqrt{R_z}\rightarrow \sqrt{R_{a_1}},\
\frac{\det\sqrt{R_z}(\sqrt{R_z})^{-1}\pm\sqrt{R_z}}{z\det\sqrt{A}}\rightarrow
\mathrm{diag}[\pm\sqrt{\frac{a_2}{a_1}(a_3-a_1)}\ \
\pm\sqrt{\frac{a_3}{a_1}(a_2-a_1)}\ \
\pm\sqrt{\frac{1}{a_1}(a_2-a_1)(a_3-a_1)}]$ (here the $\pm$ signs
depend in a certain fashion on $a_1,\ a_2,\ a_3$). If $\ti
a_2:=a_2^{-1},\ \ti a_j:=\ti a_2+(a_j-a_2)^{-1},\ j=1,3$, then
$\frac{\ti a_2}{\ti a_1}(\ti a_3-\ti
a_1)=\frac{-1}{(a_3-a_2)(a_2-a_1)}
\frac{1}{a_1}(a_2-a_1)(a_3-a_1),\ \frac{\ti a_3}{\ti a_1}(\ti
a_2-\ti a_1)
=\frac{-1}{(a_3-a_2)(a_2-a_1)}\frac{a_3}{a_1}(a_2-a_1),\
\frac{1}{\ti a_1}(\ti a_2-\ti a_1) (\ti a_3-\ti
a_1)=\frac{-1}{(a_3-a_2)(a_2-a_1)}\frac{a_2}{a_1}(a_3-a_1)$, so
$B_{\ti a_1}$ is a conjugation of $B_{a_1}$ with an H
transformation (Bianchi). A natural question is what would be the
corresponding B transformation conjugate of $B_{a_2}$ with a
Hazzidakis transformation; Calapso provided the natural answer:
since $\ti a_1$ corresponds to $a_1$ and $\infty$ corresponds to
$a_2$, the $B_{\infty}$ is the correct answer. If $\ti
a_1:=a_1^{-1},\ \ti a_j:=\ti a_1+(a_j-a_1)^{-1},\ j=2,3$, then
$\frac{\ti a_2}{\ti a_1}(\ti a_3-\ti
a_1)=\frac{1}{(a_3-a_1)(a_2-a_1)}a_2,\ \frac{\ti a_3}{\ti a_1}(\ti
a_2-\ti a_1)=\frac{1}{(a_3-a_1)(a_2-a_1)}a_3,\ \frac{1}{\ti
a_1}(\ti a_2-\ti a_1)(\ti a_3-\ti
a_1)=\frac{1}{(a_3-a_1)(a_2-a_1)}a_1$.

$B_{\infty}$ transformations of QC (\ref{eq:qc1}) were introduced
and studied by Calapso in \cite{C1}; he reduced the $B_{\infty}$
transformation via a conjugation by an H transformation to the
singular transformation $B_{a_j}$, thus naturally completing a
formula of Bianchi \cite{B1} involving only finite singular
transformations $B_{a_j}$'s.

Calapso studied $B$ transformations of QC (\ref{eq:qc1}) from an
intrinsic point of view; although Bianchi's simple geometric
approach is replaced with complicated analytic computations, it is
Calapso's approach which provides the correct setting in order to
study $B_{\infty}$ transformations, since as $z\rightarrow\infty$
finite quantities (leafs, rollings of leafs) do not have limit.
For example $B_{z}$ transforms of the seed $x^0$ tend to $\infty$,
$(R_0^1,t_0^1)$ does not have a rigid motion as limit and the
ACPIA and the W congruence disappears.

However Calapso realized that for QC (\ref{eq:qc1}) the TC (at the
analytic level), all infinitesimal quantities and all differential
equations have a definite valid limit (after proper scaling) when
$z\rightarrow\infty$. The direction fields $\hat
m_0^1:=\lim_{z\rightarrow\infty}z^{-1}m_0^1 =iY(v_1),\ \hat
m'^1_0:=\lim_{z\rightarrow\infty}z^{-1}m'^1_0=iY(u_1)$ are well
defined, finite and isotropic; they should be well defined finite
isotropic direction fields since for $z=a_j$, the facets with
normal fields $m_0^1, m'^1_0$ are tangent to the singular set
(conic) of the singular quadric $x_{a_j}$, so $m_0^1,\ m'^1_0$
must be tangent to $C(\infty)$ when $z=\infty$. A simple argument
in favor of the existence of the $B_{\infty}$ transformation for
all QC is that these are conjugate in deformation to certain
(I)Q(W)C, so $B_{\infty}$ transformations of QC are conjugation by
an H transformation of finite singular B transformations of
certain (I)Q(W)C. Conversely $B_{\infty}$ transformations for
(I)QWC do not exist, since they should be conjugation by an H
transformation of finite singular $B_{a_j}$ transformations with
$a_j^{-1}$ having isotropic eigenspace. Taking the limit
$z\rightarrow\infty$ in the scaled equations
(\ref{eq:dify2})-(\ref{eq:new1}) we get the intrinsic
characterization of the $B_{\infty}$ transformation. It is natural
to assume that the inversion of the $B_{\infty}$ transformation is
realized if in the inversion of the $B_z$ transformation one makes
$z\rightarrow\infty$ and one performs proper scalings, by means of
which we can find the second fundamental form of $x^1$ (so we
shall infinitesimally know $x^1=B_{\infty}(x^0)$).

As $z\rightarrow\infty$ the TC becomes $X_0^T\sqrt{A}X_1=0$;
(\ref{eq:dify2}) becomes the completely integrable Ricatti
equation $Y(v_1)^T\om_0=2idv_1$ and we have the similar equation
$Y(u_1)^T\om'_0=2idu_1$. Since (\ref{eq:del+-+-}), (\ref{eq:du1}),
(\ref{eq:alg}) (and its versions under symmetries) have a definite
valid limit when $z\rightarrow\infty$, the analytic version of the
inversion of the $B_{\infty}$ transformation holds.

For finite $z$ since $x^0, x^1$ are the focal surfaces of a W
congruence, there exists $s_1:=r_1R_0\frac{m_0^1}{z}$
infinitesimal deformation of $x^0$; as $z\rightarrow\infty$ the
limit of $r_1$ satisfies $d(\log r_1)-\frac{i}{2}Y'(v_1)^T\om_0=0$
and $ds_1=0,\ s_1:=r_1R_0Y(v_1)$. Conversely if we have an
isotropic vector $s_1$ and define $r_1, v_1$ by
$r_1Y(v_1):=R_0^{-1}s_1$, then
$0=\frac{-i}{2}Y'(v_1)^Tr_1^{-1}R_0^{-1}ds_1=Y(v_1)^T\om_0-2idv_1,\
0=r_1^{-1}R_0^{-1}ds_1=(d(\log
r_1)-\frac{i}{2}Y'(v_1)^T\om_1)Y(v_1)$ and we get a $B_{\infty}$
transformation of $x^0$. Therefore all $B_{\infty}$ transforms of
$x^0$ are determined by isotropic translations (this property was
communicated to Calapso by Bianchi). For the BPT if
$z_1\rightarrow\infty,\ z_2\in\mathbb{C}$, then (\ref{eq:cad}),
(\ref{eq:dv3dv2}), (\ref{eq:cross}), (\ref{eq:phi2c}) remain
valid. Therefore if $x^1=B_{\infty}(x^0),\ x^2=B_{z_2}(x^0)$, then
$x^3=B_{z_2}(x^1)=B_{\infty}(x^2)$ can be algebraically recovered
from $x^j,\ j=0,1,2$. If $z_1,\ z_2\rightarrow\infty$, then
a-priori (\ref{eq:phi2c}) has no definite limit. We can consider
(\ref{eq:cad}) in this case as the formula for the Permutability
Theorem; it becomes $X(v_0,v_3)^T\sqrt{A}X(v_1,v_2)=0$ and similar
formulae for $(u_1,u_2)$ replacing $(v_1,v_2)$, but we only know
$x^3$ infinitesimally from $x^j,\ j=0,1,2$.

For an isotropic vector $c:=mf_1+n\bar f_1+pe_3,\ 2mn+p^2=0,\
n\neq 0$ we have $c=\frac{n}{2}Y(\frac{p}{n})$. For $n=0$ we can
consider $c=mf_1=-mt^2Y(t^{-1})$ for $t\rightarrow 0$. We have
$c=rR_0Y(v)$, where $v:=\frac{c^TR_0e_3}{c^TR_0f_1},\
r:=\frac{c^TR_0f_1}{2}$. If $x^1,\ x^2$ are $B_{\infty}$
transforms of $x^0$, then $c_1:=r_1R_0Y(v_1),\ c_2:=r_2R_0Y(v_2)$
are constant isotropic vectors, so $R_0X(v_1,v_2)=\frac{\pm
c_1\times c_2}{|c_1\times c_2|}$ is a constant unit vector. Using
this fact we can reduce the determination of the transforms
$x^1=B_{\infty}(x^0)$ and of their iterates to quadratures. If we
fix $x^0,\ x^1$ and let $x^2,\ x^3$ vary in
$X(v_0,v_3)^T\sqrt{A}X(v_1,v_2)=0$, or equivalently
$X(v_0,v_3)=\frac{-iY(v_0)\times\sqrt{A}X(v_1,v_2)}{Y(v_0)^T\sqrt{A}X(v_1,v_2)}$,
then a homography will be established between $v_2,\ v_3$. $R_1$
can be algebraically recovered from $\frac{\pm c_0\times
c_3}{|c_0\times
c_3|}=R_1X(v_0,v_3)=R_1\frac{-iY(v_0)\times\sqrt{A}R_0^{-1}(c_1\times
c_2)}{Y(v_0)^T\sqrt{A}R_0^{-1}(c_1\times c_2)}$, where
$c_0:=r_0R_1Y(v_0),\ c_3:=r_3R_1Y(v_3)$ are constant isotropic
vectors; then we can recover $x^1$ from $dx^1=R_1dx_0^1$ with a
quadrature. Note that while we do not know the constant unit
vector $\frac{\pm c_0\times c_3}{|c_0\times c_3|}$, we can still
find $R_1$ up to a constant rotation. As a consequence the
$B_{\infty}$ transforms can be found with a quadrature; thus the
finite singular $B_{a_j}$ transforms can be found with a
quadrature once a certain H transform is known. Note again that
the previous argument is the same for
$(R_1,c_1)\leftrightarrow(R_2,c_2)$, so $R_1$ is the rotation of
the rolling of all $B_{\infty}$ transforms of $x^0$.

Note that while for QWC the $B_{\infty}$ does not exist, we come
close to having such a transformation: $\hat
u_1:=\lim_{z\rightarrow\infty}z^{-\frac{1}{2}}u_1,\ \hat
v_1:=\lim_{z\rightarrow\infty}z^{-\frac{1}{2}}v_1,\ \hat
m_0^1:=\lim_{z\rightarrow\infty}z^{-\frac{3}{2}}m_0^1=\hat
v_1^2Y(-i\hat v_1^{-1}),\ \hat
m'^1_0:=\lim_{z\rightarrow\infty}z^{-\frac{3}{2}}m'^1_0=-\frac{1}{2}Y(-2i\hat
u_1)$; $\hat u_1,\hat v_1$ satisfy the completely integrable
Ricatti equations $\\Y(-i\hat v_1^{-1})^T\om_0=2id(-i\hat
v_1^{-1}),\  Y(-2i\hat u_1)^T\om'_0=2id(-2i\hat u_1)$ (scaled
versions of (\ref{eq:dify2})) and the (non-symmetric) scaled
version of the TC $iZ_0^T\sqrt{A}\hat Z_1-e_3^T\hat
Z_1-\frac{1}{2}=0$. Again equations (\ref{eq:del+-+-}),
(\ref{eq:du1}), (\ref{eq:alg}) (and its versions under
symmetries), (\ref{eq:dify2})-(\ref{eq:new1}) have a definite
valid limit when $z\rightarrow\infty$ and after proper scaling, by
means of which one can prove the same equations with
$(u_0,v_0,\hat\Del^+,\hat\Del'^-,\mathcal{N}_0)\leftrightarrow
(\hat u_1,\hat
v_1,\hat\Del'^-,\hat\Del^+,\widehat{\mathcal{N}_1})$, but since
the scaled $\widehat{\mathcal{N}_1}$ is different than
$\mathcal{N}_{\hat 1}$ as a function of the scaled $\hat u_1,\hat
v_1$, the version of (\ref{eq:flatcon}) thus obtained loses its
geometric significance. The generation by isotropic translation
property remains, so we can find $\hat u_1,\hat v_1$ only by
algebraic computations.

Similarly we have the case for $a_j\rightarrow 0$. For example for
QC (\ref{eq:qc1}) and $a_1$ we replace $X_0$ with $[0\
\frac{i(2+u_0v_0)}{\sqrt{2}(u_0-v_0)}\ \frac{u_0+v_0}{u_0-v_0}]^T$
and $x_0^0$ becomes the tangent surface of a conic
$\mathcal{X}_0^0$; thus $x^0$ becomes the tangent surface of a
curve with the same curvature as that of $\mathcal{X}_0^0$ and the
rolling of $\mathcal{X}_0^0$ on $x^0$ becomes the rolling of two
developables with the osculating spaces of the lines of regression
coinciding. In this case $\om_0$ depends only on $u_0 (v_0)$ and
(\ref{eq:dify2}) becomes a Ricatti ODE upon whose integration we
find developables whose lines of regression have the same
curvature as that of $\mathcal{X}_0^0$. Therefore we can state
this geometrically: as the tangent surface $x_0^0$ of
$\mathcal{X}_0^0$ rolls with correspondence of lines of
regressions on an applicable developable the ruling families on a
confocal quadric will describe two line congruences. The confocal
conic (which is the intersection of the confocal quadric with the
plane of the singular quadric) rigidly moved in this rolling will
describe a focal surface for both these congruences. The
developables of these congruences have lines of regression on this
focal surface and will be applicable to $x_0^0$ (with
correspondence of lines of regression).

\subsection{(Totally real) deformations of totally real quadrics}
\label{subsec:deformations6} \noindent

\noindent We are interested in deformations $x^0$ of surfaces
 $x_0^0\subseteq x_0$ with real linear element $|dx_0^0|^2$.

Let $d\bar x_0^0=Rdx_0^0,\ R\subset\mathbf{O}_3(\mathbb{C}),\ \bar
R^T=R$. Since $R\hat N_0^0$ is a multiple of $\bar{\hat N}_0^0$
and keeping account of $K(x_0^0)\subset\mathbb{R}$, we get $R\hat
N_0^0=e^{ic}\bar{\hat N}_0^0,\ c\in\mathbb{R}$ constant. Applying
to this last relation $R^{-1}d$ we get $\om\times\hat N_0^0=
(e^{ic}R^{-1}\bar AR-A)dx_0^0$, so $(e^{ic}R^{-1}\bar AR-A)\hat
N_0^0=\la\hat N_0^0,\ \la\in\mathbb{C}$; by conjugation we obtain
$\bar\la=-e^{-ic}\la$. Assume $e^{ic}R^{-1}\bar AR-A\neq 0$; if
$v\in\mathbb{C}^3\setminus\{0\}$ is one of the remaining
eigenvectors of $e^{ic}R^{-1}\bar AR-A$ and having eigenvalue
$\mu\in\mathbb{C}$, then $v^T\hat N_0^0=0$ and $R^{-1}\bar v$ is
also eigenvector with eigenvalue $-e^{ic}\bar\mu$.

Note that $\bar{\hat N}_0^0\times\bar A\bar{\hat
N}_0^0=e^{-2ic}R\hat N_0^0\times\bar AR\hat
N_0^0=e^{-ic}\det(R)R^{-1}(\hat N_0^0\times A\hat N_0^0)$ combined
with (I)$\&$(II) of \S\ \ref{subsec:algpre22} give
$e^{2ic}=\frac{\ep\bar{\mathcal{A}}}{\mathcal{A}}$ and
$\mathcal{A}^2\in\mathbb{R},\ e^{4ic}=1$.

Since $\bar{\sqrt{a}}=\si_{a}\sqrt{\bar a},\ \si_{a}=1$ for $a$
not negative real and $-1$ otherwise, we have
$\overline{\sqrt{M}}=\si_{M}\sqrt{\bar M}$ where $M$ is SJ with
non-isotropic kernel and the signature $\si_M$ is diagonal with
$\pm 1$ entries (the sign is same for any SJ block). All
computations are preserved by conjugation; for example for QC
$(\bar R_0,\ \bar t_0)(\bar x_0^0,d\bar x_0^0)=(\bar x^0,d\bar
x^0),\ \bar x_0^0=\si_{A}(\sqrt{\bar A^{-1}})^{-1}\bar X_0
=(\sqrt{\bar A^{-1}})^{-1}X_0(\si_{X_0}(\bar u_0,\bar v_0)),\ \bar
x_z^1=...=\sqrt{\bar R_{\bar z}}(\sqrt{\bar
A^{-1}})^{-1}X_1(\si_{X_1}(\bar u_1, \bar v_1))$, where
$\si_{X_0},\ \si_{X_1}$ are compositions of some of the
involutions $\si_1^2,\si_2,\si_3$ of \S\ \ref{subsec:algpre21};
thus deforming the quadric $x_0$ is equivalent to deforming the
quadric $\bar x_0$. Since the condition $(R_0m_0^1)^Tdx^1=0$ holds
for $x^1=B_z(x^0)$, the condition $(\bar R_0\bar m_0^1)^Td\bar
x^1=0$ holds for $\bar x^1=\overline{B_z(x^0)}=B_{\bar z}(\bar
x^0)\ (\mathrm{or}\ =B'_{\bar z}(\bar x^0))$, so the $B_{\bar z},\
B'_{\bar z}$ transforms of $\bar x^0$ are conjugates of the $B_z,\
B'_z$ transforms of $x^0$.

If $x_0^0$ is totally real (the seed $x^0$ may be non-totally
real), then the B transformation behaves well wrt to these
surfaces (due mainly to its algebraic character). For $z$ and
initial value of the Ricatti equation (\ref{eq:dify2}) chosen to
satisfy certain totally real conditions the corresponding surface
$x_0^1\subseteq x_0$ of the leaf $x^1$ is totally real (not
necessarily of the same signature as that of $x_0^0$). However
most $B_z$ transforms of seeds $x^0$ with totally real linear
element do not obey these totally real requirements, in which case
the leafs $x^1$ and their corresponding surfaces $x_0^1\subseteq
x_0$ although being real $2$-dimensional may have complex linear
element. In this case the leaves are not interesting from a
totally real point of view, but they are still useful to generate
via conjugation and an application of the BPT surfaces with
$x_0^3$ or both $x_0^3,x^3$ being totally real respectively of the
same signature as that of $x^0,x_0^0$ and $x^3$ in the same
surrounding space as that of $x^0$.

If $\bar x^0=\ep x^0,\ \ep:=\mathrm{diag}[\ep_1\ \ep_2\ \ep_3]$
and $x^1=B_z(x^0)$, then $\bar x^1$ is a $B_{\bar z}\ (B'_{\bar
z})$ transform of $\bar x^0$. We would like to have $x^2:=\ep \bar
x^1$ a $B_{\bar z}\ (B'_{\bar z})$ transform of $x^0$. Once this
is accomplished, if $x^2\neq x^1$, then one can find $x^3$ from
the BPT applied to $(x^0,(x^1,z),(x^2,\bar z))$. Since the first
three vertices of the quadrilateral $(x^0,(x^1,z),(x^2,\bar
z)),x^3)$ are preserved under $((x^1,z)\leftrightarrow (x^2,\bar
z))\circ\ep\circ\bar{\cdot}$, we conclude that also the fourth one
has to be preserved: $\bar x^3=\ep x^3$. If $x_0^0$ is given by
$\begin{bmatrix}x_0^0\\1\end{bmatrix}^T
\begin{bmatrix}A&B\\B^T&C\end{bmatrix}\begin{bmatrix}x_0^0\\1\end{bmatrix}=0$,
then, since with $R:=\bar R_0^T\ep R_0,\ t:=\ep t_0-\bar R_0^T\bar
t_0$ we have $\begin{bmatrix}\bar
x_0^0\\1\end{bmatrix}=\begin{bmatrix}R&t\\0&1\end{bmatrix}
\begin{bmatrix}x_0^0\\1\end{bmatrix}$, we have $\\\begin{bmatrix}R&t\\0&1\end{bmatrix}^T
\begin{bmatrix}\bar A&\bar B\\\bar B^T&\bar C\end{bmatrix}\begin{bmatrix}R&t\\0&1\end{bmatrix}
=\begin{bmatrix}A&B\\B^T&C\end{bmatrix}$ (a quadric is uniquely
determined by the requirement that $3$ facets in general position
are among its tangent planes among the tangent planes of $x_0^0$
one can certainly find such facets; thus the above equality is
true, a-priori, modulo a constant $e^{ic}$, but it turns out that
$c\in\mathbb{R}$ and one can absorb it in $A,B,C$). This ensures
that $\begin{bmatrix}x_{\bar
z}^2\\1\end{bmatrix}:=\begin{bmatrix}R_0&t_0\\0&1\end{bmatrix}^{-1}
\begin{bmatrix}x^2\\1\end{bmatrix}=\\=\begin{bmatrix}R&t\\0&1\end{bmatrix}^{-1}
\begin{bmatrix}\bar x_z^1\\1\end{bmatrix}$ satisfies
$\begin{bmatrix}x_{\bar
z}^2\\1\end{bmatrix}^T(\begin{bmatrix}A&B\\B^T&C\end{bmatrix}^{-1}-
\bar
z\begin{bmatrix}I_3&0\\0^T&0\end{bmatrix})^{-1}\begin{bmatrix}x_{\bar
z}^2\\1\end{bmatrix}=0$. We need now only prove that if $w_0^1$ is
a ruling of $x_0$ at $x_0^1$ (that is
$(w_0^1)^TAw_0^1=(w_0^1)^T\hat N_0^1=0$), then
$w_0^2:=(\sqrt{R_{\bar z}})^{-1}R^T \overline{\sqrt{R_z}}\bar
w_0^1$ is a ruling of $x_0$ at $x_0^2$. This is certainly true,
since $RAR^T=\bar A$ implies $\\((\sqrt{R_{\bar
z}})^{-1}R^T\overline{\sqrt{R_z}})^TA (\sqrt{R_{\bar
z}})^{-1}R^T\overline{\sqrt{R_z}}=\bar A$; also $((\sqrt{R_{\bar
z}})^{-1}R^T\overline{\sqrt{R_z}})^T\hat N_0^2=\bar{\hat N}_0^1$
boils down to the above, $\bar At+\bar B=RB$ and $AC(\bar
z)+(I_3-\sqrt{R_{\bar z}})B= \bar
A\overline{C(z)}+(I_3-\overline{\sqrt{R_z}})\bar B=0$.

If $\bar x_0^0=\ep x_0^0,\ x^1=B_z(x^0)$, then as above $x_{\bar
z}^2:=\ep\bar x_z^1\subseteq x_{\bar z}$ satisfies the TC
requirement $(N_0^0)^T(x_{\bar z}^2-x_0^0)=0$. If $x_z^1\neq
x_{\bar z}^2$, then one can complete $(x_0^0,(x_z^1,z),(x_{\bar
z}^2,\bar z))$ to an SITC quadrilateral $(x_0^0,(x_z^1,z),(x_{\bar
z}^2,\bar z),x_0^3)$. As above we conclude that $\bar x_0^3=\ep
x_0^3$ we only need to show that $x^2:=(R_0,t_0)x_{\bar z}^2$ is a
$B_{\bar z}\ (B'_{\bar z})$ transform of $x^0$.

If $|dx_0^0|^2$ is real, then $d\bar x_0^0=Rdx_0^0,\
R\in\mathbf{O}_3(\mathbb{C}),\ \bar R^T=R$ and $d\bar x^0=\bar
R_0RR_0^Tdx^0$

Sometimes one gets deformations for free: if $|dx^0|^2$ is real,
the conjugation becomes a rigid motion for $x^0$, but $|dx^1|^2$
may not be real, in which case we get a new deformation $\bar
x^1$. If $x^0$ is real, then the composition of two distinct
$B_z,\ B_{\bar z}$ transformations (that is $x^1$ not real; we
always assume that $B_{\bar z}= \overline{B_z}$ is the conjugate
of the leaf in discussion) gives, by the BPT, a surface $x^3$;
applying the conjugation to the quadrilateral $(x^0,\ x^1,\ \bar
x^1,\ x^3)$ we get $\bar x^3=x^3$, so $x^3$ is real. We shall
prove later that if $x_0^0$ is totally real, then $x^3$ is
applicable to $x^0$.

This result can be extended to the case when there exists a rigid
motion $(R,t)$ such that $x^0=(R,t)\bar x^0$; then, if $x^1\neq
(R,t)\bar x^1$, the surface $x^3$ given by the BPT satisfies
$x^3=(R,t)\bar x^3$. In particular, if $x^0$ is totally real of
signature $\mathrm{diag}[\ep_1\ \ep_2\ \ep_3],\ \ep_j=\pm 1$, we
can take $(R,t):=(\mathrm{diag}[\ep_1\ \ep_2\ \ep_3],0)$; if $x^1$
is not totally real of the same signature as $x^0$, then $x^3$ is
totally real with the same signature as $x^0$; if furthermore
$x_0^0$ is totally real then $x^3$ is applicable to $x^0$. By
iteration we obtain M\"{o}bius configurations with $z_1,\ z_2,\
...,z_n$ and their conjugates; the composition $x^{2^{2n}-1}$ of
all the transformations $B_{z_j},\ B_{\bar z_j},\ j=1,...,n$ will
have the same character as $x^0$; moreover $x^{2^{2n}-1}$ will be
applicable to $x^0$ if $x_0^0$ is totally real.

If we are in the previous conditions ($x^0=(R,t)\bar x^0,\ x^1\neq
(R,t)\bar x^1,\ x^3=(R,t)\bar x^3$) and there is a rigid motion
$(R',t')$ such that $x_0^0=(R',t')\bar x_0^0$, then
$(R_0,t_0)(R',t')=(R,t)(\bar R_0,\bar t_0)$. We have
$\\x_0^3=(R_3,t_3)^{-1}(R,t)(\bar R_3,\bar t_3)^{-1}\bar x_0^3$;
if $(R_3,t_3)^{-1}(R,t)(\bar R_3,\bar t_3)^{-1}$ is a rigid
motion, then $x_0^3$ has a character similar to $x_0^0$. This is
the case for $x_0^0$ totally real: we shall see later that we have
$\\(\bar u_0,\bar v_0)=\si(u_0,v_0)$, where $\si$ is the
composition of some of $\si_1,\ \si_2,\ \si_3$ and uniquely
describes the totally real quadric $x_0^0$. Note that the SITC
provides homographies between $v_0,v_1,v_2,v_3$ and
$u_0,v_1,v_2,u_3$; using these two homographies and their
conjugates we get $(\bar u_3,\bar v_3)=\si(u_3,v_3)$, so $x_0^3$
coincides with $x_0^0$. For example for the real ellipsoid
$\si_{A}=I_3$ and we have $u_0\bar v_0=u_3\bar v_3=-2$.

The composition of complex conjugate B transformations (with
certain rationality conditions) of the axis of the tractrix
(vacuum soliton) gives a surface called a breather (see
Rogers-Schieff \cite{RS1}).

With $\si_1^2,\ \si_2,\ \si_3,\ \si_4$ the involutions of \S\ 2.1
and $\ep=:(-1)^{s(\ep)}$ the totally real unit spheres $X(u,v)$
have signatures as follows: $\mathrm{diag}[\ep\ \ep\ 1],\
\mathrm{diag}[\ep\ -\ep\ 1],\ \mathrm{diag}[\ep\ -\ep\ -1],\
\mathrm{diag}[1\ 1\ -1]$ respectively for
$\si_4=\si_1^{2(s(\ep)+1)}\si_2\si_3,\ \si_1^{2s(\ep)},\
\si_1^{2(s(\ep)+1)}\si_2,\ \si_3$. The totally real equilateral
paraboloids $Z(u,v)$ have signatures $\mathrm{diag}[\ep\ \ep\ 1],\
\mathrm{diag}[\ep\ -\ep\ 1]$ for $\si_4=\si_1^{2s(\ep)}\si_2,\
\si_1^{2s(\ep)}$. Thus all totally real QC (\ref{eq:qc1}) and QWC
(\ref{eq:qwc1}) with $a_j\in\mathbb{R}$ are obtained respectively
from totally real spheres and equilateral paraboloids by changing
the signature ot the $j^{\mathrm{th}}$ coordinate if $a_j<0$. We
need now only describe the totally real QC (\ref{eq:qc1}), QWC
(\ref{eq:qwc1}) with $\sqrt{a_1^{-1}}=\sqrt{\bar
a_2^{-1}}=b_1+ib_2\notin \mathbb{R}\cup i\mathbb{R}$, the
remaining quadrics (for which only $a_j\in\mathbb{R}$ is possible)
and some (I)QWC for which the canonical form does not admit
totally real realizations with the coordinates being real or
purely imaginary.

\subsection{Discrete deformations of quadrics}\label{subsec:deformations7}
\noindent

\noindent The study of M\"{o}bius configurations by Bianchi
(\cite{B2},Vol {\bf 5},(117)) and discrete CGC and isothermic
surfaces by Bobenko-Pinkall \cite{BP1},\cite{BP2} motivates the
study of DDQ.

A $\mathbb{Z}\times\mathbb{Z}$ lattice of B transforms $x^{(j,k)}$
of a deformation $x^{(0,0)}$ of a surface $x_0^{(0,0)}\subset
x_0$, $x_0$ being a quadric
($x^{(j+1,k+1)}=B_{z_j}(x^{(j,k+1)})=B_{z'_k}(x^{(j+1,k)}),\
z_j,z'_k\in\mathbb{C},\ j,k\in\mathbb{Z}$) forms a special
$2$-dimensional family of DDQ; a DDQ is formed by points
corresponding under the B transformations. All surfaces
$x^{(j,k)}$ can be algebraically obtained once we know the
surfaces on a connected {\it initial value tree} (IVT) which
contains once each transformation $B_{z_j},\ B_{z'_k},\
j,k\in\mathbb{Z}$ (for example the cross $x^{(j,0)},\ x^{(0,k)},\
j,k\in\mathbb{Z}$ or the staircase $x^{(j,j)},\ x^{(j+1,j)},\
j\in\mathbb{Z}$); more generally from a connected tree containing
once each transformation $B_{z_j},\ B_{z'_k},\ j,k\in\mathbb{Z}$
one can recover the corresponding double infinite M\"{o}bius
moving configuration and the DDQ as a subset of this (for example
the connected tree $x^{(0,0)},\ B_{z_j}(x^{(0,0)}),\
B_{z'_k}(x_{(0,0)}),\ j,k\in\mathbb{Z}$). The four segments
joining a point of a DDQ with its neighbors reside in the tangent
plane of the subjacent deformation, so a DDQ naturally admits a
discrete Gau\ss\ map.

From (\ref{eq:dre}) we get:
\begin{eqnarray}\label{eq:deriv}
(R_{(j+1,k)},t_{(j+1,k)})^{-1}(R_{(j,k)},t_{(j,k)})=
(R_{(j,k)}^{(j+1,k)},t_{(j,k)}^{(j+1,k)}),\nonumber\\
(R_{(j,k+1)},t_{(j,k+1)})^{-1}(R_{(j,k)},t_{(j,k)})=
(R_{(j,k)}^{(j,k+1)},t_{(j,k)}^{(j,k+1)})
\end{eqnarray}
and the permutability theorem allows us to recover all the
surfaces $x_{(j,k)}$ from an IVT using only algebraic
computations:
\begin{eqnarray}\label{eq:comm}
(R_{(j+1,k+1)},t_{(j+1,k+1)})^{-1}(R_{(j,k)},t_{(j,k)})
=(R_{(j+1,k)}^{(j,k+1)},t_{(j+1,k)}^{(j,k+1)}).
\end{eqnarray}
One can remove the scaffolding of subjacent surfaces; thus the
notion of DDQ can be parted with the notion of deformations of
quadrics.

A DDQ is a $\mathbb{Z}\times\mathbb{Z}$ lattice of facets
$(x^{(j,k)},T^{(j,k)})$ such that

{\it (i) $T^{(j,k)}$ contains the points of its $4$ neighboring
facets.

(ii) There is a  $\mathbb{Z}\times\mathbb{Z}$ lattice of facets
$(x_0^{(j,k)},T_{x_0^{(j,k)}}x_0)$ in the quadric $x_0$ and a
$\mathbb{Z}\times\mathbb{Z}$ lattice of rigid motions
$(R_{(j,k)},t_{(j,k)})$ (discrete rolling) satisfying
(\ref{eq:deriv}) and
$(R_{(j,k)},t_{(j,k)})(x_0^{(j,k)},T_{x_0^{(j,k)}}x_0)=(x^{(j,k)},T^{(j,k)})$.}

In requiring that (\ref{eq:deriv}) is satisfied we assume that
$(R_{(j,k)},t_{(j,k)})^{-1}(x^{(j,k+1)},T^{(j,k+1)})$ is the facet
centered at $x_{z'_k}^{(j,k+1)}$ and spanned by
$V_{(j,k)}^{(j,k+1)}$ and one of the rulings of $x_{z'_k}$ at
$x_{z'_k}^{(j,k+1)}$, etc.

The discrete rolling
$(x^{(j,k)},T^{(j,k)})=(R_{(j,k)},t_{(j,k)})(x_0^{(j,k)},T_{x_0^{(j,k)}}x_0)$
is a natural discrete version of the rolling
$(x,dx)=(R,t)(x_0,dx_0)$.

As the rotation of the rolling must have constant determinant, all
discrete rotations $R_{(j,k)}$ must have same determinant (we can
take it to be $1$), so the rulings involved in a homography of a
Bianchi quadrilateral must be the same; we can take them to be
$v$'s. However DDQ can have discrete rollings with discrete
rotations of different determinants, as long as the co-cycle
requirements are satisfied.

The discrete version of the compatibility condition
$d\wedge\om_0+\frac{1}{2}\om_0\times\wedge\om_0=0$ is
$\\R_{(j,k)}^{-1}R_{(j+1,k)}R_{(j+1,k)}^{-1}R_{(j+1,k+1)}=
R_{(j,k)}^{-1}R_{(j,k+1)}R_{(j,k+1)}^{-1}R_{(j+1,k+1)}$, or the
co-cycle relation
$\\(R_{(j+1,k)}^{(j,k)},t_{(j+1,k)}^{(j,k)})(R_{(j+1,k+1)}^{(j+1,k)},t_{(j+1,k+1)}^{(j+1,k)})
=(R_{(j,k+1)}^{(j,k)},t_{(j,k+1)}^{(j,k)})(R_{(j+1,k+1)}^{(j,k+1)},t_{(j+1,k+1)}^{(j,k+1)})$.

A DDQ is recovered only by algebraic computations from the
knowledge of the facets $(x^{(j,k)},T^{(j,k)})$ along an IVT and a
particular point $x_0^{(j_0,k_0)}$ on $x_0$ corresponding to a
point $x^{(j_0,k_0)}$ of a facet on that tree. That is so because
we know a rigid motion $(R_{(j_0,k_0)},t_{(j_0,k_0)})$ (modulo a
rotation of the facet into itself), then the confocal quadrics
upon which the points of the neighboring facets are taken by this
rigid motion, so by use of the RMPIA's we know the facets on $x_0$
corresponding to the neighboring facets (of the initial facet),
etc. Note that the facets along an IVT cannot be arbitrarily
prescribed: if we roll a facet to its corresponding point on $x_0$
by a rigid motion, then the neighboring rolled facets must contain
one of the rulings of the confocal quadric they are situated on
(there are three choices of such quadrics so 6 choices of such
facets), so if we prescribe the lengths of the segments of an IVT
and the angles between neighboring segments of an IVT, then the
DDQ is uniquely determined, modulo countable many choices.

While all discrete versions of the fundamental forms, Gau\ss\
curvature, etc are good as long as one replaces each derivative
with an edge, the correct discrete version is the one where
computations are exact and not approximate to order two. The
fundamental forms can be found at the level of a single Bianchi
quadrilateral; we consider the one from the discussion of the SITC
with the edges and the discrete projective Gau\ss\ map
highlighted:
\begin{center}
$\xymatrix{R_2^0\hat N_0^2\ar@{<-}[d]_{\mathcal{A}(R_2^0\hat
N_0^2)\times\hat N_0^0}
\ar@{-}[r]&\ar@{-}[r]^{\mathcal{A}(R_2^1\hat N_0^3)\times
R_2^0\hat N_0^2}&\ar@{->}[r]&
R_2^1\hat N_0^3\ar@{<-}[d]^{\mathcal{A}(R_2^1\hat N_0^3)\times R_1^0\hat N_0^1}\\
\hat N_0^0\ar@{-}[r]&\ar@{-}[r]^{\mathcal{A}(R_1^0\hat
N_0^1)\times\hat N_0^0} &\ar@{->}[r]&R_1^0\hat N_0^1}$
\end{center}

From (\ref{eq:V21c}) we have $x^0_u\rightsquigarrow
V_0^1=\mathcal{A}R_1^0\hat N_0^1\times\hat N_0^0,\
x^0_v\rightsquigarrow V_0^2=\mathcal{A}R_2^0\hat N_0^2\times\hat
N_0^0$. It is convenient to work with the projective Gau\ss\ map
$\hat N$ instead of $N$: $\hat N^0_u\rightsquigarrow R_1^0\hat
N_0^1-\hat N_0^0,\ \hat N^0_v\rightsquigarrow R_2^0\hat N_0^2-\hat
N_0^0$. From the existence of the Bianchi quadrilateral we get
$\mathcal{A}(R_2^0\hat N_0^2)\times\hat N_0^0+
\mathcal{A}(R_2^1\hat N_0^3)\times R_2^0\hat N_0^2=
\mathcal{A}(R_1^0\hat N_0^1)\times\hat N_0^0+
\mathcal{A}(R_2^1\hat N_0^3)\times R_1^0\hat N_0^1$; equivalently
\begin{eqnarray}\label{eq:ljk}
(R_2^1\hat N_0^3-\hat N_0^0)\times (R_1^0\hat N_0^1-R_2^0\hat
N_0^2)=0.
\end{eqnarray}
A similar relation for CGC $-1$ appears in Bobenko-Pinkall
\cite{BP2}. Thus the existence of the Bianchi quadrilateral is the
discrete version of $x^0_{uv}=x^0_{vu}$ (equivalently of the
compatibility condition $\om_0\times\wedge dx_0^0=0$).

Since $(x_u^0)^T\hat N^0_u\rightsquigarrow\mathcal{A}(R_1^0\hat
N_0^1\times\hat N_0^0)^T (R_1^0\hat N_0^1-\hat N_0^0)=0$, it is
clear that the piecewise linear curves $j,k=$const on the DDQ are
discrete versions of asymptotes. The discrete second fundamental
form is $(x_{uv}^0)^T\frac{\hat N^0}{|\hat N^0|}\rightsquigarrow
\mathcal{A}\frac{(\hat N_0^0)^T}{|\hat N_0^0|}((R_2^1\hat
N_0^3)\times R_2^0\hat N_0^2) =\mathcal{A}\frac{(\hat
N_0^0)^T}{|\hat N_0^0|}((R_2^1\hat N_0^3)\times R_1^0\hat N_0^1)$.
Using (\ref{eq:ljk}) the condition that the discrete Gau\ss\
curvature coincides with the Gau\ss\ curvature
$K(x_0^0)=\frac{-1}{\mathcal{A}^2|\hat N_0^0|^4}$ becomes $(\hat
N_0^0)^T((R_1^0\hat N_0^1)\times R_2^0\hat N_0^2)=0$, which is
straightforward.

Bobenko-Pinkall \cite{BP2} discuss many other issues, including
discrete B transformations and the corresponding discrete
permutability theorem. Most of their results should remain valid
for DDQ as long as the correct metric-projective equivalents are
found.

Now the existence of the $\mathcal{M}_3$ configuration has to do
with commuting of third order derivatives, so the algebraic
relations needed for the existence of this configuration may be
the correct discretization of the GCM equations.

\subsection{Delaunay's, Guichard's and Weingarten's results}\label{subsec:deformations8}
\noindent

\noindent
\subsubsection{Delaunay's result and its converse}

All diagonal parabolas are homothetic to $\\p_0=[2\sinh (u)\ \sinh
^2(u)]^T$ with focus $f=[0\ 1]^T$ we have the rolling of $p_0$ on
the line $p=[0\ \cosh(u)\sinh(u)+u)]^T$: $(p,dp)=(R,t)(p_0,dp_0)$,
where $R$ is
the rotation $\begin{bmatrix} \tanh(u)&-\frac{1}{\cosh(u)}\\
\frac{1}{\cosh(u)}&\tanh(u)
\end{bmatrix}$.

The curve described by the focus through the rolling of $p_0$ on
$p$ is the catenary $p+R(f-p_0)=[\cosh(u),u]^T$, so we get the
statement for minimal surfaces of revolution.

All diagonal conics with center are homothetic to $e_0=[\cos(u)\
b\sin(u)]^T$ with excess $E,\ E^2:=b^2-1$ and foci $f=[0\ E]^T$
(there are two choices for the sign of $E$). We have the rolling
of $e_0$ on the line $e=[0\ \int \sqrt{1+E^2\cos^2(u)}du]^T$:
$(e,de)=(R,t)(e_0,de_0)$, where $R$ is the rotation
$\frac{1}{\sqrt{1+E^2\cos^2(u)}}\begin{bmatrix}
b\cos(u)&\sin(u)\\
-\sin(u)&b\cos(u)
\end{bmatrix}$.
The curves described by the foci through the rolling of $e_0$ on
$e$ are
$c=e+R(f-e_0)=e+\frac{E\sin(u)-b}{\sqrt{1+E^2\cos^2(u)}}[1\ -
E\cos(u)]^T$ and if we denote $[\cos(r)\
\sin(r)]^T:=\frac{1}{\sqrt{1+E^2\cos^2(u)}}[1\ - E\cos(u)]^T$,
then we have
$dr=-\cot(r)d\log(\cos(r))=\frac{E\sin(u)du}{1+E^2\cos^2(u)}$ and
$dc=\frac{1}{\sqrt{1+E^2\cos^2(u)}}[E\cos(u)\
1]^Tdu+(E\sin(u)-b)d(\frac{1}{\sqrt{1+E^2\cos^2(u)}}[1\ -
E\cos(u)]^T)=(du+(E\sin(u)-b)dr)[-\sin(r)\
\cos(r)]^T=\frac{b(b-E\sin(u))du}{1+E^2\cos^2(u)}[-\sin(r)\
\cos(r)]^T$.

We have $dc\times
d^2c=\frac{E\sin(u)[b(b-E\sin(u))]^2du^3}{[1+E^2\cos^2(u)]^3}[0\
0\ 1]^T$ and the curvature $k_c$ of $c$ is $\frac{|dc\times
d^2c|}{|dc|^3}=\frac{E\sin(u)}{b(b-E\sin(u))}$. The mean curvature
of the surface obtained by the rotation of $c$ around the $y$ axis
is the sum of $k_c$  and $\cos(r)$ of the inverse of the $x$
coordinate (the projection of the curvature of the parallel of the
surface on the normal to the surface), that is
$\frac{E\sin(u)}{b(b-E\sin(u))}
+\frac{1}{\sqrt{1+E^2\cos^2(u)}}\frac{\sqrt{1+E^2\cos^2(u)}}
{E\sin(u)-b}=-\frac{1}{b}$, so we obtain CMC surface of
revolution.

For the general conic with center ($e_0=[a\cos(u)\ b\sin(u)]^T$)
we get the mean curvature $-\frac{1}{b}$ since the mean curvature
scales inversely. Also because CMC surfaces of revolution are
described in arc-length parametrization of the meridian by the
$2^{\mathrm{nd}}$ order ODE
$H=\sqrt{1-x'^2}(\frac{1}{x}-\frac{x''}{1-x'^2})$ all CMC surface
of revolution are obtained in this way (pick initial conditions at
a point where the $x$ coordinate is critical, so we have only two
parameters).

For $u,a,b\in\mathbb{R},\ b>a>0\
(i(u-\frac{\pi}{2}),ia,b\in\mathbb{R})$ we have real ellipse
(hyperbola), real foci and real CMC (the real ellipse $[a\cos(u)\
b\sin(u)]^T,\ a,b,u\in\mathbb{R},\ b>a>0$ (hyperbola $[a\sinh(u)\
b\cosh(u)]^T,\\ a,b,u\in\mathbb{R}$) has also imaginary foci $[\pm
i\sqrt{b^2-a^2}\ 0)]^T$ ($[\pm i\sqrt{a^2+b^2}\ 0]^T$)).

Conversely we have the rolling of the curve $c_0$ on the line
$c=[0\ \int|dc_0|]^T:\ (c,dc)=(R,t)(c_0,dc_0),\\ R^T=[N_0\
\frac{dc_0}{|dc_0|}],\ t=c-Rc_0$. Then $t$ is the curve described
by a point (which can be considered the origin) rigidly moved with
$c_0$ as it rolls on $c$. We have:
$R^{-1}dt=\frac{N_0^TR^{-1}dRdc_0}{|dc_0|^2}(N_0\times dc_0)\times
c_0,\ |dt|^2=k_{c_0}^2|c_0|^2|dc_0|^2,\ N_t=\frac{1}{|c_0|}Rc_0,\
k_t=\frac{N_t^Td^2t}{|dt|^2}=\frac{1}{|c_0||dt|^2}(|dt|^2-N_0^Tc_0N_0^TR^{-1}dRdc_0)=
\frac{1}{|c_0|}+\frac{N_0^Tc_0}{k_{c_0}|c_0|^3}$. Since the
surface of revolution described by the rotation of $t$ around the
$y$ axis has CMC, $t$ coincides with the curve described by a
(the) focus of a conic as it rolls on the $y$ axis is uniquely
determined modulo rigid motions by $|dt|,\ k_t$. But these
quantities uniquely determine $|c_0|,\ |dc_0|$, which in turn
uniquely determine modulo rigid motions $c_0$ (since
$(N_0^Tc_0)^2=|c_0|^2-\frac{1}{4|dc_0|^2}(d|c_0|^2)^2$ and
$k_{c_0}N_0^tc_0|dc_0|=c_0^Td\frac{dc_0}{|dc_0|}=d\frac{d|c_0|^2}{2|dc_0|}-|dc_0|$).
Thus $c_0$ is uniquely determined and must be a conic.

\subsubsection{Guichard's result and its converse} (see Bianchi (\cite{B1}, ch XX))

For this subsubsection we use the invariant notation established
at the beginning of \ref{subsec:rolling2}.

Bianchi proposed and solved the inversion of the Guichard result:
find all surfaces $x_0$ and fixed points (which can be considered
the origin) such that all the surfaces $t$ described by this point
through the rolling of $x_0$ on an applicable surface $x$ have the
same CMC $H(t)$ (including the statement in space forms and CGC
$K(t)$).

From (\ref{eq:cmc}) we thus need:
\begin{eqnarray}\label{eq:secon}
x_{0(j+1)(k+1)}=(-1)^{j+k+1}n^{-1}B^{jk}g,\
2K(x_0)n+x_{0jk}B^{jk}=0;
\end{eqnarray}
using Gau\ss\ theorem for $x_0$ we get from the first three
relations $K(x_0)n^2=(N_0^Tx_0)^2|x_0|^2$, so the last relation is
a consequence of the first three.

We are interested in the Gau\ss\ equation of $x_0$ for the square
of the distance from a point and distance from a plane (Darboux
equations).

Let $\rho:=|x_0|^2,\ x_0=(N_0N_0^T+dx_0\llcorner\na x_0^T)x_0$, so
$4(N_0^Tx_0)^2=4\rho-|\na \rho|^2$.

Also
$\na^2\rho=2(x_0^T\na^2x_0+|dx_0|^2)=2(x_0^TN_0N_0^Td^2x_0+|dx_0|^2)$
and thus the Gau\ss\ equation becomes, with the area form
$da:=\sqrt{\det(|dx_0|^2)}=\sqrt{g}du^1\wedge du^2$:
\begin{eqnarray}\label{eq:da1}
\det(\na^2\rho-2|dx_0|^2)=(4\rho-|\na\rho|^2)K(x_0)da^2.
\end{eqnarray}
For $e$ constant unit vector let $\tau:=e^Tx_0$; then
$e=(N_0N_0^T+dx_0\llcorner\na x_0^T)e$, so
$(e^TN_0)^2=1-|\na\tau|^2,\
\na^2\tau=e^T\na^2x_0=e^TN_0N_0^Td^2x_0$, so the Gau\ss\ equation
becomes:
\begin{eqnarray}\label{eq:da2}
\det(\na^2\tau)=(1-|\na\tau|^2)K(x_0)da^2.
\end{eqnarray}
Once either version of the Gau\ss\ equation is satisfied the
Codazzi-Mainardi equations are automatically satisfied.

Also from (\ref{eq:da1}) one can easily read Darboux's result
about cyclic systems:  circles of cyclic systems have centers at
$x_0+\frac{\na x_0(\mu)}{2}$ and radii $\mu+\frac{|\na\mu|^2}{4}$
where $-\mu$ satisfies the Gau\ss\ equation for the square of the
distance from a point. With $\mu^2:=\rho$ we have from
(\ref{eq:secon})
$\frac{\na^2(\mu^2)-2|dx_0|^2}{2\mu\sqrt{1-|\na\mu|^2}}=N_0^Td^2x_0
=-\frac{(da\llcorner\na\mu)^2+(1-|\na\mu|^2)|dx_0|^2}
{2(H(t)\mu+1)\mu\sqrt{1-|\na\mu|^2}}$, so we get
\begin{eqnarray}\label{eq:mupa}
\mu\na^2\mu+d\mu^2-|dx_0|^2=
-\frac{(da\llcorner\na\mu)^2+(1-|\na\mu|^2)|dx_0|^2}
{2(H(t)\mu+1)}.
\end{eqnarray}
Applying to this $2\llcorner\na \mu$ and using $d(|\na
\mu|^2)=2\na^2\mu\llcorner\na \mu$ we get $\mu
d(|\na\mu|^2)+2|\na\mu|^2d\mu-2d\mu=
-\frac{(1-|\na\mu|^2)d\mu}{H(t)\mu+1}$, or
$\frac{-d(|\na\mu|^2)}{1-|\na\mu|^2}=\frac{(2H(t)\mu+1)d\mu}{\mu(H(t)\mu+1)}$,
so $d\log(1-|\na\mu|^2)=d\log(\frac{1}{\mu(H(t)\mu+1)})$; thus
$|\na\mu|^2=1-\frac{c}{\mu(H(t)\mu+1)}$ for some constant $c$.
Applying $\llcorner\llcorner\pa x_0^2$ to (\ref{eq:mupa}) we get
$\mu\Del\mu+|\na\mu|^2-2=-\frac{|\na\mu|^2+2(1-|\na\mu|^2)}{2(H(t)\mu+1)},\\
\Del\mu:=\na^2\mu\llcorner\llcorner\pa x_0^2$. Applying again
$\llcorner\llcorner(d\mu\llcorner\pa a)^2,\ \pa
a:=(\sqrt{g})^{-1}\pa_{u^1}\wedge\pa_{u^2}$ to (\ref{eq:mupa}) we
get $\\\mu\na^2\mu\llcorner\llcorner(d\mu\llcorner\pa
a)^2-|\na\mu|^2=
-\frac{|\na\mu|^4+(1-|\na\mu|^2)|\na\mu|^2}{2(H(t)\mu+1)}$.
Equation (\ref{eq:da1}) can be rewritten as:
$\mu^2\det(\na^2\mu)+(\mu\na^2\mu\llcorner\llcorner(d\mu\llcorner\pa
a)^2 -\mu\Del\mu-|\na\mu|^2+1)da^2=\mu^2(1-|\na\mu|^2)K(x_0)da^2$,
or
$\det(\na^2\mu)=(1-|\na\mu|^2)(K(x_0)+\frac{H(t)}{\mu(H(t)\mu+1)})da^2$.
If we set $\tau:=k\mu,\ k:=\frac{1}{\sqrt{1+4cH(t)}}$, then
$1-|\na\tau|^2=k^2(1-|\na\mu|^2)(1+2H(t)\mu)^2,\\
K(x_0)=\frac{1}{4\mu^2(H(t)\mu+1)^2}$ and $\tau$ satisfies
(\ref{eq:da2}).

Thus $x_0$ has the property that the ratio of the distances to a
point and to a plane is a constant, so $x_0$ is a quadric of
revolution around the focal axis (with center if $k\neq 1$ and
without center if $k=1$).

Bianchi (\cite{B1}, ch XX) proved that all CMC surfaces occur as
translations of rolling of quadrics of revolution (around the
focal axis) in a $2$-dimensional fashion for a fixed quadric; this
is how the $2$-dimensionality of the B transformation for quadrics
of revolution was first established. In particular all real
surfaces applicable to real quadrics of revolution around the
focal axis are real analytic in the common conjugate system of
coordinates (a conformal parametrization of the second fundamental
form), which corresponds to the real analytic lines of curvature
parametrization on the real CMC surface.

The reality conditions are not necessary for the validity of the
computations, so we get statements about surfaces in
$\mathbb{C}^3$ applicable to quadrics of revolution and CMC
surfaces in $\mathbb{C}^3$.

\subsubsection{Darboux's geometrization of Weingarten's result}
(see Bianchi (\cite{B2},Vol {\bf 4},(202)))

\subsection{Conic systems}\label{subsec:deformations9}
\noindent

\noindent Motivated by Darboux's and Ribacour's result about
cyclic systems, consider the integrable $3$-dimensional
distribution of facets obtained as follows: intersect the tangent
planes of a surface $y_0$ with a quadric $x_z$ (thus along
conics); for each point of the conic take the facet determined by
the segment joining the point of the conic with the origin of the
tangent plane and one of the rulings passing through that point.
As we roll $y_0$ on any applicable surface $y:\
(y,dy)=(R,t)(y_0,dy_0)$ while rigidly moving this distribution in
the process, does this distribution always remain integrable? The
case when $x_z$ is non-isotropic singular quadric (tangent bundle
of a conic) was classified by Bianchi (\cite{B2},Vol {\bf
4},(173)): $y_0$ must be a quadric confocal to $x_z$ and we
recover the finite singular B transformations of $x_0$ (Bianchi in
fact showed more: a non-isotropic plane must also have the
projective structure of a singular quadric, as a-priori the choice
of rulings is not well defined in a plane).

Again it is convenient to derive certain identities before
imposing the TC: consider $x_z=x_z(u,v),\ |N_0|=1,\ N_0^Tdy_0=0,\
V:=x_z-y_0$ and (\ref{eq:baseb}) becomes
\begin{eqnarray}\label{eq:basebc}
(V_vV_u^T+V_uV_v^T)N_0=
V_{uv}V^TN_0+\frac{2}{\mathcal{B}}(A^{\circ -1}-zI_3+\Ga)N_0,\nonumber\\
\Ga:=-(V+y_0)y_0^T\ \mathrm{for\ QC},\
:=(V+y_0)(Le_3)^T+Le_3y_0^T\ \mathrm{for\ (I)QWC}.
\end{eqnarray}
Again multiplying (\ref{eq:basebc}) on the left respectively with
$N_0^T,(V\times N_0)^T$ in the TC assumption $V^TN_0=0$ and with
$n:=(A^{\circ -1}-y_0y_0^T)N_0$ for QC, $:=(A^{\circ
-1}+y_0(Le_3)^T+Le_3y_0^T)N_0$ for (I)QWC,
$m:=\mathcal{B}x_{zu}\times V,\ m':=\mathcal{B}x_{zv}\times V$ we
get the analogous of (\ref{eq:simpas}):
\begin{eqnarray}\label{eq:simp3}
N_0^Tx_{zu}N_0^Tx_{zv}\mathcal{B}=n^TN_0-z,\nonumber\\
x_{zv}^T(I_3-2N_0N_0^T)m=-2V^T(N_0\times n).
\end{eqnarray}
The condition (\ref{eq:dify1}) remains the same; equation
(\ref{eq:dify2}) becomes
\begin{eqnarray}\label{eq:difyc}
m^T\om m^TN_0 +2m^T(zN_0-n)dv=0,\ \om:=\al^{-1}(R^{-1}dR).
\end{eqnarray}
Since applying directly (\ref{eq:disa}) is complicated (we have
$V=\frac{m\times N_0}{N_0^T\mathcal{B}x_{zu}}$), direct
computations again must be applied as in the B transformation
case. If we let $z\rightarrow\infty$ for $x_z$ QC, then we get
$N_0^TX=0,\ Y(v)^T\om=2idv$ again isotropic translations, just as
for $B_{\infty}$. This can also be considered as infinite cyclic
systems: we intersect the tangent spaces of a surface with a null
cone with center at $\infty$; the Darboux equation for the
distance from a point to a surface becomes the equation for the
distance from a point at $\infty$ to the surface, that is the
distance from a plane to the surface. This can be easier seen by
considering the surrounding space of constant curvature $k$ and
letting $k\rightarrow 0$.

Thus applying $d\wedge$ to (\ref{eq:difyc}), using the equation
itself, $dm=m_vdv+dy_0\times\mathcal{B}x_{zu}$, (\ref{eq:om}) and
$dm^T\wedge\om=dv\wedge m_v^T\om$ we get $(m^TN_0)^2(N_0^T(m\times
m_v+2zm)-2m^Tn)N_0^T(\frac{1}{2}\om\times\wedge\om)
+2m^T\om\wedge(-(dy_0\times\mathcal{B}x_{zu})^T(N_0\times n)\times
m+m^TN_0m^T(zdN_0-dn)-m^T(zN_0-n)m^TdN_0)=0$. Note $((N_0\times
n)\times m)\times\mathcal{B}x_{zu}=-(N_0\times
n)^T\mathcal{B}x_{zu}m$ and the remaining two terms in the second
parenthesis give $m^T(ndN_0^T-N_0dn^T)m$. Similarly to
(\ref{eq:intalg}) we can use (\ref{eq:simp3}) to get
$N_0^T(m\times m_v+2zm)=\frac{2zV^T(N_0\times
n)-n^TN_0V^Tm_v}{N_0^Tx_{zv}}=2m^Tn$ (use $V^Tm_v=-m^Tx_{zv}$ and
further (\ref{eq:simp3}) again). Thus we get the restriction
$2m^T\om\wedge(dy_0^Tm(\mathcal{B}x_{zu})^T(N_0\times
n)+m^T(ndN_0^T-N_0dn^T)m)=0$; this is a polynomial of degree $\le
6$ in $v$ which must be identically $0$ for any $v$ and choice of
$\om$ satisfying (\ref{eq:om}) and thus the restrictions on the
geometry of $y_0$ revealed themselves: $m^Tdy_0=0$ (in general
this is over-determined; the infinitesimal B transformation may be
considered an example) or $dy_0^Tm(\mathcal{B}x_{zu})^T(N_0\times
n)+m^T(ndN_0^T-N_0dn^T)m=0$.

\section{Particular deformations in $\mathbb{C}^{n+p}$ of quadrics in
$\mathbb{C}^{n+1}$}\label{sec:particular} \setcounter{equation}{0}

\subsection{Overview and particular deformations in $\mathbb{C}^{2+p}$ of quadrics in
$\mathbb{C}^{3}$}
\noindent

\noindent Some results about surfaces do not generalize for
sub-manifolds in scalar product spaces. If the co-dimension is
lower, the GCMR equation will impose higher restrictions. However
the Ricci equations (partially) disappear if we only consider
deformations with (partially) flat normal bundle; further the
existence of imaginary numbers relieves part of restrictions (for
example according to Cartan \cite{C2} $\mathbf{H}^n(\mathbb{R})$
cannot be locally immersed in $\mathbb{R}^{n+p},\ p<n-1$, but it
can be embedded in $\mathbb{R}^n\times i\mathbb{R}$ as a
pseudo-sphere; also Tenenblat and Terng in \cite{TT1} have
developed the B transformation for local immersions of
$\mathbf{H}^n(\mathbb{R})$ in $\mathbb{R}^{2n-1}$; more generally
Tenenblat, Terng, Uhlenbeck and collaborators for space forms in
space forms). Beginning with a seed sub-manifold the same
formalism of B transformation can be applied and we can replace
the GCRM equations of the B transform with a Ricatti equation.

Note that (\ref{eq:dx0}) is true for sub-manifolds in
$\mathbb{C}^{n+p}$: if $(x,dx)=(R,t)(x_0,dx_0),\
(R,t)\subset\\\mathbb{O}_{n+p}(\mathbb{C}^{n+p})
\ltimes\mathbb{C}^{n+p},\
a_{jkl}:=\pa_{u^j}x_0^TR^{-1}\pa_{u^k}R\pa_{u^l}x_0$, then
$a_{jkl}=a_{jlk},\ a_{jkl}=-a_{lkj}$ so
$a_{jlk}=-a_{klj}=-a_{kjl}=a_{ljk}=a_{lkj}=-a_{jkl}=-a_{jlk}$. Let
$N_0$ be orthonormal frame of $x_0,\ N:=RN_0,\ \mathbb{C}^{n+p}\ni
a=a^{\top}+a^{\bot}=(I-N_0N_0^T)a+N_0N_0^Ta,\
R^{-1}dRa=(N_0N_0^TR^{-1}dR(I-N_0N_0^T)+R^{-1}dRN_0N_0^T)a=([N_0,R^{-1}dRN_0]-
N_0N_0^TR^{-1}dRN_0N_0^T)a$, where $[a,b]:=ba^T-ab^T$. The part
$[N_0,R^{-1}dRN_0]$ will encode the difference of the second
fundamental forms of $x,x_0:\
[N_0,R^{-1}dRN_0]dx_0=R^{-1}NN^Td^2x- N_0N_0^Td^2x_0$, while the
part $N_0N_0^TR^{-1}dRN_0N_0^T$ the difference of the normal
connections of $x,x_0:\
N_0N_0^TR^{-1}dRN_0=R^{-1}NN^TdN-N_0N_0^TdN_0$. The rotation of
the rolling can be changed so as to take $N_0$ to any other normal
frame $N_a$ of $x:\ N_a:=aN,\
a\subset\mathbb{O}_{n+p}(\mathbb{C}^{n+p}),\ adx=dx,\ R_a:=aR$.

We can now complete the only remaining part of the proof when
$n=2$: the complete integrability of the Ricatti equation; in
doing so the necessary conditions on the seed $x^0$ reveal
themselves. Let $N^0:=R_0N_0^0$ be the first vector of a normal
frame of $x^0$.

We have
$\om_0:=N_0^0\times(R_0^{-1}dR_0N_0^0)^{\bot_{\mathbb{C}^3}}$ and
we need only prove
$0=d\wedge\om_0+\frac{1}{2}\om_0\times\wedge\om_0
=\frac{1}{2}((R_0^{-1}dN^0)^{\bot_{\mathbb{C}^3}}
\times\wedge(R_0^{-1}dN^0)^{\bot_{\mathbb{C}^3}}
-dN_0^0\times\wedge dN_0^0)+N_0^0\times
d\wedge(R_0^{-1}dR_0N_0^0)^{\bot_{\mathbb{C}^3}}$. The first term
is normal, while the second one is tangent, so both must be $0$.
We have
$(R_0^{-1}dN^0)^{\bot_{\mathbb{C}^3}}=\\(dx_0^0)^TR_0^{-1}dN^0\llcorner\na
x_0^0=-(d^2x^0)^TN^0\llcorner\na x_0^0,\
dN_0^0=-(d^2x_0^0)^TN_0^0\llcorner\na x_0^0,\
-d\wedge(R_0^{-1}dR_0N_0^0)^{\bot_{\mathbb{C}^3}}
=d\wedge((N_0^0)^TR_0^{-1}dR_0(dx_0^0\llcorner\na
x_0^0))=d\wedge(((N^0)^Td^2x^0-(N_0^0)^Td^2x_0^0)\llcorner\na
x_0^0)$. But $0=d\wedge((N_0^0)^Td^2x_0^0\llcorner\na x_0^0)$ are
the Codazzi-Mainardi equations of $x_0^0$, so the two terms being
$0$ is equivalent to the statement that $x^0$ is a sub-manifold of
a ruled $3$-dimensional flat $x^0+wR_0N_0^0,\ w\in \mathbb{C}$ in
$\mathbb{C}^{p+2}$. The B transformation thus obtained is more
general than the one for $p=1$ because the flats of the seed and
the leaf are in general distinct (there is the question if this B
transformation can be extended to B transformation of the
corresponding flats, just as the B transformation for quadrics in
$\mathbb{C}^3$ can be extended the the $B$ transformation of
triply conjugate systems). This result is similar to the B
transformation of singular $2$-dimensional quadrics, but is not in
the spirit of Tenenblat-Terng B transformation, so to get a higher
dimensional generalization of deformation of quadrics a different
point of view with $p=n-1$ must be established. We reproduce in
the next subsection the Tenenblat-Terng B transformation, in the
hope that it will provide us with insight as to what is the
correct setting.

\subsection{Tenenblat-Terng's B\"{a}cklund transformation of $\mathbf{H}^n(\mathbb{R})$
in $\mathbb{R}^{2n-1}$}
\noindent

\noindent We need the notion of two $(n-1)$-dimensional subspaces
$V_1,\ V_2$ of $(\mathbb{R}^{2n-2},\ <,>)$ being {\it isoclinic}
(having constant inclination $\si$). Define a symmetric bilinear
inner product $<,>'$ on $V_2$ by
$<v_1,v_2>':=<(v_1)^{\bot_{V_1}},(v_2)^{\bot_{V_1}}>$. There is a
self-adjoint operator $A$ on $V_2$ such that
$<v_1,v_2>'=<Av_1,v_2>$; the angles $\si_1,...,\si_{n-1}$ between
$V_1,\ V_2$ are defined such that $\cos^2(\si_j),\ j=1,...,n-1$
are the eigenvalues of $A$. Geometrically these are obtained as
follows: $\cos\si_1=<v_1,(v_1)^{\bot_{V_1}}>:=\max_{v\in V_2,\
|v|=1}<v,(v)^{\bot_{V_1}}>,\
\cos\si_2=<v_2,(v_2)^{\bot_{V_1}}>:=\max_{v\in V_2,\ |v|=1,\
<v,v_1>=0}<v,(v)^{\bot_{V_1}}>$ and so on. This accounts for
orthonormal bases $E^j=[E'^j\ E''^j],\ E'^j=[e^j_1...e^j_{n-1}],\
E''^j=[e^j_n...e^j_{2n-2}],\ j=1,2$ of $\mathbb{R}^{2n-2}$ such
that ${E'}^j$ is basis of $V_j,\ j=1,2,\ AE'^2=E'^2D'^2$ and
\begin{eqnarray}\label{eq:isoc}
E^2=E^1\begin{bmatrix}D'&-D''\\D''&D'\end{bmatrix},
\end{eqnarray}
where $D':=\mathrm{diag}[\cos\si_1...\cos\si_{n-1}],\
D'':=\mathrm{diag}[\sin\si_1...\sin\si_{n-1}]$. If
$\si_1=...=\si_{n-1}=:\si$, then $V_1,\ V_2$ have constant
inclination $\si$.

Consider $x$ an $n$-dimensional sub-manifold of $\mathbb{R}^N$ and
$X=[X'\ X''],\ X'=[e_1...e_n],\ X''=[e_{n+1}...e_N]$ a moving
orthonormal frame along $x$ such that $\{e_j\}_{j=1,...,n}$ are
tangent and $\\\{e_{\al}\}_{\al=n+1,...,N}$ normal. Then
$\om':=X'^TdX',\ \om'':=X'^TdX''=-(X''^TdX')^T,\
\ti\om:=X''^TdX''$ will respectively encode the Levi-Civita
connection on $x$, the second fundamental form of $x$ and the
normal connection of $x$. We have $\om'^T=-\om',\
\ti\om^T=-\ti\om,\ dx=X'\om,\ \om:=X'^Tdx$, so
$d\wedge\om=dX'^T\wedge X'\om=-\om'\wedge\om$. Since $0=X''^Tdx$
we get $0=d\wedge(X''^Tdx)=dX''^T\wedge X'\om=\om''^T\wedge\om$.
Using Cartan's lemma we find $\om''^T=H(\om)$, where $H$ is linear
symmetric (if $\om=\{\om^j\}_{j=1,...,n},\
\om''=\{{\om''}_{\al}^j\}_{j=1,...,n,\ \al=n+1,...,N},\ H=\{h_{\al
i}^j\}_{i,j=1,...,n,\ \al=n+1,...,N}$, then
${\om''}_{\al}^i=\sum_jh_{\al j}^i\om^j,\ h_{\al j}^i=h_{\al
i}^j$). If $\Omega':=d\wedge\om'+\om'\wedge\om',\
\ti\Omega:=d\wedge\ti\om+\ti\om\wedge\ti\om$ are the Gau\ss\ and
normal curvature forms, then imposing the compatibility condition
$d\wedge(X^TdX)+(X^TdX)\wedge(X^TdX)=0,\
X^TdX=\begin{bmatrix}\om'&\om''\\-\om''^T&\ti\om\end{bmatrix}$ we
obtain the GCMR equations:

$$\Omega'-\om''\wedge\om''^T=0,\ d\wedge\om''+\om'\wedge\om''+\om''\wedge\ti\om=0,\
\ti\Omega-\om''^T\wedge\om''=0.$$

The first fundamental form of $x$ is
$|dx|^2=|\om|^2:=\sum_j(\om^j)^2$ and the second one is
$X''X''^Td^2x=X''\om''^T\om=\sum_{i,\al}\om^i{\om''}_{\al}^ie_{\al}=\sum_{i,j,\al}h_{\al
j}^i\om^i\om^je_{\al}$. $x$ has CGC $-1$ iff
$\Omega'=-\om\wedge\om^T$; if it has flat normal bundle
$\ti\Omega=0$, then we can choose $X'$ formed by principal
directions $h_{\al j}^i=\del_j^ib_{\al}^i$ and $X''$ giving normal
connection $\ti\om=0$. We have the next theorem due to Cartan and
Moore:

{\it (Cartan-Moore)
If $x$ is an $n$-dimensional CGC $-1$ in $\mathbb{R}^{2n-1}$, then
it has flat normal bundle. Let $X''$ be a normal frame such that
$\ti\om=0$. Then $x$ admits a parametrization by its lines of
curvature (Tchebyshef coordinates):

$$|dx|^2=\sum_ia_i^2(du^i)^2,\ a_i\neq 0,\ \sum_ia_i^2=1,\ X''X''^Td^2x=
\sum_{i,\al}b_{\al}^ia_i^2(du^i)^2e_{\al}.$$

Moreover $\sum_i\pa_{u_i}$ is the unique asymptotic vector in the
first orthant}.

For example if we have $\{e_j\}_{j=1,...,2n-1}$ the standard basis
of $\mathbb{R}^{2n-1}$ and 'rotate' the tractrix
$0<t\rightarrow\int_0^t\sqrt{1-e^{-2u}}due_1+e^{-t}\sum_{i=2}^n\la^ie_{2i-2},\
\la^i>0,\ \sum_{i=2}^{n}(\la^i)^2=1$, then we obtain
$x=\int_0^{u^1}\sqrt{1-e^{-2u}}due_1+e^{-u^1}\sum_{i=2}^n\la^i(\cos\frac{u^i}{\la^i}e_{2i-2}
+\sin\frac{u^i}{\la^i}e_{2i-1})$ and
$X'=[\frac{\pa_{u^1}x}{|\pa_{u^1}x|}...\frac{\pa_{u^n}x}{|\pa_{u^n}x|}],\
{\om'}_i^j=0,\ i,j>1,\ {\om'}_1^i=-e^{-u^1}du^i,\ i>1,\ \om=[du^1\
e^{-u^1}du^2...e^{-u^1}du^n]^T$. One can easily check that
$\Omega'=-\om\wedge\om^T$, so $x$ has CGC $-1$. For
$u^1=\mathrm{const}$ we obtain Clifford tori (deformations of
$\mathbb{R}^{n-1}$ in $\mathbb{S}^{2n-3}$), so the correct
generalization of the $2$-dimensional real pseudo-sphere in
$\mathbb{R}^3$ is to replace the homothetic parallels with
homothetic Clifford tori (under the same homothety). These flats
are parallel horospheres in $\mathbf{H}^n(\mathbb{R})$, so the
$n$-dimensional real pseudo-sphere is applicable to a region
bounded by a horosphere of radius half of the Poincar\'{e}'s unit
ball model of $\mathbf{H}^n(\mathbb{R})$; the orthogonal
trajectories of these flats are meridians and thus one can
establish the applicability correspondence at the level of linear
elements without need of further computing the curvature $\Omega'$
(the parallel horospheres on the pseudo-sphere
$\subset\mathbb{R}^n\times(i\mathbb{R})$ are intersections of this
with the pencil of $n$-dimensional isotropic spaces tangent to
$C(\infty)$ at one point). To get the parametrization as in the
theorem further apply the change $u^1=\ln(\cosh(v^1)),\ v^1>0$;
the normal frame is
$\frac{1}{\sqrt{1-e^{-2u^1}(1-(\la^i)^2)}}\{e^{-u^1}\la^ie_1
+\sqrt{1-e^{-2u^1}}(\cos\frac{u^i}{\la^i}e_{2i-2}
+\sin\frac{u^i}{\la^i}e_{2i-1})\}_{i=2,...,n}$. Note that we never
essentially used the fact that the particular deformation of
$\mathbb{R}^{n-1}$ in $\mathbb{S}^{2n-3}$ is a Clifford torus, so
from any such deformation we get one of $\mathbf{H}^n(\mathbb{R})$
in $\mathbb{R}^{2n-1}$.

Since $-\{\om^i\wedge\om^j\}_{i,j=1,...,n}=-\om\wedge\om^T=\Omega'
=\om''\wedge\om''^T=\{\sum_{\al}{\om''}_{\al}^i\wedge
{\om''}_{\al}^j\}_{i,j=1,...,n}
=\{\sum_{\al}b_{\al}^ib_{\al}^j\om^i\wedge\om^j\}_{i,j=1,...,n}$
we get
\begin{eqnarray}\label{eq:1}
\sum_{\al}b_{\al}^ib_{\al}^j=-1,\ i\neq j.
\end{eqnarray}
Since $\sum_i\pa_{u^i}$ is an asymptotic vector we get
$\sum_{i,\al}b_{\al}^ia_i^2e_{\al}=0$, or
\begin{eqnarray}\label{eq:2}
\sum_ib_{\al}^ia_i^2=0\ \forall\ \al=n+1,...,2n-1.
\end{eqnarray}

Multiplying this by $b_{\al}^j$ and summing after $\al$ we get
$\sum_{i,\al}b_{\al}^jb_{\al}^ia_i^2=0$ using (\ref{eq:1}) we
finally get $(\sum_{\al}a_jb_{\al}^j)^2+a_j^2=1$. Thus the
$n$-dimensional sub-manifold $A:=\{a^i_j\}_{i,j=1,...,n},\
a^1_i:=a_i,\ i=1,...,n,\ a^i_j:=a_jb_{n+i-1}^j,\ i=2,...,n,\
j=1,...,n$ is orthogonal and has nonzero entries on the first row.
The GCMR equations for $x$ become the {\it generalized sine-Gordon
equation} (GSGE) for $A$. Consider on $x$ the frame $X'$ formed by
the principal vectors $v_i:=\frac{1}{a^1_i}\pa_{u^i}x$ and $X''$
for which $\ti\om=0$. Then $\om^i=a^1_idu^i,\
{\om'}_j^i=\frac{\pa_{u^j}a^1_i}{a^1_j}du^i-\frac{\pa_{u^i}a^1_j}{a^1_i}du^j,\
{\om''}_{n+j-1}^i=a^j_idu^i,\ i=1,...,n,\ j=2,...,n$; with
$\del:=\mathrm{diag}[du^1...du^n]$ we have $\om=\del A^Te_1$ (here
$e_1\in\mathbb{R}^n$ is the first vector of the standard basis and
in the future it may be in $\mathbb{R}^m$ for different values of
$m$, as suitable) and the GSGE becomes the completely integrable
\begin{eqnarray}\label{eq:3}
d\wedge\om'+\om'\wedge\om'=-\del A^Te_1\wedge(\del A^Te_1)^T,\
\del\wedge\om'+A^TdA\wedge\del=0.
\end{eqnarray}
Conversely, if the orthogonal $n$-dimensional sub-manifold $A$
satisfies the GSGE and has nonzero entries on the first row, then
we know the first and second fundamental forms of a CGC $-1$
$n$-dimensional sub-manifold in the lines of curvature
parametrization and we can recover it by the integration of a
Ricatti equation and a quadrature. Note however that the GSGE
cannot apply to the pseudo-sphere itself since the shape is
imaginary degenerate.

Note that in order to get Terng's orthogonal $A$ of the GSGE we
only used the shape (extrinsic) part of the Gau\ss\ equation; the
remaining intrinsic part of the Gau\ss\ equation and the
Codazzi-Mainardi-Ricci equations followed immediately and did not
further impose functional restrictions on $A$; in fact if one
d-wedges the intrinsic Gau\ss\ equation (the first equation of
(\ref{eq:3})), then one naturally gets the second equation of
(\ref{eq:3}) (namely the Codazzi-Mainardi equations) and
conversely. In fact this is the method employed by Berger, Bryant
and Griffiths \cite{BB} by using only the Gau\ss\ equation and its
consequences obtained by d-wedging (because of its algebraic
character one does mostly Algebraic Geometry); once the
obstructions imposed by the Gau\ss\ equation and its higher
d-wedges are satisfied, everything else falls into place. If one
makes a detailed analysis of (\ref{eq:3}) and \`{a} la
Cartan-K\"{a}hler, then one gets the Berger, Bryant and Griffiths
\cite{BB} dependence of deformations of $\mathbf{H}^n(\mathbb{R})$
in $\mathbb{R}^{2n-1}$ on $n(n-1)$ functions of one variable. Note
that Cartan's argument according to which the co-dimension $n-1$
cannot be lowered follows from linear algebra and his exteriorly
orthogonal symmetric $2$-forms: basically from the extrinsic part
of the Gau\ss\ equation we have
$\sum_{\al=1}^rH_\al[e_i,e_j]H_\al+[e_i,e_j]=0$, where
$\{e_i\}_{i=1,...,n}$ is the standard basis of $\mathbb{R}^n,\
[e_i,e_j]:=e_je_i^T-e_ie_j^T$ and the symmetric shape operator is
$\{H_\al\}_{\al=1,...,r}=\{h_{\al
i}^j\}_{i,j=1,...,n,\al=1,...,r}$, so $H_1,...,H_r,I_n$ are
exteriorly orthogonal symmetric matrices. This exteriorly
orthogonal property remains invariant under linear transformations
of $\mathbb{R}^n$ (the exteriorly part) and under orthogonal
transformations of $\mathbb{R}^{r+1}$ (the orthogonal part).
Cartan proves in (\cite{C2}, ch II, \S\ 20) that if $r\le n-1$ and
one completes $H_1,...,H_r,I_n$ with $0$'s up to $n$ exteriorly
orthogonal symmetric $2$-forms, then the resulting forms are
linearly independent (so $r=n-1$) and can be found from an
orthogonal substitution of $2$-forms which are perfect squares:
this is exactly the above theorem if the squares are
$(\sqrt{|a_i|}du^i)^2$ and Terng's completely integrable GSGE
shows that the infinitesimal picture (purely algebraic in
character) integrates to the finite picture with no other
restrictions.

In the construction of the B transformation and the BPT two types
of frames will be useful: principal vectors completed with a
normal frame with $0$ normal connection and frames in which the
normal spaces of $x_0$ and $x_1:=B_{\si}(x_0)$ clearly appear
isoclinic. The change of tangent frames of $x_0$ will be an
orthogonal $n$-dimensional sub-manifold. An important observation
is that this is $A_1$; by symmetry the change of tangent frames of
$x_1$ will be $JA_0,\ J:=\mathrm{diag}[-1\ 1...1]$.

The orthogonal $A_1$ (and thus the B transform
$x_1=B_{\si}(x_0)$) will be given by the integration of the
Ricatti equation:
\begin{eqnarray}\label{eq:4}
dA_1=A_1\om'_0+A_1\del A_0^TDA_1-DA_0\del,
\end{eqnarray}
where $D:=\mathrm{diag}[\csc\si\ \cot\si...\cot\si]$. Imposing the
compatibility condition $d\wedge$ on (\ref{eq:4}) and using the
equation itself we thus need: $0=(A_1\om'_0+A_1\del
A_0^TDA_1-DA_0\del)\wedge\om'_0+A_1d\wedge\om'_0+(A_1\om'_0+A_1\del
A_0^TDA_1-DA_0\del)\wedge\del A_0^TDA_1-A_1\del\wedge
dA_0^TDA_1-A_1\del A_0^TD\wedge(A_1\om'_0+A_1\del
A_0^TDA_1-DA_0\del)-DdA_0\wedge\del=A_1(\om'_0\wedge\om'_0
+d\wedge\om'_0+(DA_0\del)^T\wedge(DA_0\del))-
DA_0(\del\wedge\om'_0+A_0^TdA_0\wedge\del)+A_1(\om'_0\wedge\del+\del\wedge
A_0^TdA_0)A_0^TDA_1$. But because $A_0$ satisfies GSGE, this
becomes $\del A_0^Te_1\wedge
e_1^TA_0\del=(DA_0\del)^T\wedge(DA_0\del)$, which is
straightforward (use $\csc^2\si-\cot^2\si=1$ and $D^2=\cot^2\si
I_n+e_1e_1^T$). Therefore (\ref{eq:4}) is completely integrable
and admits solution for any initial value of $A_1$. If the initial
value is orthogonal, we would like the solution to remain
orthogonal, or $0=dA_1^TA_1+A_1^TdA_1=(-\om'_0A_1^T+A_1^TDA_0\del
A_1^T-\del A_0^TD)A_1+A_1^T(A_1\om'_0+A_1\del A_0^TDA_1-DA_0\del)$
which is straightforward.

It is useful for the simplicity of the computations to add to the
normal frames with $n-1$ columns the zero column vector as the
first column and still call the extended 'frames' thus obtained
frames; all computations remain valid; for example $\om_0'':=\del
A_0^T\frac{J+I}{2}$ is augmented by a first column $0\
(\om''_{0in}:=0)$, which does not change the GCMR equations or any
other equation. Putting together
$d\wedge\om_0=-\om_0'\wedge\om_0,\
d\wedge\om_0''+\om_0'\wedge\om_0''=0$, or $-\del\wedge
dA_0^Te_1=-\om_0'\wedge\del A_0^Te_1,\ -\del\wedge
dA_0^T\frac{J+I}{2}+\om_0'\wedge\del A_0^T\frac{J+I}{2}=0$ we get
the second equation of (\ref{eq:3}). With
$i(A):=\begin{bmatrix}A&0\\0&I\end{bmatrix}$ the change of tangent
frames extends to frames as $[X_0'A_1^T\ \ X_0'']=X_0i(A_1^T)$.
Since the normal spaces of $x_0,\ x_1$ in the frames
$X_0i(A_1^T),\ X_1i(A_0^TJ)$ must have equal inclination $\si$ and
$x_1:=x_0+\sin\si X_0'A_1^Te_1$ is the formula for the B
transform, we must extend (\ref{eq:isoc}) by inserting a row and a
column of zeros in the middle and a block diagonal matrix of
dimension $1$ to the upper left corner and having entry $-1$. Thus
\begin{eqnarray}\label{eq:5}
X_1i(A_0^TJ)=X_0i(A_1^T)R_{\si},
\end{eqnarray}
where $R_{\si}:=\begin{bmatrix}-1&0&|&0&0\\0&\cos\si
I_{n-1}&|&0&\sin\si
I_{n-1}\\\_\_&\_\_&\_\_&\_\_&\_\_\\0&0&|&0_{1,1}&0\\0&-\sin\si
I_{n-1}&|&0&\cos\si I_{n-1}\end{bmatrix}=\sin\si\begin{bmatrix}
DJ&\frac{J+I}{2}\\-\frac{J+I}{2}&D\frac{J+I}{2}\end{bmatrix}$
(multiplication on the right (left) by $\frac{J+I}{2}$
(respectively $\frac{I-J}{2}$) replaces the first column (row)
with $0$ or respectively the remaining columns (rows) with $0$).
Therefore $X_1'=\sin\si(X_0'A_1^TD-X_0''\frac{J+I}{2})A_0,\
X_1''=\sin\si(X_0'A_1^T+X_0''D)\frac{J+I}{2}$. Since
$\om_0'={X_0'}^TdX_0'$, (\ref{eq:4}) can be written as
$A_1{X_0'}^Td(X_0'A_1^T)=DA_0\del A_1^T-A_1\del A_0^TD$. We have
$dx_1=X_0'\om_0+\sin\si d(X_0'A_1^T)e_1=(X_0'\del A_0^T+\sin\si
d(X_0'A_1^T))e_1$, so ${X_1''}^Tdx_1=\sin\si\frac{J+I}{2}(A_1\del
A_0^T+\sin\si A_1{X_0'}^Td(X_0'A_1^T)-\sin\si
D{\om_0''}^TA_1^T)e_1=\sin\si\frac{J+I}{2}(A_1\del
A_0^T+\sin\si(DA_0\del A_1^T-A_1\del A_0^TD)-
\sin\si\frac{J+I}{2}DA_0\del A_1^T)e_1=0$, since $\sin\si
De_1=e_1$. Also ${X_1'}^Tdx_1=\\=\sin\si
A_0^T(\sin\si\frac{J+I}{2}{\om_0''}^TA_1^T+DA_1\del A_0^T+\sin\si
DA_1{X_0'}^Td(X_0'A_1^T))e_1=\sin\si
A_0^T(\sin\si\frac{J+I}{2}A_0\del A_1^T+DA_1\del A_0^T+\sin\si
D(DA_0\del A_1^T-A_1\del A_0^TD))e_1=\sin^2\si
A_0^T(\frac{J+I}{2}+D^2)A_0\del A_1^Te_1=\del A_1^Te_1$, so
$\om_1:={X_1'}^Tdx_1=\del A_1^Te_1$ and thus $x_1$ is
$n$-dimensional if $A_1$ has nonzero entries on the first row. We
have ${X_1'}^TdX_1''=\sin^2\si
A_0^T(\frac{J+I}{2}{\om_0''}^TA_1^T+DA_1{X_0'}^Td(X_0'A_1^T)+
DA_1\om_0''D)\frac{J+I}{2}=\sin^2\si A_0^T(\frac{J+I}{2}A_0\del
A_1^T+D(DA_0\del A_1^T-A_1\del A_0^TD)+DA_1\del
A_0^TD)\frac{J+I}{2}=\del A_1^T\frac{J+I}{2}$, so
$\om_1'':={X_1'}^TdX_1''=\del A_1^T\frac{J+I}{2}$. Therefore $x_1$
has CGC $-1$:
$\Omega_1'=\om_1''\wedge{\om_1''}^T=-\om_1\wedge\om_1^T$ and zero
normal connection:
${X_1''}^TdX_1''=\frac{J+I}{2}\sin^2\si(A_1{X_0'}^Td(X_0'A_1^T)
+A_1\om_0''D-D{\om_0''}^TA_1)\frac{J+I}{2}
=\frac{J+I}{2}\sin^2\si(DA_0\del A_1^T-A_1\del A_0^TD+A_1\del
A_0^TD-DA_0\del A_1^T)\frac{J+I}{2}=0$. $A_1$ satisfies GSGE, all
computations admit the symmetry $(0,\si)\leftrightarrow (1,-\si)$
and lines of curvature (asymptotic directions) correspond on
$x_0,\ x_1$.

Note that from the isoclinic spaces property we always assumed
that all angles $\si$ are positive; in fact some of them can be
negative, so all computations involving diagonal matrices with
entries $\sin\si$ are valid if one attaches a sign to these
entries (which may differ from entry to entry). The signs can be
absorbed in the rows of $A_1$; since these changes do not change
the linear element and the shape of $x_1$, they do not change
$x_1$. For example in equation (\ref{eq:4}) the signs of the
entries in the diagonal $D$ are transferred by multiplication on
the left to the rows of $A_1$ and in $x_1=x_0+\sin\si
X'_0A_1^Te_1$ one can multiply the rows of $A_1$ with same signs
and keep $\sin\si>0$ for all frames in the isoclinic picture.

\subsubsection{Bianchi Permutability Theorem}

Let $x_i=B_{\si_i}(x_0),\ i=1,2$. We need to find just by
algebraic computations $B_{\si_2}(x_1)=x_3=B_{\si_1}(x_2)$. We
assume the theorem to be true, derive an algebraic formula and
post-priori show that the formula works. From (\ref{eq:5}) we thus
have:
$X_2i(A_0^T)R_{-\si_2}i(JA_2)i(A_1^T)R_{\si_1}i(JA_0)=X_0i(A_1^T)R_{\si_1}i(JA_0)=X_1=
X_3i(A_1^T)R_{-\si_2}i(JA_3)=\\X_2i(A_3^T)R_{\si_1}i(JA_2)i(A_1^T)R_{-\si_2}i(JA_3)$,
or with $X:=A_3A_0^T,\ C:=A_2A_1^T:\
i(X)R_{-\si_2}i(JC)R_{\si_1}i(J)=R_{\si_1}i(JC)R_{-\si_2}i(JX)$,
or $X(D_2CD_1+\frac{J+I}{2})=(D_1CD_2+\frac{J+I}{2})X,\ (X(D_1
-D_2C)-(D_1C-D_2))\frac{J+I}{2}=\frac{J+I}{2}((D_1-CD_2)X-(CD_1-D_2))=0$.
Since
$\sin\si_1(X_0'A_1^T+\sin\si_2(X_0'A_1^TD_1-X_0''\frac{J+I}{2})A_0A_3^T)e_1=
x_1-x_0+\sin\si_2X_1'A_3^Te_1=x_3-x_0=x_2-x_0+\sin\si_1X_2'A_3^Te_1
=\sin\si_2(X_0'A_2^T+\sin\si_1(X_0'A_2^TD_2-X_0''\frac{J+I}{2})A_0A_3^T)e_1$,
we get $X_0'A_1^T(D_2+D_1X^T-C^T(D_1+D_2X^T))e_1=0$, or
$e_1^T(X(D_1-D_2C)-(D_1C-D_2))=0$. If we denote
$a:=X(D_1-D_2C)-(D_1C-D_2)$, then we have
$e_1^Ta=a\frac{J+I}{2}=0,\ a=a\frac{I-J}{2}$. But
$0=(X(D_2CD_1+\frac{J+I}{2})-(D_1CD_2+\frac{J+I}{2})X)(D_1-D_2C)=(D_1C-D_2)D_2a+aD_2(D_1C-D_2)$,
or $a(\frac{I-J}{2}D_2(D_1C-D_2)\frac{J+I}{2})=0$, so $a=0$ and
\begin{eqnarray}\label{eq:permh}
X(D_1-D_2C)=D_1C-D_2.
\end{eqnarray}
Therefore we only need to prove that $A_3$ given by
(\ref{eq:permh}) satisfies (\ref{eq:4}) for $(A_0,\ \si)$ replaced
by $(A_1,\ \si_2),\ (A_2,\si_1)$; by symmetry it is enough to
prove only one relation. Since $dA_1=A_1\om'_0+A_1\del
A_0^TD_1A_1-D_1A_0\del,\ dA_2=A_2\om'_0+A_2\del
A_0^TD_2A_2-D_2A_0\del$, we get $d(A_2A_1^T)=\\(A_2\om'_0+A_2\del
A_0^TD_2A_2-D_2A_0\del)A_1^T+A_2(-\om'_0A_1^T+A_1^TD_1A_0\del
A_1^T-\del A_0^TD_1)=A_2(\del
A_0^T(D_2A_2-D_1A_1)-(A_2^TD_2-A_1^TD_1)A_0\del)A_1^T$. Thus if we
prove the similar relation $d(A_3A_0^T)=A_3(\del
A_1^T(D_2A_3+D_1A_0)-(A_3^TD_2+A_0^TD_1)A_1\del)A_0^T$, then since
$dA_0=A_0\om'_1-A_0\del A_1^TD_1A_0+D_1A_1\del$ we obtain what we
want: $dA_3=A_3\om'_1+A_3\del A_1^TD_2A_3-D_2A_1\del$.

But the formulae for $d(A_2A_1^T),\ d(A_3A_0^T)$ can be written
as: $dC=CA_1\del A_0^T(D_2C-D_1)-(D_2-CD_1)A_0\del A_1^T,\
dX=XA_0\del A_1^T(D_2X+D_1)-(D_2+XD_1)A_1\del A_0^T$ and
differentiating (\ref{eq:permh}) we get
$dX(D_1-D_2C)=(XD_2+D_1)dC$; thus we need to prove $(XA_0\del
A_1^T(D_2X+D_1)-(D_2+XD_1)A_1\del
A_0^T)(D_1-D_2C)=(XD_2+D_1)(CA_1\del
A_0^T(D_2C-D_1)-(D_2-CD_1)A_0\del A_1^T)$. The terms containing
$A_1\del A_0^T$ become: $(-(D_2+XD_1)+(XD_2+D_1)C)A_1\del
A_0^T(D_1-D_2C)$, which is $0$ because of (\ref{eq:permh}). Using
(\ref{eq:permh}) we have $(D_2X+D_1)(D_1-D_2C)=D_1^2-D_2^2,\
(XD_2+D_1)(D_2-CD_1)=X(D_2^2-D_1^2)$ and the terms containing
$A_0\del A_1^T$ cancel because
$D_2^2-D_1^2=(\cot^2\si_2-\cot^2\si_1)I$. Finally we need
$I=X^TX=((D_1-D_2C)^{-1})^T(D_1C-D_2)^T(D_1C-D_2)(D_1-D_2C)^{-1}$,
or $0=(D_1C-D_2)^T(D_1C-D_2)-(D_1-D_2C)^T(D_1-D_2C)
=C^T(D_1^2-D_2^2)C+D_2^2-D_1^2=
(\cot^2\si_2-\cot^2\si_1)(I-C^TC)$.

Again to prove the existence of the $\mathcal{M}_p$ configuration
we need only prove the existence of the $\mathcal{M}_3$
configuration (discrete deformations of $\mathbf{H}^n(\mathbb{R})$
will be obtained by considering $\mathbb{Z}^n$ lattices of B
transformations which are subsets of infinite M\"{o}bius
configurations).

Consider $(D_1D_2)^{-1}((D_1^2-D_2^2)
D_1A_1(D_1A_1-D_2A_2)^{-1}-D_1^2)=A_3A_0^{-1}=(D_1D_2)^{-1}((D_1^2-D_2^2)
D_2A_2(D_1A_1-D_2A_2)^{-1}-D_2^2),\
A_5A_0^{-1}=(D_1D_3)^{-1}((D_3^2-D_1^2)
D_1A_1(D_3A_4-D_1A_1)^{-1}-D_1^2),\
A_6A_0^{-1}=(D_2D_3)^{-1}((D_2^2-D_3^2)
D_2A_2(D_2A_2-D_3A_4)^{-1}-D_2^2)$; thus with
$\Box:=(D_2^2-D_3^2)D_1A_1+(D_3^2-D_1^2)D_2A_2+(D_1^2-D_2^2)D_3A_4$
we have $(D_2A_3A_0^{-1}-D_3A_5A_0^{-1})^{-1}A_1=((D_1^2-D_2^2)
(D_1A_1-D_2A_2)^{-1}-(D_3^2-D_1^2)(D_3A_4-D_1A_1)^{-1})^{-1}
=(D_1A_1-D_2A_2)\Box^{-1}(D_3A_4-D_1A_1)$ and similarly
$(D_3A_6A_0^{-1}-D_1A_3A_0^{-1})^{-1}A_2=(D_1A_1-D_2A_2)\Box^{-1}(D_2A_2-D_3A_4)$.
Now $D_1((D_2^2-D_3^2)D_2A_3A_0^{-1}
(D_2A_3A_0^{-1}-D_3A_5A_0^{-1})^{-1}A_1-D_2^2A_1)=
(D_1^2D_2A_2-D_2^2D_1A_1)\Box^{-1}(D_2^2-D_3^2)(D_3A_4-D_1A_1)-D_2^2D_1A_1
=(D_1^2D_2A_2-D_2^2D_1A_1)\Box^{-1}(D_3^2-D_1^2)(D_2A_2-D_3A_4)-D_1^2D_2A_2
=D_2((D_3^2-D_1^2)D_1A_3A_0^{-1}
(D_3A_6A_0^{-1}-D_1A_3A_0^{-1})^{-1}A_2-D_1^2A_2)$, so the very
lhs and rhs provide the good definition of and afford themselves
the name $D_1D_2D_3A_7$.

\subsection{The Bianchi-Lie approach}
\noindent

\noindent Following the Bianchi-Lie approach in generalizing the B
transformation of the $2$-dimensional pseudo-sphere to B
transformation of an arbitrary $2$-dimensional quadric, it is
natural to ask what are the leaves of the considered distribution
of isoclinic subspaces when the CGC $-1$ seed $x^0$ is in the
particular position of being the actual $n$-dimensional
pseudo-sphere. Once the correct metric-projective formulation is
found, all that is left to do is the analytic confirmation of its
generalization to arbitrary quadrics. If the CGC $-1$ seed is in
the particular position of actually being the pseudo-sphere, then
the facets of the considered distribution are situated on confocal
pseudo-spheres; this constitutes a strong indication if favor of
the Bianchi-Lie approach. However, the Tenenblat-Terng approach
has to be slightly modified, as for the $n$-dimensional
pseudo-sphere $|x_0|^2=-1$ in the Lorentz space
$\mathbb{R}^n\times i\mathbb{R}$ the shape is highly degenerated
and one cannot prescribe the orthogonal sub-manifold $A$.

Let $x_0:=-ie_0+\rho(ie_0+u),\ \rho:=\frac{2}{1-|u|^2},\
u:=\sum_{j=1}^nu^je_j\in\mathbb{C}^n$ be a parametrization of the
pseudo-sphere $x_0\subset \mathbb{C}^{n+1},\ |x_0|^2=-1$ with
asymptotic isotropic cone $|u|^2=1$. Note that $u$ may not vary in
all of $\mathbb{C}^n$; it may depend on $k$ complex and $n-k$ real
parameters such that $du$ is $n$-dimensional;
$\mathbf{H}^n(\mathbb{R})$ is obtained for $u\in\mathbb{R}^n,\
|u|^2<1$. We have $d\rho=\rho^2u^Tdu$; the unit normal is $ix_0$,
$dx_0=(ie_0+u)d\rho+\rho du,\
|dx_0|^2=\rho^2|du|^2=-ix_0^Td^2x_0$.

Again it is convenient to introduce a $0$ column vector as the
first one of an orthonormal frame in the normal bundle; all normal
frames will be of the form
$X''\in\mathbf{M}_{2n-1,n}(\mathbb{C}),\ X''e_1=0,\
X''^TX''=I_{2,n}$ and all computations with orthonormal frames
remain valid as long as one uses the transpose instead of inverse.

We have $|\pa_{u^j}x_0|=\rho,\ \pa_{u^j}x_0=\rho(e_j+\rho
u^j(ie_0+u))$,
$X'_0:=\rho^{-1}[\pa_{u^1}x_0...\pa_{u^n}x_0]\in\mathbf{M}_{2n-1,n}(\mathbb{C})$
moving tangent orthonormal frame along $x_0$ and $X''_0:=[0\ ix_0\
e_{n+1}...e_{2n-2}] \in\mathbf{M}_{2n-1,n}(\mathbb{C})$ moving
orthonormal frame in the normal bundle of $x_0$, $X_0:=[X'_0\
X''_0]$. Further $\om_0:=X'^T_0dx_0=\rho du,\
\om'_0:=X'^T_0dX'_0=\rho[u,du],\
\om''_0:=X'^T_0dX''_0=-(X''^T_0dX'_0)^T=\rho[0\ idu\ 0\ ...\ 0],\
\ti\om_0:=X''^T_0dX''_0=0,\ \Om'_0:=d\wedge
\om'_0+\om'_0\wedge\om'_0=-\rho^2du\wedge
du^T=-\om_0\wedge\om_0^T,\ d\wedge\om_0=-\om'_0\wedge\om_0,\
\om''^T_0\wedge\om_0=0,\ \om'^T_0=-\om'_0,\ \om''_0e_1=0$ and we
have the GCMR equations
$$\Om'_0-\om''_0\wedge\om''^T_0=0,\ d\wedge\om''_0
+\om'_0\wedge\om''_0=0,\ \om''^T_0\wedge\om''_0=0.$$ Consider the
rolling $(R,t)(x_0,dx_0)=(x,dx),\ x\subset\mathbb{C}^{2n-1}$ with
flat normal bundle; $R$ is chosen so that the normal connection of
$x$ is not only flat, but further trivial: $X''^T_0R^TdRX''_0=0$.

The corresponding omegas for $x$ have same definition and
properties as those of $x_0$, but with the index $0$ removed; the
frames of $x$ are those of $x_0$ multiplied on the left by $R$. We
have $\om=\om_0$; from (\ref{eq:dx0}) we have
$\om'=(RX'_0)^Td(RX'_0)=\om'_0$, so
$(RX_0)^TdRX_0=(RX_0)^Td(RX_0)-X_0^TdX_0
=\begin{bmatrix}0&\om''-\om''_0\\-(\om''-\om''_0)^T&0\end{bmatrix}$.

The distribution of $n$-dimensional facets centered at
$x_0+\sin\si X'_0A^Te_1,\ A\in\mathbf{O}_n(\mathbb{C})$, spanned
by frames $\sin\si(X'_0A^TDJ-X''_0)$ and with normal frames
$\sin\si(X'_0A^TI_{2,n}+X''_0D)$ should remain integrable if we
roll it as we roll $x_0$ on $x$ and for any position of $x$.

. Should such leaves exist, we would have
$0=(X'_0A^TI_{2,n}+X''_0D)^TR^Td(x+\sin\si RX'_0A^Te_1)=0$, or
$-D\om''^TA^Te_1+I_{2,n}A(\csc\si\om_0+\om'_0A^Te_1-A^TdAA^Te_1)=0$.
With $K:=\frac{1}{2}A\om''D,\
L:=\frac{1}{2}I_{2,n}A\om_0e_1^T=Le_1e_1^T$ we have
$Ke_1=L^Te_1=0,\ LI_{2,n}=0$ and the last relation becomes
\begin{eqnarray}\label{eq:ric}
dAA^T=A\om'_0A^T+K-K^T+\csc\si(L-L^T)+N,\nonumber\\
N^T=-N,\ Ne_1=(-K^T+\csc\si L)e_1.
\end{eqnarray}
From (\ref{eq:ric}) and the GW, GCMR equations of $x$ we get
\begin{eqnarray}\label{eq:KLM}
K^T\wedge K=L\wedge L=L^T\wedge L=L\wedge K^T=0,\ K^T\wedge A\om_0=0,\
K\wedge K^T=-\frac{\cot^2\si}{4}(A\om_0)\wedge(A\om_0)^T,\nonumber\\
d\wedge K=(K+\csc\si(L-L^T)+N)\wedge K,\
d\wedge L=\csc\si L\wedge A\om_0e_1^T+(K+K^T-e_1e_1^TK+N)\wedge L.\nonumber\\
\end{eqnarray}
Imposing the $d\wedge$ condition on (\ref{eq:ric})$A$, using
(\ref{eq:ric}) itself, (\ref{eq:KLM}) and the intrinsic part of
the Gau\ss\ equation we get
\begin{eqnarray}\label{eq:dM}
d\wedge N=-(K+\csc\si L+N)\wedge(K+\csc\si L+N)^T+(A\om_0)\wedge(A\om_0)^T+\nonumber\\
\csc\si((e_1e_1^TK)\wedge L+L^T\wedge(e_1e_1^T K)^T)
-\csc^2\si(L\wedge(A\om_0e_1^T)+(A\om_0e_1^T)^T\wedge L^T).
\end{eqnarray}
Further imposing the $d\wedge$ condition on (\ref{eq:dM}) and
using the equation itself, (\ref{eq:ric}) and (\ref{eq:KLM}) we
obtain only further (possibly vacuous) algebraic restrictions on
$N$. Now those conditions may make (\ref{eq:ric}) into a totally
integrable Ricatti equation in $A$, in which case the leaves in
the particular position of the Lie ansatz reveal themselves by
letting $\om'':=\om''_0$ and integrating that equation (hopefully
such a venue is possible).

\subsection{Particular deformations in $\mathbb{C}^{2n-1}$ of quadrics in
$\mathbb{C}^{n+1}$}
\noindent

\noindent

\newpage

\section{Appendix 1: A gist of deformations of quadrics}\label{sec:gist} \setcounter{equation}{0}
\subsection{Introduction}
\noindent

\noindent Quadratic laws and conserved quantities (like the energy
or the linear element of a surface) are metric-projective notions
(differentials involve tangency, a projective notion). It is thus
no accident that confocal quadrics (having a metric-projective
definition) behave well with respect to quadratic laws and
conserved quantities.

We shall present a partial outline of the theory of deformations
(through bending) of quadrics, namely the existence of the
B\"{a}cklund transformation for deformations of quadrics and the
applicability (isometric) correspondence provided by the Ivory
affinity. The analysis naturally splits into {\it static}
(algebraic) computations on quadrics confocal to the given one and
{\it moving} (differential) computations in its tangent bundle;
the flat connection form is provided by rolling. All necessary
identities of the moving part boil down to valid relevant
algebraic identities of the static part; in turn these algebraic
identities naturally appear within the static picture if one
applies the the method of Archimedes-Bianchi-Lie in the moving
picture.

\subsection{Confocal quadrics}
\noindent

\noindent Consider $\mathbb{R}^3$ with the standard basis
$e_1,e_2,e_3$ and the usual Euclidean scalar product
$<x,y>:=x^Ty,\ |x|^2:=x^Tx$. The family of quadrics $\{x_z\}_z$
confocal ('with same foci') to a given one
$x_0\subset\mathbb{R}^3$ is by definition the adjugate of the
pencil of quadrics generated by the Cayley's absolute
$\mathbb{CP}^2\supset C(\infty)\ni [y\ 0]:\
\begin{bmatrix}y\\0\end{bmatrix}^T\begin{bmatrix}I_3&0\\0&0\end{bmatrix}
\begin{bmatrix}y\\0\end{bmatrix}=0$ and the adjugate of the initial quadric
$x_0$:
$\begin{bmatrix}x_z\\1\end{bmatrix}^T(\begin{bmatrix}A&0\\0&-1\end{bmatrix}^{-1}
-z\begin{bmatrix}I_3&0\\0&0\end{bmatrix})^{-1}\begin{bmatrix}x_z\\1\end{bmatrix}=0$
for $A:=\mathrm{diag}[a_1^{-1}\ a_2^{-1}\ a_3^{-1}]$ (ellipsoids
and hyperboloids with one or two sheets) or
 $\begin{bmatrix}x_z\\1\end{bmatrix}^T(\begin{bmatrix}A&-e_3\\-e_3^T&0\end{bmatrix}^{-1}
-z\begin{bmatrix}I_3&0\\0&0\end{bmatrix})^{-1}\begin{bmatrix}x_z\\1\end{bmatrix}=0$
for $A:=\mathrm{diag}[a_1^{-1}\ a_2^{-1}\ 0]$ (elliptic or
hyperbolic paraboloids).

Since rigid motions and homotheties are the only projective
transformations (homographies) which preserve $C(\infty)$, it
encodes the metric structure of $\mathbb{R}^3$ and the mixed
metric-projective character of the definition of confocal quadrics
becomes clear.

Through any point $x\in\mathbb{R}^3$ pass $3$ quadrics of the
confocal family: $x\in x_{z_j},\ j=1,2,3$; replacing the cartesian
coordinates of $x$ with $(z_j)_{j=1,2,3}$ gives elliptic
coordinates on $\mathbb{R}^3$ (for quadrics of revolution we have
less than three $z$'s; we replace in this case one $z$ with the
coordinate induced by rotation or two $z$'s with the spherical
coordinates for confocal spheres).

The classical theorems about confocal quadrics are as follows:

{\it (Lam\'{e}) The normal to $x_z$ is proportional to $\pa_zx_z$
and the vectors $(\pa_z|_{z=z_j}x_{z_j})_{j=1,2,3}$ are
orthogonal; thus the family of quadrics confocal to a given
general quadric forms an orthogonal system.}

{\it (Dupin) Confocal quadrics cut each other along lines of
curvature (they are given by the elliptic coordinates).}

{\it (Ivory) The orthogonal trajectory of a point on $x_z$ as $z$
varies is a conic, and the correspondence $x_0\rightarrow x_z$
thus established is affine.

This affine transformation (henceforth called the Ivory affinity)
preserves the lengths of segments between confocal quadrics: with
$V_0^1:=x_z^1-x_0^0,\ V_1^0:=x_z^0-x_0^1$ we have
$|V_0^1|^2=|V_1^0|^2$ for pairs of points $(x_0^0,x_z^0),\
(x_0^1,x_z^1)$ corresponding on $(x_0,x_z)$ under the Ivory
affinity}.
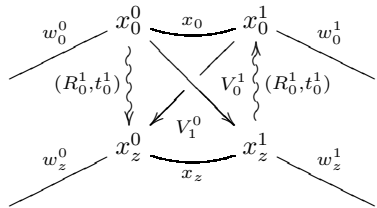
\begin{center}$\xymatrix@!0{&&x_0^0\ar@{-}[drdr]\ar@/_/@{-}[rr]^{x_0}\ar@{~>}[dd]_{(R_0^1,t_0^1)}&&
x_0^1\ar@{<~}[dd]^{(R_0^1,t_0^1)}&&\\
\ar@{-}[urr]^{w_0^0}&&&\ar[dl]^>>>>{V_1^0}\ar[dr]^>>>>>{V_0^1}&&&\ar@{-}[ull]_{w_0^1}\\&&
x_z^0\ar@{-}'[ur][urur]\ar@/_/@{-}[rr]_{x_z}&&x_z^1&\\\ar@{-}[urr]^{w_z^0}&&&
&&&\ar@{-}[ull]_{w_z^1}}$
\end{center}
{\it (Bianchi) If we have the rulings $w_0^0,\ w_0^1$ at the
points $x_0^0,\ x_0^1\in x_0$ and by use of the Ivory affinity we
get the rulings $w_z^0,\ w_z^1$ at the points $x_z^0,\ x_z^1\in
x_z$, then $[V_0^1\ \ w_0^0\ \ w_z^1]^T[V_0^1\ \ w_0^0\ \
w_z^1]=[-V_1^0\ \ w_z^0\ \ w_0^1]^T[-V_1^0\ \ w_z^0\ \ w_0^1]$, so
there exists a rigid motion
$(R_0^1,t_0^1)\in\mathbf{O}_3(\mathbb{R})\ltimes\mathbb{R}^3$ with
\begin{eqnarray}\label{eq:Qivo}
(R_0^1,t_0^1)(x_0^0,x_z^1,w_0^0,w_z^1)=(x_z^0,x_0^1,w_z^0,w_0^1).
\end{eqnarray}
Moreover $(V_0^1)^T\pa_z|_{z=0}x_z^0=(V_1^0)^T\pa_z|_{z=0}x_z^1$,
so the Ivory affinity has a nice projective property: the symmetry
of the tangency configuration
\begin{eqnarray}\label{eq:Qtc}
x_z^1\in T_{x_0^0}x_0\Leftrightarrow x_z^0\in T_{x_0^1}x_0.
\end{eqnarray}}

{\it (Jacobi) The tangent lines to a geodesic on $x_0$ remain
tangent to another confocal quadric.}

{\it (Chasles) The common tangents to two confocal quadrics form a
normal congruence and envelope geodesics on the two confocal
quadrics.}

There are many transformations of $\mathbb{R}^3$ which take a
family of confocal quadrics to another one, but a transformation
which preserves the previous properties must take lines to lines,
so it must be a homography; conversely:

{\it (Bianchi) Most homographies take a family of confocal
quadrics to another one.}

All quadrics are doubly ruled because they are equivalent from a
complex projective point of view, but the only doubly ruled real
quadrics are the hyperboloid with $1$ sheet
\begin{eqnarray}\label{eq:Qhyp}
x_z(u,v):=\sqrt{a_1-z}\frac{1-uv}{u-v}e_1+\sqrt{z-a_2}\frac{1+uv}{u-v}e_2+
\sqrt{a_3-z}\frac{u+v}{u-v}e_3,\ a_2<0,z<a_1,a_3,\nonumber\\
\end{eqnarray}
when the Ivory affinity is given by $x_z(u,v)=\sqrt{I_3-zA}\
x_0(u,v)$ and the hyperbolic paraboloid
\begin{eqnarray}\label{eq:Qpar}
x_z(u,v):=\sqrt{a_1-z}(u+v)e_1+\sqrt{z-a_2}(u-v)e_2+(2uv+\frac{z}{2})e_3,\
a_2<0,z<a_1,
\end{eqnarray}
when the Ivory affinity is given by $x_z(u,v)=\sqrt{I_3-zA}\
x_0(u,v)+\frac{z}{2}e_3$.

\subsection{The B\"{a}cklund transformation}
\noindent

\noindent B\"{a}cklund constructed in 1883 a transformation for
constant Gau\ss\ curvature $-1$ surfaces (this point of view is
due to Lie): the tangent spaces to the unit {\it pseudo-sphere}
$x_0$ (the space-like unit sphere in the Lorentz space
$\mathbb{R}^2\times i\mathbb{R}$; it is seen as time-like from the
origin) cut a confocal pseudo-sphere $x_z$ along circles, thus
highlighting a circle in each tangent space of $x_0$. Each point
of the circle, the segment joining it with the origin of the
tangent space and one of the (imaginary) rulings on $x_z$ passing
through that point determine a {\it facet} (pair of a point and a
plane passing through that point). We have thus highlighted a 3
dimensional integrable distribution of facets: its leaves are the
ruling families on $x_z$. Consider a surface $x$ (called {\it
seed}) non-rigidly isometric (applicable) to the pseudo-sphere
$x_0$. One can roll $x_0$ on $x$ such that for each position of
the rolled $x_0$ the rolled $x_0$ and $x$ will meet tangentially
and with the same differential at the tangency point.

If we roll the distribution while rolling $x_0$ on $x$ (the rolled
distribution is obtained as follows: each facet of the original
distribution corresponds to a point on $x_0$; we act on that facet
with the rigid motion of the rolling corresponding to the
highlighted point of $x_0$ in order to obtain the corresponding
facet of the rolled distribution), then it turns out that the
integrability condition of the rolled distribution is always
satisfied (we have complete integrability), so the integrability
of the rolled distribution does not depend on the shape of the
seed. Thus the rolled distribution is integrable and its leaves
(called the B\"{a}cklund transforms of $x$, denoted $B_z(x)$ and
whose determination requires the integration of a Ricatti
equation) will be applicable to the pseudo-sphere. Moreover the
seed and any leaf are the focal surfaces of a Weingarten
congruence (their second fundamental forms are proportional).

A-priori one would think that the B\"{a}cklund transformation for
the pseudo-sphere exists because the pseudo-sphere has a big group
of symmetries and thus a general quadric might not admit such a
transformation. But the fact that the B\"{a}cklund transformation
is of a general nature (independent of the shape of the seed)
means that its true nature lies not in the differential picture
with the seed and the leaf, but in the algebraic picture with
confocal pseudo-spheres where the B\"{a}cklund transformation
being independent of the shape of the seed is still valid
(although the leaves degenerate from surfaces to rulings of
$x_z$). Thus the vanishing of the shape of the seed $x$ in the
total integrability condition occurs not because of the symmetries
of the pseudo-sphere, but because of general equations (the
Gau\ss-Codazzi-Mainardi equations of the seed and of the
pseudo-sphere) coupled with algebraic consequences of the tangency
configuration and apparent at the level of confocal
pseudo-spheres: these are valid for a general quadrics.

Thus Bianchi generalized in 1906 the picture of the B\"{a}cklund
transformation of the pseudo-sphere to arbitrary quadrics: the
same statement remains true if in Lie's interpretation
'pseudo-sphere' is replaced with 'quadric' and 'circle' with
'conic'. However, the applicability correspondence on constant
Gau\ss\ curvature surfaces (which can be found in $\infty ^3$
ways) does not have an easy generalization to general quadrics.
According to Bianchi's own account (\cite{B2},(122),IV) should the
applicability correspondence exist in general, it can be presumed
to be of a general nature and thus independent of the shape of the
seed: therefore the answer lies not in the differential picture
with the seed and the leaf, but in the algebraic picture with
confocal quadrics so we have to roll back. Thus to find the
applicability correspondence one must look for a natural
correspondence between confocal quadrics. The Ivory affinity
provides such a correspondence and proving the validity of the
applicability correspondence via the Ivory affinity is very easy
once the complete integrability of the rolled distribution is
checked: if we roll the seed $x^0$ on $x_0^0$, then the tangent
space of the leaf $x^1$ corresponding to the point of tangency of
the rolled $x_0^0$ and the seed $x^0$ will be applied to the facet
centered at $x_z^1$ and spanned by $V_0^1,\ x_{zu_1}^1$; further
applying the rigid motion $(R_0^1,t_0^1)$ provided by the Ivory
affinity it will be applied to $T_{x_0^1}x_0$. In this process
$x^1_{u_1}$ is taken to $x_{zu_1}^1$ and further to $x_{0u_1}^1$,
so actually $(x^1,dx^1)$ is taken to $(x_0^1,dx_0^1)$; moreover
because of the symmetry of the tangency configuration the seed
becomes leaf and the leaf becomes seed.

Although the applicability correspondence is restricted, the
structure becomes richer (a quadric with less symmetries has less
degeneracies; for example a segment of B\"{a}cklund
transformations degenerates to the complementary transformation
for a quadric of revolution, or the geodesic flow on an ellipsoid
degenerates to the one on a sphere).

\subsection{Two algebraic consequences of the tangency
configuration}
\noindent

\noindent Consider now the tangency configuration
$(V_0^1)^TN_0^0=0$; let $\mathcal{B}:=(u-v)^2$ for (\ref{eq:Qhyp})
and $\mathcal{B}:=1$ for (\ref{eq:Qpar}),
$m_0^1:=\mathcal{B}_1x_{zu_1}^1\times V_0^1$ a normal field of the
distribution $\mathcal{D}^1$ of facets $\mathcal{F}^1$ passing
through $x_z^1$ and spanned by $V_0^1,\ x_{zu_1}^1$ (and similarly
$m'^1_0:=\mathcal{B}_1x_{zv_1}^1\times V_0^1$ by considering the
other ruling family on $x_z^1$); note that $m_0^1$ depends only on
$u_0,v_0,v_1$ and quadratically in $v_1$; this will make the
integrability condition a Ricatti equation in $v_1$.

Consider the rulings $w_0^0:=x_{0u_0}^0,\ w_0^1:=x_{0u_1}^1$ at
$x_0^0:=x_0(u_0,v_0),\ x_0^1:=x_0(u_1,v_1)$; we thus get a rigid
motion $(R_0^1,t_0^1)$ provided by the Ivory affinity. If we
change the ruling family on $x_0^1$ the action of the new rigid
motion on the facet $T_{x_0^0}x_0$ does not change, so its new
rotation must be the old rotation composed with a reflection in
$T_{x_0^0}x_0$, because of which the facets $\mathcal{F}^1,\
\mathcal{F}'^1$ reflect in $T_{x_0^0}x_0$ (thus the distributions
$\mathcal{D}^1,\ \mathcal{D}'^1$ reflect in $Tx_0$):
\begin{eqnarray}\label{eq:Qref}
(x_{zv_1}^1)^T(I_3-2N_0^0(N_0^0)^T)m_0^1=0,
\end{eqnarray}
and $x_{0v_1}^1=R_0^1(I_3-2N_0^0(N_0^0)^T)x_{zv_1}^1$; multiplying
this on the left by $(x_{0u_1}^1)^T$ and using the preservation of
lengths of rulings under the Ivory affinity we get
\begin{eqnarray}\label{eq:Qdir}
4(x_{zu_1}^1)^TN_0^0(N_0^0)^Tx_{zv_1}^1du_1dv_1=
|dx_z^1|^2-|dx_0^1|^2=-\frac{4z}{\mathcal{B}_1}du_1dv_1.
\end{eqnarray}
Thus we have the next result, essentially due to Bianchi (he uses
equivalent computations):

{\it If $x_z^1\in T_{x_0^0}x_0$, then:

I The change in the linear element from $x_z^1$ to $x_0^1$ is four
times the product of the orthogonal projections of the
differentials of the rulings of $x_z^1$ on the normal of $x_0$ at
$x_0^0$.

II The facets at $x_z^1$ spanned by $V_0^1$ and one of the rulings
of $x_z^1$ reflect in $T_{x_0^0}x_0$; therefore the distributions
$\mathcal{D}^1,\ \mathcal{D'}^1$ reflect in $Tx_0^0$.}

The algebraic relation
\begin{eqnarray}\label{eq:Qintalg}
(N_0^0)^T(2zm_0^1+m_0^1\times m_{0v_1}^1)=0
\end{eqnarray}
will appear as the total integrability condition. Using
(\ref{eq:Qref}), (\ref{eq:Qdir}) this becomes:
$0=\frac{z(m_0^1)^Tx_{zv_1}^1}{(N_0^0)^Tx_{zv_1}^1}-\mathcal{B}_1(x_{zu_1}^1)^TN_0^0(V_0^1)^T
m_{0v_1}^1=\frac{z(V_0^1)^T(\mathcal{B}_1x_{zv_1}^1\times
x_{zu_1}^1+(\mathcal{B}_1x_{zu_1}^1\times
V_0^1)_{v_1})}{(N_0^0)^Tx_{zv_1}^1}$, which is straightforward.
Replacing $(m_0^1,v_1)$ with $(m'^1_0,u_1)$ we get a similar
relation.

\subsection{Rolling surfaces and distributions}
\noindent

\noindent Let $(R_0,t_0)(x_0^0,dx_0^0)=(x^0,dx^0)$ be the rolling
of the quadric $x_0^0=x_0(u_0,v_0)$ on the applicable (rollable)
surface $x^0\subset\mathbb{R}^3$ (in general
$(R_0,t_0)\subset\mathbf{O}_3(\mathbb{R})\ltimes\mathbb{R}^3$ is a
surface, but it degenerates to a curve if $x^0$ is ruled); the
facets of the rolled distribution $(R_0,t_0)\mathcal{D}^1$ will
become tangent spaces to surfaces (leaves)
$x^1:=(R_0,t_0)x_z^1=(R_0,x^0)V_0^1$ iff the integrability
condition $0=(R_0m_0^1)^Tdx^1$ holds. Under the identification
$(\mathbf{o}_3(\mathbb{R}),[,])\simeq(\mathbb{R}^3,\times)$ we
have $R_0^{-1}dR_0\simeq N_0^0\times R_0^{-1}dR_0N_0^0=:\om_0$;
further imposing the compatibility condition $d\wedge$ to
$R_0dx_0^0=dx^0$ (we shall use the notation $d\wedge$ for exterior
derivative and $d$ for tensorial derivative) we get
\begin{eqnarray}\label{eq:Qom}
d\wedge\om_0+\frac{1}{2}\om_0\times\wedge\om_0=0,\
\om_0\times\wedge dx_0^0=0
\end{eqnarray}

and thus $\om_0$ is a flat connection form in $Tx_0^0$ (it encodes
the difference between the Gau\ss\ -Codazzi-Mainardi equations for
$x^0, x_0^0$).

We have
$R_0^{-1}dx_1=d(x_0^0+V_0^1)+R_0^{-1}dR_0V_0^1=dx_z^1+\om_0\times
V_0^1$. But $(\om_0)^{\bot}=0$ and
$dx_z^1=x_{zv_1}^1dv_1+x_{zu_1}^1du_1$, so the integrability
condition becomes $-(V_0^1)^T\om_0\times
N_0^0(m_0^1)^TN_0^0+(m_0^1)^Tx_{zv_1}^1dv_1=0$; using
(\ref{eq:Qref}) this becomes $-(V_0^1)^T\om_0\times
N_0^0+2(N_0^0)^Tx_{zv_1}^1dv_1=0$; multiplying it by
$\mathcal{B}_1(N_0^0)^Tx_{zu_1}^1$, using (\ref{eq:Qdir}) and
$-\mathcal{B}_1(N_0^0)^Tx_{zu_1}^1V_0^1=\mathcal{B}_1(V_0^1\times
x_{zu_1}^1)\times N_0^0=-m_0^1\times N_0^1$ we finally get the
Ricatti equation:
\begin{eqnarray}\label{eq:Qric}
(m_0^1)^T\om_0+2zdv_1=0.
\end{eqnarray}
We have $dm_0^1=m_{0v_1}^1dv_1+\mathcal{B}_1dx_0^0\times
x_{zu_1}^1$, so
$(dm_0^1)^T\wedge\om_0=dv_1\wedge(m_{0v_1})^T\om_0$; imposing the
total integrability condition $d\wedge$ on (\ref{eq:Qric}) and
using the equation itself we need
$-(m_0^1)^T\om_0\wedge(m_{0v_1}^1)^T\om_0+2z(m_0^1)^Td\wedge\om_0=0$,
or, using (\ref{eq:Qom}): $(N_0^0)^T(2zm_0^1+m_0^1\times
m_{0v_1}^1)(N_0^0)^T\om_0\times\wedge\om_0=0$; therefore the total
integrability is equivalent to (\ref{eq:Qintalg})*.\footnote{*
Most of \S\ 7.1-\S\ 7.5 was a 10 minutes presentation at a general
MAA session in San Antonio, TX January 12 2006. Now in August 2005
when I registered for this general session I did not see a
differential geometry AMS session; only algebraic geometry
sessions (it turned out I was wrong; on the flyer I have seen
several such differential geometry sessions). But the
simplification from late December 2005-early January 2006
fortunately made it appropriate for such a MAA session (in fact
while trying to present a simpler point of view for the MAA
session I managed to implement the simplest point of view at that
date); on top of that the link to Archimedes' method makes it even
more so, since Archimedes' ideas are currently subject to further
research in {\it The American Mathematical Monthly}. Thus the
apparently wrong denomination stuck and fulfilled its purpose.}

\subsection{ The method of Archimedes-Bianchi-Lie}
\noindent

\noindent If we roll $x_0^0$ on different sides of the seed $x^0$,
then we get the B\"{a}cklund transformation for the other ruling
family, so the rolled distributions reflect in the tangent bundle
of the seed $x^0$. Thus  (\ref{eq:Qref}) is obtained if one makes
the ansatz $x_0^0=x^0$; the same ansatz for the applicability
correspondence provided by the Ivory affinity and the inversion of
the B\"{a}cklund transformation (two focal surfaces of a line
congruence are in a symmetric relationship) implies Bianchi's
result about the existence of $(R_0^1,t_0^1)$ and the symmetry of
the tangency configuration; now (\ref{eq:Qdir}) is obtained as
previously described.

In the famous letter {\it The Method}, sent to Eratosthenes of the
University of Alexandria and lost from the XIII$^{\mathrm{th}}$
century until 1906 (the year of Bianchi's discovery, so it was
unknown to Bianchi and Lie), Archimedes states: '{\it ... certain
things first became clear to me by a mechanical method, although
they had to be proved by geometry afterwards because their
investigation by the said method did not furnish an actual proof.
But it is of course easier, when we have previously acquired, by
the method, some knowledge of the questions, to supply the proof
than it is to find it without any previous knowledge}'.

Note that the method at the level of points of facets was known to
Bianchi and Lie; however, they had never used the full method, at
the level of the planes of the facets too (for this reason the
{\it 'of course easier'} ingredient is missing from Bianchi's
proofs). Thus if one assumes Theorem I of Bianchi's theory of
deformations of quadrics a-priori to be true and to be the
metric-projective generalization of Lie's approach, then one
naturally {\it geometrically} gets the necessary algebraic
identities needed to prove Theorem I.

For surfaces this method is just a fancy way of reformulating
already known identities and which appear naturally enough at the
analytic level. But keep in mind that (\ref{eq:Qtc}),
(\ref{eq:Qref}) and (\ref{eq:Qdir}) are equivalent from an
analytic point of view and this is not the case in higher
dimensions: thus it is very difficult to find the necessary
algebraic identities of the static picture from an analytic point
of view. Therefore the method of Archimedes-Bianchi-Lie, due to
its geometric naturalness is useful in the study of higher
dimensional problems. Note that although Archimedes and
Bianchi-Lie have dealt with different problems, their approach was
the same: $P\Rightarrow Q$ with $P$ being either {\it 'The area of
a segment of a parabola is an infinite sum of areas of lines'} or
{\it 'A line is a deformation of a quadric'}. Such sentences $P$
were not accepted as true according to the standard of proof of
the times, but they had valid relevant consequences $Q$ which
elegantly solved problems not solvable with other methods of those
times; if only for this fact the standard of proof should change
so as to accept these statements as true. Note that {\it The
Method} of Archimedes went a little closer to the B
transformation: in an a-priori intuitive elementary non-rigorous
geometric argument he transferred and stuck all lines (slices of
the segment of the parabola) with their centers at the left end of
the balance (thus with $\infty$ multiplicity).

Note that the facets of the $3$-dimensional rolled distribution
are differently re-distributed into $2$-dimensional families of
facets as tangent planes to leaves when the shape of the seed
changes, but principles and properties independent of the shape of
the seed remain valid for facets even in the singular picture:
this is Archimedes' contribution to the Bianchi-Lie ansatz.

Thus I can honestly say that the {\it B transformation (for
quadrics)} is {\it one and the same with} the {\it Archimedes
balance} as {\it principles of a general nature}: they are valid
at the infinitesimal level and survive integration and conversely,
being principles of a general nature both induce by
differentiation and by particular singular configurations the
infinitesimal picture where the simplest explanation of these
principles reveals itself.

The natural question thus appears wether Google concurs with this
finding or not and indeed Google points to
{\href{http://slashwrestling.com/guests/llakor16.html} {the site}}
where the truth ({\it 'The truth is out there!'}) is spelled out:
{\it 'In the same way that Graham was a Man of Ego who looked like
a champion, Backlund was a Man of Skill who looked like something
completely different. Bob Backlund didn't even look like a
wrestler, he looked like your shop teacher from school, or Opie
Cunningham grown up. Backlund's motto was that 'any move can be
countered, any hold can be broken.' Like Archimedes, Backlund
declared, 'Give me where to stand, and I will move the earth,'
making of his own body, his lever, and of the wrestling ring, his
fulcrum.'}. All B\"{a}cklund needed was a seed and he produced
lots of leafs by means of his method. Note also that on
{\href{http://en.wikipedia.org/wiki/Geometer} {this site}} Albert
Victor B\"{a}cklund clearly appears as {\it an immediate follower}
of Archimedes; note also that B\"{a}cklund got his undergraduate
education at Lund, then travelled (including in Europe, where he
met with his conational Lie's ideas about the contact
transformation) and then he went {\it back to Lund} for the
remainder of his life. While Opie Cunningham was a friend of
Fonzie's, it seems to me that apparently Bob ain't no friend of
B\"{a}cklund's I ever heard of; only that Robert Miura is aka as
Bob (not necessarily English Bob) and he declared in (\cite{MM1},
pg 429) in 1976 as to regards to virtually the B transformation of
the KdV equation: {\it '.... However, several simplifications have
been achieved, not least of which is that all the equations to be
solved are linear!... It is instructive to point out two
simplifying features that make the method work. First, we were
able to show that the discrete eigenvalues are constant. Second,
we could compute the time dependence of $c_n$ and $b$ before
knowing the time-dependence of the solution $u$, requiring only
that we knew $u$ at $t=0$.... (We liken this situation to the
spirit of Archimedes' statement 'Give me but one firm spot on
which to stand and I will move the earth.')'}. I can honestly say
that Bob is actually pointing out Archimedes' method: first some
useful Euclidean identities are derived (corresponding to constant
eigenvalues), then we assume the method to work before we know the
method actually works and requiring in the process nothing more
than some partial initial knowledge of the method working (simple
consequences of the method) and {\it that is all there is to it};
the checking of the fact that the method works brings nothing more
besides tedious computations. Post-priori the partial knowledge
needed to prove the general configuration clearly appears from the
singular configuration by the application of the same method. Note
however that Bob (incarnated both as the Wrestling Bob ({\it 'His
name is Robert Paulson! His name is Robert Paulsen!'}) and as the
Mathematician Bob) was right about the balance: {\it 'That hit the
spot!'}. I can also honestly say that if the Mathematician Bob
were familiar with Archimedes' method (for example it being taught
in high-school), then he would have probably made the connection
in 1976 or earlier. Note also that Bianchi in (122,\S\ 10) was on
spot, as he uses angles and moments of pairs of lines to explain
the RMPIA; thus his moments provides the physical link to
Archimedes' moments of the balance.

Jules's {\it last minute} interpretation of {\it Ezekiel 12:57}
({\it 'What is significant is that I felt the touch of God. God
got involved. ... First, I'm going to deliver this case to
Marsellus. Then, basically I'm just going to walk The Earth. ...
Like Caine in 'Kung Fu'. Going from place to place, meeting new
people, getting in adventures. ... replaying the incident in my
head when I had what alcoholics refer to as a moment of clarity.
... Blessed is he who in the name of charity and good will
shepherds the weak through the Valley of Darkness, for he is truly
his brother's keeper and the finder of lost children. ... Or maybe
it is the world who is evil and selfish. ... I am the tyranny of
an evil man, but I am trying Ringo, I am trying really hard to be
the shepherd'}) reveals itself: nothing is more powerful than a
simple stick drawing a picture on the sand of a forgotten beach
only to be found by the survivors of a nearby unfortunate
shipwreck: {\it 'Bene speremus, Hominum enim vestigia video!'.

'Qui la Geometria vive!' Rock'n Roll The Rock $\&$ Peace out!}

\end{document}